\documentclass{article}
\usepackage{latexsym, amscd, amsfonts, eucal, mathrsfs, amsmath, amssymb, amsthm, xypic,xr, makeidx, stmaryrd}
\def\Z{\mathbf{Z}}
\def\C{\mathbf{C}}

\DeclareMathOperator{\Surj}{Surj} 
\DeclareMathOperator{\FamOp}{FOp} 
\newcommand{ \FamOpRed}{ \FamOp^{r}} 
\DeclareMathOperator{\Assem}{Assem} 
\DeclareMathOperator{\Emb}{Emb} 
\DeclareMathOperator{\Homeo}{Homeo} 
\DeclareMathOperator{\Germ}{Germ}
\DeclareMathOperator{\SO}{SO}
\DeclareMathOperator{\BSO}{BSO}
\DeclareMathOperator{\Smooth}{Sm}
\DeclareMathOperator{\acti}{act}

\DeclareMathOperator{\Mul}{Mul}

\DeclareMathOperator{\Spectra}{Sp}

\newcommand{\BigM}{{\mathfrak M}}

\DeclareMathOperator{\MMod}{\overline{{\calM}od}}
\DeclareMathOperator{\PrAlg}{ \overline{\Alg}}

\newcommand{\Biod}[3]{_{#1}{{\mathcal B}imod}_{#2}(#3)}

\DeclareMathOperator{\smooth}{sm}

\newcommand{\hyp}{\wedge}

\DeclareMathOperator{\Conf}{Conf}

\newcommand{\TopE}[1]{\widetilde{ \mathbb E}[#1]}

\DeclareMathOperator{\sd}{sd}
\DeclareMathOperator{\Rect}{Rect}

\DeclareMathOperator{\bfU}{{ \mathbf U}}

\newcommand{\Cube}[1]{{\Box}^{#1}}

\DeclareMathOperator{\Comm}{CRing}
\DeclareMathOperator{\CommOp}{Comm}
\DeclareMathOperator{\Top}{Top}

\DeclareMathOperator{\loc}{loc}

\newcommand{\caten}{ \star} 

\DeclareMathOperator{\Pres}{\mathcal{P}r}
\DeclareMathOperator{\LPrest}{L\mathcal{P}r}

\DeclareMathOperator{\LPress}{\LPrest}

\DeclareMathOperator{\Mor}{Mor}

\newcommand{\Sect}{\Gamma}
\DeclareMathOperator{\Fact}{Fact}
\newcommand{ \EMB}{ {\bf Man}} 
\newcommand{ \DISJ}[1]{ {\bf D}(#1)}
\DeclareMathOperator{\Ring}{Ring}

\DeclareMathOperator{\Disc}{Disc}

\DeclareMathOperator{\qunit}{qu}
\DeclareMathOperator{\CUP}{Disj}
\newcommand{\Disj}[1]{ \CUP(#1)}
\DeclareMathOperator{\Dis}{Disk}
\newcommand{\Disk}[1]{\Dis(#1)}
\DeclareMathOperator{\nunit}{nu} 
\DeclareMathOperator{\unit}{u}

\DeclareMathOperator{\CatMod}{ \Cat_{\infty}^{\Mod}}
\DeclareMathOperator{\LCatMod}{ \Pres^{\Mod}}
\DeclareMathOperator{\LCatModM}{ \Pres^{\Mod}_{\BigM/}}
\DeclareMathOperator{\lax}{lax}

\DeclareMathOperator{\FinSeg}{\Gamma}

\newcommand{ \CatAlg}{ \Cat_{\infty}^{\Alg}} 
\newcommand{ \hatCatAlg}{ \widehat{\Cat}_{\infty}^{\Alg}}
\newcommand{\LCatAlg}{ \Pres^{\Alg}}
\newcommand{ \hatCatMod}{ \widehat{\Cat}_{\infty}^{\Mod}}
\DeclareMathOperator{\LinF}{LinFun}
\newcommand{\LinFunc}[4]{\LinF^{#1}_{#2}( #3,#4)}

\DeclareMathOperator{\Envv}{Env}
\DeclareMathOperator{\PreOp}{{\mathcal P}Op_{\infty}}

\DeclareMathOperator{\aug}{aug}

\DeclareMathOperator{\glike}{gp}

\newcommand{\Andre}{Andr\'{e}}

\newcommand{\Kahler}{K\"{a}hler\,}
\DeclareMathOperator{\Der}{Der}

\DeclareMathOperator{\Free}{Free}
\DeclareMathOperator{\Baar}{Bar}

\newcommand{\OpE}[1]{{\mathbb E}[#1]}

\DeclareMathOperator{\OO}{O}
\DeclareMathOperator{\Isom}{Isom}
\DeclareMathOperator{\BTop}{BTop}

\newcommand{\seg}[1]{{\langle #1 \rangle}}
\newcommand{\Colp}[1]{\rho^{#1}}

\newcommand{\nostar}[1]{{{\langle #1 \rangle}^{\circ}}}

\DeclareMathOperator{\Sing}{Sing}
\newcommand{\Un}{Un}

\DeclareMathOperator{\Group}{\mathcal{G}rp}

\DeclareMathOperator{\calH}{\mathcal{H}}
\DeclareMathOperator{\Mon}{Mon}

\DeclareMathOperator{\Alg}{Alg}
\DeclareMathOperator{\CAlg}{CAlg}

\newcommand{\toposref}[1]{T.\ref{HTT-#1}}
\newcommand{\stableref}[1]{S.\ref{STA-#1}}
\newcommand{\monoidref}[1]{M.\ref{MON-#1}}
\newcommand{\symmetricref}[1]{C.\ref{SYM-#1}}
\newcommand{\deformationref}[1]{D.\ref{DEF-#1}}

\newcommand{\bicatref}[1]{B.\ref{BIC-#1}}

\newcommand{\degree}{\text{o}}
\newcommand{\bfA}{{\mathbf A}}

\newcommand{\wreath}{\wr}
\DeclareMathOperator{\Kan}{\mathcal{K}an}
\DeclareMathOperator{\cDelta}{{\bf \Delta}}

\DeclareMathOperator{\Set}{\mathcal{S}et}
\DeclareMathOperator{\Vect}{Vect}
\DeclareMathOperator{\sSet}{\mathcal{S}et_{\Delta}}
\DeclareMathOperator{\sCoNerve}{\mathfrak{C}}
\DeclareMathOperator{\Nerve}{N}
\DeclareMathOperator{\Subd}{{\mathcal S}ub}
\DeclareMathOperator{\Cat}{\mathcal{C}at}
\newcommand{\h}[1]{\rm{h} \! #1}

\newcommand{\Adjoint}[4]{\xymatrix@1{#2 \ar@<.4ex>[r]^-{#1} & #3 \ar@<.4ex>[l]^-{#4}}}

\renewcommand{\boxtimes}{\odot}

\DeclareMathOperator{\coker}{coker}
 
\newcommand{\etale}{{\'{e}tale}\,\,}

\newcommand{\et}{\'{e}t}

\newcommand{\bigdot}{\bullet}

\DeclareMathOperator{\Stab}{Stab}

\DeclareMathOperator{\bHom}{Map}
\DeclareMathOperator{\Sh}{Sh}

\DeclareMathOperator{\Shv}{{\mathcal S}hv}

\newcommand{\Cech}{\v{C}ech\,}

 \DeclareMathOperator{\End}{End}
 
\DeclareMathOperator{\Mod}{Mod}
\DeclareMathOperator{\CMod}{Mod}

\DeclareMathOperator{\calM}{\mathcal{M}}

\DeclareMathOperator{\GL}{GL} \DeclareMathOperator{\G}{\mathbf{G}}
\DeclareMathOperator{\calZ}{\mathcal{Z}}

\DeclareMathOperator{\Ext}{Ext} 
 
\DeclareMathOperator{\DerRing}{\mathcal{SCR}}

\DeclareMathOperator{\colim}{colim}

\DeclareMathOperator{\calA}{\mathcal{A}}
\DeclareMathOperator{\bd}{\partial}

\DeclareMathOperator{\calU}{\mathcal{U}}

\DeclareMathOperator{\calV}{\mathcal{V}}
\DeclareMathOperator{\calE}{\mathcal{E}}
\DeclareMathOperator{\R}{\mathbf{R}}
\DeclareMathOperator{\calO}{\mathcal{O}}
\DeclareMathOperator{\calT}{\mathcal{T}}

\DeclareMathOperator{\calB}{\mathcal{B}}
\DeclareMathOperator{\calK}{\mathcal{K}}
\DeclareMathOperator{\Spec}{{\bf Spec}}

\DeclareMathOperator{\Triv}{{\mathcal T}riv}
\DeclareMathOperator{\Ass}{{\mathcal A}ss}
\DeclareMathOperator{\CatAss}{{\bf Ass}}
\DeclareMathOperator{\calF}{\mathcal{F}}
\DeclareMathOperator{\calG}{\mathcal{G}}
\DeclareMathOperator{\Hom}{Hom} 
\DeclareMathOperator{\HH}{H} 
\DeclareMathOperator{\id}{id} \DeclareMathOperator{\Fun}{Fun}
\DeclareMathOperator{\calC}{\mathcal{C}}
\DeclareMathOperator{\calQ}{\mathcal{Q}}
\DeclareMathOperator{\calI}{\mathcal{I}}
\DeclareMathOperator{\calN}{\mathcal{N}}
\DeclareMathOperator{\calJ}{\mathcal{J}}
\DeclareMathOperator{\Ran}{Ran}

\DeclareMathOperator{\HC}{HC} 
\DeclareMathOperator{\calR}{\mathcal{R}}

\newcommand{\Cent}[1]{{\mathfrak Z}(#1)}
\newcommand{\Centt}[2]{{\mathfrak Z}_{#1}(#2)}
\DeclareMathOperator{\SSet}{\mathcal{S}}

\DeclareMathOperator{\Aut}{Aut}
\DeclareMathOperator{\MM}{{\mathfrak m}}
\DeclareMathOperator{\AAA}{{\mathfrak a}}
\DeclareMathOperator{\rk}{rk}
\DeclareMathOperator{\supp}{Supp}
\DeclareMathOperator{\CatLMod}{LMod}
\DeclareMathOperator{\LMod}{{\mathcal L}Mod}

\DeclareMathOperator{\calX}{\mathcal{X}}

\DeclareMathOperator{\hocolim}{hocolim}
\DeclareMathOperator{\calY}{\mathcal{Y}}

\DeclareMathOperator{\op}{op}
\DeclareMathOperator{\calD}{\mathcal{D}}
\DeclareMathOperator{\Ind}{Ind} 
\DeclareMathOperator{\calP}{\mathcal{P}} \topmargin=0in

\newtheorem{theorem}{Theorem}[subsection]
\newtheorem{lemma}[theorem]{Lemma}
\newtheorem{proposition}[theorem]{Proposition}

\newtheorem{corollary}[theorem]{Corollary}

\theoremstyle{definition}
\newtheorem{convention}[theorem]{Convention}
\newtheorem{definition}[theorem]{Definition}
\newtheorem{construction}[theorem]{Construction}
\newtheorem{example}[theorem]{Example}

\newtheorem{notation}[theorem]{Notation}

\newtheorem{warning}[theorem]{Warning}
\newtheorem{remark}[theorem]{Remark}
\newtheorem{variant}[theorem]{Variant}

\makeindex

\begin{document}

\title{Derived Algebraic Geometry VI: $\OpE{k}$-Algebras}

\maketitle
\tableofcontents

\section*{Introduction}

Let $X$ be a topological space equipped with a base point $\ast$. We let
$\Omega X$ denote the loop space of $X$, which we will identify with the space
of continuous map $p: [-1,1] \rightarrow X$ such that $f(-1)= \ast = f(1)$. 
Given a pair of loops $p, q \in \Omega X$, we can define a composite loop
$p \circ q$ by concatenating $p$ with $q$: that is, we define
$p \circ q: [-1,1] \rightarrow X$ by the formula
$$ (p \circ q)(t) = \begin{cases} q(2t+1) & \text{ if } -1 \leq t \leq 0 \\
p(2t-1) & \text{ if } 0 \leq t \leq 1.\end{cases}$$
This composition operation is associative up to homotopy, and endows
the set of path components $\pi_0 \Omega X$ with the structure of a group: namely, the fundamental group $\pi_1(X, \ast)$. However, composition of paths is not strictly associative:
given a triple of paths $p,q,r \in \Omega X$, we have
$$ (p \circ (q \circ r))(t) = \begin{cases} r(4t+3) & \text{ if } -1 \leq t \leq \frac{-1}{2} \\
q(4t + 1) & \text{ if } \frac{-1}{2} \leq t \leq 0 \\
p(2t-1) & \text{ if } 0 \leq t \leq 1. \end{cases} \quad \quad
((p \circ q) \circ r)(t) = \begin{cases} r(2t+1) & \text{ if } -1 \leq t \leq 0 \\
q(4t-1) & \text{ if } 0 \leq t \leq \frac{1}{2} \\
p(4t-3) & \text{ if } \frac{1}{2} \leq t \leq 1. \end{cases}$$
The paths $p \circ (q \circ r)$ and $(p \circ q) \circ r$ follow the same trajectories but are parametrized differently; they are homotopic but not identical.

One way to compensate for the failure of strict associativity is to consider not one composition operation but several. For every finite set $S$, let $\Rect( (-1,1) \times S, (-1,1) )$ denote the collection
of finite sequences of maps $\{ f_S: (-1,1) \rightarrow (-1,1) \}_{s \in S}$ with the following properties:
\begin{itemize}
\item[$(a)$] For $s \neq t$, the maps $f_{s}$ and $f_{t}$ have disjoint images.
\item[$(b)$] For each $s \in S$, the map $f_s$ is given by a formula
$f_{s}(t) = a t  + b$ where $a > 0$.
\end{itemize}
If $X$ is any pointed topological space, then there is an evident map
$$ \theta: (\Omega X)^{S} \times \Rect( (-1,1) \times S, (-1,1) ) \rightarrow \Omega X,$$
given by the formula
$$ \theta( \{ p_s \}_{s \in S}, \{ f_s \}_{s \in S})(t) = \begin{cases}
p_s(t') & \text{ if } t = f_s(t') \\
\ast & \text{ otherwise. }\end{cases}$$
Each of the spaces $\Rect( (-1,1) \times S, (-1,1))$ is equipped with a natural topology
(with respect to which the map $\theta$ is continuous), and the collection of spaces
$\{ \Rect( (-1,1) \times S, (-1,1) ) \}_{I}$ can be organized into a topological operad,
which we will denote by $\calC_1$. We can summarize the situation as follows:

\begin{itemize}
\item[$(\ast)$] For every pointed topological space $X$, the loop space $\Omega X$ carries an action of the topological operad $\calC_1$.
\end{itemize}

Every point of $\Rect( (-1,1) \times S, (-1,1))$ determines a linear ordering of the finite set $S$.
Conversely, if we fix a linear ordering of $S$, then the corresponding subspace of
$\Rect( (-1,1) \times S, (-1,1))$ is contractible. In other words, there is a canonical homotopy equivalence
of $\Rect( (-1,1) \times S, (-1,1))$ with the (discrete) set of linear orderings of $S$. Together, these homotopy equivalences determine a weak equivalence of the topological operad $\calC_1$ with the associative operad (see Example \ref{sulta}). Consequently, an action of the operad $\calC_1$ can be regarded as a homotopy-theoretic substitute for an associative algebra structure. In other words,
assertion $(\ast)$ articulates the idea that the loop space $\Omega X$ is equipped with a multiplication which associative up to coherent homotopy.

If $X$ is a pointed space, then we can consider also the $k$-fold loop space $\Omega^{k} X$, which
we will identify with the space of all maps $f: [-1,1]^{k} \rightarrow X$ which carry the boundary of the cube $[-1,1]^{k}$ to the base point of $X$. If $k > 0$, then we can identify $\Omega^{k} X$ with
$\Omega( \Omega^{k-1} X)$, so that $\Omega^{k} X$ is equipped with a coherently associative multiplication given by concatenation of loops. However, if $k > 1$, then the structure of
$\Omega^{k} X$ is much richer. To investigate this structure, it is convenient to introduce a higher-dimensional version of the topological $\calC_1$, called the {\it little $k$-cubes operad}. 
We begin by introducing a bit of terminology.

\begin{definition}\label{defcube}
Let $\Cube{k} = (-1,1)^{k}$ denote an open cube of dimension $k$. We will say that a map
$f: \Cube{k} \rightarrow \Cube{k}$ is a {\it rectilinear embedding} if it is given by the
formula 
$$ f( x_1, \ldots, x_k ) = ( a_1 x_1 + b_1, \ldots, a_k x_k + b_k )$$
for some real constants $a_i$ and $b_i$, with $a_i > 0$.
More generally, if $S$ is a finite set, then we will say that a map
$\Cube{k} \times S \rightarrow \Cube{k}$ is a {\it rectilinear embedding} if 
it is an open embedding whose restriction to each connected component of
$\Cube{k} \times S$ is rectilinear. Let $\Rect( \Cube{k} \times S, \Cube{k})$ denote the collection of
all rectitlinear embeddings from $\Cube{k} \times S$ into $\Cube{k}$. 
We will regard $\Rect( \Cube{k} \times S, \Cube{k})$ as a topological space
(it can be identified with an open subset of $(\R^{2k})^{I}$).
\end{definition}

The spaces $\Rect( \Cube{k} \times S, \Cube{k})$ determine a topological operad, called the
{\it little $k$-cubes operad}; we will (temporarily) denote this operad by $\calC_k$.
Assertion $(\ast)$ has an evident generalization:

\begin{itemize}
\item[$(\ast')$] Let $X$ be a pointed topological space. Then the $k$-fold loop space
$\Omega^{k} X$ carries an action of the topological operad $\calC_k$.
\end{itemize}

Assertion $(\ast')$ admits the following converse, which highlights the importance of
the little cubes operad $\calC_{k}$ in algebraic topology:

\begin{theorem}[May]\label{recog}
Let $Y$ be a topological space equipped with an action of the little cubes operad
$\calC_{k}$. Suppose that $Y$ is grouplike (Definition \ref{ungwar}). Then
$Y$ is weakly homotopy equivalent to $\Omega^{k} X$, for some pointed topological space
$X$.
\end{theorem}

A proof of Theorem \ref{recog} is given in \cite{iterated} (we will prove another version of this result as
Theorem \ref{slage}, using a similar argument). Theorem \ref{recog} can be interpreted as saying that, in some sense, the topological operad $\calC_{k}$ encodes precisely the structure that a $k$-fold loop space should be expected to possess. In the case $k=1$, the structure consists of a coherently associative multiplication. This structure turns out to be useful and interesting outside the context of topological spaces. For example,
we can consider associative algebras in the category of abelian groups (regarded as a symmetric monoidal category with respect to the tensor product of abelian groups) to recover the theory of
associative rings. 

A basic observation in the theory of structured ring spectra is that a similar phenomenon occurs for larger values of $k$: that is, it is interesting to study algebras over the topological operads $\calC_{k}$ in categories other than that of topological spaces. Our goal in this paper is to lay the general foundations for such a study,
using the formalism of $\infty$-operads developed in \cite{symmetric}. More precisely, we will define (for each integer $k \geq 0$) an $\infty$-operad $\OpE{k}$ of {\it little $k$-cubes} by applying the construction of Notation \symmetricref{capin} to the topological operad $\calC_{k}$ (see Definition \ref{defE}). We will study these operads (and algebras over them) in \S \ref{founder}.

In \cite{monoidal}, we develop a general theory of associative algebras and (right or left) modules over them. This theory has a number of applications to the study of $\OpE{k}$-algebras, which we will describe in \S \ref{sec2}. For example, we prove a version of the (generalized) Deligne conjecture in \S \ref{delconj}.

The description of the $\infty$-operads $\OpE{k}$ in terms of rectilinear embeddings between cubes suggests
that the the theory of $\OpE{k}$-algebras is closely connected with geometry. In \S \ref{sec3}, we will make
this connection more explicit by developing the theories of {\it factorizable (co)sheaves} and {\it topological chiral homology}, following ideas introduced by Beilinson and Drinfeld in the algebro-geometric context.

We conclude this paper with two appendices. The first, \S \ref{secA}, reviews a number of ideas from topology which are relevant to the subject of this paper. A large portion of this appendix, concerning the theory of constructible sheaves and $\infty$-categories of exit paths, is not really needed. However, it can be used to provide an alternative description of the $\infty$-category of constructible (co)sheaves
studied in \S \ref{stupem}. Our second appendix, \S \ref{secB}, develops several aspects of the general theory of $\infty$-operads which are used in the body of this paper, but are not covered in \cite{symmetric}.

\begin{remark}
We should emphasize that the little cubes operads and their algebras are not new objects of study. Consequently, many of the ideas treated in this paper have appeared elsewhere, though often in a rather different language. It should not be assumed that uncredited results are due to the author; we apologize in advance to any whose work does not receive the proper attribution.
\end{remark}

\subsection*{Notation and Terminology}

For an introduction to the language of higher category theory (from the point of view taken in this paper), we refer the reader to \cite{topoi}. For convenience, we will adopt the following conventions concerning references to \cite{topoi} and to the other papers in this series:

\begin{itemize}
\item[$(T)$] We will indicate references to \cite{topoi} using the letter T.
\item[$(S)$] We will indicate references to \cite{DAGStable} using the letter S.
\item[$(M)$] We will indicate references to \cite{monoidal} using the letter M.
\item[$(C)$] We will indicate references to \cite{symmetric} using the letter C.
\item[$(D)$] We will indicate references to \cite{deformation} using the letter D.
\item[$(B)$] We will indicate references to \cite{derivative} using the letter B.
\end{itemize}

We will use the notation of \cite{symmetric} throughout this paper. In particular, for every integer $n \geq 0$, we let $\nostar{n}$ denote the finite set $\{ 1, \ldots, n \}$, and $\seg{n} = \{ 1, \ldots, n, \ast \}$ the finite set obtained by adjoining a base point to $\nostar{n}$. We let $\FinSeg$ denote the category whose objects
are the finite sets $\seg{n}$, and whose morphisms are maps of finite sets $f: \seg{m} \rightarrow \seg{n}$ such that $f(\ast) = \ast$.

\subsection*{Acknowledgements}

I would like to thank John Francis, Dennis Gaitsgory, and Mike Hopkins for numerous conversations concerning the subject matter of this paper. I would also like to thank the National Science Foundation for supporting this research (via NSF grants DMS-0757293 and DMS-0943108).

\section{Foundations}\label{founder}

Our goal in this section is to introduce the $\infty$-operads $\{ \OpE{k} \}_{k \geq 0}$ of
{\it little $k$-cubes}, and to verify their basic properties. We begin in \S \ref{kloo} with a review of the relevant definitions. In \S \ref{sass1}, we will prove an important additivity result (Theorem \ref{slide}), which asserts that 
endowing an object $A$ of a symmetric monoidal $\infty$-category $\calC$ with the structure of an
$\OpE{k+k'}$-algebra is equivalent to equipping $A$ with $\OpE{k}$-algebra and $\OpE{k'}$-algebra structures, which are compatible with one another in a suitable sense. In particular, since $\OpE{1}$ is equivalent to the associative $\infty$-operad (Example \ref{sulta}), we can think of an $\OpE{k}$-algebra
$A \in \Alg_{ \OpE{k}}(\calC)$ as an object of $\calC$ which is equipped with $k$ compatible associative multiplications.

One of the original applications of the little cubes operads in topology is to the study of $k$-fold loop spaces. In \S \ref{kloop}, we will return to this setting by considering $\OpE{k}$-algebra objects in the $\infty$-category $\SSet$ of spaces. In particular, we will revisit a classical result of May, which establishes the equivalence between a suitable $\infty$-category of $\OpE{k}$-spaces and the $\infty$-category of $k$-fold loop spaces (Theorem \ref{slage}). 

Outside of topology, the little cubes operads have found a number of algebraic applications. For these purposes, it is important to know that there is a good theory not only of $\OpE{k}$-algebras $A$, but of {\em modules} over such algebras. In \S \ref{slabba} we will prove the existence of such a theory by verifying that the operads $\OpE{k}$ are {\it coherent} in the sense of Definition \symmetricref{koopa}. The proof makes use of
a general coherence criterion which we establish in \S \ref{cohcrit}.

In \S \ref{tensor2}, we will study the formation of tensor products of $\OpE{k}$-algebras. More precisely, we will give a simple ``generators and relations'' description of the tensor product $A \otimes B$ in the case where
$A$ and $B$ are free (Theorem \ref{juke}). 

In \S \ref{sec3}, we will give a geometric reformulation of the theory of $\OpE{k}$-algebras, using the formalism of factorizable (co)sheaves. In this context, it is convenient to work with {\em nonunital} algebras. In \S \ref{pluy}, we will show that there is not much difference between the theories of unital and nonunital $\OpE{k}$-algebras. More precisely, we will show that for any symmetric monoidal $\infty$-category $\calC^{\otimes}$, the $\infty$-category $\Alg_{ \OpE{k}}(\calC)$ is equivalent to a subcategory of $\Alg_{ \OpE{k} }^{\nunit}(\calC)$  which can be explicitly described (Theorem \ref{quas}).

\subsection{The $\OpE{k}$-Operads}\label{kloo}

Our goal in this section is to define the little cubes $\infty$-operads $\OpE{k}$ for $k \geq 0$, which are our main object of study throughout this paper. We begin by unfolding the definition of the topological operads $\calC_k$ described in the introduction. 

\begin{definition}\label{defE}
We define a topological category $\TopE{k}$ as follows:
\begin{itemize}
\item[$(1)$] The objects of $\TopE{k}$ are the objects $\seg{n} \in \FinSeg$.
\item[$(2)$] Given a pair of objects $\seg{m}, \seg{n} \in \TopE{k}$, a morphism from
$\seg{m}$ to $\seg{n}$ in $\TopE{k}$ consists of the following data:
\begin{itemize}
\item A morphism $\alpha: \seg{m} \rightarrow \seg{n}$ in $\FinSeg$.
\item For each $j \in \nostar{n}$ a rectilinear embedding
$\Cube{k} \times \alpha^{-1} \{j\} \rightarrow \Cube{k}$.
\end{itemize}
\item[$(3)$] For every pair of objects $\seg{m}, \seg{n} \in \TopE{k}$, we regard
$\Hom_{ \TopE{k}}( \seg{m}, \seg{n})$ as endowed with the topology induced by the presentation
$$ \Hom_{ \TopE{k}}( \seg{m}, \seg{n}) = \coprod_{ f: \seg{m} \rightarrow \seg{n} }
\prod_{ 1 \leq j \leq n} \Rect( \Cube{k} \times f^{-1} \{j\}, \Cube{k} ).$$
\item[$(4)$] Composition of morphisms in $\TopE{k}$ is defined in the obvious way.
\end{itemize}
We let $\OpE{k}$ denote the nerve of the topological category $\TopE{k}$.
\end{definition}

Corollary \toposref{tooky} implies that $\OpE{k}$ is an $\infty$-category.
There is an evident forgetful functor from $\TopE{k}$ to the (discrete) category $\FinSeg$,
which induces a functor $\OpE{k} \rightarrow \Nerve(\FinSeg)$. 

\begin{proposition}
The functor $\OpE{k} \rightarrow \Nerve(\FinSeg)$ exhibits $\OpE{k}$ as an $\infty$-operad.
\end{proposition}

\begin{proof}
We have a canonical isomorphism $\OpE{k} \simeq \Nerve^{\otimes}(\calO)$, where
$\calO$ denotes the simplicial colored operad having a single object $\Cube{k}$ with
$\Mul_{\calO}( \{ C^{k} \}_{i \in I}, C^{k} ) = \Sing X$ for $X$ the space
consisting of all rectilinear embeddings $f: \coprod_{i \in I} \Cube{k} \hookrightarrow \Cube{k}$. Since $\calO$ is a fibrant simplicial colored operad,
$\OpE{k}$ is an $\infty$-operad by virtue of Proposition \symmetricref{calp}. 
\end{proof}

\begin{definition}
We will refer to the $\infty$-operad $\OpE{k}$ as the {\it $\infty$-operad of little $k$-cubes}.
\end{definition}

\begin{remark}
Let $\Envv( \OpE{k})$ be the symmetric monoidal envelope of
$\OpE{k}$, as defined in \S \symmetricref{monenv}. We can describe the $\infty$-category
$\Envv( \OpE{k})$ informally as follows: its objects are topological space which are given
a finite unions $\coprod_{i \in I} \Cube{k}$ of the standard cube $\Cube{k}$, and its morphisms
are given by embeddings which are rectilinear on each component. The symmetric monoidal
structure on $\Envv( \OpE{k})$ is given by disjoint union.
\end{remark}

\begin{remark}\label{sove}
The mapping spaces in the topological category $\TopE{k}$ are closely related
to {\em configuration spaces}. If $I$ is a finite set and $M$ is any manifold, we let
$\Conf(I;M)$ denote the space of all {\em injective} maps from $I$ into
$M$ (regarded as an open subset of $M^I$). 
Evaluation at the origin $0 \in \Cube{k}$ induces a map
$\theta: \Rect( \Cube{k} \times I, \Cube{k}) \rightarrow \Conf(I; \Cube{k})$. We will prove that this
map is a homotopy equivalence.

Let $\overline{ \Rect}( \Cube{k} \times I, \Cube{k})$ denote the 
collection of all maps $\Cube{k} \times I \rightarrow \Cube{k}$ which
are either rectilinear embeddings, or factor as a composition
$$ \Cube{k} \times I \rightarrow I \hookrightarrow \Cube{k}.$$
Then $\theta$ factors as a composition
$$ \Rect( \Cube{k} \times I, \Cube{k})
\stackrel{ \theta'}{\rightarrow} \overline{ \Rect}( \Cube{k} \times I, \Cube{k})
\stackrel{\theta''}{\rightarrow} \Conf(I; \Cube{k})$$
where $\theta'$ is the open inclusion and $\theta''$ is given by evaluation
at the origin $0 \in \Cube{k}$. We claim that both $\theta'$ and $\theta''$ are homotopy equivalences:

\begin{itemize}
\item[$(i)$] For every map $f \in \overline{\Rect}( \Cube{k} \times I, \Cube{k})$,
let $\epsilon(f)$ denote the infumum over $i,j \in I$ of the distance from
$f( \Cube{k} \times \{i\})$ to $f( \Cube{k} \times \{j\})$ and the distance
from $f( \Cube{k} \times \{i\})$ to the boundary of $\Cube{k}$. We then
define a family of maps $\{ f_t \}_{t \in [0,1]}$ by the formula
$$ f_{t}( x_1, \ldots, x_k, i ) = f( x_1, \ldots, x_k, i) + \frac{t \epsilon(f)}{2k} (x_1, \ldots, x_k).$$
This construction determines a map $H: \overline{\Rect}( \Cube{k} \times I, \Cube{k})
\times [0,1] \rightarrow \overline{\Rect}( \Cube{k} \times I, \Cube{k})$ such that
$H | \overline{\Rect}( \Cube{k} \times I, \Cube{k}) \times \{0\}$ is the identity map
and $H$ carries $(\overline{ \Rect}( \Cube{k} \times I, \Cube{k}) \times (0,1] )$ into the open subset
$\Rect( \Cube{k} \times I, \Cube{k}) \subseteq \overline{\Rect}( \Cube{k} \times I, \Cube{k})$. 
It follows that the inclusion $\theta': \Rect(\Cube{k} \times I, \Cube{k}) \subseteq
\overline{\Rect}( \Cube{k} \times I, \Cube{k})$ is a homotopy equivalence.

\item[$(ii)$] The inclusion $j: \Conf(I; \Cube{k}) \subseteq
\overline{\Rect}( \Cube{k} \times I, \Cube{k})$ is a homotopy equivalence.
Indeed, there is a deformation retraction of $\overline{ \Rect}( \Cube{k} \times I, \Cube{k})$, which
carries a map $f: \Cube{k} \times I \rightarrow \Cube{k}$ to the family of maps
$\{ f_{t}: \Cube{k} \times I \rightarrow \Cube{k} \}_{t \in [0,1]}$ given by the formula
$$ f_{t}( x_1, \ldots, x_k, i) = f( tx_1, \ldots, tx_k, i).$$
Since $\theta''$ is a left inverse to $j$, it follows that $\theta''$ is a homotopy equivalence.
\end{itemize}
\end{remark}

\begin{example}
Suppose that $k=0$. Then $\Cube{k}$ consists of a single point, and the only rectilinear embedding
from $\Cube{k}$ to itself is the identity map. A finite collection $\{ f_i: \Cube{k} \rightarrow \Cube{k} \}_{i \in I}$ of rectlinear embeddings have disjoint images if and only if the index set $I$ has at most one element.
It follows that $\TopE{k}$ is isomorphic (as a topological category) to the subcategory
of $\FinSeg$ spanned by the injective morphisms in $\FinSeg$. We conclude that
$\OpE{k}$ is the subcategory of $\Nerve(\FinSeg)$ spanned by the injective morphisms.
\end{example}

\begin{example}\label{sulta}
Suppose that $k = 1$, so that we can identify the cube $\Cube{k}$ with the interval $(-1,1)$.
Every rectangular embedding $(-1, 1) \times I \rightarrow (-1,1)$ determines a linear
ordering of the set $I$, where $i < j$ if and only if
$f( t, i) < f(t', j)$ for all $t, t' \in (-1,1)$. This construction determines a composition of the space $\Rect( (-1,1) \times I, (-1,1))$ 
into components $\Rect_{<}( (-1,1) \times I, (-1,1))$, where $<$ ranges over all linear orderings on $I$. 
Each of the spaces $\Rect_{<}( (-1,1) \times I, (-1,1))$ is a nonempty convex set
and therefore contractible. It follows that $\Rect( (-1,1) \times I, (-1,1) )$ is homotopy equivalent to the discrete set of all linear orderings on $I$.

Using these homotopy equivalences, we obtain a weak equivalence of topological categories
$\TopE{1} \rightarrow \CatAss$, where $\CatAss$ is the category defined in Example \symmetricref{xe2}.
Passing to the homotopy coherent nerves, we obtain an equivalence of $\infty$-operads
$\OpE{1} \simeq \Ass$.
\end{example}

We can use the relationship between rectilinear embedding spaces and configuration spaces to
establish some basic connectivity properties of the $\infty$-operads $\OpE{k}$:

\begin{proposition}\label{slobe}
Let $k \geq 0$. For every pair of integers $m, n \geq 0$, the map of topological spaces
$\bHom_{ \TopE{k}}( \seg{m}, \seg{n}) \rightarrow \Hom_{ \FinSeg}( \seg{m}, \seg{n} )$
is $(k-1)$-connective. 
\end{proposition}

\begin{proof}
Unwinding the definitions, this is equivalent to the requirement that
for every finite set $I$, the space of rectinilinear embeddings $\Rect( \Cube{k} \times I, \Cube{k})$ is $(k-1)$-connective. This space is homotopy equivalent (via evaluation at the origin) to the configuration space $\Conf(I; \Cube{k} )$ of injective maps $I \hookrightarrow \Cube{k}$ (Remark \ref{sove}). We will prove more generally that $\Conf( J, \Cube{k} - F)$ is $(k-1)$-connective, where
$J$ is any finite set and $F$ is any finite subset of $\Cube{k}$. The proof proceeds by induction
on the number of elements of $J$. If $J = \emptyset$, then $\Conf( J, \Cube{k} - F)$ consists of a single point and there is nothing to prove. Otherwise, choose an element $j \in J$. Evaluation at
$j$ determines a Serre fibration $\Conf(J, \Cube{k} - F) \rightarrow \Cube{k} - F$, whose fiber over a point $x$ is the space $\Conf( J - \{j\}, \Cube{k} - ( F \cup \{x\} ))$. The inductive hypothesis guarantees
that these fibers are $(k-1)$-connective. Consequently, to show that $\Conf(J, \Cube{k} - F)$ is $(k-1)$-connective, it suffices to show that $\Cube{k} - F$ is $(k-1)$-connective. In other words, we must show that
for $m < k$, every map
$g_0: S^{m-1} \rightarrow \Cube{k} - F$ can be extended to a map
$g: D^{m} \rightarrow \Cube{k} - F$, where $D^{m}$ denotes the unit disk of dimension
$m$ and $S^{m-1}$ its boundary sphere. Without loss of generality, we may assume that
$g_0$ is smooth. Since $\Cube{k}$ is contractible, we can extend
$g_0$ to a map $g: D^{m} \rightarrow \Cube{k}$, which we may also assume to be smooth
and transverse to the submanifold $F \subseteq \Cube{k}$. Since $F$ has codimension $k$
in $\Cube{k}$, $g^{-1} F$ has codimension $k$ in $D^{m}$, so that $g^{-1} F = \emptyset$
(since $m < k$) and $g$ factors through $\Cube{k} - F$, as desired.
\end{proof}

For each $k \geq 0$, there is a stabilization functor $\TopE{k} \rightarrow \TopE{k+1}$
which is the identity on objects and is given on morphisms by taking the product with the
interval $(-1,1)$. This functor induces a map of $\infty$-operads
$\OpE{k} \rightarrow \OpE{k+1}$. Proposition \ref{slobe} immediately implies the following:

\begin{corollary}\label{infec}
Let $\OpE{\infty}$ denote the colimit of the sequence of $\infty$-operads
$$ \OpE{0} \rightarrow \OpE{1} \rightarrow \OpE{2} \rightarrow \ldots$$
Then the canonical map $\OpE{\infty} \rightarrow \Nerve(\FinSeg)$ is an equivalence
of $\infty$-operads.
\end{corollary}

Consequently, if $\calC^{\otimes}$ is a symmetric monoidal $\infty$-category, then
the $\infty$-category $\CAlg(\calC)$ of commutative algebra objects of $\calC$ can be identified with the homotopy limit of the tower of $\infty$-categories $\{ \Alg_{ \OpE{k}}(\calC) \}_{k \geq 0}$. In many situations, this tower actually stabilizes at some finite stage:

\begin{corollary}\label{manik}
Let $\calC^{\otimes}$ be a symmetric monoidal $\infty$-category. Let $n \geq 1$, and assume that the underlying $\infty$-category $\calC$ is equivalent to an $n$-category (that is, the mapping spaces
$\bHom_{\calC}(X,Y)$ are $(n-1)$-truncated for every pair of objects $X,Y \in \calC$; see \S \toposref{ncats}). Then the map $\OpE{k} \rightarrow \Nerve(\FinSeg)$ induces an
equivalence of $\infty$-categories $\CAlg(\calC) \rightarrow \Alg_{ \OpE{k}}(\calC)$ for
$k > n$.
\end{corollary}

\begin{proof}
Let $C$ and $D$ be objects of $\calC^{\otimes}$, corresponding to finite sequences of
objects $( X_1, \ldots, X_m)$ and $(Y_1, \ldots, Y_{m'})$ of objects of $\calC$. Then
$\bHom_{\calC^{\otimes}}(C,D)$ can be identified with the space
$$ \coprod_{ \alpha: \seg{m} \rightarrow \seg{n} } \prod_{ 1 \leq j \leq m'}
\bHom_{\calC}( \otimes_{ \alpha(i) = j} X_i, Y_j ),$$
and is therefore also $(n-1)$-truncated. Consequently, $\calC^{\otimes}$ is equivalent to an $n$-category. Proposition \ref{slobe} implies that the forgetful functor
$\OpE{k} \rightarrow \Nerve(\FinSeg)$ induces an equivalence of the underlying homotopy $n$-categories, and therefore induces an equivalence $\theta: \Fun_{ \Nerve(\FinSeg)}( \Nerve(\FinSeg), \calC^{\otimes}) \rightarrow \Fun_{ \Nerve(\FinSeg)}( \OpE{k}, \calC^{\otimes})$. The desired result now follows from the observation that a map $A \in \Fun_{ \Nerve(\FinSeg)}( \Nerve(\FinSeg), \calC^{\otimes})$ is a commutative algebra object of $\calC$ if and only if $\theta(A)$ is an
$\OpE{k}$-algebra object of $\calC$.
\end{proof}





\subsection{The Additivity Theorem}\label{sass1}

If $K$ is a pointed topological space, then the $k$-fold loop space
$\Omega^{k}(K)$ carries an action of the (topological) little cubes operad
$\calC_k$ of the introduction. Passing to singular complexes, we deduce that if
$X \in \SSet_{\ast}$, then the $k$-fold loop space $\Omega^{k}(X)$ can be promoted
to an $\OpE{k}$-algebra object of the $\infty$-category $\SSet$ of spaces.
The work of May provides a converse to this observation:
if $Z$ is a {\it grouplike} $\OpE{k}$-algebra object of $\SSet$ (see Definition \ref{ungwar}), then
$Z$ is equivalent to $\Omega^{k}(Y)$ for some pointed space $Y \in \SSet_{\ast}$ (see
Theorem \ref{preslage} for a precise statement). The delooping process $Z \mapsto Y$ is compatible with products in $Z$. Consequently, if $Z$ is equipped with a {\em second} action of the operad $\OpE{k'}$, which is suitable compatible with the $\OpE{k}$ action on $Z$, then we should expect that the space $Y$ again carries an action of $\OpE{k'}$, and is therefore {\em itself} homotopy equivalent to $\Omega^{k'}(X)$ for some pointed space $X \in \SSet_{\ast}$. Then $Z \simeq \Omega^{k+k'}(X)$ carries an action of the $\infty$-operad
$\OpE{k+k'}$. Our goal in this section is to show that this phenomenon is quite general, and applies to
algebra objects of an arbitrary symmetric monoidal $\infty$-category $\calC^{\otimes}$: namely,
giving an $\OpE{k+k'}$-algebra object of $\calC$ is equivalent to giving an object $A \in \calC$ which is equipped with {\em commuting} actions of the $\infty$-operads $\OpE{k}$ and $\OpE{k'}$. More precisely, we have a canonical equivalence $\Alg_{ \OpE{k+k'}}(\calC) \simeq \Alg_{ \OpE{k}}( \Alg_{\OpE{k'}}(\calC) )$
(Theorem \ref{slide}). Equivalently, we can identify $\OpE{k+k'}$ with the {\em tensor product} of
the $\infty$-operads $\OpE{k}$ and $\OpE{k'}$ (see Definition \ref{ahmad}). We first describe the bifunctor
$\OpE{k} \times \OpE{k'} \rightarrow \OpE{k+k'}$ which gives rise to this identification.

\begin{construction}\label{laxose}
Choose nonnegative integers $k, k'$. We define a topological functor
$\times: \TopE{k} \times \TopE{k'} \rightarrow \TopE{k+k'}$
as follows:
\begin{itemize}
\item[$(1)$] The diagram of functors
$$ \xymatrix{ \TopE{k} \times \TopE{k'} \ar[r]^{ \times} \ar[d] & \TopE{k+k'} \ar[d] \\
\Nerve(\FinSeg) \times \Nerve(\FinSeg) \ar[r]^{ \wedge} & \Nerve(\FinSeg) }$$
commutes, where $\wedge$ denotes the smash product functor on pointed finite sets (Notation \symmetricref{crumb}). In particular, the functor $\times$ is given on objects by the formula
$\seg{m} \wedge \seg{n} = \seg{mn}$. 
\item[$(2)$] Suppose we are given a pair of morphisms $\overline{\alpha} : \seg{m} \rightarrow \seg{n}$ in $\TopE{k}$ and $\overline{\beta}: \seg{m'} \rightarrow \seg{n'}$ in $\TopE{k'}$.
Write $\overline{\alpha} = ( \alpha, \{ f_j: \Cube{k} \times \alpha^{-1} \{j \}, \Cube{k} \}_{j \in \nostar{n} })$
and $\overline{\beta} = ( \beta, \{ f'_{j}: \Cube{k} \times \beta^{-1} \{j\}, \Cube{k} \}_{j' \in \nostar{n'} })$. 
We then define $\overline{\alpha} \times \overline{\beta}: \seg{mm'} \rightarrow \seg{nn'}$ to be given by the pair $( \alpha \wedge \beta, \{ f_{j} \times f_{j'}: \Cube{k+k'} \times \alpha^{-1} \{j\} \times \beta^{-1} \{j' \}, \Cube{k+k'} \}_{j \in \nostar{n}, j' \in \nostar{n'} })$.
\end{itemize}

Passing to homotopy coherent nerves, we obtain a bifunctor of $\infty$-operads (see Definition \symmetricref{rocks}) $\OpE{k} \times \OpE{k'} \rightarrow \OpE{k+k'}$. 
\end{construction}

A version of the following fundamental result was proven by Dunn (see \cite{dunn}): 

\begin{theorem}[Dunn]\label{slide}
Let $k, k' \geq 0$ be nonnegative integers. Then the bifunctor
$\OpE{k} \times \OpE{k'} \rightarrow \OpE{k+k'}$ of Construction \ref{laxose}
exhibits the $\infty$-operad $\OpE{k+k'}$ as
a tensor product of the $\infty$-operads $\OpE{k}$ with $\OpE{k'}$ (see Definition \ref{ahmad}).
\end{theorem}

\begin{example}[Baez-Dolan Stabilization Hypothesis]\label{swee}
Theorem \ref{slide} implies that supplying an $\OpE{k}$-monoidal structure on an $\infty$-category
$\calC$ is equivalent to supplying $k$ compatible monoidal structures on $\calC$.
Fix an integer $n \geq 1$, and let $\Cat_{\infty}^{\leq n}$ denote the full subcategory of
$\Cat_{\infty}$ spanned by those $\infty$-categories which are equivalent to $n$-categories.
For $\calC, \calD \in \Cat_{\infty}^{\leq n}$, the mapping space $\bHom_{ \Cat_{\infty} }(\calC, \calD)$
is the underlying Kan complex of $\Fun(\calC, \calD)$, which is equivalent to an $n$-category (Corollary \toposref{zooka}). It follows that $\Cat_{\infty}^{\leq n}$ is equivalent to an $(n+1)$-category.
Let us regard $\Cat_{\infty}^{\leq n}$ as endowed with the Cartesian monoidal structure.
Corollary \ref{manik} implies that $\CAlg( \Cat_{\infty}^{\leq n}) \simeq \Alg_{ \OpE{k}}( \Cat_{\infty}^{\leq n} )$ for $k \geq n+2$. Combining this observation with Corollary \symmetricref{co1}, we deduce that if $\calC$ is
an $n$-category, then supplying an $\OpE{k}$-monoidal structure on $\calC$ is equivalent to supplying
a symmetric monoidal structure on $\calC$. This can be regarded as a version of the ``stabilization hypothesis'' proposed in \cite{baezdolan} (the formulation above applies to $n$-categories where all
$k$-morphisms are invertible for $k > 1$, but the argument can be applied more generally.)
\end{example}

\begin{example}[Braided Monoidal Categories]\label{sweeta}
Let $\calC$ be an ordinary category. According to Example \ref{swee}, supplying an 
$\OpE{k}$-monoidal structure on $\Nerve(\calC)$ is equivalent to supplying a symmetric monoidal
structure on the $\infty$-category $\calC$ if $k \geq 3$. If $k=1$, then supplying
an $\OpE{k}$-monoidal structure on $\Nerve(\calC)$ is equivalent to supplying a
monoidal structure on $\calC$ (combine Example \ref{sulta}, Proposition \symmetricref{algass}, and Remark \monoidref{otherlander}). Let us therefore focus our attention on the case $n=2$.
In view of Corollary \symmetricref{co1}, giving an $\OpE{2}$-monoidal structure on $\Nerve(\calC)$ is equivalent
to exhibiting $\Nerve(\calC)$ as an $\OpE{2}$-algebra object of
$\Cat_{\infty}$. Theorem \ref{slide} provides an equivalence $\Alg_{ \OpE{2}}(\Cat_{\infty})
\simeq \Alg_{ \OpE{1}}( \Alg_{\OpE{1}}(\Cat_{\infty}) )$. Combining this with Example
\ref{sulta}, Proposition \symmetricref{algass}, and Remark \monoidref{otherlander}, we can
view $\Nerve(\calC)$ as an (associative) monoid object in the $\infty$-category $\Cat_{\infty}^{\Mon}$ of monoidal $\infty$-categories. This structure allows us to view $\calC$ as a monoidal category with respect to some tensor product $\otimes$, together with a {\em second} multiplication given by a monoidal functor
$$\boxtimes: (\calC, \otimes) \times (\calC, \otimes) \rightarrow (\calC, \otimes).$$ 
This second multiplication also has a unit, which is a functor from
the one-object category $[0]$ into $\calC$. Since this functor is required to be monoidal, it carries
the unique object of $[0]$ to the unit object ${\bf 1} \in \calC$, up to canonical isomorphism.
It follows that ${\bf 1}$ can be regarded as a unit with respect to both tensor product operations
$\otimes$ and $\boxtimes$.

We can now exploit the classical Eckmann-Hilton argument to show that the tensor product
functors $\otimes, \boxtimes: \calC \times \calC \rightarrow \calC$ are isomorphic.
Namely, our assumption that $\boxtimes$ is a monoidal functor gives a chain of isomorphisms
\begin{eqnarray} X \boxtimes Y & \simeq & (X \otimes {\bf 1}) \boxtimes ( {\bf 1} \boxtimes Y) \\
& \simeq & (X \boxtimes {\bf 1}) \otimes ({\bf 1} \boxtimes Y) \\
& \simeq & X \otimes Y
\end{eqnarray}
depending naturally on $X$ and $Y$. Consequently, $\boxtimes$ is determined by
$\otimes$ as a functor from $\calC \times \calC$ into $\calC$. However, it gives rise to additional
data when viewed as a {\em monoidal} functor: a monoidal structure on the tensor product
functor $\otimes: \calC \times \calC \rightarrow \calC$ supplies a canonical isomorphism
$$ (W \otimes X) \otimes (Y \otimes Z) \simeq (W \otimes Y) \otimes (X \otimes Z).$$
Taking $W$ and $Z$ to be the unit object, we get a canonical isomorphism
$\sigma_{X,Y}: X \otimes Y \rightarrow Y \otimes X$. Conversely, if we are given 
a collection of isomorphisms $\sigma_{X,Y}: X \otimes Y \rightarrow Y \otimes X$, we can
try to endow $\otimes: \calC \times \calC \rightarrow \calC$ with the structure of a monoidal
functor by supplying the isomorphisms
$$(W \otimes X) \otimes (Y \otimes Z) \simeq W \otimes (X \otimes Y) \otimes Z
\stackrel{ \sigma_{X,Y}}{\simeq} W \otimes (Y \otimes X) \otimes Z \simeq (W \otimes Y) \otimes (X \otimes Z)$$
together with the evident isomorphism ${\bf 1} \otimes {\bf 1} \simeq {\bf 1}$. Unwinding
the definitions, we see that these isomorphisms supply a monoidal structure on the
functor $\otimes$ if and only if the following condition is satisfied:

\begin{itemize}
\item[$(1)$] For every triple of objects $X,Y,Z \in \calC$, the isomorphism
$\sigma_{X,Y \otimes Z}$ is given by the composition
$$ X \otimes (Y \otimes Z) \simeq (X \otimes Y) \otimes Z
\stackrel{ \sigma_{X,Y} }{\simeq} (Y \otimes X) \otimes Z
\simeq Y \otimes (X \otimes Z) \stackrel{\sigma_{X,Z} }{\rightarrow} Y \otimes (Z \otimes X)
\simeq (Y \otimes Z) \otimes X.$$
\end{itemize}

In this case, we have a diagram of monoidal functors
$$ \xymatrix{ \calC \times \calC \times \calC \ar[r]^{ \otimes \times \id} \ar[d]^{\id \times \otimes} & \calC \times \calC \ar[d]^{\otimes} \\
\calC \times \calC \ar[r]^{\otimes} & \calC }$$
such that the underlying diagram of categories commutes up to canonical isomorphism $\alpha$
(supplied by the monoidal structure on $\calC$). Unwinding the definitions, we see that the natural transformation $\alpha$ is an isomorphism of {\em monoidal} functors if and only if the following additional condition is satisfied:
\begin{itemize}
\item[$(2)$] For every triple of objects $X,Y, Z \in \calC$, the isomorphism $\sigma_{X \otimes Y,Z}$ is given by the composition $$ (X \otimes Y) \otimes Z \simeq X \otimes (Y \otimes Z)
\stackrel{\sigma_{Y,Z}}{\simeq} X \otimes (Z \otimes Y)
\simeq (X \otimes Z) \otimes Y \stackrel{ \sigma_{X,Z}}{\simeq}
(Z \otimes X) \otimes Y \simeq Z \otimes (X \otimes Y).$$
(Equivalently, the inverse maps $\sigma_{X,Y}^{-1}: Y \otimes X \simeq X \otimes Y$ satisfy
condition $(1)$.)
\end{itemize}
A natural isomorphism $\sigma_{X,Y}: X \otimes Y \simeq Y \otimes X$ is called a
{\it braiding} on the monoidal category $(\calC, \otimes)$ if it satisfies conditions $(1)$ and $(2)$.
A {\it braided monoidal category} is a monoidal category equipped with a braiding. We can summarize our discussion as follows: if $\calC$ is an ordinary category, then endowing $\calC$ with the structure of a braided monoidal category is equivalent to endowing the nerve $\Nerve(\calC)$ with the structure of an $\OpE{2}$-monoidal $\infty$-category.
\end{example}

\begin{remark}
It follows from Example \ref{sweeta} that if $\calC$ is a monoidal category containing
a sequence of objects $X_1, \ldots, X_n$, then the tensor product
$X_1 \otimes \cdots \otimes X_n$ is the fiber of a local system of objects of $\calC$
over the space $\Rect( \Cube{2} \times \{1, \ldots, n \}, \Cube{2})$. In other words, the tensor product
$X_1 \otimes \cdots \otimes X_n$ is endowed with an action of the fundamental group
$\pi_1 \Conf( \{1, \ldots, n\}, \R^2)$ of configurations of $n$ distinct points in the plane $\R^{2}$ (Remark \ref{sove}). The group $\pi_1 \Conf( \{1, \ldots, n\}, \R^{2})$ is the Artin pure braid group on $n$ strands.
The action of $\pi_1 \Conf( \{1, \ldots, n\}, \R^{2})$ on $X_1 \otimes \cdots \otimes X_n$ can be constructed by purely combinatorial means, by matching the standard generators of the Artin 
braid group with the isomorphisms $\sigma_{ X_{i}, X_{j} }$. However, Theorem \ref{slide} provides a much more illuminating geometric explanation of this phenomenon.
\end{remark}

The proof of Theorem \ref{slide} will occupy our attention for the remainder of this section. In
what follows, we will assume that the reader is familiar with the notation introduced in \S \ref{wreath}.
We begin by observing that the bifunctor $\OpE{k} \times \OpE{k'} \rightarrow \OpE{k+k'}$
factors as a composition
$$ \OpE{k} \times \OpE{k'} \stackrel{\theta'}{\rightarrow} \OpE{k} \wreath \OpE{k'}
\stackrel{\theta}{\rightarrow} \OpE{k+k'},$$
where $\theta'$ is the map described in Remark \ref{hulker}. In the special case
where $k=1$, we will denote the functor $\theta$ by
$\Psi_k: \OpE{1} \wreath \OpE{k} \rightarrow \OpE{k+1}$. The key step in the proof
of Theorem \ref{slide} is the following:

\begin{proposition}\label{inku}
Let $k \geq 0$ be a nonnegative integer. Then the map
$\theta: \OpE{1} \wreath \OpE{k} \rightarrow \OpE{k+1}$ induces
a weak equivalence of $\infty$-preoperads $( \OpE{1} \wreath \OpE{k}, M) \rightarrow \OpE{k+1}^{\natural}$. Here $M$ denotes the collection of all inert morphisms in
$\OpE{1} \wreath \OpE{k}$.
\end{proposition}

Assuming Proposition \ref{inku} for the moment, we can give the proof of Theorem \ref{slide}.

\begin{proof}[Proof of Theorem \ref{slide}]
We proceed by induction on $k$. If $k = 0$, the desired result follows from Proposition
\symmetricref{skymall}, since the $\infty$-operad $\OpE{k}$ is unital. If $k = 1$, we consider the factorization
$$ \OpE{1}^{\natural} \odot \OpE{k'}^{\natural} \rightarrow
( \OpE{1} \wreath \OpE{k'}, M) \rightarrow \OpE{1+k'}^{\natural}$$
and apply Proposition \ref{inku} together with Theorem \ref{kuj}. If
$k > 1$, we have a commutative diagram
$$ \xymatrix{ \OpE{1}^{\natural} \odot \OpE{k-1}^{\natural} \odot \OpE{k'}^{\natural} \ar[r] \ar[d] & \OpE{k}^{\natural} \odot \OpE{k'}^{\natural} \ar[d] \\
\OpE{1}^{\natural} \odot \OpE{k+k'-1}^{\natural} \ar[r] & \OpE{k+k'}^{\natural}. }$$
The inductive hypothesis guarantees that the horizontal maps and the left vertical map are
weak equivalences of $\infty$-preoperads, so that the right vertical map is a weak equivalence as well.
\end{proof}

We now turn to the proof of Proposition \ref{inku}. Our first step is to establish the following lemma.

\begin{lemma}\label{cruss}
Let $\seg{m}$ be an object of $\OpE{k+1}$, and let
$\calJ = ( \OpE{1} \wreath \OpE{k}) \times_{ \OpE{k+1}} (\OpE{k+1})_{\seg{m}/}$,
and let $\calJ^{0}$ denote the full subcategory of $\calJ$ spanned by objects
which correspond to inert morphisms $\seg{m} \rightarrow \Psi_k(X)$ in
$\OpE{k+1}$. Then the inclusion $\calJ^{0} \subseteq \calJ$ is left anodyne.
\end{lemma}

To prove this, we will need to introduce some notation.

\begin{construction}\label{tupus}
Fix an integer $m \geq 0$. Let $P$ denote the partially ordered set defined as follows:
\begin{itemize}
\item[$(i)$] An element of $P$ consists of an subset $I \subseteq (-1,1)$ which
can be written as a finite disjoint union of closed intervals, together with
a surjection of finite sets $\chi: \nostar{m} \rightarrow \pi_0 I$
is a surjection of finite sets.

\item[$(ii)$] We have $(I, \chi) \leq (I', \chi')$ if and only if $I \subseteq I'$ and
the diagram
$$ \xymatrix{ & \pi_0 I \ar[dr] & \\
\nostar{m} \ar[ur]^{\chi} \ar[rr]^{\chi'} & & \pi_0 I' }$$
commutes.
\end{itemize}

Fix an active morphism $\alpha: \seg{m} \rightarrow \seg{n}$ in $\FinSeg$.
We can lift $\seg{n}$ to an object $( \seg{n} ) \in \OpE{1} \wreath \Nerve(\FinSeg)$,
allowing us to identify $\alpha$ with an object $D$ of the $\infty$-category
$$ \calD = (\OpE{1} \wreath \Nerve(\FinSeg) ) \times_{ \Nerve(\FinSeg) }
\Nerve(\FinSeg)^{\seg{m} / }.$$
Let $\calD' = \calD^{/D}$ denote the $\infty$-category whose objects
consist of maps $\gamma: ( \seg{m_{1}}, \ldots, \seg{m_{b}}) \rightarrow
( \seg{n} )$ in $\OpE{1} \wreath \Nerve(\FinSeg)$ together with a commutative diagram
$$ \xymatrix{ \seg{m} \ar[r]^-{\beta} \ar[d]^{\id} & \seg{m_1 + \ldots + m_{b} } \ar[d] \\
\seg{m} \ar[r]^{\alpha} & \seg{n} }$$
in $\FinSeg$, and let $\calC_0$ denote the full subcategory spanned by those diagrams where
$\beta$ is an isomorphism and each of the integers $m_{i}$ is positive.

Given an element $(I, \chi) \in P$, write
$I$ as a disjoint union of intervals $I_{1} \cup \ldots \cup I_{b}$, where
$x < y$ whenever $x \in I_{i}$, $y \in I_{j}$, and $i < j$. For $1 \leq i \leq b$, there
is a unique order-preserving bijection $\chi^{-1} \{ I_{i} \} \simeq \nostar{m_i}$ for some
$m_i > 0$; these bijections together determine an isomorphism
$\beta: \seg{m} \rightarrow \seg{m_1 + \cdots + m_{b} }$ in $\FinSeg$.
There is a unique map $\gamma_0: ( \seg{m_1}, \ldots, \seg{m_b}) \rightarrow
( \seg{n} )$ in $\Nerve(\FinSeg)^{\amalg}$ such that $\alpha = \overline{\gamma}_0 \circ \beta$, where
$\overline{\gamma}_0$ denotes the image of $\gamma_0$ in $\FinSeg$. The morphism
$\gamma_0$ lifts to a morphism $\gamma$ in $\OpE{1} \wreath \Nerve(\FinSeg)$ by
specifying rectilinear embeddings $\{ e_{j}: (-1,1) \rightarrow (-1,1) \}_{1 \leq j \leq b}$ such that
$e_{j} (-1,1)$ is the interior of $I_{j}$. The pair $(\beta, \gamma)$ determines an object of
$\calC_0$, and the construction $(I, \chi) \mapsto (\beta, \gamma)$ extends naturally
to a functor $\theta: \Nerve(P) \rightarrow \calC_0$.
\end{construction}

\begin{lemma}\label{buz}
Let $q: X \rightarrow S$ be a Cartesian fibration of simplicial sets. If each
fiber of $q$ is weakly contractible, then $q$ is a weak homotopy equivalence.
\end{lemma}

\begin{proof}
We will prove that for any map of simplicial sets $S' \rightarrow S$, the induced map
$q_{S'}: X \times_{S} S' \rightarrow S'$ is a weak homotopy equivalence.
Since the collection of weak homotopy equivalences is stable under filtered colimits,
we can reduce to the case where the simplicial set $S'$ is finite. We now work by induction
on the dimension of $n$ of $S'$. If $S'$ is empty, the result is obvious; otherwise,
let $T$ be the set of nondegenerate $n$-simplices of $S'$ so that we have a pushout
diagram
$$ \xymatrix{ T \times \bd \Delta^n \ar[r] \ar[d] & T \times \Delta^n \ar[d] \\
S'' \ar[r] & S'. }$$
The inductive hypothesis guarantees that $q_{T \times \bd \Delta^n}$ and
$q_{S''}$ are weak homotopy equivalences. Since the usual model structure on
the category of simplicial sets is left proper, it will suffice to show that
$q_{T \times \Delta^n}$ is a weak homotopy equivalence. Since the collection
of weak homotopy equivalences is stable under coproducts, we can reduce to the case
where $S' \simeq \Delta^n$ is an $n$-simplex.

We wish to show that $X' = X \times_{S} \Delta^n$ is weakly contractible.
Note that $X'$ is an $\infty$-category. Since the map $X' \rightarrow \Delta^n$
is a Cartesian fibration, the inclusion $X'_0 = X' \times_{ \Delta^n} \{0\} \subseteq X'$
admits a right adjoint. It follows that $X'$ is weakly homotopy equivalent to
$X'_0$, which is a fiber of the map $q$ and therefore weakly contractible by assumption.
\end{proof}

\begin{lemma}\label{cornall}
Fix an integer $b \geq 0$, and let $Q_b$ denote the set of sequences
$(I_1, \ldots, I_b)$, where each $I_{j} \subseteq (-1,1)$ is a closed interval,
and we have $x < y$ whenever $x \in I_{i}$, $y \in I_{j}$, and $i < j$.
We regard $Q$ as a partially ordered set, where
$(I_1, \ldots, I_b) \leq (I'_1, \ldots, I'_{b})$ if $I_j \subseteq I'_{j}$ for $1 \leq i \leq j$.
Then the nerve $\Nerve(Q_b)$ is weakly contractible.
\end{lemma}

\begin{proof}
The proof proceeds by induction on $b$. If $b = 0$, then $Q_{b}$ has a single element
and there is nothing to prove. Otherwise, we observe that ``forgetting'' the last coordinate
induces a Cartesian fibration $q: \Nerve(Q_b) \rightarrow \Nerve(Q_{b-1})$. 
We will prove that the fibers of $q$ are weakly contractible, so that
$q$ is a weak homotopy equivalence (Lemma \ref{buz}). Fix
an element $x = ( [t_1, t'_1], [t_2, t'_2], \ldots, [t_{b-1}. t'_{b-1}]) \in Q_{b-1}$. 
Then $q^{-1}\{x\}$ can be identified with the nerve of the partially ordered set
$Q' = \{ (t_b, t'_b): t'_{b-1} < t_b < t'_b < 1 \}$, where
$(t_b, t'_b) \leq (s_{b}, s'_{b})$ if $t_b \geq s_{b}$ and $t'_{b} \leq s'_{b}$.

The map $(t_b, t'_b) \mapsto t'_{b}$ is a monotone map from
$Q'$ to the open interval $(t'_{b-1}, 1)$. This map determines a coCartesian
fibration $q': \Nerve(Q') \rightarrow \Nerve(t'_{b-1}, 1)$. The
fiber of $q'$ over a point $s$ can be identified with the
{\em opposite} of the nerve of the interval $(t'_{b-1}, s)$, and is therefore
weakly contractible. Applying Lemma \ref{buz}, we deduce that
$q'$ is a weak homotopy equivalence, so that $\Nerve(Q')$ is weakly contractible as desired.
\end{proof}

\begin{lemma}\label{colman}
Let $\theta: \Nerve(P) \rightarrow \calC_0$ be the functor of Construction \ref{tupus}.
Then $\theta$ induces a cofinal map $\Nerve(P)^{op} \rightarrow \calC_{0}^{op}$.
\end{lemma}

\begin{proof}
Let $\overline{\calD} = (\Ass \wreath \Nerve(\FinSeg)) \times_{ \Nerve(\FinSeg)}
\Nerve(\FinSeg)^{\seg{m}/ }$, let $\overline{D}$ denote the image of
$D$ in $\overline{\calD}$ under the map induced by the trivial Kan fibration
$\OpE{1} \rightarrow \Ass$ of Example \ref{sulta}. Let 
$\overline{\calD}' = \overline{\calD}^{/ \overline{D}}$ and let
$\overline{\calC}_0$ be the essential image of $\calC_0$ in
$\overline{\calD'}$. Finally, let $\overline{\theta}$ denote the composition
of $\theta$ with the trivial Kan fibration $\calC_0 \rightarrow \overline{\calC}_0$.
We are now reduced to the (purely combinatorial) problem of showing that
$\overline{\theta}$ induces a cofinal map $\Nerve(P)^{op} \rightarrow \overline{\calC}^{op}_{0}$.

According to Theorem \toposref{hollowtt}, it will suffice to prove that for every object
$C \in \overline{\calC}_0$, the $\infty$-category
$\Nerve(P)_{/C} = \Nerve(P) \times_{ \overline{\calC}_0} ( \overline{\calC_0})_{/C}$
is weakly contractible. Note that $\infty$-category $\overline{\calC}_0$ is actually
the nerve of a small category $\calJ$ whose objects consist of isomorphisms
$$ \beta: \nostar{m} \simeq \coprod_{1 \leq i \leq b} \nostar{m_{i}}$$
where each $m_i > 0$, together with a linear ordering $<$ on the set of indices
$\{ 1, \ldots, b \}$. Without loss of generality, we may assume that
$C = ( \beta, <)$, where $<$ is the standard linear ordering on the set
$\{ 1, \ldots, b \}$. Then $\Nerve(P)_{/C}$ can be identified with the nerve of the
partially ordered set $P'$ of triples $(I, \chi, f)$, where $I \subseteq (-1,1)$ is a disjoint
union of closed intervals, $\chi: \nostar{m} \rightarrow \pi_0 I$ is a surjection of finite
sets, and $f: \pi_0 I \rightarrow \{1, \ldots, b\}$ is a nonstrictly increasing map
such that $\beta^{-1} \nostar{m_i} = (\chi \circ f)^{-1} \{i\}$. 
Let $P'_0$ denote the partially ordered subset consisting of those pairs for which
$f$ is a bijection. We observe that the inclusion $\Nerve(P'_0) \subseteq \Nerve(P')$
has a left adjoint. Consequently, it will suffice to show that $\Nerve(P'_0)$ is weakly contractible, which follows from Lemma \ref{cornall}.
\end{proof}

\begin{construction}\label{swellman}
Let $\alpha: \seg{m} \rightarrow \seg{n}$ be an active morphism in $\FinSeg$, and let
$P$ be the partially ordered set of Construction \ref{tupus}. We let
$P^{\triangleleft}$ denote the partially ordered set obtained from $P$ by adjoining
a new smallest element, which we will denote by $- \infty$. We define a contravariant functor $T$
from $P^{\triangleleft}$ to the category of topological spaces as follows:
\begin{itemize}
\item[$(i)$] Let $(I, \chi)$ be an element of $P$. Then
$$ T( I, \chi) = \prod_{i \in \nostar{n}} \prod_{j \in \pi_0 I} \Rect( \Cube{k} 
\times ( \alpha^{-1} \{i\} \cap \chi^{-1}(j)), \Cube{k}).$$
\item[$(ii)$] We set $T( - \infty ) = \prod_{ i \in \nostar{n}} \Rect( \Cube{k+1}
\times \alpha^{-1} \{i\}, \Cube{k+1}).$
\item[$(iii)$] For $(I, \chi) \in P$, the map $T( I, \chi) \rightarrow T( - \infty)$ carries
an element $\{ e_{i,j} \in \Rect( \Cube{k} \times ( \alpha^{-1} \{i\}, \chi^{-1} \{j\}), \Cube{k} )$
to the collection of rectilinear embeddings $\{ e'_{i} \in \Rect( \Cube{k+1} \times \alpha^{-1} \{i\}, \Cube{k+1}) \}$ characterized by the property that for $a \in \nostar{m}$, the image
of $e'_{ \alpha(a)} | \Cube{k+1} \times \{ a\}$ is the product of
the image of $J \times e_{ \alpha(a), \chi(j) }$, where $J$ is the interior of the connected
component of $I$ given by $J$.
\end{itemize}
\end{construction}

\begin{lemma}\label{cisgn}
The contravariant functor $T$ of Construction \ref{swellman} exhibits the simplicial set
$\Sing_{\bigdot} T( - \infty)$ as a homotopy colimit of the diagram of simplicial sets
$\{ \Sing_{\bigdot} T(x) \}_{x \in P}$.
\end{lemma}

\begin{proof}
Given a finite set $S$ and an open interval $J$, we let $\Conf'(S, J \times \Cube{k})$
denote the space of all maps $S \rightarrow J \times \Cube{k}$ such that the composite map
$S \rightarrow J \times \Cube{k} \rightarrow \Cube{k}$ is injective. Since $J$ is contractible, we
deduce:
\begin{itemize}
\item[$(\ast)$] The projection map $J \times \Cube{k} \rightarrow \Cube{k}$ induces a homotopy equivalence $\Conf'( S, J \times \Cube{k}) \rightarrow \Conf(S, \Cube{k})$.
\end{itemize}

We define a contravariant functor $T_0$ from $P^{\triangleleft}$ to the category of topological spaces as follows:
\begin{itemize}
\item[$(i)$] Let $(I, \chi)$ be an element of $P$. Then $T_0(I, \chi)$ is the space
$$ \prod_{ i \in \nostar{n} } \prod_{j \in \pi_0 I} \Conf'( \alpha^{-1} \{i\} \cap \chi^{-1} \{j\}, I_{j} \times \Cube{k} ),$$ 
where $I_{j}$ denotes the interior of the connected component of $I$ corresponding to
$j \in \pi_0 I$. 
\item[$(ii)$] Set $T_0( - \infty ) = \prod_{ i \in \nostar{n} } \Conf( \alpha^{-1} \{i\}, \Cube{k+1})$.
\item[$(iii)$] For $(I, \chi) \in P$, the map $T_0(I, \chi) \rightarrow T_0(- \infty)$ is the canonical inclusion.
\end{itemize}
There is a natural transformation of functors $\gamma: T \rightarrow T_0$, which is uniquely
determined by the requirement that it induces the map $T( - \infty ) \rightarrow T_0( - \infty)$ given
by the product of the maps
$$ \Rect( \Cube{k+1} \times \alpha^{-1} \{i\}, \Cube{k+1}) \rightarrow
\Conf( \alpha^{-1} \{i\}, \Cube{k+1} )$$
given by evaluation at the origin of $\Cube{k+1}$. It follows from
Remark \ref{sove} that the map $T( - \infty) \rightarrow T_0( - \infty)$ is a homotopy equivalence.
For $(I, \chi) \in P$, we have a commutative diagram
$$ \xymatrix{ T(I, \chi) \ar[dr]^{\phi} \ar[rr]^{\gamma_{(I,\chi)}} & & T_0(I, \chi) \ar[dl]^{\phi'} \\
& \prod_{i,j} \Conf( \alpha^{-1} \{i\} \times \chi^{-1} \{j\}, \Cube{k}) & }$$
The map $\phi'$ is a homotopy equivalence by virtue of $(\ast)$, and the map
$\phi$ is a homotopy equivalence by virtue of Remark \ref{sove}. It follows that
$\gamma_{(I, \chi)}$ is a homotopy equivalence, so that $\gamma$ is a weak equivalence of functors.

To prove that $T$ exhibits $T( - \infty)$ has a homotopy colimit of $T | P$, it will suffice
to show that $T_0$ exhibits $T_0( - \infty)$ as a homotopy colimit of $T_0 | P$. We will deduce
this by applying Theorem \ref{vankamp}. Fix a point $x \in T_0( - \infty )$, which we will
identify with a map $x: \nostar{m} \rightarrow \Cube{k+1}$ such that
$x(a) \neq x(b)$ if $\alpha(a) = \alpha(b)$. Let $x_0: \nostar{m} \rightarrow (-1,1)$
be the composition of $x$ with the projection to the first coordinate, and let
$x_1: \nostar{m} \rightarrow \Cube{k}$ be given by the projection to the remaining coordinates.
To apply Theorem \ref{vankamp}, we must
show that the partially ordered set $P_{x}$ has weakly contractible nerve, where
$P_{x}$ denotes the subset of $P$ spanned by those pairs $(I, \chi)$ satisfying
the following condition:
\begin{itemize}
\item[$(a)$] The subset $I \subseteq (-1,1)$ contains the set $x_0( \nostar{m} )$ in its interior.
\item[$(b)$] The map $\chi$ is characterized by the property that
$\chi(i) \in \pi_0 I$ is the connected component containing $x_0( \nostar{m})$.
\item[$(c)$] Given a pair of elements $a, b \in \nostar{m}$ such that
$\alpha(a) = \alpha(b)$ and $x_1(a) = x_1(b)$, the points
$x_0(a)$ and $x_0(b)$ belong to different connected components of $I$.
\end{itemize}

The contractibility of $\Nerve( P_x)$ follows from the observation that
$P_{x}^{op}$ is filtered: for any finite set $S \subset P_{x}$, there
exists another element $(I, \chi) \in P_{x}$ such that
$(I, \chi) \leq s$ for each $s \in S$. Indeed, we can take $I$ to be 
the union of closed intervals $\bigcup_{ a \in \nostar{m}} [ x_0(a) - \epsilon, x_0(a) + \epsilon ]$
for any sufficiently small positive real number $\epsilon$ (the map 
$\chi: \seg{m} \rightarrow \pi_0 I$ is uniquely determined by requirement $(b)$).
\end{proof}

\begin{lemma}\label{tople}
Let $X \rightarrow S$ be Cartesian fibration of simplicial sets, let
$S' \rightarrow S$ be a map such that $f^{op}$ is cofinal. Then the induced map $f': X' = X \times_{S} S' \rightarrow X$ is a weak homotopy equivalence.
\end{lemma}

\begin{proof}
It follows from Propositions \toposref{huneff} and \toposref{strokhop} that
${f'}^{op}$ is cofinal, and therefore a weak homotopy equivalence (Proposition \toposref{cofbasic}).
\end{proof}

\begin{proof}[Proof of Lemma \ref{cruss}]
We wish to show that the inclusion $(\calJ^{0})^{op} \rightarrow \calJ^{op}$ is cofinal (Proposition \toposref{cofbasic}). Let $\calI = \OpE{1} \wreath \OpE{k}$, and let $\overline{\calJ} = \calI \times_{ \OpE{k+1}} \OpE{k+1}^{\seg{m}/}$. It follows from Proposition \toposref{certs} that the natural map
$\calJ \rightarrow \overline{\calJ}$ is a categorical equivalence; we let $\overline{\calJ}^{0}$ denote the essential image of $\calJ^{0}$ in $\overline{\calJ}$. It will therefore suffice to show that the inclusion
$( \overline{\calJ}^{0})^{op} \subseteq \overline{\calJ}^{op}$ is cofinal. 
Using Theorem \toposref{hollowtt} and Proposition \toposref{certs}, we can reduce to showing
that for every object $\overline{X} \in \overline{\calJ}$, the $\infty$-category
$\overline{\calJ}^{0} \times_{ \overline{\calJ}} \overline{\calJ}^{/\overline{X}}$ is weakly contractible.

Let $X = ( \seg{n_1}, \seg{n_2}, \ldots, \seg{n_b})$ denote the image of $\overline{X}$ in $\calI$. 
We can identify $\overline{X}$ with a morphism $\alpha: \seg{m} \rightarrow \Psi_{k}(X)
= \seg{n_1 + \ldots + n_b}$ in $\OpE{k+1}$. Let $K(\alpha)$ be the $\infty$-category
$\overline{\calJ}^{0} \times_{ \overline{\calJ} } \overline{\calJ}^{/\overline{X}}$ by $\calK(\alpha)$.
We can identify objects of $\calK(\alpha)$ with triples $(Y, \gamma, \sigma)$, where
where $Y = ( \seg{n'_{1}}, \ldots, \seg{n'_{b'}} )$ is another object of $\calI$,
$\gamma: Y \rightarrow X$ is a morphism in $\calI$, and $\sigma$
is a commutative diagram
$$ \xymatrix{ \seg{m} \ar[d]^{\id} \ar[r]^{\beta} & \Psi_{k}(Y) \ar[d]^{ \Psi_{k}(\gamma)} \\
\seg{m} \ar[r]^{\alpha} & \Psi_{k}(X) }$$
in $\OpE{k+1}$, where $\beta$ is inert. The map $\gamma$ can in turn be identified
with a morphism $\gamma_0: \seg{b'} \rightarrow \seg{b}$ in $\OpE{1}$ together with
a collection of maps $\{ \gamma_i: \seg{b'_i} \rightarrow \seg{b_j} \}_{ i \in \gamma_0^{-1} \{j\} }$
in $\OpE{k}$. Let $\calK'(\alpha)$ denote the full subcategory of $\calK(\alpha)$ spanned
by those diagrams for which $\gamma_0$ and each of the morphisms $\gamma_i$ are active.
The inclusion $\calK'(\alpha) \subseteq \calK(\alpha)$ admits a left adjoint, and is therefore a weak homotopy equivalence. Let $\calK''(\alpha)$ denote the full subcategory of $\calK'(\alpha)$ spanned
by those objects for which each of the integers $n'_{i}$ is positive. The inclusion
$\calK''(\alpha) \subseteq \calK'(\alpha)$ admits a right adjoint, and is therefore a weak homotopy equivalence. It will therefore suffice to show that $\calK''(\alpha)$ is weakly contractible.

The map $\alpha$ factors as a composition
$$ \seg{m} \stackrel{ \alpha'}{\rightarrow} \seg{m'} \stackrel{\alpha''}{\rightarrow} \Psi_{k}(X)$$
where $\alpha'$ is inert and $\alpha''$ is active. Composition with $\alpha'$ induces
an equivalence of $\infty$-categories $\calK''(\alpha'') \rightarrow \calK''(\alpha)$.
We may therefore replace $\alpha$ by $\alpha''$ and thereby reduce to the case where
$\alpha$ is active. Then $\alpha$ determines a partition $m = m_1 + \ldots + m_b$ of
$m$ and a collection of active maps $\{ \alpha_{i}: \seg{m_i} \rightarrow \Psi_k(\seg{n_i}) \}_{1 \leq i \leq b}$ in $\OpE{k+1}$. The $\infty$-category $\calK''(\alpha)$ is equivalent to the direct
product $\prod_{ 1 \leq i \leq b} \calK''( \alpha_i)$. It now suffices to prove that each
$\calK''(\alpha_i)$ is weakly contractible: in other words, we may replace $\alpha$ by $\alpha_i$ and thereby reduce to the case where $X = ( \seg{n} )$ lies in the image of the embedding
$\OpE{k} \hookrightarrow \OpE{1} \wreath \OpE{k}$.

Let $\calA$ denote the $\infty$-category $\Nerve(\FinSeg)_{ \seg{m}/} \times_{ \Nerve(\FinSeg)} \Nerve(\FinSeg)^{\amalg}$,
whose objects are pairs consisting of an object $( \seg{a_1}, \ldots, \seg{a_{p}} ) \in \Nerve(\FinSeg)^{\amalg}$ together with a morphism $\beta: \seg{m} \rightarrow
\seg{ a_1 + \cdots + a_p}$ of pointed finite sets, and let $A = ( ( \seg{n}), \seg{m} \rightarrow \seg{n})$ denote the object of $\calA$ determined by $\alpha$. Let $\calA^{0}$ denote the full subcategory
of $\calA$ spanned by those objects for which each $a_i$ is positive and $\beta$ is an equivalence.
Let $\calI' = \calI \times_{ \Nerve(\FinSeg)^{\amalg} } \calA$, and let
$X'$ denote the object $(X, A)$ of $\calI'$. Note that the map
$\overline{\calJ} \rightarrow \calI'$ is a pullback of
$(\OpE{k+1})^{ \seg{m}/ } \rightarrow \OpE{k+1} \times_{ \Nerve(\FinSeg)} \Nerve(\FinSeg)^{ \seg{m}/}$, and therefore a left fibration. It follows that the fiber $\overline{\calJ}_{X'} = \overline{\calJ} \times_{ \calI'} \{X'\}$ is a Kan complex. 
Let $$\calC = \calA^{0} \times_{ \Fun( \{0\}, \calA)} \Fun( \Delta^1, \overline{\calJ})
\times_{ \Fun( \{1\}, \calI') } \{X' \}.$$
Evaluation at $\{1\}$ induces a categorical fibration
$$ \phi: \calC \rightarrow \overline{\calJ}_{X'}$$ 
such that $\phi^{-1} \{ \overline{X} \} = \calK(\alpha)$. 
Choose a contractible Kan complex $K$ and a Kan fibration
$K \rightarrow \overline{\calJ}_{X'}$ whose image
contains $\overline{X}$, and consider the diagram
$$ \xymatrix{ \calK''(\alpha) \ar[r] \ar[d] & K \times_{ \overline{\calJ}_{X'} } \calC \ar[r] \ar[d]^{\phi'} & \calC \ar[d]^{\phi} \\
\{ \overline{X} \} \ar[r] & K \ar[r] & \overline{\calJ}_{X'}. }$$
Since $\phi$ is a categorical fibration, every square in this diagram is a homotopy
pullback (with respect to the Joyal model structure). Because $K$ is
contractible, the left horizontal maps are categorical equivalences. Consequently, to prove
that $\calK''(\alpha)$ is weakly contractible, it suffices to show that the map $\phi'$ is a weak
homotopy equivalence. Because the usual model structure on simplicial sets
is right proper (and the map $K \rightarrow \overline{\calJ}_{X'}$ is a Kan fibration), we may reduce
to the problem of showing that $\phi: \calC \rightarrow \overline{\calJ}_{X'}$ is a weak homotopy equivalence. To prove this, we will need to make some auxiliary constructions.

Let $P$ be the partially ordered set of Construction \ref{tupus}, and let
$T: P^{\triangleleft} \rightarrow \Top$ be the contravariant functor described in
Construction \ref{swellman}. We define a topological category $\overline{P}$ as follows:
\begin{itemize}
\item[$(i)$] The set of objects of $\overline{P}$ is
$P \cup \{-\infty, \infty \} = P^{\triangleleft} \cup \{ \infty \}$. 
\item[$(ii)$] The mapping spaces in $\overline{P}$ are given by the formula
$$ \bHom_{ \overline{P} }(x,y) = \begin{cases} \ast & \text{ if } x = y = \infty \\
\ast & \text { if } x,y \in P^{\triangleleft}, x \leq y \\
T(x) & \text{ if } x \in P^{\triangleleft}, y = \infty \\
\emptyset & \text{ otherwise.} \end{cases}$$
\end{itemize}
Let $\overline{\calD}$ denote the fiber product
$\Nerve( \overline{P})^{/ \infty} \times_{ \Nerve(\overline{P})} \Nerve(P^{\triangleleft})$, and
let $\calD = \overline{\calD} \times_{ \Nerve( P^{\triangleleft} )} \Nerve(P)$. 

Let $\overline{P}_0$ be the full subcategory of $\overline{P}$ spanned by the objects
$\pm \infty$. There is an evident retraction $r$ of $\overline{P}$ onto $\overline{P}_0$, given on objects by the formula $$r(x) = \begin{cases} \infty & \text{ if } x = \infty \\
- \infty & \text{ otherwise. }\end{cases}$$
We also have a topological functor $\rho_0: \overline{P}_0 \rightarrow \TopE{k+1}$, given on objects by $$ \rho_0( - \infty ) = \seg{m} \quad \quad \rho_0( \infty ) = \seg{n}. $$
The composition $\rho_0 \circ r$ induces a map of $\infty$-categories
$\Nerve( \overline{P}) \rightarrow \OpE{k+1}$, which determines in turn a map
$$\overline{\calD} \rightarrow \{ \seg{m} \} \times_{ \Fun( \{0\}, \OpE{k+1}) }
\Fun( \Delta^1, \OpE{k+1}) \times_{ \Fun( \{1\}, \OpE{k+1})} \{ \seg{n} \}.$$
This map factors through $\overline{\calJ}_{X'} \subseteq \{ \seg{m} \} \times_{ \Fun( \{0\}, \OpE{k+1}) }
\Fun( \Delta^1, \OpE{k+1}) \times_{ \Fun( \{1\}, \OpE{k+1})} \{ \seg{n} \}$, and therefore
determines a map $\rho: \overline{\calD} \rightarrow \overline{\calJ}_{X'}$.

Note that we can identify objects of $\calD$ with triples $(I, \chi, \eta)$, where $(I, \chi) \in P$ and
$\eta \in T( I, \chi)$. Given such an object, write $I$ as a disjoint union of closed intervals
$I_1 \cup \ldots \cup I_b$, and choose an order-preserving bijection $\chi^{-1} [I_j] \simeq \nostar{m_j}$ for $1 \leq j \leq b$. The point $\eta$ can be regarded as a morphism $\gamma$ from
the object $Y = ( \seg{m_1}, \ldots, \seg{m_b})$ to $( \seg{n} )$ in $\calI$. This morphism fits
into a commutative diagram $\sigma$
$$ \xymatrix{ \seg{m} \ar[r] \ar[d]^{\id} & \Psi_k(Y) \ar[d]^{ \Psi_k(\gamma)} \\
\seg{m} \ar[r]^{\eta} & \seg{n} }$$
in the $\infty$-category $\OpE{k+1}$. The construction $(I, \chi, \eta) \mapsto (Y, \gamma, \sigma)$
extends to a functor $\rho': \calD \rightarrow \calC$. We have a commutative diagram
$$ \xymatrix{ \calD \ar[r] \ar[d]^{\rho' } & \overline{\calD} \ar[d]^{\rho} \\
\calC \ar[r]^{\phi} & \overline{\calJ}_{X'}. }$$
We wish to prove that $\phi$ is a weak homotopy equivalence. We will complete the proof by verifying the following:
\begin{itemize}
\item[$(a)$] The inclusion $\calD \subseteq \overline{\calD}$ is a weak homotopy equivalence.
\item[$(b)$] The map $\rho$ is a weak homotopy equivalence.
\item[$(c)$] The map $\rho'$ is a weak homotopy equivalence.
\end{itemize}

To prove $(a)$, we observe that the right fibration
$q: \overline{\calD} \rightarrow \Nerve(P^{\triangleleft}) \simeq \Nerve(P)^{\triangleleft}$, which is equivalent to the right fibration
$\Nerve( \overline{P})_{/ \infty} \times_{ \Nerve( \overline{P})} \Nerve( P^{\triangleleft}) \rightarrow \Nerve(P^{\triangleleft})$ (Proposition \toposref{certs}) obtained by applying the 
unstraightening functor
$\Un$ of \S \toposref{rightstraight} to the functor $(\Sing_{\bigdot} \circ T): (P^{\triangleleft})^{\op} \rightarrow \sSet$. In other words, $q$ is the right fibration associated to the functor
$\Nerve(P^{op})^{\triangleright} \rightarrow \SSet$ given by the nerve of
$\Sing_{\bigdot} \circ T$. Lemma \ref{cisgn} implies that $\Sing_{\bigdot} \circ T$ is a homotopy
colimit diagram, so that $\Nerve(P^{op})^{\triangleright} \rightarrow \SSet$ is a colimit diagram
(Theorem \toposref{colimcomparee}). Applying Proposition \toposref{charspacecolimit}, we deduce that the inclusion $\calD \subseteq \overline{\calD}$ is a weak homotopy equivalence.

To prove $(b)$, let $\overline{\calD}_{ - \infty }$ denote the fiber product $\overline{\calD} \times_{ \Nerve( P^{\triangleleft}) } \{ - \infty \}$, and consider the diagram
$$ \xymatrix{ & \overline{\calD} \ar[dr]^{\rho} & \\
\overline{\calD}_{ - \infty} \ar[ur]^{j} \ar[rr]^{\rho''} & & \overline{\calJ}_{X'}. }$$
Since $\overline{\calD} \rightarrow \Nerve( P^{\triangleleft})$ is a right fibration and
$- \infty$ is an initial object of $\Nerve( P^{\triangleleft})$, the inclusion $j$ is a weak homotopy equivalence (Lemma \ref{tople}). The map $\rho''$ is a homotopy equivalence because both its domain and codomain can be identified with the summand $$\bHom_{ \OpE{k+1}}( \seg{m}, \seg{n} )
\times_{ \Hom_{\FinSeg}( \seg{m}, \seg{n})} \{ \alpha_0 \} \subseteq \bHom_{ \OpE{k+1}}( \seg{m}, \seg{n}),$$
where $\alpha_0$ denotes the morphism $\seg{m} \rightarrow \seg{n}$ in $\FinSeg$ determined by $\alpha$. It follows from the two-out-of-three property that $\rho$ is a weak homotopy equivalence as required.

To prove $(c)$, let $\calC_0$ be defined as in Lemma \ref{tupus}, and observe that
we have a commutative diagram
$$ \xymatrix{ \calD \ar[r]^{ \rho'} \ar[d] & \calC \ar[d]^{\psi} \\
\Nerve(P) \ar[r]^{\theta} & \calC_0. }$$
We will prove:
\begin{itemize}
\item[$(c')$] The map $\psi$ is a right fibration.
\item[$(c'')$] For every object $x \in P$, the induced map
$\calD \times_{ \Nerve(P) } \{x\} \rightarrow \calC \times_{ \calC_0} \{\theta(x)\}$ is a weak homotopy equivalence.
\end{itemize}
Then condition $(c'')$ and Lemma \ref{buz} guarantee that the map
$\calD \rightarrow \calC \times_{ \calC_0} \Nerve(P)$ is a weak homotopy equivalence,
while condition $(c')$, Lemma \ref{tople}, and Lemma \ref{colman} guarantee that
$\calC \times_{ \calC_0} \Nerve(P) \rightarrow \calC$ is a weak homotopy equivalence. It follows that the composition $\rho': \calD \rightarrow \calC$ is a weak homotopy equivalence, as required by $(c)$.

We now prove $(c')$.
Let $\calI_0 = \OpE{1} \wreath \Nerve(\FinSeg)$, let $\calI'_0 = \calI_0 \times_{ \Nerve(\FinSeg)^{\amalg} } \calA$, and let $X'_0$ denote the image of
$X'$ in $\calI'_0$. We note that $\psi$ factors as a composition
$$ \calC \stackrel{\psi_0}{\rightarrow} (\calA^{0} \times_{\calA} \overline{\calJ})
\times_{ \calI'} (\calI')^{/X_0} \stackrel{\psi_1}{\rightarrow}
\calA^{0} \times_{\calA} (\calI')^{/X'} \stackrel{\psi_2}{\rightarrow}
\calA^{0} \times_{\calA} (\calI'_0)^{/ X'_0}.$$
Since $\overline{\calJ} \rightarrow \calI'$ is a left fibration and the inclusion $\{0\} \subseteq \Delta^1$
is left anodyne, the map $\psi_0$ is a trivial Kan fibration. The map $\psi_1$ is also a trivial
Kan fibration, because the forgetful functor $\OpE{k+1} \rightarrow \Nerve(\FinSeg)$
induces a trivial Kan fibration between the underlying Kan complexes.
The map $\psi_2$ factors as a composition
$$ \calA^{0} \times_{\calA} (\calI')^{/X'}
\stackrel{\psi'_2}{\rightarrow} \calA^{0} \times_{\calA}
( \calI' \times_{ \calI'_0} (\calI'_0)^{/X'_0}) \stackrel{\psi''_{2}}{\rightarrow}
\calC_0.$$
Here $\psi'_2$ is a pullback of 
$\Fun( \Delta^1, \calI') \rightarrow \Fun( \Delta^1, \calI'_0) \times_{ \Fun( \{1\}, \calI'_0) }
\Fun( \{1\}, \calI')$ (which is a right fibration since the inclusion $\{1\} \subseteq \Delta^1$ is right anodyne and the map $\calI' \rightarrow \calI'_0$ is a categorical fibration), and
$\psi''_2$ is a pullback of $\calA^{0} \times_{ \Nerve(\FinSeg)^{\amalg} } \OpE{k}^{\amalg}
\rightarrow \calA^{0}$ (which is a trivial Kan fibration, since the map
$\OpE{k} \rightarrow \Nerve(\FinSeg)$ induces a trivial Kan fibration over the full subcategory
of $\FinSeg$ spanned by the injective maps).

To prove $(c)$, consider an object $x = ( I_1 \cup \ldots \cup I_b, \chi)$, and let
$Y = ( \seg{m_1}, \ldots, \seg{m_b})$ be the corresponding object of
$\calI$. Both $\calD_{x}$ and $\calC_{ \theta(x)}$ are homotopy equivalent to
the summand $\bHom_{ \calI}( Y, \seg{n} ) \times_{ \Hom_{ \FinSeg}}( \seg{m}, \seg{n})
\{ \alpha_0 \} \subseteq \bHom_{ \calI}( Y, \seg{n}).$
\end{proof}

\begin{proof}[Proof of Proposition \ref{inku}]
We use the strategy of Proposition \symmetricref{algass}. Fix a fibration of
$\infty$-operads $\calC^{\otimes} \rightarrow \OpE{k+1}$, and let
$\calZ$ denote the full subcategory of $\Fun_{ \OpE{k+1}}( \OpE{1} \wreath
\OpE{k}, \calC^{\otimes})$ spanned by those functors which carry inert morphisms
to inert morphisms; we will prove that composition with $\Psi_{k}$ induces
an equivalence of $\infty$-categories $\theta: \Alg_{ \OpE{k+1}}(\calC) \rightarrow \calZ$.

Let $\calK$ denote the mapping cylinder of the functor
$\Psi_k: ( \OpE{1} \wreath \OpE{k}) \rightarrow \OpE{k+1}$: that is,
we let $\calK$ be a simplicial set equipped with a map $\calK \rightarrow \Delta^1$
satisfying the following universal mapping property: given any map of simplicial
sets $K \rightarrow \Delta^1$, the set $\Hom_{ \Delta^1}( K, \calK)$ can be identified
with the set of commutative diagrams
$$ \xymatrix{ K \times_{ \Delta^1} \{0\} \ar[d] \ar[r] & \OpE{1} \wreath \OpE{k} \ar[d]^{\Psi_{k}} \\
K \ar[r] & \OpE{k+1}. }$$
Then $\calK \rightarrow \Delta^1$ is a correspondence from
$\calK_{0} \simeq \OpE{k+1}$ to $\calK_{1} \simeq \OpE{1} \wreath \OpE{k}$.
There is a canonical retraction $r$ of $\calK$ onto $\OpE{k+1}$.

Let $\calX$ denote the full subcategory of $\Fun_{ \OpE{k}}( \calK, \calC^{\otimes})$ spanned by those functors $A$ which satisfy the following pair of conditions:
\begin{itemize}
\item[$(i)$] The restriction of $A | \OpE{k}$ belongs to $\Alg_{ \OpE{k}}(\calC)$.
\item[$(ii)$] For every object $X \in \OpE{1} \wreath \OpE{k}$, the canonical map
$A( \Psi_k X) \rightarrow A(X)$ is an equivalence in $\calC^{\otimes}$.
\end{itemize}
Note that condition $(ii)$ is equivalent to the requirement that $A$ is a left Kan extension of
$A | \OpE{k}$. It follows from Proposition \toposref{lklk} that the restriction map
$\beta: \calX \rightarrow \Alg_{ \OpE{k} }( \calC)$ is a trivial Kan fibration.
The functor $\theta: \Alg_{ \OpE{k}}(\calC) \rightarrow \calZ$ factors as a composition
$$ \Alg_{ \OpE{k}}(\calC) \stackrel{\theta'}{\rightarrow} \calX
\stackrel{ \theta''}{\rightarrow} \calZ,$$
where $\theta'$ is a section of the trivial Kan fibration $\beta$ (and therefore a categorical
equivalence). It will therefore suffice to show that the map $\theta''$ is a trivial Kan fibration.
In view of Proposition \toposref{lklk}, it will suffice to verify the following:
\begin{itemize}
\item[$(a)$] Let $A \in \Fun_{ \OpE{k}}(\calK, \calC^{\otimes})$ be a functor
such that $A_0 = A | ( \OpE{1} \wreath \OpE{k})$ belongs to
$\calZ$. Then $A \in \calX$ if and only if and $A$ is a $p$-right Kan extension of $A_0$.
\item[$(b)$] Every object $A_0 \in \calZ$ admits a $p$-right
Kan extension $A \in \Fun_{ \OpE{k}}( \calK, \calC^{\otimes})$.
\end{itemize}

To prove $(a)$, fix an object $A \in \Fun_{ \OpE{k}}( \calK, \calC^{\otimes})$
such that $A_0 \in \calZ$.
Let $\seg{m} \in \OpE{k}$, and let $\calJ^{0} \subseteq \calJ$ be the inclusion
described in Lemma \ref{cruss}. Using Lemma \ref{cruss}, we deduce:

\begin{itemize}
\item[$(a_0)$] The functor $A$ is a $p$-right Kan extension of $A_0$ at
$\seg{m}$ if and only if the map
$$ \overline{f}^0: ( \calJ^{0})^{\triangleleft} \rightarrow
\calK \stackrel{A}{\rightarrow} \calC^{\otimes}$$
is a $p$-limit diagram.
\end{itemize}

Let $\calJ^{1}$ denote the full subcategory of
$\calJ^{0}$ spanned by those inert morphisms
$\seg{m} \rightarrow \Psi_{k}(X)$ in $\OpE{k+1}$ for which $X$ lies over the object
$( \seg{1} )$ of $\Nerve(\FinSeg)^{\amalg}$. We next claim:

\begin{itemize}
\item[$(a_1)$] The map $f^0 = \overline{f}^0 | \calJ^{0}$
is a $p$-right Kan extension of $f^1 = \overline{f}^0 | \calJ^1$. 
Consequently (by virtue of Lemma \toposref{kan0} and $(a_0)$)
the functor $A$ is a $p$-right Kan extension of $A_0$ at $\seg{m}$ if and only if
the functor $\overline{f}^1 = \overline{f}^0 | (\calJ^1)^{\triangleleft}$
is a $p$-limit diagram.
\end{itemize}

To prove $(a_1)$, consider an arbitrary object $\alpha: \seg{m} \rightarrow \Psi_{k}(X)$
of $\calJ^0$, where $X$ lies over the object $( \seg{n_1}, \ldots, \seg{n_j}) \in 
\Nerve(\FinSeg)^{\amalg}$. Let $\calD$ denote the $\infty$-category
$\calJ^{0}_{\alpha/} \times_{\calJ^{0}} \calJ^{1}$. Every object
$D \in \calD$ determines an inert morphism $\seg{n_1 + \cdots + n_j} \rightarrow
\seg{1}$ in $\FinSeg$, which we can identify with an element
$i_{D} \in \nostar{ n_1 + \cdots + n_j }$. The assignment $D \mapsto i(D)$ determines
a decomposition of $\calD$ as a disjoint union
$\calD \simeq \coprod_{ i \in \nostar{ n_1 + \cdots + n_j}} \calD_{i}$,
where each $\calD_{i}$ is a contractible Kan complex containing a vertex
$D_i$, which induces an inert morphism $\beta_{i}: X \rightarrow X_i$ in $\OpE{1} \wreath \OpE{k}$. 
It follows that $f^0$ is a $p$-right Kan extension of $f^1$ at
$\alpha$ if and only if $A_0$ exhibits $A_0(X)$ as a $p$-product of the objects
$\{ A_0(X_i) \}_{1 \leq i \leq n_1 + \cdots + n_j}$, which follows from our assumption
that $A_0 \in \calZ$.

Every object of $\calJ^1$ determines an inert morphism $\seg{m} \rightarrow \seg{1}$
in $\Nerve(\FinSeg)$, which we can identify with an element $i \in \nostar{m}$.
This assignment determines a decomposition of $\calJ^1$ as a disjoint union
$\calJ^1 \simeq \coprod_{1 \leq i \leq m} \calJ^{1}_{i}$.
Let $X_0$ denote the unique vertex of $\OpE{1} \wreath \OpE{k}$ lying
over the object $( \seg{1} ) \in \Nerve(\FinSeg)^{\amalg}$. 
Each of the $\infty$-categories $\calJ^{1}_{i}$ is a contractible Kan complex
containing a vertex $\alpha_{i}: \seg{m} \rightarrow \Psi_k( X_0)$. Combining
this observation with $(a_1)$, we deduce:

\begin{itemize}
\item[$(a_2)$] The functor $A$ is a $p$-right Kan extension of $A_0$
at $\seg{m}$ if and only if the morphisms $A(\alpha_{i})$ exhibit
$A( \seg{m} )$ as a $p$-product of the objects $\{ A_0( X_0) \}_{1 \leq i \leq m}$.
\end{itemize}

Using the fact that $p$ is a fibration of $\infty$-operads and allowing the integer $m$ to vary, we
deduce the following version of $(a)$:

\begin{itemize}
\item[$(a_3)$] The functor $A$ is a $p$-right Kan extension of $A_0$ if and only if,
for every pair of integers $1 \leq i \leq m$, the morphism $\alpha_{i}$ described
above determines an inert morphism $A(\alpha_i):
A( \seg{m}) \rightarrow A_0(X_0)$. 
\end{itemize}

We now prove $(a)$. Assume first that $A$ is a $p$-right Kan extension of $A_0$;
we will show that $A$ satisfies conditions $(i)$ and $(ii)$. To prove
$(i)$, we must show that if $\gamma: \seg{m} \rightarrow \seg{m'}$ is an
inert morphism in $\OpE{k}$, then $A(\gamma)$ is an inert morphism in $\calC^{\otimes}$.
Using Remark \symmetricref{casper}, we can reduce to the case where $m' = 1$.
We have a commutative diagram
$$ \xymatrix{ A( \seg{m} ) \ar[rr]^{ A(\gamma)} \ar[dr]^{A(\gamma')} & & A( \seg{1}) \ar[dl]^{A(\gamma'')} \\
& A_0(X_0). & }$$
Applying $(a_3)$, we deduce that $A(\gamma')$ and $A(\gamma'')$ are inert; in particular,
$A( \gamma'')$ is an equivalence so that $A(\gamma)$ is inert as desired.

We now prove $(ii)$. Fix an object $X \in \OpE{1} \wreath \OpE{k}$, and let 
$\seg{m} = \Psi_k(X) \in \OpE{k+1}$. We wish to
show that the canonical map $\gamma: A( \seg{m} ) \rightarrow A_0(X)$ is
an equivalence. For $1 \leq i \leq m$, let $\beta_i: X \rightarrow X_i$ be defined as above,
so that our assumption that $A_0 \in \calZ$ guarantees that $A_0( \beta_i)$ is inert.
Since $\calC^{\otimes}$ is an $\infty$-operad, it will suffice to show that
each composition $A( \beta_i \circ \gamma)$ is inert, which follows immediately
from $(a_3)$, since $\beta_{i} \circ \gamma$ is equivalent to the morphism
$\alpha_i: \seg{m} \rightarrow \Phi_{k}(X_0)$.

We now prove the converse: suppose that the functor $A$ satisfies conditions
$(i)$ and $(ii)$; we wish to prove that $A$ is a $p$-right Kan extension of
$A_0$. In view of $(a_3)$, it suffices to show that for $1 \leq i \leq m$, the map
$A(\alpha_i): A( \seg{m}) \rightarrow A_0( X_0)$ is inert.
This map factors as a composition
$$ A( \seg{m}) \stackrel{ \rho}{\rightarrow} A( \seg{1}) \stackrel{ \rho'}{\rightarrow}
A_0(X_0)$$
where $\rho$ is inert by virtue of $(i)$ and $\rho'$ is an equivalence by virtue of $(ii)$.
This completes the proof of $(a)$.

We now prove $(b)$. Fix $A_0 \in \calZ$. To prove that
$A_0$ admits a $p$-right Kan extension $A \in \Fun_{ \OpE{k}}( \calK, \calC^{\otimes})$,
it suffices to show that for every object $\seg{m} \in \OpE{k}$, the composite diagram
$$f: \calJ \rightarrow \OpE{1} \wreath \OpE{k} \stackrel{A_0}{\rightarrow}
\calC^{\otimes}$$
can be extended to a $p$-limit diagram $\overline{f} \in \Fun_{ \OpE{k}}(
\calJ^{\triangleleft}, \calC^{\otimes})$. Since
the inclusion $\calJ^{0} \subseteq \calJ$
is left anodyne, it suffices to extend $f^0 = f | \calJ^{0}$ to
a $p$-limit diagram. Since $f^0$ is a $p$-right Kan extension of $f^1 =
f| \calJ^1$, we can reduce to showing that
$f^1$ can be extended to a $p$-limit diagram (Lemma \toposref{kan0}).
Since $\calJ^{1}$ is equivalent to the discrete simplicial set $\nostar{m}$, we are reduced
to showing that a map $\nostar{m} \rightarrow \calC$ can be extended to a $p$-limit diagram
in $\Fun_{ \OpE{k+1}}( {\nostar{m}}^{\triangleleft}, \calC^{\otimes})$, which follows immediately from our assumption that $p$ is a fibration of $\infty$-operads. 
\end{proof}

\subsection{Iterated Loop Spaces}\label{kloop}

Let $X$ be a topological space equipped with a base point $\ast$, and let
$\Top$ denote the category of topological spaces. For each
$n \geq 0$, let $\theta_X( \seg{n} ) \simeq (\Omega^{k} X)^n$ denote the collection of 
$n$-tuples of maps $f_1, \ldots, f_n: [-1,1]^{k} \rightarrow X$
such that each $f_{i}$ carries the boundary of $[-1,1]^{k}$ to the base point
$\ast \in X$. The construction $\seg{n} \rightarrow \Sing_{\bigdot} \theta_X(\seg{n})$ determines a simplicial functor $\theta_{X}: \Sing \TopE{k} \rightarrow \Kan$ (which encodes the idea that the iterated loop space $\Omega^{k} X$ is acted on by the little cubes operad $\calC_{k}$). This construction depends functorially on $X$. Restricting our attention to the case where $X = |K|$, where $K$ is a
(pointed) Kan complex, we obtain a simplicial functor
$$ \Kan_{\ast/} \times \Sing \TopE{k} \rightarrow \Kan.$$
Passing to nerves and using the evident equivalence $\Nerve(\Kan_{\ast/}) \rightarrow
\SSet_{\ast}$, we obtain a functor
$$ \Nerve(\theta): \SSet_{\ast} \times \OpE{k} \rightarrow \SSet.$$
For every pointed space $K$, the resulting map $\OpE{k} \rightarrow \SSet$ is evidently
an $\OpE{k}$-monoid object of $\SSet$ (in the sense of Definition \symmetricref{gammaobj}).
Consequently, $\Nerve(\theta)$ is adjoint to a functor
$\beta: \SSet_{\ast} \rightarrow \Mon_{ \OpE{k}}( \SSet).$
We will refer to $\OpE{k}$-monoid objects of $\SSet$ simply as {\it $\OpE{k}$-spaces}, and
$\Mon_{ \OpE{k}}(\SSet)$ as the {\it $\infty$-category of $\OpE{k}$-spaces}.
The functor $\beta$ implements the observation that for every pointed space $X$,
the $k$-fold loop space of $X$ is an $\OpE{k}$-space. This observation has a converse:
the functor $\beta$ is {\em almost} an equivalence of $\infty$-categories. However, it fails to be an equivalence for two reasons:

\begin{itemize}
\item[$(a)$] If $X$ is a pointed space, then the $k$-fold loop space
$\Omega^{k} X$ contains no information about the homotopy groups
$\pi_{i} X$ for $i < k$. More precisely, if $f: X \rightarrow Y$ is a map of pointed spaces
which induces isomorphisms $\pi_{i} X \rightarrow \pi_{i} Y$ for $i \geq k > 0$, then
the induced map $\Omega^{k} X \rightarrow \Omega^{k} Y$ is a weak homotopy equivalence
of spaces (which underlies a weak homotopy equivalence of $\OpE{k}$-monoids). Consequently,
the functor $\beta: \SSet_{\ast} \rightarrow \Mon_{ \OpE{k}}( \SSet)$ fails to be conservative.
To correct this problem, we need to restrict our attention to {\it $k$-connective} spaces:
that is, pointed spaces $X$ such that $\pi_{i} X \simeq \ast$ for $i < k$; in this case, there is no information about low-dimensional homotopy groups to be lost.

\item[$(b)$] Suppose that $k > 0$ and let $Y \in \Mon_{ \OpE{k}}(\SSet)$; we will abuse notation by identifying $Y$ with the space $Y(\seg{1})$. Then $Y$ carries an action of the $\infty$-operad $\OpE{1}$: in particular, there is a multiplication map $Y \times Y \rightarrow Y$ which is unital and associative up to homotopy. 
This multiplication endows the set of connected components $\pi_0 Y$ with the structure of a monoid
(which is commutative if $k > 1$). If $Y \simeq \Omega^{k} X$ lies in the image of
the functor $\beta$, then we have a canonical isomorphism $\pi_0 Y \simeq \pi_{k} X$
(compatible with the monoid structures on each side). In particular, we deduce that the monoid
$\pi_0 Y$ is actually a group (that is, $Y$ is {\it grouplike} in the sense of Definition \ref{ungwar} below).
\end{itemize}

\begin{remark}
In the case $k=0$, issues $(a)$ and $(b)$ do not arise: in fact, we have canonical
equivalences of $\infty$-categories
$$ \SSet_{\ast} \simeq \Alg_{ \OpE{0}}( \SSet) \simeq \Mon_{ \OpE{0}}( \SSet)$$
(here we regard $\SSet$ as endowed with the Cartesian monoidal structure).
The first equivalence results from Proposition \symmetricref{ezalg}, and the second from
Proposition \symmetricref{ungbatt}; the composition of these equivalences agrees with the map $\beta$ defined above. For this reason, we will confine our attention to the case $k > 0$ in what follows.
\end{remark}

We now introduce some terminology to address objection $(b)$.

\begin{definition}\label{ungwar}
Let $\calX$ be an $\infty$-topos, and let $\phi: \Nerve(\cDelta)^{op} \rightarrow \Ass$ be defined as in
Construction \symmetricref{urpas}. We will say that an $\Ass$-monoid object $X: \Ass \rightarrow \calX$ is
{\it grouplike} if the composition
$$ \Nerve(\cDelta)^{op} \rightarrow \Ass \stackrel{X}{\rightarrow} \calX$$
is a groupoid object of $\calX$ (see \S \toposref{gengroup}). Let $\Mon_{\Ass}^{\glike}(\calX)$ be the full subcategory of $\Mon_{\Ass}(\calX)$ spanned by the grouplike $\Ass$-monoid objects of $\calX$.

We will say that an $\OpE{1}$-monoid object $X: \OpE{1} \rightarrow \calX$ is {\it grouplike}
if it belongs to the essential image of $\Mon_{\Ass}^{\glike}(\calX)$ under the equivalence of
$\infty$-categories $\Mon_{\Ass}(\calX) \rightarrow \Mon_{\OpE{1}}(\calX)$ induced by
the equivalence $\OpE{1} \rightarrow \Ass$. We let $\Mon_{ \OpE{1}}^{\glike}(\calX)
\subseteq \Mon_{\OpE{1}}(\calX)$ denote the full subcategory spanned by the grouplike
$\OpE{1}$-monoid objects of $\calX$.

If $k > 0$, then we will say that an $\OpE{k}$-monoid object $X: \OpE{k} \rightarrow \calX$
is {\it grouplike} if the composite map $\OpE{1} \hookrightarrow \OpE{k} \stackrel{X}{\rightarrow} \calX$
is an grouplike $\OpE{1}$-monoid object of $\calX$. We let $\Mon_{ \OpE{k}}^{\glike}(\calX) \subseteq
\Mon_{ \OpE{k}}(\calX)$ denote the full subcategory spanned by the grouplike $\OpE{k}$-monoid objects. 
\end{definition}

\begin{remark}
Let $X$ be an $\OpE{k}$-monoid object of an $\infty$-topos $\calX$, for $k > 0$, and let
us abuse notation by identifying $X$ with the underlying object $X( \seg{1}) \in \calX$.
Then $X$ is equipped with a multiplication map $X \times X \rightarrow X$, which is associative
up to homotopy. Using Lemma \toposref{slurpy}, we deduce that $\tau_{\leq 0} X$ is
an associative monoid in the ordinary topos $\h{(\tau_{\leq 0} \calX)}$. Unwinding the definitions,
we see that $X$ is grouplike if and only if $\tau_{\leq 0} X$ is a group object of
$\h{ \tau_{\leq 0} \calX}$. In particular, the condition that $X$ is grouplike does not depend on which of the natural embeddings $\OpE{1} \hookrightarrow \OpE{k}$ is chosen.
\end{remark}

\begin{remark}\label{splish}
Let $X$ be an $\OpE{k}$-monoid object of an $\infty$-topos $\calX$ for $k > 0$. If
$X$ is $1$-connective, then $\tau_{\leq 0} X$ is a final object of $\h{ \tau_{\leq 0} \calX}$, so that $X$ is grouplike.
\end{remark}

\begin{remark}\label{jazwind}
An $\OpE{k}$-monoid object $X$ of the $\infty$-category $\SSet$ of spaces is {\it grouplike} if and only if
the monoid $\pi_0 X$ is a group. Note that the truncation functor $X \rightarrow \pi_0 X$
preserves colimits. Since the category of (commutative) groups is stable under colimits in the larger category of (commutative) monoids, we deduce that $\Mon_{ \OpE{k}}^{\glike}(\SSet)$ is stable
under small colimits in $\Mon_{ \OpE{k}}( \SSet)$.
\end{remark}

We now prove an abstract version of our main result:

\begin{theorem}\label{preslage}
Let $k > 0$, let $\calX$ be an $\infty$-topos, and let
$\calX_{\ast}^{\geq k}$ denote the full subcategory of $\calX_{\ast}$ spanned by those pointed spaces which are $k$-connective. Then there is a canonical equivalence of $\infty$-categories $\alpha: \calX_{\ast}^{\geq k} \simeq \Mon_{ \OpE{k}}^{\glike}( \calX)$.
\end{theorem}

\begin{proof}
We first observe that the forgetful functor $\Mon_{ \OpE{k}}( \calX_{\ast}) \rightarrow
\Mon_{ \OpE{k}}( \calX)$ is an equivalence of $\infty$-categories.
This map fits into a commutative diagram
$$ \xymatrix{ \Alg_{ \OpE{k}}(\calX_{\ast}) \ar[r] \ar[d] & \Alg_{ \OpE{k}}( \calX) \ar[d] \\
\Mon_{ \OpE{k}}( \calX_{\ast}) \ar[r] & \Mon_{ \OpE{k}}( \calX). }$$
The vertical maps are categorical equivalences by Proposition \symmetricref{ungbatt}, so
we are reduced to proving that the upper horizontal map is a categorical equivalence.
The $\infty$-category $\Alg_{ \OpE{k}}( \calX_{\ast})$ is equivalent to
$\Alg_{ \OpE{k}}( \calX)_{ {\bf 1}/}$, where ${\bf 1}$ is a trivial $\OpE{k}$-algebra in $\calX$;
since $\OpE{k}$ is unital, it follows from Proposition \symmetricref{gargle1} that the
forgetful functor $\Alg_{ \OpE{k}}( \calX)_{ {\bf 1}/} \rightarrow \Alg_{ \OpE{k}}( \calX)$
is an equivalence of $\infty$-categories. It will therefore suffice to construct
an equivalence $\alpha': \calX^{\geq k}_{\ast} \rightarrow \Mon_{ \OpE{k}}( \calX)$.
The construction proceeds by recursion on $k$. Suppose first that $k > 0$, so we can write
$k = k_{-} + k_{+}$ where $0 < k_{-}, k_{+} < k$.
The inductive hypothesis guarantees the
existence of equivalences $\alpha'_{-}: \calX_{\ast}^{\geq k_{-}} \rightarrow
\Mon_{ \OpE{k_{-}}}^{\glike}( \calX_{\ast})$ and $\alpha'_{+}: \calX_{\ast}^{\geq k_{+}} \rightarrow
\Mon_{ \OpE{k_{+}}}^{\glike}( \calX_{\ast})$. Note that a pointed object $X$ of $\calX$ is
$k$-connective if and only if $\Omega^{k''} X$ is $k'$-connective; moreover, any $k'$-connective
$\OpE{k''}$-monoid object of $\calX_{\ast}$ is automatically grouplike by Remark \ref{splish}.
It follows that $\alpha'_{+}$ and $\alpha'_{-}$ induce an equivalence
$$\gamma: \calX_{\ast}^{\geq k} \simeq \Mon_{ \OpE{k_{+}}}( \Mon_{ \OpE{k_{-}}}( \calX_{\ast})).$$
Let $\delta: \Mon_{ \OpE{k}}(\calX_{\ast}) \rightarrow \Mon_{ \OpE{k_{+}}}( \Mon_{\OpE{k_{-}}}(\calX_{\ast}))$ be the map induced by the $\infty$-operad bifunctor
$\OpE{k_{-}} \times \OpE{k_{+}} \rightarrow \OpE{k}$. Then $\delta$ fits into a commutative diagram
$$ \xymatrix{ \Alg_{ \OpE{k}}( \calX_{\ast}) \ar[r]^-{\delta'} \ar[d] & \Alg_{ \OpE{k_{+}}}( \Alg_{\OpE{k_{-}}}(\calX_{\ast})) \ar[d] \\
\Mon_{ \OpE{k}}( \calX_{\ast}) \ar[r]^-{\delta} & \Mon_{ \OpE{k_{+}}}( \Mon_{ \OpE{k_{-}}}( \calX_{\ast})). }$$
The vertical maps are categorical equivalences by Proposition \symmetricref{ungbatt}, and the
map $\delta'$ is a categorical equivalence by virtue of Theorem \ref{slide}; it follows that
$\delta$ is likewise a categorical equivalence. Let $\delta^{-1}$ be a homotopy inverse to
$\delta$. We now complete the proof by setting $\alpha' = \delta^{-1} \circ \gamma$.

It remains to treat the case where $k=1$.
Let $\calC$ denote the full subcategory of $\Fun( \Nerve(\cDelta_{+}^{op}, \calX)$
spanned by those augmented simplicial objects $X_{\bigdot}$ satisfying the following
conditions:
\begin{itemize}
\item[$(i)$] The underlying map $f: X_0 \rightarrow X_{-1}$ is an effective epimorphism in $\calX$.
\item[$(ii)$] The augmented simplicial object $X_{\bigdot}$ is a \Cech nerve of $f$.
\item[$(iii)$] The object $X_0 \in \calX$ is final.
\end{itemize}
Using Proposition \toposref{lklk}, we deduce that the construction $X_{\bigdot} \mapsto f$
determines a trivial Kan fibration from $\calC$ to the full subcategory of $\Fun( \Delta^1, \calX)$
spanned by those morphisms $f: X_0 \rightarrow X_{-1}$ where $X_0$ is a final object of
$\calX$ and $f$ is an effective epimorphism. Let $\phi_0: \calX_{\ast}^{\geq 1} \rightarrow \calC$
be a section of this trivial Kan fibration. Let $\Mon(\calX) \subseteq \Fun( \Nerve(\cDelta^{op}), \calX)$
be the $\infty$-category of monoid objects of $\calX$, and $\Mon^{\glike}(\calX) \subseteq \Mon(\calX)$ the full subcategory spanned by the grouplike monoid objects. Since $\calX$ is an $\infty$-topos,
the restriction map $\phi_1: \calC \rightarrow \Mon^{\glike}(\calX)$ is an equivalence of $\infty$-categories. Using Propositions \monoidref{ungbat}, \symmetricref{ungbatt}, and \symmetricref{algass}, we deduce that the restriction functor $\Mon_{\Ass}^{\glike}(\calX) \rightarrow \Mon^{\glike}(\calX)$ is
an equivalence of $\infty$-categories which admits a homotopy inverse $\phi_2$. Let
$\phi_3: \Mon_{ \Ass}^{\glike}(\calX) \rightarrow \Mon_{ \OpE{1}}^{\glike}(\calX)$ be the equivalence of
$\infty$-categories induced by the categorical equivalence $\OpE{1} \rightarrow \Ass$ of Example \ref{sulta}. We now define $\alpha$ to be the composite equivalence
$$ \calX_{\ast}^{\geq 1} \stackrel{\phi_0}{\rightarrow} \calC
\stackrel{ \phi_1}{\rightarrow} \Mon^{\glike}(\calX) \stackrel{\phi_2}{\rightarrow}
\Mon^{\glike}_{\Ass}(\calX) \stackrel{\phi_3}{\rightarrow} \Mon^{\glike}_{\OpE{1}}(\calX).$$
\end{proof}

\begin{corollary}\label{sant}
The loop functor $\Omega: \SSet^{\geq 1}_{\ast} \rightarrow \SSet$
is conservative and preserves sifted colimits.
\end{corollary}

\begin{proof}
Using Theorem \ref{preslage}, we may reduce to the problem of showing that the forgetful
functor $\theta: \Mon_{ \OpE{1}}^{\glike}(\SSet) \rightarrow \SSet$ is conservative and preserves sifted colimits. Since $\Mon_{\OpE{1}}^{\glike}(\SSet)$ is stable under colimits in
$\Mon_{ \OpE{1}}(\SSet)$, it suffices to show that the forgetful functor
$\Mon_{ \OpE{1}}(\SSet) \rightarrow \SSet$ is conservative and preserves sifted colimits. This follows from Proposition \symmetricref{ungbatt}, Proposition \symmetricref{fillfemme}, and Corollary \symmetricref{jumunj22}.
\end{proof}

\begin{corollary}\label{stunnek}
For every integer $n \geq 0$, the loop functor $\Omega: \SSet^{\geq n+1}_{\ast} \rightarrow \SSet^{\geq n}_{\ast}$ is conservative
and preserves sifted colimits.
\end{corollary}

\begin{corollary}\label{koso}
Let $\Spectra_{\geq 0}$ denote the $\infty$-category of connective spectra.
Then the functor $\Omega^{\infty}_{\ast}: \Spectra_{\geq 0} \rightarrow \SSet_{\ast}$
is conservative and preserves sifted colimits.
\end{corollary}

\begin{proof}
Write $\Spectra_{\geq 0}$ as the limit of the tower 
$\cdots \rightarrow \SSet^{\geq 1}_{\ast} \stackrel{ \Omega}{\rightarrow} \SSet^{\geq 0}_{\ast}$
and apply Corollary \ref{stunnek}.
\end{proof}

\begin{remark}\label{monoidchar}
Let $\calX$ be an $\infty$-topos, and regard $\calX$ as endowed with the Cartesian symmetric monoidal structure. Theorem \ref{preslage} guarantees the existence of
an equivalence $\theta: \calX^{\geq 1}_{\ast} \simeq \Mon^{\glike}( \calX) \simeq \Alg^{\glike}(\calX)$, where $\Alg^{\glike}( \calX)$ denotes the essential image of $\Mon^{\glike}(\calX)$ under the equivalence of $\infty$-categories $\Mon(\calX) \simeq \Alg(\calX)$ of Proposition \monoidref{ungbat}. 
This equivalence fits into a commutative diagram
$$ \xymatrix{ \Fun( \Delta^1, \calX) \times_{ \Fun( \{1\}, \calX) } \calX^{\geq 1}_{\ast} \ar[r]^-{ \overline{\theta} } \ar[d] & \Mod^{\glike}(\calX) \ar[d] \\
\calX^{\geq 1}_{\ast} \ar[r]^{\theta} & \Alg^{\glike}(\calX), }$$
where $\Mod^{\glike}(\calX)$ denotes the fiber product $\Mod(\calX) \times_{ \Alg(\calX) } \Alg^{\glike}(\calX)$ and $\overline{\theta}$ is an equivalence of $\infty$-categories. In other words, if
$X \in \calX$ is a pointed connected object, then there is a canonical equivalence between the $\infty$-topos $\calX_{/X}$ and the $\infty$-category $\Mod_{ \theta(X) }(\calX)$ of $\theta(X)$-module objects of $\calX$. 

To prove this, we let $\calD$ denote the full subcategory of $\Fun( \Delta^1 \times \Nerve( \cDelta_{+})^{op}, \calX)$ spanned by those functors $F$ with the following properties:
\begin{itemize}
\item[$(i)$] The functor $F$ is a right Kan extension of its restriction to the full subcategory 
$\calK \subseteq \Delta^1 \times \Nerve( \cDelta_{+})^{op}$ spanned by the objects $(0, [-1])$, $(1, [-1])$, and $(1, [0])$.
\item[$(ii)$] The object $F(1, [0]) \in \calX$ is final.
\item[$(iii)$] The augmentation map $F( 1, [0] ) \rightarrow F(1, [-1])$ is an effective epimorphism
(equivalently, the object $F( 1, [-1]) \in \calX$ is $1$-connective).
\end{itemize}
It follows from Proposition \toposref{lklk} that the restriction map $F \mapsto F | \calK$ determines
a trivial Kan fibration $\calD \rightarrow \Fun( \Delta^1, \calX) \times_{ \Fun( \{1\}, \calX) } \calX^{\geq 1}_{\ast}$. To construct the functor $\overline{\theta}$, it will suffice to show that the restriction functor
$F \mapsto F | ( \Delta^1 \times \Nerve( \cDelta)^{op} )$ is a trivial Kan fibration from
$\calD$ onto $\Mon^{L}(\calX) \times_{ \Mon(\calX) } \Mon^{\glike}(\calX)$, where
$\Mon^{L}(\calX)$ is described in Definition \monoidref{monll}. Using Proposition \toposref{acekan}, we see that $(i)$ is equivalent to the following pair of assertions:
\begin{itemize}
\item[$(i_0)$] The restriction $F | (\{1\} \times \Nerve( \cDelta_{+})^{op})$ is a right
Kan extension of its restriction to $\{1\} \times \Nerve( \cDelta^{\leq 0}_{+})^{op}$.
\item[$(i_1)$] The functor $F$ determines a Cartesian natural transformation from
$F_0 = F | ( \{0\} \times \Nerve( \cDelta_{+})^{op})$ to $F_1 = F| ( \{1\} \times \Nerve( \cDelta_{+})^{op} )$. 
\end{itemize}
Assertions $(i_0)$, $(ii)$, and $(iii)$ are equivalent to requirement that the functor $F_1$ belongs to
the full subcategory $\calC \subseteq \Fun( \Nerve( \cDelta_+)^{op}, \calX)$ appearing in the proof of Theorem \ref{preslage}. In particular, these conditions guarantee that $F_1$ is a colimit diagram.
Combining this observation with Theorem \toposref{charleschar} allows us to replace $(i_1)$ by the following pair of conditions:
\begin{itemize}
\item[$(i'_1)$] The functor $F_0$ is a colimit diagram.
\item[$(i''_1)$] The restriction $F | ( \Delta^1 \times \Nerve( \cDelta)^{op} )$ is a Cartesian transformation from $F_0 | \Nerve( \cDelta)^{op}$ to $F_1 | \Nerve( \cDelta)^{op}$. 
\end{itemize}
It follows that $\calY$ can be identified with the full subcategory of $\Fun( \Delta^1, \Nerve( \cDelta_{+})^{op})$ spanned by those functors $F$ such that 
$F' = F| ( \Delta^1 \times \Nerve(\cDelta)^{op} )$ belongs to $\Mon^{L}(\calX) \times_{ \Mon(\calX) } \Mon^{\glike}(\calX)$ and $F$ is a left Kan extension of $F'$. The desired result now follows from
Proposition \toposref{lklk}.
\end{remark}

\begin{remark}\label{monoidchur}
In the situation of Remark \ref{monoidchar}, let $X$ be a pointed $1$-connective object of the
$\infty$-topos $\calX$. Under the equivalence $\calX_{/X} \simeq \Alg_{ \theta(X)}(\calX)$,
the forgetful functor $\Mod_{ \theta X}(\calX) \rightarrow \calX$ corresponds to the functor
$(Y \rightarrow X) \mapsto (Y \times_{X} {\bf 1} )$ given by passing to the fiber over the base point $\eta: {\bf 1} \rightarrow X$ (here ${\bf 1}$ denotes the final object of $\calX$). 
It follows that the free module functor $\calX \rightarrow \Mod_{\theta(X)}(\calX)$ corresponds to
the functor $\calX \simeq \calX_{ / {\bf 1}} \rightarrow \calX_{/X}$ given by composition with
$\eta$. 
\end{remark}

We note that the loop functor $\Omega: \SSet^{\geq 1}_{\ast} \rightarrow \SSet$ is corepresentable
by the $1$-sphere $S^1 \in \SSet^{\geq 1}_{\ast}$. It follows from Corollary \ref{sant} that $S^1$ is
a compact projective object of $\SSet^{\geq 1}_{\ast}$. Since the collection of compact projective objects of $\SSet^{\geq 1}_{\ast}$ is stable under finite coproducts, we deduce the following:

\begin{corollary}\label{wanton}
Let $F$ be a finitely generated free group, and $BF$ its classifying space. Then
$BF$ is a compact projective object of $\SSet^{\geq 1}_{\ast}$.
\end{corollary}

For each $n \geq 0$, let $F(n)$ denote the free group on $n$ generators, and $BF(n)$ a classifying space for $F(n)$. Let $\calF$ denote the full subcategory of the category of groups spanned by the objects $\{ F(n) \}_{n \geq 0}$. We observe that the construction $F(n) \mapsto BF(n)$ determines a fully faithful embedding $i: \Nerve(\calF) \rightarrow \SSet^{\geq 1}_{\ast}$.
Let $\calP_{\Sigma}( \Nerve(\calF))$ be defined as in \S \toposref{stable11}
(that is, $\calP_{\Sigma}(\Nerve(\calF))$ is the $\infty$-category freely generated by $\Nerve(\calF)$ under sifted colimits).

\begin{remark}\label{studious}
According to Corollary \toposref{smokerr}, the $\infty$-category $\calP_{\Sigma}( \Nerve(\calF))$
is equivalent to the underlying $\infty$-category of the simplicial model category $\bfA$ of simplicial groups.
\end{remark}

It follows from Proposition \toposref{surottt} that the fully faithful embedding $i$ is
equivalent to a composition
$$ \Nerve(\calF) \stackrel{j}{\rightarrow} \calP_{\Sigma}( \Nerve(\calF)) \stackrel{F}{\rightarrow} \SSet^{\geq 1}_{\ast},$$
where $F$ is a functor which preserves sifted colimits (moreover, the functor $F$ is essentially unique).

\begin{corollary}\label{cansa}
The functor $F: \calP_{\Sigma}( \Nerve(\calF)) \rightarrow \SSet^{\geq 1}_{\ast}$ is an equivalence of $\infty$-categories.
\end{corollary}

\begin{remark}\label{jurm}
Combining Corollary \ref{cansa} and Remark \ref{studious}, we recover the following classical fact:
the homotopy theory of pointed connected spaces is equivalent to the homotopy theory of simplicial groups.
\end{remark}

\begin{proof}[Proof of Corollary \ref{cansa}]
Since $i: \Nerve(\calF) \rightarrow \SSet^{\geq 1}_{\ast}$ is fully faithful and its essential image
consists of compact projective objects (Corollary \ref{wanton}), Proposition \toposref{smearof} implies that $F$ is fully faithful. We observe that the functor $i$ preserves finite coproducts, so that $F$ preserves small colimits by virtue of Proposition \toposref{surottt}. Using Corollary \toposref{adjointfunctor}, we deduce
that $F$ admits a right adjoint $G$. Since $F$ is fully faithful, $G$ is a colocalization functor; to complete the proof, it will suffice to show that $G$ is conservative. 

Let $f: X \rightarrow Y$ be a morphism in $\SSet^{\geq 1}_{\ast}$ such that $G(f)$ is an equivalence; we wish to prove that $f$ is an equivalence. Let $\Z$ be the free group on one generator, and
$j\Z$ its image in $\calP_{\Sigma}( \Nerve(\calF))$. Then $f$ induces a homotopy equivalence
$$ \bHom_{ \SSet^{\geq 1}_{\ast}}( S^1, X) \simeq \bHom_{ \calP_{\Sigma}( \Nerve(\calF))}( j\Z, GX)
\rightarrow \bHom_{ \calP_{\Sigma}( \Nerve(\calF))}( j \Z, GY) \simeq \bHom_{ \SSet^{\geq 1}_{\ast}}( S^1, Y).$$
It follows that $\Omega(f): \Omega X \rightarrow \Omega Y$ is a homotopy
equivalence, so that $f$ is a homotopy equivalence by virtue of Corollary \ref{sant}.
\end{proof}

We are now ready to prove a more precise version of Theorem \ref{preslage}:

\begin{theorem}\label{slage}
Let $k > 0$, and let $\beta: \SSet_{\ast} \rightarrow \Mon_{ \OpE{k}}( \SSet)$ be
the functor described at the beginning of this section. Then $\beta| \SSet_{\ast}^{\geq k}$ is equivalent
to the functor $\alpha$ constructed in the proof of Theorem \ref{preslage}, so that
$\beta$ induces an equivalence of $\infty$-categories
$\SSet_{\ast}^{\geq k} \rightarrow \Mon_{ \OpE{k}}^{\glike}( \SSet).$
\end{theorem}

\begin{proof}
As before, we work by induction on $k$. Suppose first that $k = 1$, and let
$\alpha^{-1}$ be a homotopy inverse to $\alpha$. We wish to show that the
composition $\theta: \alpha^{-1} \circ \beta$ is equivalent to the identity functor
from $\SSet_{\ast}^{\geq 1}$ to itself. Let $F: \calP_{\Sigma}( \Nerve(\calF)) \rightarrow
\SSet_{\ast}^{\geq 1}$ be the equivalence of $\infty$-categories of Corollary \ref{cansa};
it will suffice to construct an equivalence $F \simeq \theta \circ F$.

Let $\Omega: \SSet^{\geq 1}_{\ast} \rightarrow \SSet$ denote the loop space functor. It is easy to see that there is an equivalence $\Omega \circ \theta \simeq \Omega$. Using Corollary \ref{sant}, we deduce that
$\theta$ commutes with sifted colimits. In view of Proposition \toposref{surottt}, it will suffice to show
that $F \circ j$ is equivalent to $\theta \circ F \circ j$ in the $\infty$-category
$\Fun( \Nerve(\calF), \SSet_{\ast}^{\geq 1})$; here $j: \Nerve(\calF) \rightarrow \calP_{\Sigma}( \Nerve(\calF))$ denotes the Yoneda embedding so that $F \circ j$ is equivalent to the classifying space functor $i: \Nerve(\calF) \rightarrow \SSet_{\ast}^{\geq 1}$.

The functor $\pi_0: \SSet \rightarrow \Nerve(\Set)$ induces a functor from
$U: \Mon_{ \OpE{1}}^{\glike}(\SSet)$ to the nerve of the category $\calG$ of groups. This functor
admits a right adjoint $T$, given by the fully faithful embedding
$$\Nerve(\calG) \simeq \Mon_{ \OpE{1}}^{\glike}( \Nerve(\Set)) \subseteq \Mon_{ \OpE{1}}^{\glike}(\SSet).$$
The composition $U \circ \beta$ can be identified with the functor which carries a pointed space
$X$ to its fundamental group $\pi_1 X$, while $\alpha^{-1} \circ T$ carries a group $G$
to a classifying space $BG \in \calA \subseteq \SSet^{\geq 1}_{\ast}$. Consequently, on
$\SSet^{\geq 1}_{\ast}$, the composition $\alpha^{-1} \circ T \circ U \circ \beta$
agrees with the truncation functor $\tau_{\leq 1}$, so there is a natural transformation of functors
$v: \id_{ \SSet^{\geq 1}_{\ast}} \rightarrow \alpha^{-1} \circ T \circ U \circ \beta$.
Since $T$ and $U$ are adjoint, we also have a unit transformation
$u: \theta \rightarrow \alpha^{-1} \circ T \circ U \circ \beta$. The natural transformation
$v$ is an equivalence when restricted to $1$-truncated spaces, and the natural
transformation $u$ is an equivalence when restricted to spaces $X$ such that $\Omega X$ is
discrete. In particular, $u$ and $v$ are both equivalences on the essential image of the fully faithful embedding $i: \Nerve(\calF) \rightarrow \SSet_{\ast}^{\geq 1}$. It follows that $u$ and $v$
determines an equivalence of functors
$$ F \circ j \simeq i \simeq \alpha^{-1} \circ T \circ U \circ \beta \circ i
\simeq \alpha^{-1} \circ \beta \circ i \simeq \theta \circ F \circ j.$$
This completes the proof in the case $k=1$.

Suppose now that $k > 1$. We observe that the functor $\beta$ factors as a composition
$$ \SSet_{\ast}^{\geq k} \stackrel{\beta'}{\rightarrow} \Mon_{ \OpE{k}}^{\glike}( \SSet_{\ast})
\stackrel{\beta''}{\rightarrow} \Mon_{ \OpE{k}}^{\glike}( \SSet),$$
where the functor $\beta''$ is an equivalence of $\infty$-categories (as in the proof of Theorem \ref{preslage}). Consequently, it will suffice to show that $\beta'$ is equivalent to the functor
$\alpha'$ constructed in the proof of Proposition \ref{preslage}. Write $k = k_{-} + k_{+}$, where
$0 < k_{-}, k_{+} < k$. By the inductive hypothesis, we may assume that the functors
$$ \beta'_{-}: \SSet^{\geq 1}_{\ast} \rightarrow \Mon_{ \OpE{k_{-}}}( \SSet_{\ast})
\quad \quad \beta'_{+} : \SSet^{\geq 1}_{\ast} \rightarrow \Mon_{ \OpE{k_{+}}}( \SSet_{\ast})$$
are equivalent to the functors $\alpha'_{-}$ and $\alpha'_{+}$ constructed in the proof of
Proposition \ref{preslage}. The equivalence of $\alpha'$ and $\beta'$ follows from the homotopy commutativity of the diagram
$$ \xymatrix{  \SSet_{\ast}^{\geq k} \ar[r]^-{\beta'_{+}} \ar[d]^{\beta'} & \Mon_{ \OpE{k_{+}}}( \SSet_{\ast}^{\geq k_{-}}) \ar[d]^{\beta'_{-}} \\
\Mon_{ \OpE{k}}^{\glike}( \SSet_{\ast}) \ar[r]^-{\delta} & \Mon_{ \OpE{k_{+}}}( \Mon_{\OpE{k_{-}}}^{\glike}(\SSet_{\ast})). }$$
\end{proof}

\subsection{Coherence of the Little Cubes Operads}\label{slabba}

In this section, we will use the coherence criterion of \S \ref{cohcrit} (more specifically, Theorem \ref{uggus}) to prove the following result, which guarantees the existence of a good theory of modules over $\OpE{k}$-algebras: 

\begin{theorem}\label{cubecoh}
Let $k \geq 0$ be a nonnegative integer. Then the little cubes $\infty$-operad $\OpE{k}$ is
coherent.
\end{theorem}

In order to prove Theorem \ref{cubecoh}, we will need to introduce a few simple constructions for passing converting information about simplicial or topological operads (such as $\TopE{k}$) into information about
their underlying $\infty$-operads (such as $\OpE{k}$). 

\begin{notation}
Let $\calO$ be a simplicial operad (that is, a simplicial colored operad having a single distinguished object), and let $\calO^{\otimes}$ be the simplicial category described in Notation \symmetricref{capin}:
the objects of $\calO^{\otimes}$ are objects $\seg{n} \in \FinSeg$, and the morphisms spaces
$\calO^{\otimes}$ are given by the formula
$$ \bHom_{ \calO^{\otimes}}( \seg{m}, \seg{n}) = \coprod_{ \alpha: \seg{m} \rightarrow \seg{n}}
\prod_{ 1 \leq i \leq n} \Mul_{\calO}( \alpha^{-1} \{i\}, \{i\})$$
where $\alpha$ ranges over all maps $\seg{m} \rightarrow \seg{n}$ in $\FinSeg$.
We will say that a morphism in $\calO^{\otimes}$ is {\it active} if its image in $\FinSeg$ is
active, and we let $\bHom_{\calO^{\otimes}}^{\acti}( \seg{m}, \seg{n})$ denote the summand of
$\bHom_{\calO^{\otimes}}( \seg{m}, \seg{n})$ spanned by the active morphisms.

We will say that $\calO$ is {\it unital} if $\Mul_{\calO}( \emptyset, \{ 0 \})$ is isomorphic to
$\Delta^0$; in this case, every semi-inert morphism $\alpha: \seg{m} \rightarrow \seg{n}$ in
$\FinSeg$ can be lifted {\em uniquely} to a morphism $\overline{\alpha}$ in $\calO^{\otimes}$.
In particular, the canonical inclusion $i: \seg{m} \rightarrow \seg{m+1}$ admits a unique lift
$\overline{i}: \seg{m} \rightarrow \seg{m+1}$ in $\calO^{\otimes}$. Composition with
$\overline{i}$ induces a map of simplicial sets
$$ \theta: \bHom^{\acti}_{ \calO^{\otimes}}( \seg{m+1}, \seg{n}) \rightarrow \bHom^{\acti}_{\calO^{\otimes}}( \seg{m}, \seg{n}).$$
For every active morphism $f: \seg{m} \rightarrow \seg{n}$ in $\calO^{\otimes}$, we will denote the simplicial set
$\theta^{-1} \{f\}$ by $\Ext_{\Delta}(f)$; we will refer to $\Ext_{\Delta}(f)$ as the 
{\it space of strict extensions of $f$}.
\end{notation}

\begin{construction}\label{sadper}
Let $\calO$ be a fibrant simplicial operad, and let $\Nerve(\calO)^{\otimes}$ be the 
underlying $\infty$-operad (Definition \ref{eli}). Suppose we are given a sequence
of active morphisms
$$ \seg{m_0} \stackrel{ f_1}{\rightarrow} \seg{m_1} \stackrel{f_2}{\rightarrow}
\ldots \stackrel{f_n}{\rightarrow} \seg{m_n}$$
in $\calO^{\otimes}$. This sequence determines an $n$-simplex $\sigma$ of
$\Nerve(\calO)^{\otimes}$. Let $S \subseteq [n]$ be a proper nonempty subset having maximal
element $j-1$. We define a map of simplicial sets
$\theta: \Ext_{\Delta}(f_j) \rightarrow \Ext( \sigma, S)$
as follows: for every $k$-simplex $\tau: \Delta^{k} \rightarrow \Ext_{\Delta}(f_{j})$, 
$\theta(\tau)$ is a $k$-simplex of $\Ext(\sigma,S)$ corresponding to a map of simplicial categories
$\psi: \sCoNerve[ \Delta^n \times \Delta^{k+1}] \rightarrow \calO^{\otimes}$, which may
be described as follows:
\begin{itemize}
\item[$(i)$] On objects, the functor $\psi$ is given by the formula
$$ \psi( n', k' ) = \begin{cases} \seg{ m_{n'} } & \text{ if } k' = 0 \text{ or } n' \notin S \\
\seg{ m_{n'} + 1} & \text{ otherwise. } \end{cases}$$
\item[$(ii)$] Fix a pair of vertices $(n', k'), (n'',k'') \in \Delta^n \times \Delta^{k+1}$.
Then $\psi$ induces a map of simplicial sets
$\phi: \bHom_{ \sCoNerve[ \Delta^n \times \Delta^{k+1}]}( (n', k'), (n'', k''))
\rightarrow \bHom_{ \calO^{\otimes}}( \psi(n',k'), \psi(n'',k''))$. The left hand side
can be identified with the nerve of the partially ordered set 
$P$ of chains 
$$(n',k') = (n_0, k_0) \leq (n_1, k_1) \leq \ldots \leq (n_p, k_p) = (n'', k'')$$
in $[n] \times [k+1]$. If $\psi( n',k') = \seg{m_{n'}}$
or $\psi( n'', k'') = \seg{m_{n''} + 1}$, then $\phi$ is given by the constant map determined by
$f_{n''} \circ \cdots \circ f_{n'+1}$. Otherwise, $\phi$ is given by composing the morphisms
$f_{j-1} \circ \cdots \circ f_{n'+1}$ and  $f_{n''} \circ \cdots \circ f_{j+1}$ with the map
$$ \Nerve(P) \stackrel{\phi_0}{\rightarrow} \Delta^{k} \stackrel{\tau}{\rightarrow}
\Ext_{\Delta}( f_j) \rightarrow \bHom_{\calO^{\otimes}}( \seg{m_{j-1}+1}, \seg{m_j}),$$
where $\phi_0$ is induced by the map of partially ordered sets $P \rightarrow [k]$ which
carries a chain $(n',k') = (n_0, k_0) \leq (n_1, k_1) \leq \ldots \leq (n_p, k_p) = (n'', k'')$
to the supremum of the set $\{ k_i - 1: n_i \in S\} \subseteq [k]$.
\end{itemize}
\end{construction}

\begin{remark}\label{saucewell}
In the situation of Construction \ref{sadper}, the simplicial set
$\Ext(\sigma,S)$ can be identified with the homotopy fiber of the map
$$ \beta: \bHom_{ \calO^{\otimes}}^{\acti}( \seg{m_{j-1} + 1}, \seg{m_j})
\rightarrow \bHom_{ \calO^{\otimes}}^{\acti}( \seg{m_{j-1}}, \seg{m_j}),$$
while $\Ext_{\Delta}( f_{j})$ can be identified with the actual fiber of
$\beta$. The map $\theta$ of Construction \ref{sadper} can be identified with the canonical
map from the actual fiber to the homotopy fiber. 
\end{remark}

\begin{proposition}\label{sendy}
Let $\calO$ be a fibrant simplicial operad, and assume that every morphism 
in the simplicial category $\calO = \calO^{\otimes}_{\seg{1}}$ admits a homotopy inverse.
Suppose that, for every pair active morphisms
$f_0: \seg{m} \rightarrow \seg{n}$ and $g_0: \seg{n} \rightarrow \seg{1}$ in $\calO^{\otimes}$,
there exist morphisms $f: \seg{m} \rightarrow \seg{n}$, $h: \seg{n} \rightarrow \seg{n}$, and
$g: \seg{n} \rightarrow \seg{1}$ satisfying the following conditions:
\begin{itemize}
\item[$(i)$] The map $f$ is homotopic to $f_0$, the map $g$ is homotopic to $g_0$, and the
map $h$ is homotopic to $\id_{\seg{n}}$.
\item[$(ii)$] Each of the sequences
$$ \Ext_{\Delta}( h) \rightarrow \bHom_{\calO^{\otimes}}^{\acti}( \seg{n+1}, \seg{n})
\rightarrow \bHom_{\calO^{\otimes}}^{\acti}( \seg{n}, \seg{n})$$
$$ \Ext_{\Delta}(g \circ h) \rightarrow \bHom_{\calO^{\otimes}}^{\acti}( \seg{n+1}, \seg{1})
\rightarrow \bHom_{\calO^{\otimes}}^{\acti}( \seg{n}, \seg{1})$$
$$ \Ext_{\Delta}(h \circ f) \rightarrow \bHom_{\calO^{\otimes}}^{\acti}( \seg{m+1}, \seg{n})
\rightarrow \bHom_{\calO^{\otimes}}^{\acti}( \seg{m}, \seg{n})$$
$$ \Ext_{\Delta}(g \circ h \circ f) \rightarrow \bHom_{\calO^{\otimes}}^{\acti}( \seg{m+1}, \seg{1})
\rightarrow \bHom_{\calO^{\otimes}}^{\acti}( \seg{m}, \seg{1})$$
is a homotopy fiber sequence.
\item[$(iii)$] The diagram
$$ \xymatrix{ \Ext_{\Delta}( h) \ar[r] \ar[d] & \Ext_{\Delta}(g \circ h) \ar[d] \\
\Ext_{\Delta}( h \circ f) \ar[r] & \Ext_{\Delta}( g \circ h \circ f) }$$
is a homotopy pushout square of simplicial sets.
\end{itemize}
Then the $\infty$-operad $\Nerve(\calO)^{\otimes}$ is coherent.
\end{proposition}

\begin{proof}
We will show that $\Nerve(\calO)^{\otimes}$ satisfies criterion $(3)$ of Theorem \ref{uggus}.
Suppose we are given a degenerate 
$3$-simplex $\sigma:$
$$ \xymatrix{ & \seg{n} \ar[dr]^{\id} \ar[rr]^{g_0} & & \seg{1} \\
\seg{m} \ar[ur]^{f_0} \ar[rr]^{f_0} & & \seg{n} \ar[ur]^{g_0}}$$
in $\Nerve(\calO)^{\otimes}$, where $f$ and $g$ are active. 
We wish to show that the diagram
$$ \xymatrix{ \Ext( \sigma, \{0,1\} ) \ar[r] \ar[d] & \Ext( \sigma | \Delta^{ \{0,1,3 \} }, \{0,1\} ) \ar[d] \\
\Ext( \sigma | \Delta^{ \{0,2,3\} }, \{0\} ) \ar[r] & \Ext( \sigma| \Delta^{ \{0,3\} }, \{0\} ) }$$
is a homotopy pushout square of Kan complexes. In proving this, we are free to replace
$\sigma$ by any equivalent diagram $\sigma': \Delta^3 \rightarrow \Nerve(\calO)^{\otimes}$.
We may therefore assume that $\sigma'$ is determined by a triple of morphisms
$f: \seg{m} \rightarrow \seg{n}$, $h: \seg{n} \rightarrow \seg{n}$, and $g: \seg{n} \rightarrow \seg{1}$
satisfying conditions $(ii)$ and $(iii)$ above. Using Remark \ref{saucewell}, we see
that Construction \ref{sadper} determines a weak homotopy equivalence between the diagrams
$$ \xymatrix{ \Ext_{\Delta}( h) \ar[r] \ar[d] & \Ext_{\Delta}(g \circ h) \ar[d] &  \Ext( \sigma, \{0,1\} ) \ar[r] \ar[d] & \Ext( \sigma | \Delta^{ \{0,1,3 \} }, \{0,1\} ) \ar[d]  \\
\Ext_{\Delta}( h \circ f) \ar[r] & \Ext_{\Delta}( g \circ h \circ f) & \Ext( \sigma | \Delta^{ \{0,2,3\} }, \{0\} ) \ar[r] & \Ext( \sigma| \Delta^{ \{0,3\} }, \{0\} ). }$$
Since the diagram on the left is a homotopy pushout square by virtue of $(iii)$, the diagram
on the right is also a homotopy pushout square.
\end{proof}

\begin{proof}[Proof of Theorem \ref{cubecoh}]
Let $\calO = \Sing \TopE{k}$ denote the simplicial operad associated to the topological operad
$\TopE{k}$. We will say that a rectilinear embedding
$f \in \Rect( \Cube{k} \times \nostar{n}, \Cube{k})$ is {\it generic} if $f$ can be extended to
an $\overline{f}: \overline{ \Cube{k}} \times \nostar{n} \rightarrow \Cube{k}$, where
$\overline{ \Cube{k} } = [-1,1]^{k}$ is a closed cube of dimension $k$. We will say that an active morphism $f: \seg{n} \rightarrow \seg{m}$ in $\calO^{\otimes}$ is {\it generic} if it corresponds to a sequence of $m$ rectlinear embeddings which are generic.

We observe the following:
\begin{itemize}
\item[$(a)$] If $f$ is generic, then the difference $\Cube{k} - \overline{f}( \overline{ \Cube{k}} \times \nostar{n})$
is homotopy equivalent to $\Cube{k} - f( \{0\} \times \nostar{n})$. It follows that the sequence
$\Ext_{\Delta}(f) \rightarrow \bHom_{ \calO^{\otimes}}( \seg{n+1}, \seg{1})
\rightarrow \bHom_{ \calO^{\otimes}}( \seg{n}, \seg{1})$ is homotopy
equivalent to the fiber sequence of configuration spaces (see Remark \ref{sove})
$$ \Cube{k} - f( \{0\} \times \nostar{n}) \rightarrow \Conf( \nostar{n+1}, \Cube{k})
\rightarrow \Conf( \nostar{n}, \Cube{k}),$$
hence also a homotopy fiber sequence. More generally, if $f: \seg{n} \rightarrow \seg{m}$ is generic, then 
$$ \Ext_{\Delta}(f) \rightarrow \bHom_{ \calO^{\otimes}}^{\acti}( \seg{n+1}, \seg{m})
\rightarrow \bHom_{\calO^{\otimes}}^{\acti}( \seg{n}, \seg{m}),$$
is a fiber sequence.
\item[$(b)$] Every rectlinear embedding $f_0 \in \Rect( \Cube{k} \times \nostar{n}, \Cube{k})$
is homotopic to a generic rectilinear embedding $f$ (for example, we can take $f$
to be the composition of $f_0$ with the ``contracting'' map
$\Cube{k} \times \nostar{n} \simeq (\frac{-1}{2}, \frac{1}{2})^{k} \times \nostar{n}
\hookrightarrow \Cube{k} \times \nostar{n}$). Similarly, every active morphism in $\calO^{\otimes}$ is homotopic to a generic morphism.
\item[$(c)$] The collection of generic morphisms in $\calO^{\otimes}$ is stable under composition.
\end{itemize}

To prove that $\OpE{k}$ is coherent, it will suffice to show that the simplicial operad $\calO$ satisfies the criteria of Proposition \ref{sendy}. It is clear that every map in $\calO$ admits a homotopy inverse
(in fact, every rectilinear embedding from $\Cube{k}$ to itself is homotopic to the identity). In
view of $(a)$, $(b)$, and $(c)$ above, it will suffice to show that the diagram
$$ \xymatrix{ \Ext_{\Delta}( h) \ar[r] \ar[d] & \Ext_{\Delta}(g \circ h) \ar[d] \\
\Ext_{\Delta}( h \circ f) \ar[r] & \Ext_{\Delta}( g \circ h \circ f) }$$
is a homotopy pushout square for every triple of active morphisms
$$ \seg{m} \stackrel{f}{\rightarrow} \seg{n} \stackrel{h}{\rightarrow} \seg{n} \stackrel{g}{\rightarrow} \seg{1}$$
in $\calO^{\otimes}$, provided that each of the underlying rectilinear embeddings is generic.

Let $U_0 \subseteq U_1 \subseteq U_2$ be the images of
$g \circ h \circ f$, $g \circ h$, and $g$, respectively. Let $\overline{U}_i$ denote the closure
of $U_{i}$. We now set
$$ V = \Cube{k} - \overline{U}_1 \quad \quad W = U_2 - \overline{U}_0.$$
Note that $V \cup W = \Cube{k} - \overline{U}_0$ and $V \cap W = U_2 - \overline{U}_1$.
The argument of Remark \ref{sove} shows that evaluation at the origin of $\Cube{k}$ determines weak homotopy equivalences
$$ \Ext_{\Delta}(h) \rightarrow \Sing( V \cap W) \quad \quad \Ext_{\Delta}( g \circ h) \rightarrow
\Sing(V)$$
$$ \Ext_{\Delta}( h \circ f) \rightarrow \Sing(W) \quad \quad \Ext_{\Delta}( g \circ h \circ f) \rightarrow \Sing(W \cup V).$$
It will therefore suffice to show that the diagram
$$ \xymatrix{ \Sing( V \cap W) \ar[r] \ar[d] & \Sing(V) \ar[d] \\
\Sing(W) \ar[r] & \Sing(W \cup V) }$$
is a homotopy pushout square of Kan complexes, which follows from Theorem \ref{vankamp}.
\end{proof}

\subsection{Tensor Products of $\OpE{k}$-Algebras}\label{tensor2}

Let $\calO^{\otimes}$ be any $\infty$-operad, and let $\calC^{\otimes}$ be a symmetric monoidal $\infty$-category. As explained in \S \symmetricref{comm1.8}, the $\infty$-category $\Alg_{\calO}(\calC)$ of $\calO$-algebras in $\calC$ inherits the structure of a symmetric monoidal $\infty$-category.
In particular, for every pair of objects $A, B \in \Alg_{\calO}(\calC)$, we have another object
$A \otimes B \in \Alg_{\calO}(\calC)$, which is given on objects by the formula
$$ (A \otimes B)(X) = A(X) \otimes B(X)$$
for $X \in \calO$.

In the special case where $\calO^{\otimes} = \Nerve(\FinSeg)$ is the commutative $\infty$-operad,
the tensor product $A \otimes B$ can be identified with the coproduct of $A$ and $B$ in the
$\infty$-category $\Alg_{\calO}(\calC) = \CAlg(\calC)$ (Proposition \symmetricref{cocarten}).
For other $\infty$-operads, this is generally not the case. Suppose, for example, that
$\calO^{\otimes}$ is the associative $\infty$-operad, and that $\calC$ is the (nerve of the) ordinary category $\Vect_{\C}$ of vector spaces over the field $\C$ of complex numbers. Then
$\Alg_{\Ass}(\calC)$ is equivalent the nerve of the category of
associative $\C$-algebras. Given a pair of associative $\C$-algebras $A$ and $B$, there is a diagram of associative algebras
$$ A \rightarrow A \otimes_{\C} B \leftarrow B,$$
but this diagram does not exhibit $A \otimes_{\C} B$ as a coproduct of $A$ and $B$. Instead,
it exhibits $A \otimes_{\C} B$ as the quotient of the coproduct $A \coprod B$ by the (two-sided) ideal generated by commutators $[a,b] = ab-ba$, where $a \in A$ and $b \in B$. In other words,
$A \otimes_{\C} B$ is freely generated by $A$ and $B$ subject to the condition that
$A$ and $B$ commute in $A \otimes_{\C} B$. 

Our goal in this section is to obtain an $\infty$-categorical generalization of the above assertion. We will replace the ordinary category $\Vect_{\C}$ by an arbitrary symmetric monoidal $\infty$-category $\calC$, and the associative $\infty$-operad $\Ass$ by a little $k$-cubes operad
$\OpE{k}$, for any $k \geq 0$ (we can recover the case of associative algebras by taking $k=1$,
by virtue of Example \ref{sulta}). Assume that $\calC$ admits small colimits, and that the tensor product of $\calC$ preserves small colimits separately in each variable. Then the forgetful functor
$\Alg_{\OpE{k}}( \calC) \rightarrow \calC$ admits a left adjoint $\Free: \calC \rightarrow \Alg_{\OpE{k}}(\calC)$ (Corollary \symmetricref{spaltwell}). Given a pair of objects $C,D \in \calC$, the tensor product $\Free(C) \otimes \Free(D)$
is generally not equivalent to the coproduct $\Free(C) \coprod \Free(D) \simeq \Free(C \coprod D)$.
To measure the difference, we note that every binary operation $f \in \Mul_{\OpE{k}}( \{ \seg{1}, \seg{1} \}, \seg{1} )$ gives rise to a map
$$ \phi_{f}: C \otimes D \rightarrow \Free(C \coprod D) \otimes \Free(C \coprod D)
\stackrel{f}{\rightarrow} \Free(C \coprod D).$$
Note that the composite map $C \otimes D \rightarrow \Free(C \coprod D) \stackrel{\psi}{\rightarrow} \Free(C) \otimes \Free(D)$ does not depend on $f$. 
The space of choices for the binary operation $f$ is homotopy equivalent to the configuration space of
pairs of points in $\R^{k}$ (Remark \ref{sove}), which is in turn homotopy equivalent to a sphere $S^{k-1}$. Allowing $f$ to vary, we obtain a map
$$\phi: (C \otimes D) \otimes S^{k-1} \rightarrow \Free(C \coprod D)$$
in $\calC$, where we regard $\calC$ as tensored over the $\infty$-category $\SSet$ of spaces
as explained in \S \toposref{quasilimit7}. Equivalently, we can view $\phi$ as a map
$$ \Free(C \otimes D \otimes S^{k-1}) \rightarrow \Free(C \coprod D),$$
which fits into a diagram
$$ \xymatrix{ \Free(C \otimes D \otimes S^{k-1} ) \ar[r] \ar[d] & \Free(C \coprod D) \ar[d] \\
\Free(C \otimes D) \ar[r] & \Free(C) \otimes \Free(D). }$$
The commutativity of this diagram encodes the fact that $\psi \circ \phi_{f}$ is independent of $f$; equivalently, it reflects the idea that $C$ and $D$ ``commute'' inside the tensor product
$\Free(C) \otimes \Free(D)$. The main result of this section can be formulated as follows:

\begin{theorem}\label{juke}
Let $k \geq 0$, let $\calC$ be a symmetric monoidal $\infty$-category which admits countable colimits, 
and assume that the tensor product on $\calC$ preserves countable colimits separately in each variable. Let $\Free: \calC \rightarrow \Alg_{ \OpE{k}}( \calC)$ be a left adjoint to the forgetful functor.
Then, for every pair of objects $C,D \in \calC$, the construction sketched above gives rise to a pushout diagram
$$ \xymatrix{ \Free(C \otimes D \otimes S^{k-1}) \ar[r] \ar[d] & \Free(C) \coprod \Free(D) \ar[d] \\
\Free(C \otimes D) \ar[r] & \Free(C) \otimes \Free(D) }$$
in $\calC$.
\end{theorem}

\begin{example}
Suppose that $k = 0$. In this case, we can identify the $\infty$-category
$\Alg_{\OpE{k}}(\calC)$ with the $\infty$-category $\calC_{ {\bf 1}/}$ (Proposition \symmetricref{ezalg}; here ${\bf 1}$ denotes the unit object of $\calC$, and the
free algebra functor $\Free: \calC \rightarrow \Alg_{\OpE{k}}(\calC)$ is given by the formula
$C \mapsto {\bf 1} \coprod C$. In this case, Theorem \ref{juke} asserts that the diagram
$$ \xymatrix{ {\bf 1} \ar[r] \ar[d] & {\bf 1} \coprod C \coprod D \ar[d] \\
{\bf 1} \coprod (C \otimes D) \ar[r] & ( {\bf 1} \coprod C) \otimes ( {\bf 1} \coprod D). }$$
This follows immediately from the calculation
$$ ( {\bf 1} \coprod C) \otimes ( {\bf 1} \coprod D) \simeq {\bf 1} \coprod C \coprod D \coprod (C \otimes D).$$
\end{example}

\begin{example}
In the case $k=1$, we can replace the operad $\OpE{1}$ with the associative $\infty$-operad
$\Ass$ (Example \ref{sulta}). In this case, Theorem \ref{juke} is equivalent to the assertion that 
the diagram
$$ \xymatrix{ \Free(C \otimes D) \ar@<.4ex>[r]^-{f} \ar@<-.4ex>[r]_-{g} & \Free(C) \coprod \Free(D) \ar[r] & 
\Free(C) \otimes \Free(D)}$$
is a coequalizer, where $f$ and $g$ are induced by the maps $C \otimes D \rightarrow \Free(C) \coprod \Free(D)$ given by multiplication on $\Free(C) \coprod \Free(D)$ in the two possible orders.
\end{example}

\begin{example}
We can also take $k = \infty$ in Theorem \ref{juke}. In this case, the sphere $S^{k-1}$ is contractible, so the left vertical map $\Free( C \otimes D \otimes S^{k-1}) \rightarrow \Free(C \otimes D)$ is an equivalence. Consequently, Theorem \ref{juke} reduces to the assertion that the right vertical map
$\Free(C) \coprod \Free(D) \rightarrow \Free(C) \otimes \Free(D)$ is an equivalence. This follows
from Proposition \symmetricref{cocarten}, since the $\infty$-operad $\OpE{k}$ is equivalent to the commutative $\infty$-operad $\Nerve(\FinSeg)$ (Corollary \ref{infec}). 
\end{example}

Let us now outline our approach to the proof of Theorem \ref{juke}. The rough idea is to construct
a functor $\calC \times \calC \rightarrow \Alg_{ \OpE{k}}( \calC)$ whose values can be computed in two different ways: the first computation will show that this functor is given by the formula
$(C,D) \mapsto \Free(C) \otimes \Free(D)$, while the second computation will show that it is given by
$$ (C,D) \mapsto \Free(C \otimes D) \coprod_{ \Free( C \otimes D) \otimes S^{k-1}} ( \Free(C) \coprod \Free(D)).$$
The construction will use the formalism of operadic left Kan extensions developed in
\S \symmetricref{iljest}, and the comparison between the two calculations rests on a transitivity
result for operadic left Kan extensions which is proven in \S \ref{tensor1}. 

For the discussion which follows, we fix $\infty$-operads $\calO^{\otimes}$ and
${\calD}^{\otimes}$. We will freely employ the notation of \S \ref{coprodinf}; in particular,
we let $\calT_{\calO}$ be the correspondence of $\infty$-operads
defined in Notation \ref{skaz}. Let $\psi: \calT_{\calO} \rightarrow
{\calD}^{\otimes}$ be a map of $\infty$-operad families and let $q: \calC^{\otimes} \rightarrow {\calD}^{\otimes}$
be a coCartesian fibration of $\infty$-operads, so that the $\infty$-categories
$\Alg^{-}_{\calO}(\calC)$, $\Alg^{+}_{\calO}(\calC)$, and $\Alg^{\pm}_{\calO}(\calC)$ are defined
as in Construction \ref{urmouth}. Similarly, we can define $\infty$-categories
$\Fun^{-}_{\calD}( \calO, \calC)$, $\Fun^{+}_{\calD}(\calO, \calC)$, and
$\Fun^{\pm}_{\calD}(\calO, \calC)$. There are evident forgetful functors
$$ \Alg^{-}_{\calO}(\calC) \rightarrow \Fun^{-}_{\calD}(\calO, \calC)$$
$$ \Alg^{+}_{\calO}(\calC) \rightarrow \Fun^{+}_{\calD}(\calO, \calC)$$
$$ \Alg^{\pm}_{\calO}(\calC) \rightarrow \Fun^{\pm}_{\calD}(\calO, \calC).$$
Under some mild hypotheses, these forgetful functors admit left adjoints, which we will
denote by $F^{-}$, $F^{+}$, and $F^{\pm}$. The construction
$(X,Y) \mapsto F^{-}(X) \otimes F^{+}(Y)$ determines a functor from $\Fun^{-}_{\calD}(\calO, \calC) \times
\Fun^{+}_{\calD}(\calO, \calC)$ to $\Alg^{\pm}_{\calO}(\calC)$. Our first step is to give a convenient description of this functor, using the theory of $q$-left Kan extensions.

\begin{remark}
In the special case where ${\calD}^{\otimes} = \FinSeg$ (so that $\calC^{\otimes}$ is a symmetric monoidal $\infty$-category), the $\infty$-categories $\Alg^{-}_{\calO}(\calC)$, 
$\Alg^{+}_{\calO}(\calC)$, and $\Alg^{\pm}_{\calO}(\calC)$ coincide and the superscripts become superfluous.
\end{remark}

\begin{notation}\label{bbess}
We let $\calT_{\calO}^{0}$ denote the subcategory of $\calT_{\calO}$ spanned
by all those morphisms $f: X \rightarrow Y$ satisfying the following condition:
\begin{itemize}
\item[$(\ast)$] If $X$ and $Y$ belong to $\calO^{\otimes} \boxplus \calO^{\otimes} \subseteq \calT_{\calO}$, then the image of $f$ in $\Nerve(\FinSeg)$ is inert.
\end{itemize}
\end{notation}

\begin{remark}
It follows easily from Lemma \ref{kwop} that $\calT^{0}_{\calO} \rightarrow \Delta^1 \rightarrow \Nerve(\FinSeg)$ is a $\Delta^1$-family of $\infty$-operads, which we view as a correspondence
of $\infty$-operads from ${\calO'}^{\otimes} \boxplus {\calO'}^{\otimes}$ to $\calO$; here
${\calO'}^{\otimes}$ denotes the fiber product $\calO^{\otimes} \times_{ \Nerve(\FinSeg)} \Triv$. 
\end{remark}

\begin{proposition}\label{onjuke}
Let $\kappa$ be an uncountable regular cardinal, let $\calO^{\otimes}$ be an $\infty$-operad which is essentially $\kappa$-small, let $q: \calC^{\otimes} \rightarrow \calD^{\otimes}$ be a coCartesian fibration of $\infty$-operads and $\psi: \calT_{\calO} \rightarrow \calD^{\otimes}$ an $\infty$-operad family map.
Assume that for each $D \in \calD$, the fiber $\calC_{D}$ admits $\kappa$-small colimits,
and that the $\calD$-monoidal structure on $\calC^{\otimes}$ is compatible with $\kappa$-small colimits.
Then:
\begin{itemize}
\item[$(1)$] For each symbol $\sigma \in \{ -, +, \pm \}$, the forgetful functor
$\Alg^{\sigma}_{\calO}(\calC) \rightarrow \Fun_{\calD}^{\sigma}(\calO, \calC)$ admits a left
adjoint $F^{\sigma}$. If we let ${\calO'}^{\otimes}$ denote the fiber product
$\calO^{\otimes} \times_{ \Nerve(\FinSeg)} \Triv$, then $F^{\sigma}$ is given by composing a homotopy
inverse to the trivial Kan fibration $\Alg^{\sigma}_{\calO'}(\calC) \rightarrow \Fun^{\sigma}_{\calD}( \calO', \calC) \simeq
\Fun^{\sigma}_{\calD}(\calO, \calC)$ with the functor of operadic $q$-left Kan extension along the inclusion ${\calO'}^{\otimes} \rightarrow \calO^{\otimes}$.
\item[$(2)$] Let $F^{(2)}: \Fun_{\calD}^{-}(\calO, \calC) \times \Fun_{\calD}^{+}( \calO, \calC) \rightarrow \Alg^{\pm}_{\calO}(\calC)$ be the functor given by the formula $X,Y \mapsto F(X) \otimes F(Y)$. Then
$F^{(2)}$ is equivalent to the composition
$$ \Fun^{-}_{\calD}(\calO, \calC) \times \Fun^{+}_{\calD}(\calO, \calC) \simeq
\Fun^{\lax}_{\calD^{\otimes}}( {\calO'}^{\otimes} \boxplus {\calO'}^{\otimes}, \calC^{\otimes})
\stackrel{f_{02}}{\rightarrow} \Alg^{\pm}_{\calO}(\calC),$$
where $f_{02}$ denotes the functor given by operadic left Kan extension along the
correspondence of $\infty$-operads $\calT^{0}_{\calO}$ of Notation \ref{bbess}.
\end{itemize}
\end{proposition}

\begin{proof}
Assertion $(1)$ is a special case of Corollary \symmetricref{spaltwell}. To prove
$(2)$, let $p: \Delta^2 \rightarrow \Delta^1$ be the map which collapses the edge
$\Delta^{ \{0,1\} } \subseteq \Delta^2$, and let $\calM^{\otimes}$ be the subcategory of
$\calT_{\calO} \times_{ \Delta^1} \Delta^2$ spanned by those morphisms
$f: X \rightarrow Y$ satisfying the following condition: if both $X$ and $Y$ belong
to $\calM^{\otimes} \times_{ \Delta^2} \{0\}$, then the image of $f$ in $\Nerve(\FinSeg)$ is inert.
The canonical map $\calM^{\otimes} \rightarrow \Delta^2 \rightarrow \Nerve(\FinSeg)$ is a
$\Delta^2$-family of $\infty$-operads. Let $f_{ij}: \Alg_{ \calM_i}(\calC) \rightarrow
\Alg_{ \calM_j}(\calC)$ be the functor given by operadic left Kan extension
along the correspondence $\calM^{\otimes} \times_{ \Delta^2 } \Delta^{ \{i,j\} }$ for
$0 \leq i < j \leq 2$. Note that for every object $X$ of $\calM^{\otimes}_{0}$, there
exists a $q$-coCartesian morphism $X \rightarrow Y$, where $Y \in \calM^{\otimes}_{1}$;
here $q$ denotes the projection $\calM^{\otimes} \rightarrow \Delta^2$. It follows from Example
\bicatref{gabbe} that $q$ is a flat inner fibration, so we have an equivalence of functors
$f_{02} \simeq f_{12} \circ f_{01}$ (Theorem \ref{transit}).

To study the functor $f_{01}$, we note that $\calM \otimes_{ \Delta^2} \Delta^{\{0,1\} }$
is the correspondence associated to the inclusion of $\infty$-operads
$$ {\calO'}^{\otimes} \boxtimes {\calO'}^{\otimes} \rightarrow \calO^{\otimes} \boxtimes \calO^{\otimes}.$$
We have a homotopy commutative diagram of $\infty$-categories
$$ \xymatrix{ \Alg_{ \calM_0}(\calC) \ar[r]^{f_{01}} \ar[d] & \Alg_{ \calM_1}(\calC) \ar[r]^{f_{12}} \ar[d] & \Alg_{\calM_2}(\calC) \ar[d] \\
\Fun_{\calD}^{-}( \calO, \calC) \times \Fun_{\calD}^{+}( \calO, \calC) \ar[r]^-{g_{01}} & \Alg^{-}_{\calO}(\calC) \times \Alg^{+}_{\calO}(\calC) \ar[r]^-{g_{12}} & \Alg^{\pm}_{\calO}(\calC) }$$
where the vertical maps are categorical equivalences (Theorem \ref{cross}), the map
$g_{01}$ can be identified with $F^{-} \times F^{+}$ (by virtue of $(1)$), and the map $g_{12}$ can
be identified with the tensor product $\otimes$ on $\Alg^{\pm}_{\calO}(\calC)$ (Proposition \ref{alug}).
Composing these identifications, we obtain the desired description of $f_{02}$.
\end{proof}

To deduce a version of Theorem \ref{juke} from Proposition \ref{onjuke}, we would like to obtain a {\em different} description of the functor $f_{02}$ given by operadic left Kan extension along $\calT^{0}_{\calO}$. This description will also be obtained from Theorem \ref{transit}, but using a more interesting factorization of the correspondence $\calT^{0}_{\calO}$.

\begin{notation}\label{ingus}
We define categories $\calJ_{0}$, $\calJ_1$, and $\calJ_2$ as follows:
\begin{itemize}
\item[$(1)$] The category $\calJ_1$ has as objects triples
$( \seg{n}, S, T)$, where $S$ and $T$ are subsets of $\seg{n}$ which contain the base
point such that $\seg{n} = S \cup T$. A morphism from $( \seg{n}, S, T)$ to
$( \seg{n'}, S', T')$ in $\calJ_1$ consists of a map $\alpha: \seg{n} \rightarrow \seg{n'}$ in
$\FinSeg$ which restricts to inert morphisms $[S] \rightarrow [S']$, $[T] \rightarrow [T']$.

\item[$(2)$] The category $\calJ_0$ is the full subcategory of $\calJ_1$ spanned by those objects $( \seg{n}, S, T)$ for which $S \cap T = \{ \ast \}$

\item[$(3)$] The category $\calJ_2$ coincides with $\FinSeg$.
\end{itemize}

Let $\phi_{01}: \calJ_0 \rightarrow \calJ_1$ be the inclusion, let
$\phi_{12}: \calJ_1 \rightarrow \calJ_2$ be the forgetful functor, and let
$\phi_{02} = \phi_{12} \circ \phi_{01}$. We can assemble the categories
$\calJ_{i}$ into one large category $\calJ$ as follows:
\begin{itemize}
\item[$(i)$] An object of $\calJ$ is a pair $(J, i)$, where $0 \leq i \leq 2$ and $X$ is an object of
$\calJ_i$.
\item[$(ii)$] Given a pair of objects $(I,i)$, $(J, j) \in \calJ$, we have
$$ \Hom_{\calJ}( (I,i), (J,j) ) = \begin{cases} \Hom_{\calJ_j}( \phi_{ij}(I), J) & \text{ if } i < j \\
\Hom_{\calJ_j}( I, J) & \text{ if } i = j \\
\emptyset & \text{ if } i > j. \end{cases}$$
\end{itemize}

Let $\Subd'$ be the category defined in Notation \ref{skaz}. The inclusion
$\calJ_0 \subseteq \Subd$ and equivalence $\calJ_2 \simeq \FinSeg$ extend to a functor
$r: \calJ \rightarrow \Subd'$, which is given on $\calJ_1$ by the formula
$r( \seg{n}, S, T) = \seg{n}$. For every $\infty$-operad $\calO^{\otimes}$, we define
$\calM[ \calO]^{\otimes}$ to be the fiber product $\Nerve(\calJ) \times_{ \Nerve( \Subd')} \calT_{\calO}$.
Let $\Phi$ denote the composite map
$$\calM[\calO]^{\otimes} \rightarrow \Nerve(\calJ) \stackrel{r'}{\rightarrow} \Delta^2 \times \Nerve(\FinSeg).$$
\end{notation}

Repeating the proof of Lemma \ref{kwop}, we obtain the following:

\begin{lemma}\label{kwop2}
Let $\calO^{\otimes}$ be an $\infty$-operad. The map
$\Phi: \calM[ \calO]^{\otimes} \rightarrow \Delta^2 \times \Nerve(\FinSeg)$ of
Notation \ref{ingus} exhibits $\calM[\calO]^{\otimes}$ as a $\Delta^2$-family of $\infty$-operads.
\end{lemma}

For $i \in \{0,1,2\}$, we let $\calM[\calO]^{\otimes}_{i}$
denote the fiber $\calM[\calO]^{\otimes} \times_{ \Delta^2 } \{i\}$. Let
${\calO'}^{\otimes}$ denote the fiber product $\calO^{\otimes} \times_{ \Nerve(\FinSeg)} \Triv$. 
The fiber $\calM[\calO]^{\otimes}_{0}$ is isomorphic to ${\calO'}^{\otimes} \boxtimes {\calO'}^{\otimes}$, while $\calM[\calO]^{\otimes}_{2}$ is isomorphic to $\calO^{\otimes}$ itself. We will denote
the inner fiber $\calM[\calO]^{\otimes}_{1}$ by $\calQ^{\otimes}$. This $\infty$-operad is more exotic: it in some sense encodes the ``quadradic part'' of the $\infty$-operad $\calO^{\otimes}$.
Note that the fiber product $\calM[\calO]^{\otimes} \times_{ \Delta^2} \Delta^{ \{0,2\} }$ is isomorphic to the correspondence $\calT^{0}_{\calO}$ of Notation \ref{bbess}.  

To proceed further with our analysis, we need the following technical result, whose proof we defer until the end of this section:

\begin{proposition}\label{corpos}
Let $q: \calO^{\otimes} \rightarrow \Nerve(\FinSeg)$ be an $\infty$-operad satisfying the following conditions:
\begin{itemize}
\item[$(1)$] The $\infty$-operad $\calO^{\otimes}$ is coherent.
\item[$(2)$] The underlying $\infty$-category $\calO$ is a Kan complex.
\end{itemize}
Then the map $\calM[\calO]^{\otimes} \rightarrow \Delta^2$ is a flat inner fibration.
\end{proposition}

Note that if $\psi: \calT_{\calO} \rightarrow \calD^{\otimes}$ is a map of $\infty$-operad families, then
composition with $\psi$ induces another map of $\infty$-operad families
$\calM[\calO]^{\otimes} \rightarrow \calD^{\otimes}$.

\begin{corollary}\label{sup}
Let $\kappa$ be an uncountable regular cardinal, let $\calO^{\otimes}$ be an essentially $\kappa$-small $\infty$-operad, let $q: \calC^{\otimes} \rightarrow \calD^{\otimes}$ be a coCartesian fibration
of $\infty$-operads, and let $\psi: \calT_{\calO} \rightarrow \calD^{\otimes}$ be an $\infty$-operad family map. Assume that each fiber $\calC_{D}$ of $q$ admits
$\kappa$-small colimits, and that the $\calD$-monoidal structure on $\calC$ is compatible with
$\kappa$-small colimits.

For $\sigma \in \{ -, +, \pm \}$, let $F^{\sigma}: \Fun^{\sigma}_{\calD}( \calO, \calC)
\rightarrow \Alg^{\sigma}_{\calO}(\calC)$ be a left adjoint to the restriction functor. Let
Let $f_{01}: \Fun^{-}_{\calD}(\calO, \calC) \times \Fun^{+}_{\calD}(\calO, \calC) \simeq
\Alg_{ \calM[\calO]_{0}}(\calC) \rightarrow \Alg_{ \calQ}(\calC)$ be given by operadic
left Kan extension along the correspondence $\calM[\calO]^{\otimes} \times_{ \Delta^2} \Delta^{ \{0,1\} }$, and let $f_{12}: \Alg_{ \calQ}(\calC) \rightarrow \Alg^{\pm}_{\calO}(\calC)$ be the functor
given by left Kan extension along $\calM[\calO]^{\otimes} \times_{ \Delta^2} \Delta^{ \{1,2\} }$. 
If $\calO^{\otimes}$ is coherent and $\calO$ is a Kan complex, then the composition $f_{12} \circ f_{01}$ can be identified with the composite functor $$ \Fun^{-}_{\calD}(\calO, \calC) \times \Fun^{+}_{\calD}(\calO, \calC) \stackrel{F^{-} \times F^{+}}{\longrightarrow}
\Alg^{-}_{\calO}(\calC) \times \Alg^{+}_{\calO}(\calC) \stackrel{ \otimes}{\longrightarrow} \Alg^{\pm}_{\calO}(\calC).$$
\end{corollary}

\begin{proof}
Combine Proposition \ref{onjuke}, Theorem \ref{transit}, and Proposition \ref{corpos}.
\end{proof}

To apply Corollary \ref{sup}, we need to understand the functors $f_{01}$ and $f_{12}$ better.
To this end, we need to introduce some additional notation. Fix an $\infty$-operad family map
$\psi: \calT_{\calO} \rightarrow \calD^{\otimes}$, so that $\psi$ induces $\infty$-operad maps
$\psi_{-}, \psi_{+}, \psi_{\pm}: \calO^{\otimes} \rightarrow \calD^{\otimes}$. We note that
$\psi$ also determines natural transformations
$$ \psi_{-} {\rightarrow} \psi_{\pm} {\leftarrow} \psi_{+}.$$
If $q: \calC^{\otimes} \rightarrow \calD^{\otimes}$ is a coCartesian fibration of $\infty$-operads, then coCartesian transport along these transformations determines a pair of functors
$$ \Fun^{-}_{\calD}( \calO, \calC) \stackrel{u_-}{\rightarrow} \Fun^{\pm}_{\calD}( \calO, \calC)
\stackrel{ u_{+}}{\leftarrow} \Fun^{+}_{\calD}( \calO, \calC).$$

Let $\calQ_0^{\otimes}$ denote the full subcategory of $\calQ^{\otimes}$ spanned by those objects
whose image in $\Nerve(\calJ_1)$ belongs to the full subcategory
$\Nerve(\calJ_0) \subseteq \Nerve(\calJ_1)$. 
There is an evident forgetful functor
$\calQ^{\otimes}_0 \rightarrow \calM[ \calO]^{\otimes}_0$.
If $\calO^{\otimes}$ is unital, then
this forgetful functor is a trivial Kan fibration. Choosing a section, we obtain a map of
$\infty$-operads $i: \calM[\calO]^{\otimes}_{0} \rightarrow \calQ^{\otimes}$.
There is also a natural transformation $\id \rightarrow i$ of functors 
$\calM[ \calO]^{\otimes}_0 \rightarrow \calM[ \calO]^{\otimes}$.

Composition with $i$ induces a forgetful functor
$$\theta: \Alg_{\calQ}(\calC) \rightarrow \Alg_{\calM[\calO]_0}(\calC)
\simeq \Fun^{\pm}_{\calD}(\calO, \calC) \times \Fun^{\pm}_{\calD}(\calO, \calC).$$
Under the hypotheses of Corollary \ref{sup}, the functor $\theta$ has a left adjoint $\theta^{L}: \Fun^{\pm}_{\calD}(\calO, \calC) \times \Fun^{\pm}_{\calD}(\calO, \calC) \rightarrow \Alg_{ \calQ}(\calC)$, given by operadic left Kan extension along $i$. Let $F_{(2)}:
\Fun^{-}_{\calD}( \calO, \calC) \times \Fun^{+}_{\calD}( \calO, \calC) \rightarrow \Alg_{\calQ}(\calC)$
be the composition of $\theta^{L}$ with $u_{-} \times u_{+}$. Since the composition 
$\theta \circ f_{01}$ is equivalent to $u_{-} \times u_{+}$, we obtain a natural
transformation of functors $\alpha: F_{(2)} \rightarrow f_{01}$.

We would like to measure the failure of $\alpha$ to be an equivalence. To this end, consider
the fully faithful embedding $\Triv \rightarrow \Nerve(\calJ_1)$ given by the formula
$\seg{n} \mapsto ( \seg{n}, \seg{n}, \seg{n})$. This embedding determines a map of
$\infty$-operads $j: {\calO'}^{\otimes} \rightarrow \calQ^{\otimes}$. Composition
with $j$ induces a forgetful functor
$j^{\ast}: \Alg_{ \calQ}(\calC) \rightarrow \Alg^{\pm}_{\calO'}(\calC) \simeq \Fun^{\pm}_{\calD}( \calO, \calC)$.
The hypotheses of Corollary \ref{sup} guarantee that $j^{\ast}$ admits a left adjoint
$j_{!}$, given by operadic left Kan extension along $j$. 

\begin{lemma}\label{sokol}
Let $\calO^{\otimes}$, $q: \calC^{\otimes} \rightarrow \calD^{\otimes}$, and
$\psi: \calT_{\calO} \rightarrow \calD^{\otimes}$ be as in Corollary \ref{sup}, and assume
that $\calO^{\otimes}$ is unital. Then the diagram
$$ \xymatrix{ j_{!} j^{\ast} F_{(2)} \ar[r]^{ j_{!} j^{\ast} \alpha} \ar[d] & j_{!} j^{\ast} f_{01} \ar[d] \\
F_{(2)} \ar[r]^{\alpha} & f_{01} }$$
is a pushout diagram of functors from $\Fun^{-}_{\calD}(\calO, \calC) \times \Fun^{+}_{\calD}(\calO, \calC)$ to $\Alg_{ \calQ}(\calC)$. 
\end{lemma}

The proof of Lemma \ref{sokol} will be given at the end of this section. Let us accept
Lemma \ref{sokol} for the moment, and see how it leads to a version of Theorem \ref{juke}.
Note that the correspondence $\calM[\calO]^{\otimes} \times_{ \Delta^2} \Delta^{ \{1,2\} }$ is
associated to the forgetful map of $\infty$-operads $k: \calQ \rightarrow \calO^{\otimes}$. Let $k^{\ast}: \Alg^{\pm}_{ \calO}(\calC) \rightarrow \Alg_{ \calQ}(\calC)$ denote the induced map, so that $k^{\ast}$ is right adjoint to the operadic left Kan extension functor $f_{12}:
\Alg_{ \calQ}(\calC) \rightarrow \Alg^{\pm}_{\calO}(\calC)$. 
Since $f_{12}$ preserves pushouts, we deduce from Lemma \ref{sokol} the existence of a pushout
diagram
$$ \xymatrix{ f_{12} j_{!} j^{\ast} F_{(2)} \ar[r] \ar[d] & f_{12} j_{!} j^{\ast} f_{01} \ar[d] \\
j_{12} F_{(2)} \ar[r] & f_{12} f_{01}  }$$
of functors from $\Fun(\calO, \calC) \times \Fun(\calO, \calC)$ to $\Alg_{\calO}(\calC)$.
Let us identify the terms in this diagram. The functor $f_{12} j_{!}$ is left adjoint to
the forgetful functor $j^{\ast} k^{\ast}: \Alg^{\pm}_{\calO}(\calC) \rightarrow \Fun^{\pm}_{\calD}(\calO, \calC)$, and is therefore equivalent to $F^{\pm}$. The composition $j^{\ast} f_{01}$ can be identified
with the tensor product functor $\otimes: \Fun^{-}_{\calD}(\calO, \calC) \times \Fun^{+}_{\calD}(\calO, \calC) \rightarrow \Fun^{\pm}_{\calD}(\calO, \calC)$ determined by the composite map
$\calT_{ \calO'} \rightarrow \calT_{\calO} \stackrel{\psi}{\rightarrow} \calD^{\otimes}$.
In the case where $\calO^{\otimes}$ is coherent and $\calO$ is a Kan complex, Corollary \ref{sup} allows us to identify
$f_{12} f_{01}$ with the composition
$$ \Fun^{-}_{\calD}(\calO, \calC) \times \Fun^{+}_{\calD}(\calO, \calC) \stackrel{F^{-} \times F^{+}}{\rightarrow}
\Alg^{-}_{\calO}(\calC) \times \Alg^{+}_{\calO}(\calC) \stackrel{\otimes}{\longrightarrow} \Alg^{\pm}_{\calO}(\calC).$$ 
Finally, $j_{12} F_{(2)}$ is the composition of $u_{-} \times u_{+}$ with a left adjoint to $\theta k^{\ast}$, which coincides with the composition $$ \Alg^{\pm}_{\calO}(\calC) \rightarrow \Fun^{\pm}_{\calD}(\calO, \calC) \rightarrow \Fun^{\pm}_{\calD}(\calO, \calC) \times \Fun^{\pm}_{\calD}(\calO, \calC).$$
It follows that $j_{12} F_{(2)}$ can be identified with the composition of $F^{\pm}$ with
the coproduct of the functors $u_{-}$ and $u_+$. Combining this identifications, we arrive at the following conclusion:

\begin{theorem}\label{pjule}
Let $\calO^{\otimes}$ be a coherent $\infty$-operad whose underlying $\infty$-category $\calO$ is a Kan complex, $\psi: \calT_{\calO} \rightarrow \calD^{\otimes}$ a map of $\infty$-operad families, and $q: \calC^{\otimes} \rightarrow \calD^{\otimes}$ a coCartesian fibration of $\infty$-operads. Assume that there exists an uncountable regular cardinal $\kappa$ such that $\calO^{\otimes}$ is essentially $\kappa$-small, each fiber $\calC_{D}$ of $q$ admits $\kappa$-small colimits, and the $\calD$-monoidal structure on $\calC$ is compatible with $\kappa$-small colimits. Then, for every pair of objects $V \in \Fun^{-}_{\calD}( \calO, \calC),
W \in \Fun^{+}_{\calD}( \calO, \calC)$, there is a canonical pushout diagram
$$ \xymatrix{ F^{\pm}(j^{\ast} F_{(2)}(X,Y)) \ar[r] \ar[d] & F^{\pm}( X \otimes Y ) \ar[d] \\
F^{\pm}( e_{-}(X) \coprod e_{+}(Y) ) \ar[r] & F^{-}(X) \otimes F^{+}(Y) }$$
in the $\infty$-category $\Alg^{\pm}_{\calO}(\calC)$. This diagram depends functorially
on $X$ and $Y$ $($in other words, it is given by a pushout diagram of functors
from $\Fun^{-}_{\calD}(\calO, \calC) \times \Fun^{+}_{\calD}(\calO, \calC)$ to $\Alg^{\pm}_{\calO}(\calC)${}$)$. 
\end{theorem}

\begin{remark}
In the special case where $\calD^{\otimes} = \Nerve(\FinSeg)$, the superscripts in Theorem \ref{pjule} are superfluous and the functors $e_{-}$ and $e_{+}$ are equivalent to the identity. In this case,
we obtain a pushout diagram
$$ \xymatrix{ F( j^{\ast} F_{(2)}(X,Y)) \ar[r] \ar[d] & F(X \otimes Y) \ar[d] \\
F(X \coprod Y) \ar[r] & F(X) \otimes F(Y). }$$
\end{remark}

We can now proceed with the proof of our main result.

\begin{proof}[Proof of Theorem \ref{juke}]
Let $\calO^{\otimes}$ be the $\infty$-operad $\OpE{k}$, and let
$\calC^{\otimes}$ be a symmetric monoidal $\infty$-category. We will assume
that $\calC$ admits countable colimits and that the tensor product on $\calC$ preserves countable colimits separately in each variable.
Since $\calO$ is a contractible Kan complex,
evaluation at $\seg{1} \in \calO$ induces a trivial Kan fibration
$e: \Fun(\calO, \calC) \rightarrow \calC$. Let $\Free: \calC \rightarrow \Alg_{\calO}(\calC)$
be the functor obtained by composing the free functor
$F: \Fun(\calO, \calC) \rightarrow \Alg_{\calO}(\calC)$
of \S \ref{tensor2} with a section $s$ of $e$, and let
$f_{(2)}: \calC \times \calC \rightarrow \calC$ be the composition
$$ \calC \times \calC \stackrel{ s \times s}{\longrightarrow}
\Fun(\calO, \calC) \times \Fun(\calO, \calC)
\stackrel{ j^{\ast} F_{(2)}}{\longrightarrow} \Fun( \calO, \calC) \stackrel{e}{\rightarrow} \calC.$$
Unwinding the definitions, we see that $f_{(2)}$ is the colimit of the functors
$p_{!}: \calC \times \calC \rightarrow \calC$ indexed by the binary operations
$$p \in \Mul_{\calO}( \{ \seg{1}, \seg{1} \}, \seg{1} ) \simeq \Rect( \nostar{2} \times \Cube{k}, \Cube{k})
\simeq S^{k-1}.$$
Because $\calC^{\otimes}$ is symmetric monoidal, this diagram of functors is constant,
and we can identify $f_{(2)}$ with the functor $(C,D) \mapsto C \otimes D \otimes S^{k-1}$.
Invoking Theorem \ref{pjule} (note that $\calO^{\otimes}$ is coherent by Theorem \ref{cubecoh}), we obtain the desired pushout diagram
$$ \xymatrix{ \Free( C \otimes D \otimes S^{k-1}) \ar[r] \ar[d] & \Free(C \otimes D) \ar[d] \\
\Free( C \coprod D) \ar[r] & \Free(C) \otimes \Free(D) }$$
in the $\infty$-category $\Alg_{\calO}(\calC)$.
\end{proof}

We conclude this section with the proofs of Lemma \ref{sokol} and Proposition \ref{corpos}.

\begin{proof}[Proof of Lemma \ref{sokol}]
Let $\overline{\calC}^{\otimes}$ denote the fiber product
$\calC^{\otimes} \times_{ \calD^{\otimes}} \calQ^{\otimes}$, and let
$q': \overline{\calC}^{\otimes} \rightarrow \calQ^{\otimes}$ denote the projection map.
We observe that $\Alg_{\calQ}(\calC)$ can be identified with a full subcategory of
$\Fun_{ \calQ^{\otimes}}( \calQ^{\otimes}, \overline{\calC}^{\otimes})$. We will prove
that for every pair of objects $X \in \Fun^{-}_{\calD}(\calO, \calC)$
and $Y \in \Fun^{+}_{\calD}( \calO, \calC)$, the diagram
$$ \xymatrix{ j_{!} j^{\ast} F_{(2)}(X,Y) \ar[r] \ar[d] & j_{!} j^{\ast} f_{01}(X,Y) \ar[d] \\
F_{(2)}(X,Y) \ar[r] & f_{01}(X,Y) }$$
is a pushout in $\Fun_{ \calQ^{\otimes}}( \calQ^{\otimes}, \overline{\calC}^{\otimes})$. 
In view of Lemma \monoidref{surtybove}, it will suffice to show that for
each object $Q \in \calQ^{\otimes}$, the diagram $\sigma_{Q}:$
$$ \xymatrix{ (j_{!} j^{\ast} F_{(2)}(X,Y))(Q) \ar[r] \ar[d] & (j_{!} j^{\ast} f_{01}(X,Y))(Q) \ar[d] \\
F_{(2)}(X,Y)(Q) \ar[r] & f_{01}(X,Y)(Q) }$$
is a $q'$-colimit diagram in $\overline{\calC}^{\otimes}$. 

Since $q'$ is a pullback of $q: \calC^{\otimes} \rightarrow \calD^{\otimes}$, it is a coCartesian fibration. For every morphism $\beta: Q \rightarrow Q'$ in $\calQ^{\otimes}$, let
$\beta_{!}: \overline{\calC}^{\otimes}_{Q} \rightarrow \overline{\calC}^{\otimes}_{Q'}$ be the induced map. To prove that $\sigma_{Q}$ is a $q'$-colimit diagram, it will suffice to show that
each of the diagrams $\beta_{!}(\sigma_{Q})$ is a pushout diagram in
$\overline{\calC}^{\otimes}_{Q'}$ (Proposition \toposref{relcolfibtest}). 
Let $\seg{n}$ denote the image of $Q'$ in $\Nerve(\FinSeg)$, and choose inert morphisms
$\gamma(i): Q' \rightarrow Q'_{i}$ lying over $\Colp{i}: \seg{n} \rightarrow \seg{1}$ for
$1 \leq i \leq n$. Since $q'$ is a coCartesian fibration of $\infty$-operads, the functors
$\gamma(i)_{!}$ induce an equivalence
$$ \overline{\calC}^{\otimes}_{Q'} \rightarrow \prod_{1 \leq i \leq n} \overline{\calC}^{\otimes}_{Q'_{i}}.$$
It follows that $\beta_{!}(\sigma_{Q})$ is a pushout diagram if and only if each
$\gamma(i)_{!} \beta_{!}(\sigma_{Q})$ is a pushout diagram in $\overline{\calC}^{\otimes}_{Q'_{i}}$. 
We may therefore replace $\beta$ by $\gamma(i) \circ \beta$ and thereby reduce to the case
$Q' \in \calQ$. The map $\beta$ factors as a composition
$$ Q \stackrel{ \beta'}{\rightarrow} Q'' \stackrel{\beta''}{\rightarrow} Q'$$
where $\beta'$ is inert and $\beta''$ is active. Note that $\beta'_{!}( \sigma_{Q})$ is equivalent
to $\sigma_{Q''}$; we may therefore replace $Q$ by $Q''$ and thereby reduce to the case where
$\beta$ is active. Together these conditions imply that $Q$ belongs to the image of
either $j: {\calO'}^{\otimes} \rightarrow \calQ^{\otimes}$ or the essential image of $i: {\calO'}^{\otimes} \boxplus {\calO'}^{\otimes} \rightarrow \calQ^{\otimes}$. We consider each case in turn.

Suppose that $Q$ belongs to the image of $j$. We claim in this case that
the vertical maps in the diagram $\sigma_{Q}$ are equivalences. It follows that
$\beta_{!}( \sigma_{Q})$ has the same property, and is therefore automatically a pushout diagram.
Our claim is a special case of the following more general assertion: let $U$ be an arbitrary
object of $\Alg_{\calQ}(\calC)$: then the counit map $j_{!} j^{\ast} U \rightarrow U$ induces an equivalence $(j_{!} j^{\ast} U)(Q) \rightarrow U(Q)$. The functor $j_{!}$ is computed by
the formation of operadic $q$-left Kan extension: consequently, $(j_{!} j^{\ast} U)(Q)$ is an operadic
$q$-colimit of the diagram 
$$ (\calM[\calO]^{\otimes}_{0})^{\acti}_{/Q} \rightarrow
\calQ^{\otimes} \stackrel{U}{\rightarrow} \calC^{\otimes}.$$
The desired assertion now follows from the observation that $Q$ belongs to the
image of $j$, so that $( \calM[ \calO]^{\otimes}_{0})^{\acti}_{/Q}$ contains $Q$ as a final object.

Suppose instead that $Q$ belongs to the essential image of $i$. We claim in this case
that the horizontal maps in the diagram $\sigma_{Q}$ are equivalences. As before, it follows
that $\beta_{!} \sigma_{Q}$ has the same property and is therefore automatically a pushout diagram.
To prove the claim, we first show that for {\em any} map $U \rightarrow V$ in
$\Alg_{ \calM[\calO]_0}(\calC) \simeq \Fun^{\pm}_{\calD}(\calO, \calC) \times \Fun^{\pm}_{\calD}(\calO, \calC)$, the induced map 
$\xi: (j_{!} U)(Q) \rightarrow (j_! V)(Q)$ is an equivalence in $\calC^{\otimes}$. To see this, we observe
that $(j_{!} U)(Q)$ and $(j_! V)(Q)$ are given by operadic $q$-colimits of diagrams
$$ \calM[\calO]^{\otimes}_0 \times_{ \calQ^{\otimes} } (\calQ^{\otimes})^{\acti}_{/Q} \rightarrow \calM[\calO]^{\otimes}_0 \rightarrow \calC^{\otimes}.$$
To prove that $\xi$ is an equivalence, it suffices to show that these diagrams are equivalent.
The assumption that $Q$ belongs to the image of $i$ guarantees that every object of
$\calM[\calO]^{\otimes}_0 \times_{ \calQ^{\otimes} } (\calQ^{\otimes})^{\acti}_{/Q}$ lies over $\seg{0} \in \Nerve(\FinSeg)$: the desired result now follows from the observation that every morphism in
$\calC^{\otimes}_{\seg{0}}$ is an equivalence, since $\calC^{\otimes}_{\seg{0}}$ is a (contractible) Kan complex.

To complete the proof, we must show that the map $\xi': F_{(2)}(X,Y)(Q) \rightarrow f_{01}(X,Y)(Q)$
is an equivalence whenever $Q$ belongs to the image of $i$. Let $U \in \Alg_{ \calM[\calO]_0}(\calC)$ be a preimage of $(X,Y)$ under the equivalence
$\Alg_{ \calM[\calO]_0}(\calC) \rightarrow \Fun^{-}_{\calD}(\calO, \calC) \times \Fun^{+}_{\calD}(\calO, \calC)$. Then $F_{(2)}(X,Y)(Q)$ and $f_{01}(X,Y)(Q)$ can be identified with the operadic $q$-colimits of diagrams
$$ \calM[\calO]^{\otimes}_0 \times_{ \calQ^{\otimes} } (\calQ^{\otimes})^{\acti}_{/Q} \rightarrow \calM[\calO]^{\otimes}_0 \stackrel{U}{\rightarrow} \calC^{\otimes}$$
$$ \calM[\calO]^{\otimes}_0 \times_{ \calM[\calO]^{\otimes}_{/Q} }
( \calM[\calO]^{\otimes})^{\acti}_{/Q} \rightarrow \calM[\calO]^{\otimes}_0
\stackrel{U}{\rightarrow} \calC^{\otimes}.$$
To prove that $\xi'$ is an equivalence, it suffices to show that the functor
$$ \epsilon: \calM[\calO]^{\otimes}_0 \times_{ \calQ^{\otimes} } (\calQ^{\otimes})^{\acti}_{/Q}
\rightarrow 
\calM[\calO]^{\otimes}_0 \times_{ \calM[\calO]^{\otimes}_{/Q} }
( \calM[\calO]^{\otimes})^{\acti}_{/Q} \rightarrow \calM[\calO]^{\otimes}_0$$
is a categorical equivalence. Both the domain and codomain of $\epsilon$ are
right-fibered over $(\calM[\calO]^{\otimes}_{0})^{\acti}$: it will therefore suffice
to show that $\epsilon$ induces a homotopy equivalence after passing to
the fiber over any object $P \in \calM[\calO]^{\otimes}_{0}$ (Corollary \toposref{usefir}).
Unwinding the definitions, we must show that the canonical map
$$\bHom_{ \calQ^{\otimes}}( i(P), Q) \rightarrow \bHom_{ \calM[\calO]^{\otimes}}(P, Q)$$
is a homotopy equivalence. This follows from a simple calculation, using our assumptions that
$\calO^{\otimes}$ is unital and that $Q$ belongs to the essential image of $i$.
\end{proof}

\begin{proof}[Proof of Proposition \ref{corpos}]
Let $\calO^{\otimes}$ be a coherent $\infty$-operad such that $\calO$ is a Kan complex; we wish to prove that the inner fibration $\calM[ \calO]^{\otimes} \rightarrow \Delta^2$ is flat.. Fix an object $X \in \calM[\calO]^{\otimes}_{0}$, corresponding to a pair of objects
$X_{-}, X_{+} \in \calO^{\otimes}$, and let $Z \in \calO^{\otimes} \simeq \calM[ \calO]^{\otimes}_{2}$.
Suppose we are given a morphism $X \rightarrow Z$ in $\calM[\calO]^{\otimes}$.
We wish to prove that the $\infty$-category $\calC = \calM[\calO]^{\otimes}_{X/ \, /Z} \times_{ \Delta^2} \{1\}$ is weakly contractible. Let $\emptyset \in \calO^{\otimes}_{\seg{0}}$ be a final object of
$\calO^{\otimes}$. Since $\calO^{\otimes}$ is unital, $\emptyset$ is also an initial object
of $\calO^{\otimes}$. We can therefore choose a diagram $\sigma:$
$$ \xymatrix{ \emptyset \ar[r] \ar[d] & X_{-} \\
X_{+} & }$$
in $\calO^{\otimes}_{/Z}$. Let $\overline{\calC}[\sigma]$ denote the full subcategory of
$\calO^{\otimes}_{\sigma/ \, /Z}$ spanned by those diagrams
$$ \xymatrix{ \emptyset \ar[r] \ar[d] & X_{-} \ar[d]^{f} \\
X_{+} \ar[r]^{g} & Q }$$
where $f$ and $g$ are both semi-inert (by virtue of $(2)$, this is equivalent to the requirement
that $q(f)$ and $q(g)$ are semi-inert morphisms in $\Nerve(\FinSeg)$) and the map
$q(X_{-}) \coprod q(X_{+}) \rightarrow q(Q)$ is a surjection. Using the fact that
$\emptyset$ is an initial object of $\calO^{\otimes}$, we obtain a trivial Kan fibration
$\overline{\calC}[\sigma] \rightarrow \calC$. It will therefore suffice to show that
$\overline{\calC}[\sigma]$ is weakly contractible.

Let $q( X_{+} ) = \seg{m}$. The proof will proceed by induction on $m$. If $m = 0$, then
$\overline{\calC}[\sigma]$ has an initial object and there is nothing to prove. Otherwise, the map
$\emptyset \rightarrow X_{+}$ factors as a composition
$$ \emptyset \rightarrow X'_{+} \stackrel{\alpha}{\rightarrow} X_{+},$$
where $q(\alpha)$ is an inclusion $\seg{m-1} \hookrightarrow \seg{m}$.
Let $\tau: \Delta^1 \coprod_{ \{0\} }\Delta^2 \rightarrow \calO^{\otimes}_{/Z}$ denote the diagram $X_{-} \leftarrow \emptyset \rightarrow X'_{+} \rightarrow X_{+}$, and let
$\tau_0 = \tau | ( \Delta^1 \coprod_{ \{0\} } \Delta^1)$. Let $\calD$ denote the full subcategory
of the fiber product $$\Fun( \Delta^1, (\calO^{\otimes}_{/Z})_{\tau_0/}) \times_{ \Fun( \{1\}, ( \calO^{\otimes}_{/Z})_{\tau_0/}}
(\calO^{\otimes}_{/Z})_{\tau/}$$ spanned by those diagrams
$$ \xymatrix{ \emptyset \ar[r] \ar[d] & X_{-} \ar[d] \\
X'_{+} \ar[r] \ar[d] & Q' \ar[d] \\
X_{+} \ar[r] & Q }$$
in $\calO^{\otimes}_{/Z}$ where the maps $X_{-} \rightarrow Q'$, $X'_{+} \rightarrow Q'$,
and $Q' \rightarrow Q$ are semi-inert and the map $q(X_{-}) \coprod_{ q(\emptyset)}
q( X'_{+}) \rightarrow q(Q')$ is surjective.
Let $\calD_0 \subseteq \calD$ be the full subcategory spanned by those diagrams for which the map
$X_{+} \rightarrow Q$ is semi-inert and the map
$q(Q') \coprod_{ q(X'_{+})} q(X_{+}) \rightarrow q(Q)$ is surjective. We have canonical maps
$$ \overline{\calC}[\sigma] \stackrel{\phi}{\leftarrow} \calD_0 \subseteq \calD \stackrel{\psi}{\rightarrow} 
\overline{\calC}[ \tau_0].$$ 
The map $\phi$ admits a right adjoint and is therefore a weak homotopy equivalence,
and the simplicial set $\overline{\calC}[\tau_0]$ is weakly contractible by the inductive hypothesis.
The inclusion $\calD_0 \subseteq \calD$ admits a right adjoint, and is therefore a weak homotopy equivalence. To complete the proof, it will suffice to show that $\psi$ is a weak homotopy equivalence. We have a homotopy pullback
diagram
$$ \xymatrix{ \calD \ar[d]^{\psi} \ar[r] & (\calK_{\calO})_{\alpha/ \, / \id_{Z} } \ar[d]^{\psi'} \\
\overline{\calC}[\tau_0] \ar[r] & \calO^{\otimes}_{ X'_{+}/ \, / Z}.} $$
The
coherence of $\calO^{\otimes}$ guarantees that the map $\calK_{\calO} \rightarrow \calO^{\otimes}$ is a flat inner fibration, so that $\psi'$ satisfies the hypotheses of Lemma \ref{clapsos} (see Example \ref{stublos}). It follows that $\psi$ is a weak homotopy equivalence, as desired.
\end{proof}

\subsection{Nonunital Algebras}\label{pluy}

Let $A$ be an abelian group equipped with a commutative and associative multiplication
$m: A \otimes A \rightarrow A$. A {\it unit} for the multiplication $m$ is an element
$1 \in A$ such that $1a = a$ for each $a \in A$. If there exists a unit for $A$, then that unit is unique
and $A$ is a commutative ring (with unit). Our goal in this section is to prove an analogous result, where
the category of abelian groups is replaced by an arbitrary symmetric monoidal $\infty$-category
$\calC$ (Corollary \ref{suz}).

An analogous result for associative algebras was proven in \S \monoidref{giddug}. Namely, we proved
that if $A$ is a nonunital associative algebra object of a monoidal $\infty$-category $\calC$ which
is {\em quasi-unital} (Definition \ref{laster}), then $A$ can be promoted (in an essentially unique fashion) to an associative algebra with unit (Theorem \monoidref{uniqueunit}). Roughly speaking, the idea is to realize $A$ as the algebra of endomorphisms of $A$, regarded as a right module over itself. This proof does not immediately generalize to the commutative context, since the endomorphism algebra of an $A$-module is noncommutative in general. We will therefore take a somewhat different approach: rather than trying to mimic the proof of Theorem \monoidref{uniqueunit}, we will combine Theorem \monoidref{uniqueunit} with Theorem \ref{slide} to deduce an analogous result for $\OpE{k}$-algebras (Theorem \ref{quas}). We then obtain the result for commutative algebras by passing to the limit $k \rightarrow \infty$. 
We begin with a discussion of nonunital algebras in general.

\begin{definition}\label{placa}
Let $\Surj$ denote the subcategory of $\Sect$ containing all objects of
$\Sect$, such that a morphism $\alpha: \seg{m} \rightarrow \seg{n}$ belongs to
$\Surj$ if and only if it is surjective. If $\calO^{\otimes}$ is an $\infty$-operad, we let $\calO^{\otimes}_{\nunit}$ denote the fiber
product $\calO^{\otimes} \times_{ \Nerve(\FinSeg)} \Nerve(\Surj)$. If $\calC^{\otimes} \rightarrow \calO^{\otimes}$ is a fibration of $\infty$-operads, we let $\Alg^{\nunit}_{\calO}(\calC)$
denote the $\infty$-category $\Alg_{ \calO_{\nunit}}(\calC)$ of $\calO_{\nunit}$-algebra
objects of $\calC$; we will refer to $\Alg^{\nunit}_{\calO}(\calC)$ as the {\it $\infty$-category
of nonunital $\calO$-algebra objects of $\calC$}.
\end{definition}

Our goal is to show that if $\calO^{\otimes}$ is a little $k$-cubes operad $\OpE{k}$
for some $k \geq 1$, then the $\infty$-category $\Alg_{\calO}^{\nunit}(\calC)$ of nonunital
$\calO$-algebra objects of $\calC$ is not very different from the $\infty$-category
$\Alg_{\calO}(\calC)$ of unital $\calO$-algebras objects of $\calC$. More precisely, we will
show that the restriction functor $\Alg_{\calO}(\calC) \rightarrow \Alg_{\calO}^{\nunit}(\calC)$
induces an equivalence of $\Alg_{\calO}(\calC)$ onto a subcategory 
$\Alg_{\calO}^{\qunit}(\calC) \subseteq \Alg_{\calO}^{\nunit}(\calC)$ whose objects are required to
admit units up to homotopy and whose morphisms are required to preserve those units
(see Definition \ref{laster} below). Our next step is to define the $\infty$-categories $\Alg_{\calO}^{\qunit}(\calC)$ more precisely.

\begin{definition}\label{laster}
Let $k \geq 1$, let $q: \calC^{\otimes} \rightarrow \OpE{k}$ be a coCartesian fibration of
$\infty$-operads, and let $A \in \Alg_{\OpE{k}}^{\nunit}( \calC)$; we will abuse notation by identifying
$A$ with its image in the underlying $\infty$-category $\calC$. 

Let ${\bf 1}$ denote a unit object of $\calC$. The multiplication map
$A \otimes A \rightarrow A$ induces an associative multiplication
$$m: \Hom_{ \h{\calC}}( {\bf 1},A) \times \Hom_{ \h{\calC}}( {\bf 1}, A) \rightarrow
\Hom_{ \h{\calC}}( {\bf 1}, A).$$
We will say that morphism $e: {\bf 1} \rightarrow A$ is a {\it quasi-unit} for $A$
if its homotopy class $[e]$ is both a left and a right unit with respect to the multiplication $m$.
We will say that $A$ is {\it quasi-unital} if it admits a quasi-unit $e: {\bf 1} \rightarrow A$.

Let $f: A \rightarrow B$ be a morphism between nonunital $\OpE{k}$-algebra objects of $\calC$,
and assume that $A$ admits a quasi-unit $e: {\bf 1} \rightarrow A$. We will say that
$f$ is {\it quasi-unital} if the composite map $f \circ e: {\bf 1} \rightarrow B$ is a quasi-unit for $B$;
in this case, $B$ is also quasi-unital. We let $\Alg_{\OpE{k}}^{\qunit}(\calC)$ denote the subcategory
of $\Alg_{ \OpE{k}}^{\nunit}(\calC)$ spanned by the quasi-unital algebras and quasi-unital morphisms between them.
\end{definition}

\begin{remark}
In the situation of Definition \ref{laster}, a map $e: {\bf 1} \rightarrow A$ is a quasi-unit for $A$
if and only if each of the composite maps
$$ A \simeq {\bf 1} \otimes A \stackrel{e \otimes \id}{\longrightarrow} A \otimes A \stackrel{m}{\rightarrow} A \quad \quad A \simeq A \otimes {\bf 1} \stackrel{ \id \otimes e}{\longrightarrow} A \otimes A \stackrel{m}{\rightarrow} A$$
is homotopic to the identity. If $k > 1$, then the multiplication on $A$ and the tensor product on $\calC$ are commutative up to homotopy, so these conditions are equivalent to one another.
\end{remark}

\begin{remark}
Let $A \in \Alg_{ \OpE{k}}^{\nunit}(\calC)$ be as in Definition \ref{laster}. Then a quasi-unit
$e: {\bf 1} \rightarrow A$ is uniquely determined up to homotopy, if it exists. Consequently,
the condition that a map of nonunital $\OpE{k}$-algebras $f: A \rightarrow B$ be quasi-unital 
is independent of the choice of $e$.
\end{remark}

\begin{example}
Let $q: \calC^{\otimes} \rightarrow \OpE{k}$ be a coCartesian fibration of $\infty$-operads,
and let $\theta: \Alg_{ \OpE{k}}(\calC) \rightarrow \Alg^{\nunit}_{\OpE{k}}(\calC)$ be the
restriction functor. Then $\theta$ carries $\OpE{k}$-algebra objects of
$\calC$ to quasi-unital objects of $\Alg^{\nunit}_{\OpE{k}}(\calC)$, and morphisms
of $\OpE{k}$-algebras to quasi-unital morphisms in $\Alg^{\nunit}_{\OpE{k}}(\calC)$.
Consequently, $\theta$ can be viewed as a functor from
$\Alg_{\OpE{k}}(\calC)$ to $\Alg^{\qunit}_{\OpE{k}}(\calC)$.
\end{example}

Our main result about nonunital algebras is the following:

\begin{theorem}\label{quas}
Let $k \geq 1$ and let $q: \calC^{\otimes} \rightarrow \OpE{k}$ be a coCartesian fibration
of $\infty$-operads. Then the forgetful functor $\theta: \Alg_{ \OpE{k}}(\calC) \rightarrow
\Alg_{ \OpE{k}}^{\qunit}(\calC)$ is an equivalence of $\infty$-categories.
\end{theorem}

The proof of Theorem \ref{quas} is somewhat elaborate, and will be given at the end of this section.

\begin{remark}\label{proquas}
In the situation of Theorem \ref{quas}, we may assume without loss of generality
that $\calC^{\otimes}$ is small (filtering $\calC^{\otimes}$ if necessary).
Using Proposition \symmetricref{skimmy}, we deduce the
existence of a presentable $\OpE{k}$-monoidal $\infty$-category
$\calD^{\otimes} \rightarrow \OpE{k}$ and a fully faithful
$\OpE{k}$-monoidal functor $\calC^{\otimes} \rightarrow \calD^{\otimes}$.
We have a commutative diagram
$$ \xymatrix{ \Alg_{ \OpE{k}}( \calC) \ar[d]^{\theta} \ar[r] & \Alg_{ \OpE{k}}( \calD) \ar[d]^{\theta'} \\
\Alg_{\OpE{k}}^{\qunit}(\calC) \ar[r] & \Alg_{ \OpE{k}}^{\qunit}(\calD) }$$
where the horizontal maps are fully faithful embeddings, whose essential images consist
of those (unital or nonunital) $\OpE{k}$-algebra objects of $\calD$ whose underlying object
belongs to the essential image of the embedding $\calC \hookrightarrow \calD$. To prove that $\theta$ is a categorical equivalence, it suffices to show that $\theta'$ is a categorical equivalence.
In other words, it suffices to prove Theorem \ref{quas} in the special case where
$\calC^{\otimes}$ is a presentable $\OpE{k}$-monoidal $\infty$-category.
\end{remark}

We can use Theorem \ref{quas} to deduce an analogous assertion regarding commutative
algebras. Let $\calC^{\otimes}$ be a symmetric monoidal $\infty$-category. We let
$\CAlg^{\nunit}(\calC)$ denote the $\infty$-category $\Alg_{\CommOp}^{\nunit}(\calC)$ of nonunital
commutative algebra objects of $\calC$. Definition \ref{laster} has an evident analogue
for nonunital commutative algebras and maps between them: we will say that a 
nonunital commutative algebra $A \in \CAlg^{\nunit}(\calC)$ is {\it quasi-unital} if there
exists a map $e: {\bf 1} \rightarrow A$ in $\calC$ such that the composition
$$ A \simeq {\bf 1} \otimes A \stackrel{e \otimes \id}{\longrightarrow} A \otimes A \rightarrow A$$
is homotopic to the identity (in the $\infty$-category $\calC$). In this case, $e$ is uniquely determined up to homotopy and we say that $e$ is a quasi-unit for $A$; a morphism $f: A \rightarrow B$
in $\CAlg^{\nunit}(\calC)$ is {\it quasi-unital} if $A$ admits a quasi-unit $e: {\bf 1} \rightarrow A$
such that $f \circ e$ is a quasi-unit for $B$. The collection of quasi-unital commutative algebras
and quasi-unital morphisms between them can be organized into a subcategory
$\CAlg^{\qunit}(\calC) \subseteq \CAlg^{\nunit}(\calC)$.

\begin{corollary}\label{suz}
Let $\calC^{\otimes}$ be a symmetric monoidal $\infty$-category. Then the forgetful functor
$\CAlg(\calC) \rightarrow \CAlg^{\qunit}(\calC)$ is an equivalence of $\infty$-categories.
\end{corollary}

\begin{proof}
In view of Corollary \ref{infec}, we have an equivalence of $\infty$-operads
$\OpE{\infty} \rightarrow \Nerve(\FinSeg)$. It will therefore suffice to show that
the forgetful functor $\Alg_{ \OpE{\infty}}(\calC) \rightarrow \Alg^{\qunit}_{\OpE{\infty}}(\calC)$
is an equivalence of $\infty$-categories. This map is the homotopy inverse limit
of a tower of forgetful functors $\theta_{k}: \Alg_{ \OpE{k}}(\calC) \rightarrow \Alg^{\qunit}_{\OpE{k}}(\calC)$, each of which is an equivalence of $\infty$-categories by Theorem \ref{quas}.
\end{proof}

As a first step toward understanding the forgetful functor
$\theta: \Alg_{\calO}(\calC) \rightarrow \Alg^{\nunit}_{\calO}(\calC)$, let us study the left adjoint
to $\theta$. In classical algebra, if $A$ is a nonunital ring, then we can canonically enlarge
$A$ to a unital ring by considering the product $A \oplus \Z$ endowed with the multiplication
$(a,m)(b,n) = (ab + mb + na, mn)$. Our next result shows that this construction works quite generally:

\begin{proposition}\label{calba}
Let $\calO^{\otimes}$ be a unital $\infty$-operad, let
$q: \calC^{\otimes} \rightarrow \calO^{\otimes}$ be a coCartesian fibration of $\infty$-operads
which is compatible with finite coproducts, and let
$\theta: \Alg_{\calO}(\calC) \rightarrow \Alg_{\calO}^{\nunit}(\calC)$ be the forgetful functor. Then:
\begin{itemize}
\item[$(1)$] For every object $A \in \Alg_{\calO}^{\nunit}(\calC)$, there exists another object
$A^{+} \in \Alg_{\calO}(\calC)$ and a map $A \rightarrow \theta( A^{+})$ which exhibits
$A^{+}$ as a free $\calO$-algebra generated by $A$.
\item[$(2)$] A morphism $f: A \rightarrow \theta(A^{+})$ in $\Alg_{\calO}^{\nunit}(\calC)$ exhibits
$A^{+}$ as a free $\calO$-algebra generated by $A$ if and only if, for every object
$X \in \calO$, the map $f_X: A(X) \rightarrow A^{+}(X)$ and the unit map
${\bf 1}_{X} \rightarrow A^{+}(X)$ exhibit $A^{+}(X)$ as a coproduct of $A(X)$ and
the unit object ${\bf 1}_{X}$ in the $\infty$-category $\calC_{X}$.
\item[$(3)$] The functor $\theta$ admits a left adjoint.
\end{itemize}
\end{proposition}

\begin{proof}
For every object $X \in \calO$, the $\infty$-category $\calD = \calO^{\otimes}_{\nunit}
\times_{\calO^{\otimes} } (\calO^{\otimes})_{/X}^{\acti}$ can be written as a disjoint
union of $\calD_0 = (\calO^{\otimes}_{\nunit})_{/X}^{\acti}$ with the full subcategory $\calD_1 \subseteq \calD$
spanned by those morphisms $X' \rightarrow X$ in $\calO^{\otimes}$ where
$X' \in \calO^{\otimes}_{\seg{0}}$. The $\infty$-category $\calD_0$ contains $\id_{X}$ as
a final object. Since $\calO^{\otimes}$ is unital, the $\infty$-category $\calD_1$ is a contractible Kan complex containing a vertex $v: X_0 \rightarrow X$. It follows that the inclusion
$\{ \id_X, v \}$ is cofinal in $\calD$. Assertions $(1)$ and $(2)$ now follow from
Proposition \symmetricref{wax1} (together with Propositions \symmetricref{chocolateoperad} and \symmetricref{optest}). Assertion $(3)$ follows from $(1)$ (Corollary \symmetricref{slavetime}).
\end{proof}

In the stable setting, there is a close relationship between nonunital algebras and augmented algebras.
To be more precise, we need to introduce a bit of terminology.

\begin{definition}
Let $q: \calC^{\otimes} \rightarrow \calO^{\otimes}$ be a coCartesian fibration of $\infty$-operads,
and assume that $\calO^{\otimes}$ is unital. An {\it augmented $\calO$-algebra object of $\calC$}
is a morphism $f: A \rightarrow A_0$ in $\Alg_{\calO}(\calC)$, where $A_0$ is a trivial algebra.
We let $\Alg_{\calO}^{\aug}(\calC)$ denote the full subcategory of $\Alg_{\calO}(\calC)$ spanned by
the augmented $\calO$-algebra objects of $\calC$.
\end{definition}

The following result will not play a role in the proof of Theorem \ref{quas}, but is of some independent interest:

\begin{proposition}\label{propers}
Let $q: \calC^{\otimes} \rightarrow \calO^{\otimes}$ be a coCartesian fibration of
$\infty$-operads. Assume that $\calO^{\otimes}$ is unital and that $q$ exhibits
$\calC$ as a stable $\calO$-monoidal $\infty$-category. Let $F: \Alg^{\nunit}_{\calO}(\calC)
\rightarrow \Alg_{\calO}(\calC)$ be a left adjoint to the forgetful functor
$\theta: \Alg_{\calO}(\calC) \rightarrow \Alg^{\nunit}_{\calO}(\calC)$. 
Let $0 \in \Alg^{\nunit}_{\calO}(\calC)$ be a final object, so that $F(0) \in \Alg_{\calO}(\calC)$
is a trivial algebra (Proposition \ref{calba}). Then $F$ induces an equivalence of
$\infty$-categories
$$ T: \Alg^{\nunit}_{\calO}(\calC) \simeq \Alg^{\nunit}_{\calO}(\calC)^{/ 0}
\rightarrow \Alg^{\aug}_{\calO}(\calC).$$
\end{proposition}

\begin{proof}
Let $p: \calM \rightarrow \Delta^1$ be a correspondence associated to the adjunction
$\Adjoint{F}{ \Alg_{\calO}^{\nunit}(\calC)}{\Alg_{\calO}(\calC)}{\theta}$. Let
$\calD$ denote the full subcategory of $\Fun_{ \Delta^1}( \Delta^1 \times \Delta^1, \calM)$ spanned
by those diagrams $\sigma$
$$ \xymatrix{ A \ar[d] \ar[r]^{f} & A^{+} \ar[d]^{g} \\
A_0 \ar[r]^{f_0} & A_0^{+} }$$
where $A_0$ is a final object of $\Alg_{\calO}^{\nunit}(\calC)$ and the maps $f$
and $f'$ are $p$-coCartesian; this (together with Proposition \ref{calba}) guarantees
that $A_0^{+} \in \Alg_{\calO}(\calC)$ is a trivial algebra so that $g$ can be regarded
as an augmented $\calO$-algebra object of $\calC$. Using Proposition \toposref{lklk}, we deduce
that the restriction functor $\sigma \mapsto A$ determines a trivial Kan fibration
$\calD \rightarrow \Alg^{\nunit}_{\calO}(\calC)$. By definition, the functor $T$ is obtained
by composing a section of this trivial Kan fibration with the restriction map
$\phi: \calD \rightarrow \Alg^{\aug}_{\calO}(\calC)$ given by $\sigma \mapsto g$.
To complete the proof, it will suffice to show that $\phi$ is a trivial Kan fibration.

Let $K$ denote the full subcategory of $\Delta^1 \times \Delta^1$ obtained by removing
the object $(0,0)$, and let $\calD_0$ be the full subcategory of
$\Fun_{\Delta^1}(K, \calM)$ spanned by those diagrams
$$ A^{+} \stackrel{g}{\rightarrow} A_0^{+} \stackrel{f_0}{\leftarrow} A_0$$
where $A_0$ is a final object of $\Alg^{\nunit}_{\calO}(\calC)$ and $A_0^{+}$ is a trivial
$\calO$-algebra object of $\calC$; note that this last condition is equivalent to the requirement
that $f_0$ be $p$-coCartesian. The functor $\phi$ factors as a composition
$$ \calD \stackrel{\phi'}{\rightarrow} \calD_0 \stackrel{\phi''}{\rightarrow} \Alg^{\aug}_{\calO}(\calC).$$
We will prove that $\phi'$ and $\phi''$ are trivial Kan fibrations.

Let $\calD_1$ be the full subcategory of $\Fun_{\Delta^1}(\Delta^1, \calM)$ spanned by
the $p$-coCartesian morphisms $f_0: A_0 \rightarrow A_0^{+}$ where
$A_0$ is a final object of $\Alg^{\nunit}_{\calO}(\calC)$. It follows from Proposition
\toposref{lklk} that the restriction map $f_0 \mapsto A_0$ determines a trivial Kan fibration
from $\calD_1$ to the contractible Kan complex of final objects in $\Alg^{\nunit}_{\calO}(\calC)$, so
that $\calD_1$ is contractible. The restriction map $f_0 \mapsto A_0^{+}$ is a categorical
fibration $\overline{\phi}''$ from $\calD_1$ onto the contractible Kan complex of initial objects of
$\Alg_{\calO}(\calC)$. It follows that $\overline{\phi}''$ is a trivial Kan fibration. The map
$\phi''$ is a pullback of $\overline{\phi}''$, and therefore also a trivial Kan fibration.

We now complete the proof by showing that $\phi'$ is a trivial Kan fibration.
In view of Proposition \toposref{lklk}, it will suffice to show that a diagram
$\sigma \in \Fun_{\Delta^1}(\Delta^1 \times \Delta^1, \calM)$ belongs to
$\calD$ if and only if $\sigma_0 = \sigma | K$ belongs to $\calD_0$ and
$\sigma$ is a $p$-right Kan extension of $\sigma_0$. Unwinding the definitions
(and using Corollary \symmetricref{slimycamp2}), we are reduced to showing that
if we are given a diagram
$$ \xymatrix{ A \ar[d] \ar[r]^{f} & A^{+} \ar[d]^{g} \\
A_0 \ar[r]^{f_0} & A_0^{+} }$$
where $A_0$ is a final object of $\Alg^{\nunit}_{\calO}(\calC)$ and $A_0^{+}$ is a trivial
algebra, then $f$ is $p$-coCartesian if and only if the induced diagram
$$ \xymatrix{ A(X) \ar[r]^{f_X} \ar[d] & A^{+}(X) \ar[d] \\
A_0(X) \ar[r] & A_0^{+}(X) }$$
is a pullback square in $\calC_{X}$, for each $X \in \calO$. Since $\calC_{X}$ is
a stable $\infty$-category, this is equivalent to the requirement that the induced map
$\psi: \coker(f_X) \rightarrow A_0^{+}(X)$ is an equivalence. The map $\psi$ fits into a commutative diagram
$$ \xymatrix{ {\bf 1}_X \ar[rr] \ar[d] & & {\bf 1}_{X} \ar[d] \\
A^{+}(X) \ar[r] & \coker(f) \ar[r] & A_0^{+}(X) }$$
where the vertical maps are given by the units for the algebra objects $A^{+}$ and $A_0^{+}$.
Since $A_0^{+}(X)$ is a trivial algebra, the unit map ${\bf 1}_{X} \rightarrow A_0^{+}(X)$
is an equivalence. Consequently, it suffices to show that $f$ is $p$-coCartesian
if and only if each of the composite maps ${\bf 1}_{X} \rightarrow A^{+}(X) \rightarrow \coker(f)$
is an equivalence. We have a pushout diagrm
$$ \xymatrix{ {\bf 1}_{X} \coprod A(X) \ar[r] \ar[d] & A^{+}(X) \ar[d] \\
{\bf 1}_X \ar[r] & \coker(f). }$$
Since $\calC_X$ is stable, the lower horizontal map is an equivalence if and only if
the upper horizontal map is an equivalence. The desired result now follows immediately
from the criterion described in Proposition \ref{calba}.
\end{proof}

Let us now return to the proof of Theorem \ref{quas}. We begin by treating the case $k=1$.
Without loss of generality, we may assume that $q: \calC^{\otimes} \rightarrow \OpE{1}$
is the pullback of a coCartesian fibration of $\infty$-operads $\calD^{\otimes} \rightarrow \Ass$
(Example \ref{sulta}). Let $\phi: \Nerve( \cDelta)^{op} \rightarrow \Ass$ be defined as in
Construction \symmetricref{urpas}, so that the pullback of $\calD^{\otimes}$ by $\phi$ determines a monoidal structure on the $\infty$-category $\calD$. The map $\phi$ restricts to a functor $\phi_0: \Nerve( \cDelta_{s})^{op} \rightarrow \Ass \times_{ \Nerve(\FinSeg)} \Nerve(\Surj)$.
Composition with $\phi_0$ determines a functor
$\Alg_{\Ass}^{\nunit}(\calD) \rightarrow \Alg^{\nunit}(\calD)$ (see \S \monoidref{digunit}).
We have the following nonunital analogue of Proposition \symmetricref{algass}:

\begin{proposition}
Let $q:\calD^{\otimes} \rightarrow \Ass$ be a coCartesian fibration of $\infty$-operads.
Then the functor $\Alg_{\Ass}^{\nunit}(\calD) \rightarrow \Alg^{\nunit}(\calD)$ constructed above is an equivalence of $\infty$-categories.
\end{proposition}

\begin{proof}
Let $\CatAss^{\nunit}$ denote the subcategory $\CatAss \times_{ \FinSeg} \Surj$.
We define a category $\calI$ as follows:
\begin{itemize}
\item[$(1)$] An object of $\calI$ is either an object of $\cDelta^{op}_s$ or an object of $\CatAss^{\nunit}$.
\item[$(2)$] Morphisms in $\calI$ are give by the formulas
$$\Hom_{\calI}( [m], [n] ) = \Hom_{ \cDelta^{op}_s }( [m], [n]) \quad \Hom_{\calI}( \seg{m}, \seg{n} ) = \Hom_{\CatAss^{\nunit}}( \seg{m}, \seg{n} )$$
$$ \Hom_{ \calI}( \seg{m}, [n] ) = \Hom_{ \CatAss^{\nunit}}( \seg{m}, \phi_0( [n] ) ) \quad
 \Hom_{ \calI}( [n], \seg{m} ) = \emptyset.$$ 
\end{itemize}
where $\phi_0: \cDelta^{op}_{s} \rightarrow \CatAss^{\nunit}$ is the functor defined above.
We observe that $\phi_0$ extends to a retraction $r: \calI \rightarrow \CatAss^{\nunit}$.
Let $\overline{\Alg}(\calD)$ denote the full subcategory of 
$\Fun_{ \Ass}( \Nerve(\calI), \calD^{\otimes})$ consisting of those functors $f: \Nerve(\calI) \rightarrow \calD^{\otimes}$ 
such that $q \circ f = \overline{\psi}$ and the following additional conditions are satisfied:
\begin{itemize}
\item[$(i)$] For each $n \geq 0$, $f$ carries the canonical map $\seg{n} \rightarrow [n]$ in
$\calI$ to an equivalence in $\calD^{\otimes}$.
\item[$(ii)$] The restriction $f | \Nerve(\cDelta_s)^{op}$ belongs to $\Alg^{\nunit}(\calD)$.
\item[$(ii')$] The restriction $f | \Nerve( \CatAss^{\nunit})$ is a nonunital $\Ass$-algebra object of $\calC$.
\end{itemize}

If $(i)$ is satisfied, then $(ii)$ and $(ii')$ are equivalent to one another. Moreover, $(i)$ is equivalent to the assertion that $f$ is a $q$-left Kan extension of $f | \Nerve( \CatAss^{\nunit})$.
Since every functor $f_0: \Nerve(\CatAss^{\nunit}) \rightarrow \calD^{\otimes}$ admits a $q$-left Kan extension (given, for example, by $f_0 \circ r$), Proposition \toposref{lklk} implies that the restriction map
$p: \overline{\Alg}(\calD) \rightarrow \Alg_{\Ass}^{\nunit}(\calD)$ is a trivial Kan fibration. The map
$\theta$ is the composition of a section to $p$ (given by composition with $r$) 
and the restriction map $p': \overline{\Alg}(\calD) \rightarrow \Alg^{\nunit}(\calD)$. It will therefore suffice to show that $p'$ is a trivial fibration. In view of Proposition \toposref{lklk}, this will follow from the following pair of assertions:

\begin{itemize}
\item[$(a)$] Every $f_0 \in \Alg^{\nunit}(\calD)$ admits a $q$-right Kan extension
$f \in \Fun_{ \Ass}( \Nerve(\calI), \calD^{\otimes})$.
\item[$(b)$] Given $f \in \Fun_{ \Ass}( \Nerve(\calI), \calD^{\otimes})$
such that $f_0 = f | \Nerve(\cDelta_s^{op})$ belongs to $\Alg^{\nunit}(\calD)$, $f$ is a $q$-right Kan extension of $f_0$ if and only if $f$ satisfies condition $(i)$ above.
\end{itemize}

To prove $(a)$, we fix an object $\seg{n} \in \CatAss^{\nunit}$. Let $\calJ$ denote the category
$\cDelta^{op}_{s} \times_{ \CatAss^{\nunit}} ( \CatAss^{\nunit})_{\seg{n}/ },$
and let $g$ denote the composition
$ \Nerve(\calJ) \rightarrow \Nerve(\cDelta^{op}_{s}) \rightarrow \calD^{\otimes}.$
According to Lemma \toposref{kan2}, it will suffice to show that $g$ admits a $q$-limit in $\calD^{\otimes}$ (for each $n \geq 0$). The objects of $\calJ$ can be identified with surjective
morphisms $\alpha: \seg{n} \rightarrow \phi_0([m])$ in $\CatAss$. Let $\calJ_0 \subseteq \calJ$
denote the full subcategory spanned by those objects for which $\alpha$ is inert.
The inclusion $\calJ_0 \subseteq \calJ$ has a right adjoint, so that $\Nerve(\calJ_0)^{op} \rightarrow \Nerve(\calJ)^{op}$ is cofinal. Consequently, it will suffice to show that $g_0 = g | \Nerve(\calJ_0)$ admits a $q$-limit in $\calD^{\otimes}$.

Let $\calJ_1$ denote the full subcategory of $\calJ_0$ spanned by the morphisms
$\Colp{j}: \seg{n} \rightarrow \phi_0( [1] )$. Using our assumption that
$f_0$ is a nonunital algebra object of $\calD$, we deduce that $g_0$ is a $q$-right Kan extension of
$g_1 = g_0 | \Nerve(\calJ_1)$. In view of Lemma \toposref{kan0}, it will suffice to show that
the map $g_1$ has a $q$-limit in $\calD^{\otimes}$. But this is clear; our assumption that $f_0$ belongs to $\Alg^{\nunit}(\calD)$ guarantees that $f_0$ exhibits $f_0([n])$ as a $q$-limit of $g_1$. This proves $(a)$. Moreover, the proof shows that $f$ is a $q$-right Kan extension of $f_0$ at $\seg{n}$ if and only if $f$ induces an equivalence $f( \seg{n} ) \rightarrow f( [n] )$; this immediately implies $(b)$ as well.
\end{proof}

Given a coCartesian fibration of $\infty$-operads $\calD^{\otimes} \rightarrow \Ass$, we let
$\Alg_{\Ass}^{\qunit}(\calD)$ denote the fiber product $\Alg_{\Ass}^{\nunit}(\calD) \times_{ \Alg^{\nunit}(\calD)} \Alg^{\qunit}(\calD)$, where $\Alg^{\qunit}(\calD)$ is defined as in \S \monoidref{digunit}. We have a commutative diagram
$$ \xymatrix{ \Alg(\calD) \ar[d]^{\theta'} & \Alg_{\Ass}(\calD) \ar[l] \ar[r] \ar[d] & \Alg_{ \OpE{1}}(\calD) \ar[d]^{\theta} \\
\Alg^{\qunit}(\calD) & \Alg_{\Ass}^{\qunit}(\calD) \ar[l] \ar[r] & \Alg_{ \OpE{1}}^{\qunit}(\calD) }$$
in which the horizontal maps are categorical equivalences. Theorem \monoidref{uniqueunit} implies that
$\theta'$ is an equivalence of $\infty$-categories, so that $\theta$ is likewise an equivalence of
$\infty$-categories. This proves Theorem \ref{quas} in the special case $k=1$.

The proof of Theorem \ref{quas} in general will proceed by induction on $k$. For the remainder
of this section, we will fix an integer $k \geq 1$, and assume that Theorem \ref{quas} has been
verified for the $\infty$-operad $\OpE{k}$. Our goal is to prove that Theorem \ref{quas} is valid
also for $\OpE{k+1}$. 
Fix a coCartesian fibration of $\infty$-operads $q: \calC^{\otimes} \rightarrow \OpE{k+1}$; we wish to show that the forgetful functor $\theta: \Alg_{\OpE{k+1}}(\calC) \rightarrow
\Alg_{ \OpE{k+1}}^{\qunit}(\calC)$ is an equivalence of $\infty$-categories. In view of
Remark \ref{proquas}, we can assume that $\calC^{\otimes}$ is a presentable
$\OpE{k+1}$-monoidal $\infty$-category.

We begin by constructing a left homotopy inverse to $\theta$. Consider the bifunctor of
$\infty$-operads $\OpE{1} \times \OpE{k} \rightarrow \OpE{k+1}$ of \S \ref{sass1}.
Using this bifunctor, we can define $\OpE{1}$-monoidal $\infty$-categories
$\Alg_{ \OpE{k}}( \calC)^{\otimes}$ and $\Alg_{\OpE{k}}^{\nunit}(\calC)^{\otimes}$.
Moreover, the collection of quasi-unital $\OpE{k}$-algebras and quasi-unital morphisms
between them are stable under tensor products, so we can also consider an $\OpE{1}$-monoidal
subcategory $\Alg_{ \OpE{k}}^{\qunit}(\calC)^{\otimes} \subseteq \Alg_{ \OpE{k}}^{\nunit}(\calC)^{\otimes}$.
Similarly, we have $\OpE{k}$-monoidal $\infty$-categories
$\Alg_{ \OpE{1}}(\calC)^{\otimes}$, $\Alg_{ \OpE{1}}^{\nunit}(\calC)^{\otimes}$, and
$\Alg^{\qunit}_{\OpE{1}}(\calC)^{\otimes}$.

There is an evident forgetful functor $\Alg_{\OpE{k+1}}^{\nunit}( \calC) \rightarrow
\Alg_{ \OpE{1}}^{\nunit}( \Alg_{ \OpE{k}}^{\nunit}(\calC))$, which obviously restricts to
a functor $\psi_0: \Alg_{\OpE{k+1}}^{\qunit}(\calC) \rightarrow \Alg_{ \OpE{1}}^{\nunit}( \Alg_{\OpE{k}}^{\qunit}(\calC))$. Using the inductive hypothesis (and Corollary \toposref{usefir}), we deduce that the
evident categorical fibration $\Alg_{ \OpE{k}}(\calC)^{\otimes} \rightarrow \Alg_{\OpE{k}}^{\qunit}(\calC)^{\otimes}$ is a categorical equivalence and therefore a trivial Kan fibration.
It follows that the induced map $\Alg_{ \OpE{1}}^{\nunit}( \Alg_{\OpE{k}}(\calC)) \rightarrow
\Alg_{ \OpE{1}}^{\nunit}( \Alg_{ \OpE{j}}^{\qunit}(\calC))$ is a trivial Kan fibration, which admits
a section $\psi_1$. Let $\psi_2$ be the evident equivalence
$\Alg_{ \OpE{1}}^{\nunit}( \Alg_{ \OpE{k}}(\calC)) \simeq \Alg_{ \OpE{k}}( \Alg^{\nunit}_{\OpE{1}}(\calC))$.
We observe that the composition $\psi_2 \circ \psi_1 \circ \psi_0$ carries
$\Alg_{ \OpE{k+1}}^{\qunit}( \calC)$ into the subcategory
$\Alg_{ \OpE{k}}( \Alg^{\qunit}_{\OpE{1}}(\calC)) \subseteq \Alg_{ \OpE{k}}( \Alg^{\nunit}_{\OpE{1}}(\calC))$.
Using the inductive hypothesis and Corollary \toposref{usefir} again, we deduce that
the forgetful functor $\Alg_{ \OpE{k}}( \Alg_{\OpE{1}}(\calC)) \rightarrow \Alg_{ \OpE{k}}( \Alg_{\OpE{1}}^{\qunit}(\calC))$ is a trivial Kan fibration, which admits a section
$\psi_3$. Finally, Theorem \ref{slide} implies that the functor
$\Alg_{\OpE{k+1}}(\calC) \rightarrow \Alg_{ \OpE{1}}( \Alg_{ \OpE{k}}(\calC))$ is an equivalence
of $\infty$-categories which admits a homotopy inverse $\psi_{4}$. Let
$\psi$ denote the composition $\psi_{4} \psi_3 \psi_2 \psi_1 \psi_0$. Then
$\psi$ is a functor from $\Alg^{\qunit}_{\OpE{k+1}}(\calC)$ to $\Alg_{\OpE{k+1}}(\calC)$.
The composition $\psi \circ \theta$ becomes homotopic to the identity after composing with
the functor $\Alg_{\OpE{k+1}}(\calC) \simeq \Alg_{ \OpE{1}}^{\qunit}( \Alg^{\qunit}_{\OpE{k}}(\calC))
\subseteq \Alg^{\nunit}_{\OpE{1}}( \Alg^{\nunit}_{\OpE{k}}(\calC))$, and is therefore homotopic to the identity
on $\Alg_{ \OpE{k+1}}(\calC)$. 

To complete the proof of Theorem \ref{quas}, it will suffice to show that
the composition $\theta \circ \psi$ is equivalent to the identity functor from
$\Alg^{\qunit}_{\OpE{k+1}}(\calC)$ to itself. This is substantially more difficult, and the proof will require a brief digression. In what follows, we will assume that the reader is familiar with the theory of centralizers of maps of
$\OpE{k}$-algebras developed in \S \ref{clapser} (see Definition \ref{saigd}). 

\begin{definition}\label{cullstone}
Let $\calC^{\otimes} \rightarrow \OpE{k}$ be a coCartesian fibration of $\infty$-operads,
let $A$ and $B$ be $\OpE{k}$-algebra objects of $\calC$, and let
$u: {\bf 1} \rightarrow A$ be a morphism in $\calC$. 
We let $\bHom_{ \Alg_{\OpE{k}}(\calC)}^{u}(A,B)$ be the summand of the mapping
space $\bHom_{ \Alg_{\OpE{k}}(\calC)}(A,B)$ given by those maps
$f: A \rightarrow B$ such that $f \circ u$ is an invertible element in the monoid
$\Hom_{ \h{\calC}}( {\bf 1}, B)$.

Let $f: A \rightarrow B$ be a morphism in $\Alg_{ \OpE{k}}(\calC)$ and let
$u: {\bf 1} \rightarrow A$ be as above. We will say that $f$ is a {\it $u$-equivalence} if,
for every object $C \in \Alg_{ \OpE{k}}(\calC)$, composition with $f$ induces
a homotopy equivalence
$$\bHom_{ \Alg_{\OpE{k}}(\calC)}^{fu}(B,C) \rightarrow \bHom_{ \Alg_{\OpE{k}}(\calC)}^{u}(A,C).$$
\end{definition}

\begin{remark}\label{soad}
Let $M$ be an associative monoid. If $x$ and $y$ are commuting elements of $M$, then the product
$xy = yx$ is invertible if and only if both $x$ and $y$ are invertible. In the situation of
Definition \ref{cullstone}, this guarantees that if $u: {\bf 1} \rightarrow A$ and
$v: {\bf 1} \rightarrow A$ are morphisms in $\calC$ such that $u$ and $v$ commute
in the monoid $\Hom_{ \h{\calC}}( {\bf 1}, A)$ and $w$ denotes the product map
${\bf 1} \simeq {\bf 1} \otimes {\bf 1} \stackrel{u \otimes v}{\longrightarrow} A \otimes A \rightarrow A$, then
we have
$\bHom_{ \Alg_{\OpE{k}}(\calC)}^{w}(A,B) = \bHom_{ \Alg_{\OpE{k}}(\calC)}^{u}(A,B)
\cap \bHom_{ \Alg_{\OpE{k}}(\calC)}^{v}(A,B)$ (where the intersection is formed in the mapping
space $\bHom_{ \Alg_{\OpE{k}}(\calC)}(A,B)$).
It follows that if $f: A \rightarrow B$ is a $u$-equivalence or a $v$-equivalence, then it is also a $w$-equivalence.
\end{remark}

\begin{remark}
Let $\calC^{\otimes} \rightarrow \OpE{k}$ be a presentable $\OpE{k}$-monoidal $\infty$-category, 
and let $e: {\bf 1} \rightarrow A$ be the unit map for an $\OpE{k}$-algebra object
$A \in \Alg_{\OpE{k}}(\calC)$. We will abuse notation by identifying
$A$ with the underlying nonunital $\OpE{k}$-algebra object, and let
$A^{+}$ be the free $\OpE{k}$-algebra generated by this nonunital
$\OpE{k}$-algebra (see Proposition \ref{calba}). Let $e^{+}$ denote the composite map
${\bf 1} \rightarrow A \rightarrow A^{+}$. Then the counit map $v: A^{+} \rightarrow A$
is an $e^{+}$-equivalence. To see this, it suffices to show that for every object
$B \in \Alg_{\OpE{k}}(\calC)$, composition with $v$ induces a homotopy equivalence
$$ \bHom_{ \Alg_{\OpE{k}}(\calC)}(A,B) = \bHom_{ \Alg_{\OpE{k}}(\calC)}^{e}(A,B)
\rightarrow \bHom_{ \Alg_{\OpE{k}}(\calC)}^{e^{+}}(A^{+},B).$$
Note that any nonunital algebra morphism $f: A \rightarrow B$ carries
$e$ to an idempotent element $[f \circ e]$ of the monoid $\Hom_{ \h{\calC}}({\bf 1} ,B)$,
so $f \circ e$ is a quasi-unit for $B$ if and only if $[f \circ e]$ is invertible.
Consequently, the homotopy equivalence $\bHom_{ \Alg_{\OpE{k}}}(A^{+},B) \simeq
\bHom_{ \Alg_{\OpE{k}}^{\nunit}}(A,B)$ induces an identification
$ \bHom_{ \Alg_{\OpE{k}}(\calC)}^{e^{+}}(A^{+},B) \simeq \bHom_{ \Alg_{\OpE{k}}^{\qunit}(\calC)}(A,B)$.
The desired result now follows from Theorem \ref{quas}.
\end{remark}

\begin{lemma}\label{tab2}
Let $q: \calC^{\otimes} \rightarrow \OpE{k+1}$ be a presentable $\OpE{k+1}$-monoidal $\infty$-category, so that $\Alg_{\OpE{k}}(\calC)$ inherits the structure of an $\OpE{1}$-monoidal
$\infty$-category. Let $f: A \rightarrow A'$ be a morphism in $\Alg_{ \OpE{k}}(\calC)$,
and let $u: {\bf 1} \rightarrow A$ be a morphism in $\calC$ such that $f$ is a $u$-equivalence.
Let $B \in \Alg_{ \OpE{k}}(\calC)$ and $v: { \bf 1} \rightarrow B$ be an arbitrary morphism in $\calC$.
Then:
\begin{itemize}
\item[$(1)$] The induced map $f \otimes \id_{B}$ is a $u \otimes v: {\bf 1} \rightarrow A \otimes B$
equivalence.
\item[$(2)$] The induced map $\id_{B} \otimes f$ is a $v \otimes u: {\bf 1} \rightarrow B \otimes A$-equivalence.
\end{itemize}
\end{lemma}

\begin{proof}
We will prove $(1)$; the proof of $(2)$ is similar. Let $e_{A}: {\bf 1} \rightarrow A$ and
$e_{B}: {\bf 1} \rightarrow B$ denote the units of $A$ and $B$, respectively. We note that
$u \otimes v$ is homotopic to the product of maps $e_{A} \otimes v$ and
$u \otimes e_{B}$ which commute in the monoid $\Hom_{ \h{\calC}}({\bf 1}, A \otimes B)$.
By virtue of Remark \ref{soad}, it will suffice to show that $f \otimes \id_{B}$ is
a $w$-equivalence, where $w = u \otimes e_{B}$. 

Let $w'$ be the composition of $w$ with $f \otimes \id_{B}$, and let
$C \in \Alg_{\OpE{k}}(\calC)$. We have a commutative diagram
$$ \xymatrix{ \bHom_{ \Alg_{\OpE{k}}(\calC)}^{w'}( A' \otimes B, C) \ar[rr] \ar[dr] & &
\bHom_{ \Alg_{\OpE{k}}(\calC}^{w}(A \otimes B, C) \ar[dl] \\
& \bHom_{ \Alg_{\OpE{k}}(\calC)}( B, C) & }$$
and we wish to show that the horizontal map is a homotopy equivalence.
It will suffice to show that this map induces a homotopy equivalence after passing
to the homotopy fibers over any map $g: B \rightarrow C$. This is equivalent to the requirement that $f$ induces a homotopy
equivalence
$\bHom_{\Alg_{\OpE{k}}(\calC)}^{fu}( A', \Centt{\OpE{k}}{g} ) \rightarrow \bHom_{ \Alg_{\OpE{k}}(\calC)}^{u}(A, \Centt{\OpE{k}}{g} )$, which follows from our assumption that $f$ is a $u$-equivalence.
\end{proof}

\begin{lemma}\label{sidew}
Let $\calC^{\otimes} \rightarrow \OpE{k}$ be a presentable $\OpE{k}$-monoidal
$\infty$-category, let $A \in \Alg_{ \OpE{k}}(\calC)$, and let $u: {\bf 1} \rightarrow A$ be a morphism in
the underlying $\infty$-category $\calC$. Then there exists a morphism
$f: A \rightarrow A[u^{-1}]$ in $\Alg_{ \OpE{k}}(\calC)$ with the following universal properties:
\begin{itemize}
\item[$(1)$] The map $f$ is a $u$-equivalence.
\item[$(2)$] The composite map $fu$ is a unit in the monoid $\Hom_{ \h{\calC}}( {\bf 1}, A[u^{-1}])$.
\end{itemize}
\end{lemma}

\begin{proof}
Let $P: \Alg_{\OpE{k}}(\calC) \rightarrow \Mon_{\OpE{k}}(\SSet)$ be the functor
described in \S \ref{kloop}. The inclusion $\Mon_{\OpE{k}}^{\glike}(\SSet) \subseteq
\Mon_{\OpE{k}}(\SSet)$ admits a right adjoint $G$ which can be described informally as follows:
$G$ carries an $\OpE{k}$-space $X$ to the subspace $X^{\glike} \subseteq X$ given by the
union of those connected components of $X$ which are invertible in $\pi_0 X$. 
Let $J: \Mon_{\OpE{k}}(\SSet) \rightarrow \SSet$ be the forgetful functor, and let
$\chi: \Alg_{\OpE{k}}(\calC) \rightarrow \SSet$ be the functor corepresented by $A$.
We can identify $u$ with a point in the space $JP(A)$, which determines natural
transformation of functors $\chi \rightarrow JP$. Let $\chi$ denote the fiber product
$\chi \times_{ JP} JGP$ in the $\infty$-category $\Fun( \Alg_{\OpE{k}}(\calC), \SSet)$.
Since $\chi$, $J$, $G$, and $P$ are all accessible functors which preserve small limits,
the functor $\chi'$ is accessible and preserves small limits, and is therefore corepresentable
by an object $A[u^{-1}] \in \Alg_{ \OpE{k}}(\calC)$ (Proposition \toposref{representableprime}). 
The evident map $\chi' \rightarrow \chi$ induces a map $f: A \rightarrow A[u^{-1}]$ which
is easily seen to have the desired properties.
\end{proof}

\begin{remark}\label{saug}
Let $\calC^{\otimes} \rightarrow \OpE{k}$ be as in Lemma \ref{sidew}, 
let $f: A \rightarrow B$ be a morphism in $\Alg_{\OpE{k}}(\calC)$ and let
$u: {\bf 1} \rightarrow A$ be a morphism in $\calC$. Then $f$ is a $u$-equivalence
if and only if it induces an equivalence $A[u^{-1}] \rightarrow B[ (fu)^{-1}]$ in the
$\infty$-category $\Alg_{\OpE{k}}(\calC)$.
\end{remark}

\begin{example}\label{treed}
Let $A \in \Alg^{\qunit}_{\OpE{k}}(\calC)$ be a nonunital algebra equipped with a quasi-unit
$e_A: {\bf 1} \rightarrow A$. Let $A^{+} \in \Alg_{\OpE{k}}(\calC)$ be an algebra equipped with
a nonunital algebra map $\beta: A \rightarrow A^{+}$ which exhibits $A^{+}$ as the free $\OpE{k}$-algebra generated by $A$. Then the composite map $\gamma_0: A \rightarrow A^{+} \rightarrow A^{+}[ (\beta e_{A})^{-1}]$ is quasi-unital, and therefore (by Theorem \ref{quas}) lifts to an
$\OpE{k}$-algebra map $\gamma: \overline{A} \rightarrow A^{+}[ (\beta e_A)^{-1}]$. Using
Theorem \ref{quas} again, we deduce that $\gamma$ is an equivalence in $\Alg_{\OpE{k}}(\calC)$,
so that $\gamma_0$ is an equivalence of nonunital algebras.
\end{example}

We now return to the proof of Theorem \ref{quas} for a presentable
$\OpE{k+1}$-monoidal $\infty$-category $\calC^{\otimes} \rightarrow \OpE{k+1}$.
We will assume that Theorem \ref{quas} holds for the $\infty$-operad $\OpE{k}$, so that
the forgetful functor $\Alg_{\OpE{k}}(\calC) \rightarrow \Alg^{\qunit}_{\OpE{k}}(\calC)$ is
an equivalence of $\infty$-categories. Consequently, all of the notions defined above
for $\OpE{k}$-algebras make sense also in the context of quasi-unital $\OpE{k}$-algebras; we will make use of this observation implicitly in what follows.

Let $\calD$ denote the fiber product
$$ \Fun( \bd \Delta^1, \Alg_{ \OpE{k}}^{\qunit}(\calC)) \times_{ \Fun( \bd \Delta^1,
\Alg_{ \OpE{k}}^{\nunit}(\calC))} \Fun( \Delta^1, \Alg_{ \OpE{k}}^{\nunit}(\calC))$$
whose objects are nonunital maps $f: A \rightarrow B$ between quasi-unital
$\OpE{k}$-algebra objects of $\calC$, and whose morphisms
are given by commutative diagrams
$$ \xymatrix{ A \ar[r]^{f} \ar[d] & B \ar[d] \\
A' \ar[r]^{f'} & B' }$$
where the vertical maps are quasi-unital. Let $\calD_0$ denote the full subcategory
$\Fun( \Delta^1, \Alg^{\qunit}_{\OpE{k}}(\calC)) \subseteq \calD$ spanned by the
quasi-unital maps $f: A \rightarrow B$. The inclusion $\calD_0 \hookrightarrow \calD$
admits a left adjoint $L$, given informally by the formula $(f: A \rightarrow B) \mapsto
(A \rightarrow B[ (fe_A)^{-1}])$, where $e_{A}: {\bf 1} \rightarrow A$ denotes the unit of $A$.
Using Remark \ref{saug}, we deduce the following:

\begin{lemma}\label{tab}
If $\alpha$ is a morphism in $\calD$ corresponding to a commutative diagram
$$ \xymatrix{ A \ar[r]^{f} \ar[d]^{g} & B \ar[d]^{g'} \\
A' \ar[r] & B', }$$
then $L(\alpha)$ is an equivalence if and only if the following pair of conditions is
satisfied:
\begin{itemize}
\item[$(i)$] The map $g$ is an equivalence.
\item[$(ii)$] The map $g'$ is an $f e_A$-equivalence, where $e_{A}: {\bf 1} \rightarrow A$
denotes a quasi-unit for $A$.
\end{itemize}
\end{lemma}

Note that the $\OpE{1}$-monoidal structure on $\Alg_{\OpE{k}}^{\nunit}(\calC)$ induces an
$\OpE{1}$-monoidal structure on the $\infty$-category $\calD$. 

\begin{lemma}\label{stor}
The localization functor $L: \calD \rightarrow \calD_0 \subseteq \calD$ is compatible with
the $\OpE{1}$-monoidal structure on $\calD$. In other words, if $\alpha: D \rightarrow D'$ is an
$L$-equivalence in $\calD$ and $E$ is any object of $\calD$, then the induced maps
$D \otimes E \rightarrow D' \otimes E$ and $E \otimes D \rightarrow E \otimes D'$ are
again $L$-equivalences.
\end{lemma}

\begin{proof}
Combine Lemmas \ref{tab} and \ref{tab2}.
\end{proof}

Combining Lemma \ref{stor} with Proposition \symmetricref{localjerk2}, we deduce
that $L$ can be promoted to an $\OpE{1}$-monoidal functor from
$\calD$ to $\calD_0$; in particular, $L$ induces a functor
$L': \Alg_{ \OpE{1}}^{\nunit}( \calD) \rightarrow \Alg_{ \OpE{1}}^{\nunit}(\calD_0)$ which
is left adjoint to the inclusion and therefore comes equipped with a natural transformation
$\alpha: \id_{ \Alg_{\OpE{1}}^{\nunit}(\calD)} \rightarrow L'$.

We are now ready to complete the proof of Theorem \ref{quas}. Let
$G: \Alg_{ \OpE{k+1}}( \calC) \rightarrow \Alg^{\nunit}_{\OpE{k+1}}(\calC)$ denote the forgetful
functor, let $F$ be a left adjoint to $G$ (Proposition \ref{propers}), and let
$\beta: \id_{ \Alg_{\OpE{k+1}}^{\nunit}(\calC)} \rightarrow G \circ F$ be a unit transformation.
Let $j: \Alg^{\qunit}_{\OpE{k+1}}(\calC) \rightarrow \Alg^{\nunit}_{\OpE{k+1}}(\calC)$ be the inclusion
functor and let $\xi: \Alg^{\nunit}_{\OpE{k+1}}( \calC) \rightarrow \Alg_{\OpE{1}}^{\nunit}( \Alg_{\OpE{k}}^{\nunit}(\calC))$ be the forgetful functor. If $A \in \Alg_{\OpE{k+1}}^{\nunit}(\calC)$ is quasi-unital,
then $GF(A)$ is likewise quasi-unital. Consequently, the construction
$A \mapsto \xi( \beta_{A})$ induces a functor
$\epsilon: \Alg^{\qunit}_{ \OpE{k+1}}(\calC) \rightarrow \Alg^{\nunit}_{\OpE{1}}(\calD)$. Let
$L'$ and $\alpha: \id \rightarrow L'$ be defined as above. The induced natural transformation
$\epsilon \rightarrow L' \epsilon$ can be regarded as a functor
from $\Alg^{\qunit}_{\OpE{k+1}}(\calC)$ to $\Fun( \Delta^1 \times \Delta^1, \Alg^{\nunit}_{\OpE{1}}( \Alg_{\OpE{k}}^{\nunit}(\calC)))$. This functor can be described informally as follows: it carries
a quasi-unital algebra $A$ to the diagram
$$ \xymatrix{ A \ar[r]^{\beta_A} \ar[d] & F(A) \ar[d] \\
A \ar[r] & F(A)[ (\beta_{A} e_A)^{-1}], }$$
where $e_{A}: {\bf 1} \rightarrow A$ denotes the quasi-unit of $A$. It follows
from Example \ref{treed} that the lower horizontal map is an equivalence.
Consequently, the above functor can be regarded as a natural transformation from
$\xi GF j$ to $\xi j$ in the $\infty$-category $\Fun( \Alg_{\OpE{k}}^{\qunit}(\calC),
\Alg_{\OpE{1}}^{\nunit}( \Alg^{\qunit}_{\OpE{k}}(\calC)))$. Composing with
$\psi_4 \circ \psi_3 \circ \psi_2 \circ \psi_1$, we obtain a natural transformation
$\delta:  \psi \theta F \rightarrow \psi$ of functors from $\Alg_{\OpE{k+1}}^{\qunit}(\calC)$
to $\Alg_{\OpE{k}}(\calC)$. Since $\psi \theta$ is homotopic to the identity, we can
view $\delta$ as a natural transformation from $F | \Alg_{\OpE{k+1}}^{\qunit}(\calC)$ to
$\psi$. This transformation is adjoint to a map of functors $\id_{ \Alg_{\OpE{k+1}}^{\qunit}(\calC)} \rightarrow \theta \circ \psi$. It is easy to see that this transformation is an equivalence
(using the fact that the forgetful functor $\Alg_{\OpE{k+1}}^{\qunit}(\calC) \rightarrow \calC$
is conservative, by Corollary \symmetricref{jumunj22}), so that $\psi$ is a right homotopy
inverse to $\theta$. This completes the proof of Theorem \ref{quas}.

\section{Applications of Left Modules}\label{sec2}

In \cite{monoidal}, we saw that there is a good theory of (associative) algebra objects and
(left and right) modules over them in an arbitrary monoidal $\infty$-category $\calC^{\otimes} \rightarrow \Nerve(\FinSeg)$. Our goal in this section is to describe some connections between this theory and the more general theory of $\infty$-operads developed in \cite{symmetric}, and to obtain some basic results about
the little cubes $\infty$-operads $\OpE{k}$ as a consequence. 

We will begin in \S \ref{conl} by showing that if $\calC^{\otimes}$ is a monoidal $\infty$-category
and $A$ is an algebra object of $\calC^{\otimes}$, then the $\infty$-category $\Mod_{A}^{R}(\calC)$ of
right $A$-module objects of $\calC$ is left-tensored over $\calC$. Moreover, the construction
$(\calC^{\otimes}, A) \rightarrow (\calC^{\otimes}, \Mod_{A}^{R}(\calC) )$ determines a functor
$\Theta$ from the $\infty$-category of monoidal $\infty$-categories equipped with an algebra object
to the $\infty$-category of monoidal $\infty$-categories equipped with a left module (Construction \ref{spaas}). 

The module category $\Mod_{A}^{R}(\calC)$ has a canonical object, given by the algebra $A$ itself.
Since right and left multiplication commute with one another, we can regard $A$ as a {\em left} $A$-module object in $\Mod_{A}^{R}(\calC)$. In fact, $\Mod_{A}^{R}(\calC)$ is in some sense {\em freely} generated by this left $A$-module object (Theorem \ref{bix2}). We will prove this result, together with the formally similar
Theorem \ref{bix1}, in \S \ref{probix}. In \S \ref{probix2}, we will apply these results to prove a number of formal properties of the functor $\Theta$, which establish a close connection between the theory of algebras objects and left module categories for a symmetric monoidal $\infty$-category $\calC^{\otimes}$.

In \S \ref{clapser}, we will explain how to formulate the theory of algebras and their left (or right) modules in terms of the general formalism of $\infty$-operads developed in \cite{symmetric}. More precisely, we will define an $\infty$-operad $\LMod$, and show that giving an $\LMod$-algebra object of a symmetric monoidal
$\infty$-category $\calC^{\otimes}$ is equivalent to giving a pair $(A,M)$, where $A$ is an (associative) algebra object of $\calC^{\otimes}$ and $M$ is a left $A$-module (Proposition \ref{hunters}).

Let $\calO^{\otimes}$ be a coherent $\infty$-operad. In \S \ref{delconj}, we will apply the theory of left modules to develop a general theory of {\em centralizers} of maps $f: A \rightarrow B$ between $\calO^{\otimes}$-algebra objects of a symmetric monoidal $\infty$-category $\calC^{\otimes}$. Taking $\calO^{\otimes}$ to be a little cubes operad $\OpE{k}$ (which is coherent by Theorem \ref{cubecoh}), we obtain a proof of the
{\it generalized Deligne conjecture} (the original statement of Deligne's conjecture corresponds to the case $k=1$). As another application, we describe the construction of the {\it Koszul dual} of an 
(augmented) $\OpE{k}$-algebra (Example \ref{kd}). 

If $M$ is any associative monoid, then its group $M^{\times}$ of invertible elements acts on $M$ by conjugation. In \S \ref{adjrep}, we will describe a higher-categorical analogue of this adjoint action, which applies to $\OpE{k}$-algebra objects of an arbitrary symmetric monoidal $\infty$-category. Our main result, Theorem \ref{machus}, will be applied in \S \ref{cotan} to obtain a convenient description of the cotangent complex of an $\OpE{k}$-algebra $A$ (Theorem \ref{curtis}). 

\subsection{Algebras and their Module Categories}\label{conl}

Let $k$ be a commutative ring, and let $A$ be an associative $k$-algebra. If $M$
is a right $A$-module and $V$ is an $k$-module, then the tensor product $V \otimes_{k} M$ 
carries a right action of $A$, given by the formula $(v \otimes m)a = v \otimes ma$.
Via this construction, we can view the category $\Mod_{A}^{R}$ of right $A$-modules as a {\em left} module over the monoidal category $\Mod_{k}$ of $k$-modules. Our goal in this section is to generalize the above situation, replacing the category $\Mod_{k}$ of $k$-modules by an arbitrary monoidal $\infty$-category. 

Our goal in this section is twofold:

\begin{itemize}
\item[$(i)$] We will introduce an $\infty$-category $\CatAlg$ whose objects are pairs
$(\calC^{\otimes}, A)$, where $\calC^{\otimes}$ is a monoidal $\infty$-category
and $A$ is an algebra object of $\calC^{\otimes}$. Roughly speaking, a morphism from
$(\calC^{\otimes}, A)$ to $(\calD^{\otimes}, B)$ in $\CatAlg$ consists of a monoidal
functor $F: \calC^{\otimes} \rightarrow \calD^{\otimes}$ together with a map of algebras
$F(A) \rightarrow B$. For a more precise description, see Definition \ref{spork} below. 

\item[$(ii)$] Let $\CatMod$ denote the $\infty$-category of pairs 
$(\calC^{\otimes}, \calM)$, where
$\calC^{\otimes}$ is a monoidal $\infty$-category and $\calM$ is an $\infty$-category
left-tensored over $\calC^{\otimes}$ (see Definition \monoidref{defcatmod}).
We will define functor $\Theta$ from a subcategory of $\CatAlg$ to $\CatMod$.
Informally, the functor $\Theta$ associates to every pair $(\calC^{\otimes}, A)$ the
$\infty$-category $\Mod_{A}^{R}(\calC)$ of right $A$-module objects of $\calC^{\otimes}$.
\end{itemize}

The relevant constructions are straightforward but somewhat tedious. The reader who does
not wish to become burdened by technicalities is invited to proceed directly to \S \ref{probix2}, where we will undertake a deeper study of the functor $\Theta$. 

We begin by introducing some terminology.

\begin{definition}\label{paswa}
Let $S$ be a simplicial set. A {\it coCartesian $S$-family of monoidal $\infty$-categories} is
a coCartesian fibration $q: \calC^{\otimes} \rightarrow \Nerve( \cDelta)^{op} \times S$
with the following property: for every vertex $s \in S$, the induced map of fibers
$\calC^{\otimes}_{s} = \calC^{\otimes} \times_{S} \{s\} \rightarrow \Nerve( \cDelta)^{op}$ is a monoidal
$\infty$-category. In this case, we will say that {\it $q$ exhibits $\calC^{\otimes}$ as a coCartesian
$S$-family of monoidal $\infty$-categories}.
\end{definition}

\begin{notation}
If $q: \calC^{\otimes} \rightarrow \Nerve( \Delta)^{op} \times S$ is a coCartesian $S$-family of monoidal $\infty$-categories, we let $\calC$ denote the fiber product $\calC^{\otimes} \times_{ \Nerve( \cDelta)^{op} } \{ [1] \}$, so that $q$ induces a coCartesian fibration
$\calC \rightarrow S$.
\end{notation}

\begin{example}
Let $\Cat_{\infty}^{\Mon}$ denote the $\infty$-category of monoidal $\infty$-categories. The $\infty$-category $\Cat_{\infty}^{\Mon}$ is equivalent to the $\infty$-category $\Mon(\Cat_{\infty})$ of monoid objects of $\Cat_{\infty}$ (see Remark \monoidref{otherlander}). In particular, there is a canonical map
$\Nerve( \cDelta^{op} ) \times \Cat_{\infty}^{\Mon} \rightarrow \Cat_{\infty}$, which classifies
a coCartesian fibration $q: \widetilde{ \Cat}_{\infty}^{\Mon} \rightarrow
\Nerve( \cDelta^{op}) \times \Cat_{\infty}^{\Mon}$. The coCartesian fibration $q$ exhibits
$\widetilde{\Cat}_{\infty}^{\Mon}$ as a coCartesian $\Cat_{\infty}^{\Mon}$-family of monoidal $\infty$-categories. Moreover, this family of monoidal $\infty$-categories is universal in the following sense:
for every simplicial set $S$, the construction $( \phi: S \rightarrow \Cat_{\infty}^{\Mon})
\mapsto S \times_{ \Cat_{\infty}^{\Mon} } \widetilde{ \Cat}_{\infty}^{\Mon}$ establishes a bijection
between the collection of equivalence classes of diagrams $S \rightarrow \Cat_{\infty}^{\Mon}$ and the collection of equivalence classes of coCartesian $S$-families of monoidal $\infty$-categories
$\calC^{\otimes} \rightarrow \Nerve( \cDelta)^{op} \times S$ (with essentially small fibers).
\end{example}

\begin{definition}
Let $K$ and $S$ be simplicial sets. We will say that a coCartesian $S$-family of
$\infty$-categories $q: \calC^{\otimes} \rightarrow \Nerve( \cDelta)^{op} \times S$ is
{\it compatible with $K$-indexed colimits} if the following conditions are satisfied:
\begin{itemize}
\item[$(i)$] For each vertex $s \in S$, the fiber $\calC_{s}$ admits $K$-indexed colimits.
\item[$(ii)$] For each vertex $s \in S$, the tensor product functor
$\calC_{s} \times \calC_{s} \rightarrow \calC_{s}$ preserves $K$-indexed colimits separately in each variable.
\item[$(iii)$] For every edge $s \rightarrow t$ in $S$, the induced functor
$\calC_{s} \rightarrow \calC_{t}$ preserves $K$-indexed colimits. 
\end{itemize}
If $\calK$ is a collection of simplicial sets, we will say that $q$ is {\it compatible with $\calK$-indexed colimits} if it is compatible with $K$-indexed colimits for each $K \in \calK$.
\end{definition}

\begin{notation}\label{kost}
Let $\calK$ be a collection of simplicial sets. We let $\Cat_{\infty}^{\Mon}(\calK)$ denote the
subcategory of $\Cat_{\infty}^{\Mon}$ whose objects are monoidal $\infty$-categories
$\calC^{\otimes}$ which are compatible with $\calK$-indexed colimits and whose morphisms
are monoidal functors $F: \calC^{\otimes} \rightarrow \calD^{\otimes}$ such that the underlying functor
$\calC \rightarrow \calD$ preserves $\calK$-indexed colimits. Let
$\widetilde{ \Cat}_{\infty}^{\Mon}(\calK)$ denote the fiber product
$\widetilde{\Cat}_{\infty}^{\Mon} \times_{ \Cat_{\infty}^{\Mon} } \Cat_{\infty}^{\Mon}(\calK)$.
The evident map $q: \widetilde{ \Cat}_{\infty}^{\Mon}(\calK) \rightarrow
\Nerve( \cDelta)^{op} \times \Cat_{\infty}^{\Mon}(\calK)$ exhibits $\widetilde{ \Cat}_{\infty}^{\Mon}(\calK)$ as a coCartesian $\Cat_{\infty}^{\Mon}(\calK)$-family of monoidal $\infty$-categories which is compatible with $\calK$-indexed colimits. Moreover, it is universal with respect to this property: for every simplicial set $S$, pullback along $q$ induces a bijection from equivalence classes of diagrams
$S \rightarrow \Cat_{\infty}^{\Mon}(\calK)$ and equivalence classes of
coCartesian $S$-families of monoidal $\infty$-categories which are compatible with $\calK$-indexed colimits.
\end{notation}

\begin{definition}\label{prespork}
We let $\CatAlg$ denote the full subcategory of the fiber product
$$\Cat_{\infty}^{\Mon} \times_{ \Fun_{ \Nerve(\cDelta)^{op}} ( \Nerve( \cDelta)^{op}, \Nerve( \cDelta)^{op} \times \Cat_{\infty}^{\Mon})} 
\Fun_{ \Nerve( \cDelta)^{op}}( \Nerve( \cDelta)^{op}, \widetilde{ \Cat}^{\Mon}_{\infty})$$
spanned by those pairs $(\calC^{\otimes}, A)$, where $\calC^{\otimes} \in \Cat_{\infty}^{\Mon}$ is a monoidal $\infty$-category and $A$ is an algebra object of the monoidal $\infty$-category 
$\widetilde{ \Cat}^{\Mon}_{\infty} \times_{ \Cat_{\infty}^{\Mon} } \{ \calC^{\otimes} \} \simeq \calC^{\otimes}$. If $\calK$ is a collection of simplicial sets, we let
$\CatAlg(\calK)$ denote the fiber product $\CatAlg \times_{ \Cat_{\infty}^{\Mon}} \Cat_{\infty}^{\Mon}(\calK)$. 
\end{definition}

\begin{remark}\label{spork}
The $\infty$-category $\CatAlg(\calK)$ is characterized up to equivalence by the following universal property: for any simplicial set $S$, there is a bijection between equivalence classes of diagrams
$S \rightarrow \CatAlg(\calK)$ and equivalence classes of diagrams
$$ \xymatrix{ & \calC^{\otimes} \ar[d]^{q} \\
\Nerve( \cDelta)^{op} \times S \ar[ur]^{A} \ar[r]^{\id} & \Nerve( \cDelta)^{op} \times S, }$$
where $q$ exhibits $\calC^{\otimes}$ as a coCartesian $S$-family of monoidal $\infty$-categories
whose fibers are essentially small, $q$ is compatible with $\calK$-indexed colimits, and
$A$ is an $S$-family of algebra objects of $\calC^{\otimes}$.
\end{remark}

\begin{definition}\label{splot}
If $q: \calC^{\otimes} \rightarrow \Nerve( \cDelta)^{op} \times S$ is a coCartesian $S$-family of monoidal $\infty$-categories, then we will say that a map $p: \calM^{\otimes} \rightarrow \calC^{\otimes}$
{\it exhibits $\calM^{\otimes}$ as a coCartesian $S$-family of $\infty$-categories left-tensored over $\calC^{\otimes}$} if the following conditions are satisfied:
\begin{itemize}
\item[$(i)$] The composition $q \circ p$ is a coCartesian fibration.
\item[$(ii)$] The map $p$ is a categorical fibration which carries $(q \circ p)$-coCartesian morphisms to $q$-coCartesian morphisms.
\item[$(iii)$] For every vertex $s \in S$, the induced map of fibers $\calM^{\otimes}_{s} \rightarrow \calC^{\otimes}_{s}$ exhibits $\calM^{\otimes}_{s}$ as an $\infty$-category left-tensored over $\calC^{\otimes}_{s}$ (see Definition \monoidref{ulult}). \end{itemize}
\end{definition}

\begin{notation}
If $\calM^{\otimes} \rightarrow \calC^{\otimes} \rightarrow \Nerve( \cDelta)^{op} \times S$
is as in Definition \ref{splot}, we let $\calM$ denote the fiber product $\calM^{\otimes} \times_{ \Nerve( \cDelta)^{op} } \{ [0] \}$.
\end{notation}

\begin{remark}\label{cay1}
Assuming conditions $(i)$ and $(ii)$ of Definition \ref{splot} are satisfied, condition $(iii)$ is equivalent to the requirement that for each $n \geq 0$ and each $s \in S$, the canonical map $\calM^{\otimes}_{([n],s)} \rightarrow \calM^{\otimes}_{([0],s)} \times \calC^{\otimes}_{([n],s)}$ is an equivalence of $\infty$-categories.
\end{remark}

\begin{remark}\label{cay2}
In the situation of Definition \ref{splot}, conditions $(i)$, $(ii)$, and $(iii)$ guarantee that the map
$\calM^{\otimes} \rightarrow \calC^{\otimes}$ is a locally coCartesian fibration (Proposition \toposref{fibertest}). Conversely, suppose that $p: \calC^{\otimes} \rightarrow \Nerve( \cDelta)^{op} \times S$ is a coCartesian $S$-family of monoidal $\infty$-categories, and let $q: \calM^{\otimes} \rightarrow \calC^{\otimes}$ be a locally coCartesian categorical fibration. Then conditions $(i)$ and
$(ii)$ of Definition \ref{splot} are equivalent to the following:
\begin{itemize}
\item[$(\ast)$] Let $e: M \rightarrow N$ be a locally $q$-coCartesian edge of $\calM^{\otimes}$ such that
$q(e)$ is $p$-coCartesian. Then $e$ is $q$-coCartesian.
\end{itemize}
To see this, let us first suppose that $q$ satisfies $(i)$ and $(ii)$, and let $e: M \rightarrow N$ be as in $(\ast)$. Then $(i)$ implies that there exists
a $(p \circ q)$-coCartesian edge $e': M \rightarrow M'$ lifting $(p \circ q)(e)$. Condition
$(ii)$ ensures that $q(e')$ is $p$-coCartesian. We may therefore assume without loss of generality
that $q(e) = q(e')$. Proposition \toposref{protohermes} guarantees that $e'$ is 
$q$-coCartesian, hence locally $q$-coCartesian and therefore homotopic to $e$; this proves
that $e$ is $q$-coCartesian. 

Conversely, suppose that $(\ast)$ is satisfied. Let $M \in \calM^{\otimes}$ and let
$\alpha: (p \circ q)(M) \rightarrow X$ be an edge of $\Nerve(\cDelta)^{op} \times S$.
We will prove $(i)$ and $(ii)$ by showing that $\alpha$ can be lifted to a $(p \circ q)$-coCartesian
edge $e: M \rightarrow N$ such that $q(e)$ is $p$-coCartesian. Using our assumption that
$p$ is a coCartesian fibration, we can choose a $p$-coCartesian edge $\overline{\alpha}: q(M) \rightarrow \overline{X}$. Since $q$ is a locally coCartesian fibration, we can choose a locally
$q$-coCartesian edge $e: M \rightarrow N$ lifting $\overline{\alpha}$. Condition $(\ast)$ guarantees that $e$ is $q$-coCartesian as required.

Using Lemma \toposref{charloccart}, we can reformulate condition $(\ast)$ as follows:
\begin{itemize}
\item[$(\ast')$] Suppose we are given a $2$-simplex
$$ \xymatrix{ & M' \ar[dr]^{e''} & \\
M \ar[ur]^{e'} \ar[rr]^{e} & & M'' }$$
in $\calM^{\otimes}$, where $e'$ and $e''$ are locally $q$-coCartesian
and $q(e)$ is $p$-coCartesian. Then $e$ is locally $q$-coCartesian.
\end{itemize}
\end{remark}

\begin{notation}\label{suspik}
Let $\calC^{\otimes} \rightarrow \Nerve( \cDelta)^{op} \times S$ be a coCartesian $S$-family of $\infty$-operads. We define a simplicial set $\widetilde{\Alg}(\calC)$ equipped with a forgetful map
$\widetilde{\Alg}(\calC) \rightarrow S$ so that the following universal property is satisfied:
for every map of simplicial sets $K \rightarrow S$, there is a canonical bijection
$$ \Hom_{S}( K, \widetilde{\Alg}(\calC)) \simeq \Hom_{ \Nerve( \cDelta)^{op} \times S}( \Nerve(\cDelta)^{op} \times K, \calC^{\otimes}).$$
We let $\Alg(\calC)$ denote the full simplicial subset of $\widetilde{\Alg}(\calC)$ spanned by those
vertices which correspond to algebra objects in the monoidal $\infty$-category $\calC^{\otimes}_{s}$, for some vertex $s \in S$.

Suppose we are given a map $\calM^{\otimes} \rightarrow \calC^{\otimes}$ satisfying condition
$(iii)$ of Definition \ref{splot}. We let $\widetilde{\Mod}^{L}(\calM)$ denote a simplicial set
with a map $\widetilde{\Mod}^{L}(\calM)$ having the following universal property: for
every map of simplicial sets $K \rightarrow S$, there is a canonical bijection
$$ \Hom_{S}( K, \widetilde{\Mod}^{L}(\calM)) \simeq \Hom_{ \Nerve( \cDelta)^{op} \times S}( \Nerve(\cDelta)^{op} \times K, \calM^{\otimes}).$$
We let $\Mod^{L}(\calM)$ denote the full simplicial subset of $\widetilde{\Mod}^{L}(\calM)$ whose vertices are left module objects of $\calM^{\otimes}_{s}$, for some vertex $s \in S$.
\end{notation}

\begin{remark}
In the special case where $S = \Delta^0$, the terminology of Notation \ref{suspik} agrees with that of
Definitions \monoidref{suskin} and \monoidref{defmond}.
\end{remark}

The following result is an easy consequence of Proposition \toposref{doog}:

\begin{lemma}\label{precane}
Let $q: \calC^{\otimes} \rightarrow \Nerve( \cDelta)^{op} \times S$ be a coCartesian $S$-family of monoidal $\infty$-categories, let $p: \calM^{\otimes} \rightarrow \calC^{\otimes}$ be a coCartesian $S$-family of $\infty$-categories left-tensored over $\calC^{\otimes}$. Then:
\begin{itemize}
\item[$(1)$] The map $q': \Alg(\calC) \rightarrow S$ is a coCartesian fibration of simplicial sets.
\item[$(2)$] A morphism $A \rightarrow A'$ in $\Alg(\calC)$ is $q'$-coCartesian if and only if
the underying map $A([1]) \rightarrow A'([1])$ is a $q$-coCartesian morphism in $\calC \subseteq \calC^{\otimes}$.
\item[$(3)$] The map $r: \Mod^{L}(\calM) \rightarrow S$ is a coCartesian fibration of simplicial sets.
\item[$(4)$] A morphism $M \rightarrow M'$ in $\Mod^{L}(\calM)$ is $r$-coCartesian if and only if
its image in $\Alg(\calC)$ is $q'$-coCartesian, and the induced map $M([0]) \rightarrow M'([0])$ is
a $(q \circ p)$-coCartesian morphism in $\calM \subseteq \calM^{\otimes}$. 
\end{itemize}
\end{lemma}

\begin{definition}
Let $\calK$ be a collection of simplicial sets, and let $\calC^{\otimes} \rightarrow \Nerve(\cDelta)^{op} \times S$ be a coCartesian $S$-family of monoidal $\infty$-categories which is compatible with $\calK$-indexed colimits. We will say that a coCartesian $S$-family $\calM^{\otimes}$ of $\infty$-categories left-tensored over $\calC^{\otimes}$ is {\it compatible with $\calK$-indexed colimits} if the following conditions are satisfied:
\begin{itemize}
\item[$(i)$] For every vertex $s \in S$ and each $K \in \calK$, the $\infty$-category
$\calM_{s}$ admits $K$-indexed colimits.
\item[$(ii)$] For every vertex $s \in S$ and each $K \in \calK$, the action map
$\calC_{s} \times \calM_{s} \rightarrow \calM_{s}$ preserves $K$-indexed colimits separately in each variable.
\item[$(iii)$] For every edge $s \rightarrow t$ in $S$ and each $K \in \calK$, the induced functor
$\calM_{s} \rightarrow \calM_{t}$ preserves $K$-indexed colimits.
\end{itemize}
\end{definition}

\begin{lemma}\label{precanus}
Let $p: \calC^{\otimes} \rightarrow \Nerve(\cDelta)^{op}$ be a coCartesian $S$-family of monoidal $\infty$-categories and let $q: \calM^{\otimes} \rightarrow \calC^{\otimes}$ 
be a coCartesian $S$-family of $\infty$-categories
left-tensored over $\calC^{\otimes}$. Assume that both $p$ and $q$ are compatible with
$\Nerve(\cDelta)^{op}$-indexed colimits. Then:
\begin{itemize}
\item[$(1)$] The forgetful functor $r: \Mod^{L}(\calM) \rightarrow \Alg(\calC)$ is a coCartesian fibration of simplicial sets. 
\item[$(2)$] Let $f: M \rightarrow N$ be an in $\Mod^L(\calM)$, lying over an edge
$f_0: A \rightarrow B$ in $\Alg(\calC)$, which in turn lies over an edge $\alpha: s \rightarrow t$
in $S$. Then $f$ is $r$-coCartesian if and only if it induces an equivalence
$B \otimes_{ \alpha_{!} A} M \rightarrow N$ in the $\infty$-category $\Mod^L(\calM)_{t}$, where
$\alpha_{!}: \Alg(\calC)_{s} \rightarrow \Alg(\calC)_{t}$ denotes the functor induced by $\alpha$.
\end{itemize}
\end{lemma}

\begin{proof}
Choose a vertex $M \in \Mod^L(\calM)$
lying over $A \in \Alg(\calC)$, and let $f_0: A \rightarrow B$ be an edge of $\Alg(\calC)$ lying
over an edge $\alpha: s \rightarrow t$ in $S$. To prove $(1)$, we will show that $f_0$ can be lifted to
an $r$-coCartesian morphism of $\Mod^L(\calM)$; assertion $(2)$ will be a consequence of our construction. Let $r': \Alg(\calC) \rightarrow S$ denote the canonical projection. Using Lemma \ref{precane}, we can lift $\alpha$ to an $(r' \circ r)$-coCartesian morphism $f': M \rightarrow M'$ in $\Mod^L(\calM)$; let $f'_0: A \rightarrow A'$ denote the image of $f'$ in $\Alg(\calC)$. Lemma
\ref{precane} guarantees that $f'_0$ is $r'$-coCartesian, so we can identify
$A'$ with $\alpha_{!} A$; moreover, there exists a $2$-simplex $\sigma$ of $\Alg(\calC)$ 
corresponding to a diagram
$$ \xymatrix{ & A' \ar[dr]^{f''_0} & \\
A \ar[rr]^{f_0} \ar[ur]^{f'_0} & & B. }$$
We will prove that $f''_0$ can be lifted to an $r$-coCartesian morphism $f''$ of $\Mod^L(\calM)$.
Using the fact that $r$ is an inner fibration, it will follow that there is a composition
$f= f'' \circ f'$ lifting $f_0$, which is also $r$-coCartesian by virtue of Proposition \toposref{protohermes}.
We may therefore replace $f_0$ by $f''_0$ and thereby reduce to the case where $s=t$ and the
edge $\alpha$ is degenerate.

Corollary \monoidref{tara} shows that $f_0$ can be lifted to a locally $r$-coCartesian
morphism $f: M \rightarrow B \otimes_{A} M$ in $\Mod^{L}(\calM)_{s}$. Since the projection
$r_s: \Mod^L(\calM)_{s} \rightarrow \Alg(\calM)_{s}$ is a Cartesian fibration (Corollary \monoidref{thetacart}), we deduce that $f$ is $r_{s}$-coCartesian (Corollary \toposref{grutt1}).
To prove that $f$ is $r$-coCartesian, it will suffice to show that for every edge
$\beta: s \rightarrow t$ in $S$, the image $\beta_{!}(f)$ is an $r_{t}$-coCartesian morphism
of the fiber $\Mod^L(\calM)_{t}$. Using the characterization of $r_{t}$-coCartesian morphisms
supplied by Corollary \monoidref{tara}, we see that this is equivalent to the requirement that the canonical map $\alpha_{!}( B) \otimes_{ \alpha_{!} A} \alpha_{!} M \rightarrow 
\alpha_{!}( B \otimes_{A} M)$ is an equivalence in $\Mod^L(\calM)_{t}$. This is clear, since
the functor $\alpha_{!}$ preserves tensor products and geometric realizations of simplicial objects.
\end{proof}

\begin{definition}\label{spull}
Let $p: \calC^{\otimes} \rightarrow \Nerve( \cDelta)^{op} \times S$ be a coCartesian $S$-family of monoidal $\infty$-categories. An {\it $S$-family of algebra objects of $\calC^{\otimes}$} is
a section of the projection map $\Alg(\calC) \rightarrow S$. 

If $q: \calM^{\otimes} \rightarrow \calC^{\otimes}$ is a map satisfying condition $(iii)$ of Definition \ref{splot} and $A$ is an $S$-family of algebra objects of $\calC^{\otimes}$,
then we let $\Mod^{L}_A(\calM)$ denote the fiber product $\Mod(\calM) \times_{ \Alg(\calC) } S$.
\end{definition}

\begin{remark}\label{swiff}
In the situation of Definition \ref{spull}, if $p$ and $q$ are compatible with $\Nerve(\cDelta)^{op}$-indexed colimits, then Lemma \ref{precanus} implies that the projection map $\Mod^{L}_{A}(\calM) \rightarrow S$ is a coCartesian fibration of simplicial sets.
\end{remark}

\begin{variant}
Let $p: \calC^{\otimes} \rightarrow \Nerve(\cDelta)^{op} \times S$ be a coCartesian $S$-family of monoidal $\infty$-categories. We will say that an inner fibration $q: \calM^{\otimes} \rightarrow \calC^{\otimes}$
is a {\it locally coCartesian $S$-family of $\infty$-categories left-tensored over $\calC^{\otimes}$} if, for every edge
$\Delta^1 \rightarrow S$, the induced map $q_{\Delta^1}: \calM^{\otimes} \times_{ S} \Delta^1 \rightarrow
\calC^{\otimes} \times_{S} \Delta^1$ is a coCartesian $\Delta^1$-family of $\infty$-categories
left-tensored over $\calC^{\otimes} \times_{S} \Delta^1$. If $K$ is a simplicial set, we will say that
$q$ is {\it compatible with $K$-indexed colimits} if each $q_{\Delta^1}$ is compatible with $K$-indexed colimits. If $p$ and $q$ are compatible with $\Nerve(\cDelta)^{op}$-indexed colimits and
$A$ is an $S$-family of algebra objects of $\calC^{\otimes}$, then Remark \ref{swiff} implies
that the map $\Mod^{L}_{A}(\calM) \rightarrow S$ is a locally coCartesian fibration of simplicial sets.
\end{variant}

\begin{variant}\label{cluin}
In the situation of Definition \ref{splot}, 
there is an evident dual notion of a {\it locally coCartesian $S$-family of $\infty$-categories $ \calM^{\otimes} \rightarrow \calC^{\otimes}$ right-tensored over $\calC^{\otimes}$}.
Given an $S$-family of algebra objects $A$ of $\calC^{\otimes}$, we can then define
a locally coCartesian fibration $\Mod^{R}_{A}(\calM) \rightarrow S$ (provided that $\calC^{\otimes}$ and $\calM^{\otimes}$ are compatible with $\Nerve(\cDelta)^{op}$-indexed colimits), whose fiber over a vertex
$s \in S$ is the $\infty$-category of right $A$-module objects of the fiber $\calM_{s}$.
\end{variant}

Let $\CatMod$ be the $\infty$-category of Definition \monoidref{defcatmod}, whose objects
are diagrams $\calM^{\otimes} \rightarrow \calC^{\otimes} \rightarrow \Nerve( \cDelta)^{op}$ which
exhibit $\calC^{\otimes}$ as a monoidal $\infty$-category and $\calM^{\otimes}$ as an $\infty$-category
which is left-tensored over $\calC^{\otimes}$. We will informally describe the objects of
$\CatMod$ as pairs $(\calC, \calM)$, where $\calC$ is an $\infty$-category equipped with a monoidal structure and $\calM$ is an $\infty$-category equipped with a left action of $\calC$.
If $\calK$ is a collection of simplicial sets,
we let $\CatMod(\calK)$ denote the subcategory of $\CatMod$ whose objects are diagrams
where $\calC$ and $\calM$ admit $\calK$-indexed colimits and the tensor product functors
$$ \calC \times \calC \rightarrow \calC \quad \quad \calC \times \calM \rightarrow \calM$$
preserve $\calK$-indexed colimits separately in each variable, and whose morphisms 
are maps $(\calC, \calM) \rightarrow (\calC', \calM')$ such that the underlying functors
$\calC \rightarrow \calC'$, $\calM \rightarrow \calM'$ preserve $\calK$-indexed colimits.

\begin{remark}\label{stla}
The $\infty$-category $\CatMod(\calK)$ is characterized by the following universal property: for every
simplicial set $S$, there is a canonical bijection between equivalence classes of maps
$S \rightarrow \CatMod(\calK)$ and equivalence classes of diagrams
$\calM^{\otimes} \rightarrow \calC^{\otimes} \rightarrow \Nerve( \cDelta)^{op} \times S$
which exhibit $\calC^{\otimes}$ as a coCartesian $S$-family of monoidal $\infty$-categories compatible with $\calK$-indexed colimits, and $\calM^{\otimes}$ as a coCartesian $S$-family of $\infty$-categories left-tensored over $\calC^{\otimes}$.
\end{remark}

We now sketch the construction of the functor $\Theta$.

\begin{construction}\label{cloops}
Let $c: \Nerve( \cDelta)^{op} \times \Nerve( \cDelta)^{op} \rightarrow \Nerve( \cDelta)^{op}$ be
the concatenation functor, given on objects by the formula
$c( [m], [n]) = [m] \star [n] \simeq [m+n + 1]$.
Let $\pi_0, \pi_1: \Nerve( \cDelta)^{op} \times \Nerve( \cDelta)^{op} \rightarrow \Nerve(\cDelta)^{op}$ denote the projection functors. The inclusions of linearly ordered sets $[m] \hookrightarrow [m] \star [n] \hookleftarrow [n]$
induce natural transformations
$\pi_0 \stackrel{\alpha}{\leftarrow} c \stackrel{\beta}{\rightarrow} \pi_1$, which we can identify
with a map $\gamma: \Lambda^2_0 \times \Nerve( \cDelta)^{op} \times \Nerve(\cDelta)^{op} \rightarrow \Nerve( \cDelta)^{op}$.

Let $q: \calC^{\otimes} \rightarrow \Nerve( \cDelta)^{op} \times S$ be a coCartesian $S$-family of
monoidal $\infty$-categories. We define an $\infty$-category $\overline{B}( \calC^{\otimes})$ 
equipped with a map $p: \overline{B} \rightarrow \Nerve( \cDelta)^{op} \times \Nerve( \cDelta)^{op} \times S$ so that the following universal property is satisfied: for every map of simplicial sets
$K \rightarrow \Nerve( \cDelta)^{op} \times \Nerve( \cDelta)^{op} \times S$, there is a canonical bijection
of $\Hom_{/\Nerve( \cDelta)^{op} \times \Nerve( \cDelta)^{op} \times S}( K, \overline{B}( \calC^{\otimes}))$ with the collection of all commutative diagrams
$$ \xymatrix{ \Lambda^2_0 \times K \ar[d] \ar[r] & \calC^{\otimes} \ar[d] \\
\Lambda^2_0 \times \Nerve( \cDelta)^{op} \times \Nerve( \cDelta)^{op} \times S \ar[r]^-{\gamma \times \id_{S}} & \Nerve(\cDelta)^{op} \times S.}$$

We can identify an object of the the fiber of the map $p$ over a triple $( [m], [n], s)$ with
a pair of morphisms $C \stackrel{f}{\leftarrow} M \stackrel{g}{\rightarrow} C'$ in the
$\infty$-category $\calC^{\otimes}_{s}$, where $C \in (\calC^{\otimes}_{s})_{[m]}$, $M \in (\calC^{\otimes}_{s})_{[m+n+1]}$, $D \in (\calC^{\otimes}_{s})_{[n]}$, and the maps
$f$ and $g$ cover the inclusions $[m] \hookrightarrow [m] \star [n] \hookleftarrow [n]$ in
$\cDelta$. Let $B( \calC^{\otimes})$ denote the full subcategory of $\overline{B}(\calC^{\otimes})$ spanned by those objects for which the morphisms $f$ and $g$ are $q$-coCartesian.
\end{construction}

\begin{remark}
Evaluation at the vertices $\{1\}, \{2\} \subseteq \Lambda^2_0$ induces a canonical map
$B( \calC^{\otimes}) \rightarrow \calC^{\otimes} \times_{S} \calC^{\otimes}$.
\end{remark}

\begin{lemma}\label{cuppose}
Let $q: \calC^{\otimes} \rightarrow \Nerve(\cDelta)^{op} \times S$ be a coCartesian $S$-family of monoidal $\infty$-categories. Then:
\begin{itemize}
\item[$(1)$] Let $\calD^{\otimes}$ denote the fiber product $\calC^{\otimes} \times_{S} \calC^{\otimes}$.
The projection $$\pi: \calD^{\otimes}
\rightarrow \calC^{\otimes} \times_{S} ( \Nerve(\cDelta)^{op} \times S) \simeq \Nerve(\cDelta)^{op} \times \calC^{\otimes}$$
exhibits $\calD^{\otimes}$ as a coCartesian $\calC^{\otimes}$-family
of monoidal $\infty$-categories.
\item[$(2)$] The map $p: B( \calC^{\otimes}) \rightarrow \calD^{\otimes}$
exhibits $B( \calC^{\otimes})$ as a locally coCartesian $\calC^{\otimes}$-family of $\infty$-categories
left-tensored over $\calD^{\otimes}$.
\item[$(3)$] 
If $\sigma: \Delta^2 \rightarrow \calC^{\otimes}$ is a $2$-simplex such that
$\sigma | \Delta^{ \{0,1\} }$ is a $q$-coCartesian edge of $\calC^{\otimes}$, then
the restriction $p_{ \sigma}: B( \calC^{\otimes}) \times_{ \calC^{\otimes} } \Delta^2 \rightarrow
\Nerve( \cDelta)^{op} \times \Delta^2$ is a coCartesian $\Delta^2$-family of $\infty$-categories
right-tensored over $\calD^{\otimes} \times_{ \calC^{\otimes} } \Delta^2$.
\end{itemize}
\end{lemma}

\begin{proof}
Assertion $(1)$ is obvious. In view of Remarks \ref{cay1} and \ref{cay2}, assertions $(2)$ and
$(3)$ are consequences of the following:
\begin{itemize}
\item[$(a)$] The map $p$ is a locally coCartesian fibration.
\item[$(b)$] Suppose we are given a $2$-simplex
$$ \xymatrix{ & M' \ar[dr]^{e''} & \\
M \ar[ur]^{e'} \ar[rr]^{e} & & M'' }$$
in $B(\calC^{\otimes})$ where $e'$ and $e''$ are locally $p$-coCartesian,
the edge $p(e')$ is $\pi$-coCartesian, and the edge $\pi( p(e') )$ is
$q$-coCartesian. Then $e$ is locally $p$-coCartesian.
\item[$(c)$] For every object $C \in \calC^{\otimes}$, the induced map
$$ B(\calC^{\otimes})_{C, [n]} \rightarrow \calD^{\otimes}_{C, [n]} \times
B( \calC^{\otimes})_{C, [0]}$$
is an equivalence of $\infty$-categories.
\end{itemize}
We first prove $(a)$. Proposition \toposref{doog} implies that
the maps $r: \calD^{\otimes} \rightarrow \Nerve( \cDelta)^{op} \times S \times \Nerve(\cDelta)^{op}$
and $r': B( \calC^{\otimes}) \rightarrow \Nerve( \cDelta)^{op} \times S \times \Nerve( \cDelta)^{op}$
are coCartesian fibrations, and that $p$ carries $r'$-coCartesian edges to $r$-coCartesian edges. It therefore suffices to prove that for every object $X = ([m], s, [n]) \in \Nerve(\cDelta)^{op} \times S \times \Nerve(\cDelta)^{op}$, the induced map of fibers $p_X: B( \calC^{\otimes})_{X} \rightarrow
\calD^{\otimes}_{X}$ is a locally coCartesian fibration (Proposition \toposref{fibertest}). 
We now observe that $p_X$ is equivalent to the projection map
$\calC^{m+n+1}_{s} \rightarrow \calC_{s}^{m+n}$ which omits the ``middle'' factor.
This proves $(a)$. Assertion $(b)$ follows from the description of the class of locally $p$-coCartesian morphisms supplied by Proposition \toposref{fibertest}, and assertion $(c)$ is easy.
\end{proof}

\begin{remark}\label{slopper}
Suppose that $\calC^{\otimes} \rightarrow \Nerve(\cDelta)^{op} \times S$ is a coCartesian
$S$-family of monoidal $\infty$-categories which is compatible with $\calK$-indexed colimits, for some collection of simplicial sets $\calK$. Then the locally coCartesian $\calC^{\otimes}$-family of $\infty$-categories $B( \calC^{\otimes}) \rightarrow \calD^{\otimes}$ appearing in Lemma \ref{cuppose} will again be compatible with $\calK$-indexed colimits.
\end{remark}

\begin{construction}\label{spays}
Let $q: \calC^{\otimes} \rightarrow \Nerve( \cDelta)^{op} \times S$ be a coCartesian $S$-family of
monoidal $\infty$-categories which is compatible with $\Nerve(\cDelta)^{op}$-indexed colimits, and let $A$ be an $S$-family of algebra objects of $\calC^{\otimes}$. Then $A$ determines a section
of the projection $\pi: \calD^{\otimes} \rightarrow \calC^{\otimes}$, which we will also denote by
$A'$. We let $\Mod_{A}^{R}(\calC)^{\otimes}$ denote the $\infty$-category
$\Mod_{A'}^{R}( B(\calC^{\otimes}) )$ described in Variant \ref{cluin}.
\end{construction}

\begin{remark}\label{spluce}
In the situation of Construction \ref{spays}, suppose that the simplicial set $S$ consists of a single vertex,
so that $\calC^{\otimes}$ is a monoidal $\infty$-category.
Fix an object $C \in \calC^{\otimes}_{[0]}$. Then the fiber product
$\Mod^{R}_{A}( \calC)^{\otimes} \times_{ \calC^{\otimes} } \{C\}$ can be identified
with the $\infty$-category $\Mod^{R}_{A}(\calC)$ of right $A$-module objects of $\calC$ (Definition \monoidref{defmond}). More generally, each fiber $\Mod^{R}_{A}(\calC)^{\otimes}_{[n]}$ is
equivalent to a product $\calC^{\otimes}_{[n]} \times \Mod^{R}_{A}(\calC)$, with projection onto the first factor induced by the forgetful functor $\Mod^{R}_{A}(\calC)^{\otimes} \rightarrow \calC^{\otimes}$.
\end{remark}

Using Lemma \ref{cuppose} and Remark \ref{slopper}, we see that the canonical projection
$p: \Mod_{A}^{R}(\calC)^{\otimes} \rightarrow \calC^{\otimes}$ is a locally coCartesian categorical fibration
of simplicial sets satisfying condition $(\ast')$ of Remark \ref{cay2}. Remarks \ref{spluce} and \ref{cay1} imply that $p$ satisfies condition $(iii)$ of Definition \ref{splot}. We can summarize our analysis as follows:

\begin{proposition}\label{kanus}
Let $\calK$ be a collection of simplicial sets which includes $\Nerve( \cDelta)^{op}$, 
let $q: \calC^{\otimes} \rightarrow \Nerve( \cDelta)^{op} \times S$ be a coCartesian $S$-family of $\infty$-categories which is compatible with $\calK$-indexed colimits, and let
$A$ be an $S$-family of algebra objects of $\calC^{\otimes}$. 
Then the forgetful functor $p: \Mod^{R}_{A}( \calC)^{\otimes} \rightarrow \calC^{\otimes}$ exhibits $\Mod^{R}_{A}(\calC)^{\otimes}$ as a coCartesian $S$-family of $\infty$-categories left-tensored over $\calC^{\otimes}$ which is compatible with $\calK$-indexed colimits.
\end{proposition}

\begin{remark}
Let $\calC^{\otimes}$ be a monoidal $\infty$-category, and let $A,B \in \Alg(\calC)$. 
Proposition \ref{kanus} implies that $\Mod_{B}^{R}(\calC)^{\otimes}$ is an $\infty$-category
left-tensored over $\calC^{\otimes}$. The $\infty$-category $\Mod_{A}^{L}( \Mod_{B}^{R}(\calC))$
of left $A$-module objects of $\Mod_{B}^{R}(\calC)$ is isomorphic to the $\infty$-category of
$\Biod{A}{B}{\calC}$ of $(A,B)$-bimodules (Definition \symmetricref{bidefn}). 
\end{remark}

\begin{construction}\label{spaas}
Fix a collection of simplicial sets $\calK$ which includes $\Nerve( \cDelta)^{op}$. 
Let $\calC^{\otimes}$ denote the fiber product $\CatAlg(\calK) \times_{ \Cat_{\infty}^{\Mon} } \widetilde{ \Cat}_{\infty}^{\Mon}$, so that we have a coCartesian $\CatAlg(\calK)$-family of
monoidal $\infty$-categories $\calC^{\otimes} \rightarrow \Nerve( \FinSeg) \times \CatAlg(\calK)$.
By construction, there is a canonical $\CatAlg(\calK)$-family of algebra objects of
$\calC^{\otimes}$, which we will denote by $A$. Let $\Mod^{R}_{A}( \calC)^{\otimes}$
be the $\infty$-category of Construction \ref{spays}, so that we have a forgetful functor
$\Mod^{R}_{A}(\calC)^{\otimes} \rightarrow \calC^{\otimes}$ which exhibits 
$\Mod^{R}_{A}(\calC)^{\otimes}$ as a coCartesian $\CatAlg(\calK)$-family of $\infty$-categories which are left-tensored over $\calC^{\otimes}$. Remark \ref{stla} implies that this family is classified by a functor
$\Theta: \CatAlg(\calK) \rightarrow \CatMod(\calK)$. Note that the composite functor
$\CatAlg(\calK) \rightarrow \CatMod(\calK) \rightarrow \Cat_{\infty}^{\Mon}(\calK)$ classifies the
coCartesian $\CatAlg(\calK)$-family of monoidal $\infty$-categories $\calC^{\otimes}$,
and is therefore equivalent to the evident forgetful functor $\CatAlg(\calK) \rightarrow \Cat_{\infty}^{\Mon}(\calK)$.
Replacing $\Theta$ by an equivalent functor if necessary, we will henceforth assume that the diagram
$$ \xymatrix{ \CatAlg(\calK) \ar[rr]^{\Theta} \ar[dr] & & \CatMod(\calK) \ar[dl] \\
& \Cat_{\infty}^{\Mon}(\calK) & }$$
is commutative.
\end{construction}

\begin{remark}
More informally, we can describe the functor $\Theta: \CatAlg(\calK) \rightarrow \CatMod(\calK)$ as follows: to every object $(\calC^{\otimes}, A)$ of $\CatAlg(\calK)$ (given by a monoidal $\infty$-category
$\calC^{\otimes}$ and an algebra object $A \in \Alg(\calC)$), it associates the
$\infty$-category $\Mod_{A}^{R}(\calC)$ of right $A$-module objects of $\calC$, viewed
as an $\infty$-category left-tensored over $\calC$.
\end{remark}

\subsection{Properties of $\Mod_{A}(\calC)$}\label{probix}

Let $\calC$ be a monoidal $\infty$-category and let $A$ be an algebra object of $\calC$.
The $\infty$-category $\Mod_{A}^{R}(\calC)$ of right $A$-module objects of $\calC$ admits a left action of the $\infty$-category $\calC$: informally speaking, if $M$ is a right $A$-module and $C \in \calC$, then $C \otimes M$ admits a right $A$-module structure given by the map
$$ (C \otimes M) \otimes A \simeq C \otimes (M \otimes A) \rightarrow C \otimes M.$$
(for a complete construction, we refer to Proposition \ref{kanus}). In this section, we will prove that (under some mild hypotheses) the $\infty$-category $\Mod_{A}^{R}(\calC)$ enjoys two important features (which will be formulated more precisely below):

\begin{itemize}
\item[$(A)$] If $\calN$ is an $\infty$-category left-tensored over $\calC$, then the $\infty$-category
of $\calC$-linear functors from $\Mod_{A}^{R}(\calC)$ to $\calN$ is equivalent to the
$\infty$-category $\Mod_{A}^{L}(\calN)$ of left $A$-module objects of $\calN$ (Theorem \ref{bix2}).

\item[$(B)$] If $\calM$ is an $\infty$-category right-tensored over $\calC$, then the tensor product
$\calM \otimes_{\calC} \Mod_{A}^{R}(\calC)$ is equivalent to the $\infty$-category
$\Mod_A^{R}(\calM)$ of right $A$-module objects of $\calM$ (Theorem \ref{bix1}).
\end{itemize}

We begin by formulating assertion $(A)$ more precisely. 
\begin{definition}
Let $q: \calC^{\otimes} \rightarrow \Nerve( \cDelta)^{op}$ be a monoidal $\infty$-category, and suppose we are given maps $p: \calM^{\otimes} \rightarrow \calC^{\otimes}$
and $p': \calN^{\otimes} \rightarrow \calC^{\otimes}$ which exhibit $\calM = \calM^{\otimes}_{[0]}$
and $\calN = \calN^{\otimes}_{[0]}$ as $\infty$-categories left-tensored over $\calC$. 
A {\it $\calC$-linear functor} from $\calM$ to $\calN$ is a functor $F: \calM^{\otimes} \rightarrow \calN^{\otimes}$ with the following properties:
\begin{itemize}
\item[$(i)$] The diagram
$$ \xymatrix{ \calM^{\otimes} \ar[dr]^{p} \ar[rr]^{F} & & \calN^{\otimes} \ar[dl] \\
& \calC^{\otimes} & }$$
is commutative.
\item[$(ii)$] The functor $F$ carries locally $p$-coCartesian morphisms of $\calM^{\otimes}$ to
locally $p'$-coCartesian morphisms of $\calN^{\otimes}$.
\end{itemize}
We let $\LinFunc{}{\calC}{\calM}{\calN}$ denote the full subcategory of
$\Fun_{\calC^{\otimes}}( \calM^{\otimes}, \calN^{\otimes})$ spanned by the $\calC$-linear functors from $\calM$ to $\calN$.
\end{definition}

\begin{remark}
Note that restriction to the fiber over $[0] \in \Nerve(\cDelta)^{op}$ induces a forgetful functor
$\theta: \LinFunc{}{\calC}{\calM}{\calN} \rightarrow \Fun( \calM, \calN)$. If $\calK$ is a collection of simplicial sets such that both $\calM$ and $\calN$ admit $\calK$-indexed colimits, we let
$\LinFunc{\calK}{\calC}{\calM}{\calN}$ denote the full subcategory of $\LinFunc{}{\calC}{\calM}{\calN}$ spanned by those functors $F$ such that $\theta(F): \calM \rightarrow \calN$ preserves $\calK$-indexed colimits.
\end{remark}

\begin{remark}\label{slayed}
Let $\calC^{\otimes}$ be a monoidal $\infty$-category, and let
$F: \calM^{\otimes} \rightarrow \calN^{\otimes}$ be a $\calC$-linear functor from
$\calM$ to $\calN$. Then composition with $F$ determines a commutative diagram
$$ \xymatrix{ \Mod^{L}( \calM) \ar[rr] \ar[dr] & & \Mod^{L}(\calN) \ar[dl] \\
& \Alg(\calC). & }$$
In particular, for every algebra object $A \in \Alg_{\calC}$, we have an induced functor
$\Mod^{L}_{A}(\calM) \rightarrow \Mod^{L}_{A}(\calN)$. This construction depends functorially
on $F$ in an obvious sense, so we get a functor
$$ \LinFunc{}{\calC}{\calM}{\calN} \rightarrow \Fun( \Mod^{L}_{A}(\calM), \Mod^{L}_{A}(\calN) ).$$
\end{remark}

We can now give a precise statement of $(A)$: 

\begin{theorem}\label{bix2}
Let $\calK$ be a collection of simplicial sets which includes $\Nerve( \cDelta)^{op}$,
let $\calC^{\otimes}$ be a monoidal $\infty$-category, and $\calM^{\otimes} \rightarrow \calC^{\otimes}$ an $\infty$-category left-tensored over $\calC$. Assume that $\calC$ and $\calM$ admit $\calK$-indexed colimits, and that the tensor product functors
$$\calC \times \calC \rightarrow \calC \quad \quad \calC \times \calM \rightarrow \calM$$
preserve $\calK$-indexed colimits separately in each variable. Let $A$ be an algebra object of $\calC$, and let $\theta$ denote the composition
$$ \LinFunc{\calK}{\calC}{\Mod^{R}_A(\calC)}{\calM} 
\subseteq \LinFunc{}{\calC}{ \Mod^{R}_A(\calC)}{\calM}
\stackrel{\theta'}{\rightarrow}
\Fun( \Mod^{L}_{A}( \Mod^{R}_{A}(\calC)), \Mod^{L}_{A}(\calM)) \stackrel{\theta''}{\rightarrow} \Mod^{L}_A(\calM),$$
where $\theta'$ is the map described in Remark \ref{slayed} and
$\theta''$ is given by evaluation at the $A$-bimodule given by $A$. Then
$\theta$ is an equivalence of $\infty$-categories.
\end{theorem}

We turn to assertion $(B)$. In order to make sense of the relative tensor product
$\calM \otimes_{ \calC} \Mod_{A}^{R}(\calC)$, we need to interpret each factor as an object of
a relevant $\infty$-category. To this end, let us recall a bit of notation. Fix a collection
of simplicial sets $\calK$. We let $\Cat_{\infty}(\calK)$ be the subcategory of
$\Cat_{\infty}$ whose objects are $\infty$-categories which admit $\calK$-indexed colimits and
whose morphisms are functors which preserve $\calK$-indexed colimits, and regard
$\Cat_{\infty}(\calK)$ as endowed with the (symmetric) monoidal structure described in \S \symmetricref{comm7}. 

The basic features of $\Cat_{\infty}(\calK)$ are summarized in the following result:

\begin{lemma}\label{spuke}
Let $\calK$ be a small collection of simplicial sets. Then
the $\infty$-category $\Cat_{\infty}(\calK)$ is presentable, and the tensor product
$\otimes: \Cat_{\infty}(\calK) \times \Cat_{\infty}(\calK) \rightarrow \Cat_{\infty}(\calK)$
preserves small colimits separately in each variable.
\end{lemma}

\begin{proof}
We first show that $\Cat_{\infty}(\calK)$ admits small colimits. Let $\calJ$ be an
$\infty$-category, and let $\chi: \calJ \rightarrow \Cat_{\infty}(\calK)$ be a diagram.
Let $\chi'$ denote the composition
$$ \calJ \stackrel{\chi}{\rightarrow} \Cat_{\infty}(\calK) \subseteq \Cat_{\infty},$$
and let $\calC$ be a colimit of the diagram $\chi'$ in $\Cat_{\infty}$. Let
$\calR$ denote the collection of all diagrams in $\calC$ given by a composition
$$ K^{\triangleright} \stackrel{p}{\rightarrow} \chi(J) \rightarrow \calC,$$
where $K \in \calK$ and $p$ is a colimit diagram. It follows from Proposition
\toposref{cupper1} that there exists a functor $F: \calC \rightarrow \calD$ with the following
properties:
\begin{itemize}
\item[$(i)$] For every diagram $q: K^{\triangleright} \rightarrow \calC$ belonging
to $\calR$, the composition $F \circ q$ is a colimit diagram.
\item[$(ii)$] The $\infty$-category $\calD$ admits $\calK$-indexed colimits.
\item[$(iii)$] For every $\infty$-category $\calE$ which admits $\calK$-indexed colimits, composition with $F$ induces an equivalence from the full subcategory of $\Fun( \calD, \calE)$ spanned by those
functors which preserve $\calK$-indexed colimits to the full subcategory of
$\Fun(\calC, \calE)$ spanned by those functors such that the composition
$\chi(J) \rightarrow \calC \rightarrow \calE$ preserves $\calK$-indexed colimits for each
$J \in \calJ$.
\end{itemize}
The map $F$ allows us to promote $\calD$ to an object of
$\overline{\calD} = ( \Cat_{\infty})_{\chi'/ }$. Using $(i)$ and $(ii)$, we deduce that $\overline{\calD}$ lies
in the subcategory $\Cat_{\infty}(\calK)_{ \chi/} \subseteq (\Cat_{\infty})_{ \chi'/}$, and $(iii)$ that this lifting exhibits $\calD$ as a colimit of the diagram $\chi$.

We next show that the tensor product $\otimes: \Cat_{\infty}(\calK) \times \Cat_{\infty}(\calK) \rightarrow \Cat_{\infty}(\calK)$ preserves small colimits separately in each variable. It will suffice to show that
for every object $\calC \in \Cat_{\infty}(\calK)$, the operation $\calD \mapsto \calC \otimes \calD$ admits a right adjoint. This right adjoint is given by the formula $\calE \mapsto \Fun^{\calK}(\calC, \calE)$, where
$\Fun^{\calK}(\calC, \calE)$ denotes the full subcategory of $\Fun( \calC, \calE)$ spanned by those functors which preserve $\calK$-indexed colimits.

We now complete the proof by showing that $\Cat_{\infty}(\calK)$ is presentable.
Fix an uncountable regular cardinal $\kappa$ so that $\calK$ is $\kappa$-small
and every simplicial set $K \in \calK$ is $\kappa$-small. Choose another regular
cardinal $\tau$ such that $\kappa < \tau$ and $\kappa \ll \tau$: that is, $\tau_0^{\kappa_0} < \tau$ whenever $\tau_0 < \tau$ and $\kappa_0 < \kappa$. Let $\Cat_{\infty}^{\tau}(\calK)$ denote the full subcategory
of $\Cat_{\infty}^{\tau}(\calK)$ spanned by those $\infty$-categories $\calC$ which are $\tau$-small
and admit $\calK$-indexed colimits. Then $\Cat_{\infty}^{\tau}(\calK)$ is an essentially small $\infty$-category; it will therefore suffice to prove that every object $\calC \in \Cat_{\infty}(\calK)$ is the colimit
(in $\Cat_{\infty}(\calK)$) of a diagram taking values in $\Cat_{\infty}^{\tau}(\calK)$.

Let $A$ be the collection of all simplicial subsets $\calC_0 \subseteq \calC$ with the following properties:
\begin{itemize}
\item[$(a)$] The simplicial set $\calC_0$ is an $\infty$-category.
\item[$(b)$] The $\infty$-category $\calC_0$ admits $\calK$-indexed colimits.
\item[$(c)$] The inclusion $\calC_0 \hookrightarrow \calC$ preserves $\calK$-indexed colimits.
\item[$(d)$] The simplicial set $\calC_0$ is $\tau$-small.
\end{itemize}
Our proof rests on the following claim:
\begin{itemize}
\item[$(\ast)$] For every $\tau$-small simplicial subset $\calC_0 \subseteq \calC$, there
exists a $\tau$-small simplicial subset $\calC'_0 \subseteq \calC$ which contains
$\calC_0$ and belongs to $A$.
\end{itemize}
Let us regard the set $A$ as partially ordered with respect to inclusions, and we
have an evident functor $\rho: A \rightarrow \sSet$.
From assertion $(\ast)$, it follows that $A$ is filtered (in fact, $\tau$-filtered) and that 
$\calC$ is the colimit of the diagram $\rho$ (in the ordinary category $\sSet$).
Since the collection of categorical equivalences in $\sSet$ is stable under filtered colimits,
we deduce that $\calC$ is the homotopy colimit of the diagram $\rho$ (with respect to the Joyal model structure), so that $\calC$ is the colimit of the induced diagram $\Nerve(\rho): \Nerve(A) \rightarrow \Cat_{\infty}$ (Theorem \toposref{colimcomparee}). Requirement $(c)$ guarantees that
every inclusion $\calC_0 \subseteq \calC_1$ between elements of $A$ is a functor which
preserves $\calK$-indexed colimits, so that $\Nerve(\rho)$ factors through $\Cat_{\infty}(\calK)$.
We claim that $\calC$ is a colimit of the diagram $\Nerve(\rho)$ in the $\infty$-category
$\Cat_{\infty}(\calK)$. Unwinding the definitions, this amounts to the following assertion:
for every $\infty$-category $\calE$ which admits $\calK$-indexed colimits, a functor
$F: \calC \rightarrow \calE$ preserves $\calK$-indexed colimits if and only if
$F| \calC_0$ preserves $\calK$-indexed colimits for each $\calC_0 \in A$. The
``only if'' direction is obvious. To prove the converse, choose $K \in \calK$ and
a colimit diagram $p: K^{\triangleright} \rightarrow \calC$. The image of
$p$ is $\tau$-small, so that $(\ast)$ guarantees that $p$ factors through
$\calC_0$ for some $\calC_0 \in A$. Requirement $(c)$ guarantees that
$p$ is also a colimit diagram in $\calC_0$, so that $F \circ p$ is a colimit diagram
provided that $F | \calC_0$ preserves $\calK$-indexed colimits.

It remains only to prove assertion $(\ast)$. Fix a $\tau$-small subset $\calC_0 \subseteq \calC$.
We define a transfinite sequence of $\tau$-small simplicial subsets $\{ \calC_{\alpha} \subseteq \calC \}_{\alpha < \kappa}$. If $\alpha$ is a nonzero limit ordinal, we take $\calC_{\alpha} = \bigcup_{\beta < \alpha} \calC_{\beta}$. If $\alpha = \beta + 1$, we define $\calC_{\alpha}$ to be any $\tau$-small
simplicial subset of $\calC$ with the following properties:
\begin{itemize}
\item Every map $\Lambda^n_i \rightarrow \calC_{\beta}$ for $0 < i < n$ extends to an $n$-simplex of $\calC_{\alpha}$.
\item For each $K \in \calK$ and each map $q: K \rightarrow \calC_{\beta}$, there
exists an extension $\overline{q}: K^{\triangleright} \rightarrow \calC_{\alpha}$
which is a colimit diagram in $\calC$.
\item Given $n > 0$ and a map $f: K \star \bd \Delta^n \rightarrow \calC_{\beta}$ such
that the restriction $f| K \star \{0\}$ is a colimit diagram in $\calC$, there exists
a map $\overline{f}: K \star \Delta^n \rightarrow \calC_{\alpha}$ extending $f$.
\end{itemize}
Our assumption that $\kappa \ll \tau$ guarantees that we can satisfy these conditions
by adjoining a $\tau$-small set of simplices to $\calC_{\beta}$. Let $\calC'_0 = \bigcup_{\alpha < \kappa} \calC_{\alpha}$. Then $\calC'_0$ contains $\calC_0$, and belongs to $A$ as desired.
\end{proof}

Now suppose that we are given a monoidal $\infty$-category $\calC^{\otimes}$
and maps $\calM^{\otimes} \rightarrow \calC^{\otimes} \leftarrow \calN^{\otimes}$ which
exhibit $\calM = \calM^{\otimes}_{[0]}$ as an $\infty$-category which is right-tensored
over $\calC$ and $\calN = \calN^{\otimes}_{[0]}$ as an $\infty$-category which is left-tensored over $\calC$. Fixing a collection of simplicial sets $\calK$, we further assume that $\calC$, $\calM$, and $\calN$ admit $\calK$-indexed colimits, and that the tensor product functors
$$ \calM \times \calC \rightarrow \calM \quad \quad \calC \times \calC \rightarrow \calC \quad \quad \calC \times \calN \rightarrow \calN$$
preserve $\calK$-indexed colimits separately in each variable. We can identify
$\calC^{\otimes}$ with an associative algebra object of $\Cat_{\infty}(\calK)$, and
the $\infty$-categories $\calM$ and $\calN$ with right and left modules over this associative algebra,
respectively (see Remark \monoidref{otherlander} and Corollary \monoidref{spg}). Consequently,
we can define the tensor product $\calM \otimes_{\calC} \calN$ by applying the construction of
\S \monoidref{balpair}: namely, we first form the two-sided bar construction
$\Baar_{\calC}(\calM, \calN)_{\bigdot}$ of Definition \monoidref{barcon}, and then take the geometric realization in the $\infty$-category $\Cat_{\infty}(\calK)$.

The bar complex $\Baar_{\calC}(\calM, \calN)_{\bigdot}$ is given informally by the formula
$[n] \mapsto \calM \otimes \calC^{(\otimes n)} \otimes \calN$, where the tensor product is
formed in the $\infty$-category $\Cat_{\infty}(\calK)$. Consequently, this bar construction is dependent on the choice of the collection $\calK$. To emphasize this dependence, we will denote
$\Baar_{\calC}(\calM, \calN)_{\bigdot}$ by $\Baar_{\calC}^{\calK}(\calM, \calN)_{\bigdot}$. If $\calK' \subseteq \calK$, then we have a forgetful functor $\Cat_{\infty}(\calK) \rightarrow \Cat_{\infty}(\calK')$ which is lax monoidal, and induces a natural transformation $\Baar_{\calC}^{\calK'}(\calM, \calN)_{\bigdot} \rightarrow \Baar_{\calC}^{\calK}( \calM, \calN)$. In particular, we have a map
$\theta: \Baar_{\calC}^{\emptyset}( \calM, \calN)_{\bigdot} \rightarrow \Baar_{\calC}^{\calK}(\calM, \calN)$. The map $\theta$ is characterized by the following universal property:

\begin{itemize}
\item[$(\ast)$] For each $n \geq 0$ and every $\infty$-category $\calE$ which admits $\calK$-indexed colimits, composition with $\theta$ induces an equivalence from the full subcategory of
$\Fun( \Baar_{\calC}^{\calK}(\calM, \calN)_{n}, \calE)$ spanned by those functors which preserve
$\calK$-indexed colimits to the full subcategory of
$\Fun( \Baar_{\calC}^{\emptyset}( \calM, \calN)_{n}, \calE) \simeq \Fun( \calM \times \calC^{n} \times \calN, \calE)$ spanned by those functors which preserve $\calK$-indexed colimits separately in each variable.
\end{itemize}

We can identify $\Baar_{\calC}^{\emptyset}( \calM, \calN)_{\bigdot}$ with a simplicial object $\Nerve(\cDelta)^{op} \rightarrow \Cat_{\infty}$. Unwinding the definitions, we see that this object is classified by the coCartesian fibration $q: \calM^{\otimes} \times_{\calC^{\otimes} } \calN^{\otimes} \rightarrow \Nerve(\cDelta)^{op}$. Proposition \toposref{charcatcolimit} allows us to identify the geometric realization $| \Baar_{\calC}^{\emptyset}( \calM, \calN)_{\bigdot} |$ with the $\infty$-category
obtained from $\calM^{\otimes} \times_{ \calC^{\otimes} } \calN^{\otimes}$ obtained by inverting all
of the $q$-coCartesian morphisms. Combining this observation with $(\ast)$, we obtain
the following concrete description of the relative tensor product in $\Cat_{\infty}(\calK)$:

\begin{lemma}\label{sider}
Let $\calK$ be a small collection of simplicial sets and let
$\calM^{\otimes} \rightarrow \calC^{\otimes} \leftarrow \calN^{\otimes}$ be as above.
For every $\infty$-category $\calE$ which admits $\calK$-indexed colimits, the natural transformation
$\Baar_{\calC}^{\emptyset}(\calM, \calN)_{\bigdot} \rightarrow \Baar_{\calC}^{\calK}(\calM, \calN)_{\bigdot}$ induces an equivalence of $\infty$-categories from the full subcategory of
$\Fun( \calM \otimes_{\calC} \calN, \calE)$ spanned by those functors which preserve $\calK$-indexed colimits to the full subcategory of $\Fun( \calM^{\otimes} \times_{\calC^{\otimes}} \calN^{\otimes}, \calE)$
spanned by those functors $F$ with the following properties:
\begin{itemize}
\item[$(i)$] The functor $F$ carries $q$-coCartesian morphisms to equivalences in $\calE$, where
$q: \calM^{\otimes} \times_{ \calC^{\otimes} } \calN^{\otimes} \rightarrow \Nerve( \cDelta)^{op}$ denotes the canonical projection.
\item[$(ii)$] For each $n \geq 0$, the functor
$$\calM \times \calC^{n} \times \calN \simeq q^{-1} \{ [n] \} \rightarrow \calE$$
preserves $\calK$-indexed colimits separately in each variable.
\end{itemize}
\end{lemma}

Our next goal is to apply Lemma \ref{sider} to construct a canonical map
$\calM \otimes_{\calC} \Mod_{A}^{R}(\calC) \rightarrow \Mod_{A}^{R}(\calM)$.
Fix a monoidal $\infty$-category $\calC^{\otimes}$ and a map
$\calM^{\otimes} \rightarrow \calC^{\otimes}$ which exhibits $\calM$ as an $\infty$-category right-tensored over $\calC$. For this, we need a variant on Construction \ref{cloops}:

\begin{construction}\label{cloopsprime}
Let  $\pi_0 \stackrel{\alpha}{\leftarrow} c \stackrel{\beta}{\rightarrow} \pi_1$ be defined as in Construction \ref{cloops}, let $p: \calC^{\otimes} \rightarrow \Nerve(\cDelta)^{op}$ be a monoidal $\infty$-category, and let
$q: \calM^{\otimes} \rightarrow \calC^{\otimes}$ be an $\infty$-category right-tensored
over $\calC$. We define an $\infty$-category $\overline{B}( \calM^{\otimes})$ 
equipped with a map 
$\overline{B} \rightarrow \Nerve( \cDelta)^{op} \times \Nerve( \cDelta)^{op}$ so that the following universal property is satisfied: for every map of simplicial sets
$K \stackrel{\phi}{\rightarrow} \Nerve( \cDelta)^{op} \times \Nerve( \cDelta)^{op}$, there is a canonical bijection
of $\Hom_{(\sSet)_{/ \Nerve( \cDelta)^{op} \times \Nerve( \cDelta)^{op} \times S}}( K, \overline{B}( \calM^{\otimes}))$ with the collection of all commutative diagrams
$$ \xymatrix{ \Nerve(\cDelta)^{op} \times \Nerve( \cDelta)^{op} \times \Delta^1 \ar[d]^{\alpha} & K \times \Delta^1 \ar[l]_-{\phi \times \id} \ar[d] & \calK \times \{0\} \ar[l] \ar[r] \ar[d] & \calK \times \Delta^1 \ar[r]^-{ \phi \times \id} \ar[d] & \Nerve( \cDelta)^{op} \times \Nerve( \cDelta)^{op} \times \Delta^1 \ar[d]^{\beta} \\
\Nerve( \cDelta)^{op} & \calM^{\otimes} \ar[l]_{pq} & \calM^{\otimes} \ar[l]_{\id} \ar[r]^{q} & \calC^{\otimes} \ar[r]^-{p} & \Nerve( \cDelta)^{op}. }$$
More informally, $\overline{B}$ is the $\infty$-category whose fiber over $([m], [n]) \in \Nerve(\cDelta)^{op} \times \Nerve( \cDelta)^{op}$ consists of a triple
$(M, f: M \rightarrow M_0, g: q(M) \rightarrow C)$, where $M \in \calM^{\otimes}_{[m] \star [n]}$, 
$f$ is a morphism in $\calM^{\otimes}$ covering the inclusion $[m] \hookrightarrow [m] \star [n]$,
and $g$ is a morphism in $\calC^{\otimes}$ covering the inclusion $[n] \hookrightarrow [m] \star [n]$.
Let $B( \calM^{\otimes})$ denote the full subcategory of $\overline{B}( \calM^{\otimes})$ spanned by
those objects for which $f$ is $(p \circ q)$-coCartesian and $g$ is $p$-coCartesian.
Composition with $q$ induces a map $B( \calM^{\otimes}) \rightarrow B( \calC^{\otimes})$, where
$B( \calC^{\otimes})$ is defined as in Construction \ref{cloops}. We therefore obtain a map
$$\theta: B( \calM^{\otimes}) \rightarrow B( \calC^{\otimes}) \times_{ \calC^{\otimes} \times \calC^{\otimes}} ( \calM^{\otimes} \times \calC^{\otimes}).$$
which is easily shown to be a trivial Kan fibration. Let $s$ denote a section of the map $\theta$.

Choose a $(p \circ q)$-coCartesian natural transformation $$
h: (\calM^{\otimes} \times_{ \Nerve( \cDelta)^{op} } ( \Nerve( \cDelta)^{op} \times \Nerve( \cDelta)^{op} )) \times \Delta^1 \rightarrow \calM^{\otimes}$$
covering the map induced by $\beta$. Then $h$ induces a map $h': \overline{B}( \calM^{\otimes}) \rightarrow \calM^{\otimes}$. By construction, the diagram
$$ \xymatrix{ B( \calM^{\otimes}) \ar[r]^{h'} \ar[d] & \calM^{\otimes} \ar[d]^{q} \\
\calM^{\otimes} \times \calC^{\otimes} \ar[r] & \calC^{\otimes} }$$
commutes up to a canonical equivalence. Since $q$ is a categorical fibration, we can
modify $h'$ by a homotopy to guarantee that this diagram is strictly commutative.
The composition $h' \circ s$ is a functor 
$h'': \calM^{\otimes} \times_{ \calC^{\otimes} } B( \calC^{\otimes} ) \rightarrow \calM^{\otimes}.$
For any algebra object $A$ in $\calC^{\otimes}$, composition with $h''$ induces a map
$\Psi: \calM^{\otimes} \times_{ \calC^{\otimes} } \Mod_{A}^{R}(\calC)^{\otimes}
\rightarrow \Mod_{A}^{R}(\calM)$. 
\end{construction}

We can think of objects of the fiber product $\calM^{\otimes} \times_{ \calC^{\otimes}} \Mod_{A}^{R}(\calC)^{\otimes}$ as finite sequences
$(M, C_1, \ldots, C_n, N)$, where $M$ is an object of $\calM$, each $C_i$ is an object of
$\calC$, and $N$ is a right $A$-module object of $\calC$. The functor $\Psi$ of Construction \ref{cloopsprime} is given informally by the formula
$(M, C_1, \ldots, C_n, N) \mapsto M \otimes C_1 \otimes \ldots \otimes C_n \otimes N$.
From this description, it is clear that $\Psi$ carries $r$-coCartesian morphisms in
$\calM^{\otimes} \times_{ \calC^{\otimes}} \Mod_{A}^{R}(\calC)^{\otimes}$ to
equivalences in $\Mod_{A}^{R}(\calM)$, where $r: \calM^{\otimes} \times_{ \calC^{\otimes}} \Mod_{A}^{R}(\calC)^{\otimes} \rightarrow \Nerve( \cDelta )^{op}$ is the projection.
Suppose furthermore that $\calK$ is a collection of simplicial sets such that $\calC$ and $\calM$ admit $\calK$-indexed colimits, and that the tensor product functors
$$ \calC \times \calC \rightarrow \calC \quad \quad \calM \times \calC \rightarrow \calM$$
preserve $\calK$-indexed colimits separately in each variable. It then follows from Corollary
\monoidref{gloop} that for each $n \geq 0$, the functor
$$ \calM \times \calC^{n} \times \Mod_{A}^{R}(\calC)
\simeq r^{-1} \{ [n] \} \stackrel{\Psi}{\rightarrow} \Mod_{A}^{R}(\calM)$$
preserves $\calK$-indexed colimits separately in each variable. It follows from Lemma \ref{sider} that
$\Psi$ determines a functor $\Phi: \calM \otimes_{ \calC} \Mod_{A}^{R}(\calC) \rightarrow \Mod_{A}^{R}(\calM)$, which is well-defined up to equivalence.

We can now formulate $(B)$ as follows:

\begin{theorem}\label{bix1}
Let $\calK$ be a small collection of simplicial sets which includes $\Nerve( \cDelta)^{op})$,
let $\calC^{\otimes}$ be a monoidal $\infty$-category, let $\calM^{\otimes} \rightarrow \calC^{\otimes}$ be an $\infty$-category right-tensored over $\calC$, and let $A \in \Alg(\calC)$ be an algebra object
of $\calC$. Suppose that the $\infty$-categories $\calC$ and $\calM$ admit $\calK$-indexed colimits, and the tensor product functors
$$ \calC \times \calC \rightarrow \calC \quad \quad \calM \times \calC \rightarrow \calM$$
preserve $\calK$-indexed colimits separately in each variable. Then the above construction
yields an equivalence of $\infty$-categories $\calM \otimes_{\calC} \Mod_{A}^{R}(\calC) \rightarrow
\Mod_{A}^{R}(\calM)$, where the tensor product is taken in $\Cat_{\infty}(\calK)$.
\end{theorem}

\begin{remark}
In the formulation of Theorems \ref{bix1} and \ref{bix2}, the $\infty$-categories $\Mod_{A}^{R}(\calM)$ and $\Mod_{A}^{L}(\calM)$ do not depend on the class of simplicial sets $\calK$. It follows that the $\infty$-categories $\calM \otimes_{\calC} \Mod_{A}^{R}(\calC)$ and 
$\LinFunc{\calK}{\calC}{ \Mod_{A}^{R}(\calC)}{\calM}$ do not depend on $\calK$, provided that $\calK$ contains $\Nerve( \cDelta)^{op}$.
\end{remark}

The proofs of Theorem \ref{bix1} and \ref{bix2} are very similar, and rest on an analysis of the
forgetful functor $\Mod^{R}_{A}(\calC) \rightarrow \calC$. We observe that this functor is
$\calC$-linear. More precisely, evaluation at the point $[0]$ in $\Nerve( \cDelta)^{op}$ induces
a $\calC$-linear functor $G: \Mod^{R}_{A}(\calC)^{\otimes} \rightarrow \calN^{\otimes}$,
where $\calN^{\otimes} \rightarrow \calC^{\otimes}$ exhibits $\calN \simeq \calC$ as left-tensored over
itself (see Example \monoidref{sumai}). We have the following fundamental observation:

\begin{lemma}\label{pax}
Let $\calC^{\otimes}$ be a monoidal $\infty$-category containing an algebra object $A$,
and let $\calN^{\otimes} \rightarrow \calC^{\otimes}$ exhibit $\calN \simeq \calC$ as left-tensored over itself as in Example \monoidref{sumai}. Consider the commutative diagram
$$ \xymatrix{ \calN^{\otimes} \ar[dr]^{q} & & \Mod^{R}_{A}(\calC)^{\otimes} \ar[dl]^{p'} \ar[ll]^{G} \\
& \calC^{\otimes}, & }$$
where $G$ is defined as above. Then there exists a functor
$F: \calN^{\otimes} \rightarrow \Mod^{R}_{A}(\calC)^{\otimes}$ and a natural transformation
$u: \id_{ \calN^{\otimes}} \rightarrow G \circ F$ which exhibits $F$ as a left adjoint to $G$
relative to $\calC^{\otimes}$ (see Definition \deformationref{coughly} and Remark \deformationref{coughco}). Moreover, $F$ is a $\calC$-linear functor from $\calC$ to $\Mod^{R}_{A}(\calC)$.
\end{lemma}

\begin{remark}
More informally, Lemma \ref{pax} asserts that the forgetful functor
$\Mod^{R}_{A}(\calC) \rightarrow \calC$ and its left adjoint $C \mapsto C \otimes A$
commute with the action of $\calC$ by left multiplication.
\end{remark}

\begin{proof}
We observe that $p$ and $p'$ are locally coCartesian fibrations (Lemma \monoidref{excel}).
Moreover, for each object $C \in \calC^{\otimes}$, the induced map on fibers
$\Mod^{R}_{A}(\calC)^{\otimes}_{C} \rightarrow \calN^{\otimes}_{C}$ is equivalent to the forgetful functor $\theta: \Mod^{R}_{A}(\calC) \rightarrow \calC$, and therefore admits a left adjoint (Proposition \monoidref{pretara}). We complete the proof by observing that the functor $G$ satisfies
hypothesis $(2)$ of Proposition \deformationref{newadjprop}. Unwinding the definitions, this
results from the observation that the canonical maps
$(C \otimes D) \otimes A \rightarrow C \otimes (D \otimes A)$ expressing the coherent
associativity of the tensor product on $\calC$ are equivalences in $\calC$.
\end{proof}

\begin{lemma}\label{costanew}
Let $p: \calC^{\otimes} \rightarrow \Nerve( \cDelta)^{op}$ be a monoidal $\infty$-category, and let $q: \calM^{\otimes} \rightarrow \calC^{\otimes}$ and $\calN^{\otimes} \rightarrow \calC^{\otimes}$ be $\infty$-categories
left-tensored over $\calC$. Let $K$ be a simplicial set such that $\calN$ admits $K$-indexed
colimits, and such that for each $C \in \calC$, the tensor product functor
$\{ C \} \times \calN \subseteq \calC \times \calN \rightarrow \calN$ preserves $K$-indexed colimits.
Then:
\begin{itemize}
\item[$(1)$] The $\infty$-category $\LinFunc{}{\calC}{\calM}{\calN}$ admits $K$-indexed colimits.
\item[$(2)$] A map $f: K^{\triangleright} \rightarrow \LinFunc{}{\calC}{\calM}{\calN}$
is a colimit diagram if and only if, for each $M \in \calM$, the induced map
$K^{\triangleright} \rightarrow \calN$ is a colimit diagram.
\end{itemize}
\end{lemma}

\begin{remark}\label{costu}
In the situation of Lemma \ref{costanew}, suppose that $\calK$ is a class of simplicial sets such that
$\calC$, $\calM$, and $\calN$ admit $\calK$-indexed colimits, and the tensor product functors
$$\calC \times \calC \rightarrow \calC \quad \quad \calC \times \calM \rightarrow \calM \quad \quad \calC \times \calN \rightarrow \calN$$
preserve $\calK$-indexed colimits separately in each variable. Then the full subcategory
$\LinFunc{\calK}{\calC}{\calM}{\calN} \subseteq \LinFunc{}{\calC}{\calM}{\calN}$ is stable
under $K$-indexed colimits: this follows from the characterization of $K$-indexed colimits
supplied by Lemma \ref{costanew} together with Lemma \toposref{limitscommute}.
\end{remark}

\begin{proof}
Assertion $(1)$ follows from Proposition \toposref{prestorkus}. Moreover, Proposition
\toposref{prestorkus} guarantees that a diagram $f: K^{\triangleright} \rightarrow \LinFunc{}{\calC}{\calM}{\calN}$ is a colimit if and only if, for every object $M \in \calM^{\otimes}$ having image
$C \in \calC^{\otimes}$ and $[n] \in \Nerve(\cDelta)^{op}$, the induced map $f_M: K^{\triangleright} \rightarrow \calN^{\otimes}_{C}$ is a colimit diagram. The necessity of condition $(2)$ is now obvious. For the sufficiency, we note that if we choose a $p$-coCartesian morphism $\alpha: C \rightarrow C_0$
in $\calC^{\otimes}$ covering the inclusion $[0] \simeq \{n\} \subseteq [n]$ in $\cDelta$
and a locally $q$-coCartesian morphism $M \rightarrow M_0$ lifting $\alpha$, then
we have a homotopy commutative diagram
$$ \xymatrix{ & K^{\triangleright} \ar[dl]^{p_M} \ar[dr]^{p_{M_0}} & \\
\calN^{\otimes}_{C} \ar[rr]^{ \alpha_{!} } & & \calN^{\otimes}_{C_0} }$$
where $\alpha_{!}$ is an equivalence of $\infty$-categories, so that $f_{M}$ is a colimit 
if and only if $f_{M_0}$ is a colimit diagram.
\end{proof}

We now turn to the proofs of Theorems \ref{bix2} and \ref{bix1}.

\begin{proof}[Proof of Theorem \ref{bix2}]
Let $\calN^{\otimes} \rightarrow \calC^{\otimes}$ exhibit the monoidal $\infty$-category $\calC$ as left-tensored over itself, as in Example \monoidref{sumai}, and let $G: \Mod_{A}^{R}(\calC) \rightarrow \calN^{\otimes}$ and $F: \calN^{\otimes} \rightarrow \Mod_{A}^{R}(\calC)$ be as in Lemma \ref{pax}.
Then $F$ and $G$ induce adjoint functors
$$ \Adjoint{f}{\LinFunc{\calK}{\calC}{\calN}{\calM}}{\LinFunc{\calK}{\calC}{\Mod_{A}^{R}(\calC)}{\calM}}{g}$$
We first claim that evaluation at the unit object ${\bf 1} \in \calN \simeq \calC$ induces an equivalence
of $\infty$-categories $\phi: \LinFunc{\calK}{\calC}{\calN}{\calM} \rightarrow \calM$. It will suffice
to show that for every simplicial set $K$, the induced map
$\Fun( K, \LinFunc{\calK}{\calC}{\calN}{\calM}) \rightarrow \Fun(K, \calM)$ induces a bijection
on equivalence classes of objects. Replacing $\calM$ by $\Fun(K, \calM)$, we are reduced to proving that $\phi$ induces a bijection on equivalence classes of objects. On the left hand side, the set of equivalence classes can be identified with $\pi_0 \bHom_{ \Mod^{L}_{ \calC}( \Cat_{\infty}(\calK)}( \calC, \calM)$. Using Proposition \monoidref{pretara}, we can identify this with the set
$\pi_0 \bHom_{ \Cat_{\infty}(\calK)}( \SSet(\calK), \calM) \simeq \pi_0 \bHom_{ \Cat_{\infty} }( \Delta^0, \calM)$, which is the set of equivalence classes of objects of $\calM$ as required.

Let $T: \LinFunc{ \calK}{\calC}{\Mod_{A}^{R}(\calC)}{\calM} \rightarrow \calM$ denote the composition
of the functor $g$ with the equivalence $\phi$. We have a homotopy commutative diagram of $\infty$-categories
$$ \xymatrix{ \LinFunc{ \calK}{\calC}{\Mod_{A}^{R}(\calC)}{\calM} \ar[rr]^{\theta} \ar[dr]^{T} & & 
\Mod^{L}_{A}(\calM) \ar[dl]^{T'} \\
& \calM, & }$$
where $T$ is the evident forgetful functor. We will prove that $\theta$ is an equivalence showing
that this diagram satisfies the hypotheses of Corollary \monoidref{littlerbeck}:

\begin{itemize}
\item[$(a)$] The $\infty$-categories $\LinFunc{\calK}{\calC}{\Mod_{A}^{R}(\calC)}{\calM}$ and
$\Mod^{L}_{A}(\calM)$ admit geometric realizations of simplicial objects. In the first case,
this follows from Lemma \ref{costanew} and Remark \ref{costu}. In the second, it follows from
Corollary \monoidref{gloop}.

\item[$(b)$] The functors $T$ and $T'$ admit left adjoints, which we will denote by
$U$ and $U'$. The left adjoint $U$ is given by composing $f$ with a homotopy inverse to the equivalence $\phi$, and the left adjoint $U'$ is supplied by Proposition \monoidref{pretara}.

\item[$(c)$] The functor $T'$ is conservative and preserves geometric realizations of simplicial objects. 
The first assertion follows from Corollary \monoidref{thetacart} and the second from Corollary \monoidref{gloop}.

\item[$(d)$] The functor $T$ is conservative and preserves geometric realizations of simplicial objects.
The second assertion follows from Lemma \ref{costanew}. To prove the first, suppose that
$\alpha: S \rightarrow S'$ is a natural transformation of $\calC$-linear functors from
$\Mod_{A}^{R}(\calC)$ to $\calM$, each of which preserves $\calK$-indexed colimits, and
that $T(\alpha)$ is an equivalence. We wish to show that $\alpha$ is an equivalence.
Let us abuse notation by identifying $S$ and $S'$ with the underlying maps
$\Mod_{A}^{R}(\calC) \rightarrow \calM$, and let $\calX$ be the full subcategory of
$\Mod_{A}^{R}(\calC)$ spanned by those objects $X$ for which $\alpha$ induces an equivalence
$S(X) \rightarrow S'(X)$ in $\calM$. We wish to show that $\calX = \Mod_{A}^{R}(\calC)$. 

Since $\phi$ is an equivalence of $\infty$-categories, we conclude that $g(\alpha)$ is an
equivalence in $\LinFunc{\calK}{\calC}{\calN}{\calM}$. In other words, the $\infty$-category
$\calX$ contains the essential image of the free module functor $\calC \simeq \calN \rightarrow
\Mod_{A}^{R}(\calC)$. Since $S$ and $S'$ preserve $\calK$-indexed colimits, $\calX$ is stable under
geometric realizations of simplicial objects. The equality $\calX = \Mod_{A}^{R}(\calC)$ now follows from Proposition \monoidref{littlebeck}.

\item[$(e)$] The natural transformation $T' \circ U' \rightarrow T \circ U$ is an equivalence
of functors from $\calM$ to itself. Unwinding the definitions, we see that both of these functors
are given by tensoring with the object $A \in \calC$.
\end{itemize} 
\end{proof}

\begin{proof}[Proof of Theorem \ref{bix1}]
The forgetful functor $\Mod_{A}^{R}(\calC) \rightarrow \calC$ can be viewed as a map between left $\calC$-module objects in $\Cat_{\infty}(\calK)$, and therefore induces a functor $G: \calM \otimes_{\calC} \Mod_{A}^{R}(\calC) \rightarrow
\calM \otimes_{\calC} \calC \simeq \calM$ (where the last equivalence is given by Proposition \monoidref{usss}). Let $G': \Mod_{A}^{R}(\calM) \rightarrow \calM$ be the evident forgetful functor. We have a diagram
$$ \xymatrix{ \calM \otimes_{\calC} \Mod_{A}^{R}(\calC) \ar[dr]^{G} \ar[rr]^{\Phi} & & \Mod_{A}^{R}(\calM) \ar[dl]^{G'} \\
& \calM, & }$$
which commutes up to canonical homotopy. To prove that $\Phi$ is an equivalence of $\infty$-categories, it will suffice to show that this diagram satisfies the hypotheses of Corollary \monoidref{littlerbeck}:
\begin{itemize}
\item[$(a)$] The $\infty$-categories $\calM \otimes_{\calC} \Mod_{A}^{R}(\calC)$ and
$\Mod_{A}^{R}(\calM)$ admit geometric realizations of simplicial objects. In the first case, this follows from our assumption that $\Nerve(\cDelta)^{op} \in \calK$; in the second case, it follows from
Corollary \monoidref{gloop} (since $\calC$ admits geometric realizations and tensor product with $A$ preserves geometric realizations).
\item[$(b)$] The functors $G$ and $G'$ admit left adjoints, which we will denote by $F$ and $F'$.
The existence of $F'$ follows from Proposition \monoidref{pretara} (which also shows that $F'$ is given informally by the formula $M \mapsto M \otimes A$). Similar reasoning shows that the forgetful functor
$\Mod_{A}^{R}(\calC) \rightarrow \calC$ admits a left adjoint $F_0$. It is not difficult to see that this left adjoint can be promoted to a map of $\infty$-categories left-tensored over $\calC$, so that it induces a functor $\id \otimes F_0: \calM \otimes_{\calC} \calC \rightarrow \calM \otimes_{\calC} \Mod_{A}^{R}(\calC)$ which is left adjoint to $G$.

\item[$(c)$] The functor $G'$ is conservative and preserves geometric realizations of simplicial objects.
The first assertion follows from Corollary \monoidref{thetacart} and the second from Corollary \monoidref{gloop} (since $\Nerve(\cDelta)^{op} \in \calK$). 

\item[$(d)$] The functor $G$ is conservative and preserves geometric realizations of simplicial objects.
The second assertion is obvious (since $G$ is a morphism in $\Cat_{\infty}(\calK)$ by construction).
The proof that $G$ is conservative is a bit more involved. Let $c: \Nerve( \cDelta)^{op}
\times \Nerve( \cDelta_{+})^{op} \times \Nerve( \cDelta)^{op} \rightarrow \Nerve(\cDelta)^{op}$ be the concatenation functor, given by 
the formula $([l], [m], [n]) \mapsto [l] \star [m] \star [n] \simeq [l+m+n - 2]$, and let
$c_0: \Nerve( \cDelta)^{op} \times \Nerve( \cDelta_{+})^{op} \times \Nerve( \cDelta)^{op} 
\rightarrow \Nerve( \cDelta)^{op}$ be given by $([l],[m],[n]) \mapsto [l] \star [n]$. The canonical inclusions $[l] \star [n] \hookrightarrow [l] \star [m] \star [n]$ induce a natural transformation of functors
$\alpha: c \rightarrow c_0$. Let $p: \calC^{\otimes} \rightarrow \Nerve( \cDelta)^{op}$ be the canonical map and $\pi: \calC^{\otimes} \times_{ \Nerve(\cDelta)^{op} } ( \Nerve( \cDelta)^{op} \times \Nerve( \cDelta_{+})^{op} \times \Nerve(\cDelta)^{op} ) ) \rightarrow \calC^{\otimes}$ the projection,
and choose a $p$-coCartesian natural transformation $\overline{\alpha}: \pi \rightarrow \pi'$ covering
$\alpha$. Adjusting $\pi'$ by a homotopy if necessary, we can assume that
$\pi'$ induces a functor $T_{\bigdot}: \Nerve( \cDelta_{+})^{op} \times \Mod_{A}^{R}( \calC)^{\otimes}
\rightarrow \Mod_{A}^{R}(\calC)^{\otimes}$. We will view the functor $T_{\bigdot}$ as an augmented simplicial
object in the category of left $\calC$-module functors from $\Mod_{A}^{R}(\calC)$ to itself, given informally by the formula $T_{n}(N) = N \otimes A^{ \otimes (n+1)}$; in particular, the functor
$T_{-1}$ is equivalent to the identity functor. For each object
$N \in \Mod_{A}^{R}(\calC)$, the canonical map $\epsilon: | T_{\bigdot} N | \rightarrow T_{-1} N \simeq N$ is an equivalence: to prove this, it suffices to show that the image of $\epsilon$ under the forgetful functor $\theta: \Mod_{A}^{R}(\calC) \rightarrow \calC$ is an equivalence (Corollary \monoidref{thetacart}): we note that $\theta | T_{\bigdot} N|$ can be identified with the relative
tensor product $N \otimes_{A} A$ and that $\theta( \epsilon)$ is the equivalence of Proposition \monoidref{usss}). 

Using Lemma \ref{sider}, we see that $T_{\bigdot}$ determines an augmented simplicial object
$U_{\bigdot}: \Nerve( \cDelta_{+})^{op} \rightarrow \Fun( \calM \otimes_{\calC} \Mod_{A}^{R}(\calC),
\calM \otimes_{\calC} \Mod_{A}^{R}(\calC) )$. Note that each $U_n$ preserves $\calK$-indexed colimits.
Let $\calX$ be the full subcategory of $\calM \otimes_{\calC} \Mod_{A}^{R}( \calC)$ spanned by those objects for which the canonical map $| U_{\bigdot} X| \rightarrow U_{-1} X \simeq X$ is an equivalence. 
Since each $U_{n}$ preserves $\calK$-indexed colimits, the full subcategory $\calX$ is stable
under $\calK$-indexed colimits in $\calM \otimes_{\calC} \Mod_{A}^{R}(\calC) )$. Since
$\calM \otimes_{\calC} \Mod_{A}^{R}(\calC) )$ is generated under $\calK$-indexed colimits by the essential image of the tensor product functor $\otimes: \calM \times \Mod_{A}^{R}(\calC) \rightarrow \calM \otimes_{\calC} \Mod_{A}^{R}(\calC)$ (which obviously belongs to $\calX$), we conclude that
$\calX = \calM \otimes_{\calC} \Mod_{A}^{R}(\calC) $. 

Now suppose that $f: X \rightarrow Y$ is a morphism in $\calM \otimes_{\calC} \Mod_{A}^{R}(\calC) )$
such that $G(f)$ is an equivalence. We wish to prove that $f$ is an equivalence. Note that
$f$ is equivalent to the geometric realization $| U_{\bigdot} f |$ (in the $\infty$-category
$\Fun( \Delta^1, \calM \otimes_{\calC} \Mod_{A}^{R}(\calC) )$); it therefore suffices to show that
$U_n(f)$ is an equivalence for $n \geq 0$. We complete the argument by observing that $U_n$ factors through $G$ (since $T_n$ factors through the forgetful functor $\Mod_{A}^{R}(\calC) \rightarrow \calC$).

\item[$(e)$] The canonical natural transformation $G' \circ F' \rightarrow G \circ F$
is an equivalence of functors from $\calM$ to itself. This follows easily from the descriptions of $F$ and $F'$ given above: both functors are given by tensor product with $A$.
\end{itemize}
\end{proof}

\subsection{Behavior of the Functor $\Theta$}\label{probix2}

In \S \ref{conl}, we saw that if $\calC^{\otimes} \rightarrow \Nerve( \cDelta)^{op}$ is a monoidal $\infty$-category and $A$ is an algebra object of $\calC$, then the $\infty$-category $\Mod_{A}^{R}(\calC)$ of right
$A$-module objects of $\calC$ has the structure of an $\infty$-category left-tensored over $\calC$.
Moreover, the construction $(\calC^{\otimes}, A) \mapsto ( \calC^{\otimes}, \Mod^{R}_{A}(\calC) )$ determines a functor
$$\Theta: \CatAlg(\calK) \rightarrow \CatMod(\calK)$$
for any collection of simplicial sets $\calK$ which contains $\Nerve(\cDelta)^{op}$ (see Construction \ref{spaas}). In this section, we will apply the main results of \S \ref{probix} (Theorems \ref{bix1} and \ref{bix2}) to establish some basic formal properties of $\Theta$. We can describe our goals more specifically as follows:

\begin{itemize}
\item[$(1)$] If $A$ is an associative ring, we can {\em almost} recover $A$ from the category $\Mod_{A}^{R}$ of
right $A$-modules. More precisely, if we let $M$ denote the ring $A$ itself, regarded as a right
$A$-module, then left multiplication by elements of $A$ determines a canonical isomorphism
$A \rightarrow \End_{A}(M)$. In other words, the data of the associative ring $A$ is equivalent to the data of the category $\Mod_{A}^{R}$ of right $A$-modules together with its distinguished object $M$. An analogous result holds if we replace the category of abelian groups by a more general monoidal $\infty$-category $\calC$:
the functor $\Theta$ induces a fully faithful embedding $\Theta_{\ast}$ from the
$\infty$-category $\CatAlg(\calK)$ to the $\infty$-category of triples $( \calC^{\otimes}, \calM, M)$,
where $\calC^{\otimes}$ is a monoidal $\infty$-category, $\calM$ is an $\infty$-category left-tensored over $\calC$, and $M \in \calM$ is a distinguished object. We refer the reader to Theorem \ref{postcur} for a precise statement. 

\item[$(2)$] If we work in the setting of {\em presentable} $\infty$-categories, then the functor
$\Theta_{\ast}$ admits a right adjoint, which carries a triple $(\calC^{\otimes}, \calM, M)$ to the pair
$(\calC^{\otimes}, A)$, where $A \in \Alg(\calC)$ is the algebra of endomorphisms of $M$ (Theorem \ref{curly}). 

\item[$(3)$] The $\infty$-categories $\CatAlg(\calK)$ and $\CatMod(\calK)$ admit symmetric monoidal
structures, and $\Theta$ can be promoted to a symmetric monoidal functor (Theorem \ref{saylime}). 

\end{itemize}

We begin by addressing a small technical point regarding the behavior of the functor
$\Theta$ with respect to base change:

\begin{proposition}\label{komb}
Let $\calK$ be a small collection of simplicial sets which includes $\Nerve( \cDelta)^{op}$, and
consider the commutative diagram
$$ \xymatrix{ \CatAlg(\calK) \ar[rr]^{\Theta} \ar[dr]^{\phi} & & \CatMod(\calK) \ar[dl]_{\psi} \\
& \Cat_{\infty}^{\Mon}(\calK). & }$$
The functors $\phi$ and $\psi$ are coCartesian fibrations, and the functor $\Theta$ carries
$\phi$-coCartesian morphisms to $\psi$-coCartesian morphisms.
\end{proposition}

\begin{proof}We first show that $\phi$ is a coCartesian fibration. Let $\widetilde{ \Cat}_{\infty}^{\Mon, \calK}$
be as defined in Notation \ref{kost}, and let $\calX = \Fun( \Nerve( \cDelta)^{op}, \widetilde{ \Cat}_{\infty}^{\Mon, \calK}) \times_{ \Fun( \Nerve( \cDelta)^{op}, \Nerve( \cDelta)^{op} \times \Cat_{\infty}^{\Mon}(\calK))} \Cat_{\infty}^{\Mon}(\calK)$.
Let us denote an object of $\calX$ by a pair $(\calC^{\otimes}, A)$, where
$\calC^{\otimes}$ is a monoidal $\infty$-category (compatible with $\calK$-indexed colimits)
and $A \in \Fun_{ \Nerve( \cDelta)^{op}}( \Nerve( \cDelta)^{op}, \calC^{\otimes})$. 
It follows from Proposition \toposref{doog} that the projection map $\phi': \calX \rightarrow \Cat_{\infty}^{\Mon}(\calK)$ is a coCartesian fibration; moreover, a morphism
$( \calC^{\otimes}, A) \rightarrow ( \calD^{\otimes}, B)$ in $\calX$ is $\phi'$-coCartesian if and only if
the underlying map $F(A) \rightarrow B$ is an equivalence, where $F: \calC^{\otimes} \rightarrow \calD^{\otimes}$ denotes the underlying monoidal functor. In this case, if $A$ is an algebra object of $\calC^{\otimes}$, then $B \simeq F(A)$ is an algebra object of $\calD^{\otimes}$. Note that
$\CatAlg(\calK)$ can be identified with the full subcategory of $\calX$ spanned by those pairs
$(\calC^{\otimes},A)$ where $A$ is an algebra object of $\calC^{\otimes}$. It follows that if
$f: X \rightarrow Y$ is a $\phi'$-coCartesian morphism of $\calX$ such that $X \in \CatAlg(\calK)$, then
$Y \in \CatAlg(\calK)$. We conclude that $\phi = \phi' | \CatAlg(\calK)$ is again a coCartesian fibration, and that a morphism in $\CatAlg(\calK)$ is $\phi$-coCartesian if and only if it is $\phi'$-coCartesian.

We next prove that $\psi$ is a coCartesian fibration. Note that the equivalence $\CatMod \simeq \Mod^{L}(\Cat_{\infty})$ of
Corollary \monoidref{spg} restricts to an equivalence $\CatMod(\calK) \simeq \Mod^{L}( \Cat_{\infty}(\calK) )$. Moreover, the functor $\psi$ can be identified with the forgetful functor
$\Mod^{L}( \Cat_{\infty}(\calK)) \rightarrow \Alg( \Cat_{\infty}(\calK) )$, and is therefore
a Cartesian fibration (Corollary \monoidref{thetacart}). Consequently, to prove that
$\psi$ is a coCartesian fibration, it will suffice to show that for every morphism
$F: \calC^{\otimes} \rightarrow \calD^{\otimes}$ in the $\infty$-category $\Alg( \Cat_{\infty}(\calK) ) \simeq \Cat_{\infty}^{\Mon}(\calK)$, the forgetful functor $\Mod^{L}_{ \calD^{\otimes}}( \Cat_{\infty}(\calK) )
\rightarrow \Mod^{L}_{\calC^{\otimes}}( \Cat_{\infty}(\calK) )$ admits a left adjoint. This follows
immediately from Lemmas \ref{spuke} and \monoidref{peacestick}.

It remains only to prove that $\Theta$ carries $\phi$-coCartesian morphisms to
$\psi$-coCartesian morphisms. Unwinding the definitions, we must show that
if $F: \calC^{\otimes} \rightarrow \calD^{\otimes}$ is a morphism in $\Alg( \Cat_{\infty}(\calK))$ and
$A$ is an algebra object of $\calC^{\otimes}$, then the canonical map 
$\rho: \calD \otimes_{\calC} \Mod^{R}_{A}(\calC) \rightarrow \Mod^{R}_{FA}(\calD)$ is an
equivalence of $\infty$-categories (each of which is left-tensored over $\calD$). 
Note that $\Mod^{R}_{FA}(\calD)$ can be identified with the $\infty$-category
$\Mod^{R}_{A}(\calD)$, where we regard $\calD$ as right-tensored over the $\infty$-category $\calC$ via the monoidal functor $F$. Under this identification, the functor $\rho$ is given by
the equivalence of Theorem \ref{bix1}.
\end{proof}

Let us now study the image of the initial object of $\CatAlg(\calK)$ under the functor $\Theta$. 

\begin{notation}
Fix a small collection of simplicial sets $\calK$. We let $\SSet(\calK)$ denote the unit object of the monoidal $\infty$-category $\Cat_{\infty}(\calK)$: it can be described concretely as the smallest full subcategory of $\SSet$ which contains the final object $\Delta^0$ and is closed under $\calK$-filtered colimits (Remark \toposref{poweryoga}). Since the formation of Cartesian products in $\SSet$ preserves small colimits in each variable, the full subcategory $\SSet(\calK) \subseteq \SSet$ is stable under finite products. We may therefore regard $\SSet(\calK)$ as equipped with the Cartesian monoidal structure, which endows it with the structure of an algebra object of $\Cat_{\infty}(\calK)$. We let
$\BigM$ denote the object of $\CatMod(\calK)$ given by the left action of
$\SSet(\calK)$ on itself. 
\end{notation}

According to Corollary \monoidref{puterry}, we can identify
$\BigM$ with $\Theta( \SSet(\calK)^{\times}, {\bf 1})$, where ${\bf 1}$ denotes the unit object
$\Delta^0 \in \SSet(\calK)$, regarded as an algebra object of $\SSet(\calK)$. 

\begin{lemma}\label{stuke}
Let $\calK$ be a small collection of simplicial sets. 
Then the pair $( \SSet(\calK)^{\times}, {\bf 1})$ is an initial object of $\CatAlg(\calK)$. 
\end{lemma}

\begin{proof}
Let $\phi: \CatAlg(\calK) \rightarrow \Cat_{\infty}^{\Mon}(\calK)$ denote the forgetful functor.
Then $\phi( \SSet(\calK)^{\times}, {\bf 1})$ is an initial object of $\Cat_{\infty}^{\Mon}(\calK) \simeq
\Alg( \Cat_{\infty}(\calK))$ (Proposition \monoidref{gurgle}). It will therefore suffice to show
that $( \SSet(\calK)^{\times}, {\bf 1})$ is a $\phi$-initial object of $\CatAlg(\calK)$ (Proposition \toposref{basrel}). Since $\phi$ is a coCartesian fibration (Proposition \ref{komb}), this is equivalent to the requirement that for every $\phi$-coCartesian morphism $\alpha: ( \SSet(\calK)^{\times}, {\bf 1}) \rightarrow ( \calC^{\otimes}, A)$, the object $A$ is initial in the fiber $\phi^{-1} \{ \calC^{\otimes} \} \simeq \Alg(\calC)$ (Proposition \toposref{relcolfibtest}). This follows immediately from Proposition \monoidref{gurgle}.
\end{proof}

It follows from Lemma \ref{stuke} that the forgetful functor
$\theta: \CatAlg(\calK)_{ ( \SSet(\calK)^{\times}, {\bf 1})/ } \rightarrow \CatAlg(\calK)$ is a trivial Kan fibration.
We let $\Theta_{\ast}$ denote the composition
$$ \CatAlg(\calK) \simeq \CatAlg(\calK)_{(\SSet(\calK)^{\times}, {\bf 1})} \stackrel{\Theta}{\rightarrow}
\CatMod(\calK)_{\BigM/},$$
where the first map is given by a section of $\theta$. 

\begin{remark}
An object of $\CatMod(\calK)_{\BigM/}$ is given by a morphism
$( \SSet(\calK)^{\times}, \SSet(\calK) ) \rightarrow (\calC^{\otimes}, \calM)$ in $\CatMod(\calK)$, given
by a monoidal functor $\SSet(\calK)^{\times} \rightarrow \calC^{\otimes}$ which preserves $\calK$-indexed colimits (which is unique up to a contractible space of choices by Proposition \monoidref{gurgle}) together with a functor $f: \SSet(\calK) \rightarrow \calM$ which preserves
$\calK$-indexed colimits. In view of Remark \toposref{poweryoga}, such a functor is determined
uniquely up to equivalence by the object $f( \Delta^0) \in \calM$. Consequently, we can informally regard $\CatMod(\calK)_{\BigM/}$ as an $\infty$-category whose objects are triples
$(\calC^{\otimes}, \calM, M)$, where $(\calC^{\otimes}, \calM) \in \CatMod$ and $M \in \calM$ is an object.
\end{remark}

\begin{theorem}\label{postcur}
Let $\calK$ be a small collection of simplicial sets which contains $\Nerve(\cDelta)^{op}$. Then the functor
$\Theta_{\ast}: \CatAlg(\calK) \rightarrow \CatMod(\calK)_{\BigM/}$ is
fully faithful.
\end{theorem}

\begin{lemma}\label{snoka}
Let $\calK$ be a small collection of simplicial sets which contains $\Nerve( \cDelta)^{op}$, 
let $( \calC^{\otimes}, \calM, M)$ be an object of $\CatMod(\calK)_{ \BigM/}$, and suppose
that there exists an algebra object $E \in \Alg(\calC)$ such that $M$ can be promoted to an object
$\overline{M} \in \Mod^{L}_{E}(\calM)$ where the action $E \otimes M \rightarrow M$ exhibits
$E$ as a morphism object $\Mor_{\calM}(M,M)$ (see 
morphism object $\Mor_{\calM}(M,M)$ (Definition \monoidref{supiner}). Then $E$ represents
the right fibration
$$ \CatAlg(\calK) \times_{ \CatMod(\calK)_{\BigM/} } ( \CatMod(\calK)_{\BigM/})_{/ (\calC^{\otimes}, \calM, M)}.$$
\end{lemma}

\begin{proof}
Theorem \ref{bix2} implies the existence of a functor
$\phi: \Mod^{R}_{E}(\calC)^{\otimes} \rightarrow \calM^{\otimes}$ of $\infty$-categories
left-tensored over $\calC^{\otimes}$ together with an identification $\alpha: M \simeq \phi(E)$ of left $E$-modules. The pair $(\phi, \alpha)$ determines an object
$\eta \in \CatAlg(\calK) \times_{ \CatMod(\calK)_{\BigM/} } ( \CatMod(\calK)_{\BigM/})_{/ (\calC^{\otimes}, \calM, M)}$ lying over $E$. We claim that this object is final. To prove this,
consider an arbitrary object $(\calD^{\otimes}, A) \in \CatAlg(\calK)$; we wish to show that
the map
$$ \bHom_{ \CatAlg(\calK)}( (\calD^{\otimes},A), (\calC^{\otimes},E) )
\rightarrow \bHom_{ \CatMod(\calK)_{\BigM/}}( \Theta_{\ast}( \calD^{\otimes}, A),
( \calC^{\otimes}, \calM, M) )$$
is a homotopy equivalence. It will suffice to prove the result after passing to the homotopy
fiber over a point of $\bHom_{ \Cat_{\infty}^{\Mon}(\calK)}( \calD^{\otimes}, \calC^{\otimes})$, corresponding to a monoidal functor $F$. Using Propositions \ref{komb} and \toposref{compspaces} and replacing $\calD^{\otimes}$ by $\calC^{\otimes}$ (and $A$ by $FA \in \Alg(\calC)$), we
are reduced to proving that the diagram
$$ \xymatrix{ \bHom_{ \Alg(\calC)}( A, E) \ar[r] \ar[d] & \bHom_{ \Mod_{ \calC}( \Cat_{\infty}(\calK))}( \Mod^{R}_{A}(\calC), \calM) \ar[d] \\
\{ M \} \ar[r] & \calM }$$
is a pullback square. Theorem \ref{bix2} allows us to identify the upper right corner of this diagram with the $\infty$-category $\Mod^{L}_{A}(\calM)$, and the desired result follows from Corollary \monoidref{usebilly}.
\end{proof}

\begin{proof}[Proof of Theorem \ref{postcur}]
Fix objects $(\calC^{\otimes}, A)$, $(\calD^{\otimes},B) \in \CatAlg(\calK)$.
We wish to show that the canonical map
$\theta:  \bHom_{ \CatAlg(\calK)}( (\calC^{\otimes},A), (\calD^{\otimes},B) ) \rightarrow \bHom_{\CatMod(\calK)_{\BigM/}}( \Theta_{\ast}( \calC^{\otimes},A),
\Theta_{\ast}( \calD^{\otimes}, B) )$ is a homotopy equivalence. Let $M \in \Mod^{R}_{B}(\calD)$
denote the right $B$-module given by the action of $B$ on itself.
In view of Lemma \ref{snoka}, it will suffice to show that the canonical map $m: B \otimes M \rightarrow M$ exhibits $B$ as a morphism object $\Mor_{ \Mod^{R}_{B}(\calD)}(M,M)$.
In other words, we must show that for every object $D \in \calD$, the multiplication map $m$ induces a homotopy equivalence
$\bHom_{ \calD}( D, B) \rightarrow \bHom_{ \Mod^{R}_{B}(\calD)}( D \otimes M, M)$. This follows from Proposition \monoidref{pretara}.
\end{proof}

Lemma \ref{snoka} can also be used to show that the fully faithful embedding $\Theta_{\ast}$
admits a right adjoint, provided that we can guarantee the existence of endomorphism objects
$\Mor_{\calM}(M,M)$. For this, it is convenient to work a setting where we require all $\infty$-categories to be presentable. For this, we need to introduce a bit of terminology.

\begin{notation}
Let $\widehat{\Cat}_{\infty}$ denote the $\infty$-category of (not necessarily small) $\infty$-categories, which contains $\Cat_{\infty}$ as a full subcategory. Similarly, we define $\infty$-categories
$\hatCatAlg \supset \CatAlg$ and $\hatCatMod \supset \CatMod$ by allowing
monoidal $\infty$-categories and left-tensored $\infty$-categories which are not small. Let $\calK$ denote the collection of {\em all} small simplicial sets, and let $\hatCatAlg(\calK)$ and
$\hatCatMod(\calK)$ be defined as in \S \ref{conl}. Construction \ref{spaas} generalizes immediately to give a functor $\widehat{ \Theta}: \hatCatAlg(\calK) \rightarrow
\hatCatMod(\calK)$.
We let
$\LCatAlg$ denote the full subcategory of $\hatCatAlg(\calK)$ spanned by those pairs
$( \calC^{\otimes}, A)$ where the $\infty$-category $\calC$ is presentable, and 
$\LCatMod$ the full subcategory of $\hatCatMod(\calK)$ spanned by those pairs
$( \calC^{\otimes}, \calM)$ where $\calC$ and $\calM$ are both presentable. 
It follows from Corollary \monoidref{underwhear} that the functor
$\widehat{\Theta}: \hatCatAlg(\calK) \rightarrow \hatCatMod(\calK)$ restricts
to a functor $\LCatAlg \rightarrow \LCatMod$, which we will also denote by $\widehat{\Theta}$.
Similarly, if we let $\BigM$ denote the object $\Theta( \SSet^{\times}, \Delta^0) \simeq ( \SSet^{\times}, \SSet) \in \LCatMod$, then we have a functor
$\widehat{\Theta}_{\ast}: \LCatAlg \rightarrow \LCatModM$. 
\end{notation}

\begin{theorem}\label{curly}
The functor $\widehat{\Theta}_{\ast}: \LCatAlg \rightarrow \LCatModM$ is fully faithful and admits a right adjoint.
\end{theorem}

\begin{proof}
The first assertion follows by applying Theorem \ref{postcur} in a larger universe.
For the second, it will suffice to show that for every object $X = ( \calC^{\otimes}, \calM, M)
\in \LCatModM$, the right fibration $\LCatAlg \times_{ \LCatModM} ( \LCatModM)_{ /X }$ is representable (Proposition \toposref{adjfuncbaby}).
In view of Lemma \ref{stuke}, it will suffice to show that there exists an algebra object
$E \in \Alg(\calC)$ such that $M$ can be promoted to a module $\overline{M} \in \Mod^{L}_{E}(\calC)$
such that the action $E \otimes M \rightarrow M$ exhibits $E$ as a morphism object
$\Mor_{\calM}(M,M)$. According to Propositions \monoidref{ugher} and \monoidref{poofer}, this is equivalent to requiring the existence of an algebra object $E \in \Alg( \calC[M])$ such that
the underlying object in $\calC[M]$ is final. According to Corollary \monoidref{firebaugh}, it will suffice
to show that $\calC[M]$ has a final object (which then admits an essentially unique algebra structure):
that is, it will suffice to show that there exists a morphism object $\Mor_{\calM}(M,M)$. This follows
from Proposition \monoidref{enterich}.
\end{proof}

\begin{remark}
Informally, the right adjoint to $\widehat{\Theta}_{\ast}$ carries an object
$(\calC^{\otimes}, \calM, M) \in \LCatModM$ to the pair $(\calC^{\otimes}, E) \in \LCatAlg$,
where $E \in \Alg(\calC)$ is the algebra of endomorphisms of the object $M \in \calM$.
\end{remark}

We next investigate the behavior of the functor $\Theta$ with respect to tensor products of
$\infty$-categories. 

\begin{notation}
The $\infty$-category $\Cat_{\infty}^{\Mon}$ admits finite products, and can therefore be regarded
as endowed with Cartesian symmetric monoidal structure. Let $\calK$ be a small collection of simplicial sets. We define a subcategory $\Cat_{\infty}^{\Mon}(\calK)^{\otimes} \subseteq \Cat_{\infty}^{\Mon, \times}$ as follows:
\begin{itemize}
\item[$(1)$] Let $C$ be an object of $\Cat_{\infty}^{\Mon, \times}$, given by a finite
sequence of monoidal $\infty$-categories $( \calC^{\otimes}_1, \ldots, \calC^{\otimes}_n)$. Then
$C \in \Cat_{\infty}^{\Mon}(\calK)^{\otimes}$ if and only if each of the underlying $\infty$-categories
$\calC_{i}$ admits $\calK$-indexed colimits, and the tensor product functors $\calC_i \times \calC_i \rightarrow \calC_i$ preserve $\calK$-indexed colimits separately in each variable.
\item[$(2)$] Let $F: ( \calC_1^{\otimes}, \ldots, \calC_m^{\otimes}) \rightarrow
( \calD_{1}^{\otimes}, \ldots, \calD_n^{\otimes})$ be a morphism in 
$\Cat_{\infty}^{\Mon, \times}$ covering a map $\alpha: \seg{m} \rightarrow \seg{n}$ in
$\FinSeg$, where the objects $( \calC_1^{\otimes}, \ldots, \calC_m^{\otimes})$ and
$( \calD_{1}^{\otimes}, \ldots, \calD_n^{\otimes})$ belong to $\Cat_{\infty}^{\Mon}(\calK)^{\otimes}$.
Then $F$ belongs to $\Cat_{\infty}^{\Mon}(\calK)^{\otimes}$ if and only if
the induced functor $\prod_{ \alpha(i) = j } \calC_i \rightarrow \calD_j$ preserves $\calK$-indexed colimits separately in each variable, for $1 \leq j \leq n$.
\end{itemize}

We let $\CatMod(\calK)^{\otimes}$ denote the subcategory of
$\CatMod^{\times}$ described as follows:
\begin{itemize}
\item[$(1')$] Let $C$ be an object of $\CatMod^{\times}$, corresponding to a finite
sequence $( ( \calC^{\otimes}_1, \calM_1), \ldots, (\calC^{\otimes}_n, \calM_n ) )$. 
Then $C \in \CatMod(\calK)^{\otimes}$ if and only if each $\calC_i$ and each $\calM_i$ admit
$\calK$-indexed colimits, and the tensor product functors 
$$ \calC_i \times \calC_i \rightarrow \calC_i \quad \quad \calC_i \times \calM_i \rightarrow \calM_i$$
preserves $\calK$-indexed colimits separately in each variable.
\item[$(2')$] 
Let $F: ( (\calC_1^{\otimes}, \calM_1), \ldots, (\calC_m^{\otimes},\calM_m)) \rightarrow
( (\calD_{1}^{\otimes}, \calN_1), \ldots, (\calD_n^{\otimes}, \calN_n))$ be a morphism in 
$\Cat_{\infty}^{\Mon, \times}$ covering a map $\alpha: \seg{m} \rightarrow \seg{n}$ in
$\FinSeg$, where the objects 
$( (\calC_1^{\otimes}, \calM_1), \ldots, (\calC_m^{\otimes},\calM_m))$ and
$( (\calD_{1}^{\otimes}, \calN_1), \ldots, (\calD_n^{\otimes}, \calN_n))$ belong to
$\CatMod(\calK)^{\otimes}$. Then $F$ belongs to 
$\CatMod(\calK)^{\otimes}$ if and only if
the induced functors 
$$\prod_{ \alpha(i) = j } \calC_i \rightarrow \calD_j \quad \quad \prod_{ \alpha(i) = j} \calM_i \rightarrow \calN_j$$ 
preserves $\calK$-indexed colimits separately in each variable, for $1 \leq j \leq n$.
\end{itemize}

We let $\CatAlg(\calK)^{\otimes}$ denote the fiber product
$(\CatAlg)^{\times} \times_{ \Cat_{\infty}^{\Mon, \times} } \Cat_{\infty}^{\Mon}(\calK)^{\otimes}$.
\end{notation}

\begin{remark}\label{upsy}
Assume that $\calK$ consists entirely of {\em sifted} simplicial sets (this is satisfied, for example, if
$\calK = \{ \Nerve(\cDelta)^{op} \}$. Then we can identify $\CatAlg(\calK)^{\otimes}$, 
$\CatMod(\calK)^{\otimes}$, and $\Cat_{\infty}^{\Mon}(\calK)^{\otimes}$ with subcategories
of $\CatAlg(\calK)^{\times}$, $\CatMod(\calK)^{\times}$, and $\Cat_{\infty}^{\Mon}(\calK)^{\times}$, respectively.
\end{remark}

\begin{theorem}\label{saylime}
Let $\calK$ be a small collection of simplicial sets. Then:
\begin{itemize}
\item[$(1)$] The map $\Cat_{\infty}^{\Mon}(\calK)^{\otimes}
\rightarrow \Nerve(\FinSeg)$ determines a symmetric monoidal structure on
$\Cat_{\infty}^{\Mon}$, and the maps
$\CatAlg(\calK)^{\otimes} \rightarrow \Cat_{\infty}^{\Mon}(\calK)^{\otimes} \leftarrow
\CatMod(\calK)^{\otimes}$ are coCartesian fibrations of symmetric monoidal $\infty$-categories.

\item[$(2)$] The functor $\Theta: \CatAlg( \{ \Nerve(\cDelta)^{op} \} ) \rightarrow \CatMod( \{ \Nerve(\cDelta)^{op} \})$ preserves products, and
therefore induces a symmetric monoidal functor $\Theta^{\times}: \CatAlg( \{ \Nerve(\cDelta)^{op} \})^{\times} \rightarrow \CatMod( \{ \Nerve(\cDelta)^{op} \})^{\times}$.

\item[$(3)$] Assume that $\Nerve( \cDelta)^{op} \in \calK$. Then the functor $\Theta^{\times}$ of $(2)$ restricts to a functor  $\Theta^{\otimes}: \CatAlg(\calK)^{\otimes} \rightarrow \CatMod(\calK)^{\otimes}$
(see Remark \ref{upsy}).

\item[$(4)$] The functor $\Theta^{\otimes}$ is symmetric monoidal.
\end{itemize}
\end{theorem}

\begin{proof}
We first prove $(1)$. Recall that $\Cat_{\infty}^{\Mon}$ can be identified with the $\infty$-category 
$\Alg_{\Ass}(\Cat_{\infty})$ of associative algebra objects of $\Cat_{\infty}$
(Remark \monoidref{otherlander} and Proposition \symmetricref{algass}). Here we regard
$\Cat_{\infty}$ as endowed with the Cartesian symmetric monoidal structure. The
$\infty$-category $\Alg_{\Ass}(\Cat_{\infty})$ inherits a symmetric monoidal structure from
that of $\Cat_{\infty}$ (see Example \symmetricref{slabber}), which is also Cartesian; we therefore have
an induced identification $\Alg_{\Ass}(\Cat_{\infty})^{\otimes} \simeq \Cat_{\infty}^{\Mon, \times}$. 
Under this equivalence, the subcategory $\Cat_{\infty}^{\Mon}(\calK)^{\otimes}$ corresponds to
the subcategory $\Alg_{\Ass}(\Cat_{\infty}(\calK))^{\otimes}$, which is a symmetric monoidal $\infty$-category (Example \symmetricref{slabber} again). This completes the proof that $\Cat_{\infty}^{\Mon}(\calK)^{\otimes}$ is a symmetric monoidal $\infty$-category. A similar argument
(using Theorem \symmetricref{podus}) shows that $\CatMod(\calK)^{\otimes} \rightarrow
\Cat_{\infty}^{\Mon}(\calK)^{\otimes}$ is a coCartesian fibration of symmetric monoidal $\infty$-categories. Finally, we observe that the functor $\CatAlg(\calK)^{\otimes} \rightarrow
\Cat_{\infty}^{\Mon}(\calK)^{\otimes}$ is a pullback of $(\CatAlg)^{\times} \rightarrow
\Cat_{\infty}^{\Mon, \times}$, which is easily seen to be a coCartesian fibration of $\infty$-operads.

Assertion $(2)$ is obvious, and assertion $(3)$ follows from Corollary \monoidref{gloop}.
We will prove $(4)$. It is easy to see that $\Theta^{\otimes}$ is a map of $\infty$-operads,
and it follows from Corollary \monoidref{puterry} that $\Theta^{\otimes}$ preserves
unit objects. Consequently, it will suffice to show that for every pair of objects
$( \calC^{\otimes}, A), (\calD^{\otimes}, B) \in \CatAlg(\calK)$, the induced map
$$ \Theta( \calC^{\otimes}, A) \otimes \Theta( \calD^{\otimes}, B) \rightarrow
\Theta( ( \calC^{\otimes}, A) \otimes (\calD^{\otimes},B) )$$
is an equivalence in $\CatMod(\calK)$. In other words, we wish to show that
$\Theta$ induces an equivalence of $\infty$-categories
$$ \theta: \Mod^{R}_{A}(\calC) \otimes \Mod^{R}_B(\calD) \rightarrow \Mod^{R}_{ A \otimes B}( \calC \otimes \calD)$$
(here the tensor products are taken in $\Cat_{\infty}(\calK)$). We have a homotopy commutative
diagram of $\infty$-categories
$$ \xymatrix{ \Mod^{R}_{A}(\calC) \otimes \Mod^{R}_{B}(\calD) \ar[rr]^{\theta} \ar[dr]^{G} & & \Mod^{R}_{A \otimes B}( \calC \otimes \calD) \ar[dl] \\
& \calC \otimes \calD. & }$$
To prove that $\theta$ is a categorical equivalence, it will suffice to show that this diagram
satisfies the hypotheses of Corollary \monoidref{littlerbeck}:
\begin{itemize}
\item[$(a)$] The $\infty$-categories $\Mod^{R}_{A}(\calC) \otimes \Mod^{R}_{B}(\calD)$ and
$\Mod^{R}_{ A \otimes B}( \calC \otimes \calD)$ admit geometric realizations of simplicial objects.
This follows from our assumption that $\Nerve(\cDelta)^{op} \in \calK$.
\item[$(b)$] The functors $G$ and $G'$ admit left adjoints, which we will denote by $F$ and $F'$.
The existence of $F'$ is guaranteed by Proposition \monoidref{pretara}, and is given informally
by the formula $X \mapsto X \otimes (A \otimes B)$. Similarly, Proposition \monoidref{pretara} guarantees that the forgetful functors $\Mod^{R}_{A}(\calC) \rightarrow \calC$ and $\Mod^{R}_{B}(\calD) \rightarrow \calD$ admit left adjoints, given by tensoring on the right with $A$ and $B$, respectively.
The tensor product of these left adjoints is a left adjoint to $G$.
\item[$(c)$] The functor $G'$ is conservative and preserves geometric realizations of simplicial objects.
The first assertion follows from Corollary \monoidref{thetacart} and the second from Corollary \monoidref{gloop}.
\item[$(d)$] The functor $G$ is conservative and preserves geometric realizations of simplicial objects.
The second assertion is obvious: $G$ is a tensor product of the forgetful functors
$\Mod^{R}_{A}(\calC) \rightarrow \calC$ and $\Mod^{R}_{B}(\calD) \rightarrow \calD$, each of which preserves geometric realizations (and can therefore be interpreted as a morphism in $\Cat_{\infty}(\calK)$) by
Corollary \monoidref{gloop}. To prove that $G$ is conservative, we factor $G$ as a composition
$$ \Mod^{R}_{A}(\calC) \otimes \Mod^{R}_{B}(\calD) \stackrel{G_0}{\rightarrow}
\calC \otimes \Mod^{R}_{B}(\calD) \stackrel{G_1}{\rightarrow} \calC \otimes \calD.$$
Using Proposition \monoidref{ussr}, we can identify $G_1$ with the forgetful functor
$$ (\calC \otimes \calD) \otimes_{\calD} \Mod^{R}_{B}(\calD) \rightarrow (\calC \otimes \calD) \otimes_{\calD} \calD \simeq \calC \otimes \calD.$$
Theorem \ref{bix1} allows us to identify the left hand side with the $\infty$-category
$\Mod^{R}_{B}(\calC \otimes \calD)$. Under this identification, $G_1$ corresponds to the forgetful functor $\Mod^{R}_{B}(\calC \otimes \calD) \rightarrow \calC \otimes \calD$, which is 
conservative by Corollary \monoidref{thetacart}. A similar arguments shows that $G_0$ is conservative, so that $G \simeq G_1 \circ G_0$ is conservative as required.

\item[$(e)$] The canonical natural transformation $G' \circ F' \rightarrow G \circ F$ is an
equivalence of functors from $\calC \otimes \calD$ to itself. This is clear from the descriptions
of $F$ and $F'$ given above: both compositions are given by right multiplication by the object
$A \otimes B \in \calC \otimes \calD$.
\end{itemize}
\end{proof}

\begin{remark}\label{kulpus}
Fix a small collection of simplicial sets $\calK$ which contains $\Nerve( \cDelta)^{op}$. 
Let $\calC^{\otimes}$ be a symmetric monoidal $\infty$-category. Assume that 
$\calC$ admits $\calK$-indexed colimits and that the tensor product
$\calC \times \calC \rightarrow \calC$ preserves $\calK$-indexed colimits separately in each variable.
Then we $\calC^{\otimes}$ as a commutative algebra object in the (symmetric monoidal) $\infty$-category $\Cat_{\infty}^{\Mon}(\calK)$. The fiber products
$$ \Nerve( \FinSeg) \times_{ \Cat_{\infty}^{\Mon}(\calK)} \CatAlg(\calC)^{\otimes} 
\quad \quad \Nerve( \FinSeg) \times_{ \Cat_{\infty}^{\Mon}(\calK)^{\otimes} } \CatMod(\calC)^{\otimes}$$
can be identified with the symmetric monoidal $\infty$-categories $\Alg(\calC)^{\otimes}$ and
$\Mod^{L}_{\calC}( \Cat_{\infty}(\calK))^{\otimes}$, respectively. It follows from Theorem
\ref{saylime} and Proposition \ref{komb} that $\Theta$ determines a symmetric monoidal functor
$\Theta^{\otimes}_{\calC}: \Alg(\calC)^{\otimes} \rightarrow \Mod^{L}_{\calC}( \Cat_{\infty}(\calK))^{\otimes}$.
\end{remark}

\begin{corollary}\label{sabis}
Let $\calK$ and $\calC^{\otimes}$ be as in Remark \ref{kulpus}, and let
$\calO^{\otimes}$ be a unital $\infty$-operad. Then the functor $\Theta^{\otimes}_{\calC}$
induces a fully faithful functor
$$ \theta: \Alg_{ \calO}( \Alg(\calC) ) \rightarrow \Alg_{\calO}( \Mod^{L}_{\calC}( \Cat_{\infty}(\calK) ) ).$$
\end{corollary}

\begin{proof}
Let $\Alg(\calC)_{\ast}^{\otimes}$ and $\Mod^{L}_{\calC}( \Cat_{\infty}(\calK))^{\otimes}_{\ast}$ be
unitalizations of $\Alg(\calC)^{\otimes}$ and $\Mod^{L}_{\calC}( \Cat_{\infty}(\calK) )^{\otimes}$,
respectively (see \S \ref{unitos}; note that $\Alg(\calC)^{\otimes}$ is already a unital $\infty$-operad, so that
$\Alg(\calC)^{\otimes}_{\ast} \simeq \Alg(\calC)^{\otimes}$). The functor $\Theta^{\otimes}_{\calC}$ induces a symmetric monoidal
functor $\Alg(\calC)_{\ast}^{\otimes} \rightarrow \Mod^{L}_{\calC}( \Cat_{\infty}(\calK))^{\otimes}_{\ast},$
and Theorem \ref{postcur} guarantees that this functor is fully faithful. We have a commutative
diagram
$$ \xymatrix{ \Alg_{\calO}( \Alg(\calC)_{\ast}) \ar[r]^-{\theta_{\ast}} \ar[d] & \Alg_{\calO}( \Mod^{L}_{\calC}( \Cat_{\infty}(\calK))_{\ast}) \ar[d] \\
\Alg_{\calO}( \Alg(\calC) ) \ar[r]^-{\theta} & \Alg_{\calO}( \Mod^{L}_{\calC}( \Cat_{\infty}(\calK) )) }$$
where $\theta_{\ast}$ is fully faithful. Since $\calO^{\otimes}$ is unital, the vertical maps are
categorical equivalences, so that $\theta$ is fully faithful as well.
\end{proof}

\begin{corollary}\label{eliac}
Let $\calK$ and $\calC^{\otimes}$ be as in Remark \ref{kulpus}. Then for $n \geq 1$, we have a fully faithful functor
$\Alg_{\OpE{n}}( \calC) \rightarrow \Alg_{ \OpE{n-1}}( \Mod^{L}_{\calC}( \Cat_{\infty}(\calK)))$.
\end{corollary}

\begin{proof}
Combine Corollary \ref{sabis}, Proposition \symmetricref{algass}, and Theorem \ref{slide}. 
\end{proof}

Corollary \ref{eliac} furnishes a convenient way of understanding the notion of an $\OpE{n}$-algebra:
giving an $\OpE{n}$-algebra object $A \in \Alg_{\OpE{n}}(\calC)$ is equivalent to giving the underlying
associative algebra object $A_0 \in \Alg(\calC)$, together with an $\OpE{n-1}$-structure on the
$\infty$-category $\Mod^{R}_{A_0}(\calC)$ of right $A_0$-modules (with unit object given by the module $A_0$ itself).

\subsection{The $\infty$-Operad $\LMod$}\label{clapser}

Proposition \symmetricref{algass} implies that giving a monoidal $\infty$-category
$\calC^{\otimes} \rightarrow \Nerve( \cDelta)^{op}$ is equivalent to giving a coCartesian fibration of $\infty$-operads $q: {\calC'}^{\otimes} \rightarrow \Ass$. Our goal in this section is to extend this equivalence to the case of modules. More precisely, we will show that giving an $\infty$-category $\calM^{\otimes} \rightarrow \calC^{\otimes}$ left-tensored over $\calC^{\otimes}$ is equivalent to extending $q$ to a coCartesian fibration of $\infty$-operads ${\calM'}^{\otimes} \rightarrow \LMod$, for a suitably defined $\infty$-operad $\LMod$
(see Remark \ref{puuff}). Moreover, giving a left module object of $\calM^{\otimes}$ is equivalent to giving
an $\LMod$-algebra object of ${\calM'}^{\otimes}$ (Proposition \ref{hunters}).

We begin by defining the $\infty$-operad $\LMod$. 

\begin{definition}
We define a category $\CatLMod$ as follows:
\begin{itemize}
\item[$(1)$] The objects of $\CatLMod$ are pairs $( \seg{n},S)$, where $\seg{n}$ is an object of
$\FinSeg$ and $S$ is a subset of $\nostar{n}$.
\item[$(2)$] A morphism from $( \seg{n}, S)$ to $( \seg{n'}, S')$ in $\CatLMod$ consists of a
morphism $\alpha: \seg{n} \rightarrow \seg{n'}$ in $\CatAss$ satisfying the following conditions:
\begin{itemize}
\item[$(i)$] If $s' \in \nostar{n} - S'$, then $\alpha^{-1} \{s'\}$ does not intersect $S$.
\item[$(ii)$] If $s' \in S'$, then $\alpha^{-1} \{s'\}$ contains exactly one element of $S$, and that
element is final with respect to the linear ordering of $\alpha^{-1} \{s'\}$.
\end{itemize}
\end{itemize}
We let $\LMod$ denote the nerve of the category $\CatLMod$.
\end{definition}

There is an evident forgetful functor $\CatLMod \rightarrow \FinSeg$. It follows from 
Example \symmetricref{xe3} that the induced map $\LMod = \Nerve( \CatLMod) \rightarrow \Nerve(\FinSeg)$ exhibits $\LMod$ as an $\infty$-operad. We will show in a moment that this
is the $\infty$-operad which governs left modules over associative algebras: that is, giving
an $\LMod$-algebra is equivalent to giving a pair $(A,M)$, where $A$ is an associative
algebra and $M$ is a left $A$-module (Proposition \ref{hunters}).

\begin{construction}\label{sass}
We will identify $\CatAss$ with the full subcategory of $\LMod$ spanned by those pairs
$(\seg{n}, S)$ where $S = \emptyset$.

Let $\phi: \cDelta^{op} \rightarrow \CatAss \subset \LMod$ be the functor described in
Construction \symmetricref{urpas}. We define another functor $\phi^{L}: \cDelta^{op} \rightarrow \LMod$
as follows:
\begin{itemize}
\item[$(1)$] For each $n \geq 0$, we have $\phi^{L}( [n] ) = (\seg{n+1}, \{n+1\})$.
\item[$(2)$] Given a morphism $\alpha: [n] \rightarrow [m]$ in $\cDelta$, the associated morphism
$\phi(\alpha): \seg{m+1} \rightarrow \seg{n+1}$ is given by the formula
$$\phi(\alpha)(i) = \begin{cases} j & \text{if } (\exists j) [\alpha(j-1) < i \leq \alpha(j)] \\
n+1 & \text{ if } (\forall j) [\alpha(j) < i] \\
\ast & \text{ otherwise. } \end{cases}$$
\end{itemize}

There is an evident natural transformation of functors $\phi^{L} \rightarrow \phi$, which determines
a map of simplicial sets $\Nerve( \cDelta)^{op} \rightarrow \Fun( \Delta^1, \LMod)$.
\end{construction}

\begin{remark}
The inclusion $\CatAss \subseteq \LMod$ induces a fully faithful embedding
of $\infty$-operads $\Ass \rightarrow \LMod$. In particular, every coCartesian fibration
of $\infty$-operads $q: \calC^{\otimes} \rightarrow \LMod$ determines an $\Ass$-monoidal
$\infty$-category $\calC^{\otimes} \times_{\LMod} \Ass$, and therefore 
a monoidal $\infty$-category $\calC^{\otimes} \times_{ \LMod} \Nerve(\cDelta)^{op}$
(see Construction \symmetricref{urpas}). We will refer to the fiber product
$\calC^{\otimes} \times_{ \LMod} \Nerve(\cDelta)^{op}$ as the {\it underlying monoidal
$\infty$-category} of $\calC^{\otimes}$, and denote it by $\calA[q]^{\otimes}$.
\end{remark}

\begin{construction}\label{stonger}
Suppose that $q: \calC^{\otimes} \rightarrow \LMod$ is a coCartesian fibration of $\infty$-operads.
We let $\calM[q]^{\otimes}$ denote the fiber product $\Fun( \Delta^1, \calC^{\otimes})
\times_{ \Fun( \Delta^1, \LMod)} \Nerve( \cDelta)^{op}$, where
$\Nerve(\cDelta)^{op}$ maps to $\Fun( \Delta^1, \LMod)$ via the functor described in
Construction \ref{sass}.
 \end{construction}

\begin{remark}\label{jumpter}
Let $q: \calC^{\otimes} \rightarrow \LMod$ be as in Construction \ref{stonger}. 
Evaluation at $\{0\} \subseteq \Delta^1$ induces a trivial Kan fibration $\calM[q]^{\otimes}
\rightarrow \calC^{\otimes} \times_{\LMod} \Nerve(\cDelta)^{op}$, where
the map $\Nerve(\cDelta)^{op} \rightarrow \LMod$ is determined by the functor
$\phi^{L}: \cDelta^{op} \rightarrow \CatLMod$ of Construction \ref{sass}.
\end{remark}

\begin{proposition}\label{slapfight}
Let $q: \calC^{\otimes} \rightarrow \LMod$ be a coCartesian fibration of $\infty$-operads.
Then evaluation at $\{1\} \subseteq \Delta^1$ induces a map
$\theta: \calM[q]^{\otimes} \rightarrow \calA[q]^{\otimes}$, which exhibits
$\calM[q]^{\otimes}$ as an $\infty$-category left-tensored over the monoidal $\infty$-category $\calA[q]^{\otimes}$ (see Definition \monoidref{ulult}).
\end{proposition}

\begin{proof}
It is easy to see that $\theta$ is a categorical fibration.
Let $p: \calA[q]^{\otimes} \rightarrow \Nerve(\cDelta)^{op}$ and
$p': \calM[q]^{\otimes} \rightarrow \Nerve(\cDelta)^{op}$ be the projection maps.
It follows from Proposition \toposref{doog} that the map $p'$ is a coCartesian fibration, and
that $\theta$ carries $p'$-coCartesian morphisms to $p$-coCartesian morphisms.
To complete the proof, it suffices to show that for every $n \geq 0$, the inclusion
$\{n\} \subseteq [n]$ induces an equivalence of $\infty$-categories
$\calM[q]^{\otimes}_{[n]} \rightarrow \calA[q]^{\otimes}_{[n]} \times \calM[q]^{\otimes}_{[0]}$.
Using Remark \ref{jumpter}, we see that this is equivalent to the requirement that
the natural map 
$$\calC^{\otimes}_{ ( \seg{n+1}, \{n+1\})} \rightarrow \calC^{\otimes}_{ ( \seg{n}, \emptyset)}
\times \calC^{\otimes}_{ (\seg{1}, \{1\}) }$$
is an equivalence. This follows from the observation that the maps
$$ ( \seg{n}, \emptyset) \leftarrow ( \seg{n+1}, \{n+1\}) \rightarrow ( \seg{1}, \{1\} )$$
determine a splitting of the object $( \seg{n+1}, \{n+1\}) \in \LMod$ (Definition \symmetricref{capper}).
\end{proof}

The main result of this section is the following analogue of Proposition \symmetricref{algass}:

\begin{proposition}\label{hunters}
Let $q: \calC^{\otimes} \rightarrow \LMod$ be a coCartesian fibration of $\infty$-operads.
Then Construction \ref{stonger} induces an equivalence of $\infty$-categories
$\theta: \Alg_{ \LMod}( \calC) \rightarrow \Mod^{L}( \calM[q] ).$
\end{proposition}

\begin{proof}
We define a category $\calI$ as follows:
\begin{itemize}
\item[$(a)$] An object of $\calI$ is either an object of $\cDelta^{op} \times [1]$ or an object of $\CatLMod$.
\item[$(b)$] Morphisms in $\calI$ are give by the formulas
$$\Hom_{\calI}( ([m],i), ([n],j) ) = \Hom_{ \cDelta^{op} \times [1]}( ([m],i), ([n],j))$$
$$ \Hom_{\calI}( (\seg{m},S), (\seg{n},T) ) = \Hom_{\CatLMod}( (\seg{m},S), (\seg{n},T) )$$
$$ \Hom_{ \calI}( (\seg{m},S), ([n],0) ) =
\Hom_{ \CatLMod}( (\seg{m},S), \phi^{L} [n] )$$
$$\Hom_{\calI}( (\seg{m},S), ([n],1))
= \Hom_{\CatLMod}( ( \seg{m},S), \phi [n]) $$
$$\Hom_{ \calI}( ([n],i), (\seg{m},T) ) = \emptyset.$$ 
\end{itemize}
where $\phi, \phi^{L}: \cDelta^{op} \rightarrow \CatLMod$ are defined as in Construction \ref{sass}.
There is a canonical retraction $r$ from $\calI$ onto $\CatLMod$, given on objects of
$\cDelta^{op} \times [1]$ by the formula $r( [n], 0) = \phi^{L}([n])$, $r([n],1) = \phi([n])$.

Let $\overline{\Alg}(\calC)$ denote the full subcategory of 
$\Fun_{ \LMod}( \Nerve(\calI), \calC^{\otimes})$ consisting of those functors $f: \Nerve(\calI) \rightarrow \calC^{\otimes}$ 
such that $q \circ f = r$ and the following additional conditions are satisfied:
\begin{itemize}
\item[$(i)$] For every object
$([n],i) \in \cDelta^{op} \times [1]$, the functor $f$ carries the canonical map
$r([n],i) \rightarrow ([n],i)$ in $\calI$ to an equivalence in $\calC^{\otimes}$.
\item[$(ii)$] The restriction $f | \Nerve(\cDelta)^{op} \times \Delta^1$ determines
an object of $\Mod^{L}( \calM[q])$. 
\item[$(iii)$] The restriction $f | \LMod$ is an $\LMod$-algebra object of $\calC$.
\end{itemize}

It is easy to see that condition $(ii)$ follows from $(i)$ and $(iii)$. Moreover, $(i)$ is equivalent to the requirement that $f$ is a $q$-left Kan extension of $f| \LMod$. Since every functor $f_0 \in \Fun_{ \LMod}( \LMod, \calC^{\otimes})$ admits a $q$-left Kan extension $f \in \Fun_{ \LMod}( \Nerve(\calI), \calC^{\otimes})$ (given, for example, by $f_0 \circ r$), Proposition \toposref{lklk} implies that the restriction map $p: \overline{\Alg}(\calC) \rightarrow \Alg_{\LMod}(\calC)$ is a trivial Kan fibration. The map $\theta$ is the composition of a section to $p$ (given by composition with $r$) 
with the restriction map $p': \overline{\Alg}(\calC) \rightarrow \Mod^{L}(\calM[q])$. It will therefore suffice to show that $p'$ is a trivial fibration. In view of Proposition \toposref{lklk}, this will follow from the following pair of assertions:

\begin{itemize}
\item[$(a)$] Every $f_0 \in \Mod^{L}( \calM[q])$ admits a $q$-right Kan extension
$f \in \Fun_{ \LMod}( \Nerve(\calI), \calC^{\otimes})$.
\item[$(b)$] Given $f \in \Fun_{ \LMod}( \Nerve(\calI), \calC^{\otimes})$
such that $f_0 = f | \Nerve(\cDelta)^{op} \times \Delta^1$ belongs to $\Mod^{L}( \calM[q])$, $f$ is a $q$-right Kan extension of $f_0$ if and only if $f$ satisfies conditions $(i)$ and $(ii)$ above.
\end{itemize}

To prove $(a)$, we fix an object $(\seg{n},S) \in \LMod$. Let $\calJ$ denote the category
$( \cDelta^{op} \times [1]) \times_{ \calI} ( \calI)_{(\seg{n},S)/ },$
and let $g$ denote the composition
$ \Nerve(\calJ) \rightarrow \Nerve(\cDelta)^{op} \times \Delta^1 \stackrel{f_0}{\rightarrow} \calC^{\otimes}.$
According to Lemma \toposref{kan2}, it will suffice to show that $g$ admits a $q$-limit in $\calC^{\otimes}$ (compatible with the evident map $\Nerve(\calJ)^{\triangleleft}) \rightarrow \LMod$). 
The objects of $\calJ$ can be identified with
morphisms $\alpha: (\seg{n},S) \rightarrow r([m],i)$ in $\LMod$. Let $\calJ_0 \subseteq \calJ$
denote the full subcategory spanned by those objects for which $\alpha$ is inert.
The inclusion $\calJ_0 \subseteq \calJ$ has a right adjoint, so that $\Nerve(\calJ_0)^{op} \rightarrow \Nerve(\calJ)^{op}$ is cofinal. Consequently, it will suffice to show that $g_0 = g | \Nerve(\calJ_0)$ admits a $q$-limit in $\calC^{\otimes}$ (compatible with the evident map
$\Nerve(\calJ_0)^{\triangleleft} \rightarrow \LMod$).

Let $\calJ_1$ denote the full subcategory of $\calJ_0$ spanned by the morphism which are either of the form $\Colp{j}: (\seg{n},S) \rightarrow r( [1],0)$ where $j \in S$ or $\Colp{j}: ( \seg{n}, S) \rightarrow r([1],1)$ when $j \notin S$. Using our assumption that
$f_0 \in \Mod^{L}( \calM[q])$, we deduce that $g_0$ is a $q$-right Kan extension of
$g_1 = g_0 | \Nerve(\calJ_1)$. In view of Lemma \toposref{kan0}, it will suffice to show that
the map $g_1$ admits a $q$-limit in $\calC$ (compatible with the map
$\Nerve(\calJ_1)^{\triangleleft} \rightarrow \LMod$). This follows immediately from our assumption that
$q$ is a fibration of $\infty$-operads, thereby proving $(a)$. Moreover, the proof shows that
$f$ is a $q$-right Kan extension of $f_0$ if and only if it satisfies the following
condition:
\begin{itemize}
\item[$(i')$] For every object $(\seg{n}, S) \in \LMod$ as above and every morphism
$\alpha: (\seg{n},S) \rightarrow ([1],i)$ in $\calI$ belonging to $\calJ_1$, the morphism $f(\alpha)$ is inert
in $\calC^{\otimes}$. 
\end{itemize}

To prove $(b)$, it will suffice to show that if 
$f \in \Fun_{ \LMod}( \Nerve(\calI), \calC^{\otimes})$ satisfies condition $(ii)$, then
it satisfies conditions $(i)$ and $(iii)$ if and only if it satisfies condition $(i')$. We first prove the
``only if'' direction. Assume that $f \in \overline{\Alg}(\calC)$, and let $\alpha: (\seg{n}, S) \rightarrow ([1],i)$ be as in $(i')$. Then $\alpha$ factors as a composition
$$ (\seg{n}, S) \stackrel{ \alpha'}{\rightarrow} (\seg{1}, T) \stackrel{\alpha''}{\rightarrow} ([1],i),$$
where $\alpha'$ is inert (so that $f(\alpha')$ is inert by virtue of $(iii)$) and $f(\alpha'')$ is an equivalence
by virtue of $(i)$.

Suppose now that $f$ satisfies $(i')$ and $(ii)$. We first show that $f$ satisfies $(i)$.
Fix an object $([n],i)$ in $\Nerve(\cDelta)^{op} \times \Delta^1$; we wish to show that
$f$ carries the canonical map $\alpha: r([n],i) \rightarrow ([n],i)$ to an equivalence in
$\calC^{\otimes}$. For $1 \leq j \leq n$, let $\beta_j: ([n],i) \rightarrow ([1], i')$ be the map
carrying $[1]$ to the interval $\{ j-1, j \} \subseteq [n]$, where $$i' = \begin{cases} 0 & \text{ if } i=0, j=n \\
1 & \text{otherwise.} \end{cases}$$
Condition $(ii)$ guarantees that each of the maps $f( \beta_{j})$ is inert in $\calC^{\otimes}$.
Since $\calC^{\otimes}$ is an $\infty$-operad, we deduce that $f( \alpha)$ is an equivalence
if and only if each of the composite maps $f( \beta_{j} \circ \alpha)$ is inert. We now conclude by invoking $(i')$. 

It remains to show that $f$ satisfies $(iii)$.
By virtue of Remark \symmetricref{casper}, it will suffice to show that $f(\alpha)$ is inert
whenever $\alpha$ is an inert morphism of $\LMod$ of the form $(\seg{n}, S) \rightarrow ( \seg{1}, T)$.
Let $\beta: (\seg{1}, T) \rightarrow ( [1], i)$ be a morphism in $\calI$ such that
$r(\beta)$ is an equivalence. Then $f(\beta)$ is an equivalence by virtue of $(i)$, so it will suffice
to show that $f( \beta \circ \alpha)$ is an inert morphism in $\calC^{\otimes}$: this follows from $(i')$.
\end{proof}

\begin{definition}
We let $\AAA$ denote the object $(\seg{1}, \emptyset) \in \LMod$ and
$\MM$ the object $( \seg{1}, \{1\} ) \in \LMod$. If
$q: \calC^{\otimes} \rightarrow \LMod$ is a fibration of $\infty$-operads, then we let
$\calC_{\AAA}$ and $\calC_{\MM}$ denote the fiber products $\calC^{\otimes} \times_{ \LMod} \{ \AAA \}$
and $\calC^{\otimes} \times_{ \LMod} \{ \MM \}$, respectively.
\end{definition}

\begin{remark}\label{puuff}
Suppose that $q: \calC^{\otimes} \rightarrow \LMod$ is a coCartesian fibration of
$\infty$-operads. According to Corollary \symmetricref{co1}, we can identify
$q$ with an $\LMod$-algebra object of $\Cat_{\infty}$ (endowed with the Cartesian symmetric monoidal structure), is in turn equivalent to giving a left module object of $\Cat_{\infty}$ (Proposition \ref{hunters}).
In other words, we can think of $q$ as providing an $\infty$-category with an associative
tensor product (namely, the fiber $\calC_{\AAA}$) together with a left action of this $\infty$-category
on another $\infty$-category (the fiber $\calC_{\MM}$). This is an equivalent way of describing the data encoded in a left-tensored $\infty$-category (see Proposition \ref{slapfight}).
\end{remark}

\subsection{Centralizers and Deligne's Conjecture}\label{delconj}

Let $A$ be an associative algebra over a field $k$ with multiplication $m$. Then the cyclic bar complex
$$ \ldots \rightarrow A \otimes_k A \otimes_k A \stackrel{ m \otimes \id - \id \otimes m}{\longrightarrow} A \otimes_k A$$ provides a resolution of $A$ by free $(A \otimes_k A^{op})$-modules; we will denote this resolution by $P_{\bigdot}$. The cochain complex $\HC^{\ast}(A) = \Hom_{A \otimes_{k} A^{op}}( P_{\bigdot}, A)$
is called the {\it Hochschild cochain complex} of the algebra $A$. The cohomologies of
the Hochschild cochain complex (which are the $\Ext$-groups $\Ext^{i}_{ A \otimes_{k} A^{op}}( A, A)$)
are called the {\it Hochschild cohomology} groups of the algebra $A$. A famous conjecture of Deligne 
asserts that the Hochschild cochain complex $\HC^{\ast}(A)$ carries an action of the little $2$-disks operad:
in other words, that we can regard $\HC^{\ast}(A)$ as an $\OpE{2}$-algebra object in the
$\infty$-category of chain complexes over $k$. This conjecture has subsequently been proven by many authors in many different ways (see, for example, \cite{mcsmith}, \cite{kontsoib}, and \cite{tamarkin}). 
In this section we will outline a proof of Deligne's conjecture, using the ideas presented in this paper. Our basic strategy can be outlined as follows:

\begin{itemize}
\item[$(1)$] Let $\calC^{\otimes}$ be an arbitrary presentable symmetric monoidal $\infty$-category, and let
$f: A \rightarrow B$ be an $\OpE{k}$-algebra object of $\calC$. We will prove that there exists another
$\OpE{k}$-algebra $\Centt{\OpE{k}}{f}$ of $\calC$, which is universal with respect to the existence of a
commutative diagram
$$ \xymatrix{ & \Centt{\OpE{k}}{f} \otimes A \ar[dr] & \\
A \ar[ur]^{u} \ar[rr]^{f} & & B }$$
in the $\infty$-category $\Alg_{ \OpE{k}}(\calC)$, where $u$ is induced by the unit map
${\bf 1} \rightarrow \Centt{\OpE{k}}{f}$ (Theorem \ref{centex}).

\item[$(2)$] In the case where $A = B$ and $f$ is the identity map, we will denote
$\Centt{\OpE{k}}{f}$ by $\Centt{\OpE{k}}{A}$. We will see that $\Centt{\OpE{k}}{A}$ has the structure of
an $\OpE{k+1}$-algebra object of $\calC$. (More generally, the centralizer construction is functorial
in the sense that there are canonical maps of $\OpE{k}$-algebras $\Centt{\OpE{k}}{f} \otimes \Centt{\OpE{k}}{g} \rightarrow
\Centt{\OpE{k}}{f \circ g}$; in the special case $f = g = \id_{A}$ this gives rise to an associative algebra
structure on the $\OpE{k}$-algebra $\Centt{\OpE{k}}{A}$, which promotes $\Centt{\OpE{k}}{A}$ to an $\OpE{k+1}$-algebra
by Theorem \ref{slide}.) 

\item[$(3)$] We will show that the image of $\Centt{\OpE{k}}{f}$ in $\calC$ can be identified with a classifying
object for morphisms from $A$ to $B$ in the $\infty$-category $\Mod^{\OpE{k}}_{A}(\calC)$ (Theorem \ref{serpses}). In the special case where $k=1$ and $f = \id_{A}$, we can identify this with a classifying object for endomorphisms of $A$ as an $A$-bimodule (see \S \symmetricref{bimid}), recovering the usual definition of Hochschild cohomology.
\end{itemize}

We begin with a very general discussion of centralizers. 

\begin{definition}\label{saigd}
Let $q: \calC^{\otimes} \rightarrow \LMod$ be a coCartesian fibration of
$\infty$-operads. Let ${\bf 1}$ denote the unit object of $\calC_{\AAA}$, and suppose we are given a morphism $f: {\bf 1} \otimes M \rightarrow N$ in $\calC_{\MM}$. A
{\it centralizer of $f$} is final object of the $\infty$-category
$$ (\calC_{\AAA})_{{\bf 1}/} \times_{ (\calC_{\MM})_{{\bf 1 \otimes M}/} } (\calC_{\MM})_{{\bf 1} \otimes M/ \, /N}.$$
We will denote such an object, if it exists, by $\Cent{f}$. We will refer to $\Cent{f}$ as the {\it centralizer} of the morphism $f$. 
\end{definition}

\begin{remark}
We will generally abuse notation by identifying $\Cent{f}$ with its image in the 
$\infty$-category $\calC$. By construction, this object is equipped with a map
$\Cent{f} \otimes A \rightarrow B$ which fits into a commutative diagram
$$ \xymatrix{ & \Cent{f} \otimes A \ar[dr] &  \\
{\bf 1} \otimes A \ar[ur] \ar[rr]^{f} & & B. }$$
\end{remark}

\begin{remark}\label{sly}
In the situation of Definition \ref{saigd}, choose an algebra object $A \in \Alg_{\LMod}(\calC)$
such that $A( \MM) = M$ and $A | \Ass$ is a trivial algebra. Then we can identify
centralizers for a morphism $f: M \rightarrow N$ with morphism objects
$\Mor_{ \calM_{M/} }(M,N)$ computed in the $\LMod$-monoidal $\infty$-category
$\calC^{\otimes}_{A_{\LMod}/}$. 
\end{remark}

\begin{definition}\label{hulkus}
Let $q: \calC^{\otimes} \rightarrow \LMod$ be a fibration of $\infty$-operads, and let
$M \in \calC_{\MM}$ be an object. A {\it center} of $M$ is a final object of the
fiber product $\Alg_{\LMod}( \calC) \times_{ \calC_{\MM}} \{M\}$. If such an object exists,
we will denote it by $\Cent{M}$.
\end{definition}

\begin{remark}
In the situation of Definition \ref{hulkus}, we will often abuse notation by identifying
the center $\Cent{M}$ with its image in the $\infty$-category $\Alg_{\Ass}(\calC)$ of associative algebra objects of the $\Ass$-monoidal $\infty$-category $\calC_{\AAA}$.
\end{remark}

Our first goal is to show that, as our notation suggests, the theory of centers is closely related to the theory of centralizers. Namely, we have the following:

\begin{proposition}\label{staffstuff}
Let $q: \calC^{\otimes} \rightarrow \LMod$ be a coCartesian fibration of $\infty$-operads,
and let $M \in \calC_{\MM}$. Suppose that there exists a centralizer of the 
canonical equivalence $e: {\bf 1} \otimes M \rightarrow M$. Then there exists a center of $M$. Moreover, an object $A \in \Alg_{\LMod}(\calC)$ with $A(\MM)= M$ is a center of $M$ if and only if the unit map of $A$ and the diagram 
$$ \xymatrix{ & A(\AAA) \otimes M \ar[dr] & \\
{\bf 1} \otimes M \ar[ur] \ar[rr]^{e} & & M }$$
exhibits $A(\AAA)$ as a centralizer of $e$.
\end{proposition}

The proof will require a few preliminaries.

\begin{proposition}\label{jukk}
Let $q: \calC^{\otimes} \rightarrow \LMod$ be a fibration of $\infty$-operads.
Assume that $\calC^{\otimes}$ admits a $\AAA$-unit object (Definition \symmetricref{cuspe}), and let
$\theta: \Alg_{\LMod}(\calC) \rightarrow \calC_{\MM}$ be the functor given by evaluation
at $\MM$. Then:
\begin{itemize}
\item[$(1)$] The functor $\theta$ admits a left adjoint $L$.
\item[$(2)$] Let $\Alg'_{\LMod}(\calC) \subseteq \Alg_{\LMod}(\calC)$ be the essential image of $L$. Then $\theta$ induces a trivial Kan fibration $\Alg'_{\LMod}(\calC) \rightarrow \calC_{\MM}$.
In particular, $L$ is fully faithful.
\item[$(3)$] An $\LMod$-algebra object $A \in \Alg_{\LMod}(\calC)$ belongs to 
$\Alg'_{\LMod}(\calC)$ if and only if $A| \Ass$ is a trivial $\Ass$-algebra (see \S \symmetricref{unitr}).
\end{itemize}
\end{proposition}

\begin{proof}
Let us identify the $\infty$-operad $\Triv$ with the full subcategory of
$\LMod$ spanned by those objects having the form $( \seg{n}, \nostar{n})$. The functor
$\theta$ factors as a composition
$$ \Alg_{\LMod}(\calC) \stackrel{\theta'}{\rightarrow} \Alg_{\Triv}(\calC) \stackrel{\theta''}{\rightarrow} \calC_{\MM},$$
where $\theta''$ is a trivial Kan fibration (Example \symmetricref{algtriv}). To prove
$(1)$, it will suffice to show that $\theta'$ admits a left adjoint. We claim that this left
adjoint exists, and is given by operadic $q$-left Kan extension along the inclusion
$\Triv \subset \LMod$. According to Proposition \symmetricref{wax1}, it suffices to verify
that for each $\overline{M} \in \Alg_{\Triv}(\calC)$ and
every object of the form $X = (\seg{1}, S) \in \LMod$, the map
$\Triv \times_{ \LMod } (\LMod)^{\acti}_{/X} \rightarrow \calC^{\otimes}$
can be extended to an operadic $q$-colimit diagram (lying over the natural map
$( \Triv \times_{ \LMod} (\LMod)^{\acti})_{/X})^{\triangleright} \rightarrow \LMod$).
If $S = \{1\}$, then $X \in \Triv$ and the result is obvious. If $S = \emptyset$, then
the desired result follows from our assumption that $\calC^{\otimes}$ admits an $\AAA$-unit.
This proves $(1)$. Moreover, we see that if $A \in \Alg_{\LMod}(\calC)$, then a map
$f: \overline{M} \rightarrow A| \Triv$ exhibits $A$ as a free $\LMod$-algebra generated
by $\overline{M}$ if and only $A | \Ass$ is a trivial algebra and $f$ is an equivalence.
It follows that the unit map $\id \rightarrow \theta \circ L$ is an equivalence, so that $L$ is a
fully faithful embedding whose essential image $\Alg'_{\LMod}(\calC)$ is as described
in assertion $(3)$. To complete the proof, we observe that $\theta | \Alg'_{\LMod}(\calC)$ is an equivalence of $\infty$-categories and also a categorical fibration, and therefore a trivial Kan fibration.
\end{proof}

\begin{proposition}\label{reductocenter}
Let $q: \calC^{\otimes} \rightarrow \LMod$ be a fibration of $\infty$-operads.
Assume that $\calC^{\otimes}$ admits a $\AAA$-unit object (see Definition \symmetricref{cuspe}), let 
$A \in \Alg_{\LMod}(\calC)$ be an algebra object such that $A | \Ass$ is trivial,
and let $M = A( \MM)$.
Let ${\calC'}^{\otimes} = \calC^{\otimes}_{ A_{\LMod}/}$ be defined as in
Notation \ref{kurplex}, and let $M'$ be the object $\id_{M} \in
\calC'_{\MM} \simeq (\calC_{\MM})_{M/}$. Then the forgetful functor
$$ \theta: \Alg_{\LMod}(\calC') \times_{ \calC'_{\MM}} \{M'\} \rightarrow \Alg_{\LMod}(\calC)
\times_{ \calC_{\MM}} \{M\}$$
is a trivial Kan fibration. In particular, $M$ has a center in
$\calC^{\otimes}$ if and only if $M'$ has a center in ${\calC'}^{\otimes}$.
\end{proposition}

\begin{proof}
Note that $A: \LMod \rightarrow \calC^{\otimes}$ is a coCartesian section of
$q$, so that ${\calC'}^{\otimes} \rightarrow \LMod$ is a coCartesian fibration
of $\infty$-operads (Theorem \ref{eli}).
Since $\theta$ is evidently a categorical fibration, it will suffice to show that $\theta$
is a trivial Kan fibration. To this end, we let $\calA$ denote the full subcategory
of $\Fun( \Delta^1, \Alg_{\LMod}(\calC))$ spanned by those morphisms $A' \rightarrow A$
which exhibit $A'$ as an $\Alg'_{\LMod}(\calC)$-colocalization of $A$, where
$\Alg'_{\LMod}(\calC)$ is the full subcategory of $\Alg_{\LMod}(\calC)$ described in Proposition
\ref{jukk} (in other words, a morphism $A' \rightarrow A$ belongs to $\calA$ if and only if
$A' | \Ass$ is a trivial algebra and the map $A'(\MM) \rightarrow A(\MM)$ is an equivalence).
Evaluation at $\{1\} \subseteq \Delta^1$ and $\MM \in \LMod$ induces a functor
$e: \calA \rightarrow \calC_{\MM}$. The map $\theta$ factors as a composition
$$ \Alg_{\LMod}( \calC') \times_{ \calC'_{\MM}} \{M'\} \stackrel{\theta'}{\rightarrow}
\calA \times_{ \calC_{\MM} } \{M\} \stackrel{\theta''}{\rightarrow} \Alg_{\LMod}(\calC) \times_{ \calC_{\MM}} \{M\},$$
where $\theta''$ is a pullback of the trivial Kan fibration $\calA \rightarrow \Alg_{\LMod}(\calC)$
given by evaluation at $\{1\}$. We conclude by observing that $\theta'$ is also an equivalence of $\infty$-categories.
\end{proof}

\begin{proposition}\label{easycenter}
Let $q: \calC^{\otimes} \rightarrow \LMod$ be a coCartesian fibration of $\infty$-operads
and let $M \in \calC_{\MM}$ be such that a morphism object $\Mor_{ \calC_{\MM}}(M,M)$ exists
in $\calC_{\AAA}$ (Definition \monoidref{supiner}). Then there exists a center $\Cent{M}$;
moreover, an algebra object $A \in \Alg_{\LMod}(\calC)$ with $A( \MM) = M$ is
a center of $M$ if and only if the canonical map $A(\AAA) \otimes M \rightarrow M$ exhibits
$A(\AAA)$ as a morphism object $\Mor_{ \calC_{\MM}}(M,M)$.
\end{proposition}

\begin{proof}
Combine Proposition \monoidref{ugher}, Corollary \monoidref{slimycomp}, and
Proposition \ref{hunters}.
\end{proof}

\begin{proof}[Proof of Proposition \ref{staffstuff}]
Combine Proposition \ref{reductocenter}, Proposition \ref{easycenter}, and
Remark \ref{sly}.
\end{proof}

We are primarily interested in studying centralizers in the setting of $\calO^{\otimes}$-algebra objects
of a symmetric monoidal $\infty$-category $\calC^{\otimes}$. To emphasize the role of $\calO^{\otimes}$, it is convenient to introduce a special notation for this situation:

\begin{definition}\label{coopin}
Let $\calO^{\otimes}$ and $\calD^{\otimes}$ be $\infty$-operads, and let $p: \calO^{\otimes} \times \LMod \rightarrow \calD^{\otimes}$ be a bifunctor of $\infty$-operads. Suppose that $q: \calC^{\otimes} \rightarrow \calD^{\otimes}$ is a coCartesian fibration of $\infty$-operads. Then we have an induced coCartesian
fibration of $\infty$-operads $q': \Alg_{\calO}(\calC)^{\otimes} \rightarrow \LMod$ (see Notation \symmetricref{scambag}). If
$f: A \rightarrow B$ is a morphism in $\Alg_{\calO}(\calC)_{\MM}$, then we let
$\Centt{\calO}{f}$ denote the centralizer of $f$ (as a morphism in $\Alg_{\calO}(\calC)^{\otimes}$), provided that this centralizer exists. If $A \in \Alg_{\calO}(\calC)_{\MM}$, we let $\Centt{\calO}{A}$ denote the center of $A$.
\end{definition}

\begin{remark}
The primary case of interest to us is that in which $\calD^{\otimes} = \Nerve( \FinSeg)$, so that
$\calC^{\otimes}$ can be regarded as a symmetric monoidal $\infty$-category and the map
$p: \calO^{\otimes} \times \LMod \rightarrow \calD^{\otimes}$ is uniquely determined.
In this case, we will denote $\Alg_{\calO}(\calC)_{\MM} \simeq \Alg_{\calO}(\calC)_{\AAA}$ simply by
$\Alg_{\calO}(\calC)$. If $A \in \Alg_{\calO}(\calC)$, we can identify the center $\Centt{\calO}{A}$ (if it exists)
with an associative algebra object of the symmetric monoidal $\infty$-category $\Alg_{\calO}(\calC)$.
If $\calO^{\otimes}$ is a little cubes operad, then Theorem \ref{slide} and Example \ref{sulta} provide equivalences of $\infty$-categories
$$ \Alg_{ \OpE{k+1}}(\calC) \rightarrow \Alg_{ \OpE{1}}( \Alg_{ \OpE{k}}(\calC) ) \leftarrow
\Alg_{ \Ass}( \Alg_{ \OpE{k}}(\calC) ),$$
so we can identify $\Centt{ \OpE{k}}{A}$ with an $\OpE{k+1}$-algebra object of $\calC$. 
\end{remark}

In the situation of Definition \ref{coopin}, it is generally not possible to prove the existence of centralizers by direct application of Proposition \ref{easycenter}: the tensor product of $\calO$-algebra objects usually does not commute with colimits in either variable, so there generally does not exist a morphism object
$\Mor_{ \Alg_{\calO}(\calC)_{\MM} }( A, B)$ for a pair of algebras $A, B \in \Alg_{\calO}(\calC)_{\MM}$. 
Nevertheless, if $\calO$ is coherent, then we will show that the centralizer $\Centt{\calO}{f}$ of a morphism $f: A \rightarrow B$ exists under very general conditions:

\begin{theorem}\label{centex}
Let $\calO^{\otimes}$ be a coherent $\infty$-operad, let $p: \calO^{\otimes} \times \LMod \rightarrow \calD^{\otimes}$ be a bifunctor of $\infty$-operads, and let $q: \calC^{\otimes} \rightarrow \calD^{\otimes}$ exhibit $\calC^{\otimes}$ as a presentable $\calD$-monoidal $\infty$-category. Then,
for every morphism $f: A \rightarrow B$ in $\Alg_{ \calO}( \calC)_{\MM}$, there exists a centralizer
$\Centt{\calO}{f}$. 
\end{theorem}

\begin{corollary}\label{centcor}
Let $k \geq 0$, and let $\calC^{\otimes}$ be a symmetric monoidal $\infty$-category. Assume
that $\calC$ is presentable and that the tensor product $\otimes: \calC \times \calC \rightarrow \calC$ preserves
small colimits separately in each variable. Then:
\begin{itemize}
\item[$(1)$] For every morphism $f: A \rightarrow B$ in $\Alg_{ \OpE{k}}(\calC)$, there exists a centralizer
$\Centt{ \OpE{k}}{f} \in \Alg_{ \OpE{k}}(\calC)$.
\item[$(2)$] For every object $A \in \Alg_{ \OpE{k}}(\calC)$, there exists a center
$$ \Centt{\OpE{k}}{A} \in \Alg_{ \Ass}( \Alg_{ \OpE{k} }(\calC) \simeq \Alg_{ \OpE{k+1} }(\calC).$$
\end{itemize}
\end{corollary}

\begin{proof}
Combine Theorems \ref{centex} and \ref{cubecoh}. 
\end{proof}

\begin{example}[Koszul Duality]\label{kd}
Let $\calC^{\otimes}$ be a symmetric monoidal $\infty$-category. Assume that $\calC$ is presentable and that the tensor product on $\calC$ preserves small colimits separately in each variable. 
An {\it augmented $\OpE{k}$-algebra} in $\calC$ is a map of $\OpE{k}$-algebras
$\epsilon: A \rightarrow {\bf 1}$, where ${\bf 1}$ is a trivial $\OpE{k}$-algebra object of $\calC$.
It follows from Theorem \ref{serpses} that $\epsilon$ admits a centralizer $\Cent{\epsilon}$.
We will refer to this centralizer as the {\it Koszul dual} of $\epsilon$, and denote it by
$A^{\vee}$. By definition, $A^{\vee}$ is universal among $\OpE{k}$-algebras such that
$A \otimes A^{\vee}$ is equipped with an augmentation $\epsilon': A \otimes A^{\vee} \rightarrow
{\bf 1}$ compatible with the augmentation $\epsilon$ on $A$. In this case,
the composite map
$$ \epsilon^{\vee}: A^{\vee} \simeq {\bf 1} \otimes A^{\vee} \rightarrow A \otimes A^{\vee} \stackrel{\epsilon'}{\rightarrow} {\bf 1}$$
is an augmentation on $A^{\vee}$, so we can again regard $A^{\vee}$ as an augmented $\OpE{k}$-algebra object of $\calC$. In many cases, the relationship between $A$ and $A^{\vee}$ is symmetric: that is, $\epsilon'$ also exhibits $A$ as a centralizer of $\epsilon^{\vee}$. We will discuss this construction in more detail elsewhere (see also Example \ref{emore}).
\end{example}

\begin{example}[Drinfeld Centers]
Let $\calC$ be an $\OpE{k}$-monoidal $\infty$-category. Using Example \symmetricref{interplace} and
Proposition \symmetricref{ungbatt}, we can view $\calC$ as an $\OpE{k}$-algebra object of
the $\infty$-category $\Cat_{\infty}$ (which we regard as endowed with the Cartesian symmetric monoidal structure). Corollary \ref{centcor} guarantees the existence of a center $\Centt{\OpE{k}}{\calC}$, which
we can view as an $\OpE{k+1}$-monoidal $\infty$-category. In the special case where $k=1$ and
$\calC$ is (the nerve of) an ordinary monoidal category, the center $\Centt{\OpE{1}}{\calC}$ is also equivalent to the nerve of an ordinary category $\calZ$. Example \ref{sweeta} guarantees that $\calZ$ admits the structure of a braided monoidal category. This braided monoidal category $\calZ$ is called the {\em Drinfeld center} of
the monoidal category underlying $\calC$ (see, for example, \cite{drinfeldcenter}). Consequently, we can view the construction $\calC \mapsto \Centt{ \OpE{k} }{\calC}$ as a higher-categorical generalization of the theory of the Drinfeld center.
\end{example}

Our goal for the remainder of this section is to provide a proof of Theorem \ref{centex}. The idea is to change $\infty$-categories to maneuver into a situation where Proposition \ref{easycenter} can be applied. To carry out this strategy, we will need to introduce a bit of notation.

\begin{definition}\label{sably}
Let $\calO^{\otimes}$ be an $\infty$-operad and $S$ an $\infty$-category. A {\it coCartesian
$S$-family of $\calO$-operads} is a map $q: \calC^{\otimes} \rightarrow \calO^{\otimes} \times S$
with the following properties:
\begin{itemize}
\item[$(i)$] The map $q$ is a categorical fibration.
\item[$(ii)$] The underlying map $\calC^{\otimes} \rightarrow \Nerve(\FinSeg) \times S$ exhibits
$\calC^{\otimes}$ as an $S$-family of $\infty$-operads, in the sense of Definition \symmetricref{sike}.
\item[$(iii)$] For every object $C \in \calC^{\otimes}$ with $q(C) = (X,s) \in \calO^{\otimes} \times S$
and every morphism $f: s \rightarrow s'$ in $S$, there exists a $q$-coCartesian morphism
$C \rightarrow C'$ in $\calC^{\otimes}$ lifting the morphism $(\id_{X}, f)$.
\end{itemize}
\end{definition}

\begin{remark}
Let $q: \calC^{\otimes} \rightarrow \calO^{\otimes} \times S$ be a coCartesian $S$-family of $\calO$-operads. Condition $(iii)$ of Definition \ref{sably} guarantees that the underlying map
$\calC^{\otimes} \rightarrow S$ is a coCartesian fibration, classified by some map
$\chi: S \rightarrow \Cat_{\infty}$. The map $q$ itself determines a natural transformation from
$\chi$ to the constant functor $\chi_0$ taking the value $\calO^{\otimes}$, so that
$\chi$ determines a functor $\overline{\chi}: S \rightarrow \Cat_{\infty}^{/\calO^{\otimes}}$. 
This construction determines a bijective correspondence between equivalences classes of
$S$-families of $\calO$-operads and equivalence classes of functors from $S$ to the
$\infty$-category $\Cat_{\infty}^{\lax, \calO}$ of $\calO$-operads (see Remark \symmetricref{slame}),
which we can identify with a subcategory of $\Cat_{\infty}^{/ \calO^{\otimes}}$. 
\end{remark}

\begin{definition}\label{scoll}
Let $\calO^{\otimes}$ be an $\infty$-operad, $S$ an $\infty$-category, and
$q: \calC^{\otimes} \rightarrow \calO^{\otimes} \times S$ be a coCartesian $S$-family of
$\calO$-operads. We define a simplicial set $\Alg^{S}_{\calO}(\calC)$ equipped with a map
$\Alg^{S}_{\calO}(\calC) \rightarrow S$ so that the following universal property is satisfied:
for every map of simplicial sets $T \rightarrow S$, there is a canonical bijection of
$\Hom_{S}( T, \Alg^{S}_{\calO}(\calC))$ with the subset of
$\Hom_{\calO^{\otimes} \times S}( \calO^{\otimes} \times T, \calC)$ spanned by those
maps with the property that for each vertex $t \in T$, the induced map
$\calO^{\otimes} \rightarrow \calC^{\otimes}_{t}$ belongs to $\Alg_{\calO}(\calC_{t})$.
\end{definition}

\begin{remark}
If $q: \calC^{\otimes} \rightarrow \calO^{\otimes} \times S$ is as in Definition \ref{scoll}, then
the induced map $q': \Alg^{S}_{\calO}(\calC) \rightarrow S$ is a coCartesian fibration.
We will refer to a section $A: S \rightarrow \Alg^{S}_{\calO}(\calC)$ of $q'$ as an
{\it $S$-family of $\calO$-algebra objects of $\calC$}. We will say that an
$S$-family of $\calO$-algebra objects of $\calC$ is {\it coCartesian} if
$A$ carries each morphism in $S$ to a $q'$-coCartesian morphism in $\Alg^{S}_{\calO}(\calC)$.
\end{remark}

\begin{definition}\label{safe}
Let $\calO^{\otimes}$ be a coherent $\infty$-operad, let
$q: \calC^{\otimes} \rightarrow \calO^{\otimes} \times S$ be a coCartesian $S$-family of
$\calO$-operads, and let $\calC_0^{\otimes}$ denote the product
$\calO^{\otimes} \times S$. If $A$ is an $S$-family of algebra objects of $\calC$, we let
$\Mod^{\calO}_{A}(\calC)$ denote the fiber product
$$ \MMod^{\calO}(\calC)^{\otimes} \times_{ \PrAlg_{\calO}(\calC)} S,$$
where $\MMod^{\calO}(\calC)$ and $\PrAlg_{\calO}(\calC)$ are defined as in \S \symmetricref{moduldef}
and the map $S \rightarrow \PrAlg_{\calO}(\calC)$ is determined by $A$.
We let $\Mod^{\calO,S}_{A}(\calC)$ denote the fiber product
$\Mod^{\calO}_{A}(\calC) \times_{ \Mod^{\calO}_{qA}(\calC_0)} \calC_0$.
\end{definition}

\begin{remark}\label{scol}
Let $q: \calC^{\otimes} \rightarrow \calO^{\otimes} \times S$ be as in Definition \ref{safe} and $A$
is an $S$-family of $\infty$-operads. Then the $\infty$-category $\Mod^{\calO,S}_{A}(\calC)^{\otimes}$ is equipped with an evident forgetful functor $\Mod^{\calO,S}_{A}(\calC)^{\otimes} \rightarrow \calO^{\otimes} \times S$. For every object $s \in S$, the fiber $\Mod^{\calO,S}_{A}(\calC)^{\otimes}_{s}
= \Mod^{\calO}_{A}(\calC)^{\otimes} \times_{S} \{s\}$ is canonically isomorphic to the
$\infty$-operad $\Mod^{\calO}_{A_s}(\calC_s)^{\otimes}$ defined in \S \symmetricref{moduldef}.
\end{remark}

We will need the following technical result, whose proof will be given at the end of this section.

\begin{proposition}\label{squisher}
Let $\calO^{\otimes}$ be a coherent $\infty$-operad, $q: \calC^{\otimes} \rightarrow \calO^{\otimes} \times S$ a coCartesian $S$-family of $\calO$-operads, and $A \in \Alg^{S}_{\calO}(\calC)$
a coCartesian $S$-family of $\calO$-algebras. Then:
\begin{itemize}
\item[$(1)$] The forgetful functor
$q': \Mod^{\calO,S}_{A}(\calC)^{\otimes} \rightarrow \calO^{\otimes} \times S$ is a coCartesian $S$-family of $\infty$-operads.
\item[$(2)$] Let $\overline{f}$ be a morphism in $\Mod^{\calO,S}_{A}(\calC)^{\otimes}$ whose image in
$\calO^{\otimes}$ is degenerate. Then $\overline{f}$ is $q'$-coCartesian if and only if its image in
$\calC^{\otimes}$ is $q$-coCartesian.
\end{itemize}
\end{proposition}

\begin{remark}\label{hurley}
In the situation of Proposition \ref{squisher}, suppose that $\calO^{\otimes}$ is the $0$-cubes operad
$\OpE{0}$. Let $\calC$ denote the fiber product $\calC^{\otimes} \times_{ \calO^{\otimes}} \calO$.
Then the forgetful functor $\theta: \Mod^{\calO,S}_{A}(\calC) \rightarrow \calC$ is a trivial Kan fibration.
To prove this, it suffices to show that $\theta$ is a categorical equivalence (since it is evidently a categorical fibration). According to Corollary \toposref{usefir}, it suffices to show that
$\theta$ induces a categorical equivalence after passing to the fiber over each vertex of $S$, which follows from Proposition \symmetricref{hir}.
\end{remark}

Suppose now that $q: \calC^{\otimes} \rightarrow \calO^{\otimes} \times S$ is
a coCartesian $S$-family of $\calO$-operads and that $A$ is a coCartesian
$S$-family of $\calO$-algebra objects of $\calC$. Then $A$ determines an
$S$-family of $\calO$-algebra objects of $\Mod_{A}^{\calO,S}(\calC)$, which we will
denote also by $A$. Note that, for each $s \in S$, $A_{s} \in \Alg_{\calO}( \Mod_{A_s}^{\calO}(\calC_s))$
is a trivial algebra and therefore initial in $\Alg_{\calO}( \Mod_{A_s}^{\calO}(\calC_s))$ (Proposition
\symmetricref{gargle1}).
Let $\Alg_{\calO}^{S}(\calC)^{A_{S}/}$ be defined
as in \S \toposref{consweet} and let $\Alg_{\calO}^{S}( \Mod_{A}^{\calO,S}(\calC))^{A_S/}$ be defined similarly. We have a commutative diagram
$$ \xymatrix{ \Alg_{\calO}^{S}( \Mod_{A}^{\calO,S}(\calC))^{A_{S}/} \ar[rr]^{\theta} \ar[dr] & & \Alg_{\calO}^{S}(\calC)^{A_S/} \ar[dl] \\
& S. & }$$
The vertical maps are coCartesian fibrations and $\theta$ preserves coCartesian morphisms.
Using Corollary \symmetricref{skoke}, we deduce that $\theta$ induces a categorical equivalence after passing to the fiber over each object of $S$. Applying Corollary \toposref{usefir}, we deduce the following:

\begin{proposition}\label{canc}
Let $\calO^{\otimes}$ be a coherent $\infty$-operad, $q: \calC^{\otimes} \rightarrow \calO^{\otimes} \times S$ a coCartesian $S$-family of $\calO$-operads, and $A$ a coCartesian $S$-family of
$\calO$-algebra objects of $\calC$. Then the forgetful functor
$$\theta: \Alg_{\calO}^{S}( \Mod_{A}^{\calO,S}(\calC))^{A_S/} \rightarrow \Alg_{\calO}^{S}(\calC)^{A_S/}$$
is an equivalence of $\infty$-categories.
\end{proposition}

Proposition \ref{canc} provides a mechanism for reducing questions about centralizers of arbitrary
algebra morphisms $f: A \rightarrow B$ to the special case where $A$ is a trivial algebra.

\begin{remark}\label{hotmother}
Let $\calA^{\otimes} \rightarrow \LMod$ be a coCartesian fibration of $\infty$-operads, and let
${\bf 1}$ denote the unit object of the $\Ass$-monoidal $\infty$-category $\calA_{\AAA}$. 
Let $f: M_0 \rightarrow M$ be a morphism in $\calA_{\MM}$, and consider the fiber product
$\calX = \calA_{\AAA} \times_{ \calA_{\MM}} (\calA_{\MM})_{/M}$, where the map $\calA_{\AAA} \rightarrow \calA_{\MM}$ is given by tensor product with $M_0$. We will identify the tensor product
${\bf 1} \otimes M_0$ with $M_0$, so that the pair $( {\bf 1}, f: M_0 \rightarrow M)$ can be
identified with an object $X \in \calX$. The undercategory $\calX_{X/}$ can be identified with
the fiber product $(\calA_{\AAA})_{{\bf 1}/} \times_{ (\calA_{\MM})_{M_0/}} (\calC_{\MM})_{M_0/ \, /M}$. 
Using Proposition \toposref{needed17}, we deduce that the forgetful functor
$\calX_{X/} \rightarrow \calX$ induces an equivalence between the full subcategories spanned by the final objects of $\calX_{X/}$ and $\calX$. In other words:
\begin{itemize}
\item[$(i)$] A map $\epsilon: {\bf 1} \rightarrow Z$ in $\calA_{\AAA}$ together with a commutative diagram
$$ \xymatrix{ & Z \otimes M_0 \ar[dr]^{g} & \\
M_0 \ar[ur]^{\epsilon \otimes \id_{M_0}} \ar[rr]^{f} & & M}$$
in $\calA_{\MM}$ is a centralizer of $f$ if and only if the underlying morphism $g$ exhibits
$Z$ as a morphism object $\Mor_{ \calA_{\MM}}(M_0, M)$. 

\item[$(ii)$] For any object $Z \in \calA_{\AAA}$ and any morphism $Z \otimes M_0 \rightarrow M$ which
exhibits $Z$ as a morphism object $\Mor_{ \calA_{\MM}}(M_0, M)$, there exists a map
${\bf 1} \rightarrow Z$ and a commutative diagram
$$ \xymatrix{ & Z \otimes M_0 \ar[dr]^{g} & \\
M_0 \ar[ur]^{\epsilon \otimes \id_{M_0}} \ar[rr]^{f} & & M}$$
satisfying the conditions of $(i)$.
\end{itemize}
\end{remark}

\begin{proposition}\label{louda}
Let $q: \calC^{\otimes} \rightarrow \calO^{\otimes} \times \LMod$ be a coCartesian
$\LMod$-family of $\calO^{\otimes}$-operads, and assume that the induced map
$\Alg^{\LMod}_{\calO}(\calC) \rightarrow \LMod$ is a coCartesian fibration of $\infty$-operads
(this is automatic if, for example, the map $\calC^{\otimes} \rightarrow \LMod$ is a coCartesian fibration of $\infty$-operads). For every object $X \in \calO^{\otimes}$, we let
$\calC_{X,\AAA}$ denote the fiber of $q$ over the vertex $(\AAA,X)$, and define
$\calC_{X,\MM}$ similarly. Let $f: A_0 \rightarrow A$ be a morphism
in $\Alg_{\calO}^{\LMod}(\calC)_{\MM}$. Assume that:
\begin{itemize}
\item[$(i)$] The $\infty$-operad $\calO^{\otimes}$ is unital.
\item[$(ii)$] The algebra object $A_0$ is trivial (see \S \symmetricref{unitr}).
\item[$(iii)$] For every object $X \in \calO^{\otimes}$, there exists a morphism object
$\Mor_{ \calC_{X,\MM}}( A_0(X), A(X) ) \in \calC_{X,\AAA}$.
\end{itemize}
Then:
\begin{itemize}
\item[$(1)$] There exists a centralizer $\Cent{f} \in \Alg_{\calO}^{\LMod}(\calC)_{\AAA}$.
\item[$(2)$] Let $Z \in \Alg_{\calO}^{\LMod}(\calC)_{\AAA}$ be an algebra object. Then a commutative diagram
$$ \xymatrix{ & Z \otimes A_0 \ar[dr]^{g} & \\
A_0 \ar[ur] \ar[rr]^{f} & & A }$$
exhibits $Z$ as a centralizer of $f$ if and only if, for every object $X \in \calO$, the induced map
$g_{X}: \Cent{X} \otimes A_0(X) \rightarrow A(X)$ exhibits $\Cent{X}$ as a morphism object
$\Mor_{ \calC_{X,\MM}}( A_0(X), A(X))$.
\end{itemize}
\end{proposition}

\begin{proof}
Let ${\bf 1} \in \Alg^{\LMod}_{\calO}(\calC)_{\AAA}$ be a trivial algebra; we will abuse notation by
identifying the tensor product ${\bf 1} \otimes A_0$ with $A_0$. To prove $(1)$, we must show that the $\infty$-category
$$\calA = (\Alg^{\Mod}_{\calO}( \calC)_{\AAA})_{{\bf 1}/} \times_{
(\Alg_{\calO}^{\LMod}(\calC)_{\MM})_{ A_0/ } } (\Alg^{\LMod}_{\calO}(\calC)_{\MM})_{A_0 / \, /A}$$
has a final object. Let $\calC^{\otimes}_{\AAA}$ denote the fiber product
$\calC^{\otimes} \times_{ \LMod } \{ \AAA \}$,
define $\calC^{\otimes}_{\MM}$ similarly, and set
$$\calE^{\otimes} = (\calC^{\otimes}_{\AAA})_{ {\bf 1}_{\calO}/}
\times_{ (\calC^{\otimes}_{\MM})_{ {A_0}_{\calO}/ }} (\calC^{\otimes}_{\MM})
_{ {A_0}_{\calO}/ \, / A_{\calO}}$$
(see \S \ref{cluper} for an explanation of this notation). Using Theorem \ref{eli} (and assumption $(ii)$), we deduce that the evident forgetful functor
$\calE^{\otimes} \rightarrow \calO^{\otimes}$ is a coCartesian fibration of $\infty$-operads;
moreover, we have a canonical isomorphism
$\calA \simeq \Alg_{\calO}( \calE)$. For each object $X \in \calO^{\otimes}$, the $\infty$-category
$\calE_{X} = \calE^{\otimes} \times_{ \calO^{\otimes} } \{X \}$ is equivalent
to the fiber product
$$ (\calC_{X,\AAA})_{ {\bf 1}(X)/} \times_{ (\calC_{X,\MM})_{ A_0(X)/}}
(\calC_{X,\MM})_{ A_0(X)/ \, /A(X)},$$
which has a final object by virtue of assumption $(iii)$ and Remark \ref{hotmother}.
It follows that $\calA$ has a final object; moreover,
an object $A \in \calA \simeq \Alg_{\calO}(\calE)$ is final if and only if each
$A(X)$ is a final object of $\calE_{X}$. This proves $(1)$, and reduces assertion
$(2)$ to the contents of Remark \ref{hotmother}.
\end{proof}

We now apply Proposition \ref{louda} to the study of centralizers in general. Fix a coherent
$\infty$-operad $\calO^{\otimes}$, a bifunctor of $\infty$-operads
$\calO^{\otimes} \times \LMod \rightarrow \calD^{\otimes}$, and a coCartesian fibration
of $\infty$-operads $q: \calC^{\otimes} \rightarrow \calD^{\otimes}$. Let
$A \in \Alg_{\calO}(\calC)_{\MM}$, and let $\overline{A} \in \Alg_{\LMod}( \Alg_{\calO}(\calC))$
be an algebra such that $\overline{A}_{\MM} = A$ and $\overline{A}_{\AAA}$ is a trivial algebra.
We can regard $\overline{A}$ as a coCartesian $\LMod$-family of
$\calO$-algebra objects of $\calC^{\otimes} \times_{ \calD^{\otimes} } ( \calO^{\otimes} \times \LMod)$.
Let $\overline{\calC}^{\otimes} = \Mod^{\calO, \LMod}_{\overline{A}}(\calC)$ be the coCartesian
$S$-family of $\calO$-operads given by Proposition \ref{squisher}. Since $\overline{A}_{\AAA}$ is trivial, the forgetful functor $\overline{\calC}^{\otimes}_{\AAA} \rightarrow \calC^{\otimes}_{\AAA} =
\calC^{\otimes} \times_{ \calD^{\otimes}} ( \calO^{\otimes} \times \{ \AAA \} )$ is an equivalence
of $\calO$-operads, and induces an equivalence of $\infty$-categories
$\Alg_{\calO}( \overline{\calC})_{\AAA} \rightarrow \Alg_{\calO}(\calC)_{\AAA}$.
It follows from Proposition \ref{canc} that every morphism
$f: A \rightarrow B$ in $\Alg_{\calO}(\calC)_{\MM}$ is equivalent to $\theta( f' )$, where
$f': A \rightarrow B'$ is a morphism in $\Alg_{\calO}(\overline{\calC})_{\MM}$; here we abuse notation
by identifying $A$ with the associated trivial $\calO$-algebra object of $\overline{\calC}^{\otimes}_{\MM}$. It follows from Proposition \ref{canc} that the forgetful functor $\theta$ induces an identification
between centralizers of $f$ in $\Alg_{\calO}(\calC)_{\AAA}$ and centralizers of
$f'$ in $\Alg_{\calO}(\overline{\calC})_{\AAA}$. Combining this observation with
Proposition \ref{louda}, we obtain the following result:

\begin{theorem}\label{serpses}
Let $\calO^{\otimes}$ be a coherent $\infty$-operad,
$\calO^{\otimes} \times \LMod \rightarrow \calD^{\otimes}$ a bifunctor of $\infty$-operads,
and $q: \calC^{\otimes} \rightarrow \calD^{\otimes}$ a coCartesian fibration of $\infty$-operads.
Let $f: A \rightarrow B$ be a morphism in $\Alg_{\calO}(\calC)_{\MM}$  and let $\overline{\calC}^{\otimes}$ and $f': A \rightarrow B'$ be defined as above.
Assume that:
\begin{itemize}
\item[$(\ast)$] For every object $X \in \calO$, there exists a morphism object
$\Mor_{ \overline{\calC}_{X,\MM}}( A(X), B'(X) ) \in \overline{\calC}_{X, \AAA}$.
\end{itemize}
Then:
\begin{itemize}
\item[$(1)$] There exists a centralizer $\Cent{f} \in \Alg_{\calO}( \calC)_{\AAA}$.
\item[$(2)$] Let $Z$ be an arbitrary object of $\Alg_{\calO}(\calC)_{\AAA}$, and let $\sigma:$
$$ \xymatrix{ & Z \otimes A \ar[dr]^{g} & \\
A \ar[rr]^{f} \ar[ur] & & B }$$
be a commutative diagram in $\Alg_{\calO}(\calC)_{\MM}$. Let $Z'$ be a preimage of
$Z$ in $\Alg_{\calO}( \overline{\calC})_{\AAA}$, so that $\sigma$ lifts (up to homotopy) to a commutative diagram
$$ \xymatrix{ & Z' \otimes A \ar[dr] & \\
A \ar[ur] \ar[rr]^{f'} & & B' }$$
in $\Alg_{\calO}(\overline{\calC})_{\MM}$. Then $\sigma$ exhibits $Z$ as a centralizer of
$f$ if and only if, for every object $X \in \calO$, the induced map
$Z'(X) \otimes A(X) \rightarrow B'(X)$ exhibits $Z'(X)$ as a morphism object
$\Mor_{ \overline{\calC}_{X,\MM}}( A(X), B'(X) ) \in \overline{\calC}_{X, \AAA}$.
\end{itemize}
\end{theorem}

\begin{corollary}\label{placet}
In the situation of Theorem \ref{serpses}, suppose that $\calO^{\otimes}$ is the $0$-cubes
$\infty$-operad $\OpE{0}$. Then we can identify centralizers of a morphism
$f: A \rightarrow B$ in $\Alg_{\calO}(\calC)_{\MM}$ with morphism objects
$\Mor_{ \calC_{\MM}}(A,B)$ in $\calC_{\AAA}$. 
\end{corollary}

\begin{proof}
Combine Theorem \ref{serpses} with Remark \ref{hurley}.
\end{proof}

\begin{remark}
More informally, we can state Theorem \ref{serpses} as follows: the centralizer of a morphism
$f: A \rightarrow B$ can be identified with the classifying object for $A$-module maps from
$A$ to $B$. In particular, the center $\Cent{A}$ can be identified with the endomorphism algebra of
$A$, regarded as a module over itself.
\end{remark}

We now return to the proof of our main result.

\begin{proof}[Proof of Theorem \ref{centex}]
Combine Theorem \ref{serpses}, Proposition \monoidref{enterich}, and Theorem \symmetricref{turkwell}.
\end{proof}

We conclude this section with the proof of Proposition \ref{squisher}. First, we need a lemma.

\begin{lemma}\label{elial}
Let $n \geq 2$, and let $\calC \rightarrow \Delta^n$ be an inner fibration of $\infty$-categories.
Let $q: \calD \rightarrow \calE$ be another inner fibration of $\infty$-categories. 
Every lifting problem of the form
$$ \xymatrix{ \Lambda^n_0 \times_{\Delta^n} \calC \ar[r]^-{g} \ar[d] & \calD \ar[d]^{q} \\
\calC \ar@{-->}[ur] \ar[r] & \calE }$$
admits a solution, provided that $g| \Delta^{ \{0,1\} } \times_{ \Delta^n } \calC$ is
a $q$-left Kan extension of $g| \{0\} \times_{ \Delta^n} \calC$.
\end{lemma}

\begin{proof}
We first define a map $r: \Delta^n \times \Delta^1 \rightarrow \Delta^n$, which is given
on vertices by the formula
$$ r(i,j) = \begin{cases} 0 & \text{ if } (i,j) = (1,0) \\
i & \text{otherwise,} \end{cases}$$
and let $j: \Delta^n \rightarrow \Delta^n \times \Delta^1$ be the map
$(\id,j_0)$, where $j_0$ carries the first two vertices of $\Delta^n$ to
$\{0\} \subseteq \Delta^1$ and the remaining vertices to
$\{1\} \subseteq \Delta^1$.

Let $K = (\Lambda^n_0 \times \Delta^1) \coprod_{ \Lambda^n_0 \times \{0\} }
( \Delta^n \times \{0\})$, let $\calC' = (\Delta^n \times \Delta^1) \times_{ \Delta^n} \calC$, and let
$\calC'_0 = K \times_{ \Delta^n} \calC$. We will show that there exists a solution to the lifting problem $$ \xymatrix{ \calC'_0 \ar[r] \ar[d]^{g'} & \calD \ar[d]^{q} \\
\calC' \ar[r] \ar@{-->}[ur] & \calE. }$$
Composing this solution with the map $\calC \rightarrow \calC'$ induced by $j$, we will obtain
the desired result.

For every simplicial subset $L \subseteq \Delta^n$, let $\calC'_{L}$ denote the fiber product
$$((L \times \Delta^1) \coprod_{ L \times \{0\} } (\Delta^n \times \{0\} )) \times_{ \Delta^n} \calC,$$
and let $X_{L}$ denote the full subcategory of 
$\Fun_{ \calE}( \calC'_{L}, \calD) \times_{ \Fun_{\calE}(\calC'_{\emptyset}, \calD)} 
\{ g' | \calC'_{\emptyset} \}$ spanned
by those functors $F$ with the following property:
for each vertex $v \in L$, the restriction of
$F$ to $(\{v\} \times \Delta^1) \times_{\Delta^n} \calC$ is a $q$-left Kan extension
of $F| ( \{ v \} \times \{0\}) \times_{ \Delta^n} \calC$. 

To complete the proof, it will suffice to show
that the restriction map $X_{ \Delta^n} \rightarrow X_{ \Lambda^n_0}$ is surjective on vertices.
We will prove the following stronger assertion:
\begin{itemize}
\item[$(\ast)$] For every inclusion $L' \subseteq L$ of simplicial subsets of $\Delta^n$, the
restriction map $\theta_{L',L}: X_{L} \rightarrow X_{L'}$ is a trivial Kan fibration.
\end{itemize}

The proof proceeds by induction on the number of nondegenerate simplices of $L$. If $L' = L$, then $\theta_{L',L}$ is an isomorphism and there is nothing to prove. Otherwise, choose a nondegenerate simplex $\sigma$ of $L$ which does not belong
to $L'$, and let $L_0$ be the simplicial subset of $L$ obtained by removing $\sigma$.
The inductive hypothesis guarantees that the map $\theta_{L', L_0}$ is a trivial Kan fibration.
Consequently, to show that $\theta_{L',L}$ is a trivial Kan fibration, it will suffice to show that
$\theta_{L_0,L}$ is a trivial Kan fibration. Note that $\theta_{L_0,L}$ is a pullback of the map
$\theta_{ \bd \sigma, \sigma}$: we may therefore assume without loss of generality that
$L = \sigma$ is a simplex of $\Delta^n$.

Since the map $\theta_{L',L}$ is evidently a categorical fibration, it will suffice to show that
each $\theta_{L',L}$ is a categorical equivalence. We may assume by the inductive hypothesis that $\theta_{\emptyset,L'}$ is a categorical equivalence. By a two-out-of-three argument, we may reduce to the problem of showing that $\theta_{\emptyset,L'} \circ \theta_{L',L} = \theta_{\emptyset,L}$ is a categorical equivalence. In other words, we may assume that $L'$ is empty. We are now reduced
to the problem of showing that the map $X_{\sigma} \rightarrow X_{\emptyset}$ is
a trivial Kan fibration, which follows from Proposition \toposref{lklk}.
\end{proof}

\begin{proof}[Proof of Proposition \ref{squisher}]
It follows from Proposition \symmetricref{supercae} and Remark \symmetricref{saef} that
$q'$ is a categorical fibration and the induced map $\Mod^{\calO,S}_{A}(\calC)^{\otimes}
\rightarrow \Nerve(\FinSeg) \times S$ exhibits $\Mod^{\calO,S}_{A}(\calC)^{\otimes}$ as
an $S$-family of $\infty$-operads (note that the projection $\Mod^{\calO,S}_{A}(\calC)^{\otimes}
\rightarrow \Mod^{\calO}_{A}(\calC)^{\otimes}$ is an equivalence of $\infty$-categories). 
To complete the proof of $(1)$, it will suffice to show
that $q'$ satisfies condition $(iii)$ of Definition \ref{sably}. That is, we must show that if
$M$ is an object of $\Mod^{\calO,S}_{A}(\calC)^{\otimes}$ having image
$(X,s)$ in $\calO^{\otimes} \times S$ and $f: s \rightarrow s'$ is a morphism in
$S$, then $(\id_X, f)$ can be lifted to a $q'$-coCartesian morphism $M \rightarrow M'$ in
$\Mod^{\calO,S}_{A}(\calC)^{\otimes}$. Let $\calA$ be the full subcategory of
$(\calO^{\otimes})^{X/}$ spanned by the semi-inert morphisms $X \rightarrow Y$ in
$\calO^{\otimes}$, and let $\calA_0$ be the full subcategory of $\calA$ spanned by the null morphisms. 
The object $M \in \Mod^{\calO,S}_{A}(\calC)^{\otimes}$ determines a functor
$F: \calA \rightarrow \calC^{\otimes}$. Let $F_0$ denote the composite map
$\calA \stackrel{F}{\rightarrow} \calC^{\otimes} \rightarrow \calO^{\otimes}$.
Since $q$ exhibit $\calC^{\otimes}$ as
a coCartesian $S$-family of $\infty$-operads, there exists a $q$-coCartesian natural transformation
$H: \calA \times \Delta^1 \rightarrow \calC^{\otimes}$ from $F$ to another map $F'$, such that
$q \circ H$ is the product map $\calA \times \Delta^1 \stackrel{F_0 \times f}{\rightarrow} \calO^{\otimes} \times S$. Let $H': \calA_0 \times \Delta^1 \rightarrow \calC^{\otimes}$ be the composition
$\calA_0 \times \Delta^1 \rightarrow \calO^{\otimes} \times S \stackrel{A}{\rightarrow} \calC^{\otimes}$.
Since $A$ is a coCartesian $S$-family of $\calO$-algebras, the functors $H | \calA_0 \times \Delta^1$ and $H'$ are equivalent; we may therefore modify $H$ by a homotopy (fixed on
$\calA \times \{0\}$) and thereby assume that $H| \calA_0 \times \Delta^1 = H'$, so that
$H$ determines a morphism $\alpha: M \rightarrow M'$ in $\Mod^{\calO,S}_{A}(\calC)^{\otimes}$
lying over $(\id_{X},f)$. To complete the proof of $(1)$, it will suffice to show that
$\alpha$ is $q'$-coCartesian. 

Let $\calC_0^{\otimes} = \calO^{\otimes} \times S$. 
We have a commutative diagram of $\infty$-categories
$$ \xymatrix{ \Mod^{\calO}_{A}(\calC)^{\otimes} \ar[r]^{r} \ar[d]^{q'} & \MMod^{\calO}(\calC)^{\otimes}
\times_{ \MMod^{\calO}(\calC_0)^{\otimes}} \calC_0^{\otimes} \ar[dr]^{p} \ar[d]^{p'} & \\
\calO^{\otimes} \times S \ar[r] & \calO^{\otimes} \times \PrAlg_{\calO}(\calC) \times_{\PrAlg_{\calO}(\calC_0)} \calC_0^{\otimes} \ar[r]_-{p''} & \calC_0^{\otimes} }$$
Since the upper square is a pullback diagram, it will suffice to show that $r(\alpha)$ is
$p'$-coCartesian. In view of Proposition \toposref{stuch}, we are reduced to showing that
$r(\alpha)$ is $p$-coCartesian and that $(p' \circ r)(\alpha)$ is $p''$-coCartesian.

To prove that $r(\alpha)$ is $p$-coCartesian, we must show that every lifting problem of the form
$$ \xymatrix{ \Lambda^n_{0} \ar[r]^-{g} \ar[d] & \MMod^{\calO}(\calC)^{\otimes} \times_{ \MMod^{\calO}(\calC)^{\otimes} } \calC_0 \ar[d] \\
\Delta^n \ar@{-->}[ur] \ar[r] & \calC_0^{\otimes} }$$
admits a solution, provided that $n \geq 2$ and that $g$ carries the initial edge of $\Lambda^n_{0}$ to the morphism determined by $H$. Unwinding the definitions, this amounts to solving a lifting problem of the form
$$ \xymatrix{ \Lambda^n_0 \times_{ \calO^{\otimes} } \calK_{\calO} \ar[d] \ar[r]^-{G} & \calC^{\otimes} \ar[d] \\
\Delta^n \times_{ \calO^{\otimes}} \calK_{\calO} \ar[r] \ar@{-->}[ur] & \calO^{\otimes} \times S. }$$
The existence of a solution to this lifting problem is guaranteed by Lemma \ref{elial}. The assertion that
$(p' \circ r)(\alpha)$ is $p''$-coCartesian can be proven in the same way. This completes the proof of $(1)$.

Let $\overline{f}: M \rightarrow M''$ be as in $(2)$, let $f$ be the image of $\overline{f}$ in $S$, and let
$\widetilde{f}: M \rightarrow M'$ be the $q'$-coCartesian map constructed above.
We have a commutative diagram
$$ \xymatrix{ & M' \ar[dr]^{g} & \\
M \ar[ur]^{ \widetilde{f}} \ar[rr]^{ \overline{f}} & & M''.}$$
Let $\theta: \Mod^{\calO}_{A}(\calC)^{\otimes} \rightarrow \calC^{\otimes}$ be the forgetful functor.
By construction, $\theta( \overline{f})$ is a $q$-coCartesian morphism in $\calC^{\otimes}$, so that
$\theta( \widetilde{f})$ is $q$-coCartesian if and only if $\theta(g)$ is an equivalence.
We note that $\overline{f}$ is $q'$-coCartesian if and only if the map $g$ is an equivalence.
The ``only if'' direction of $(2)$ is now obvious, and the converse follows from Remark
\ref{scol} together with Corollary \symmetricref{postsabel}.
\end{proof}

\subsection{The Adjoint Representation}\label{adjrep}

Let $A$ be an associative ring, and let $A^{\times}$ be the collection of units in $A$. 
Then $A^{\times}$ forms a group, which acts on $A$ by conjugation. This action is given by a group homomorphism $\phi: A^{\times} \rightarrow \Aut(A)$ whose kernel is the subgroup of $A^{\times}$ consisting of units which belong to the center: this group can be identified with the group of units
of the center $\Cent{A}$. In other words, we have an exact sequence of groups
$$ 0 \rightarrow \Cent{A}^{\times} \rightarrow A^{\times} \rightarrow \Aut(A).$$
Our goal in this section is to prove a result which generalizes this statement in the following ways:
\begin{itemize}
\item[$(a)$] In place of a single associative ring $A$, we will consider instead a map of
algebras $f: A \rightarrow B$. In this setting, we will replace the automorphism group
$\Aut(A)$ by the set $\Hom(A,B)$ of algebra homomorphisms from $A$ to $B$. This set
is acted on (via conjugation) by the group $B^{\times}$ of units in $B$. Moreover, the
stabilizer of the element $f \in \Hom(A,B)$ can be identified with the group of units
$\Cent{f}^{\times}$ of the centralizer of the image of $f$. In particular, we have an exact
sequence of pointed sets
$$ \Cent{f}^{\times} \hookrightarrow B^{\times} \rightarrow \Hom(A,B).$$

\item[$(b)$] Rather than considering rings (which are associative algebra objects of the category of abelian groups), we will consider algebra objects in an arbitrary symmetric monoidal $\infty$-category $\calC$. In this setting, we need to determine appropriate analogues of the sets
$\Cent{f}^{\times}$, $B^{\times}$, and $\Hom(A,B)$ considered above. In the last case this is straightforward: the analogue of the set $\Hom(A,B)$ of ring homomorphisms from $A$ to $B$
is the {\em space} $\bHom_{ \Alg(\calC)}(A,B)$ of morphisms in the
$\infty$-category $\Alg(\calC)$. In \S \ref{kloop}, we will define unit subspaces
$\Centt{\Ass}{f}^{\times} \subseteq \bHom_{ \calC}( {\bf 1}, \Centt{\Ass}{f})$ and $B^{\times}
\subseteq \bHom_{\calC}( {\bf 1}, B)$. 

The collection of units in an associative ring $R$ is equipped with the structure of group (with respect to multiplication). We will see that there is an analogous structure on the space of units
$B^{\times}$ for an associative algebra object $B$ of an arbitrary symmetric monoidal
$\infty$-category $\calC$: namely, $B^{\times}$ is a loop space. That is, there exists
a pointed space $X(B)$ and a homotopy equivalence $B^{\times} \simeq \Omega X(B)$.
There is an ``action'' of the loop space $B^{\times}$ on the mapping space $\bHom_{ \Alg(\calC}(A,B)$. This action is encoded by a fibration $X(A,B) \rightarrow X(B)$, whose homotopy fiber
(over the base point of $X(B)$) can be identified with $\bHom_{ \Alg(\calC)}(A,B)$.
In particular, a morphism of associative algebra objects $f: A \rightarrow B$
determines a base point of $X(A,B)$, and we will see that the loop space $\Omega X(A,B)$
can be identified with the the space of units $\Centt{\Ass}{f}^{\times}$. In other words, we have a fiber sequence of spaces
$$ \bHom_{ \Alg(\calC)}( A,B) \rightarrow X(A,B) \rightarrow X(B)$$
which, after looping the base and total space, yields a fiber sequence
$$ \Centt{\Ass}{f}^{\times} \rightarrow B^{\times} \rightarrow \bHom_{ \Alg(\calC)}(A,B)$$
analogous to the exact sequence of sets described in $(a)$.

\item[$(c)$] Instead of considering only associative algebras, we will consider algebras over an arbitrary little cubes operad $\OpE{k}$ (according to Example \ref{sulta}, we can recover the case
of associative algebras by setting $k=1$). If $B$ is an $\OpE{k}$-algebra object of a symmetric monoidal
$\infty$-category $\calC$, then we can again define a space of units $B^{\times} \subseteq
\bHom_{ \calC}( {\bf 1}, B)$. The space $B^{\times}$ has the structure of a $k$-fold loop space:
that is, one can define a pointed space $X(B)$ and a homotopy equivalence
$B^{\times} \simeq \Omega^{k} X(B)$. If $A$ is another $\OpE{k}$-algebra object of
$\calC$, then there exists a fibration $X(A,B) \rightarrow X(B)$ whose fiber (over a well-chosen point of $X(B)$)
can be identified with $\bHom_{ \Alg_{ \OpE{k}}(\calC)}(A,B)$. In particular, every $\OpE{k}$-algebra map $f: A \rightarrow B$ determines a base point of the total space $X(A,B)$, and the
$k$-fold loop space $\Omega^{k} X(A,B)$ can be identified with the space of units
$\Cent{\OpE{k}}{f}^{\times}$ (see Definition \ref{clupus} below). We therefore have a fiber sequence of spaces
$$ \bHom_{ \Alg_{\OpE{k}}(\calC)}( A,B) \rightarrow X(A,B) \rightarrow X(B)$$
which yields, after passing to loop spaces repeatedly, a fiber sequence
$$ \Centt{\OpE{k}}{f}^{\times} \rightarrow B^{\times} \stackrel{\phi}{\rightarrow} \Omega^{k-1} \bHom_{ \Alg_{\OpE{k}}(\calC)}(A,B).$$
We should regard the map $\phi$ as a $k$-dimensional analogue of the adjoint action
of the unit group $B^{\times}$ of an associative ring $B$ on the set of maps
$\Hom(A,B)$.
\end{itemize}

Our first step is to {\em define} the spaces of units appearing in the above discussion. This requires a bit of a digression.

\begin{definition}\label{huntrex}
Let $\calO^{\otimes}$ be an $\infty$-operad, and let $\calO^{\otimes}_{\ast} \subseteq
\Fun( \Delta^1, \calO^{\otimes})$ be the $\infty$-category of pointed objects of $\calO^{\otimes}$.
The forgetful functor $q: \calO^{\otimes}_{\ast} \rightarrow \calO^{\otimes}$ is a left fibration
of simplicial sets. We let $\chi_{\calO}: \calO^{\otimes} \rightarrow \SSet$ denote a functor which classifies $q$.
\end{definition}

\begin{proposition}
Let $q: \calO^{\otimes} \rightarrow \Nerve(\FinSeg)$ be an $\infty$-operad and let $\chi_{\calO}: \calO^{\otimes} \rightarrow \SSet$
be as in Definition \ref{huntrex}. Then $\chi_{\calO}$ is a $\calO$-monoid object of $\SSet$.
\end{proposition}

\begin{proof}
We must show that if $X \in \calO^{\otimes}_{\seg{n}}$, and if $\alpha_i: X \rightarrow X_i$ are a collection of inert morphisms in $\calO^{\otimes}$ lifting the maps $\Colp{i}: \seg{n} \rightarrow \seg{1}$ for $1 \leq i \leq n$, then the induced map $\chi_{\calO}( X) \rightarrow \prod_{1 \leq i \leq n} \chi_{\calO}(X_i)$ is a homotopy equivalence. Let $0$ denote a final object of $\calO^{\otimes}$; then
the left hand side is homotopy equivalent to $\bHom_{ \calO^{\otimes}}(0,X)$, while the right hand side is homotopy equivalent to $\prod_{1 \leq i \leq n} \bHom_{ \calO^{\otimes}}( 0, X_i)$. The desired result now follows from the observation that the maps $\alpha_i$ exhibit $X$ as a $q$-product of the objects
$\{ X_i \}_{1 \leq i \leq n}$.
\end{proof}

\begin{remark}
An $\infty$-operad $\calO^{\otimes}$ is unital if and only if the functor
$\chi_{\calO}: \calO^{\otimes} \rightarrow \SSet$ is equivalent to the constant functor taking
the value $\Delta^0$.
\end{remark}

\begin{definition}\label{clupus}
Let $q: \calC^{\otimes} \rightarrow \calO^{\otimes}$ be a fibration of $\infty$-operads, where
$\calO^{\otimes}$ is unital, and let $\chi_{\calC}: \calC^{\otimes} \rightarrow
\SSet$ be as in Definition \ref{huntrex}. Composition with $\chi_{\calC}$ determines a functor
$\Alg_{\calO}(\calC) \rightarrow \Mon_{ \calO}( \SSet)$. 

Suppose that $\calO^{\otimes} = \OpE{k}$, where $k > 0$. Since the collection of grouplike
$\OpE{k}$-spaces is stable under colimits in $\Mon_{ \OpE{k}}(\SSet)$ (Remark \ref{jazwind})
the inclusion $i: \Mon_{ \OpE{k}}^{\glike}(\SSet) \subseteq \Mon_{ \OpE{k}}(\SSet)$ preserves small colimits. It follows from Proposition \ref{preslage} that $\Mon_{\OpE{k}}^{\glike}(\SSet)$ is equivalent
to $\SSet^{\geq k}_{\ast}$, and therefore presentable. Using Corollary \toposref{adjointfunctor}, we deduce that the inclusion functor $i$ admits a right adjoint $G$. We let $\GL_1: \Alg_{\OpE{k}}(\calC) \rightarrow \Mon_{ \OpE{k}}^{\glike}(\SSet)$ denote the composite functor
$$ \Alg_{\calO}(\calC) \stackrel{ \chi_{\calC} \circ}{\rightarrow} \Mon_{ \OpE{k}}(\calC)
\stackrel{G}{\rightarrow} \Mon_{ \OpE{k}}^{\glike}(\calC).$$
If $A \in \Alg_{\OpE{k}}(\calC)$, we will often write $A^{\times}$ in place of
$\GL_1(A)$; we will refer to $A^{\times}$ as the {\it $\OpE{k}$-space of units in $A$}.

In the special case $k=0$, we let $\GL_1: \Alg_{ \OpE{k}}(\calC) \rightarrow \Mon_{ \OpE{k}}(\SSet)
\simeq \SSet_{\ast}$ be the functor defined by composition with $\chi_{\calC}$; we will also denote
this functor by $A \mapsto A^{\times}$.
\end{definition}

We are now prepared to state our main result:

\begin{theorem}\label{machus}
Let $\calC^{\otimes}$ be a symmetric monoidal $\infty$-category. Assume that the underlying
$\infty$-category $\calC$ is presentable and that the tensor product $\otimes: \calC \times \calC \rightarrow \calC$ preserves small colimits separately in each variable. Fix an integer $k \geq 0$, and let
$\bHom: \Alg_{ \OpE{k}}(\calC)^{op} \times \Alg_{ \OpE{k}}(\calC) \rightarrow \SSet$ be the adjoint of the Yoneda embedding $\Alg_{ \OpE{k}}(\calC) \rightarrow
\Fun( \Alg_{ \OpE{k}}(\calC)^{op}, \SSet)$. There exists another functor
$X: \Alg_{ \OpE{k}}( \calC)^{op} \times \Alg_{ \OpE{k}}(\calC) \rightarrow \SSet$ and a natural transformation $\alpha: \bHom \rightarrow X$ with the following properties:
\begin{itemize}
\item[$(1)$] For every object $B \in \Alg_{\OpE{k}}(\calC)$ and every morphism
$f: A' \rightarrow A$ in $\Alg_{\OpE{k}}(\calC)$, the diagram
$$ \xymatrix{ \bHom(A,B) \ar[r] \ar[d] & \bHom(A',B) \ar[d] \\
X(A,B) \ar[r] & X(A',B) }$$
is a pullback square.
\item[$(2)$] Let $f: A \rightarrow B$ be a morphism in $\Alg_{\OpE{k}}(\calC)$, so
that the map $f$ determines a base point of the space $X(A,B)$ (via $\alpha$). Then there is a canonical
homotopy equivalence $\Omega^{k} X(A,B) \simeq \Centt{\OpE{k}}{f}^{\times}$.
\end{itemize}
\end{theorem}

\begin{remark}
In the situation of Theorem \ref{machus}, it suffices to prove assertion $(1)$ in the case
where $A'$ is the initial object ${\bf 1} \in \Alg_{\OpE{k}}( \calC)$. This follows
by applying Lemma \toposref{transplantt} to the diagram
$$ \xymatrix{ \bHom(A,B) \ar[r] \ar[d] & \bHom(A',B) \ar[r] \ar[d] & \bHom( {\bf 1}, B) \ar[d] \\
X(A,B) \ar[r] & X(A',B) \ar[r] & X( {\bf 1},B). }$$
\end{remark}

\begin{remark}
In the special case where $A'$ is the initial object ${\bf 1} \in \Alg_{\OpE{k}}(\calC)$,
the space $\bHom(A',B)$ is contractible, so that part $(1)$ of Theorem \ref{machus} asserts
the existence of a fiber sequence
$$ \bHom(A,B) \rightarrow X(A,B) \rightarrow X({ \bf 1}, B).$$
Fixing a base point $(f: A \rightarrow B) \in \bHom(A,B)$ and taking loop spaces repeatedly, we have a fiber sequence
$$ \Omega^{k} X(A,B) \rightarrow \Omega^{k} X( {\bf 1},B) \rightarrow
\Omega^{k-1} \bHom_{ \Alg_{\OpE{k}}(\calC)}(A,B)$$
We observe that there is a canonical natural transformation $\beta: \Centt{\OpE{k}}{f_0} \rightarrow B$ of
functors $\OpE{k} \rightarrow \calC^{\otimes}$. The natural transformation $\beta$ induces
an equivalence of $\OpE{k}$-spaces $\Centt{\OpE{k}}{f_0}^{\times} \rightarrow B^{\times}$. Invoking
part $(2)$ of Theorem \ref{machus}, we obtain the fiber sequence
$$ \Centt{\OpE{k}}{f}^{\times} \rightarrow B^{\times} \rightarrow \Omega^{k-1} \bHom_{\Alg_{\OpE{k}}(\calC)}(A,B)$$
described in $(c)$.
\end{remark}

An $\OpE{k}$-algebra object $A$ of
a symmetric monoidal $\infty$-category $\calC$ determines an $(\infty,n)$-category $C(A)$ enriched over $\calC$. One approach to the proof of Theorem \ref{machus} would be to define $X(A,B)$ to be the space of functors from $C(A)$ into $C(B)$. Since we do not wish to develop the theory of enriched
$(\infty,n)$-categories in this paper, our proof will proceed along somewhat different lines: we will use an inductive approach, which iteratively replaces the $\infty$-category $\calC$ by the $\infty$-category
$\Mod^{L}_{\calC}$ of $\infty$-categories left-tensored over $\calC$. To guarantee that this replacement does not destroy our hypothesis that $\calC$ is presentable, we need to introduce a few restrictions on the $\calC$-modules that we allow.

\begin{notation}
Let $\kappa$ be a regular cardinal. Recall that a presentable $\infty$-category
$\calC$ is {\it $\kappa$-compactly generated} if $\calC$ is generated by its $\kappa$-compact objects
under the formation of small, $\kappa$-filtered colimits (see \S \toposref{compactgen}). 
If $\calC$ and $\calD$ are $\kappa$-compactly generated $\infty$-categories, then we will say that
a functor $F: \calC \rightarrow \calD$ is {\it $\kappa$-good} if $F$ preserves small colimits and carries $\kappa$-compact objects of $\calC$ to $\kappa$-compact objects of $\calD$. Equivalently,
$F$ is $\kappa$-good if $F$ admits a right adjoint $G$ which commutes with $\kappa$-filtered colimits
(Proposition \toposref{comppress}). 

Let $\LPress$ denote the $\infty$-category of presentable
$\infty$-categories and colimit-preserving functors. We let $\LPress_{\kappa}$ denote the subcategory of
the $\infty$-category $\LPress$ whose objects are $\kappa$-compactly generated presentable
$\infty$-categories and whose morphisms are $\kappa$-good functors.
\end{notation}

\begin{lemma}\label{clopus}
Let $\kappa$ be an uncountable regular cardinal. Then:
\begin{itemize}
\item[$(1)$] Let $\calK$ denote the collection of all $\kappa$-small simplicial sets.
Then the functor $\calC \mapsto \Ind_{\kappa}(\calC)$ determines an equivalence of
$\infty$-categories from $\Cat_{\infty}(\calK)$ to $\LPress_{\kappa}$.
\item[$(2)$] The $\infty$-category $\LPress_{\kappa}$ is presentable.
\end{itemize}
\end{lemma}

\begin{proof}
Note that assertion $(2)$ follows immediately from $(1)$ and Lemma \ref{spuke}.
It is clear that the functor $\Ind_{\kappa}: \Cat_{\infty}(\calK) \rightarrow \LPress_{\kappa}$ is
essentially surjective. To prove that it is fully faithful, it will suffice to show that for every
pair of $\infty$-categories $\calC, \calD \in \Cat_{\infty}(\calK)$, the canonical map
$\theta: \Fun( \calC, \calD) \rightarrow \Fun( \Ind_{\kappa}(\calC), \Ind_{\kappa}(\calD))$ induces
an equivalence of $\infty$-categories from the full subcategory $\Fun'(\calC, \calD)$ of $\Fun( \calC, \calD)$ spanned by the
those functors which preserve $\calK$-indexed colimits to the full subcategory $\Fun'( \Ind_{\kappa}(\calC), \Ind_{\kappa}(\calD) )$ of $\Fun( \Ind_{\kappa}(\calC), \Ind_{\kappa}(\calD) )$ spanned by the $\kappa$-good functors. Let $\Fun'( \calC, \Ind_{\kappa}(\calD))$ denote the full subcategory of
$\Fun( \calC, \Ind_{\kappa}(\calD) )$ consisting of those functors which preserve $\calK$-indexed colimits and carry $\calC$ into the full subcategory of $\Ind_{\kappa}(\calD)$ spanned by the $\kappa$-compact objects. We have a homotopy commutative diagram of $\infty$-categories
$$ \xymatrix{ \Fun'(\calC, \calD) \ar[rr]^{\theta} \ar[dr]^{\theta'} & & \Fun'( \Ind_{\kappa}(\calC), \Ind_{\kappa}(\calD) ) \ar[dl]^{\theta''} \\
& \Fun'( \calC, \Ind_{\kappa}(\calD) ), & }$$
where $\theta'$ and $\theta''$ are given by composing with the Yoneda embeddings for $\calD$ and
$\calC$, respectively. To complete the proof, it will suffice to show that $\theta'$ and $\theta''$ are categorical equivalences.

To show that $\theta'$ is a categorical equivalence, let $\calD'$ denote the collection of all
$\kappa$-compact objects of $\Ind_{\kappa}(\calD)$. Since $\calD'$ is stable under $\kappa$-small colimits in $\Ind_{\kappa}(\calD)$, $\Fun'( \calC, \Ind_{\kappa}(\calD) )$ is isomorphic to the full subcategory of $\Fun( \calC, \calD')$ spanned by those functors which preserve $\kappa$-small colimits. It will therefore suffice to show that the Yoneda embedding induces an equivalence $\calD \rightarrow \calD'$. Lemma \toposref{stylus} guarantees that $\calD'$ is an idempotent completion of $\calD$; it will therefore suffice to show that $\calD$ is idempotent complete. This follows from Proposition \toposref{slanger}, since $\kappa$ is assumed to be uncountable so that $\calD$ admits sequential colimits.

Repeating the previous argument with $\calC$ in place of $\calD$, we see that an object of
$\Ind_{\kappa}(\calC)$ is $\kappa$-compact if and only if it lies in the image of the Yoneda embedding
$j: \calC \rightarrow \Ind_{\kappa}(\calC)$. Consequently, to prove that $\theta''$ is a categorical equivalence, it suffices to show that composition with $j$ induces an equivalence from the
full subcategory of $\Fun( \Ind_{\kappa}(\calC), \Ind_{\kappa}(\calC) )$ spanned by those functors which preserve small colimits to the full subcategory of $\Fun( \calC, \Ind_{\kappa}(\calC) )$ spanned by those functors which preserve $\kappa$-small colimits; this follows from Proposition \toposref{sumatch}.
\end{proof}

\begin{remark}
In the statement of Lemma \ref{clopus}, we can drop the requirement that the cardinal $\kappa$ is uncountable if we are willing to restrict our attention to idempotent complete $\infty$-categories.
\end{remark}

We now study the interaction between the subcategory $\LPress_{\kappa} \subseteq \LPress$
and the symmetric monoidal structure $\LPress^{\otimes}$ on $\LPress$ constructed in 
\S \symmetricref{comm7}. Let $\LPress^{\otimes}_{\kappa}$ denote the subcategory of
$\LPress^{\otimes}$ whose objects are finite sequences $(\calC_1, \ldots, \calC_n)$ where
each of the $\infty$-categories $\calC_i$ is $\kappa$-compactly generated, and whose morphisms
are given by maps $( \calC_1, \ldots, \calC_m) \rightarrow ( \calD_1, \ldots, \calD_n)$ covering a map
$\alpha: \seg{m} \rightarrow \seg{n}$ in $\FinSeg$ such that the functor
$\prod_{ \alpha(i) = j } \calC_i \rightarrow \calD_j$ preserves $\kappa$-compact objects for
$1 \leq j \leq n$.

\begin{lemma}\label{clopus2}
Let $\kappa$ be an uncountable regular cardinal. Then:
\begin{itemize}
\item[$(1)$] If $\calC$ and $\calD$ are $\kappa$-compactly generated presentable monoidal $\infty$-categories, then $\calC \otimes \calD$ is $\kappa$-compactly presented. Moreover, the collection
of $\kappa$-compact objects of $\calC \otimes \calD$ is generated under $\kappa$-small colimits
by tensor products of the form $C \otimes D$, where $C \in \calC$ and $D \in \calD$ are $\kappa$-compact.
\item[$(2)$] The composite map $\LPress_{\kappa}^{\otimes} \subseteq \LPress^{\otimes} \rightarrow \Nerve(\FinSeg)$ exhibits $\LPress_{\kappa}^{\otimes}$ as a symmetric monoidal $\infty$-category, and the inclusion $\LPress_{\kappa}^{\otimes} \subseteq \LPress^{\otimes}$ is a symmetric monoidal functor.
\item[$(3)$] Let $\calK$ denote the collection of all $\kappa$-small simplicial sets. The functor $\Ind_{\kappa}$ induces a symmetric monoidal equivalence
$\Cat_{\infty}(\calK)^{\otimes} \rightarrow \LPress_{\kappa}^{\otimes}$. 
\item[$(4)$] The tensor product $\otimes: \LPress_{\kappa} \times \LPress_{\kappa} \rightarrow \LPress_{\kappa}$ preserves colimits separately in each variable.
\end{itemize}
\end{lemma}

\begin{proof}
Recall (see \S \symmetricref{comm7}) that $\Ind_{\kappa}$ determines a symmetric monoidal
functor from $\Cat_{\infty}(\calK)$ to $\LPress$. To prove $(1)$, we note that if $\calC \simeq \Ind_{\kappa}(\calC_0)$ and $\calD \simeq \Ind_{\kappa}(\calD_0)$, then
$\calC \otimes \calD \simeq \Ind_{\kappa}( \calC_0 \otimes \calD_0 )$ is a $\kappa$-compactly generated $\infty$-category. To prove the second assertion of $(1)$, it suffices to show that
$\calC_0 \otimes \calD_0$ is generated under $\kappa$-small colimits by the essential image of the functor $\calC_0 \times \calD_0 \rightarrow \calC_0 \otimes \calD_0$, which is clear. 
Assertion $(2)$ follows immediately from $(1)$. Assertion $(3)$ follows from Lemma \ref{clopus}, and assertion $(4)$ follows from $(3)$ together with Lemma \ref{spuke}. 
\end{proof}

\begin{lemma}\label{silper}
Let $\calC^{\otimes}$ be a symmetric monoidal $\infty$-category. Assume that $\calC$ is presentable and that the tensor product $\otimes: \calC \times \calC \rightarrow \calC$ preserves small colimits separately in each variable. Then there exists an uncountable regular cardinal $\kappa$ with the following properties:
\begin{itemize}
\item[$(1)$] The $\infty$-category $\calC$ is $\kappa$-compactly generated.
\item[$(2)$] The tensor product $\otimes: \calC \times \calC \rightarrow \calC$ preserves 
$\kappa$-compact objects, and the unit object ${ \bf 1} \in \calC$ is $\kappa$-compact.
\item[$(3)$] For every algebra object $A \in \Alg(\calC)$, the $\infty$-category
$\Mod^{R}_{A}(\calC)$ is $\kappa$-compactly generated. 
\item[$(4)$] For every algebra object $A \in \Alg(\calC)$, the action functor
$\otimes: \calC \times \Mod^{R}_{A}(\calC) \rightarrow \Mod^{R}_{A}(\calC)$ preserves $\kappa$-compact objects.
\end{itemize}
\end{lemma}

\begin{proof}
Choose an regular cardinal $\kappa_0$ such that $\calC$ is $\kappa_0$-compactly generated. Let $\calC_0$ be the full subcategory of $\calC$ spanned by the $\kappa_0$-compact objects, and let $\calC_1$ denote the smallest full subcategory of $\calC$ which contains
$\calC_0$, the unit object of $\calC$, and the essential image of the tensor product functor
$\otimes: \calC_0 \times \calC_0 \rightarrow \calC$. Since $\calC_1$ is essentially small, there
exists a regular cardinal $\kappa > \kappa_0$ such that every object in $\calC_1$ is $\kappa$-small.
We claim that $\kappa$ has the desired properties. It is clear that $\kappa$ is uncountable and that $(1)$ is satisfied. 

To prove $(2)$, choose $\kappa$-compact objects $C, D \in \calC$. Then $C$ and $D$ can be written as $\kappa$-small colimits $\varinjlim( C_{\alpha})$ and $\varinjlim( D_{\beta})$, where the objects
$C_{\alpha}$ and $D_{\beta}$ are $\kappa_0$-compact. Then $C \otimes D \simeq \varinjlim( C_{\alpha} \otimes D_{\beta})$ is a $\kappa$-small colimit of objects belonging to
$\calC_1$, and is therefore $\kappa$-compact.

We now prove $(3)$. According to Corollary \monoidref{gloop}, the forgetful functor
$G: \Mod_{A}^{R}(\calC) \rightarrow \calC$ preserves $\kappa$-filtered colimits (in fact, all small colimits). It follows from Proposition \toposref{comppress} that the left adjoint $F$ to
$G$ preserves $\kappa$-compact objects. Let $\calX$ denote the full subcategory of
$\Mod^{R}_{A}(\calC)$ generated under small colimits by objects of the form $F(C)$, where
$C \in \calC$ is $\kappa$-compact; we will show that $\calX = \Mod^{R}_{A}(\calC)$. 
For each $M \in \Mod^{R}_{A}(\calC)$, we can write $M \simeq A \otimes_{A} M = | \Baar_{A}(A,M)_{\bigdot} |$ (see Proposition \monoidref{usss}). Consequently, to show that $M \in \calX$, it will suffice to show that $X$ contains $F( A^{\otimes n-1} \otimes G(M) )$ for each $n \geq 1$. We are therefore reduced to proving that $F(C) \in \calX$ for each $C \in \calC$, which is clear (the functor
$F$ preserves small colimits and $C$ can be written as a colimit of $\kappa$-compact objects of $\calC$ by $(1)$).

We now prove $(4)$. Let $\calY$ denote the full subcategory of $\Mod^{R}_A(\calC)$ spanned by those objects $M$ such that $C \otimes M \in \Mod^{R}_{A}(\calC)$ is $\kappa$-compact for every
$\kappa$-compact object $C \in \calC$. The $\infty$-category $\calY$ is evidently closed under
$\kappa$-small colimits in $\Mod^{R}_{A}(\calC)$. Since $C \otimes F(D) \simeq F(C \otimes D)$, it follows from $(2)$ that $\calY$ contains $F(D)$ for every $\kappa$-compact object $D \in \calC$.
Since every object of $\calY$ is $\kappa$-compact in $\Mod^{R}_{A}(\calC)$, we have a fully faithful embedding $f: \Ind_{\kappa}(\calY) \rightarrow \Mod^{R}_{A}(\calC)$, which preserves small colimits by Proposition \toposref{sumatch}. The essential image $\calY'$ of $f$ is stable under small colimits and contains $F(D)$ for every $\kappa$-compact object $D \in \calC$, so that $\calX \subseteq \calY'$.
It follows that $f$ is essentially surjective and therefore an equivalence of $\infty$-categories.
Lemma \toposref{stylus} now guarantees that the collection of $\kappa$-compact objects of
$\Mod^{R}_{A}(\calC)$ is an idempotent completion of $\calY$. Since $\kappa$ is uncountable,
$\calY$ is stable under sequential colimits and therefore idempotent complete. It follows that
$\calY$ contains every $\kappa$-compact object of $\Mod^{R}_{A}(\calC)$, as desired.
\end{proof}

We now proceed with the proof of our main result.

\begin{proof}[Proof of Theorem \ref{machus}]
We proceed by induction on $k$. Assume first that $k=0$.
Let $X$ denote the composite
functor
$$ \Alg_{ \OpE{k}}( \calC)^{op} \times \Alg_{ \OpE{k}}(\calC)
\rightarrow \calC^{op} \times \calC \stackrel{H'}{\rightarrow} \SSet,$$
where $H$ is the adjoint of the Yoneda embedding for $\calC$ (given informally by 
$H(C,C') = \bHom_{\calC}( C,C')$). The forgetful functor $\theta: \Alg_{\OpE{k}}(\calC) \rightarrow \calC$ determines a natural transformation of functors $\bHom \rightarrow X$. We claim that
this functor satisfies conditions $(1)$ and $(2)$ of Theorem \ref{machus}. 

Suppose we are given a morphism $A' \rightarrow A$ in
$\Alg_{\OpE{k}}( \calC)$ and an object $B \in \Alg_{\OpE{k}}(\calC)$.
Let ${\bf 1}$ denote the unit object of $\calC_{\MM}$. Proposition \symmetricref{ezalg} implies that
$\Alg_{ \OpE{k}}(\calC)$ is equivalent to $(\calC)^{ {\bf 1}/}$. It follows that we have a natural transformation of fiber sequences
$$ \xymatrix{ \bHom(A,B) \ar[r] \ar[d] & \bHom(A',B) \ar[d] \\
X(A,B) \ar[r] \ar[d] & X(A',B) \ar[d] \\
\bHom_{\calC}( {\bf 1}, \theta(B)) \ar[r] & \bHom_{\calC}( {\bf 1}, \theta(B) ).}$$
Since the bottom horizontal map is a homotopy equivalence, the upper square is a homotopy
pullback square. This proves $(1)$. To prove $(2)$, we invoke Corollary \ref{placet} to identify
$\Centt{\OpE{k}}{f}^{\times} = \bHom_{ \calC}( {\bf 1}, \Centt{\OpE{k}}{f}^{\times})$ with the mapping space
$\bHom_{ \calC}( \theta(A), \theta(B)) = X(A,B)$. 

We now treat the case where $k > 0$. Applying Corollary \ref{eliac} (in the setting of
$\infty$-categories which are not necessarily small, which admit small colimits) we
obtain a fully faithful embedding
$\psi: \Alg_{ \OpE{k}}(\calC) \rightarrow \Alg_{ \OpE{k-1}}( \Mod^{L}_{\calC}( \LPress))$.
Let $\kappa$ be an uncountable regular cardinal satisfying the conditions of Lemma \ref{silper}
and let $\calC' = \CMod_{\calC}( \LPress_{\kappa})$. Using Corollary \monoidref{underwhear},
Lemma \ref{clopus}, and Lemma \ref{clopus2}, we deduce that $\calC'$ is a presentable
$\infty$-category equipped with a symmetric monoidal structure, such that the tensor product
$\otimes: \calC' \times \calC' \rightarrow \calC'$ preserves colimits separately in each variable.
The functor $\psi$ induces a fully faithful embedding $\Alg_{ \OpE{k}}(\calC) \rightarrow \Alg_{ \OpE{k-1}}(\calC')$, which we will also denote by $\psi$. 

Let $\bHom': \Alg_{ \OpE{k-1}}(\calC')^{op} \times \Alg_{ \OpE{k-1}}( \calC') \rightarrow \SSet$ be the adjoint to the Yoneda embedding. Invoking the inductive hypothesis, we deduce that there
exists another functor $X': \Alg_{ \OpE{k-1}}( \calC')^{op} \times \Alg_{ \OpE{k-1}}(\calC') \rightarrow \SSet$ and a natural transformation $\alpha': \bHom' \rightarrow X'$ satisfying hypotheses
$(1)$ and $(2)$ for the $\infty$-category $\calC'$. Let $X$ denote the composition
$$ \Alg_{ \OpE{k}}(\calC)^{op} \times \Alg_{ \OpE{k} }(\calC) \stackrel{\psi \times \psi}{\rightarrow}
\Alg_{ \OpE{k-1}}( \calC')^{op} \times \Alg_{ \OpE{k-1}}(\calC') \stackrel{X'}{\rightarrow} \SSet.$$
Since $\psi$ is fully faithful, the composition
$$ \Alg_{ \OpE{k}}(\calC)^{op} \times \Alg_{ \OpE{k} }(\calC) \stackrel{\psi \times \psi}{\rightarrow}
\Alg_{ \OpE{k-1}}( \calC')^{op} \times \Alg_{ \OpE{k-1}}(\calC') \stackrel{\bHom'}{\rightarrow} \SSet$$
is equivalent to $\bHom$, so that $\alpha'$ induces a natural transformation of functors
$\alpha: \bHom \rightarrow X$. 

It is clear from the inductive hypothesis that the natural transformation $\alpha$ satisfies
condition $(1)$. We will prove that $\alpha$ satisfies $(2)$. Let $f: A \rightarrow B$ be a morphism
in $\Alg_{ \OpE{k}}( \calC)$, and let $\psi(f): \calA \rightarrow \calB$ be the induced morphism
in $\Alg_{ \OpE{k-1}}(\calC')$. Let $Z = \Centt{\OpE{k}}{f}$, so that we have a commutative diagram
$$ \xymatrix{ & Z \otimes A \ar[dr] & \\
A \ar[ur] \ar[rr]^{f} & & B.}$$
Applying the (symmetric monoidal) functor $\psi$, we obtain a diagram
$$ \xymatrix{ & \psi(Z) \otimes \calA \ar[dr] & \\
\calA \ar[ur] \ar[rr]^{\psi(f)} & & \calB, }$$
which is classified by a map $\beta: \psi(Z) \rightarrow \Centt{\OpE{k-1}}{\psi(f)}$. 
The inductive hypothesis guarantees a homotopy equivalence
$\Centt{\OpE{k-1}}{\psi(f)}^{\times} \simeq \Omega^{n-1} X'( \calA, \calB) \simeq \Omega^{n-1} X(A,B)$.
Passing to loop spaces, we get an homotopy equivalence
$\Omega \Centt{\OpE{k-1}}{\psi(f)}^{\times} \simeq \Omega^{n} X(A,B)$. We will complete the proof
by showing the following:
\begin{itemize}
\item[$(a)$] There is a canonical homotopy equivalence $Z^{\times} \simeq \Omega \psi(Z)^{\times}$.

\item[$(b)$] The map $\beta$ induces a homotopy equivalence $\Omega \psi(Z)^{\times}
\rightarrow \Omega \Centt{\OpE{k-1}}{\psi(f)}^{\times}$.
\end{itemize}

Assertion $(a)$ is easy: the space $\Omega \psi(Z)^{\times}$ can be identified with
the summand of the mapping space $\bHom_{ \Mod^{R}_{Z}(\calC)}( Z, Z)$ spanned by
the equivalences from $Z$ to itself. Corollary \monoidref{ingurtime} furnishes an identification $\bHom_{ \Mod^{R}_{Z}(\calC)}( Z, Z) \simeq \bHom_{ \calC}( {\bf 1}, Z)$, under which the summand
$\Omega \psi(Z)^{\times} \subseteq \bHom_{ \Mod^{R}_Z(\calC)}(Z,Z)$ corresponds to the space of units $Z^{\times}$.

The proof of $(b)$ is slightly more involved. We wish to show that $\beta$ induces a homotopy equivalence $$\phi: \Omega \bHom_{ \CMod_{\calC}( \LPress_{\kappa})}( \calC, \Mod^{R}_{Z}(\calC) )
\rightarrow \Omega \bHom_{ \CMod_{ \calC}( \LPress_{\kappa})}( \calC, \Centt{\OpE{k-1}}{\psi(f)}).$$
Let $\calD^{\otimes}$ be a unitalization of the symmetric monoidal $\infty$-category $\CMod_{\calC}( \LPress)^{\otimes}$, so that the underlying $\infty$-category of $\calD$ is equivalent to
$\CMod_{\calC}(\LPress)_{\calC/}$. Since $\OpE{k-1}$ is unital, we can regard
$\Mor^{R}_Z(\calC)$ and $\Centt{\OpE{k-1}}{\psi(f)}$ as $\OpE{k-1}$-algebra objects of $\calD$. 
Regard the $\infty$-category $\CMod_{\calC}(\LPress)$ as
tensored over spaces, and let $D = \calC \otimes S^1$ (see \S \toposref{quasilimit7}), regarded
as an object of $\calD$ by choosing a base point $\ast \in S^1$. Then we can identify
$\phi$ with the morphism
$\bHom_{ \calD}( D, \Mod^{R}_Z(\calC) ) \rightarrow \bHom_{\calD}( D, \Centt{\OpE{k-1}}{\psi(f)})$. 

Theorem \ref{postcur} guarantees that the construction $C \mapsto \Mod^{R}_{C}( \calC)$ determines a fully faithful embedding of symmetric monoidal $\infty$-categories $F: \Alg(\calC)^{\otimes} \rightarrow \calD^{\otimes}$. Theorem \ref{curly} guarantees that the underlying functor
$f: \Alg(\calC) \rightarrow \calD$ admits a right adjoint $g$, so that $f$ exhibits
$\Alg(\calC)$ as a colocalization of $\calD$ which is stable under tensor products in $\calD$.
Using Proposition \symmetricref{yak2}, we see that $g$ can be regarded as a lax symmetric monoidal functor, and induces a map $\gamma: \Alg_{ \OpE{k} }(\calC) \simeq \Alg_{ \OpE{k-1} }( \Alg(\calC) ) \rightarrow \Alg_{ \OpE{k-1} }( \calD)$ which is right adjoint to the functor given by composition with
$F$. Using the fact that $\psi$ is a fully faithful symmetric monoidal functor, we deduce that
$\gamma( \beta)$ is an equivalence in $\Alg_{ \OpE{k}}(\calC)$. Consequently, to prove
that $\phi$ induces an equivalence from $\bHom_{ \calD}( D, \Mod^{R}_Z(\calC) )$
to $\bHom_{ \calD}( D, \Centt{\OpE{k-1}}{\psi(f)})$, it will suffice to show that the object $D \in \calD$ lies in
the essential image of the functor $f$. In other words, we must show that there exists
an algebra object $K \in \Alg(\calC)$ such that $\calC \otimes S^1$ is equivalent to
$\Mod^{R}_{K}(\calC)$ in the $\infty$-category $\CMod_{\calC}( \LPress)_{\calC/ }$.
Choosing a symmetric monoidal functor $\SSet^{\times} \rightarrow \calC^{\otimes}$
(which is well-defined up to a contractible space of choices), we can reduce to the case
where $\calC = \SSet$, endowed with the Cartesian symmetric monoidal structure. In this case,
$\CMod_{\calC}( \LPress)$ is equivalent to the $\infty$-category $\LPress$ of symmetric monoidal
$\infty$-categories, and the tensor product $\calC \otimes S^1$ can be identified with
the $\infty$-category $( \SSet)_{/ S^1}$ of spaces lying over the circle. In this case, we can take
$K = \Z \simeq \Omega(S^1) \in \Mon( \SSet) \simeq \Alg(\SSet)$ to be the group of integers:
the equivalence $\SSet_{ / S^1} \simeq \Alg_{ K}( \SSet)$ is provided by Remark \ref{monoidchar},
and the free module functor $\SSet \rightarrow \Alg_{K}(\SSet)$ corresponds to the map given by the base point on $S^1$ by virtue of Remark \ref{monoidchur}.
\end{proof}

\begin{warning}
The spaces $X(A,B)$ constructed in the proof of Theorem \ref{machus} depend on the regular cardinals
$\kappa$ that are chosen at each stage of the construction. We can eliminate this dependence by
replacing the functor $X$ by the essential image of the natural transformation $\alpha: \bHom \rightarrow X$ at each step.
\end{warning}

\begin{remark}
With a bit more effort, one can show that the homotopy equivalence $\Omega^{k} X(A,B) \simeq \Centt{\OpE{k}}{f}^{\times}$ appearing in Theorem \ref{machus} is an equivalence of $k$-fold loop spaces, which depends functorially on $A$ and $B$.
\end{remark}

\subsection{The Cotangent Complex of an $\OpE{k}$-Algebra}\label{cotan}

Let $k \rightarrow A$ be a map of commutative rings. The multiplication map
$A \otimes_{k} A \rightarrow A$ is a surjection whose kernel is an ideal
$I \subseteq A \otimes_{k} A$. The quotient $I/ I^2$ is an $A$-module, and there is
a canonical $k$-linear derivation $d: A \rightarrow I/ I^2$, which carries an element $a \in A$ to the image
of $(a \otimes 1 - 1 \otimes a) \in I$. In fact, this derivation is {\em universal}: for any
$A$-module $M$, composition with $d$ induces a bijection
$\Hom_{A}( I/I^2, M) \rightarrow \Der_{k}(A,M).$ In other words, the quotient $I/I^2$ can be identified with the module of \Kahler differentials $\Omega_{A/k}$.

The above analysis generalizes in a straightforward way to the setting of associative algebras.
Assume that $k$ is a commutative ring and that $A$ is an associative $k$-algebra.
Let $M$ be an $A$-bimodule (in the category of $k$-modules: that is, we require
$\lambda m = m \lambda$ for $m \in M$ and $\lambda \in k$). 
A {\it $k$-linear derivation} from $A$ into $M$ is a $k$-linear map
$d: A \rightarrow M$ satisfying the Leibniz formula $d(ab) = d(a) b + a d(b)$.
If we let $I$ denote the kernel of the multiplication map $A \otimes_{k} A \rightarrow A$, then
$I$ has the structure of an $A$-bimodule, and the formula $d(a) = a \otimes 1 - 1 \otimes a$
defines a derivation from $A$ into $M$. This derivation is again universal in the following sense:
\begin{itemize}
\item[$(\ast)$] For any bimodule $M$, composition with $d$ induces a bijection of $\Hom( I, M)$ with the set of $k$-linear derivations from $I$ into $M$. 
\end{itemize}
If $A$ is commutative, then $I/I^2$ is the universal $A$-module map which receives
an $A$-bimodule homomorphism from $I$. Consequently, $(\ast)$ can be regarded as a generalization
of the formula $\Omega_{A/k} \simeq I/I^2$. 

Our goal in this section is to obtain an $\infty$-categorical analogue of assertion $(\ast)$. Rather than than working in the ordinary abelian category of $k$-modules, we will work with a symmetric monoidal stable $\infty$-category $\calC$. In this case, we can consider
algebra objects $A \in \Alg_{\calO}(\calC)$ for any coherent $\infty$-operad $\calO^{\otimes}$.
According to Theorem \deformationref{scummer2}, we can identify the stabilization
$\Stab( \Alg_{\calO}(\calC)_{/A})$ with the stable $\infty$-category $\Mod_{A}^{\calO}(\calC)$ of
$\calO$-algebra objects of $\calC$. In particular, the absolute cotangent complex $L_{A}$ can
be identified an object of $\Mod_{A}^{\calO}(\calC)$. Our goal is to obtain a concrete description of
$L_{A}$ in the special case where $\calO^{\otimes} = \OpE{k}$ is the $\infty$-operad of little $k$-cubes.

To motivate the description, let us consider first the case where $k=1$. In this case,
we can identify $\Mod_{A}^{\calO}(\calC)$ with the $\infty$-category of $A$-bimodule objects
of $\calC$ (see \S \symmetricref{bimid}). Motivated by assertion $(\ast)$, we might suppose that
$L_{A}$ can be identified with the fiber of the multiplication map
$A \otimes A \rightarrow A$ (regarded as a map of $A$-bimodules). The domain of this map
is the {\em free} $A$-bimodule, characterized up to equivalence by the existence of a morphism
$e: {\bf 1} \rightarrow A \otimes A$ with the property that it induces homotopy equivalences
$$ \bHom_{ \Mod^{\calO}_{A}(\calC)}( A \otimes A, M) \rightarrow \bHom_{ \calC}( {\bf 1}, M)$$
(here and in what follows, we will identify $A$-module objects of $\calC$ with their images in $\calC$).

Assume now that $k \geq 0$ is arbitrary, that $\calC$ is presentable, and that
the tensor product on $\calC$ preserves colimits separately in each variable.
The forgetful functor $\Mod^{\calO}_{A}(\calC) \rightarrow \calC$
preserves small limits and colimits (Corollaries \symmetricref{sus} and \symmetricref{sabus}),
and therefore admits a left adjoint $F: \calC \rightarrow \Mod^{\calO}_{A}(\calC)$
(Corollary \toposref{adjointfunctor}). We can formulate our main result as follows:

\begin{theorem}\label{curtis}
Let $\calC^{\otimes}$ be a stable symmetric monoidal $\infty$-category and let $k \geq 0$. Assume that $\calC$ is presentable and that the tensor product operation on $\calC$ preserves colimits separately in each variable. For every $\OpE{k}$-algebra object $A \in \Alg_{\OpE{k}}(\calC)$, there is a
canonical fiber sequence
$$ F( {\bf 1}) \rightarrow A \rightarrow L_{A}[k]$$
in the stable $\infty$-category $\Mod_{A}^{\calO}(\calC)$. Here $F: \calC \rightarrow \Mod_{A}^{\calO}(\calC)$ denotes the free functor described above, and the map of $A$-modules $F({ \bf 1}) \rightarrow A$ is determines by the unit map ${\bf 1} \rightarrow A$ in the $\infty$-category $\calC$.
\end{theorem}

\begin{remark}
A version of Theorem \ref{curtis} is proven in \cite{johnthesis}.
\end{remark}

\begin{remark}
If $A$ is an $\OpE{k}$-algebra object of $\calC$, then we can think of an $A$-module
$M \in \Mod^{\OpE{k}}_{A}(\calC)$ as an object of $\calC$ equipped with a commuting family of
(left) actions of $A$ parametrized by the $(k-1)$-sphere of rays in the Euclidean space $\R^{k}$ which emanate from the origin. This is equivalent to the action of a single associative algebra object of $\calC$: namely, the topological chiral homology $\int_{ S^{k-1} } A$ (see the discussion at the end of \S \ref{klimb}). The free module $F( {\bf 1})$ can be identified with $\int_{ S^{k-1}} A$ itself.

An equivalent formulation of Theorem \ref{curtis} asserts the existence of a fiber sequence of
$A$-modules
$$ L_{A} \stackrel{\theta}{\rightarrow} \Omega^{k-1} F( {\bf 1}) \stackrel{\theta'}{\rightarrow} \Omega^{k-1} A.$$
In particular, the map $\theta$ classifies a derivation $d$ of $A$ into
$\Omega^{k-1} F( {\bf 1})$. Informally, this derivation is determined by pairing the
canonical $S^{k-1}$-parameter family of maps $A \rightarrow \int_{S^{k-1}} A$ with
the fundamental class of $S^{k-1}$. Because the induced family of composite maps
$A \rightarrow \int_{S^{k-1}} A \rightarrow A$ is constant, this derivation lands in the fiber of the map $\theta'$. When $k=1$, we can identify $F({\bf 1})$ with the tensor product $A \otimes A$, and our heuristic recovers the classical formula $d(a) = a \otimes 1 - 1 \otimes a$.
\end{remark}

\begin{remark}
Our formulation of Theorem \ref{curtis} is designed to emphasize the maximal amount of symmetry.
The shift $L_{A}[k]$ can be identified with the tensor product of $L_{A}$ with the pointed space
$S^{k}$, regarded as the one-point compactification of the Euclidean space $\R^{k}$. With respect to this identification, the fiber sequence of Theorem \ref{curtis} can be constructed so as to be equivariant
with respect to the group of self-homeomorphisms of $\R^{k}$ (which acts on the $\infty$-operad
$\OpE{k}$ up to coherent homotopy, as explained in \S \ref{colee}). However, this equivariance is not apparent from the construction we present below.
\end{remark}

We now explain how to deduce Theorem \ref{curtis} from Theorem \ref{machus}. 
Fix an $\OpE{k}$-algebra $A \in \Alg_{\OpE{k}}(\calC)$, and let
$\calE = \Alg_{\OpE{k}}(\calC)_{A/}$. Consider the functors $X,Y,Z: \calE \rightarrow \SSet_{\ast}$
given informally by the formulas
$$ X(f: A \rightarrow B) = \Omega^{n} \bHom_{\Alg_{\OpE{k}}}( A, B)
\quad Y(f: A \rightarrow B) = \Centt{\OpE{k}}{f}^{\times} \quad Z(f: A \rightarrow B) = B^{\times}.$$
Theorem \ref{machus} implies that these functors fit into a pullback diagram
$$ \xymatrix{ X \ar[r] \ar[d] & Y \ar[d] \\
\ast \ar[r] & Z, }$$
where $\ast: \calE \rightarrow \SSet_{\ast}$ is the constant diagram taking the value $\ast$.
(In fact, we have a pullback diagram in the $\infty$-category of functors from $\calE$ to the $\infty$-category $\Mon_{\OpE{k}}(\SSet)$ of $\OpE{k}$-spaces, but we will not need this).

Let $\calE' = \Alg_{ \OpE{k}}(\calC)_{A/ \, /A}$. Let $X': \calE' \rightarrow \SSet_{\ast}$ be the functor which assigns to a diagram $$ \xymatrix{ & B \ar[dr] & \\
A \ar[ur]^{f} \ar[rr]^{\id_A}  & & A }$$
the fiber of the induced map $X(f) \rightarrow X(\id_A)$, and let $Y'$ and $Z'$ be defined similarly.
Using Lemma \toposref{limitscommute}, we deduce the existence of a pullback diagram of functors
$$ \xymatrix{ X' \ar[r] \ar[d] & Y' \ar[d] \\
\ast \ar[r] & Z'. }$$

Let $\phi: \Mod_{A}^{\calO}(\calC) \rightarrow \Alg_{\OpE{k}}(\calC)_{A/ \, /A}$ be the functor
given informally by the formula $M \mapsto A \oplus M$ (that is, $\phi$ is the composition
of the identification $\Mod_{A}^{\calO}(\calC) \simeq \Stab( \Alg_{\OpE{k}}(\calC)_{A/ \, /A})$
with the functor $\Omega^{\infty}: \Stab( \Alg_{\OpE{k}}(\calC)_{A/ \, /A}) \rightarrow
\Alg_{\OpE{k}}(\calC)_{A/ \, /A}$). Let $X'' = X' \circ \phi$, and define $Y''$ and $Z''$ similarly.
We have a pullback diagram of functors
$$ \xymatrix{ X'' \ar[r] \ar[d] & Y'' \ar[d] \\
\ast \ar[r] & Z'' }$$
from $\Mod_{A}^{\calO}(\calC)$ to $\SSet_{\ast}$. 

The functor $Z''$ carries an $A$-module $M$ to the fiber of the map
$(A \oplus M)^{\times} \rightarrow A^{\times}$, which can be identified with
$\bHom_{\calC}( {\bf 1}, M) \simeq \bHom_{ \Mod^{\OpE{k}}_{A}(\calC)}( F( {\bf 1}), M)$.
In other words, the functor $Z''$ is corepresentable by the object $F( {\bf 1}) \in \Mod^{\OpE{k}}_{A}(\calC)$. Similarly, Theorem \ref{serpses} implies that the functor $Y''$ is corepresentable by
the object $A \in \Mod^{\OpE{k}}_{A}(\calC)$. By definition, the functor $X''$ is corepresentable
by the shifted cotangent complex $L_A[k]$. Since the Yoneda embedding for
$\Mod^{\OpE{k}}_{A}(\calC)$ is fully faithful, we deduce the existence of a commutative diagram of
representing objects
$$ \xymatrix{ L_A[k] & A \ar[l] \\
0 \ar[u] & F( {\bf 1}) \ar[l] \ar[u] }$$
which is evidently a pushout square. This yields the desired fiber sequence
$$ F( { \bf 1}) \rightarrow A \rightarrow L_A[k].$$
in $\Mod^{\OpE{k}}_{A}(\calC)$. 

\begin{remark}
The fiber sequence of Theorem \ref{curtis} depends functorially on $A$ (this follows from a more careful version of the construction above). We leave the details of the formulation to the reader's imagination.
\end{remark}

\begin{remark}
Let $A$ be a commutative algebra object of $\calC$. Then $A$ can be regarded as an
$\OpE{k}$-algebra object of $\calC$ for every nonnegative integer $k$. When regarded
as an $\OpE{k}$-algebra object, $A$ has a cotangent complex which we will denote
by $L_{A}^{(k)}$, to emphasize the dependence on $k$. The topological chiral homology (see \S \ref{sec3})
$\int_{ S^{k-1}} A$ can be identified with the tensor product $A \otimes S^{k-1}$ (Theorem \ref{hunefer}),
which is the $(k-1)$-fold (unreduced) suspension $\Sigma^{k-1} A$ of $A$, regarded as an object of $\Alg_{\Comm}(\calC)_{/A}$. According to Theorem \ref{curtis}, we have a canonical identification
$L_{A}^{(k)} \simeq (\ker \Omega^{k-1} \Sigma^{k-1}(A) \rightarrow A)$ in the
$\infty$-category $\calC$. Since the $\infty$-operad $\Comm$ is equivalent to the colimit of the
$\infty$-operads $\OpE{k}$ (see Corollary \ref{infec}), we conclude that the {\em commutative}
algebra cotangent complex $L_{A}$ can be computed as the colimit
$\varinjlim_{k} L_{A}^{(k)}$. Combining this observation with the above identification, we obtain
an alternative ``derivation'' of the formula $\Omega^{\infty} \Sigma^{\infty} \simeq \varinjlim_{k} \Omega^{k} \Sigma^{k}$.
\end{remark}

\begin{example}\label{emore}
Let $\calC$ be as in Theorem \ref{curtis} and let $\epsilon: A \rightarrow {\bf 1}$ be
an augmented $\OpE{k}$-algebra object of $\calC$ (see Example \ref{kd}). Theorem
\ref{serpses} guarantees the existence of a Koszul dual $A^{\vee} = \Centt{\OpE{k}}{\epsilon}$. Moreover,
as an object of the underlying $\infty$-category $\calC$, $A^{\vee}$ can be identified
with a morphism object $\Mor_{ \Mod^{\OpE{k}}_{A}(\calC)}( A, { \bf 1} )$. Combining this observation with
the fiber sequence of Theorem \ref{curtis} (and observing that the morphism object
$\Mor_{ \Mod^{\OpE{k}}_{A}(\calC)}( F({\bf 1}), {\bf 1})$ is equivalent to ${\bf 1}$), we obtain a fiber sequence
$$ \Mor_{ \Mod^{\OpE{k}}_{A}(\calC)}( L_{A}[k], {\bf 1}) \rightarrow A^{\vee} \stackrel{\theta}{\rightarrow}{\bf 1}$$
in $\calC$. The map $\theta$ underlies the augmentation on $A^{\vee}$ described in
Example \ref{kd}; we may therefore view $\Mor_{ \Mod^{\OpE{k}}_{A}(\calC)}( L_A[k], {\bf 1})$ as 
the ``augmentation ideal'' of the Koszul dual $A^{\vee}$.

In heuristic terms, we can view the $\OpE{k}$-algebra $A$ as determining a ``noncommutative scheme'' $\Spec{A}$, which is equipped with a point given by the augmentation $\epsilon$. We can think
of $L_{A}$ as a version of the cotangent bundle of $\Spec{A}$, and
$\Mor_{\Mod^{\OpE{k}}_{A}(\calC)}( L_A, {\bf 1})$ as a version of the tangent space to $\Spec{A}$
at the point determined by $\epsilon$. The above analysis shows that, up to a shift by $k$, this
``tangent space'' itself is the augmentation ideal in a {\em different} augmented $\OpE{k}$-algebra object of $\calC$ (namely, the Koszul dual algebra $A^{\vee}$).
\end{example}

\section{Factorizable Sheaves}\label{sec3}

Fix an integer $k \geq 0$. In \S \ref{founder}, we introduced the $\infty$-operad
$\OpE{k}$ of little $k$-cubes. The underlying $\infty$-category of $\OpE{k}$ has a unique object, which we
can think of as an abstract open cube $\Cube{k}$ of dimension $k$. There are a number of variations on this theme, where we consider cubes (or, equivalently, open disks) endowed with additional structures of various types. In \S \ref{colee}, we will review some of these variations and study their relationship with one another.
For our purposes, the main case of interest is that in which we require all of our cubes to be equipped with an open embedding into a topological manifold $M$ of dimension $k$. The collection of such cubes can be organized into an $\infty$-operad, which we will denote by $\OpE{M}$; we refer the reader to \ref{hullby} for a precise definition.

Roughly speaking, we can think of an $\OpE{M}$-algebra object of a symmetric monoidal $\infty$-category
$\calC^{\otimes}$ as a family of $\OpE{k}$-algebras $A_x$ parametrized by the points $x \in M$ (more accurately, one should think of this family as ``twisted'' by the tangent bundle of $M$: that is, for every
point $x \in M$ we should think of $A_x$ as an algebra over an $\infty$-operad whose objects are little disks in the tangent space $T_{M,x}$ to $M$ at $x$). There is a convenient geometric way to encode this information.
Following Beilinson and Drinfeld (see \cite{beilinson}), we define the {\it Ran space} $\Ran(M)$ of $M$ to be the collection of all nonempty finite subsets of $M$ (for a more detailed discussion of $\Ran(M)$, together with a definition of $\Ran(M)$ as a topological space, we refer the reader to \ref{secran}). To every point $S \in \Ran(M)$, the tensor product $A_{S} = \bigotimes_{s \in S} A_s$ is an object of $\calC$. We will see that the objects are the stalks of a $\calC$-valued cosheaf $\calF$ on the Ran space. We can regard $\calF$ as a constructible cosheaf which is obtained by gluing together locally constant cosheaves along the locally closet subsets $\Ran^{n}(M) = \{ S \in \Ran(M): |S| = n \} \subseteq \Ran(M)$ for $n \geq 1$; the ``gluing'' data for these restrictions reflects the multiplicative structure of the algebras $\{ A_x \}_{x \in M}$. In \S \ref{stupem}, we will
see that the construction $A \mapsto \calF$ determines an equivalence of $\infty$-categories
from the $\infty$-category of (nonunital) $\OpE{M}$-algebras in $\calC$ to a suitable $\infty$-category of {\em factorizable} $\calC$-valued cosheaves on $\Ran(M)$,
which are constructible with respect to the above stratification (Theorem \ref{manust}). 

The description of an $\OpE{M}$-algebra object $A$ of $\calC$ as a factorizable $\calC$-valued cosheaf $\calF$ on $\Ran(M)$ suggests an interesting invariant of $A$: namely, the object $\calF( \Ran(M) ) \in \calC$ given by global sections of $\calF$. In the case where $M$ is connected, we will refer to the global sections
$\calF( \Ran(M) )$ as the {\it topological chiral homology} of $M$ with coefficients in $A$, which we will denote by $\int_{M} A$. We will give an independent definition of $\int_{M} A$ (which does not require the assumption that $M$ is connected) in \S \ref{tchtch}, and verify that it is equivalent to $\calF( \Ran(M) )$ for connected $M$ in
\S \ref{stupem} (Theorem \ref{selfmaid}). The construction $A \mapsto \int_{M} A$ can be regarded as
a generalization of Hochschild homology (Theorem \ref{junl}) and has a number of excellent formal properties, which we will verify in \S \ref{klimb}. In \S \ref{nonpn}, will use the theory of topological chiral homology to formulate and prove a nonabelian version of the Poincare duality theorem (Theorem \ref{stager}). The proof
relies on a technical compatibility result between fiber products and sifted colimits, which we verify
in \S \ref{digdig}.

\begin{convention}\label{koson}
Unless otherwise specified, the word {\it manifold} will refer to a paracompact Hausdorff topological manifold of some fixed dimension $k$.
\end{convention}

\subsection{Variations on the Little Cubes Operads}\label{colee}

Fix an integer $k \geq 0$. In \S \ref{kloo}, we introduced a topological operad $\TopE{k}$ whose $n$-ary
operations are given by {\em rectilinear} open embeddings from $\Cube{k} \times \nostar{n}$ into
$\Cube{k}$. Our goal in this section is to introduce some variations on this construction, where we
drop the requirement that our embeddings be rectilinear (or replace the condition of rectilinearity by some other condition). The main observation is that the resulting $\infty$-operads are closely related to
the $\infty$-operad $\OpE{k}$ studied in \S \ref{founder} (see Proposition \ref{luge} below).  

\begin{notation}
If $M$ and $N$ are manifolds of the same dimension, we let $\Emb(M, N)$ denote the space
of all open embeddings from $M$ into $N$ (for a more detailed discussion of these embedding spaces, we refer the reader to \S \ref{cowpy}).
\end{notation}

\begin{definition}\label{defF}
Fix an integer $k \geq 0$. We define a topological category
$\TopE{ \BTop(k) }$ as follows:
\begin{itemize}
\item[$(1)$] The objects of $\TopE{ \BTop(k) }$ are the objects $\seg{n} \in \FinSeg$.
\item[$(2)$] Given a pair of objects $\seg{m}, \seg{n} \in \TopE{ \BTop(k) }$, the mapping
space $\bHom_{ \TopE{ \BTop(k) }}( \seg{m}, \seg{n})$ is given by the disjoint union
$$ \coprod_{ \alpha} \prod_{ 1 \leq i \leq n} \Emb( \R^{k} \times \alpha^{-1} \{i\}, \R^{k} )$$
taken over all morphisms $\alpha: \seg{m} \rightarrow \seg{n}$ in $\FinSeg$.
\end{itemize}
We let $\OpE{ \BTop(k) }$ denote the $\infty$-category given by the topological nerve
$\Nerve( \TopE{ \BTop(k) })$.
\end{definition}

\begin{remark}
It follows from Proposition \symmetricref{calp} that $\OpE{ \BTop(k) }$ is an $\infty$-operad.
\end{remark}

\begin{remark}
Definition \ref{defF} is a close relative of Definition \ref{defE}. In fact, choosing a homeomorphism
$\R^{k} \simeq \Cube{k}$, we obtain an inclusion of $\infty$-operads
$\OpE{k} \rightarrow \OpE{ \BTop(k) }$.
\end{remark}

\begin{remark}
The object $\seg{0}$ is initial in $\TopE{ \BTop(k) }$. It follows that $\OpE{ \BTop(k) }$ is a unital $\infty$-operad.
\end{remark}

\begin{example}\label{spant}
Suppose that $k = 1$. Every open embedding $j: \Cube{k} \times S \hookrightarrow \Cube{k}$
determines a pair $( <, \epsilon)$, where $<$ is an element of the set of
linear orderings of $S$ (given by $s < s'$ if $j(0,s) < j(0,s')$) and
$\epsilon: S \rightarrow \{ \pm 1 \}$ is a function defined so that
$\epsilon(s) = 1$ if $j | \Cube{k} \times \{s\}$ is orientation preserving, and $\epsilon(s) = -1$
otherwise. This construction determines a homotopy equivalence
$\Emb( \Cube{k} \times S, \Cube{k}) \rightarrow L(S) \times \{ \pm 1 \}^{S}$, where
$L(S)$ denotes the set of linear orderings of $S$. It follows that
$\OpE{ \BTop(k) }$ is equivalent to the nerve of its homotopy category and therefore arises from an operad in the category of sets via Construction \symmetricref{coco}. In fact, this is the operad which controls
associative algebras with involution, as described in \S \ref{justac}.
\end{example}

\begin{definition}
For each integer $k \geq 0$, we let $\BTop(k)$ denote the fiber product
$\OpE{ \BTop(k) } \times_{ \Nerve( \FinSeg) } \{ \seg{1} \}$. Then $\BTop(k)$ can be identified
with the nerve of the topological category having a single object whose endomorphism monoid
is the space $\Emb( \R^{k}, \R^{k} )$ of open embeddings from $\R^{k}$ to itself. It follows
from the Kister-Mazur theorem (Theorem \ref{scen}) that $\Emb( \R^{k}, \R^{k})$ is a grouplike topological monoid, so that $\BTop(k)$ is a Kan complex. In fact, Theorem \ref{scen} shows that $\BTop(k)$ can be identified with a classifying space for the topological group $\Top(k)$ of homeomorphisms from $\R^{k}$ to itself.
\end{definition}

\begin{remark}\label{kullman}
We can modify Definition \ref{defF} by replacing the embedding spaces
$\Emb( \R^{k} \times S, \R^{k})$ by the products $\prod_{s \in S} \Emb( \R^{k}, \R^{k})$. 
This yields another $\infty$-operad, which is canonically isomorphic to
$\BTop(k)^{\amalg}$. The evident inclusions 
$\Emb( \R^{k} \times S, \R^{k}) \hookrightarrow \prod_{s \in S} \Emb( \R^{k}, \R^{k})$
induce an inclusion of $\infty$-operads $$\TopE{ \BTop(k) } \hookrightarrow \BTop(k)^{\amalg}.$$
\end{remark}

If $k > 0$, then the Kan complex $\BTop(k)$ is not contractible (nor even simply-connected, since an orientation-reversing homeomorphisms from $\R^{k}$ to itself cannot be isotopic to the identity), so the
$\infty$-operad $\TopE{ \BTop(k) }$ is not reduced. Consequently, we can apply Theorem \ref{assum}
to decompose $\TopE{ \BTop(k) }$ as the assembly of a family of reduced $\infty$-operads. The key to understanding this decomposition is the following observation:

\begin{proposition}\label{luge}
Let $k$ be a nonnegative integer, and choose a homeomorphism $\R^{k} \simeq \Cube{k}$. 
The induced inclusion $f: \OpE{k} \rightarrow \OpE{ \BTop(k) }$ is an ornamental map of $\infty$-operads (see \S \ref{justac}). 
\end{proposition}

\begin{proof}
We will employ the notation
introduced in \S \ref{cowpy}. It will suffice to show that $f$ satisfies criterion $(3)$ of Lemma \ref{snoke}.
Unwinding
the definitions, we are reduced to showing that for every finite set $S$, the diagram
$$ \xymatrix{ \Sing(\Rect( \Cube{k} \times S, \Cube{k} )) \ar[r] \ar[d] & (\Sing \Rect( \Cube{k}, \Cube{k}))^{S} \ar[d] \\
\Sing( \Emb( \R^{k} \times S), \R^k) \ar[r] & \Sing( \Emb( \R^k, \R^k) )^{S} }$$
is a homotopy pullback square of Kan complexes. Consider the larger diagram
$$ \xymatrix{ \Sing(\Rect( \Cube{k} \times S, \Cube{k} )) \ar[r] \ar[d] & (\Sing \Rect( \Cube{k}, \Cube{k}))^{S} \ar[d] \\
\Sing( \Emb( \R^k \times S, \R^k)) \ar[r] \ar[d] & \Sing( \Emb( \R^{k}, \R^k) )^{S} \ar[d] \\
\Germ(S, \R^{k}) \ar[r] \ar[d] & \prod_{s \in S} \Germ( \{s\}, \R^{k}) \ar[d] \\
\Conf( S, \R^{k}) \ar[r] & \prod_{s \in S} \Conf( \{s\}, \R^{k}). }$$
The lower square is a pullback diagram in which the vertical maps are Kan fibrations, and
therefore a homotopy pullback diagram. The middle square is a homotopy pullback diagram
because the middle vertical maps are homotopy equivalences (Proposition \ref{splay}). The outer
rectangle is a homotopy pullback diagram because the vertical compositions are homotopy
equivalences (Remark \ref{sove}). The desired result now follows from a diagram chase.
\end{proof}

\begin{remark}\label{cola}
Fix a nonnegative integer $k$. The $\infty$-operad $\OpE{ \BTop(k) }$ is unital and its
underlying $\infty$-category is a Kan complex $\BTop(k)$. According to Theorem \ref{assum}, there exists a reduced family of $\infty$-operads $\calO^{\otimes}$ and an assembly map
$\calO^{\otimes} \rightarrow \OpE{ \BTop(k) }$. Then $\calO^{\otimes}_{\seg{0}} \simeq
\calO \simeq \BTop(k)$; we may therefore assume without loss of generality that
$\calO^{\otimes} \rightarrow \BTop(k) \times \Nerve(\FinSeg)$ is a $\BTop(k)$-family of $\infty$-operads.
Since $\OpE{k}$ is reduced, Theorem \ref{assum} guarantees that the inclusion $\OpE{k} \rightarrow \OpE{\BTop(k)}$ factors (up to homotopy) through $\calO^{\otimes}$. Without loss of generality,
this map factors through $\calO^{\otimes}_{x}$ for some vertex $x \in \BTop(k)$. The resulting map
$\OpE{k} \rightarrow \calO^{\otimes}_{x}$ is an ornamental map between reduced $\infty$-operads
(Proposition \ref{luge}), and therefore an equivalence (Lemma \ref{stoag}). We
can summarize the situation as follows: the $\infty$-operad $\OpE{ \BTop(k) }$ is obtained by assembling
a reduced $\BTop(k)$-family of $\infty$-operads, each of which is equivalent to $\OpE{k}$. More informally,
we can regard this $\BTop(k)$-family as encoding an action of the loop space $\Omega \BTop(k) \simeq
\Top(k)$ on the $\infty$-operad $\OpE{k}$, so that $\OpE{ \BTop(k) }$ can be regarded as a semidirect product of $\OpE{k}$ by the action of the topological group $\Top(k)$ of homeomorphisms of
$\R^{k}$ with itself. \end{remark}

We can summarize Remark \ref{cola} informally as follows:
if $\calC^{\otimes}$ is a symmetric monoidal $\infty$-category, then the $\infty$-category
$\Alg_{ \OpE{ \BTop(k)}}( \calC)$ can be identified with the $\infty$-category of
$\OpE{k}$-algebra objects of $\calC$ which are equipped with a compatible action of the
topological group $\Top(k)$. The requirement that $\Top(k)$ act on an $\OpE{k}$-algebra is rather strong: in practice, we often encounter situations where an algebra $A \in \Alg_{ \OpE{k}}( \calC)$ is acted on not by the whole of $\Top(k)$, but by some subgroup. Our next definition gives a convenient formulation of this situation.

\begin{definition}\label{custa2}
Let $B$ be a Kan complex equipped with a Kan fibration $B \rightarrow \BTop(k)$. 
We let $\OpE{B}$ denote the fiber product
$$ \OpE{ \BTop(k) } \times_{ \BTop(k)^{\amalg} } B^{\amalg}.$$
\end{definition}

\begin{remark}
It follows immediately from the definitions that $\OpE{B}$ is a unital $\infty$-operad, equipped with
an ornamental map $\OpE{B} \rightarrow \OpE{ \BTop(k) }$. 
\end{remark}

\begin{warning}
Our notation is slightly abusive. The $\infty$-operad $\OpE{B}$ depends not only on the Kan complex $B$, but also the integer $k$ and the map $\theta: B \rightarrow \BTop(k)$. We can think of $\theta$ as classifying a fiber
bundle over the geometric realization $|B|$, whose fibers are Euclidean spaces.
\end{warning}

\begin{remark}\label{bluw}
Let $\calO^{\otimes} \rightarrow \BTop(k) \times \Nerve(\FinSeg)$ be the $\infty$-operad family of Remark \ref{cola}. If $\theta: B \rightarrow \BTop(k)$ is any map of Kan complexes, then the fiber product
$\calO^{\otimes} \times_{ \BTop(k)} B$ is a $B$-family of reduced unital $\infty$-operads. When
$\theta$ is a Kan fibration (which we may assume without loss of generality), then
this $B$-family of $\infty$-operads assembles to the unital $\infty$-operad $\OpE{B}$
(see \S \ref{justac}). We can informally describe the situation as follows: an $\OpE{B}$-algebra object
of a symmetric monoidal $\infty$-category $\calC$ is a (twisted) family of $\OpE{k}$-algebra objects of $\calC$,
parametrized by Kan complex $B$ (the nature of the twisting is specified by the map $\theta$). 
\end{remark}

We conclude this section by illustrating Definition \ref{custa2} with some examples. Another general class of examples will be discussed in \S \ref{hullby}.

\begin{example}\label{custa}
Let $B$ be a contractible Kan complex equipped with a Kan fibration $B \rightarrow \BTop(k)$. Then
$\OpE{B}$ is equivalent to the $\infty$-operad $\OpE{k}$.
\end{example}

\begin{example}\label{kample}
Fix $k \geq 0$, and choose a homeomorphism of $\R^{k}$ with the unit ball $B(1) \subseteq \R^{k}$.
We will say that a map $f: B(1) \rightarrow B(1)$ is a {\it projective isometry}
if there exists an element $\gamma$ in the orthogonal group $\OO(k)$, a positive
real number $\lambda$, and a vector $v_0 \in B(1)$ such that $f$ is given by the formula
$f(w) = v_0 + \lambda \gamma(w)$. For every finite set $S$, we let
$\Isom^{+}( B(1) \times S, B(1))$ denote the (closed) subspace of $\Emb( B(1) \times S, B(1))$
consisting of those open embeddings whose restriction to each ball $B(1) \times \{s\}$ is 
an orientation-preserving projective isometry. Let $\TopE{\SO(k)}$ be the subcategory of $\TopE{ \BTop(k) }$ having the same objects, with morphism spaces given by
$$ \bHom_{ \TopE{\SO(k)}}( \seg{m}, \seg{n}) = \coprod_{ \alpha} \prod_{ 1 \leq i \leq n}
\Isom^{+}( B(1) \times \alpha^{-1} \{i\}, B(1) ).$$
Then $\calO^{\otimes} = \Nerve( \TopE{\SO(k)})$ is a unital $\infty$-operad. 
The inclusion $\calO^{\otimes} \hookrightarrow \OpE{ \BTop(k) }$ is an ornamental map which
induces an identification of $\calO^{\otimes}$ with the $\infty$-operad $\OpE{B}$, where $B$ is a Kan complex which plays the role of a classifying space $\BSO(k)$ for the special orthogonal group
$\SO(k)$ (and we arrange that the inclusion of topological groups $\SO(k) \rightarrow \Top(k)$ induces a
Kan fibration $\BSO(k) \rightarrow \BTop(k)$). This recovers the operad of {\em framed disks} described, for example, in \cite{salvatorewahl}. 
\end{example}

\begin{variant}\label{klam}
In Example \ref{kample}, there is no need to restrict our attention to orientation preserving
maps. If we instead allow {\em all} projective isometries, then we get another $\infty$-operad
$\calO^{\otimes} \simeq \OpE{B}$, where $B$ is a classifying space for the orthogonal group $\OO(k)$.
\end{variant}

\begin{example}\label{kample2}
In the definition of $\TopE{ \BTop(k) }$, we have allowed arbitrary open embeddings between
Euclidean spaces $\R^k$. We could instead restrict our attention to spaces of {\em smooth} open embeddings (which we regard as equipped with the Whitney topology, where convergence is given by uniform convergence of all derivatives on compact sets) to obtain an $\infty$-operad $\OpE{\Smooth}$. This can be identified with the $\infty$-operad $\OpE{B}$, where
$B$ is a classifying space for the monoid of smooth embeddings from the open ball $B(1)$ to itself.
Since every projective isometry is smooth, there is an obvious map $\calO^{\otimes} \rightarrow \OpE{\Smooth}$, where $\calO^{\otimes}$ is defined as in Variant \ref{klam}. In fact, this
map is an equivalence of $\infty$-operads: this follows from the fact that the inclusion from
the orthogonal group $\OO(k)$ into the space $\Emb^{\smooth}( B(1), B(1))$ of smooth embeddings of
$B(1)$ to itself is a homotopy equivalence (it has a homotopy inverse given by the composition
$\Emb^{\smooth}(B(1), B(1)) \rightarrow \GL_{k}( \R) \rightarrow \OO(k),$
where the first map is given by taking the derivative at the origin and the second is a homotopy
inverse to the inclusion $\OO(k) \hookrightarrow \GL_{k}(\R)$).
\end{example}

\begin{variant}\label{pklam}
In Example \ref{kample2}, we can use piecewise linear manifolds in place of smooth manifolds.
We can also consider manifolds which are equipped with additional structures, such as orientations. We leave the details to the reader.
\end{variant}

\subsection{Little Cubes in a Manifold}\label{hullby}

Let $M$ be a topological space equipped with an $\R^{k}$-bundle $\zeta \rightarrow M$. 
Assuming that $M$ is sufficiently nice, we can choose a Kan complex $B$ such that
$X$ is homotopy equivalent to the geometric realization $|B|$, and the bundle $\zeta$ is classified by a Kan fibration of simplicial sets $\theta: B \rightarrow \BTop(k)$. In this case, we can apply the construction of Definition \ref{custa2} to obtain an $\infty$-operad $\OpE{B}$. In the special case where $M$ is a topological manifold of dimension $k$ and $\zeta$ is the tangent bundle of $M$, we will denote this $\infty$-operad by $\OpE{M}$ (see Definition \ref{farsy} below for a precise definition). We can think of $\OpE{M}$ as a variation
on the $\infty$-operad $\OpE{k}$ whose objects are cubes $\Cube{k}$ equipped with an open embedding
into $M$, and whose morphisms are required to be compatible with these open embeddings (up to specified isotopy). We will also consider a more rigid version of the $\infty$-operad $\OpE{M}$, where the morphisms are required to be {\em strictly} compatible with the embeddings into $M$ (rather than merely up to isotopy);
this $\infty$-operad will be denoted by $\Nerve( \Disk{M})^{\otimes}$ (Definition \ref{kaper}). The main result of this section is Theorem \ref{sazz}, which asserts that theory of $\OpE{M}$-algebras is closely related to the 
more rigid theory of $\Nerve( \Disk{M})^{\otimes}$-algebras. 

Our first step is to describe the $\infty$-operad $\OpE{M}$ more precisely.

\begin{definition}\label{farsy}
Let $M$ be a topological manifold of dimension $k$. We define a topological category
$\calC_{M}$ having two objects, which we will denote by $M$ and $\R^{k}$, with
mapping spaces given by the formulas
$$ \bHom_{ \calC_M}( \R^{k}, \R^{k}) = \Emb( \R^{k}, \R^{k}) \quad \quad
\bHom_{ \calC_M}( \R^{k},M) = \Emb( \R^{k}, M) $$
$$ \bHom_{ \calC_{M}}( M, \R^{k}) = \emptyset \quad \quad \bHom_{ \calC_M}(M,M) = \{ \id_M \}.$$
We identify the Kan complex $\BTop(k)$ with a full subcategory of the nerve
$\Nerve(\calC_M)$. Let $B_M$ denote the Kan complex
$\BTop(k) \times_{ \Nerve( \calC_M)} \Nerve(\calC_M)_{/ M}$. We let
$\OpE{M}$ denote the $\infty$-operad $\OpE{ \BTop(k) } \times_{ \BTop(k)^{\amalg} } B_{M}^{\amalg}$
In other words, we let $\OpE{M}$ denote the $\infty$-operad $\OpE{B_M}$ introduced in Definition \ref{custa2}.
\end{definition}


\begin{remark}\label{talrod}
Let $M$ be a topological manifold of dimension $k$, and let $B_{M}$ be defined as in Definition \ref{farsy}. 
Then $\OpE{M}$ can be obtained as the assembly of a $B_M$-family of $\infty$-operads, each of which is
equivalent to $\OpE{k}$ (Remark \ref{bluw}). To justify our notation, we will show that the Kan complex
$B_M$ is canonically homotopy equivalent to the (singular complex of) $M$. More precisely, we
will construct a canonical chain of homotopy equivalences
$$ B_M \leftarrow B'_{M} \rightarrow B''_{M} \leftarrow \Sing(M).$$ 
To this end, we define topological categories $\calC'_M$ and $\calC''_{M}$, each of which
consists of a pair of objects 
$\{ \R^{k}, M \}$ with morphism spaces given by the formulas
$$ \bHom_{ \calC'_M}( \R^{k}, \R^{k}) = \Emb_0( \R^{k}, \R^{k}) \quad \quad \quad
\bHom_{ \calC'_M}( \R^{k},M) = \Emb( \R^{k}, M) $$
$$ \bHom_{\calC''_{M}}( \R^{k}, \R^{k}) = \{ 0 \} \quad \quad \quad \bHom_{ \calC''_{M}}( \R^{k}, M) = M $$
$$ \bHom_{ \calC'_{M}}( M, \R^{k}) = \emptyset = \bHom_{\calC''_M}(M, \R^{k})
\quad \quad \bHom_{ \calC'_M}(M,M) = \{ \id_M \} = \bHom_{\calC''_M}(M,M).$$
Here we let $\Emb_0( \R^{k}, \R^{k})$ denote the closed subset of $\Emb(\R^{k}, \R^{k} )$ spanned by those open embeddings $f: \R^{k} \rightarrow \R^{k}$ such that $f(0) = 0$. 

Let $\BTop'(k)$ denote the full subcategory of $\Nerve( \calC'_M)$ spanned by the object
$\R^{k}$, let $B'_{M}$ denote the fiber product
$\BTop(k) \times_{ \Nerve( \calC'_M)} \Nerve(\calC'_M)_{/M}$, and let
$B''_{M}$ denote the fiber product $\{ \R^{k} \} \times_{ \Nerve(\calC''_M)} \Nerve(\calC''_{M})_{/M}$.
We have maps of topological categories $\calC_M \stackrel{\theta}{\leftarrow} \calC'_{M} \stackrel{\theta''}{\rightarrow} \calC''_{M}$. The map $\theta$ is a weak equivalence of topological categories, and so induces a homotopy equivalence $B'_{M} \rightarrow B_{M}$. We claim that the induced map
$\psi: B'_{M} \rightarrow B''_{M}$ is also a homotopy equivalence. We can identify vertices
of $B'_{M}$ with open embeddings $\R^{k} \rightarrow M$ and vertices of $B''_{M}$ with points of
$M$; since $M$ is a $k$-manifold, the map $\psi$ is surjective on vertices. Fix a vertex
$(j: \R^{k} \hookrightarrow M) \in B'_{M}$. We have a map of homotopy fiber sequences
$$ \xymatrix{ \bHom_{ \Nerve(\calC'_{M})}( \R^{k}, \R^{k}) \ar[r] \ar[d] & \bHom_{ \Nerve(\calC'_M)}( \R^{k}, M) \ar[r]^-{\phi} \ar[d] & B'_M \ar[d] \\
\ast \ar[r] & \bHom_{ \Nerve(\calC''_{M})}( \R^{k}, M) \ar[r] & B''_M. }$$
It follows from Remark \ref{gec} that the left square is a homotopy pullback. It follows
that the map of path spaces $\bHom_{ B'_M}( j, j') \rightarrow \bHom_{ B''_{M}}( \psi(j), \psi(j') )$
is a homotopy equivalence for every $j'$ lying in the essential image of $\phi$. Since
the space $\BTop'(k)$ is connected, the map $\phi$ is essentially surjective, so that
$\psi$ is a homotopy equivalence as desired.

We note there is a canonical homotopy equivalence $\Sing(M) \rightarrow B''_M$
(adjoint to the weak homotopy equivalence appearing in Proposition \toposref{babyy}). 
Consequently, we obtain
a canonical isomorphism $B_{M} \simeq B'_{M} \simeq B''_{M} \simeq \Sing(M)$ in the homotopy
category $\calH$. It follows that
$\OpE{M}$ can be identified with the colimit of a family of $\infty$-operads parametrized by $M$,
each of which is equivalent to $\OpE{k}$. This family is generally not constant: instead, it is twisted by the principal $\Top(k)$-bundle given by the tangent bundle of $M$. In other words, if $\calO^{\otimes}$
is an $\infty$-operad, then we can think of an object of $\Alg_{ \OpE{M}}( \calO)$ as a family of
$\OpE{k}$-algebra objects of $\calO^{\otimes}$, parametrized by the points of $M$.
\end{remark}

\begin{example}
Let $M$ be the Euclidean space $\R^{k}$. Then the space $B_M$ is contractible, so that
$\OpE{M}$ is equivalent to the littles cubes operad $\OpE{k}$ (see Example \ref{custa}).
Since the $\infty$-operad $\OpE{M}$ depends functorially on $M$, we obtain another description of the ``action up to homotopy'' of the homeomorphism group $\Top(k)$ on $\OpE{k}$ (at least if we view $\Top(k)$ as a discrete group). 
\end{example}

We now introduce a more rigid variant of the $\infty$-operad $\OpE{M}$. 

\begin{definition}\label{kaper}
Let $M$ be a topological manifold of dimension $k$. Let $\Disk{M}$ denote the collection of
all open subsets $U \subseteq M$ which are homeomorphic to Euclidean space $\R^{k}$.
We regard $\Disk{M}$ as a partially ordered set (with respect to inclusions of open sets),
and let $\Nerve(\Disk{M})$ denote its nerve. Let $\Nerve( \Disk{M})^{\otimes}$ denote the subcategory
subset of $\Nerve( \Disk{M})^{\amalg}$ spanned by those morphisms
$( U_1, \ldots, U_m) \rightarrow (V_1, \ldots, V_n)$ with the following property:
for every pair of distinct integers $1 \leq i,j \leq m$ having the same image $k \in \nostar{n}$, the open subsets $U_{i}, U_j \subseteq V_{k}$ are disjoint. 
\end{definition}

\begin{remark}\label{anpus}
Let $M$ be a manifold of dimension. Then $\Nerve( \Disk{M})^{\otimes}$
is the $\infty$-operad associated to the ordinary colored operad $\calO$ whose objects
are elements of $\Disk{M}$, with morphisms given by 
$$ \Mul_{\calO}( \{ U_1, \ldots, U_n \}, V ) = \begin{cases} \ast & \text{ if } U_1 \cup \ldots \cup U_n \subseteq V \text{ and } U_i \cap U_j = \emptyset \text{ for } i \neq j \\
\emptyset & \text{ otherwise. } \end{cases}$$
In particular, $\Nerve( \Disk{M})^{\otimes}$ is an $\infty$-operad
(see Example \symmetricref{xe3}).
\end{remark}

\begin{remark}\label{lapus}
Let $\Disk{M}'$ denote the category whose objects are open embeddings
$\R^{k} \hookrightarrow M$, and whose morphisms are commutative diagrams
$$ \xymatrix{ \R^{k} \ar[rr]^{f} \ar[dr] & & \R^{k} \ar[dl] \\
& M & }$$
where $f$ is an open embedding. Then the forgetful functor
$(j: \R^{k} \hookrightarrow M) \mapsto j( \R^{k})$ determines an equivalence of categories from
$\Disk{M}'$ to $\Disk{M}$. If we regard $\Disk{M}$ as a colored operad via the construction
of Remark \ref{anpus}, then $\Disk{M}'$ inherits the structure of a colored operad, to which
we can associate an $\infty$-operad $\Nerve( \Disk{M}')^{\otimes}$ equipped with
an equivalence $\phi: \Nerve( \Disk{M}')^{\otimes} \rightarrow \Nerve( \Disk{M})^{\otimes}$.
The forgetful functor $(j: \R^{k} \hookrightarrow M) \mapsto \R^{k}$ determines a map
of colored operads from $\Disk{M}'$ to $\TopE{ \BTop(k) }$. Passing to nerves, we obtain a
map of $\infty$-operads $\Nerve( \Disk{M}')^{\otimes} \rightarrow \OpE{ \BTop(k) }$, which
naturally factors through the map $\OpE{M} \rightarrow \OpE{ \BTop(k) }$. Composing with a homotopy inverse to $\phi$, we get a map of $\infty$-operads $\Nerve( \Disk{M})^{\otimes} \rightarrow \OpE{M}$.

We can describe the situation roughly as follows: the objects of the $\infty$-operads
$\Nerve( \Disk{M})^{\otimes}$ and $\OpE{M}$ are the same: copies of Euclidean space
$\R^{k}$ equipped with an embedding in $M$. However, the morphisms are slightly different:
an $n$-ary operation in $\OpE{M}$ is a diagram of open embeddings
$$ \xymatrix{ \coprod_{1 \leq i \leq n} \R^{k} \ar[rr] \ar[dr] & & \R^{k} \ar[dl] \\
& M & }$$
which commutes up to (specified) isotopy, while an $n$-ary operation in
$\Nerve( \Disk{M})^{\otimes}$ is given by a diagram as above which commutes on the nose.
\end{remark}

The map of $\infty$-operads $\psi: \Nerve( \Disk{M})^{\otimes} \rightarrow \OpE{M}$ appearing in Remark \ref{lapus} is not an equivalence. For example, the underlying $\infty$-category of $\OpE{M}$ is the Kan complex
$B_M \simeq \Sing(M)$, while the underlying $\infty$-category of $\Nerve( \Disk{M})^{\otimes}$ is the nerve of the partially ordered set $\Disk{M}$, which is certainly not a Kan complex. However, this is essentially the only difference: the map $\psi$ exhibits $\OpE{M}$ as the $\infty$-operad obtained from
$\Nerve( \Disk{M})^{\otimes}$ by inverting each of the morphisms in $\Disk{M}$.
More precisely, we have the following result:

\begin{theorem}\label{sazz}
Let $M$ be a manifold and let $\calC^{\otimes}$ be an $\infty$-operad.
Composition with the map $$\Nerve( \Disk{M})^{\otimes} \rightarrow \OpE{M}$$ of Remark \ref{lapus}
induces a fully faithful embedding
$\theta: \Alg_{ \OpE{M} }(\calC) \rightarrow \Alg_{ \Nerve(\Disk{M})}(\calC)$. The essential image
of $\theta$ is the full subcategory of $\Alg_{ \Nerve(\Disk{M})}(\calC)$ spanned by the
locally constant $\Nerve(\Disk{M})^{\otimes}$-algebra objects of $\calC$ (see Definition \ref{plux}).  
\end{theorem}

Theorem \ref{sazz} is an immediate consequence of Proposition \ref{loam}, together with the following pair of lemmas:

\begin{lemma}\label{cate}
Let $M$ be a manifold of dimension $k$. Then the map 
$\Nerve( \Disk{M})^{\otimes} \rightarrow \OpE{M}$ induces a weak homotopy equivalence
$\psi: \Nerve(\Disk{M}) \rightarrow B_M$.
\end{lemma}

\begin{lemma}\label{postcate}
The map of $\infty$-operads $\Disk{M}^{\otimes} \rightarrow \OpE{M}$ is ornamental.
\end{lemma}

\begin{proof}[Proof of Lemma \ref{cate}]
The construction $U \mapsto B_U$ determines a functor $\chi$ from the category
$\Disk{M}$ to the category of simplicial sets. Let $X$ denote the relative nerve $\Nerve_{\chi}( \Disk{M} )$ (see \S \toposref{altstr}), so that we have a coCartesian fibration $\theta: X \rightarrow \Nerve( \Disk{M})$ whose over an object $U \in \Disk{M}$ is the Kan complex $B_U$. Remark \ref{talrod} implies
that the fibers of $\theta$ are contractible, so that $\theta$ is a trivial Kan fibration.
The projection map $\theta$ has a section $s$, which carries an object $U \in \Disk{M}$ to a chart
$\R^{k} \simeq U$ in $B_U$. The map $\psi$ is obtained by composing the section $s$
with the evident map $\psi': X \rightarrow B_{M}$. Consequently, it will suffice to show that
the map $\psi'$ is a weak homotopy equivalence. According to Proposition \toposref{charspacecolimit},
this is equivalent to the requirement that $B_{M}$ be a colimit of the diagram
$\{ U \mapsto B_U \}_{ U \in \Disk{M} }$ in the $\infty$-category of spaces $\SSet$. Using Remark \ref{talrod} again, we may reduce to showing that $\Sing M$ is a colimit of the diagram
$\{ U \mapsto \Sing U \}_{U \in \Disk{M}}$. In view of Theorem \ref{vankamp}, we need only show
that for every point $x \in M$, the partially ordered set $P: \{ U \in \Disk{M}: x \in U \}$ is weakly
contractible. In fact, $P^{op}$ is filtered: for every finite collection of open disks $U_i \subseteq M$
containing $x$, the intersection $\bigcap_{i} U_i$ is an open neighborhood of $x$ which contains
a smaller open neighborhood $V \simeq \R^{k}$ of $x$ (because $M$ is a topological manifold).
\end{proof}

\begin{proof}[Proof of Lemma \ref{postcate}]
In view of Remark \ref{kuz}, it is sufficient to show that the composite map
$$ \gamma: \Nerve( \Disk{M})^{\otimes} \rightarrow \OpE{M} \rightarrow \OpE{ \BTop(k) }$$
is ornamental. To this end, fix an object $U \in \Disk{M}$ and an
integer $m \geq 0$; wish to prove that the map
$$ \psi: \Nerve( \Disk{M})^{\otimes}_{/U} \times_{ \Nerve(\FinSeg)_{/ \seg{1} } } \{ \seg{m} \} \rightarrow \OpE{ \BTop(k) }_{/U} \times_{ \Nerve( \FinSeg)_{/ \seg{1}} } \{ \seg{m} \}$$
is a weak homotopy equivalence. We can identify the domain of $\psi$ with
the nerve $\Nerve(A)$, where $A \subseteq \Disk{M}^{m}$ denotes the partially
ordered set of sequences
$(V_1, \ldots, V_m) \in \Disk{M}^{m}$ such that $\bigcup V_{i} \subseteq U$ and
$V_{i} \cap V_{j} = \emptyset$ for $i \neq j$. 

It will now suffice to show that $\psi$ induces a homotopy equivalence after passing to the homotopy
fiber over some point of the (connected) Kan complex $\BTop(k)^{m}$. Unwinding the definitions, we must show that the canonical map
$$ \hocolim_{ (V_1, \ldots, V_m) \in A} \prod_{1 \leq i \leq m} \Sing \Emb( \R^{k}, V_i)
\rightarrow \Sing \Emb( \R^{k} \times \nostar{m}, U )$$
is a weak homotopy equivalence. Using Proposition \ref{splay}, we can reduce to showing
instead that the map
$$ \hocolim_{ (V_1, \ldots, V_m) \in A} \prod_{1 \leq i \leq m} \Germ(V_i)
\rightarrow \Germ( \nostar{m}, U)$$
is a homotopy equivalence. Both sides are acted on freely by the simplicial group
$\Germ_0( \R^{k})$. Consequently, it will suffice to show that we obtain a weak homotopy
equivalence of quotients
$$ \hocolim_{ ( V_1, \ldots, V_m) \in A} \prod_{1 \leq i \leq m} \Conf( \{i\}, V_i)
\rightarrow \Conf( \nostar{m}, U).$$
In view of Theorem \ref{vankamp}, it will suffice to show that for every injective map
$\phi: \nostar{m} \rightarrow U$, the partially ordered set $A_{\phi} = \{ (V_1, \ldots, V_m) \in A: \phi(i) \in V_i \}$ has weakly contractible nerve. This is clear, since $A_{\phi}^{op}$ is filtered
(because each point $\phi(i)$ has arbitrarily small neighborhoods homeomorphic to Euclidean space
$\R^{k}$).
\end{proof}

We can summarize Theorem \ref{sazz} informally as follows. To give an $\OpE{M}$-algebra object $A$ of
a symmetric monoidal $\infty$-category $\calC$, we need to specify the following data:
\begin{itemize}
\item[$(i)$] For every open disk $U \subseteq M$, an object $A(U) \in \calC$.
\item[$(ii)$] For every collection of disjoint open disks $V_1, \ldots, V_n$ contained in an open disk
$U \subseteq M$, a map $A(V_1) \otimes \ldots \otimes A(V_n) \rightarrow A(U)$, which is an equivalence when $n=1$.
\end{itemize}

In \S \ref{secran}, we will explain how to describe this data in another way: namely, as a cosheaf on the
Ran space of $M$ (see Definition \ref{coors}). However, in the setting of the Ran space, it is much more convenient to work with a {\em nonunital} version of the theory of $\OpE{M}$-algebras. Consequently, we will spend the remainder of this section explaining how to adapt the above ideas to the nonunital case.

\begin{definition}\label{psaab}
For every $k$-manifold $M$, we let $\OpE{M}_{\nunit}$ denote the $\infty$-operad
$\OpE{M} \times_{ \Nerve(\FinSeg)} \Nerve( \Surj)$. It follows from Remark \ref{spazzy}
and Proposition \ref{sorpi} that $\OpE{M}_{\nunit}$ is the assembly of the
$B_M$-family of $\infty$-operads $(B_M \times \Nerve(\FinSeg)) \times_{ B_M^{\amalg} } \OpE{M}_{\nunit}$,
each fiber of which is equivalent to the nonunital little cubes operad $\OpE{k}_{\nunit}$. 
\end{definition}

If $\calC^{\otimes}$ is a symmetric monoidal $\infty$-category, we let
$\Alg^{\nunit}_{\OpE{M}}(\calC)$ denote the $\infty$-category $\Alg_{ \OpE{M}_{\nunit}}(\calC)$
of nonunital $\OpE{M}$-algebra objects of $\calC$. 
Our next goal is to show that the results of \S \ref{pluy} can be generalized to the present setting: that is, for any symmetric monoidal $\infty$-category $\calC$, we can identify
$\Alg_{ \OpE{M}}(\calC)$ with a subcategory of $\Alg^{\nunit}_{ \OpE{M}}(\calC)$ (Proposition \ref{dont}). Our first step is to identify the relevant subcategory more precisely.

\begin{definition}
If $\calC^{\otimes}$ is a symmetric monoidal $\infty$-category and $M$ is a manifold of dimension $k > 0$, we will say that an $\OpE{M}_{\nunit}$-algebra object
$A \in \Alg_{ \OpE{M}_{\nunit}}( \calC)$ is {\it quasi-unital} if, for every point $U \in B_M$, the restriction of $A$ to the fiber $(\{ U \} \times \Nerve(\FinSeg)) \times_{ B_M^{\amalg} } \OpE{M}_{\nunit} \simeq \OpE{k}_{\nunit}$ determines a quasi-unital $\OpE{k}_{\nunit}$-algebra object of $\calC$, in the sense of Definition \ref{laster}.
Similarly, we will say that a map $f: A \rightarrow B$ of quasi-unital $\OpE{M}_{\nunit}$-algebra
objects of $\calC$ is {\it quasi-unital} if its restriction to each fiber $(\{ U \} \times \Nerve(\FinSeg)) \times_{ B_M^{\amalg} } \OpE{M}_{\nunit}$ determines a quasi-unital map of $\OpE{k}_{\nunit}$-algebras.
We let $\Alg_{ \OpE{M} }^{\qunit}(\calC)$ denote the subcategory fo $\Alg_{ \OpE{M}_{\nunit} }(\calC)$
spanned by the quasi-unital $\OpE{M}_{\nunit}$-algebra objects of $\calC$ and quasi-unital morphisms between them.
\end{definition}

\begin{remark}\label{elim}
Let $M$ be a manifold of dimension $k > 0$ and let $A$ be a $\OpE{M}_{\nunit}$-algebra object of a symmetric monoidal
$\infty$-category $\calC^{\otimes}$. Fix a point $U \in B_M$, corresponding to an open embedding
$\psi: \R^{k} \hookrightarrow M$. We will say that a map $u: {\bf 1} \rightarrow A(U)$ in
$\calC$ is a {\it quasi-unit} for $A$ if, for every pair of objects $V,W \in B_M$ and every morphism
$\phi: U \oplus V \rightarrow W$, the composite map
$$ A(V) \simeq {\bf 1} \otimes A(V) \stackrel{u}{\rightarrow} A(U) \otimes A(V) \rightarrow A(W)$$
is homotopic to the map induced by the composition $U \rightarrow U \oplus V \stackrel{\phi}{\rightarrow} W$ in $\OpE{M}_{\nunit}$. Note that it suffices to check this condition in the special case
where $V = W = U$ and, if $k > 1$, where $\phi$ is a single map (arbitrarily chosen). 
Unwinding the definition, we see that $A$ is quasi-unital if and only if
there exists a quasi-unit $u: {\bf 1} \rightarrow A(U)$ for each $U \in B_M$. Similarly,
a map $A \rightarrow B$ between quasi-unital $\OpE{M}_{\nunit}$-algebra objects is
quasi-unital if, for every quasi-unit $u: {\bf 1} \rightarrow A(U)$, the composite map
${\bf 1} \stackrel{u}{\rightarrow} A(U) \rightarrow B(U)$ is a quasi-unit for $B$.
Moreover, if $M$ is connected, then it suffices to check these conditions for a single
$U \in B_M$.
\end{remark}

\begin{proposition}\label{dont}
Let $M$ be a manifold of dimension $k > 0$ and let $\calC^{\otimes}$ be a symmetric monoidal $\infty$-category. Then the restriction functor $\Alg_{ \OpE{M}}( \calC) \rightarrow \Alg_{ \OpE{M}}^{\qunit}(\calC)$
is an equivalence of $\infty$-categories.
\end{proposition}

\begin{proof}
For every map of simplicial sets $K \rightarrow B_M$, let $\calO^{\otimes}_{K}$ denote
the $K$-family of $\infty$-operads $(K \times \Nerve(\FinSeg)) \times_{ B_M^{\amalg}} \OpE{M}$,
let ${\calO'}^{\otimes}_{K} = (K \times \Nerve(\FinSeg)) \times_{ B_M^{\amalg}} \OpE{M}_{\nunit}$.
Note that the projection map $q: \calO^{\otimes}_{K} \rightarrow K$ is a coCartesian fibration.
Let $\Alg'_{\calO_K}( \calC)$ denote the full subcategory of $\Alg_{\calO_K}(\calC)$ spanned
by those $\infty$-operad maps which carry $q$-coCartesian morphisms to equivalences in $\calC$, let
$\Alg'_{\calO'_{K}}(\calC)$ be defined similarly, and let
$\Alg_{ \calO_{K}}^{\qunit}(\calC)$ denote the subcategory of $\Alg'_{ \calO'_{K}}(\calC)$
spanned by those objects which restrict to quasi-unital ${\calO'}^{\otimes}_{ \{v\} } \simeq \OpE{k}_{\nunit}$-algebra objects of $\calC$ and those morphisms which restrict to quasi-unital ${\calO'}^{\otimes}_{ \{v\} } \simeq \OpE{k}_{\nunit}$-algebra maps for every vertex $v \in K$. There is an evident
restriction map $\theta_{K}: \Alg'_{\calO_K}(\calC) \rightarrow \Alg^{\qunit}_{\calO'_{K}}(\calC)$
fitting into a commutative diagram
$$ \xymatrix{ \Alg_{ \OpE{M}}(\calC) \ar[r] \ar[d] & \Alg_{ \OpE{M}}^{\qunit}(\calC) \ar[d] \\
\Alg'_{ \calO_{ K}} \ar[r]^{\theta_K} & \Alg^{\qunit}_{\calO'_{K}}( \calC). }$$
If $K = B_M$, then the vertical maps are categorical equivalences. Consequently, it will suffice to prove
that $\theta_{K}$ is an equivalence for every map of simplicial sets $K \rightarrow B_M$.
The collection of simplicial sets $K$ which satisfy this condition is clearly
stable under homotopy colimits; we can therefore reduce to the case where $K$
is a simplex, in which case the desired result follows from Theorem \ref{quas}.
\end{proof}

It follows from Lemma \ref{postcate} and Remark \ref{spazzy} that for every manifold $M$, the map
$\Disk{M}^{\otimes}_{\nunit} \rightarrow \OpE{M}_{\nunit}$ is ornamental. Combining this with Lemma \ref{cate} and Proposition \ref{loam}, we deduce the following nonunital variant of Theorem \ref{sazz}:

\begin{proposition}\label{saz}
Let $M$ be a manifold and let $\calC^{\otimes}$ be an $\infty$-operad.
Then composition with map $\Nerve( \Disk{M})^{\otimes} \rightarrow \OpE{M}$ of Remark \ref{lapus}
induces a fully faithful embedding
$\theta: \Alg^{\nunit}_{ \OpE{M} }(\calC) \rightarrow \Alg^{\nunit}_{ \Nerve(\Disk{M})}(\calC)$. The essential image
of $\theta$ is the full subcategory of $\Alg^{\nunit}_{ \Nerve(\Disk{M})}(\calC)$ spanned by the
locally constant objects.
\end{proposition}

\begin{definition}
Let $M$ be a manifold of dimension $k > 0$ and let $\calC^{\otimes}$ be a symmetric monoidal $\infty$-category. We will say that a locally constant $\Disk{M}^{\otimes}_{\nunit}$-algebra object of $\calC$
is {\it quasi-unital} if it corresponds to a quasi-unital $\OpE{M}_{\nunit}$-algebra object of $\calC$
under the equivalence of Proposition \ref{saz}. Similarly, we will say that a map
$f: A \rightarrow B$ between locally constant quasi-unital $\Disk{M}^{\otimes}_{\nunit}$-algebra objects of $\calC$ is {\it quasi-unital} if it corresponds to a quasi-unital morphism in 
$\Alg_{ \OpE{M}}^{\nunit}(\calC)$ under the equivalence of Proposition \ref{saz}.
We let $\Alg^{\qunit, \loc}_{ \Disk{M}}( \calC)$ denote the subcategory of
$\Alg_{ \Disk{M}}(\calC)$ spanned by the quasi-unital, locally constant $\Disk{M}^{\otimes}$-algebra objects of $\calC$ and quasi-unital morphisms between them.
\end{definition}

\begin{remark}
Let $A \in \Alg_{ \Disk{M}}^{\nunit}(\calC)$, let $W \in \Disk{M}$ be an open disk in $M$, and let
$U \subseteq W$ be an open disk with compact closure in $W$. We say that a map
${\bf 1} \rightarrow A(U)$ in $\calC$ is a {\it quasi-unit} for $A$ if, for every disk
$V \in \Disk{M}$ such that $V \subseteq W$ and $V \cap U = \emptyset$, the diagram
$$ \xymatrix{ {\bf 1} \otimes A(V) \ar[d] \ar[r]^-{u \otimes \id} & A(U) \otimes A(V) \ar[d] \\
A(V) \ar[r] & A(W) }$$
commutes up to homotopy. Note that if $M$ has dimension at least $2$, it suffices to check
this condition for a single open disk $V$. Unwinding the definition, we see that 
$A$ is quasi-unital if and only if there exists a quasi-unit $u: {\bf 1} \rightarrow A(U)$ for
every pair $U \subseteq W$ as above, and a map $f: A \rightarrow B$ in
$\Alg_{ \Disk{M}}^{\nunit}(\calC)$ is quasi-unital if and only if composition with $f$ carries
every quasi-unit ${\bf 1} \rightarrow A(U)$ to a quasi-unit ${\bf 1} \rightarrow B(U)$ (see
Remark \ref{elim}). In fact, it suffices to check these conditions for a single pair $U \subseteq W$ in each connected component of $M$. 
\end{remark}

Combining Proposition \ref{dont}, Theorem \ref{sazz}, and Proposition \ref{saz}, we arrive
at the following:

\begin{proposition}
Let $M$ be a manifold of dimension $k > 0$ and $\calC^{\otimes}$ a symmetric monoidal
$\infty$-category. Then the restriction functor
$\Alg_{ \Disk{M}}(\calC) \rightarrow \Alg_{ \Disk{M}}^{\qunit}(\calC)$
induces an equivalence between the full subcategories spanned by the locally constant algebras.
\end{proposition}

In other words, there is no essential loss of information in passing from unital
$\Disk{M}^{\otimes}$-algebras to nonunital $\Disk{M}^{\otimes}$-algebras, at least in the locally constant case. For this reason, we will confine our attention to nonunital algebras in \S \ref{secran}.

\subsection{The Ran Space}\label{secran}

\begin{definition}\label{coors}
Let $M$ be a manifold. We let $\Ran(M)$ denote the collection of nonempty finite subsets $S \subseteq M$
which have nonempty intersection with each connected component of $M$.
We will refer to $\Ran(M)$ as the {\it Ran space} of $M$.
\end{definition}

The Ran space $\Ran(M)$ admits a natural topology, which we will define in a moment. 
Our goal in this section is to study the basic properties of $\Ran(M)$ as a topological space.
Our principal results are Theorem \ref{adley2}, which asserts that $\Ran(M)$ is weakly contractible
(provided that $M$ is connected), and Proposition \ref{scamm}, which
characterizes sheaves on $\Ran(M)$ which are constructible with respect to the natural filtration
of $\Ran(M)$ by cardinality of finite sets.

Our first step is to define the topology on $\Ran(M)$. First, we need to introduce a bit of notation.
Suppose that $\{ U_i \}_{1 \leq i \leq n}$ is a nonempty collection of pairwise disjoint subsets of $M$.
We let $\Ran( \{ U_i \}) \subseteq \Ran(M)$ denote the collection of finite sets $S \subseteq M$ such that $S \subseteq \bigcup U_i$ and $S \cap U_i$ is nonempty for $1 \leq i \leq n$.

\begin{definition}
Let $M$ be a manifold. We will regard the Ran space $\Ran(M)$ as equipped with the coarsest topology for which the subsets $\Ran( \{ U_i \}) \subseteq \Ran(M)$ are open, for every nonempty finite collection of pairwise disjoint open sets $\{ U_i \}$ of $M$.
\end{definition}

\begin{remark}
If $\{ U_i \}$ is a nonempty finite collection of pairwise disjoint open subsets of a manifold $M$, then the open subset $\Ran( \{ U_i \}) \subseteq \Ran(M)$ is homeomorphic to a product
$\prod_{i} \Ran( U_i )$, via the map $\{ S_i \subseteq U_i \}) \mapsto ( \bigcup_i S_i \subseteq M)$. 
\end{remark}

\begin{remark}\label{spow}
Let $M$ be a manifold, and let $S = \{ x_1, \ldots, x_n \}$ be a point of
$\Ran(M)$. Then $S$ has a basis of open neighborhoods in $\Ran(M)$ of the form
$\Ran( \{ U_i \})$, where the $U_i$ range over all collections of disjoint open neighborhoods of the points
$x_i$ in $M$.
Since $M$ is a manifold, we may further assume that that each $U_i$ is homeomorphic to Euclidean space. 
\end{remark}

\begin{remark}\label{slashe}
If we choose a metric $d$ on on the manifold $M$, then the topology on $\Ran(M)$ is described
by a metric $D$, where
$$D(S,T) = \sup_{s \in S} \inf_{t \in T} d(s,t) + \sup_{t \in T} \inf_{s \in S} d(s,t).$$
It follows that $\Ran(M)$ is paracompact.
\end{remark}

Our first main result in this section is the following observation of Beilinson and Drinfeld:

\begin{theorem}[Beilinson-Drinfeld]\label{adley2}
Let $M$ be a connected manifold. Then $\Ran(M)$ is weakly contractible.
\end{theorem}

We first formulate a relative version of Theorem \ref{adley} which is slightly easier to prove.

\begin{notation}
Let $M$ be a manifold and $S$ a finite subset of $M$. We let 
$\Ran(M)_S$ denote the closed subset of $\Ran(M)$ consisting of those nonempty
finite subsets $T \subseteq \Ran(M)$ such that $S \subseteq T$.
\end{notation}

\begin{lemma}[Beilinson-Drinfeld]\label{adley}
Let $M$ be a connected manifold and let $S$ be a nonempty finite subset of $M$. Then $\Ran_S(M)$ is weakly contractible.
\end{lemma}

\begin{proof}
We first prove that $\Ran(M)_{S}$ is path connected. Let $T$ be a subset of $M$ containing
$S$. For each $t \in T$, choose a path $p_t: [0,1] \rightarrow M$ such that $p_t(0) = t$ and
$p_t(1) \in S$ (this is possible since $M$ is connected and $S$ is nonempty). Then
the map $r \mapsto S \cup \{ p_t(r) \}_{t \in T}$ determines a continuous path in
$\Ran(M)_S$ joining $T$ with $S$.
We will complete the proof by showing that for each $n > 0$, every element
$\eta \in \pi_n \Ran(M)_S$ is trivial; here we compute the homotopy group
$\pi_n$ with respect to the base point given by $S \in \Ran(M)_S$.

The topological space $\Ran(M)_{S}$ admits a continuous product
$U: \Ran(M)_S \times \Ran(M)_S \rightarrow \Ran(M)_{S}$, given by the formula
$U(T,T') = T \cup T'$. This product induces a map of homotopy groups
$$\phi: \pi_n \Ran(M)_S \times \pi_n \Ran(M)_S \rightarrow \pi_n \Ran(M)_{S}$$
(here the homotopy groups are taken with respect to the base point $S \in \Ran(M)_S$).
Since $S$ is a unit with respect to the multiplication on $\Ran(M)_S$, we conclude that
$\phi(\eta, 1) = \eta = \phi(1, \eta)$ (where we let $1$ denote the unit element of the homotopy
group $\pi_n \Ran(M)_S$). Because the composition of the diagonal embedding
$\Ran(M)_S \rightarrow \Ran(M)_S \times \Ran(M)_S$ with $U$ is the identity from $\Ran(M)_S$ to itself, we have also $\phi( \eta, \eta) = \eta$. It follows that
$$ \eta = \phi( \eta,\eta) = \phi( \eta,1) \phi( 1, \eta) = \eta^2$$
so that $\eta = 1$ as desired. 
\end{proof}

\begin{proof}[Proof of Theorem \ref{adley2}]
For every point $x \in M$, choose an open embedding $j_x: \R^{k} \hookrightarrow M$ such
that $j_x(0) = x$. Let $U_x = j_x( B(1) )$ be the image under $j_x$ of the unit ball in
$\R^{k}$, and let $V_x$ be the open subset of $\Ran(M)$ consisting of those nonempty
finite subsets $S \subseteq M$ such that $S \cap U_x \neq \emptyset$.
Let $\calJ$ be the partially ordered set of all nonempty finite subsets of $M$
(that is, $\calJ$ is the Ran space $\Ran(M)$, but viewed as a partially ordered set).
We define a functor from $\calJ^{op}$ to the category of open subsets of $\Ran(M)$ by
the formula
$$ T \mapsto V_T = \bigcap_{x \in T} V_x. $$
For each $S \in \Ran(M)$, the partially ordered set $\{ T \in \calJ: S \in V_{T} \}$
is nonempty and stable under finite unions, and therefore has weakly contractible nerve. 
It follows that $\Sing \Ran(M)$ is equivalent to the homotopy colimit of the diagram
$\{ \Sing V_{T} \}_{T \in \calJ^{op} }$. We will prove that each of the spaces $V_{T}$ is
weakly contractible, so that this homotopy colimit is weakly homotopy equivalent to
$\Nerve( \calJ^{op} )$ and is therefore weakly contractible.

Fix $T \in \calJ$, and choose a continuous family of maps $\{ h_r: \R^{k} \rightarrow \R^{k} \}_{0 \leq r \leq 1}$ with
the following properties:
\begin{itemize}
\item[$(i)$] For $0 \leq r \leq 1$, the map $h_r$ is the identity outside of a ball
$B(2) \subseteq \R^{k}$ of radius $2$.
\item[$(ii)$] The map $h_0$ is the identity.
\item[$(iii)$] The map $h_1$ carries $B(1) \subseteq \R^{k}$ to the origin.
\end{itemize}
We now define a homotopy $\phi_{T}: \Ran(M) \times [0,1] \rightarrow \Ran(M)$ by the formula
$$ \phi_{T}( S, r) = S \cup \bigcup_{ x \in T} j_x h_{t} j_x^{-1}(S).$$
The homotopy $\phi_{T}$ leaves $V_{T}$ and $\Ran(M)_{T}$ setwise fixed, and
carries $V_{T} \times \{1\}$ into $\Ran(M)_{T}$. It follows that the inclusion
$\Ran(M)_{T} \subseteq V_{T}$ is a homotopy equivalence, so that
$V_{T}$ is weakly contractible by Lemma \ref{adley}.
\end{proof}

We now discuss a natural stratification of the Ran space.

\begin{definition}
Let $M$ be a manifold. We let $\Ran^{\leq n}(M)$ denote the subspace of
$\Ran(M)$ consisting of those subsets $S \subseteq M$ having cardinality $\leq n$,
and $\Ran^{n}(M)$ the subspace of $\Ran^{\leq n}(M)$ consisting of those subsets
$S \subseteq M$ having cardinality exactly $n$.
\end{definition}

\begin{remark}
The set $\Ran^{\leq n}(M)$ is closed in $\Ran(M)$, and $\Ran^{n}(M)$ is open in 
$\Ran^{\leq n}(M)$.
\end{remark}

\begin{definition}\label{plumb}
Let $M$ be a manifold and let $\calF \in \Shv( \Ran(M) )$ be a sheaf on $\Ran(M)$.
For each $n \geq 0$, let $i(n): \Ran^{\leq n}(M) \rightarrow \Ran(M)$ denote the inclusion map.
We will say that $\calF$ is {\it constructible} if the following conditions are satisfied:
\begin{itemize}
\item[$(1)$] The canonical map $\calF \rightarrow \varprojlim_{n} i(n)_{\ast} i(n)^{\ast} \calF$
is an equivalence.
\item[$(2)$] For each $n$, the restriction of $i(n)^{\ast} \calF$ to the open subset
$\Ran^{n}(M) \subseteq \Ran^{\leq n}(M)$ is locally constant.
\end{itemize}
\end{definition}

\begin{remark}
Condition $(2)$ of Definition \ref{plumb} is equivalent to the requirement that
$\calF$ be $\Z_{\geq 0}$-constructible, where we regard $\Ran(M)$ as
$\Z_{\geq 0}$-stratified via the map $\Ran(M) \rightarrow \Z_{\geq 0}$ given by
$S \mapsto |S|$. We refer the reader to \S \ref{coofer} for a general review of the theory of constructible sheaves. Here we are required to impose condition $(1)$ because the partially ordered set
$\Z_{\geq 0}$ does not satisfy the ascending chain condition.
\end{remark}

\begin{remark}
We can endow the topological space $\Ran(M)$ with another topology, where a set
$U \subseteq \Ran(M)$ is open if and only if its intersection with each $\Ran^{\leq n}(M)$ is open
(with respect to the topology of Definition \ref{coors}). If $\calF$ is a sheaf on $\Ran(M)$ with respect
to this second topology, then condition $(1)$ of Definition \ref{plumb} is automatic: this follows
from Proposition \toposref{siegland}.
\end{remark}

The following result gives a convenient characterization of constructible sheaves on the Ran space:

\begin{proposition}\label{scamm}
Let $M$ be a manifold and $\calF \in \Shv( \Ran(M) )$. Then $\calF$ is constructible
if and only if it is hypercomplete and satisfies the following additional condition:
\begin{itemize}
\item[$(\ast)$] For every nonempty finite collection of disjoint disks $U_1, \ldots, U_n \subseteq M$
containing open subdisks $V_1 \subseteq U_1, \ldots, V_n \subseteq U_n$, the restriction map
$\calF( \Ran( \{ U_i \} ) )
\rightarrow \calF( \Ran( \{ V_i \}) )$ is a homotopy equivalence.
\end{itemize}
\end{proposition}

\begin{proof}
We first prove the ``only if'' direction. Suppose that $\calF$ is constructible. To show that
$\calF$ is hypercomplete, we write $\calF$ as a limit $\varprojlim i(n)_{\ast} i(n)^{\ast}
\calF$ as in Definition \ref{plumb}. It therefore suffices to show that each
$i(n)^{\ast} \calF$ is hypercomplete. This follows from the observation that
$\Ran^{\leq n}(M)$ is a paracompact topological space of finite covering dimension
(Corollary \toposref{fdfd}). 

We now prove every constructible sheaf $\calF \in \Shv( \Ran(M) )$ satisfies $(\ast)$. For $1 \leq i \leq n$, we invoke Theorem \ref{scen} to choose an isotopy $\{ h^{t}_i: V_i  \rightarrow U_i \}_{t \in \R}$
such that $h^0_{i}$ is the inclusion of $V_i$ into $U_i$ and $h^1_{i}$ is a homeomorphism.
These isotopies determine an open embedding
$$H: \Ran( \{ V_i \}) \times \R \rightarrow
\Ran( \{ U_i \}) \times \R.$$
Let $\calF' \in \Shv( \Ran( \{  U_i \} )\times \R )$
be the pullback of $\calF$, so that $\calF'$ is hypercomplete (see Lemma \ref{sceen} and
Example \ref{tsooch}). It follows that $H^{\ast} \calF'$ is hypercomplete. Since $\calF$
is constructible, we deduce that $\calF'$ is foliated. For $t \in \R$, let
$\calF'_{t}$ denote the restriction of $\calF$ to $\Ran( \{ V_i \}) \times \{t\}$.
We have a commutative diagram of spaces
$$ \xymatrix{ \calF( \Ran( \{ U_i \}) ) \ar[rr]^{\theta} \ar[dr]^{\theta'} & &  \calF'( 
( \Ran( \{ V_i \} ) \times \R) \ar[dl]^{\theta''}  \\
& \calF'_{1}( \Ran( \{ V_i \}) ) & }$$
Since each $h_i^{1}$ is a homeomorphism, we deduce that $\theta'$ is a homotopy equivalence.
Proposition \ref{postcanca} guarantees that $\theta''$ is a homotopy equivalence, so that $\theta$
is a homotopy equivalence by the two-out-of-three property. 
Applying Proposition \ref{postcanca} again, we deduce that the composite map
$\calF( \Ran(\{ U_i \} ) ) \rightarrow \calF'_0( \Ran( \{ V_i \}) ) \simeq \calF( \Ran( \{ V_i \} ))$ is a homotopy equivalence as desired.

We now prove the ``if'' direction of the proposition. Assume that $\calF$ is hypercomplete
and that $\calF$ satisfies $(\ast)$; we wish to prove that $\calF$ is constructible. 
We first show that the restriction of $\calF$ to each $\Ran^{n}(M)$ is locally constant.
Choose a point $S \in \Ran^{n}(M)$; we will show that $\calF | \Ran^{n}(M)$ is constant in a neighborhood of $S$. Let $S = \{ x_1, \ldots, x_n \}$, and choose disjoint open disks
$U_1, \ldots, U_n \subseteq X$ such that $x_{i} \in U_i$. Let $W \subseteq \Ran^{n}(M)$
denote the collection of all subsets $S \subseteq M$ which contain exactly one point from
each $U_i$. We will prove that $\calF | \Ran^{n}(M)$ is constant on $W$. Let
$X = \calF( \Ran( \{ U_i \} ) )$. Since
$W \subseteq \Ran( \{ U_i \})$, there is a canonical map
from the constant sheaf on $W$ taking the value $X$ to $\calF | W$; we will show that this map
is an equivalence. Since $W \simeq U_1 \times \ldots \times U_n$ is a manifold,
it has finite covering dimension so that $\Shv(W)$ is hypercomplete. Consequently, to show that a morphism in $\Shv(W)$ is an equivalence, it suffices to check after passing to the stalk at each point
$\{ y_1, \ldots, y_n \} \in W$. This stalk is given by
$\varinjlim_{V} \calF(V)$, where the colimit is taken over all open subsets $V \subseteq \Ran(M)$ containing $\{ y_1, \ldots, y_n \}$. It follows from Remark \ref{spow} that it suffices to take the
colimit over those open sets $V$ of the form $\Ran( \{ V_i \} )$, where
each $V_i \subseteq U_i$ is an open neighborhood of $y_i$. Condition $(\ast)$ guarantees
that each of the maps $X \rightarrow \calF(V)$ is a homotopy equivalence, so after passing
to the filtered colimit we obtain a homotopy equivalence $X \rightarrow \varinjlim_{V} \calF(V)$ as desired.

Let $\calG = \varprojlim_{n} i(n)_{\ast} i(n)^{\ast} \calF$ (using the notation of Definition \ref{plumb}).
To complete the proof, it will suffice to show that the canonical map $\alpha: \calF \rightarrow \calG$
is an equivalence. Since each $i(n)^{\ast} \calF$ is automatically hypercomplete (because
$\Ran^{\leq n}(M)$ is a paracompact space of finite covering dimension), we see that
$\calG$ is hypercomplete. Using the results of \S \toposref{hcovh}, we deduce that the collection
of those open sets $U \subseteq \Ran(M)$ such that $\alpha$ induces a homotopy equivalence
$\alpha_U: \calF(U) \rightarrow \calG(U)$ is stable under the formation of unions of hypercoverings.
It therefore suffices to show that $\alpha_{U}$ is an homotopy equivalence for some collection of open
sets $U$ which forms a basis for the topology of $\Ran(M)$. By virtue of Remark \ref{spow}, we may
assume that $U = \Ran( \{ U_i \} )$ for some collection of disjoint open
disks $U_1, \ldots, U_n$ meeting every connected component of $M$.

For each integer $m$, let $\calF^{\leq m} = i(m)^{\ast} \calF$.
We wish to prove that the map $\calF(U) \rightarrow \varprojlim_{m} \calF^{\leq m}( U \cap \Ran^{\leq m} (M))$ is a homotopy equivalence. In fact, we will prove that the individual maps
$\calF(U) \rightarrow \calF^{\leq m}( U \cap \Ran^{\leq m}(M))$ are homotopy equivalences
for $m \geq n$. Choose a point $x_i$ in each disk $U_i$, and let $S = \{ x_1, \ldots, x_n \}$.
Let $\calF_{S}$ denote the stalk of $\calF$ at the point $S$. We have a commutative diagram of restriction maps
$$ \xymatrix{ \calF(U) \ar[dr]^{\phi} \ar[rr] & & \calF^{\leq m}( U \cap \Ran^{\leq m}(M) ) \ar[dl]^{\phi'} \\
& \calF_{S} & }$$ 
where $\phi$ is a homotopy equivalence by the argument given above. By the two-out-of-three property, we are reduced to proving that $\phi'$ is a homotopy equivalence.

The set $U \cap \Ran^{\leq m}(M)$ admits a stratification by the linearly ordered set $[m]$, which
carries a point $T \in \Ran(M)$ to the cardinality of $T$. Let $\calC = \Sing^{[m]} (U \cap \Ran^{\leq m}(M) )$. Since $\calF$ is constructible, the sheaf $\calF | ( U \cap \Ran^{\leq m}(M))$ corresponds to some left fibration $q:\widetilde{\calC} \rightarrow \calC$ under the equivalence of $\infty$-categories
provided by Theorem \ref{mainstrat}. Under this equivalence, we can identify
$\calF^{\leq m}( U \cap \Ran^{\leq m}(M))$ with the $\infty$-category
$\Fun_{\calC}( \calC, \widetilde{\calC})$ of sections of $q$, while 
$\calF_{S}$ corresponds to the fiber of $\widetilde{\calC}_{S}$ of $q$ over the point
$S \in \calC$. To prove that $\theta'$ is an equivalence, it suffices to show that
$S$ is an initial object of $\calC$. To this end, choose homeomorphisms
$\psi_{i}: \R^{k} \rightarrow U_i$ for $1 \leq i \leq n$ such that $\psi_i( 0 ) = x_i$.
We then have a map
$$c: [0,1] \times (U \cap \Ran^{\leq m}(M)) \rightarrow (U \cap \Ran^{\leq m}(M) )$$
given by the formula $c(t,T) = \{ \psi_i( t v): \psi_i(v) \in T \}$. The continuous map
$c$ induces a natural transformation from the inclusion
$\{ S \} \hookrightarrow \calC$ to the identity functor from $\calC$ to itself, thereby
proving that $S \in \calC$ is initial as desired.
\end{proof}

To apply Proposition \ref{scamm}, it is convenient to have the following characterization of
hypercompleteness:

\begin{proposition}\label{coolman}
Let $X$ be a topological space, $\calU(X)$ the collection of open subsets of $X$, and
$\calF: \Nerve( \calU(X)^{op}) \rightarrow \SSet$ a presheaf on $X$. The following conditions
are equivalent:
\begin{itemize}
\item[$(1)$] The presheaf $\calF$ is a hypercomplete sheaf on $X$.
\item[$(2)$] Let $U$ be an open subset of $X$, $\calC$ be a category, and $f: \calC \rightarrow \calU(U)$ a functor. Suppose that, for every point $x \in U$, the full subcategory $\calC_{x} = \{ C \in \calC: x \in f(C) \} \subseteq \calC$ has weakly contractible nerve. Then $\calF$ exhibits $\calF(U)$ as a limit of the diagram
$\Nerve( \calC)^{op} \rightarrow \Nerve( \calU(X)^{op} ) \stackrel{\calF}{\rightarrow} \SSet$.
\end{itemize}
\end{proposition}

\begin{lemma}\label{coolrod}
Let $X$ be a topological space, and let $\calF \in \Shv(X)$ be an $\infty$-connective sheaf
satisfying the following condition:
\begin{itemize}
\item[$(\ast)$] Let $A$ be a partially ordered set and $f: A \rightarrow \calU(X)^{op}$ an order-preserving map such that, for every point $x \in X$, the full subcategory $A_{x} = \{ a \in A: x \in f(a) \} \subseteq A$
is filtered. Then $\calF$ exhibits $\calF(X)$ as a limit of the diagram
$\Nerve( \calA) \rightarrow \Nerve( \calU(X)^{op} ) \stackrel{\calF}{\rightarrow} \SSet$.
\end{itemize}
Then the space $\calF(X)$ is nonempty.
\end{lemma}

\begin{proof}
The functor $\calF: \Nerve(\calU(X)^{op}) \rightarrow \SSet$ classifies a left fibration
$q: \calE \rightarrow \Nerve( \calU(X)^{op} )$. We will construct a partially ordered set
$A$ and a map $\psi: \Nerve(A) \rightarrow \calE$ such that the composite map
$\Nerve(A) \rightarrow \Nerve( \calU(X)^{op})$ and each subset $A_x$ is filterd.
According to Corollary \toposref{charspacelimit}, we can identify the limit $\varprojlim_{a \in A} \calF( f(a) )$ with the Kan complex $\Fun_{ \Nerve( \calU(X)^{op})}( \Nerve(A), \calE)$, which is nonempty by construction.

We will construct a sequence of partially ordered sets
$$ \emptyset = A(0) \subseteq A(1) \subseteq \ldots $$
and compatible maps $\psi(n): \Nerve( A(n) ) \rightarrow \calE$ with the following properties:
\begin{itemize}
\item[$(i)$] For every element $a \in A(n)$, the set $\{ b \in A(n): b < a \}$ is a finite subset of $A(n-1)$.
\item[$(ii)$]
For every point $x \in X$ and every finite subset $S \subseteq A(n-1)_{x}$, there exists
an upper bound for $S$ in $A(n)_{x}$.
\end{itemize}
Assuming that this can be done, we can complete the proof by taking $A = \bigcup_{n} A(n)$ and $\psi$ be the amalgamation of the maps $\psi(n)$.

The construction now proceeds by induction on $n$. Assume that $n > 0$ and that the map
$\psi(n-1): \Nerve( A(n-1) ) \rightarrow \calE$ has already been constructed. Let
$K$ be the set of pairs $(x, S)$, where $x \in X$ and $S$ is a finite subset of
$A(n-1)_{x}$ which is closed-downwards (that is, $a \leq a'$ and $a' \in S$ implies $a \in S$). 
We define $A(n)$ to be the disjoint union $A(n-1) \coprod K$.
We regard $A(n)$ as a partially ordered set, where
$a < b$ in $A(n)$ if and only $a,b \in A(n-1)$ and $a < b$ in $A(n-1)$, or
$a \in A(n-1)$, $b = (x,S) \in K$, and $a \in S$. It is clear that $A(n)$ satisfies condition $(i)$.
It remains only to construct a map $\psi(n): \Nerve( A(n) ) \rightarrow \calE$ which extends
$\psi(n-1)$ and satisfies $(ii)$. Unwinding the definitions, we must show that for every pair
$(x,S) \in K$, the extension problem
$$ \xymatrix{ \Nerve(S) \ar[r]^{\psi'} \ar[d] & \calE \\
\Nerve(S)^{\triangleright} \ar@{-->}[ur]^{ \phi} & }$$
admits a solution, where $\psi'$ denotes the restriction
$\psi(n-1) | \Nerve(S)$ and $\phi$ carries the cone point of $\Nerve(S)^{\triangleright}$ to
an object $E \in \calE$ such that $x \in q(E) \in \calU(X)$. 

Since $S$ is finite, the subset $U = \bigcap_{ s \in S} q( \psi'(s) )$ is an open subset of
$X$ containing the point $x$. The map $\psi'$ determines a diagram
$\alpha: \Nerve(S) \rightarrow \calE \times_{ \Nerve( \calU(X)^{op} ) } \{ U \} \simeq \calF(U)$.
To prove the existence of $\phi$, it suffices to show that there exists a smaller open
subset $V \subseteq U$ containing $x$ such that the composite map
$\Nerve(S) \stackrel{\alpha}{\rightarrow} \calF(U) \rightarrow \calF(V)$ is nullhomotopic.
Since $\Nerve(S)$ is finite, it suffices to show $\alpha$ induces a nullhomotopic map
from $\Nerve(S)$ into the stalk $\calF_{x} = \varinjlim_{x \in V} \calF(V)$. We conclude by
observing that $\calF_{x}$ is contractible (since $\calF$ is assumed to be $\infty$-connective).
\end{proof}

\begin{proof}[Proof of Proposition \ref{coolman}]
Suppose first that $(1)$ is satisfied; we will verify $(2)$. Let $\chi: \calU(X) \rightarrow \Shv(X)$
be the functor which carries an open set $U$ to the sheaf $\chi_{U}$ given by the formula
$$\chi_{U}(V) = \begin{cases} \Delta^0 & \text{ if } V \subseteq U \\
\emptyset & \text{ otherwise. }\end{cases}$$
Let $\calG = \varinjlim_{ C \in \calC} \chi_{f(C)}$. For every point $x \in U$, the stalk
$\calG_{x}$ is weakly homotopy equivalent to the nerve of the category $\calC_{x}$,
and for $x \notin U$ the stalk $\calG_{x}$ is empty. If each $\calC_{x}$ has weakly contractible nerve, then we conclude that the canonical map $\calG \rightarrow \chi_{U}$ is $\infty$-connective, so that
$$\calF(U) \simeq \bHom_{ \Shv(X)}( \chi_{U}, \calF) \simeq
\bHom_{ \Shv(X)}(\calG, \calF) \simeq \varprojlim_{C \in \calC} \bHom_{ \Shv(X)}(\chi_{f(C)}, \calF)
= \varprojlim_{C \in \calC} \calF( f(C) ).$$

Now suppose that $(2)$ is satisfied. Let $S \subseteq \calU(X)$ be a covering sieve on
an open set $U \subseteq X$. Then for each $x \in U$, the partially ordered set $S_x = \{ V \in S: x \in V \}$ is nonempty and stable under finite intersections, so that $\Nerve( S_x)^{op}$ is filtered and therefore weakly contractible. It follows from $(2)$ that the map $\calF(U) \rightarrow \varprojlim_{V \in S} \calF(V)$
is a homotopy equivalence, so that $\calF$ is a sheaf. It remains to show that $\calF$ is hypercomplete.
Choose an $\infty$-connective morphism $\alpha: \calF \rightarrow \calF'$, where $\calF'$ is hypercomplete; we wish to show that $\alpha$ is an equivalence. The first part of the proof shows
that $\calF'$ also satisfies the condition stated in $(2)$. Consequently, it will suffice to prove the following:
\begin{itemize}
\item[$(\ast)$] Let $\alpha: \calF \rightarrow \calG$ be an $\infty$-connective morphism in
$\Shv(X)$, where $\calF$ and $\calG$ both satisfy $(2)$. Then $\alpha$ is an equivalence.
\end{itemize}
To prove $(\ast)$, it suffices to show
that for each open set $U \subseteq X$, $\alpha$ induces a homotopy equivalence $\alpha_U: \calF(U) \rightarrow \calG(U)$. We will show that $\alpha_{U}$ is $n$-connective for each $n \geq 0$, using
induction on $n$. If $n > 0$, then we can conclude by applying the inductive hypothesis to the diagonal
map $\beta: \calF \rightarrow \calF \times_{ \calG} \calF$. It remains to consider the case $n=0$: that is, to
show that the map $\alpha_{U}$ is surjective on connected components. In other words, we must show that every map $\chi_{U} \rightarrow \calG$ factors through $\alpha$. This follows by applying
Lemma \ref{coolrod} to the fiber product $\chi_U \times_{ \calG} \calF$ (and restricting to the open set $U$).
\end{proof}

\subsection{Topological Chiral Homology}\label{tchtch}

Let $M$ be a $k$-manifold and $\calC^{\otimes}$ a symmetric monoidal $\infty$-category.
We can think of an $\OpE{M}$-algebra $A \in \Alg_{ \OpE{M}}(\calC)$ as a family of
$\OpE{k}$-algebras $A_x \in \Alg_{ \OpE{k}}(\calC)$, parametrized by the points $x \in M$.
(More precisely, we should think of this family as {\em twisted}: each $A_x$ should really be viewed
as an algebra over the operad of little $k$-cubes in the tangent space $T_{M,x}$ to $M$ at the point $x$.) In this section, we will explain how to extract form $A$ a global invariant
$\int_{M} A$, which we call the {\it topological chiral homology} of $M$ (with coefficients in $A$).
Our construction is a homotopy-theoretic analogue of the Beilinson-Drinfeld theory of chiral homology
described in \cite{beilinson}. It should be closely related to the theory of {\it blob homology} studied by Morrison and Walker.

The basic idea of the construction is simple. According to Theorem \ref{sazz}, we can think of an
$\OpE{M}$-algebra object $A$ of a symmetric monoidal $\infty$-category $\calC$ as a functor
which assigns to every disjoint union of open disks $U \subseteq M$ an object $A(U) \in \calC$, which carries disjoint unions to tensor products. Our goal is to formally extend the definition of $A$ to {\em all} open subsets of $M$. Before we can give the definition, we need to establish a bit of terminology.

\begin{definition}\label{saefer}
Let $M$ be a manifold and $\calU(M)$ the partially ordered set of
all open subsets of $M$. We can identify objects of the $\infty$-category
$\Nerve( \calU(M))^{\amalg}$ with finite sequences $(U_1, \ldots, U_n)$
of open subsets of $M$. We let $\Nerve( \calU(M) )^{\otimes}$ denote the subcategory
of $\Nerve( \calU(M) )^{\amalg}$ spanned by those morphisms
$(U_1, \ldots, U_n) \rightarrow (V_1, \ldots, V_m)$ which cover a map
$\alpha: \seg{n} \rightarrow \seg{m}$ in $\FinSeg$ and possess the following property:
for $1 \leq j \leq m$, the sets $\{ U_i \}_{\alpha(i) = j }$ are disjoint open subsets of $V_j$.
\end{definition}

For every manifold $M$, the nerve $\Nerve( \calU(M)^{\otimes})$ is an $\infty$-operad which
contains $\Nerve( \Disk{M})^{\otimes}$ as a full subcategory. 

\begin{definition}
We will say that a symmetric monoidal $\infty$-category $\calC^{\otimes}$ is
{\it sifted-complete} if the underlying $\infty$-category $\calC$ admits small sifted colimits and the
tensor product functor $\calC \times \calC \rightarrow \calC$ preserves small sifted colimits.
\end{definition}

\begin{remark}
If a simplicial set $K$ is sifted, then the requirement that the tensor product $\calC \times \calC \rightarrow \calC$ preserve sifted colimits is equivalent to the requirement that it preserve sifted colimits separately in each variable.
\end{remark}

\begin{example}
Let $\calC^{\otimes}$ be a symmetric monoidal $\infty$-category. Assume that the underlying
$\infty$-category $\calC$ admits small colimits, and that the tensor product on $\calC$ preserves small colimits separately in each variable. Let $\calO^{\otimes}$ be an arbitrary small $\infty$-operad,
so that $\Alg_{\calO}(\calC)$ inherits a symmetric monoidal structure
(given by pointwise tensor product). The $\infty$-category $\Alg_{\calO}(\calC)$ itself admits small
colimits (Corollary \symmetricref{smurf}), but the tensor product on $\Alg_{\calO}(\calC)$ generally does not preserve colimits in each variable. However, it does preserve {\em sifted} colimits separately in each variable: this follows from Proposition \symmetricref{fillfemme}. Consequently, $\Alg_{\calO}(\calC)$ is a sifted-complete symmetric monoidal $\infty$-category.
\end{example}

The main existence result we will need is the following:

\begin{theorem}\label{cookydef}
Let $M$ be a manifold and let $q: \calC^{\otimes} \rightarrow \Nerve(\FinSeg)$ be a sifted-complete symmetric monoidal $\infty$-category. For every algebra object $A \in \Alg_{\OpE{M}}(\calC)$, the restriction
$A | \Nerve( \Disk{M})^{\otimes}$ admits an operadic left Kan extension to
$\Nerve( \calU(M)^{\otimes})$.
\end{theorem}

Assuming Theorem \ref{cookydef} for the moment, we can give the definition of topological chircal homology.

\begin{definition}\label{tch}
Let $M$ be a manifold and let $\calC^{\otimes}$ be a sifted-complete symmetric monoidal $\infty$-category.
We let $\int: \Alg_{ \OpE{M}}(\calC) \rightarrow \Alg_{ \Nerve( \calU(M))}(\calC)$ be the functor
given by restriction to $\Nerve( \Disk{M})^{\otimes}$ followed by operadic left Kan extension along
the inclusion $\Nerve(\Disk{M})^{\otimes} \rightarrow \Nerve( \calU(M)^{\otimes})$. 
If $A \in \Alg_{\OpE{M}}(\calC)$ and
$U$ is an open subset of $M$, we will denote the value of $\int(A)$ on the open set $U \subseteq M$
by $\int_{U} A \in \calC$. We will refer to $\int_{U} A$ as the {\it topological chiral homology of $U$
with coefficients in $A$}.
\end{definition}

\begin{remark}\label{tww}
To describe the content of Definition \ref{tch} more concretely, it is useful to introduce a bit of notation.
If $M$ is a manifold, we let $\Disj{M}$ denote the partially ordered subset of
$\calU(M)$ spanned by those open subsets $U \subseteq M$ which are homeomorphic to
$S \times \R^{k}$ for some finite set $S$. 
In the situation of Definition \ref{tch}, the algebra object $A$ determines a functor
$\theta: \Nerve( \Disj{M} ) \rightarrow \calC$, given informally by the formula
$$ V_1 \cup \cdots \cup V_n \mapsto A(V_1) \otimes \cdots \otimes A(V_n)$$
(here the $V_i$ denote pairwise disjoint open disks in $M$). The topological chiral homology
$\int_{M} A \in \calC$ is then given by the colimit of the diagram $\theta$.
\end{remark}

\begin{example}
Let $U \subseteq M$ be an open subset homeomorphic to Euclidean space. Then
there is a canonical equivalence $A(U) \simeq \int_{U} A$. 
\end{example}

\begin{remark}\label{spacer}
Suppose that we have a map of $\infty$-operads $\psi: \OpE{M} \rightarrow \calO^{\otimes}$, where $\calO^{\otimes}$ is some other $\infty$-operad. Let
$A \in \Alg_{\calO}(\calC)$. Then we will abuse notation by denoting the topological chiral homology $\int_{M} (\psi \circ A)$ simply by $\int_{M} A$. This abuse is consistent with the notation of
Definition \ref{kost} in the following sense: if $A \in \Alg_{ \OpE{M}}( \calC)$, then the topological
chiral homology $\int_{U} A$ of $U$ with coefficients in $A$ is equivalent to the topological
chiral homology $\int_{U} (A| \OpE{U})$ of $U$ with coefficients in
the induced $\OpE{U}$-algebra. 
\end{remark}

\begin{example}
Let $A \in \Alg_{ \OpE{ \BTop(k) }}( \calC)$. Then Remark \ref{spacer} allows us to define
the topological chiral homology $\int_{M} A$ of {\em any} $k$-manifold with coefficients in $A$.
Similarly, if $A \in \Alg_{\OpE{\Smooth}}(\calC)$ (see Example \ref{kample2}), then 
$\int_{M} A$ is defined for any smooth manifold $M$. Many other variations on this theme are possible:
roughly speaking, if $A$ is an $\OpE{k}$-algebra object of $\calC$ equipped with a compatible
action of some group $G$ mapping to $\Top(k)$, then $\int_{M} A$ is well-defined if we are provided with a reduction of the structure group of $M$ to $G$.
\end{example}

In order to prove Theorem \ref{cookydef} (and to establish the basic formal properties of topological chiral homology), we need to have good control over colimits indexed by partially ordered sets of the form
$\Disj{M}$, where $M$ is a manifold (see Remark \ref{tww}). We will obtain this control by introducing a less rigid version of the $\infty$-category $\Nerve( \Disj{M} )$, where we allow open disks in $M$ to ``move''.

\begin{definition}
Fix an integer $k \geq 0$. We let $\EMB(k)$ denote the topological category whose objects
are $k$-manifolds, with morphism spaces given by $\bHom_{ \EMB(k)}(N,M) = \Emb(N,M)$.
If $M$ is a $k$-manifold, we let $\DISJ{M}$ denote the full subcategory of
the $\infty$-category $\Nerve( \EMB(k) )_{/M}$ spanned by those objects of the form
$j: N \rightarrow M$, where $N$ is homeomorphic to $S \times \R^{k}$ for some
finite set $S$.
\end{definition}

\begin{remark}\label{saba}
An object of the $\infty$-category $\DISJ{M}$ can be identified with a finite collection
of open embeddings $\{ \psi_{i}: \R^{k} \hookrightarrow M \}_{1 \leq i \leq n}$ having disjoint images.
Up to equivalence, this object depends only on the sequence of images $( \psi_1( \R^{k}), \ldots, \psi_n( \R^{k}) )$, which we can identify with an object of the category $\Disj{M}$. However,
the morphisms in these two categories are somewhat different: a morphism in $\DISJ{M}$ is given by a diagram
$$ \xymatrix{ \coprod_{1 \leq i \leq m} \R^{k} \ar[dr]^{ \{ \phi_i \} } \ar[rr] & & \coprod_{1 \leq j \leq n} \R^{k} \ar[dl]^{ \{ \psi_j \} } \\ 
& M & }$$
which commutes up to (specified) isotopy, which does not guarantee an inclusion of images
$\bigcup \phi_i( \R^{k} ) \subseteq \bigcup \psi_{j}( \R^{k} )$.

Nevertheless, there is an evident functor $\gamma: \Nerve( \Disj{M} ) \rightarrow \DISJ{M}$, defined by choosing a parametrization of each open disk in $M$ (up to equivalence, the functor $\gamma$ is independent of these choices).
\end{remark}

The fundamental result we will need is the following:

\begin{proposition}\label{koio}
Let $M$ be a $k$-manifold. Then:
\begin{itemize}
\item[$(1)$] The functor $\gamma: \Nerve( \Disj{M}) \rightarrow \DISJ{M}$, described in Remark \ref{saba}, is cofinal.
\item[$(2)$] Let $\Disj{M}_{\nunit}$ denote the subcategory of
$\Disj{M}$ whose objects are nonempty open sets $U \in \Disj{M}$ and whose morphisms are inclusions 
$U \hookrightarrow V$ such that the induced map $\pi_0 U \rightarrow \pi_0 V$
is surjective. If $M$ is connected, then the induced functor $\Nerve( \Disj{M}^{\nunit}) \rightarrow \DISJ{M}$ is cofinal.
\end{itemize}
\end{proposition}

The second assertion of Proposition \ref{koio} will require the following technical result, which
employs the notation introduced in \S \ref{secran}. 

\begin{lemma}\label{badus}
Let $M$ be a connected manifold, let $S$ be a finite subset of $M$, and let
$\Disj{M}^{\nunit}_{S}$ denote the full subcategory of $\Disj{M}^{\nunit}$ spanned by those objects
$V \in \Disj{M}^{\nunit}$ such that $S \subseteq V$. Then the simplicial set
$\Nerve( \Disj{M}^{\nunit} )$ is weakly contractible.
\end{lemma}

\begin{proof}
For every object $V \in \Disj{M}^{\nunit}_{S}$, let $\psi(V)$ denote the subset of
$\Ran(M)$ consisting of those subsets $T$ with the following properties:
\begin{itemize}
\item[$(i)$] We have inclusions $S \subseteq T \subseteq V$.
\item[$(ii)$] The map $T \rightarrow \pi_0 V$ is surjective.
\end{itemize}

For every point $T \in \Ran(M)_{S}$, let $\calC_{T}$ denote the full subcategory of
$\Disj{M}^{\nunit}_{S}$ spanned by those objects $V$ such that $T \in \psi(V)$. 
Each of the category $\calC_{T}^{op}$ is filtered (for every finite collection 
$V_1, \ldots, V_n \in \calC_{T}$, we can choose $V \in \calC_{T}$ such that
$V \subseteq \bigcap V_{i}$ and each of the maps $\pi_0 V \rightarrow \pi_0 V_i$
is surjective: namely, take $V$ to be a union of sufficiently small open disks containing the points of $T$). It follows from Theorem \ref{vankamp} that the Kan complex
$\Sing \Ran(M)_S$ is equivalent to the homotopy colimit of the diagram
$\{ \psi(V) \}_{V \in \Disk{M}^{\nunit}_{S} }$. For each
$V \in \Disj{M}^{\nunit}_{S}$, write $V$ as a disjoint union of 
open disks $U_{1} \cup \ldots \cup U_{m}$.
Then $\psi(V)$ is homeomorphic to a product $\prod_{1 \leq i \leq m} \Ran(U_m)_{ S \cap U_m}$,
and is therefore weakly contractible by Lemmas \ref{adley} and \ref{adley2}.
It follows that the Kan complex $\Sing(\Ran(M)_S)$ is weakly homotopy equivalent to the nerve of the category
$\Disj{M}^{\nunit}$. The desired result now follows from the weak contractiblity of
$\Sing \Ran(M)_S$ (Lemmas \ref{adley} and \ref{adley2}).
\end{proof}

\begin{proof}[Proof of Proposition \ref{koio}]
We first give the proof of $(1)$. Let $S = \{ 1, \ldots, n \}$, let $U = S \times \R^{k}$, and let $\psi: U \rightarrow M$ be an open embedding
corresponding to an object of $\DISJ{M}$. According to Theorem \toposref{hollowtt}, it will suffice
to show that the $\infty$-category $\calC = \Nerve( \Disj{M} ) \times_{ \DISJ{M} } \DISJ{M}_{\psi/}$
is weakly contractible. We observe that the projection map $\calC \rightarrow \Nerve( \Disj{M} )$ is
a left fibration, associated to a functor $\chi: \Nerve( \Disj{M} ) \rightarrow \SSet$ which carries
each object $V \in \Disj{M}$ to the homotopy fiber of the map of Kan complexes $\Sing \Emb(U, V) \rightarrow \Sing \Emb(U, M)$. According to Proposition \toposref{charspacecolimit}, it will suffice
to show that the colimit $\varinjlim(\chi)$ is contractible. Since colimits in $\SSet$ are universal,
it will suffice to show that $\Sing \Emb(U,M)$ is a colimit of the diagram
$\{ \Sing \Emb(U,V) \}_{V \in \Disj{M} }$. Using Theorem \toposref{charleschar} and
Remark \ref{gec}, we are reduced to showing that $\Sing \Conf(S,M)$ is a colimit of the
diagram $\{ \Sing \Conf(S,V) \}_{V \in \Disj{M}}$. According to Theorem \ref{vankamp}, it will suffice
to show that for every injective map $j: S \hookrightarrow M$, the partially ordered set
$\Disj{M}_S = \{ V \in \Disj{M}: j(S) \subseteq V \}$ has weakly contractible nerve. This is clear, since
$\Disj{M}^{op}_S$ is filtered: every open neighborhood of $j(S)$ contains a union of sufficiently
small open disks around the points $\{ j(s) \}_{s \in S}$. 

The proof of $(2)$ is identical except for the last step: we must instead show that that for every
injective map $j: S \hookrightarrow M$, the category $\Disj{M}_{S}^{\nunit} = \{ V \in \Disj{M}^{\nunit}: j(S) \subseteq V \}$ has weakly contractible nerve, which follows from Lemma \ref{badus}.
\end{proof}

The advantage of the $\infty$-category $\DISJ{M}$ over the more rigid $\infty$-category
$\Nerve( \Disj{M} )$ is summarized in the following result:

\begin{proposition}\label{scun}
For every manifold $M$, the $\infty$-category $\DISJ{M}$ is sifted.
\end{proposition}

\begin{proof}
We wish to prove that the diagonal map $\delta: \DISJ{M} \rightarrow \DISJ{M} \times \DISJ{M}$ is cofinal. 
We have a commutative diagram
$$ \xymatrix{ & \Nerve(\Disj{M}) \ar[dl]^{\gamma} \ar[dr]^{\theta} & \\
\DISJ{M} \ar[rr]^{\delta} & & \DISJ{M} \times \DISJ{M}, }$$
where $\gamma$ is cofinal by virtue of Proposition \ref{koio}. It will therefore suffice to show that
$\theta$ is cofinal (Proposition \toposref{cofbasic}). Fix a pair of objects $\phi: U \hookrightarrow M$,
$\psi: V \hookrightarrow M$ in $\DISJ{M}$. According to Theorem \toposref{hollowtt}, it will suffice to show that the $\infty$-category $\calC = \DISJ{M}_{\phi/} \times_{ \DISJ{M} } \Nerve( \Disj{M}) 
\times_{ \DISJ{M} } \DISJ{M}_{\psi/}$ is weakly contractible. There is an evident left fibration
$\calC \rightarrow \Nerve( \Disj{M} )$, classified by a functor $\chi: \Nerve(\Disj{M}) \rightarrow \SSet$
which carries an object $W \in \Disj{M}$ to the homotopy fiber of the map
$$ \Sing(\Emb(U,W) \times \Emb(V,W)) \rightarrow \Sing( \Emb(U,M) \times \Emb(V,M) )$$
over the vertex given by $(\phi,\psi)$. Using Proposition \toposref{charspacecolimit}, we can identify the weak homotopy type of $\calC$ with the colimit $\varinjlim(\chi) \in \SSet$. Consequently, it will
suffice to show that $\varinjlim(\chi)$ is contractible. Since colimits in $\SSet$ are universal, it will suffice to show that $\Sing( \Emb(U,M) \times \Emb(V,M) )$ is a colimit of the diagram $\chi': \Nerve( \Disj{M}) \rightarrow \SSet$ given by the formula $\chi'(W) = \Sing( \Emb(U,W) \times \Emb(V,W) )$. 
Let $S \subseteq U$, $T \subseteq V$ be subsets containing one point from each connected component of $U$ and $V$, respectively, and let $\chi'': \Nerve( \Disj{M}) \rightarrow \SSet$ be the functor given by the formula $W \mapsto \Sing( \Conf(S, W) \times \Conf(T,W) )$. There is an evident
restriction functor of diagrams $\chi' \rightarrow \chi''$. Using Remark \ref{gec} and Theorem \toposref{charleschar}, we are reduced to proving that the canonical map
$\varinjlim( \chi'') \rightarrow \Sing( \Conf(S,M) \times \Conf(T,M) )$ is a homotopy equivalence. In view of
Theorem \ref{vankamp}, it will suffice to show that for every point $(j,j') \in \Conf(S,M) \times \Conf(T,M)$,
the full subcategory $\Disj{M}_{(j,j')}$ of $\Disj{M}$ spanned by those objects $W \in \Disj{M}$ such that
$j(S), j'(T) \subseteq W$ is weakly contractible. This is clear, since $\Disj{M}_{(j,j')}^{op}$ is filtered.
\end{proof}

Armed with Proposition \ref{scun}, we are ready to prove that topological chiral homology is well-defined.

\begin{proof}[Proof of Theorem\ref{cookydef}]
According to Theorem \symmetricref{oplk}, it will suffice to show that for each open set $U \subseteq M$ the induced diagram
$$\Nerve( \Disj{U} ) \stackrel{\theta}{\rightarrow} \DISJ{U} \stackrel{\beta}{\rightarrow} \OpE{M} \stackrel{A}{\rightarrow} \calC^{\otimes}$$
can be extended to an operadic colimit diagram in $\calC^{\otimes}$. Since $\theta$ is cofinal
(Proposition \ref{koio}), it suffices to show that $A \circ \beta$ can be extended to an operadic
colimit diagram in $\calC^{\otimes}$. Choose a $q$-coCartesian natural transformation from
$A \circ \beta$ to a functor $\chi: \DISJ(U) \rightarrow \calC$, given informally by the formula
$\chi( \{ \psi_{i}: V_i \hookrightarrow U \}_{1 \leq i \leq n}) = A( \psi_1) \otimes \cdots \otimes A( \psi_n)$.
In view of Proposition \symmetricref{chocolateoperad}, it will suffice to show that
$\chi$ can be extended to an operadic colimit diagram in $\calC$. Since $\DISJ{U}$ is sifted
(Proposition \ref{scun}) and the tensor product on $\calC$ preserves sifted colimits separately in each variable, it suffices to show that $\chi$ can be extended to a colimit diagram in $\calC$ (Proposition
\symmetricref{optest}). This colimit exists because $\calC$ admits sifted colimits and $\DISJ{U}$ is sifted.
\end{proof}

We close this section with the following result concerning the functorial behavior of topological chiral homology:

\begin{proposition}\label{skoo}
Let $M$ be a manifold, and let $F: \calC^{\otimes} \rightarrow \calD^{\otimes}$ be a symmetric
monoidal functor. Assume $\calC^{\otimes}$ and $\calD^{\otimes}$ are sifted-complete and that the underlying functor $F: \calC \rightarrow \calD$ preserves sifted colimits.
Then:
\begin{itemize}
\item[$(1)$] If $A \in \Alg_{ \Nerve(\calU(M))}( \calC)$ has the property that $A_0 = A| \Nerve( \Disk{M})^{\otimes}$
is locally constant and $A$ is an operadic left Kan extension of $A_0$, then
$FA$ is an operadic left Kan extension of $FA_0$.
\item[$(2)$] For any locally constant algebra $A \in \Alg_{ \OpE{M} }(\calC)$,
the canonical map
$\int_{ M} FA \rightarrow F(\int_{M} A)$ is an equivalence in $\calC$.
\end{itemize}
\end{proposition}

\begin{proof}
We first prove $(1)$. Since $A_0$ is locally constant, we can assume that $A_0$ factors as a composition
$\Nerve( \Disk{M})^{\otimes} \rightarrow \OpE{M} \stackrel{A'_0}{\rightarrow} \calC^{\otimes}$
(Theorem \ref{sazz}). We wish to prove that for every object $U \in \calU(M)$, the diagram
$FA$ exhibits $FA(U) \in \calD$ as an operadic colimit of the composite diagram
$$\Disj{M} \stackrel{\alpha}{\rightarrow} \DISJ{M} \stackrel{\beta}{\rightarrow} \OpE{M} \stackrel{A'_0}{\rightarrow} \calC^{\otimes} \stackrel{F}{\rightarrow} \calD^{\otimes}.$$
Since $\alpha$ is cofinal (Proposition \ref{koio}), it will suffice to show that $FA$ exhibits
$FA(U)$ as an operadic colimit of $F \circ A'_0 \circ \beta$. 

Let $p: \calC^{\otimes} \rightarrow \Nerve(\FinSeg)$ exhibit $\calC^{\otimes}$ as a symmetric
monoidal $\infty$-category, and let $q: \calD^{\otimes} \rightarrow \Nerve(\FinSeg)$ exhibit
$\calD^{\otimes}$ as a symmetric monoidal $\infty$-category. 
Choose a $p$-coCartesian natural transformation $\alpha$
from $A'_0 \circ \beta$ to a map $\phi: \DISJ{M} \rightarrow \calC$. Since $F$ is a symmetric monoidal
functor, $F(\alpha)$ is a $q$-coCartesian natural transformation from $F \circ A'_0 \circ \beta$ to
$F \circ \phi$. It will therefore suffice to show that $FA$ exhibits $F A(U)$ as a colimit of
the diagram $F \circ \phi$ in the $\infty$-category $\calD$ (Propositions \symmetricref{chocolateoperad} and \symmetricref{optest}). Since $F | \calC$ preserves sifted colimits and
the $\infty$-category $\DISJ{M}$ is sifted (Proposition \ref{scun}), it suffices to show that
$A(U)$ is a colimit of the diagram $\phi$. Using Propositions \symmetricref{chocolateoperad} and
\symmetricref{optest} again, we are reduced to proving that
$A(U)$ is an operadic colimit of the diagram $A'_0 \circ \beta$, which (since $\alpha$ is cofinal)
follows from our assumption that $A$ is an operadic left Kan extension of $A_0$. This completes the proof of $(1)$. Assertion $(2)$ is an immediate consequence.
\end{proof}

\subsection{Properties of Topological Chiral Homology}\label{klimb}

Our goal in this section is to establish four basic facts about the theory of topological chiral homology. In what follows, we will assume that $\calC^{\otimes}$ is a sifted-complete symmetric monoidal $\infty$-category and $M$ a topological manifold of dimension $k$.

\begin{itemize}
\item[$(1)$] For a fixed algebra $A \in \Alg_{ \OpE{M}}(\calC)$, the construction $U \mapsto \int_{U} A$ carries disjoint unions of open subsets of $M$
to tensor products in the $\infty$-category $\calC$ (Theorem \ref{clapse}).
\item[$(2)$] For a fixed open set $U \subseteq M$, the construction
$A \mapsto \int_{U} A$ carries tensor products of $\OpE{M}$-algebra objects of $\calC$ to tensor products in
$\calC$ (Theorem \ref{coopsa}).
\item[$(3)$] If $A \in \Alg_{ \OpE{M}}( \calC)$ arises from a family $\{ A_x \}_{x \in M}$ of
{\it commutative} algebra objects of $\calC$, then $\int_{U} A$ can be identified with image
in $\calC$ of the colimit $\varinjlim_{x \in U}(A_x) \in \CAlg(\calC)$ (Theorem \ref{hunefer}). 
\item[$(4)$] If $k=1$ and $M$ is the circle $S^1$, then we can view an algebra object $A \in \Alg_{ \OpE{M}}(\calC)$ as an associative algebra object
of $\calC$ (equipped with an automorphism $\theta$ given by monodromy around the circle). In this case,
the topological chiral homology $\int_{M} A$ can be identified with the ($\theta$-twisted) Hochschild homology of $A$, which is computed by an analogue of the usual cyclic bar complex (Theorem \ref{junl}).
\end{itemize}

We begin with assertion $(1)$. The functor $\int$ of Definition \ref{tch} carries $\Alg_{ \OpE{M} }(\calC)$ into
$\Alg_{ \Nerve( \calU(M) )}( \calC)$. Consequently, whenever $U_1, \ldots, U_m$ are disjoint
open subsets of $U \subseteq M$, we have a multiplication map
$$ \int_{U_1} A \otimes \cdots \otimes \int_{U_m} A \rightarrow \int_{U} A.$$

\begin{theorem}\label{clapse}
Let $M$ be a manifold and $\calC^{\otimes}$ a sifted-complete symmetric monoidal $\infty$-category.
Then for every object
$A \in \Alg_{ \OpE{M} }(\calC)$ and every collection of pairwise disjoint open subsets
$U_1, \ldots, U_m \subseteq M$, the map
$$ \int_{U_1} A \otimes \cdots \otimes \int_{U_m} A \rightarrow \int_{ \bigcup U_i} A$$
is an equivalence in $\calC$.
\end{theorem}

\begin{proof}
It follows from Proposition \ref{koio} that for each open set $U \subseteq M$, the topological chiral homology $\int_{U} A$ is the colimit of a diagram $\psi_U: \DISJ(U) \rightarrow \calC$ given
informally by the formula $\psi_U( V_1 \cup \ldots \cup V_n) = A(V_1) \otimes \cdots \otimes A(V_n)$.
Since each $\DISJ{U_i}$ is sifted (Proposition \ref{scun}) and the tensor product on $\calC$ preserves sifted colimits separately in each variable, we can identify the tensor product $\int_{U_1} A \otimes \cdots \otimes \int_{U_m} A$ with the colimit $\varinjlim_{ \DISJ(U_1) \times \ldots \times \DISJ(U_n) } ( \psi_{U_1} \otimes \cdots \otimes \psi_{U_m})$. Let $W = \bigcup U_i$. The tensor product functor $\psi_{U_1} \otimes \cdots \otimes \psi_{U_m}$ can be identified with the pullback of $\psi_W$ along the evident map
$$\alpha: \DISJ{U_1} \times \cdots \times \DISJ{U_m} \rightarrow \DISJ{W}$$
$$ (V_1 \subseteq U_1, \ldots, V_m \subseteq U_m ) \mapsto V_1 \cup \ldots \cup V_m.$$
Consequently, we are reduced to proving that the $\alpha$ induces an equivalence
$$ \varinjlim (\alpha \circ \psi_{W})
\rightarrow \varinjlim \psi_{W}.$$
It will suffice to show that $\alpha$ is cofinal. 
This follows by applying Proposition \toposref{cofbasic} to the commutative diagram
$$ \xymatrix{ \Nerve( \Disj{U_1} \times \ldots \times \Disj{U_n} ) \ar[r] \ar[d] & \DISJ{U_1} \times \ldots \times \DISJ{U_m} \ar[d]^{\alpha} \\
\Nerve( \Disj{W} ) \ar[r] & \DISJ{W}; }$$
note that the horizontal maps are cofinal by Proposition \ref{koio}, and the map $\beta$ is an isomorphism of simplicial sets.
\end{proof}

To formulate assertion $(2)$ more precisely, suppose
we are given a pair of algebras
$A,B \in \Alg_{ \OpE{M}}(\calC)$. Let $\int(A), \int(B) \in \Alg_{ \Nerve( \calU(M) )}(\calC)$ be given by operadic left Kan extension. Then $(\int(A) \otimes \int(B)) | \Nerve( \Disk{M})^{\otimes}$ is an
extension of $(A \otimes B)| \Nerve( \Disk{M})^{\otimes}$, so we have a canonical map
$\int( A \otimes B) \rightarrow \int(A) \otimes \int(B)$. 
We then have the following:

\begin{theorem}\label{coopsa}
Let $M$ be a manifold and $\calC^{\otimes}$ a sifted-complete symmetric monoidal $\infty$-category.
Then for every pair of locally constant algebras $A,B \in \Alg_{ \OpE{M}}(\calC)$, the canonical map $\theta: \int_{M} (A \otimes B) \rightarrow \int_{M} A \otimes \int_{M} B$
is an equivalence in $\calC$.
\end{theorem}

\begin{proof}
Proposition \ref{koio} allows us to identify $\int_{M} A$ with the colimit of a diagram
$\phi: \DISJ{M} \rightarrow \OpE{M} \stackrel{A}{\rightarrow} \calC$ and
$\int_{M} B$ with the colimit of a diagram
$\psi: \DISJ{M} \rightarrow \OpE{M} \stackrel{B}{\rightarrow} \calC$.
Since the tensor product on $\calC$ preserves sifted colimits, we deduce that
$\int_{M} A \otimes \int_{M} B$ is given by the colimit of the functor
$(\phi \otimes \psi): \DISJ{M} \times \DISJ{M} \rightarrow \calC$. On the other hand,
the topological chiral homology is given by the colimit of the diagram
$\delta \circ (\phi \otimes \psi): \DISJ{M} \rightarrow \calC$, where $\delta:
\DISJ{M} \rightarrow \DISJ{M} \times \DISJ{M}$ is the diagonal map. 
The map $\theta$ is induced by the $\delta$, and is an equivalence since
$\delta$ is cofinal (Proposition \ref{scun}).
\end{proof}

The proof of assertion $(3)$ is based on the following simple observation:

\begin{lemma}\label{huuk}
Let $M$ be a manifold and $\calC$ an $\infty$-category which admits small colimits.
Regard $\calC$ as endowed with the coCartesian symmetric monoidal structure (see \S \symmetricref{jilun}). Then, for every object $A \in \Alg_{\OpE{M}}(\calC)$, 
the functor $\int A$ exhibits the topological chiral homology $\int_{M} A$ as
the colimit of the diagram $A | \Nerve( \Disk{M}): \Nerve( \Disk{M}) \rightarrow \calC$.
\end{lemma}

\begin{proof}
Let $\chi: \Nerve(\Disj{M}) \rightarrow \calC$ be the functor given informally by
the formula $\chi( U_1 \cup \ldots \cup U_n) = A(U_1) \coprod \cdots \coprod A(U_n)$,
where the $U_i$ are disjoint open disks in $M$. We observe that $\chi$ is a left
Kan extension of $\chi | \Nerve( \Disk{M} )$, so that $\int_{M} A \simeq
\colim \chi \simeq \colim (\chi | \Nerve( \Disk{M} ))$ (see Lemma \toposref{kan0}). 
\end{proof}

\begin{theorem}\label{hunefer}
Let $M$ be a manifold and $\calC^{\otimes}$ a sifted-complete symmetric monoidal $\infty$-category.
Regard the Kan complex $B_M$ as the underlying
$\infty$-category of the $\infty$-operad $B_M^{\amalg}$, and let $A \in \Alg_{ B_M}( \calC)$
so that $\int_{M} A$ is well-defined (see Remark \ref{spacer}). Composing $A$ with
the diagonal map $B_M \times \Nerve(\FinSeg) \rightarrow B_M^{\amalg}$, we obtain a
functor $\psi: B_M \rightarrow \CAlg(\calC)$. Let $A' = \colim(\psi) \in \CAlg(\calC)$. Then
there is a canonical equivalence $\int_{M} A \simeq A'( \seg{1})$ in the $\infty$-category $\calC$. 
\end{theorem}

\begin{remark}
Let $A$ be as in the statement of Theorem \ref{hunefer}. It follows from Theorem \symmetricref{kinema} that $A$ is determined by the functor $\psi$, up to canonical equivalence. In other words, we may identify $A \in \Alg_{B_M}(\calC)$ with a family
of commutative algebra objects of $\calC$ parametrized by the Kan complex $B_M$
(which is homotopy equivalent to $\Sing(M)$, by virtue of Remark \ref{talrod}). Theorem
\ref{hunefer} asserts that in this case, the colimit of this family of commutative algebras is computed by the formalism of topological chiral homology.
\end{remark}

\begin{proof}[Proof of Theorem \ref{hunefer}]
Let $\phi: \Disk{M}^{\otimes} \times \FinSeg \rightarrow \Disk{M}^{\otimes}$ be the functor given by the construction
$$( (U_1, \ldots, U_m), \seg{n}) \mapsto (U'_1, \ldots, U'_{mn}),$$ where $U'_{mi+j} = U_{j}$. 
Composing $\phi$ with the map $\Nerve( \Disk{M})^{\otimes} \rightarrow B_{M}^{\amalg}
\stackrel{A}{\rightarrow} \calC^{\otimes}$, we obtain a locally constant algebra object 
$\overline{A} \in \Alg_{ \Nerve( \Disk{M})}( \CAlg(\calC) )$, where $\CAlg(\calC)$ is endowed
with the symmetric monoidal structure given by pointwise tensor product (see Example \symmetricref{slabber}).
Since the symmetric monoidal structure on $\CAlg(\calC)$ is coCartesian (Proposition \symmetricref{cocarten}), the colimit $\varinjlim( \psi)$ can be identified with the
topological chiral homology $\int_{M} \overline{A} \in \CAlg(\calC)$. Let $\theta: \CAlg(\calC)^{\otimes} \rightarrow \calC^{\otimes}$ denote the forgetful functor. We wish to prove the existence of
a canonical equivalence $\theta( \int_{M} \overline{A})  \simeq \int_{M} \theta( \overline{A})$. 
In view of Proposition \ref{skoo}, it suffices to observe that $\theta$ is a symmetric monoidal
functor and that the underlying functor $\CAlg(\calC) \rightarrow \calC$ preserves sifted
colimits (Proposition \symmetricref{fillfemme}).
\end{proof}

If $M$ is an arbitrary $k$-manifold, we can view an $\OpE{M}$-algebra object of a symmetric monoidal
$\infty$-category $\calC$ as a family of $\OpE{k}$-algebras $\{ A_x \}_{x \in M}$ parametrized by the points of $M$. In general, this family is ``twisted'' by the tangent bundle of $M$. In the special case
where $M = S^1$, the tangent bundle $T_M$ is trivial, so we can think of an $\OpE{M}$-algebra as
a family of associative algebras parametrized by the circle: that is, as an associative algebra
$A$ equipped with an automorphism $\sigma$ (given by monodromy around the circle). Our final goal in this section is to show that in this case, the topological chiral homology $\int_{S^1} A$ coincides with
the Hochschild homology of the $A$-bimodule corresponding to $\sigma$. 

Fix an object of $\DISJ{S^1}$ corresponding to a single disk
$\psi: \R \hookrightarrow S^1$. An object of $\DISJ{S^1}_{\psi/}$ is given by a diagram
$$ \xymatrix{ \R \ar[dr]^{\psi} \ar[rr]^{j} & & U \ar[dl]^{\psi'} \\
& S^1 & }$$
which commutes up to isotopy, where $U$ is a finite union of disks. The set of components
$\pi_0( S^1 - \psi'(U) )$ is finite (equal to the number of components of $U$). Fix an orientation
of the circle. We define a linear ordering $\leq$ on $\pi_0(S^1 - \psi'(U) )$ as follows:
if $x, y \in S^1$ belong to different components of $S^1 - \psi'(U)$, then we write
$x < y$ if the three points $(x, y, \psi'( j(0) )$ are arranged in a clockwise order around the circle,
and $y < x$ otherwise. This construction determines a functor from
$\DISJ{S^1}_{\psi/}$ to (the nerve of) the category of nonempty finite linearly ordered sets,
which is equivalent to $\cDelta^{op}$. A simple calculation yields the following:

\begin{lemma}\label{snake}
Let $M = S^1$, and let $\psi: \R \hookrightarrow S^1$ be any open embedding. Then the above construction determines an equivalence of $\infty$-categories $\theta: \DISJ{M}_{\psi/}
\rightarrow \Nerve( \cDelta^{op} )$.
\end{lemma}

We can now formulate the relationship between Hochschild homology and topological chiral homology precisely as follows:

\begin{theorem}\label{junl}
Let $q: \calC^{\otimes} \rightarrow \Nerve(\FinSeg)$ be a sifted-complete symmetric monoidal category. Let $A \in \Alg_{ \OpE{S^1}}( \calC)$ be an algebra determining a diagram $\chi: \DISJ{S^1} \rightarrow \calC$
whose colimit is $\int_{S^1} A$. Choose an open embedding $\psi: \R \hookrightarrow S^1$. Then
the restriction $\chi | \DISJ{S^1}_{\psi/}$ is equivalent to a composition
$$ \DISJ{S^1}_{\psi/} \stackrel{\theta}{\rightarrow} \Nerve( \cDelta^{op}) \stackrel{B_{\bigdot}}{\rightarrow} \calC,$$
where $\theta$ is the equivalence of Lemma \ref{snake} and $B_{\bigdot}$ is a simplicial object
of $\calC$. Moreover, there is a canonical equivalence
$\int_{S^1} A \simeq | B_{\bigdot} |$.
\end{theorem}

\begin{lemma}\label{snoker}
Let $\calC$ be a nonempty $\infty$-category. Then $\calC$ is sifted if and only if, for each
object $C \in \calC$, the projection map $\theta_{C}: \calC_{C/} \rightarrow \calC$ is cofinal.
\end{lemma}

\begin{proof}
According to Theorem \toposref{hollowtt}, the projection map $\theta_{C}$ is cofinal if and only if,
for every object $D \in \calC$, the $\infty$-category $\calC_{C/} \times_{ \calC} \calC_{D/}$ is weakly
contractible. Using the evident isomorphism $\calC_{C/} \times_{\calC} \calC_{D/} \simeq
\calC \times_{ (\calC \times \calC)} (\calC \times \calC)_{(C,D)/}$, we see that this is equivalent
to the cofinality of the diagonal map $\calC \rightarrow \calC \times \calC$ (Theorem \toposref{hollowtt}).
\end{proof}

\begin{proof}[Proof of Theorem \ref{junl}]
The first assertion follows from Lemma \ref{snake}. The second follows from the
observation that $\DISJ{S^1}_{\psi/} \rightarrow \DISJ{S^1}$ is a cofinal map,
by virtue of Lemma \ref{snoker} and Proposition \ref{scun}.
\end{proof}

\begin{remark}
In the situation of Theorem \ref{junl}, let us view $A$ as an associative algebra object of $\calC$ equipped with an automorphism $\sigma$. We can describe the simplicial object $B_{\bigdot}$ 
informally as follows. For each $n \geq 0$, the object $B_n \in \calC$ can be identified with
the tensor power $A^{\otimes (n+1)}$. For $0 \leq i < n$, the $i$th face map from
$B_n$ to $B_{n-1}$ is given by the composition
$$ B_n \simeq A^{\otimes i} \otimes (A \otimes A) \otimes A^{ \otimes (n-1-i)}
\rightarrow A^{\otimes i} \otimes A \otimes A^{ \otimes (n-1-i)} \simeq B_{n-1},$$
where the middle map involves the multiplication on $A$. The $n$th face map is given instead by
the composition
$$ B_{n} \simeq (A \otimes A^{ \otimes (n-1)}) \otimes A
\simeq A \otimes (A \otimes A^{\otimes n-1}) \simeq 
( A \otimes A) \otimes A^{ \otimes n-1} \rightarrow A \otimes A^{\otimes^{n-1}} \simeq B_{n-1}.$$
\end{remark}

\begin{example}\label{staza}
Let $\calE$ denote the homotopy category of the $\infty$-operad
$\OpE{ \BTop(1)}$, so that $\Nerve( \calE)$ is the $\infty$-operad describing
associative algebras with involution (see Example \ref{spant}). Then
$\Nerve( \calE)$ contains a subcategory equivalent to the associative
$\infty$-operad $\Ass$. Since the circle $S^1$ is orientable, the canonical map
$ \OpE{S^1} \rightarrow \OpE{ \BTop(1) } \rightarrow \Nerve(\calE)$ factors through
this subcategory. We obtain by composition a functor
$\Alg_{ \Ass}(\calC) \rightarrow \Alg_{ \OpE{S^1}}(\calC)$ for any symmetric monoidal 
$\infty$-category $\calC$. If $\calC$ admits sifted colimits and the tensor product on $\calC$ preserves sifted colimits, we can then define the topological chiral homology $\int_{S^1} A$. It follows from
Theorem \ref{junl} that this topological chiral homology can be computed in very simple terms:
namely, it is given by the geometric realization of a simplicial object $B_{\bigdot}$ of
$\calC$ consisting of iterated tensor powers of the algebra $A$. In fact, in this case,
we can say more: the simplicial object $B_{\bigdot}$ can be canonically promoted to a {\it cyclic object}
of $\calC$. The geometric realization of this cyclic object provides the usual bar resolution for computing the Hochschild homology of $A$.
\end{example}

\subsection{Factorizable Cosheaves and Ran Integration}\label{stupem}

Let $M$ be a manifold and let $A$ be an $\OpE{M}$-algebra object of a sifted-complete
symmetric monoidal $\infty$-category $\calC^{\otimes}$. We refer to the 
object $\int_{U} A \in \calC$ introduced in Definition \ref{tch} as the topological chiral homology of $U$ with coefficients in $A$, which is intended to suggest that (like ordinary homology) it enjoys some form
of codescent with respect to open coverings in $M$. However, the situation is more subtle: the functor
$U \mapsto \int_{U} A$ is not generally a cosheaf on the manifold $M$ itself (except in the situation described in Lemma \ref{huuk}). However, it can be used to construct a cosheaf on the Ran space $\Ran(M)$ introduced in \S \ref{secran}. In other words, we can view topological chiral homology as given by the procedure of integration over the Ran space (Theorem \ref{selfmaid}).

We begin with a review of the theory of cosheaves.

\begin{definition}\label{plink}
Let $\calC$ be an $\infty$-category, $X$ a topological space, and $\calU(X)$ the partially
ordered set of open subsets of $X$. We will say that a functor $\calF: \Nerve( \calU(X) ) \rightarrow
\calC$ is a {\it cosheaf} on $X$ if, for every object $C \in \calC$, the induced map
$$ \calF_{C}: \Nerve( \calU(X) )^{op} \stackrel{\calF}{\rightarrow} \calC^{op} \stackrel{e_{C}}{\rightarrow} \SSet$$
is a sheaf on $X$, where $e_{C}: \calC^{op} \rightarrow \SSet$ denotes the functor represented by
$\calC$. We will say that a cosheaf
$\calF: \Nerve( \calU(X) ) \rightarrow \calC$ is {\it hypercomplete} if each
of the sheaves $\calF_{C} \in \Shv(X)$ is hypercomplete. If $X$ is the Ran space of a manifold
$M$, we will say that $\calF$ is {\it constructible}
if each of the sheaves $\calF_{C}$ is constructible in the sense of Definition \ref{plumb}. 
\end{definition}

\begin{remark}\label{splait}
Let $X$ be a topological space. It follows from Proposition \ref{coolman} that a functor
$\calF: \Nerve( \calU(X) ) \rightarrow \calC$ is a hypercomplete cosheaf on $X$
if and only if, for every open set $U \subseteq X$ and every functor
$f: \calJ \rightarrow \calU(U)$ with the property that
$\calJ_{x} = \{ J \in \calJ: x \in f(J) \}$ has weakly contractible nerve for each $x \in U$, 
the functor $\calF$ exhibits $\calF(U)$ as a colimit of the diagram
$\{ \calF( f(J) ) \}_{J \in \calJ}$.

In particular, if $g: \calC \rightarrow \calD$ is a functor which preserves small colimits, then
composition with $g$ carries hypercomplete cosheaves to hypercomplete cosheaves.
Similarly, if $\calC = \calP( \calE)$ for some small $\infty$-category $\calE$, a
functor $\calF: \Nerve( \calU(X) ) \rightarrow \calC$ is a hypercomplete cosheaf if and only if,
for every $E \in \calE$, the functor $U \mapsto \calF(U)(E)$ determines a cosheaf of spaces
$\Nerve( \calU(X)) \rightarrow \SSet$. 
\end{remark}

Our first goal in this section is to show that, if $M$ is a manifold, then we can identify
$\OpE{M}$-algebras with a suitable class of cosheaves on the Ran space $\Ran(M)$.
To describe this class more precisely, we need to introduce a bit of terminology.

\begin{definition}
Let $M$ be a manifold, and let $U$ be a subset of $\Ran(M)$. 
The {\it support} $\supp U$ of $U$ is the union $\bigcup_{S \in U} S \subseteq M$.
We will say that a pair of subsets $U, V \subseteq \Ran(M)$ are {\it independent} if
$\supp U \cap \supp V = \emptyset$.
\end{definition}

\begin{definition}
If $U$ and $V$ are subsets in $\Ran(M)$, we let $U \star V$
denote the set $\{ S \cup T: S \in U, T \in V \} \subseteq \Ran(M)$. 
\end{definition}

\begin{remark}
If $U$ is an open subset of $\Ran(X)$, then $\supp U$ is an open subset of $X$.
\end{remark}

\begin{example}
If $\{ U_i \}_{1 \leq i \leq n}$ is a nonempty finite collection of disjoint open subsets of a manifold $M$, then
the open set $\Ran( \{ U_i \} ) \subseteq \Ran(M)$ defined in \S \ref{secran} can be identified with
$\Ran(U_1) \caten \Ran(U_2) \caten \cdots \caten \Ran(U_n)$. 
\end{example}

\begin{remark}
If $U$ and $V$ are open in $\Ran(M)$, then $U \caten V$ is also open in $\Ran(M)$.
\end{remark}

\begin{remark}
We will generally consider the set $U \caten V$ only in the case where $U$ and $V$
are independent subsets of $\Ran(M)$. In this case, the canonical map $U \times V
\rightarrow U \caten V$ given by the formula
$(S, T) \mapsto S \cup T$ is a homeomorphism.
\end{remark}

\begin{definition}
Let $M$ be a manifold. We define a category $\Fact(M)^{\otimes}$ as follows:
\begin{itemize}
\item[$(1)$] The objects of $\Fact(M)^{\otimes}$ are finite sequences
$( U_1, \ldots, U_n )$ of open subsets $U_i \subseteq \Ran(M)$.
\item[$(2)$] A morphism from $(U_1, \ldots, U_m)$ to $(V_1, \ldots, V_n)$
in $\Fact{M}$ is a surjective map $\alpha: \seg{m} \rightarrow \seg{n}$ in $\FinSeg$ with
the following property: for $1 \leq i \leq n$, the sets $\{ U_j \}_{\alpha(j) = i }$
are pairwise independent and $\caten_{ \alpha(j) = i} U_j \subseteq V_{i}$.
\end{itemize}
We let $\Fact(M) \subseteq \Fact(M)^{\otimes}$ denote the fiber
product $\Fact(M)^{\otimes} \times_{ \FinSeg} \{ \seg{1} \}$, so that
$\Fact(M)$ is the category whose objects are open subsets of $\Ran(M)$ and whose morphisms
are inclusions of open sets.
\end{definition}

The $\infty$-category $\Nerve( \Fact(M)^{\otimes})$ is an $\infty$-operad.
Moreover, there is a canonical map of $\infty$-operads
$\Psi: \Nerve(\Disk{M})^{\otimes}_{\nunit} \rightarrow \Nerve( \Fact(M)^{\otimes})$, given on objects
by the formula $(U_1, \ldots, U_n) \mapsto ( \Ran(U_1), \ldots, \Ran(U_n) )$. 

We can now state our main result:

\begin{theorem}\label{manust}
Let $M$ be a manifold and let $\calC^{\otimes}$ be a symmetric monoidal $\infty$-category.
Assume that $\calC$ admits small colimits and that the tensor product on $\calC$
preserves small colimits separately in each variable. Then the operation of
operadic left Kan extension along the inclusion
$\Psi: \Nerve( \Disk{M})_{\nunit}^{\otimes} \rightarrow \Nerve( \Fact(M)^{\otimes})$ determines a fully faithful embedding
$F: \Alg_{ \Nerve(\Disk{M}) }^{\nunit}( \calC) \rightarrow \Alg_{ \Nerve(\Fact(M))}( \calC)$. Moreover, the
essential image of the full subcategory $\Alg_{ \Disk{M}}^{\nunit, \loc}(\calC)$ spanned by the locally constant objects of $\Alg_{ \Nerve( \Disk{M})}^{\unit}(\calC)$ is the
full subcategory of $\Alg_{ \Fact(M)}(\calC)$ spanned by those objects $A$
satisfying the following conditions:
\begin{itemize}
\item[$(1)$] The restriction of $A$ to $\Nerve( \Fact(M) )$ is a constructible cosheaf on
$\Ran(M)$, in the sense of Definition \ref{plink}.
\item[$(2)$] Let $U,V \subseteq \Ran(M)$ be independent open sets. Then the induced map
$A(U) \otimes A(V) \rightarrow A( U \caten V)$ is an equivalence in $\calC$.
\end{itemize}
\end{theorem}

\begin{remark}
In view of Proposition \ref{saz}, we can formulate Theorem \ref{manust} more informally as follows:
giving a nonunital $\OpE{M}$-algebra object of the $\infty$-category $\calC$ is equivalent to giving
a constructible $\calC$-valued cosheaf $\calF$ on the Ran space $\Ran(M)$, with the additional feature
that $\calF( U \caten V) \simeq \calF(U) \otimes \calF(V)$ when $U$ and $V$ are independent subsets of $\Ran(M)$. Following Beilinson and Drinfeld, we may sometimes refer a cosheaf with this property as
a {\it factorizable} cosheaf on $\Ran(M)$.
\end{remark}

The proof of Theorem \ref{manust} rests on the following basic calculation:

\begin{lemma}\label{corplam}
Let $M$ be a $k$-manifold, let $D \in \OpE{M}_{\nunit}$ be an object (corresponding to a nonempty finite
collection of open embeddings $\{ \psi_i: \R^{k} \rightarrow M \}_{1 \leq i \leq m}$), let
$\chi: \Nerve(\Disj{M}_{\nunit}) \rightarrow \SSet$ be a functor classified by the left fibration
$\Nerve(\Disj{M}_{\nunit}) \times_{ \OpE{M}_{\nunit} } ( \OpE{M}_{\nunit})^{\acti}_{D/}$
(here $\Disj{M}_{\nunit}$ is defined as in Proposition \ref{koio}), and let
$\overline{\chi}: \Nerve( \Fact(M) ) \rightarrow \SSet$ be a left Kan extension of $\chi$. Then
$\overline{\chi}$ is a hypercomplete $\SSet$-valued cosheaf on $\Ran(M)$.
\end{lemma}

\begin{proof}
Recall that a natural transformation of functors $\alpha: F \rightarrow G$ from
an $\infty$-category $\calC$ to $\SSet$ is said to be {\it Cartesian} if, for every
morphism $C \rightarrow D$ in $\calC$, the induced diagram
$$ \xymatrix{ F(C) \ar[r] \ar[d] & F(D) \ar[d] \\
G(C) \ar[r] & G(D) }$$
is a pullback square in $\SSet$. Let $D'$ be the image of $D$ in $\OpE{ \BTop(k) }$, and let $\chi': \Nerve(\Disj{M}_{\nunit}) \rightarrow \SSet$ be a functor classified by the left fibration $\Nerve(\Disj{M}_{\nunit}) \times_{ \OpE{ \BTop(k) } } ( \OpE{ \BTop(k) })^{\acti}_{/D'}$. There is an evident natural transformation of functors
$\beta: \chi \rightarrow \chi'$, which induces a natural transformation
$\overline{\beta}: \overline{\chi} \rightarrow \overline{\chi}'$. It is easy to see that
$\beta$ is a Cartesian natural transformation. 
Let $S = \{ 1, \ldots, m \}$, so that we can identify $\chi'$ with the functor which assigns
to $V \in \Disj{M}_{\nunit}$ the summand
$\Sing \Emb'( S \times \R^{k}, V) \subseteq
\Sing \Emb( S \times \R^{k}, V)$ consisting of those open embeddings
$j: S \times \R^{k} \rightarrow V$ which are surjective on connected components.
Let $\chi'': \Nerve(\Disj{M}_{\nunit}) \rightarrow \SSet$ be the functor given by the formula
$V \mapsto \Sing \Conf'( S, V)$, where
$\Sing \Conf'(S, V) \subseteq \Sing \Conf( S, V )$ is
the summand consisting of injective maps $i: S \rightarrow V$ which are
surjective on connected components. We have
an evident natural transformation of functors $\gamma: \chi' \rightarrow \chi''$. Using
Remark \ref{geck}, we deduce that $\gamma$ is Cartesian, so that $\alpha = \gamma \circ \beta$ is a Cartesian natural transformation from $\chi$ to $\chi''$. 

Let $\phi: \Conf(S, M) \rightarrow \Ran(M)$ be the continuous map which assigns
to each configuration $i: S \rightarrow M$ its image $i(S) \subseteq M$
(so that $\phi$ exhibits $\Conf(S,M)$ as a finite covering space of $\Ran^{m}(M) \subseteq \Ran(M)$).
Let $\overline{\chi}'': \Nerve( \Fact(M) ) \rightarrow \SSet$ be the functor given by the formula
$U \mapsto \Sing(\phi^{-1} U)$. We observe that $\overline{\chi}''$ is canonically equivalent
to $\chi''$. We claim that $\overline{\chi}''$ is a left Kan extension of $\chi''$. To prove this,
it suffices to show that for every open subset $U \subseteq \Ran(M)$, the map
$\overline{\chi}''$ exhibits $\Sing(\phi^{-1} U)$ as a colimit of the diagram
$\{ \chi''(V) \}_{V \in \calJ}$,
where $\calJ \subseteq \Disj{M}_{\nunit}$ is the full subcategory spanned by those
unions of disks $V = U_1 \cup \ldots \cup U_n$ such that
$\Ran( \{U_i \} ) \subseteq U$. For each
$x \in \phi^{-1}(U)$, let $\calJ_x$ denote the full subcategory of $\calJ$ spanned by those
open sets $V$ such that the map $x: S \rightarrow M$ factors through a map
$S \rightarrow V$ which is surjective on connected components.
In view of Theorem \ref{vankamp}, it will suffice to show that $\calJ_x$ has weakly contractible nerve.
In fact, we claim that $\calJ_x^{op}$ is filtered: this follows from the observation that every
open neighborhood of $x(S)$ contains an open set of the form $U_1 \cup \ldots \cup U_m$, where
the $U_i$ are a collection of small disjoint disks containing the elements of $x(S)$.

The map $\alpha$ induces a natural transformation
$\overline{\alpha}: \overline{\chi} \rightarrow \overline{\chi}''$.
Using Theorem \toposref{charleschar}, we deduce
that $\overline{\alpha}$ is also a Cartesian natural transformation. We wish to show
that $\overline{\chi}$ satisfies the criterion of Remark \ref{splait}. In other words, we wish
to show that if $U \subseteq \Ran(M)$ is an open subset and
$f: \calI \rightarrow \Fact(M)$ is a diagram such that each $f(I) \subseteq U$ and
the full subcategory $\calI_{x} = \{ I \in \calI: x \in f(I) \}$ has weakly contractible nerve
for each $x \in U$, then $\overline{\chi}$ exhibits
$\overline{\chi}(U)$ as a colimit of the diagram $\{ \overline{\chi}(f(I)) \}_{I \in \calI}$.
By virtue of Theorem \toposref{charleschar}, it will suffice to show that
$\overline{\chi}''$ exhibits $\overline{\chi}''(U)$ as a colimit of the diagram
$\{ \overline{\chi}''( f(I) ) \}_{I \in \calI}$. This is an immediate consequence
of Theorem \ref{vankamp}.
\end{proof}

\begin{proof}[Proof of Theorem \ref{manust}] 
The existence of the functor $F$ follows from Corollary \symmetricref{spaltwell}.
Let $A_0$ be
a nonunital $\Nerve(\Disk{M})^{\otimes}$-algebra object of $\calC$. Using Corollary \symmetricref{spaltwell},
Proposition \symmetricref{chocolateoperad}, and Proposition \symmetricref{optest}, we see
that $A = F(A_0)$ can be described as an algebra which assigns to each $U \subseteq \Ran(M)$ a colimit of the diagram
$$ \chi_{U}: \Nerve(\Disk{M})^{\otimes} \times_{ \Nerve( \Fact(M)^{\otimes}) }
\Nerve(\Fact(M)^{\otimes})^{\acti}_{/U} \rightarrow \calC.$$
The domain of this functor can be identified with the
nerve of the category $\calC_{U}$ whose objects are finite
collections of disjoint disks $V_1, \ldots, V_n \subseteq M$ such
that $\Ran( \{ V_i \}) \subseteq U$.
In particular, if $U = \Ran(U')$ for some open disk $U' \subseteq M$, then
the one-element sequence $(U')$ is a final object of $\calC_{U}$.
It follows that the canonical map $A_0 \rightarrow A | \Nerve( \Disk{M} )^{\otimes}$
is an equivalence, so that the functor $F$ is fully faithful.

We next show that if $A= F(A_0)$ for some $A_0 \in \Alg^{\loc,\nunit}_{\Nerve( \Disk{M})}(\calC)$,
then $A$ satisfies conditions $(1)$ and $(2)$. To prove that $A$ satisfies $(2)$, we observe that if $U, V \subseteq \Ran(M)$ are independent then we have a canonical equivalence $\calC_{U \caten V} \simeq \calC_{U} \caten \calC_{V}$.
Under this equivalence, the functor $\chi_{U \caten V}$ is given by the tensor product of
the functors $\chi_{U}$ and $\chi_{V}$. The map $A(U) \otimes A(V) \rightarrow A( U \caten V)$
is a homotopy inverse to the equivalence $$\varinjlim_{ \Nerve(\calC_{U \caten V})} \chi_{U \caten V}
\simeq \varinjlim_{ \Nerve(\calC_U) \times \Nerve(\calC_V)} \chi_U \otimes \chi_V
\rightarrow (\varinjlim_{ \Nerve(\calC_{U})} \chi_U) \otimes (\varinjlim_{ \Nerve(\calC_V)} \chi_V )$$
provided by our assumption that the tensor product on $\calC$ preserves small colimits separately in each variable.

We next show that $A| \Nerve( \Fact(M) )$ is a hypercomplete cosheaf on
$\Ran(M)$. By virtue of Proposition \ref{saz}, we can assume that $A_0$ factors as a composition
$$ \Disk{M}^{\otimes}_{\nunit} \rightarrow \OpE{M}_{\nunit} \stackrel{A'_0}{\rightarrow} \calC^{\otimes}.$$
Let $\calD$ be the subcategory of $\OpE{M}_{\nunit}$ spanned by the active morphisms.
As explained in \S \symmetricref{monenv}, the $\infty$-category $\calD$ admits a symmetric
monoidal structure and we may assume that
$A'_0$ factors as a composition
$$ \OpE{M}_{\nunit} \rightarrow \calD^{\otimes} \stackrel{A''_0}{\rightarrow} \calC^{\otimes},$$
where $A''_0$ is a symmetric monoidal functor. Corollary \symmetricref{banpan} implies
that the $\calP( \calD)$ inherits a symmetric monoidal structure, and that $A''_0$ factors
(up to homotopy) as a composition
$$ \calD^{\otimes} \rightarrow \calP(\calD)^{\otimes} \stackrel{T}{\rightarrow} \calC^{\otimes}$$
where $T$ is a symmetric monoidal functor such that the underlying
functor $T_{\seg{1}}:  \calP(\calD) \rightarrow \calC$ preserves small colimits. Let $B_0$ denote the composite map
$$ \Disk{M}^{\otimes}_{\nunit} \rightarrow \OpE{M}_{\nunit} \rightarrow \calD^{\otimes} \rightarrow \calP(\calD)^{\otimes},$$
and let $B \in \Alg_{ \Fact(M)}( \calP(\calD) )$ be an operadic left Kan extension of $B_0$,
so that $A_0 \simeq T \circ B_0$ and $A \simeq T \circ B$. Since $T_{\seg{1}}$ preserves
small colimits, it will suffice to show that $B | \Nerve( \Fact(M) )$ is a hypercomplete
$\calP(\calD)$-valued cosheaf on $\Ran(M)$ (Remark \ref{splait}). Fix an object
$D \in \calD$, and let $e_{D}: \calP(\calD) \rightarrow \SSet$ be the functor given by
evaluation on $D$. In view of Remark \ref{splait}, it will suffice to show that
$e_D \circ (B | \Nerve( \Fact(M) ))$ is a hypercomplete $\SSet$-valued cosheaf
on $\Ran(M)$. The desired result is now a translation of Lemma \ref{corplam}.

To complete the proof that $A$ satisfies $(1)$, it suffices to show that
for each $C \in \calC$, the functor $U \mapsto \bHom_{\calC}( A(U), C)$
satisfies condition $(\ast)$ of Proposition \ref{scamm}. Let $U_1, \ldots, U_n \subseteq M$
be disjoint disks containing smaller disks $V_1, \ldots, V_n \subseteq M$; it will
suffice to show that the corestriction map
$ A( \Ran( \{ V_i \} ) )
\rightarrow A( \Ran( \{ U_i \} ) )$
is an equivalence in $\calC$. Since $A$ satisfies $(2)$, we can reduce to the
case where $n=1$. In this case, we have a commutative diagram
$$ \xymatrix{ A_0(V_1) \ar[r]^{\beta} \ar[d] & A_0(U_1) \ar[d] \\
A( \Ran(V_1) ) \ar[r]^{\beta'} & A( \Ran(U_1) ).}$$
The vertical maps are equivalences (since $F$ is fully faithful), and the map
$\beta$ is an equivalence because $A_0$ is locally constant.

Now suppose that $A \in \Alg_{ \Nerve( \Fact(M) )}(\calC)$ satisfies conditions $(1)$ and
$(2)$; we wish to prove that $A$ lies in the essential image of $F |  \Alg^{\loc,\nunit}_{\Nerve( \Disk{M})}(\calC)$. Let $A_0 = A | 
\Nerve( \Disk{M})^{\otimes}$. Since $A$ satisfies $(1)$, Proposition \ref{scamm} guarantees that
$A_0$ is locally constant; it will therefore suffice to show that the canonical map
$F(A_0) \rightarrow A$ is an equivalence in the $\infty$-category
$\Alg_{ \Nerve(\Fact(M) )}(\calC)$. It will suffice to show that for every open set
$U \subseteq \Ran(M)$ and every object $C \in \calC$, the induced map
$\alpha_U: \bHom_{ \calC}( A(U), C) \rightarrow \bHom_{ \calC}( F(A_0)(U), C)$
is a homotopy equivalence of spaces. Since $A$ and $F(A_0)$ both satisfy
condition $(1)$, the collection of open sets $U$ such that $\alpha_U$ is a homotopy
equivalence is stable under unions of hypercovers. Consequently, Remark \ref{spow}
allows us to assume that $U = \Ran(V_1) \caten \cdots \caten \Ran(V_n)$ for some
collection of disjoint open disks $V_1, \ldots, V_n \subseteq M$. We claim that
$\beta: F(A_0)(U) \rightarrow A(U)$ is an equivalence. Since $A$ and $F(A_0)$ both satisfy
$(2)$, it suffices to prove this result after replacing $U$ by $\Ran(V_i)$ for $1 \leq i \leq n$.
We may therefore assume that $U = \Ran(V)$ for some open disk $V \subseteq M$.
In this case, we have a commutative diagram
$$ \xymatrix{ & A_0(V) \ar[dr]^{\beta''} & \\
F(A_0)(U) \ar[rr]^{\beta} \ar[ur]^{\beta'} & & A(U). }$$
The map $\beta'$ is an equivalence by the first part of the proof, and
$\beta''$ is an equivalence by construction. The two-out-of-three property shows that
$\beta$ is also an equivalence, as desired.
\end{proof}

The construction of topological chiral homology is quite closely related to the left
Kan extension functor $F$ studied in Theorem \ref{manust}. Let $M$ be a manifold,
let $A \in \Alg_{ \Nerve(\Disk{M})}(\calC)$, and let $A_0 = A | \Nerve( \Disk{M})^{\otimes}_{\nunit}$.
Evaluating $\Psi(A_0)$ on the Ran space $\Ran(M)$, we obtain an object of
$\calC$ which we will denote by $\int^{\nunit}_{M} A$. Unwinding the definition, we see that
$\int^{\nunit}_{M} A$ can be identified with the colimit $\varinjlim_{ V \in \Disj{M}_{\nunit} } \chi(V)$,
where $\chi: \Nerve(\Disj{M}) \rightarrow \calC$ is the functor given informally by the formula
$\chi( U_1 \cup \ldots \cup U_n) = A(U_1) \otimes \cdots \otimes A(U_n)$. The topological
chiral homology $\int_{M} A$ is given by the colimit $\varinjlim_{ V \in \Disj{M} } \chi(V)$.
The inclusion of $\Disj{M}_{\nunit}$ into $\Disj{M}$ induces a map
$\int^{\nunit}_{M} A \rightarrow \int_{M} A$. We now have the following result:

\begin{theorem}\label{selfmaid}
Let $M$ be a manifold and $\calC^{\otimes}$ a symmetric monoidal $\infty$-category. Assume that $\calC$ admits small colimits and that the tensor product on $\calC$ preserves colimits separately in each variable, and let $A \in \Alg_{ \Nerve( \Disk{M} )}(\calC)$. Suppose that $M$ is connected and
that $A$ is locally constant. Then the canonical map
$\int^{\nunit}_{M} A \rightarrow \int_{M} A$ is an equivalence in $\calC$.
\end{theorem}

\begin{proof}
The map $A$ determines a diagram $\psi: \Nerve( \Disj{M}) \rightarrow \calC$, given informally
by the formula $\psi( U_1 \cup \ldots \cup U_n) = A(U_1) \otimes \cdots \otimes A(U_n)$.
We wish to prove that the canonical map $\theta: \varinjlim( \psi | \Nerve( \Disj{M}^{\nunit}) )
\rightarrow \varinjlim(\psi)$ is an equivalence.
Since $A$ is locally constant, we can use Theorem \ref{sazz} to reduce to the case
where $A$ factors as a composition
$\Disk{M}^{\otimes} \rightarrow \OpE{M} \stackrel{A'}{\rightarrow} \calC^{\otimes}.$
In this case, $\psi$ factors as a composition
$\Nerve(\Disj{M}) \rightarrow \DISJ{M} \stackrel{\psi'}{\rightarrow} \calC,$
so we have a commutative diagram
$$ \xymatrix{ \varinjlim( \psi | \Nerve( \Disj{M}^{\nunit}) ) \ar[dr]^{\theta'} \ar[rr]^{\theta} & & \varinjlim( \psi) \ar[dl]^{\theta''} \\
& \varinjlim( \psi'). & }$$
Proposition \ref{koio} guarantees that $\theta'$ and $\theta''$ are equivalences in $\calC$, so that
$\theta$ is an equivalence by the two-out-of-three property.
\end{proof}

Theorem \ref{selfmaid} can be regarded as making the functor $\Psi$ of Theorem \ref{manust}
more explicit: if $A_0$ is a locally constant quasi-unital $\Nerve( \Disk{M})^{\otimes}$-algebra and $M$ is connected, then the global sections of the associated factorizable cosheaf can be computed by the topological chiral homology construction of Definition \ref{tch}. We can also read this theorem in the other direction. If $A$ is a locally constant $\Nerve(\Disk{M})^{\otimes}$-algebra, the the functor
$U \mapsto \int_{U} A$ does not determine a cosheaf $\Nerve( \calU(M) ) \rightarrow \calC$
in the sense of Definition \ref{plink}. However, when $U$ is connected, the topological chiral
homology $\int_{U} A$ {\em can} be computed as the global sections of a sheaf on the Ran space $\Ran(U)$. This is a reflection of a more subtle sense in which the construction $U \mapsto \int_{U} A$
behaves ``locally in $U$.'' We close this section with a brief informal discussion.

Let $M$ be a manifold of dimension $k$, and let $N \subseteq M$ be a submanifold
of dimension $k-d$ which has a trivial neighborhood of the form
$N \times \R^{d}$. Let $A \in \Alg_{ \OpE{M} }(\calC)$ and let
$\int(A)$ denote the associated $\Nerve( \calU(M) )^{\otimes}$-algebra object of $\calC$.
Restricting $\int(A)$ to open subsets of $M$ of the form $N \times V$, where $V$ is a union of finitely
many open disks in $\R^{d}$, we obtain another algebra $A_{N} \in \Alg_{ \Nerve( \Disk{\R^{d}})}(\calC)$. This algebra is locally constant, and can therefore
be identified with an $\OpE{d}$-algebra object of $\calC$ (Theorem \ref{sazz}). We will denote
this algebra by $\int_{N} A$. 

\begin{warning}
This notation is slightly abusive: the $\OpE{d}$-algebra $\int_{N} A$ depends not only
on the closed submanifold $N \subseteq M$ but also on a trivialization of a neighborhood
of $N$.
\end{warning}

Suppose now that $d=1$, and that $N \subseteq M$ is a hypersurface which separates
the connected manifold $M$ into two components. Let $M_{+}$ denote the union of one of these
components with the neighborhood $N \times \R$ of $N$, and $M_{-}$ the union of the other component with $N \times \R$ of $N$. After choosing appropriate conventions regarding the orientation of
$\R$, we can endow the topological chiral homology $\int_{M_{+}} A$ with the structure of a
right module over $\int_{N} A$ (which we will identify with an associative algebra object of $\calC$), and $\int_{M_{-}} A$ with the structure of a left module over $\int_{N} A$. There is a canonical map
$$ (\int_{M_{+}} A ) \otimes_{ \int_{N} A} ( \int_{M_{-}} A) \rightarrow \int_{M} A,$$
which can be shown to be an equivalence. In other words, we can recover the topological chiral
homology $\int_{M} A$ of the entire manifold $M$ if we understand the topological chiral homologies
of $M_{+}$ and $M_{-}$, together with their interface along the hypersurface $N$. 

Using more elaborate versions of this analysis, one can compute $\int_{M} A$ using any sufficiently nice decomposition of $M$ into manifolds with corners (for example, from a triangulation of $M$). This can be made precise using the formalism of extended topological quantum field theories (see \cite{cobordism} for a sketch). 

\begin{example}
Let $M = \R^{k}$, so that the $\infty$-operad $\OpE{M}$ is equivalent to $\OpE{k}$. Let
$N = S^{k-1}$ denote the unit sphere in $\R^{k}$. We choose
a trivialization of the normal bundle to $N$ in $M$, which assigns to each point $x \in S^{k-1} \subseteq \R^{k}$ the ``inward pointing'' normal vector given by $-x$ itself. According to the above discussion, we can associate to any algebra object $A \in \Alg_{ \OpE{k} }(\calC)$ an $\OpE{1}$-algebra object
of $\calC$, which we will denote by $B = \int_{ S^{k-1} } A$. Using Example \ref{sulta} and Proposition \symmetricref{algass}, we can identify $B$ with an associative algebra object of $\calC$.
One can show that this associative algebra has the following property: there is an equivalence
of $\infty$-categories $\theta: \Mod^{\OpE{k}}_{A}( \calC) \simeq \Mod^{L}_{B}(\calC)$
which fits into a commutative diagram of $\infty$-categories
$$ \xymatrix{ \Mod^{\OpE{k}}_{A}(\calC) \ar[rr]^{\theta} \ar[dr] & & \Mod^{L}_{B}(\calC) \ar[dl] \\
& \calC & }$$
which are right-tensored over $\calC$ (in view of Theorem \ref{postcur}, the existence of such a diagram characterizes the object $B \in \Alg(\calC)$ up to canonical equivalence). Under the equivalence $\theta$, the left $B$-module $B$ corresponds to the object
$F( {\bf 1}) \in \Mod^{\OpE{k}}_{A}(\calC)$ appearing in the statement of Theorem \ref{curtis}.
\end{example}

\subsection{Digression: Colimits of Fiber Products}\label{digdig}

If $\calX$ is an $\infty$-topos, then colimits in $\calX$ are universal: that is,
for every morphism $f: X \rightarrow Y$ in $\calX$, the fiber product construction
$Z \mapsto X \times_{Y} Z$ determines a colimit-preserving functor from
$\calX_{/Y}$ to $\calX_{/X}$. In other words, the fiber product
$X \times_{Y} Z$ is a colimit-preserving functor of $Z$. The same argument shows
that $X \times_{Y} Z$ is a colimit-preserving functor of $X$. However, the dependence of the fiber product $X \times_{Y} Z$ on $Y$ is more subtle. In \S \ref{nonpn}, we will need the following result, which asserts that
the construction $Y \mapsto X \times_{Y} Z$ commutes with colimits in many situations:

\begin{theorem}\label{kanter}
Let $\calX$ be an $\infty$-topos. Let $\calC$ denote the $\infty$-category
$\Fun( \Lambda^2_2, \calX) \times_{ \Fun( \{2\}, \calX)} \calX^{\geq 1}_{\ast}$
whose objects are diagrams $X \rightarrow Z \leftarrow Y$ in $\calX$, where
$Z$ is a pointed connected object of $\calC$. Let $F: \calC \rightarrow \calX$
be the functor
$$ \calC \rightarrow \Fun( \Lambda^2_2, \calX) \stackrel{\lim}{\rightarrow} \calX$$
given informally by the formula $(X \rightarrow Z \leftarrow Y) \mapsto X \times_{Z} Y$.
The $F$ preserves sifted colimits.
\end{theorem}

\begin{proof}[Proof of Theorem \ref{kanter}]
Let $\calC'$ denote the full subcategory of $\Fun( \Delta^1 \times \Delta^1 \times \Nerve( \cDelta^{op}_{+}), \calX)$
spanned by those functors $G$ which corresponding to diagrams of augmented simplicial objects 
$$ \xymatrix{ W_{\bigdot} \ar[r] \ar[d] & X_{\bigdot} \ar[d] \\
Y_{\bigdot} \ar[r] & Z_{\bigdot} }$$
which satisfy the following conditions:
\begin{itemize}
\item[$(i)$] The object $Z_{0}$ is final.
\item[$(ii)$] The augmentation map $Z_{0} \rightarrow Z_{-1}$ is an effective epimorphism
(equivalently, $Z_{-1}$ is a connected object of $\calX$).
\item[$(iii)$] Let $K$ denote the full subcategory of $\Delta^1 \times \Delta^1 \times \Nerve( \cDelta^{op}_{+})$ spanned by the objects
$(1,0,[-1])$, $(0,1,[-1])$, $(1,1,[-1])$, and $(1,1, [0])$. Then $G$ is a right Kan extension
of $G | K$. In particular, the diagram
$$ \xymatrix{ W_{-1} \ar[r] \ar[d] & X_{-1} \ar[d] \\
Y_{-1} \ar[r] & Z_{-1} }$$
is a pullback square.
\end{itemize}
It follows from Proposition \toposref{lklk} that the restriction map
$G \mapsto G|K$ induces a trivial Kan fibration $q: \calC' \rightarrow \calC$. Note that the functor
$F$ is given by composing a section of $q$ with the evaluation functor
$G \mapsto G(0,0,[-1])$. To prove that $F$ commutes with sifted colimits, it will suffice to show
that $\calC'$ is stable under sifted colimits in $\Fun( \Delta^1 \times \Delta^1 \times \Nerve( \cDelta^{op}_{+}), \calX)$.

Let $\calD$ be the full subcategory of $\Fun( \Delta^1 \times \Delta^1 \times \Nerve( \cDelta^{op}), \calX)$ spanned by those diagrams of simplicial objects
$$ \xymatrix{ W_{\bigdot} \ar[r] \ar[d] & X_{\bigdot} \ar[d] \\
Y_{\bigdot} \ar[r] & Z_{\bigdot} }$$
satisfying the following conditions:
\begin{itemize}
\item[$(i')$] The simplicial object $Z_{\bigdot}$ is a group object of $\calX$ (that is,
$Z_{\bigdot}$ is a groupoid object of $\calX$ and $Z_{0}$ is final in $\calX$; equivalently,
for each $n \geq 0$ the natural map $Z_{n} \rightarrow Z_1^{n}$ is an equivalence).
\item[$(ii')$] For each integer $n$ and each inclusion $[0] \hookrightarrow [n]$, the induced maps
$$X_{n} \rightarrow X_0 \times Z_n \quad Y_{n} \rightarrow Y_0 \times Z_{n}
\quad W_{n} \rightarrow X_0 \times Y_0 \times Z_{n}$$
are equivalences.
\end{itemize}

Since the product functor $\calX \times \calX \rightarrow \calX$ commutes with sifted colimits
(Proposition \toposref{urbil}), we deduce that $\calD$ is stable under sifted colimits
 in $\Fun( \Delta^1 \times \Delta^1 \times \Nerve( \cDelta^{op}), \calX)$. Let
 $\calD' \subseteq \Fun( \Delta^1 \times \Delta^1 \times \Nerve( \cDelta^{op}_{+}), \calX)$
 be the full subcategory spanned by those functors $G$ such that
 $G$ is a left Kan extension of $G_0 = G | ( \Delta^1 \times \Delta^1 \times \Nerve( \cDelta^{op}_{+}))$
 and $G_0 \in \calD$. Then $\calD'$ is stable under sifted colimits in
 $\Fun( \Delta^1 \times \Delta^1 \times \Nerve( \cDelta^{op}_{+}), \calX)$. We will complete
 the proof by showing that $\calD' = \calC'$.
 
 Suppose first that $G \in \calC'$, corresponding to a commutative diagram of augmented
 simplicial objects
 $$ \xymatrix{ W_{\bigdot} \ar[r] \ar[d] & X_{\bigdot} \ar[d] \\
Y_{\bigdot} \ar[r] & Z_{\bigdot}. }$$
Condition $(iii)$ guarantees that $Z_{\bigdot}$ is a \Cech nerve of the augmentation map
$Z_0 \rightarrow Z_{-1}$. Since this augmentation map is an effective epimorphism
(by virtue of $(ii)$), we deduce that the augmented simplicial object $Z_{\bigdot}$ is a colimit
diagram. Condition $(iii)$ guarantees that the natural maps $X_{n} \rightarrow Z_{n} \times_{ Z_{-1} } X_{-1}$ is are equivalences. Since colimits in $\calX$ are universal, we deduce that
$X_{\bigdot}$ is also a colimit diagram. The same argument shows that $Y_{\bigdot}$ and
$W_{\bigdot}$ are colimit diagrams, so that $G$ is a left Kan extension of
$G_0 = G| (\Delta^1 \times \Delta^1 \times \Nerve(\cDelta^{op}))$. To complete the proof
that $G \in \calD'$, it suffices to show that $G_0$ satisfies conditions $(i')$ and $(ii')$. Condition
$(i')$ follows easily from $(i)$ and $(iii)$, and condition $(ii')$ follows from $(iii)$.

Conversely, suppose that $G \in \calD'$; we wish to show that $G$ satisfies
conditions $(i)$, $(ii)$, and $(iii)$. Condition $(i)$ follows immediately from $(i')$,
and condition $(ii)$ from the fact that $Z_{\bigdot}$ is a colimit diagram. It remains
to prove $(iii)$. Let $K'$ denote the full subcategory of $\Delta^1 \times \Delta^1 \times \Nerve( \cDelta^{op}_{+})$ spanned by the objects $(0, 1, [-1)$, $(1,0, [-1])$, and
$\{ (1,1,[n] )\}_{n \geq -1}$. Since $\calX$ is an $\infty$-topos and $Z_{\bigdot}$ is
the colimit of a groupoid object of $\calX$, it is a \Cech nerve of the augmentation
map $Z_0 \rightarrow Z_{-1}$. This immediately implies that $G|K'$ is a right Kan extension
of $G|K$. To complete the proof, it will suffice to show that $G$ is a right Kan extension
of $G|K'$ (Proposition \toposref{acekan}).

We first claim that $G$ is a right Kan extension of $G|K'$ at $(0,1,[n])$ for each $n \geq 0$.
Equivalently, we claim that each of the maps
$$ \xymatrix{ X_{n} \ar[r] \ar[d] & X_{-1} \ar[d] \\
Z_{n} \ar[r] & Z_{-1} }$$
is a pullback diagram. Since $X_{\bigdot}$ and $Z_{\bigdot}$ are both colimit diagrams,
it will suffice to show that the map $X_{\bigdot} \rightarrow Z_{\bigdot}$ is a Cartesian transformation of simplicial objects (Theorem \toposref{charleschar}): in other words, 
it will suffice to show that for every morphism $[m] \rightarrow [n]$ in $\cDelta$, the analogous
diagram
$$ \xymatrix{ X_{n} \ar[r] \ar[d] & X_{m} \ar[d] \\
Z_{n} \ar[r] & Z_{m} }$$
is a pullback square. Choosing a map $[0] \hookrightarrow [m]$,
we obtain a larger diagram
$$ \xymatrix{ X_{n} \ar[r] \ar[d] & X_{m} \ar[r] \ar[d] & X_0 \ar[d] \\
Z_{n} \ar[r] & Z_{m} \ar[r] & Z_0. }$$
Since $Z_0$ is a final object of $\calX$, condition $(ii')$ implies that the right square and the outer
rectangle are pullback diagrams, so that the left square is a pullback diagram as well.
A similar argument shows that $Y_{\bigdot} \rightarrow Z_{\bigdot}$ and $W_{\bigdot} \rightarrow
Z_{\bigdot}$ are Cartesian transformations, so that $G$ is a right Kan extension of $G|K'$
at $(1,0,[n])$ and $(0,0,[n])$ for each $n \geq 0$.

To complete the proof, we must show that $G$ is a right Kan extension of $G|K'$ at
$(0,0,[-1])$: in other words, that the diagram $\sigma:$
$$ \xymatrix{ W_{-1} \ar[r] \ar[d] & X_{-1} \ar[d] \\
Y_{-1} \ar[r] & Z_{-1} }$$
is a pullback square. Since the map $\epsilon: Z_0 \rightarrow Z_{-1}$ is an effective epimorphism,
it suffices to show that the diagram $\sigma$ becomes a pullback square after base change
along $\epsilon$. In other words, we need only show that the diagram
$$ \xymatrix{ W_0 \ar[r] \ar[d] & X_0 \ar[d] \\
Y_0 \ar[r] & Z_0 }$$
is a pullback square, which follows immediately from $(ii')$.
\end{proof}

The remainder of this section is devoted to describing some applications of Theorem \ref{kanter}. The results here are not used elsewhere in this paper, and may be safely skipped by the reader.

\begin{corollary}\label{staru}
Let $\calX$ be an $\infty$-topos, $\Group(\calX)$ the $\infty$-category of group objects of $\calX$, and $K$ a sifted simplicial set. Suppose
we are given a pullback diagram
$$ \xymatrix{ W \ar[r] \ar[d] & X \ar[d] \\
Y \ar[r] & Z }$$
in the $\infty$-category $\Fun(K^{\triangleright}, \Group(\calX) )$ satisfying the following conditions:
\begin{itemize}
\item[$(i)$] The functors $X$, $Y$, and $Z$ are colimit diagrams.
\item[$(ii)$] For every vertex $v$ of $K$, the map $Y(k) \rightarrow Z(k)$ induces an
effective epimorphism in $\calX$.
\end{itemize}
Then $W$ is a colimit diagram in $\Group(\calX)$.
\end{corollary}

\begin{proof}
We observe that condition $(ii)$ is also satisfied when $v$ is the cone point of $K^{\triangleleft}$, since
the collection of effective epimorphisms in $\calX$ is stable under colimits.

Since $\calX$ is an $\infty$-topos, the formation of colimits determines an
equivalence of $\infty$-categories from $\Group(\calX)$ to the $\infty$-category $\calX_{\ast}^{\geq 1}$
of pointed connected objects of $\calX$. Applying this equivalence, we have a commutative
diagram $\sigma:$
$$ \xymatrix{ W' \ar[r] \ar[d] & X' \ar[d] \\
Y' \ar[r] & Z' }$$
of functors from $K^{\triangleleft}$ to $\calX_{\ast}$. Since the forgetful functor
$\Group(\calX) \rightarrow \calX$ is conservative and preserves sifted colimits, we deduce
that $X'$, $Y'$, and $Z'$ are colimit diagrams, and we wish to prove that $W'$ is a colimit diagram.
This follows from Theorem \ref{kanter}, provided that we can show that $\sigma$ is a pullback
square. The diagram $\sigma$ is evidently a pullback square in $\Fun(K^{\triangleleft}, \calX^{\geq 1}_{\ast})$, so it will suffice to show that the fiber product $X' \times_{Z'} Y'$ (formed in the
larger $\infty$-category $\Fun(K^{\triangleleft},\calX_{\ast})$) belongs to
$\Fun(K^{\triangleleft}, \calX^{\geq 1}_{\ast})$. In other words, we wish to show that for every
vertex $v \in K^{\triangleleft}$, the fiber product $X'(v) \times_{ Z'(v)} Y'(v)$ is a connected
object of $\calX$. Since the map $Y(v) \rightarrow Z(v)$ is an effective epimorphism, 
we deduce that its delooping $Y'(v) \rightarrow Z'(v)$ is $1$-connective. It
follows that the projection map $X'(v) \times_{ Z'(v)} Y'(v) \rightarrow X'(v)$ is $1$-connective.
The desired result now follows from the observation that $X'(v)$ is connected.
\end{proof}

\begin{corollary}\label{lym}
Let $\calX$ be an $\infty$-topos, let $X_{\bigdot}$ be a simplicial object in
the $\infty$-category $\Group(\calX)$. Then $X_{\bigdot}$ is a hypercovering
of its geometric realization $| X_{\bigdot} |$.
\end{corollary}

\begin{proof}
Without loss of generality, we may suppose that $\calX$ is the essential image of a left
exact localization functor $L: \calP(\calC) \rightarrow \calP(\calC)$, for some small $\infty$-category
$\calC$. We may assume without loss of generality that $X_{\bigdot} \simeq L Y_{\bigdot}$, for some
simplicial object $Y_{\bigdot}$ of $\Group( \calP(\calC) )$ (for example, we can take $Y_{\bigdot} = X_{\bigdot}$). Since $L: \calP(\calC) \rightarrow \calX$ preserves colimits, we have an equivalence
$L | Y_{\bigdot} | \simeq | X_{\bigdot} |$. Since $L$ preserves hypercoverings, it will suffice to show
that $Y_{\bigdot}$ is a hypercovering of $| Y_{\bigdot} |$. For this, we need only show that
$Y_{\bigdot}(C)$ is a hypercovering of $| Y_{\bigdot}(C) |$ in $\SSet$, for each object $C \in \calC$. 
In other words, we may assume that $\calX$ is the $\infty$-topos $\SSet$ of spaces.

Let $\bfA$ denote the category of simplicial groups, regarded as a simplicial model category;
we then have a canonical equivalence of $\infty$-categories $\Nerve( \bfA^{\degree}) \rightarrow
\Group(\SSet)$ (see Remark \ref{studious}). Let $\overline{X}_{\bigdot} \in \Fun( \Nerve(\cDelta_{+}^{op}), \Group(\SSet))$ be a colimit of $X_{\bigdot}$. Using Proposition \toposref{gumby444}, we may assume
that $\overline{X}_{\bigdot}$ is image of an augmented simplicial object $G: \cDelta_{+}^{op} \rightarrow \bfA$. We will identify $G$ with a simplicial object in the category $\bfA_{/ G([-1])}$. For every
simplicial set $K$, let $G(K) \in \bfA$ denote the limit $\varprojlim_{ \sigma \in \Hom_{\sSet}(\Delta^n,K)}
G( [n] )$, computed in the category $\bfA_{/ G([-1])}$. Without loss of generality, we may assume that
$G$ is Reedy fibrant. Then the map from $X_{n}$ to the matching object $M_n(X_{\bigdot})$
(computed in the $\infty$-category $\SSet_{/ | X_{\bigdot} |}$ can be identified with the map
$\theta: G( \Delta^n) \rightarrow G( \bd \Delta^n )$. Consequently, to prove that $X_{\bigdot}$ is a hypercovering $| X_{\bigdot} |$, it will suffice to show that the map
$\pi_0 G( \Delta^n) \rightarrow \pi_0 G( \bd \Delta^n)$ is surjective. Since $\theta$ is a Kan fibration
(by virtue of our assumption that $G$ is Reedy fibrant), this is equivalent to the requirement that
$\theta: G( \Delta^n)_0 \rightarrow G( \bd \Delta^n)_0$ is a surjection of groups.

Given an inclusion of simplicial sets $A \subseteq B$, we let $G(B,A)$ denote the kernel of the
restriction map $G(B) \rightarrow G(A)$. Since the map $\overline{X}_0 \rightarrow \overline{X}_{-1}$
is an effective epimorphism, the fibration $G(\Delta^0) \rightarrow G(\emptyset)$ is surjective on connected components and therefore induces a surjection $G(\Delta^0)_0 \rightarrow G(\emptyset)_0$.
Every nonempty simplicial set $K$ contains $\Delta^0$ as a retract, so that the map
$G(K)_0 \rightarrow G(\emptyset)_0$ is likewise surjective (this is evidently true also if $K = \emptyset$). We have a commutative diagram
$$ \xymatrix{ G( \Delta^n, \emptyset)_0 \ar[r] \ar[d]^{\theta'} & G( \Delta^n)_0 \ar[r] \ar[d]^{\theta} & G( \emptyset)_0 \ar[r] \ar[d] & 0 \ar[d]  \\
G( \bd \Delta^n, \emptyset)_0 \ar[r] & G( \bd \Delta^n)_0 \ar[r] & G( \emptyset)_0 \ar[r] & 0 }$$
with exact rows. Consequently, to prove that $\theta$ is surjective, it will suffice to show that
$\theta'$ is surjective. 

Let $H_{\bigdot}$ denote the simplicial group given by the formula
$H_{n} = G( \Delta^n, \emptyset)_0$. Then $H_{\bigdot}$ is automatically fibrant. 
Consequently, the map $G( \Delta^n, \emptyset)_0 \rightarrow G( \Lambda^n_0, \emptyset)_0$ is
surjective. We have a commutative diagram
$$ \xymatrix{ G( \Delta^n, \Lambda^n_0)_0 \ar[r] \ar[d]^{\theta''} & G( \Delta^n, \emptyset)_0 \ar[r] \ar[d]^{\theta'} & G( \Lambda^n_0, \emptyset)_0 \ar[r] \ar[d] & 0 \ar[d] \\
G( \bd \Delta^n, \Lambda^n_0)_0 \ar[r] & G( \bd \Delta^n, \emptyset) \ar[r] & G( \Lambda^n_0, \emptyset)_0 \ar[r] & 0 }$$
with exact rows. Consequently, to prove that $\theta'$ is surjective, it will suffice to show that
the map $$\theta'': G( \Delta^n, \Lambda^n_0)_0 \rightarrow
G( \bd \Delta^n, \Lambda^n_0)_0 \simeq G( \Delta^{ \{1, \ldots, n \}}, \bd \Delta^{ \{1, \ldots, n\}})_0$$ is surjective.
To complete the proof, it will suffice to verify the following:

\begin{itemize}
\item[$(\ast)$] Let $G: \cDelta_{+}^{op} \rightarrow \bfA$ be an augmented simplicial object
of the category $\bfA$ of simplicial groups. Assume that $G$ is Reedy fibrant and is a homotopy colimit diagram in $\bfA$. Then the map $\theta'': G( \Delta^n, \Lambda^n_0)_0 \rightarrow
G( \Delta^{ \{1, \ldots, n\} }, \bd \Delta^{ \{1, \ldots, n \}})_0$ is surjective. 
\end{itemize}

We will prove $(\ast)$ by induction on $n$. The case $n=0$ is obvious, since the group
$G( \Delta^{ \{1, \ldots, n \}}, \bd \Delta^{ \{1, \ldots, n\}})_0 \simeq G( \bd \Delta^n, \Lambda^n_0)_0$ is trivial. To handle the inductive step, let $TG$ denote the augmented simplicial group given by the formula
$TG([m]) = G([m] \star [0]) = G([m+1])$, and form a pullback diagram (in the category of augmented
simplicial objects of $\bfA$)
$$ \xymatrix{ G' \ar[r] \ar[d] & \ast \ar[d] \\
TG \ar[r] & G.}$$ 
Since each of the face maps $TG([m]) \simeq G([m+1]) \rightarrow G([m])$ is a fibration,
the above diagram is a homotopy pullback square. Note that $TG$ is a split augmented
simplicial object of $\bfA$, and therefore automatically a homotopy colimit diagram. 
For $n \geq 0$, the face map $TG([m]) \rightarrow G([m])$ admits a section, and therefore
determines an effective epimorphism in $\SSet$. Invoking Corollary \ref{staru}, we deduce
that $G'$ is a homotopy colimit diagram in $\bfA$. We have a commutative diagram
$$ \xymatrix{ G( \Delta^n, \Lambda^n_0) \ar[r] \ar[d] & G( \Delta^{ \{1, \ldots, n\}}, \bd \Delta^{ \{1, \ldots, n \}}) \ar[d] \\
G'( \Delta^{n-1}, \Lambda^{n-1}_0) \ar[r] & G'( \Delta^{ \{1, \ldots, n-1\}}, \bd \Delta^{ \{ 1, \ldots, n-1\} }) }$$
in which the vertical maps are isomorphisms of simplicial groups. The inductive hypothesis
guarantees that
$G'( \Delta^{n-1}, \Lambda^{n-1}_0)_0 \rightarrow G'( \Delta^{ \{1, \ldots, n-1\}}, \bd \Delta^{ \{ 1, \ldots, n-1\}})_0$ is surjective. This implies that the map
$\theta'': G( \Delta^n, \Lambda^n_0)_0 \rightarrow G(  \Delta^{ \{1, \ldots, n\}}, \bd \Delta^{ \{1, \ldots, n \}})_0$
is also surjective as required.
\end{proof}

\begin{remark}\label{postlym}
Let $\calC$ be the full subcategory of $\Fun( \Nerve( \cDelta_{+}^{op}, \SSet)$ spanned by those augmented simplicial spaces $X_{\bigdot}$ whose underlying simplicial space is a hypercovering
of $X_{-1}$. Then $\calC$ is stable under products in $\Fun( \Nerve( \cDelta_{+}^{op}, \SSet)$: this
follows from the observation that the collection of effective epimorphisms in $\SSet$ is stable under products.
\end{remark}

\begin{corollary}\label{commut}
Let $\Group(\SSet)$ denote the $\infty$-category of group objects of $\SSet$.
Let $F: \Fun( \Nerve(\cDelta)^{op}, \Group(\SSet))
\rightarrow \SSet$ denote the composition of the forgetful functor
$\Fun( \Nerve( \cDelta)^{op}, \Group(\SSet)) \rightarrow \Fun( \Nerve(\cDelta)^{op}, \SSet)$ with
the geometric realization functor $\Fun( \Nerve(\cDelta)^{op}, \SSet) \rightarrow \SSet$. Then
$F$ commutes with small products.
\end{corollary}

\begin{proof}
It suffices to show that the collection of augmented simplicial objects
of $\Group(\SSet)$ which determine colimit diagrams in $\SSet$ is stable
under products. This follows immediately from Corollary \ref{lym} together
with Remark \ref{postlym}.
\end{proof}

\begin{definition}
Let $\calC$ be an $\infty$-category which admits geometric realizations. We will say that
an object $P \in \calC$ is {\it strongly projective} if $P$ corepresents a functor
$e: \calC \rightarrow \SSet$ with the following property:
\begin{itemize}
\item[$(\ast)$] For every simplicial object $X_{\bigdot}$ in $\calC$, the
simplicial space $e( X_{\bigdot})$ is a hypercovering of $e( | X_{\bigdot} |)$.
\end{itemize}
\end{definition} 

\begin{remark}\label{jutly}
Every strongly projective object of an $\infty$-category $\calC$ is projective.
\end{remark}

\begin{example}\label{hurc}
Let $\calC$ be an $\infty$-category which admits geometric realizations, and let
$P$ be a cogroup object of $\calC$ (that is, a group object of the opposite $\infty$-category
$\calC^{op}$). Then $P$ is projective if and only if it is strongly projective. The
``if'' direction is obvious (Remark \ref{jutly}). For the converse, we observe that because
$P$ is a cogroup object, the functor $e: \calC \rightarrow \SSet$ corepresented by $P$ can be lifted
to a functor $\overline{e}: \calC \rightarrow \Group(\SSet)$. It follows from
Corollary \ref{lym} that $e$ carries every simplicial object $X_{\bigdot}$ of
$\calC$ to a hypercovering of $| e( X_{\bigdot})|$. If $e$ is projective, then this geometric realization
can be identified with $e | X_{\bigdot} |$.
\end{example}

\begin{proposition}
Let $\calC$ be an $\infty$-category which admits geometric realizations. Then the collection of strongly projective objects of $\calC$ is stable under all coproducts which exist in $\calC$.
\end{proposition}

\begin{proof}
This is an immediate consequence of Remark \ref{postlym}.
\end{proof}

We close this section by describing an application of Example \ref{hurc}: namely, we will show that it is possible to construct an analogue of the theory of projectively generated $\infty$-categories without assuming the generators to be compact.

\begin{definition}
Let $\kappa$ be a regular cardinal and let $\calC$ be a small $\infty$-category which admits $\kappa$-small coproducts. We let $\calP_{\Sigma}^{\kappa}$ denote the full subcategory of $\calP_{\Sigma}$
spanned by those presheaves $\calF: \calC^{op} \rightarrow \SSet$ which carry
$\kappa$-small coproducts in $\calC$ to products in $\SSet$.
\end{definition}

\begin{remark}
In the special case $\kappa = \omega$, we have $\calP_{\Sigma}^{\kappa}(\calC) = \calP_{\Sigma}(\calC)$.
\end{remark}

\begin{proposition}\label{sadir}
Let $\kappa$ be a regular cardinal. Let $\calC$ be a small $\infty$-category which admits
$\kappa$-small coproducts and satisfies the following additional condition:
\begin{itemize}
\item[$(\ast)$] Every object of $\calC$ can be regarded
as a cogroup object of $\calC$ (that is, a group object of the opposite $\infty$-category
$\calC^{op}$).
\end{itemize}
Then the full subcategory $\calP_{\Sigma}^{\kappa}( \calC) \subseteq \calP( \calC)$ is closed under the formation of geometric realizations of simplicial objects.
\end{proposition}

\begin{example}
Let $\calC$ be a pointed $\infty$-category which admits finite colimits. Then for every
object $X \in \calC$, the suspension $\Sigma(X)$ is a cogroup object of $\calC$.
In particular, if the suspension functor $\Sigma: \calC \rightarrow \calC$ is essentially
surjective (for example, if $\calC$ is stable), then $\calC$ satisfies condition $(\ast)$ of Proposition \ref{sadir}.
\end{example}

\begin{proof}
Let $X_{\bigdot}$ be a simplicial object of $\calP_{\Sigma}^{\kappa}(\calC)$, and let
$X$ denote the geometric realization $| X_{\bigdot} |$ formed in the $\infty$-category
$\calP(\calC)$. We wish to prove that $X \in \calP_{\Sigma}^{\kappa}(\calC)$. In other words,
we wish to show that if $\{ C_{\alpha} \}$ is a $\kappa$-small collection of objects of $\calC$
having a coproduct $C \in \calC$, then the natural map
$X(C) \rightarrow \prod_{\alpha} X(C_{\alpha})$ is an equivalence. In other words, we must show
that the map $\phi: | X_{\bigdot}(C) | \rightarrow \prod_{\alpha} | X_{\bigdot}(C_{\alpha})|$ is an equivalence. Since each $X_{n}$ belongs to $\calP_{\Sigma}^{\kappa}(\calC)$, we can identify
$\phi$ with the natural map
$$ | \prod_{\alpha} X_{\bigdot}(C_{\alpha})| \rightarrow \prod_{\alpha} | X_{\bigdot}(C_{\alpha}) |.$$
Since each $C_{\alpha}$ is a cogroup object of $\calC$ and each $X_{n}$ carries
finite coproducts to finite products, we deduce that each $X_{\bigdot}(C_{\alpha})$ can be
identified with a simplicial object in the $\infty$-category $\Group(\SSet)$ of group objects of $\SSet$.
The desired result now follows from Corollary \ref{commut}.
\end{proof}

\begin{proposition}\label{snaffer}
Let $\kappa$ be a regular cardinal, and let $\calC$ be a small $\infty$-category which admits
$\kappa$-small coproducts. Assume that every object of $\calC$ has the structure of a cogroup object of $\calC$. Then:
\begin{itemize}
\item[$(1)$] The $\infty$-category $\calP^{\kappa}_{\Sigma}(\calC)$ is an accessible
localization of $\calP_{\Sigma}(\calC)$. In particular, $\calP^{\kappa}_{\Sigma}(\calC)$ is a presentable $\infty$-category.
\item[$(2)$] Let $j: \calC \rightarrow \calP(\calC)$ be the Yoneda embedding.
Then $j$ factors through $\calP^{\kappa}_{\Sigma}(\calC)$.
\item[$(3)$] The functor $j: \calC \rightarrow \calP^{\kappa}_{\Sigma}(\calC)$ preserves $\kappa$-small coproducts.
\item[$(4)$] The essential image of $j$ consists of projective cogroup objects of
$\calP^{\kappa}_{\Sigma}(\calC)$ (and, in particular, of strongly projective objects
of $\calP^{\kappa}_{\Sigma}(\calC)$: Example \ref{hurc}).
\item[$(5)$] Let $X$ be an arbitrary object of $\calP^{\kappa}_{\Sigma}(\calC)$. Then
$X$ can be written as the geometric realization of a simplicial object
$X_{\bigdot}$ of $\calP^{\kappa}_{\Sigma}(\calC)$, where each $X_{n}$ is a small
coproduct (in $\calP^{\kappa}_{\Sigma}(\calC)$) of objects lying in the essential image of
$j$. In particular, each $X_{n}$ is a projective cogroup object of $\calC$.
\item[$(6)$] An object $X \in \calP^{\kappa}_{\Sigma}(\calC)$ is projective if and only if
$X$ can be obtained as a retract of some coproduct $\coprod_{\alpha} j( C_{\alpha})$.
\end{itemize}
\end{proposition}

\begin{proof}
Assertion $(1)$, $(2)$, and $(3)$ are obvious. Let $C \in \calC$. The projectivity
of the object $j(C) \in \calP^{\kappa}_{\Sigma}(\calC)$ follows from Proposition \ref{sadir}.
By assumption, $C$ has the structure of a cogroup object of $\calC$.
Because $j$ preserves finite coproducts (by virtue of $(3)$), we conclude that $j(C)$ has the structure of a cogroup in $\calP^{\kappa}_{\Sigma}(\calC)$. This proves $(4)$. Assertion $(5)$ follows immediately from $(1)$ together with Lemma \toposref{presubato}. The ``if'' direction of $(6)$ follows from
$(5)$ (since the collection of projective objects of $\calP^{\kappa}_{\Sigma}(\calC)$ is stable
under retracts). Conversely, suppose that $X$ is a projective object of $\calP^{\kappa}_{\Sigma}(\calC)$, and let $X_{\bigdot}$ be as in assertion $(5)$. Since $X$ is projective, the identity map
$\id_{X}: X \rightarrow X \simeq | X_{\bigdot} |$ factors (up to homotopy) through some map
$X \rightarrow X_0$. This exhibits $X$ as a retract of $X_0 \simeq \coprod_{\alpha} j(C_{\alpha})$.
\end{proof}

\begin{remark}
In the situation of Proposition \ref{snaffer}, the subcategory $\calP^{\kappa}_{\Sigma}(\calC)$ is stable
under small $\kappa$-filtered colimits in $\calP(\calC)$. It follows that the essential image of the Yoneda embedding $j: \calC \rightarrow \calP^{\kappa}_{\Sigma}(\calC)$ consists of $\kappa$-compact
objects of $\calP^{\kappa}_{\Sigma}(\calC)$. Conversely, suppose that $X$ is a $\kappa$-compact projective object of $\calP^{\kappa}_{\Sigma}(\calC)$. It follows from Proposition \ref{snaffer} that
$X$ is a retract of a small coproduct $\coprod_{\alpha \in A} j( C_{\alpha})$. We can write this
coproduct as a $\kappa$-filtered colimit of coproducts $\coprod_{ \alpha \in A_0} j( C_{\alpha})$, where
$A_0$ ranges over the $\kappa$-small subsets of $A$. Since $X$ is $\kappa$-compact, it follows that
$X$ is a retract of some coproduct $\coprod_{\alpha \in A_0} j(C_{\alpha}) \simeq j( \coprod_{\alpha \in A_0} C_{\alpha})$: in other words, $X$ belongs to the idempotent completion of the essential image of $j$. In particular, every projective object of $\calP^{\kappa}_{\Sigma}(\calC)$ is strongly projective.
\end{remark}

\subsection{Nonabelian Poincare Duality}\label{nonpn}

Let $M$ be an oriented $k$-manifold. Poincare duality provides a canonical isomorphism
$$ \HH_{c}^{m}(M; A) \simeq \HH_{k-m}(M; A)$$
for any abelian group $A$ (or, more generally, for any local system of abelian groups on $M$).
Our goal in this section is to establish an analogue of this statement for {\em nonabelian} cohomology: that is, cohomology with coefficients in a local system of {\em spaces} on $M$. To formulate this analogue, we will need to replace the right hand side by the topological chiral homology
$\int_{M} A$ of $M$ with coefficients in an approparite $\OpE{M}$-algebra.

\begin{remark}
The ideas described in this section are closely related to results of Segal, McDuff, and Salvatore on configuration spaces (see \cite{segalconfig}, \cite{mcduffconfig}, and \cite{salvatoreconfig}). In particular, a special case of our main result (Theorem \ref{stager}) can be found in \cite{salvatoreconfig}.
\end{remark}

\begin{definition}\label{julka}
Let $M$ be a manifold, and let $p: E \rightarrow M$ be a Serre fibration
equipped with a distinguished section $s: M \rightarrow E$.
Given a commutative diagram
$$ \xymatrix{ | \Delta^n | \times M \ar[dr] \ar[rr]^{f} & & E \ar[dl]^{p} \\
& M, & }$$
we will say that $f$ is {\it trivial} on an open set $U \subseteq M$ if
the restriction $f| (| \Delta^n | \times U)$ is given by the composition
$$| \Delta^n | \times U \rightarrow U \subseteq M \stackrel{s}{\rightarrow} E.$$
We define the {\it support} of $f$ to be the smallest closed set $K$ such that
$f$ is trivial on $M - K$. 
Given an open set $U \subseteq M$, we let
$\Sect(U;E)$ denote the simplicial set whose $n$-simplices are maps
$f$ as above, and $\Sect_{c}(U;E)$ the simplicial subset spanned by those simplices
such that the support of $f$ is a compact subset of $U$
(in this case, $f$ is determined by its restriction $f | (| \Delta^n| \times U)$

The construction $(U_1, \ldots, U_n) \mapsto \Sect_c(U_1;E) \times \ldots
\times \Sect_c(U_n; E)$ determines a functor from $\calU(M)^{\otimes}$ to the simplicial category of Kan complexes. Passing to nerves, we obtain a functor $\Nerve( \calU{M}^{\otimes}) \rightarrow \SSet$,
which we view as a $\Nerve( \calU(M)^{\otimes})$-monoid object of $\SSet$. Let us regard
the $\infty$-category $\SSet$ as endowed with the Cartesian monoidal structure, so that
this monoid object lifts in an essentially unique way to a $\Nerve( \calU(M)^{\otimes})$-algebra object
of $\SSet$ (Proposition \symmetricref{ungbatt}). We will denote this algebra by
$E_{!}$.
\end{definition}

\begin{remark}\label{staag}
Let $p: E \rightarrow M$ be as in Definition \ref{julka}. Every inclusion of
open disks $U \subseteq V$ in $M$ is isotopic to a homeomorphism (Theorem \ref{scen}), so the
inclusion $\Sect_{c}(U;E) \rightarrow \Sect_{c}(V;E)$ is a homotopy equivalence.
It follows that the restriction $E_{!} | \Nerve( \Disk{M})^{\otimes}$ is a locally constant object of
$\Alg_{ \Nerve( \Disk{M})}( \SSet)$, and is therefore equivalent to the restriction
$E^{!} | \Nerve( \Disk{M})^{\otimes}$ for some essentially unique
$\OpE{M}$-algebra $E^{!} \in \Alg_{ \OpE{M}}(\SSet)$ (Theorem \ref{sazz}).
\end{remark}

\begin{remark}
Let $M$ be a manifold and let $p: E \rightarrow M$ be a Serre fibration equipped with a section
$s$. Then the functor $U \mapsto \Sect(U;E)$ determines a sheaf $\calF$ on $M$ with
values in the $\infty$-category $\SSet_{\ast}$ of pointed spaces (Proposition \toposref{aese}). 
Using Remark \ref{kor} and Lemma \ref{kora}, we can identify the functor
$U \mapsto \Gamma_c(U; \calF)$ of Definition \ref{spose} with the functor $U \mapsto \Sect_{c}(U; E)$.
\end{remark}

\begin{remark}
Let $p: E \rightarrow M$ be as in Definition \ref{julka}. Since $p$ is a Serre fibration,  the inverse image
$U \times_{M} E$ is weakly homotopy equivalent to a product $U \times K$ for every
open disk $\R^{k} \simeq U \subseteq M$, for some pointed topological space $K$. For every positive
real number $r$, let $X_{r}$ denote the simplicial subset of $\Sect_c(U; E)$ whose $n$-simplices
correspond to maps which are supported in the closed ball $\overline{ B(r)} \subseteq \R^{k} \simeq U$.
Then each $X_{r}$ is homotopy equivalent to the iterated loop space $\Sing( \Omega^{k} K)$. Since there exist compactly supported isotopies of
$\R^{k}$ carrying $B(r)$ to $B(s)$ for $0 < r < s$, we deduce that the inclusion
$X_{r} \subseteq X_{s}$ is a homotopy equivalence for each $r < s$. It follows that
$\Sect_{c}(U; E) = \varinjlim_{r} X_{r}$ is weakly homotopy equivalent to
$X_{r}$ for every real number $r$. 

In other words, we can think of $E^{!}: \OpE{M} \rightarrow \SSet^{\times}$ as an
algebra which assigns to each open disk $j: U \hookrightarrow M$ the $k$-fold
loop space of $F$, where $F$ is the fiber of the Serre fibration $p: E \rightarrow M$ over any point
in the image of $j$.
\end{remark}

We can now state our main result as follows:

\begin{theorem}[Nonabelian Poincare Duality]\label{stager}
Let $M$ be a $k$-manifold, and let $p: E \rightarrow M$ be a Serre fibration whose fibers
are $k$-connective, which is equipped with a section $s: M \rightarrow E$. Then
$E_{!}$ exhibits $\Sect_{c}(M;E)$ as the colimit of the diagram
$E_{!} | \Nerve(\Disj{M})$. In other words, $\Sect_c(M;E)$ is the topological chiral homology
$\int_{M} E^!$, where $E^{!} \in \Alg_{ \OpE{M}}(\SSet)$ is the algebra described in
Remark \ref{staag}.
\end{theorem}

\begin{remark}
The assumption that $p: E \rightarrow M$ have $k$-connective fibers is essential.
For example, suppose that $E = M \coprod M$ and that the section
$s: M \rightarrow E$ is given by the inclusion of the second factor. If
$M$ is compact, then the inclusion of the second factor determines a vertex
$\eta \in \Sect_c(M; E)$. The support of $\eta$ is the whole of the manifold $M$:
in particular, $\eta$ does not lie in the essential image of any of the extension maps
$i: \Sect_c(U;E) \rightarrow \Sect_c(M;E)$ where $U$ is a proper open subset of $M$.
In particular, if $U$ is a disjoint union of open disks, then $\eta$ cannot lie in the essential image of $i$ unless $k=0$ or $M$ is empty.
\end{remark}

\begin{remark}
Theorem \ref{stager} implies in particular that any compactly supported section $s'$ of 
$p: E \rightarrow M$ is homotopic to a section whose support is contained in the union
of disjoint disks in $M$. It is easy to see this directly,
at least when $M$ admits a triangulation. Indeed, let $M_0 \subseteq M$ be the $(k-1)$-skeleton of this triangulation, so that the open set $M - M_0$ consists of the interiors of the $k$-simplices of the triangulation and is thus a union of disjoint open disks in $M$. Since the fibers of
$p$ are $k$-connective, the space of sections of $p$ over the $(k-1)$-dimensional space
$M_0$ is connected. Consequently, we can adjust $s'$ by a homotopy so that it agrees with
$s$ on a small neighborhood of $M_0$ in $M$, and is therefore supported in $M - M_0$.
\end{remark}

\begin{example}
Let $M$ be the circle $S^1$, let $X$ be a connected pointed space, and let $E = X \times S^1$,
equipped with the projection map $p: E \rightarrow M$.
Then $E^! \in \Alg_{ \OpE{S^1}}(\SSet)$ is the $\OpE{S^1}$-algebra determined by
the associative algebra object $\Sing(\Omega X) \in \Alg_{\Ass}(\SSet)$. Since $M$ is compact, we can identify
$\Sect_c(S^1; E)$ with the singular complex of the space $LX = \bHom(S^1,X)$ of {\em all} sections of $p$. In view of Example \ref{staza}, Theorem \ref{stager} recovers the following classical observation: the free loop space $LX$ is equivalent to the Hochschild homology of
the based loop space $\Omega X$ (regarded as an associative algebra with respect to composition of loops).
\end{example}

\begin{remark}
Let $M$ be a $k$-manifold. We will say than an algebra $A \in \Alg_{ \OpE{M}}( \SSet)$ is
{\it grouplike} if, for every open disk $U \subseteq M$, the restriction
$A | \OpE{U} \in \Alg_{ \OpE{U} }(\SSet) \simeq \Alg_{ \OpE{k}}(\SSet)$ is
grouplike in the sense of Definition \ref{ungwar} (by convention, this condition is vacuous if $k=0$). 
For every fibration $E \rightarrow M$, the associated algebra $E^{!} \in \Alg_{ \OpE{M}}( \SSet)$ is grouplike. In fact, the converse holds as well: every grouplike object of $\Alg_{ \OpE{M}}( \SSet)$
has the form $E^{!}$, for an essentially unique Serre fibration $E \rightarrow \OpE{M}$ with
$k$-connective fibers. 

To prove this, we need to introduce a bit of notation. For each open set $U \subseteq M$, let
$\calA_{U}$ denote the simplicial category whose objects are Serre fibrations
$p: E \rightarrow U$ equipped with a section $s$, where the pair
$(U,E)$ is a relative CW complex and the fibers of $p$ are $k$-connective; an $n$-simplex of $\bHom_{\calC_{U}}( E, E')$ is a commutative diagram
$$ \xymatrix{ E \times \Delta^n \ar[dr]^{p} \ar[rr]^{f} & & E' \ar[dl]^{p'} \\
& U, & }$$
such that $f$ respects the preferred sections of $p$ and $p'$.
Let $\calB_{U}$ denote the full subcategory of $\Alg_{ \OpE{U}}( \SSet)$ spanned by the grouplike objects. The construction $E \mapsto E^{!}$ determines a functor $\theta_U: \Nerve( \calA_{U} )
\rightarrow \calB_{U}$, which we claim is an equivalence of $\infty$-categories.
If $U \simeq \R^{k}$ is an open disk in $M$, then this assertion follows from Theorem \ref{preslage}
(at least if $k > 0$; the case $k=0$ is trivial). Let
$\calJ$ denote the collection of all open subsets $U \subseteq M$ which are homeomorphic
to $\R^{k}$, partially ordered by inclusion. This collection of open sets satisfies the following condition:
\begin{itemize}
\item[$(\ast)$] For every point $x \in M$, the subset
$\calJ_x = \{ U \in \calJ: x \in U \}$ has weakly contractible nerve (in fact,
$\calJ_x^{op}$ is filtered, since every open subset of $M$ containing $x$ contains an open
disk around $x$). 
\end{itemize}
We have a commutative diagram of $\infty$-categories
$$ \xymatrix{ \Nerve( \calA_M) \ar[r]^{\theta_M} \ar[d]^{\phi} & \calB_M \ar[d]^{\psi} \\
\varprojlim_{U \in \calJ^{op}} \Nerve( \calA_U) \ar[r] & \varprojlim_{U \in \calJ^{op}} \calB_U }$$
(here the limits are taken in the $\infty$-category $\Cat_{\infty}$). Here the lower horizontal map is an equivalence of $\infty$-categories. Consequently, to prove that $\theta_M$ is an equivalence of $\infty$ categories, it suffices to show that the vertical maps are equivalences of $\infty$-categories. We consider each in turn.

For each $U \subseteq M$, let $\calC_U$ denote the simplicial category whose objects
are Kan fibrations $p: X \rightarrow \Sing(U)$. The functor
$E \mapsto \Sing(E)$ determines an equivalence of $\infty$-categories
$\Nerve( \calA_U) \rightarrow \Nerve( \calC_U)_{\ast}$.
Consequently, to show that $\phi$ is a categorical equivalence, it will suffice to show that the associated map $\Nerve( \calC_M) \rightarrow \varprojlim_{U \in \calJ^{op} } \Nerve( \calC_U)$ is a categorical
equivalence. This is equivalent to the requirement that $\Sing(M)$ is a colimit of the diagram
$\{ \Sing(U) \}_{U \in \calJ}$ in the $\infty$-category $\SSet$, which follows from $(\ast)$ and
Theorem \ref{vankamp}.

To prove that $\psi$ is a categorical equivalence, it suffices to show that $\Alg_{ \OpE{M}}(\SSet)$
is a limit of the diagram $\{ \Alg_{ \OpE{U}}( \SSet) \}_{U \in \calJ^{op}}$. For each
$U \subseteq M$, let $\calD_{U}$ denote the $\infty$-category
$\Alg_{\calO_U}(\SSet)$, where $\calO_U^{\otimes}$ denotes the $\infty$-operad family $\OpE{ \BTop(k) } \times_{ \BTop(k)^{\amalg} } (B_U \times \Nerve(\FinSeg))$. It
follows that the restriction functor $A \mapsto A | \calO^{\otimes}$ determines
an equivalence of $\infty$-categories $\Alg_{\OpE{U}}(\SSet) \rightarrow \calD_U$. It will
therefore suffices to show that $\calD_M$ is a limit of the diagram of $\infty$-categories
$\{ \calD_U \}_{ U \in \calJ^{op} }$. To prove this, we show that the functor
$U \mapsto \calO_{U}^{\otimes}$ exhibits the $\infty$-operad family $\calO^{\otimes}_{M}$
as a homotopy colimit of the $\infty$-operad families $\{ \calO^{\otimes}_{U} \}_{U \in \calJ}$.
For this, it is sufficient to show that the Kan complex $B_M$ is a homotopy colimit of the diagram
$\{ B_U \}_{U \in \calJ}$, which follows from Remark \ref{talrod}, $(\ast)$, and Theorem \ref{vankamp}.
\end{remark}

\begin{remark}\label{conswe}
In proving Theorem \ref{stager}, it is sufficient to treat the case where the manifold $M$ is connected.
To see this, we note that for every open set $U \subseteq M$, we have a map
$\theta_U: \int_{U} E^{!} \rightarrow \Sect_c(U; E)$. Assume that $\theta_U$ is a homotopy equivalence whenever $U$ is connected. We will prove that $\theta_U$ is a homotopy equivalence
whenever the set of connected components $\pi_0(U)$ is finite. It will then follow that $\theta_{U}$ is an equivalence for every open set $U \subseteq M$, since the construction $U \mapsto \theta_{U}$ commutes with filtered colimits; in particular, it will follow that $\theta_{M}$ is a homotopy equivalence.

To carry out the argument, let $U \subseteq M$ be an open set with finitely many connected components $U_1, \ldots, U_n$, so that we have a commutative diagram
$$ \xymatrix{ \prod_{1 \leq i \leq n} \int_{U_{i}} E^{!} \ar[rr]^-{ \theta_{U_1} \times \cdots \times \theta_{U_{n}}} \ar[d]^{\phi} & & \prod_{1 \leq i \leq n} \Sect_c( U_i, E) \ar[d]^{\psi} \\
 \int_{U} E^{!} \ar[rr]^-{\theta_{U} } & & \Sect_c(U,E). }$$
The map $\theta_{U_{1}} \times \cdots \times \theta_{U_{n}}$ is a homotopy equivalence since each $U_i$ is connected, the map $\phi$ is a homotopy equivalence by Theorem \ref{clapse}, and the
map $\psi$ is an isomorphism of Kan complexes; it follows that $\theta_{U}$ is a homotopy equivalence
as desired.
\end{remark}

\begin{notation}
Let $p: E \rightarrow M$ be as in Definition \ref{julka}. Given a compact set $K \subseteq M$, 
we let $\Sect_K(M; E)$ denote the simplicial set whose
$n$-simplices are commutative diagrams
$$ \xymatrix{ & (| \Delta^n | \times M) \coprod_{ | \Delta^n | \times (M-K) \times \{0\} }
( |\Delta^n| \times (M-K) \times [0,1]) \ar[dr] \ar[dl]^{f} & \\
E \ar[rr]^{p} & & M }$$
such that $f | ( | \Delta^n | \times (M-K) \times \{1\} )$ is given by the composition
$$ | \Delta^n| \times V \times \{1\} \rightarrow (M-K) \subseteq M \stackrel{s}{\rightarrow} E.$$
In other words, an $n$-simplex of $\Sect_K(M;E)$ is an $n$-parameter family of
sections of $E$, together with a nullhomotopy of this family of sections on the open set $M-K$.

Note that any $n$-simplex of $\Sect_c(M;E)$ which is trivial on $M-K$ extends canonically
to an $n$-simplex of $\Sect_{K}(M; E)$, by choosing the nullhomotopy to be constant.
In particular, if $U \subseteq M$ is any open set, then we obtain a canonical map
$$ \Sect_c(U;E) \rightarrow \varinjlim_{K \subseteq U} \Sect_{K}(M; E),$$
where the colimit is taken over the (filtered) collection of all compact subsets of $U$.
\end{notation}

\begin{remark}\label{kor}
The simplicial set $\Sect_{K}(M;E)$ can be identified with the homotopy fiber of the
restriction map $\calF(M) \rightarrow \calF(M-K)$, where $\calF \in \Shv(M)$ is the sheaf
associated to the fibration $p: E \rightarrow M$.
\end{remark}

\begin{lemma}\label{kora}
Let $p: E \rightarrow M$ be a Serre fibration equipped with a section $s$ (as in Definition \ref{julka}), let $U \subseteq M$ be an open set.
Then the canonical map
$$ \Sect_c(U;E) \rightarrow \varinjlim_{K \subseteq U} \Sect_{K}(M; E)$$
is a homotopy equivalence.
\end{lemma}

\begin{proof}
It will suffice to show that if $A \subseteq B$ is an inclusion of finite simplicial sets and
we are given a commutative diagram
$$ \xymatrix{ A \ar[r] \ar[d] & \Sect_c(U;E) \ar[d] \\
B \ar[r]_-{f} \ar@{-->}[ur]^{f'} & \varinjlim_{K \subseteq U} \Sect_{K}(M;E), }$$
then, after modifying $f$ by a homotopy that is constant on $A$, there
exists a dotted arrow $f'$ as indicated in the diagram (automatically unique, since the right vertical map is a monomorphism). Since $B$ is finite, we may assume that $f$ factors through
$\Sect_{K}(M;E)$ for some compact subset $K \subseteq U$. Such a factorization determines a pair
$(F, h)$, where $F: |B| \times M \rightarrow E$ is a map of spaces over $M$
and $h: |B| \times (M-K) \times [0,1] \rightarrow E$ is a fiberwise homotopy of
$F | (|B| \times (M-K))$ to the composite map $|B| \times (M-K) \rightarrow M \stackrel{s}{\rightarrow} E$. 
Choose a continuous map $\lambda: M \rightarrow [0,1]$ which is supported in a compact subset $K'$ of $U$ and takes the value $1$ in a neighborhood of $K$. Let $F': |B| \times M \rightarrow E$ be the map
defined by the formula
$$ F'( b, x) = \begin{cases} F(b,x) & \text{ if } x \in K \\
h( b,x, 1-\lambda(x) ) & \text{ if } x \notin K. \end{cases}$$
Then $F'$ determines a map $B \rightarrow \Sect_{c}(U;E)$ such that the composite map
$B \rightarrow \Sect_{c}(U;E) \rightarrow \varinjlim_{K \subseteq U} \Sect_{K}(M;E)$ is
homotopic to $f$ relative to $A$, as desired.
\end{proof}

We now proceed with the proof of Theorem \ref{stager}. If $M$ is homeomorphic to Euclidean space
$\R^{k}$, then $\Disj{M}$ contains $M$ as a final object and Theorem \ref{stager} is obvious. Combining this observation with Remark \ref{conswe}, we obtain an immediate proof in the case $k=0$. 
If $k=1$, then we may assume (by virtue of Remark \ref{conswe}) that $M$ is homeomorphic to either
an open interval (in which case there is nothing to prove) or to the circle $S^1$. The latter case requires some argument:

\begin{proof}[Proof of Theorem \ref{stager} for $M = S^1$]
Choose a small open disk $U \subseteq S^1$ and a parametrization $\psi: \R \simeq U$, and
let $\chi: \DISJ{S^1}_{/ \psi} \rightarrow \SSet$ be the diagram determined by $E^!$. According to Theorem \ref{junl}, the functor $\chi$ is equivalent to a composition $\DISJ{S^1}_{/ \psi} \rightarrow \Nerve( \cDelta^{op}) \stackrel{B_{\bigdot}}{\rightarrow} \SSet$ for some simplicial object $B_{\bigdot}$ of $\SSet$, and the topological chiral homology $\int_{S^1} E^!$ can be identified with the geometric
realization $| B_{\bigdot} |$. We wish to show that the canonical map
$\theta: | B_{\bigdot} | \rightarrow \Sect_c(S^1; E)$ is an equivalence in $\SSet$. Since
$\SSet$ is an $\infty$-topos, it will suffice to verify the following pair of assertions:
\begin{itemize}
\item[$(a)$] The map $\theta_0: B_0 \rightarrow \Sect_c(S^1; E)$ is an effective epimorphism.
In other words, $\theta_0$ induces a surjection $\pi_0 B_0 = \pi_0 \Sect_c(U; E) \rightarrow
\pi_0 \Sect_c(S^1; E)$.
\item[$(b)$] The map $\theta$ exhibits $B_{\bigdot}$ as a \Cech nerve of $\theta_0$. That is,
for each $n \geq 0$, the canonical map
$$ B_{n} \rightarrow B_0 \times_{\Sect_c( S^1; E) } \cdots \times_{ \Sect_c(S^1;E)} B_0$$ is a homotopy equivalence (here the fiber products are taken in the $\infty$-category $\SSet$).
\end{itemize}

To prove $(a)$, let $s: S^1 \rightarrow E$ denote our given section of the Serre fibration $p: E \rightarrow S^1$, and let $f: S^1 \rightarrow E$ denote any other section of $p$. Choose a point $x \in U$.
Since $S^1 - \{x \}$ is contractible and the fibers of $p$ are connected, there exists a (fiberwise) homotopy $h: (S^1 - \{x\}) \times [0,1] \rightarrow E$ from $f | (S^1 - \{x\})$ to
$s | ( S^1 - \{ x \})$. Let $\lambda: S^1 \rightarrow [0,1]$ be a continuous function which vanishes
in a neighborhood of $x$, and takes the value $1$ outside a compact subset of $U$.
Let $h': S^1 \times [0,1] \rightarrow E$
be the map defined by
$$ h'( y, t) = \begin{cases} f(x) & \text{ if } y = x \\
h(y, t \lambda(y) ) & \text{ if } y \neq x. \end{cases}$$
Then $h'$ determines a homotopy from $f$ to another section $f' = h' | (S^1 \times \{1\})$, whose
support is a compact subset of $U$.

We now prove $(b)$. Choose a collection of open disks $U_1, \ldots, U_n \subseteq S^1$ which
are disjoint from one another and from $U$. Then the closed set
$S^1 - (U \cup U_1 \cup \ldots \cup U_n)$ is a disjoint union of connected components
$A_0, \ldots, A_n$. Unwinding the definitions, we are required to show that the simplicial set
$\Sect_c( U \cup U_1 \cup \ldots \cup U_n; E)$ is a homotopy product of
the simplicial sets $\Sect_c( S^1 - A_{i}; E)$ in the model category
$(\sSet)_{/ \Sect_c(S^1; E)}$. For each index $i$, let $\calU_i$ denote the collection of all
open subsets of $S^1$ that contain $A_{i}$, and let $\calU = \bigcap \calU_{i}$. It follows
from Lemma \ref{kora} that we have canonical homotopy equivalences
$$ \Sect_{c}( S^1 - A_i; E) \rightarrow \varinjlim_{ V \in \calU_i} \Sect_{S^1 - V}(S^1; E)$$
$$ \Sect_{c}( U \cup U_1 \cup \ldots \cup U_n; E) \rightarrow \varinjlim_{ V \in \calU} \Sect_{S^1 - V}( S^1; E).$$
Note that for each $V \in \calU_{i}$, the forgetful map $\Sect_{S^1-V}(S^1; E) \rightarrow
\Sect_c(S^1; E)$ is a Kan fibration. It follows that each
$\varinjlim_{V \in \calU_i} \Sect_{S^1 - V}( S^1; E)$ is a fibrant object of
$(\sSet)_{/ \Sect_c(S^1; E)}$, so the relevant homotopy product coincides with the actual
product $\prod_{0 \leq i \leq n} \varinjlim_{ V_i \in \calU_i} \Sect_{S^1 - V_i}(S^1; E)$
(formed in the category $(\sSet)_{/ \Sect_c(S^1;E)}$. 
Let $\calV$ denote the partially ordered set of sequences $(V_0, \ldots, V_n) \in \calU_0 \times \cdots \times \calU_n$ such that $V_{i} \cap V_j = \emptyset$ for $i \neq j$. We observe that the inclusion
$\calV^{op} \subseteq ( \calU_0 \times \cdots \times \calU_n)^{op}$ is cofinal, and the construction
$(V_0, \ldots, V_n) \mapsto \bigcup V_{i}$ is a cofinal map from $\calV^{op}$ to $\calU^{op}$.
Consequently, we obtain isomorphisms
$$ \varinjlim_{ V \in \calU} \Sect_{S^1 - V}(S^1; E) \simeq \varinjlim_{ (V_0, \ldots, V_n)} \Sect_{S^1 - \bigcup V_i}(S^1; E)$$
$$ \prod_{0 \leq i \leq n} \varinjlim_{ V_i \in \calU_i} \Sect_{S^1 - V_i}(S^1; E)
\simeq \varinjlim_{ (V_0, \ldots, V_n) \in \calV} \prod_{0 \leq i \leq n} \Sect_{ S^1 - V_i}(S^1;E);$$
here the product is taken in the category $(\sSet)_{/ \Sect_c(S^1; E)}$. 
To complete the proof, it suffices to show that for each $(V_0, V_1, \ldots, V_n) \in \calV$, the canonical map
$$ \theta: \Sect_{ S^1 - \bigcup V_i}( S^1;E) \rightarrow \prod_{0 \leq i \leq n} \Sect_{S^1 - V_i}(S^1;E)$$
is a homotopy equivalence.

We now complete the proof by observing that $\theta$ is an isomorphism
(since the open sets $V_i$ are assumed to be pairwise disjoint).
\end{proof}

Our proof of Theorem \ref{stager} in higher dimensions will use a rather different method. We first consider the following {\em linear} version of Theorem \ref{stager}, which is an easy consequence of the version of Verdier duality presented in \S \ref{kopo}.

\begin{proposition}\label{scone}
Let $M$ be a $k$-manifold, let $\calF \in \Shv(M; \Spectra)$ be a locally constant
$\Spectra$-valued sheaf on $M$, and let $\calF' \in \Shv(M; \SSet_{\ast})$ be the
sheaf of pointed spaces given by the formula $\calF'(U) = \Omega^{\infty}_{\ast} \calF(U)$.
Assume that for every open disk $U \subseteq M$, the spectrum $\calF(U)$ is $k$-connective.
Then $\calF'$ exhibits $\Gamma_c( M; \calF')$ as a colimit of the diagram
$\{ \Gamma_c(U; \calF') \}_{U \in \Disj{M} }$ in the $\infty$-category $\SSet_{\ast}$.
\end{proposition}

\begin{proof}
It follows from Corollary \ref{staffer} that $\calF$ exhibits
$\Gamma_c( M; \calF)$ as a colimit of the diagram
$$\{ \Gamma_c(U; \calF) \}_{U \in \Disj{M} }$$ in the $\infty$-category $\Spectra$ of spectra.
It will therefore suffice to show that the functor $\Omega^{\infty}_{\ast}$ preserves
the colimit of the diagram $\{ \Gamma_c(U; \calF) \}_{U \in \Disj{M}}$. 

Let us regard the $\infty$-category $\Spectra$ as endowed with its Cartesian symmetric monoidal structure, which (by virtue of Proposition \symmetricref{fingertone}) is also the coCartesian symmetric monoidal structure. The functor $U \mapsto \Gamma_c(U; \calF)$ determines a functor
$\Nerve( \Disk{M}) \rightarrow \Spectra$, which extends to a map of $\infty$-operads
$\Nerve( \Disk{M})^{\amalg} \rightarrow \Spectra^{\amalg}$ and therefore determines
an algebra $A \in \Alg_{ \Nerve( \Disk{M})}(\Spectra)$. Since $\calF$ is locally constant, the
algebra $A$ is locally constant and is therefore equivalent to a composition
$$ \Nerve( \Disk{M})^{\otimes} \rightarrow \OpE{M} \stackrel{B}{\rightarrow} \Spectra^{\times}.$$
Let $A': \Nerve( \Disk{M})^{\otimes} \rightarrow \Spectra$ and $B': \OpE{M} \rightarrow \Spectra$
be the associated monoid objects of $\Spectra$ (see Proposition \symmetricref{ungbatt}).
We wish to show that $\Omega^{\infty}_{\ast}$ preserves the colimit of the diagram
$\{ A'(U) \}_{U \in \Disj{M} }$. In view of Proposition \ref{koio}, it will suffice to prove
that $\Omega^{\infty}_{\ast}$ preserves the colimit of the diagram
$B' | \DISJ{M}$. For every open set $U = U_1 \cup \ldots \cup U_n$ of $\DISJ{M}$, the
spectrum $B'(U) \simeq \prod_{1 \leq i \leq n} B'(U_i) \simeq
\prod_{1 \leq i \leq n} \Omega^{k} \calF(U_i)$ is connective. Since
the $\infty$-category $\DISJ{M}$ is sifted (Proposition \ref{scun}), the desired result follows from
Corollary \ref{koso}.
\end{proof}

\begin{proof}[Proof of Theorem \ref{stager} for $k \geq 2$]
Replacing $E$ by $| \Sing(E) |$, we can assume without loss of generality
that $E$ is the geometric realization of a simplicial set $X$ equipped with a Kan fibration
$X \rightarrow \Sing(M)$. We wish to prove that the canonical map
$\int_{M} E^! \rightarrow \Sect_{c}( M; E)$ is a homotopy equivalence. 
For this, it suffices to show that $\tau_{\leq m} (\int_{M} E^{!}) \rightarrow \tau_{\leq m} \Sect_{c}(M;E)$
is a homotopy equivalence for every integer $m \geq 0$. Since the truncation functor
$\tau_{\leq m}: \SSet \rightarrow \tau_{\leq m} \SSet$ preserves small colimits and finite products,
Proposition \ref{skoo} allows us to identify the left hand side with the topological chiral homology
$\int_{M} ( \tau_{\leq m} E^!)$ in the $\infty$-category $\tau_{\leq m} \SSet$. 

Regard $X$ as an object of the $\infty$-topos $\calX = \SSet_{/ \Sing(M)}$, let $X'$ be an $(m+k)$-truncation of $X$, and let $E' = | X' |$. The map $X \rightarrow X'$ induces a map
$E^! \rightarrow {E'}^{!}$ which is an equivalence on $m$-truncations, and therefore
induces an equivalence $\tau_{ \leq m} ( \int_M E^!) \rightarrow \tau_{\leq m} ( \int_{M} {E'}^{!})$. 
This equivalence fits into a commutative diagram
$$ \xymatrix{ \tau_{\leq m} \int_{M} {E}^{!} \ar[r]^{\alpha} \ar[d] & \tau_{\leq m} \Sect_{c}(M; E) \ar[d]^{\beta} \\
\tau_{\leq m} \int_{M} {E'}^{!} \ar[r]^{\alpha'} &  \tau_{\leq m} \Sect_c(M; E'), }$$
where $\beta$ is also an equivalence (since $M$ has dimension $k$). Consequently, to prove
that $\alpha$ is an equivalence, it suffices to prove that $\alpha'$ is an equivalence. We may therefore replace $X$ by $X'$ and thereby reduce to the case where $X$ is an $n$-truncated object
of $\calX$ for some $n \gg 0$. 

The proof now proceeds by induction on $n$. If $n < k$, then $X$ is both $k$-connective and $(k-1)$-truncated, and is therefore equivalent to the final object of $\calX$. In this case, both
$\int_{M} E^{!}$ and $\Sect_{c}(M; E)$ are contractible and there is nothing to prove. Assume
therefore that $n \geq k \geq 2$. Let $A = \pi_n X$, regarded as an object of the
topos of discrete objects $\Disc{\calX_{/X} }$. Since $X$ is a $2$-connective object of
$\calX$, this topos is equivalent to the topos of discrete objects $\Disc{\calX}$ of local systems of
sets on the manifold $M$. We will abuse notation by identifying $A$ with its image under this equivalence; let $K(A,n+1)$ denote the associated Eilenberg-MacLane objects of $\calX$. Let $Y = \tau_{\leq n-1} X$, so that $X$ is an $n$-gerbe over $Y$ banded by $A$ and therefore fits into a pullback square
$$ \xymatrix{ X \ar[r] \ar[d] & {\bf 1} \ar[d] \\
Y \ar[r] & K(A,n+1) }$$
Let $E_0 = |Y|$ and $E_1 = | K(A,n+1)|$, so that we have a fiber sequence 
$E \rightarrow E_0 \rightarrow E_1$
of Serre fibrations over $M$. We then have a commutative diagram
$$ \xymatrix{ \int_{M} E^! \ar[r]^{\alpha} \ar[d] & \int_{M} E_0^{!} \ar[r] \ar[d]^{\alpha_0} & \int_{M} E_1^{!} \ar[d]^{\alpha_1} \\
\Sect_{c}(M; E) \ar[r] & \Sect_{c}( M; E_0) \ar[r] & \Sect_{c}(M; E_1) }$$
where $\alpha_0$ is a homotopy equivalence by the inductive hypothesis, and $\alpha_1$ is a homotopy equivalence by Proposition \ref{scone}. Consequently, to prove that $\alpha$ is a homotopy
equivalence, it suffices to prove that the upper line is a fiber sequence. The algebras
$E^{!}$, $E_0^{!}$, and $E_{1}^!$ determine functors $\chi, \chi_0, \chi_1: \DISJ{M} \rightarrow \SSet_{\ast}$, which fit into a pullback square
$$ \xymatrix{ \chi \ar[d] \ar[r] & \ast \ar[d] \\
\chi_0 \ar[r] & \chi_1. }$$
To complete the proof, it suffices to show that the induced square of colimits
$$ \xymatrix{ \varinjlim(\chi) \ar[r] \ar[d] & \ast \ar[d] \\
\varinjlim(\chi_0) \ar[r] & \varinjlim(\chi_1) }$$
is again a pullback diagram. Since $n \geq k$, the object $K(A,n+1)$ is $(k+1)$-connective, so that
$\chi_1$ takes values in connected spaces. The desired result now follows from Theorem \ref{kanter}, since $\DISJ{M}$ is sifted (Proposition \ref{scun}).
\end{proof}

\appendix

\section{Background on Topology}\label{secA}

In this appendix, we collect together some results in topology which are relevant (directly or indirectly) to the body of this paper. We begin in \S \ref{sec5sub1} by describing a ``higher'' version of the Seifert-van Kampen theorem, which permits us to reconstruct the weak homotopy type of a topological space $X$ from any covering (or hypercovering) of $X$. Our principal application is given in \S \ref{lubla}, where we show that if $X$ is a sufficiently nice topological space, then the $\infty$-category of locally constant ($\SSet$-valued) sheaves on $X$ is equivalent to the $\infty$-category of functors from the Kan complex $\Sing(X)$ into $\SSet$ (Theorem \ref{squ}). The proof relies on having developed a good theory of locally constant sheaves (which we describe in \S \ref{clome}) and on the homotopy invariance of this theory, which we prove in \S \ref{hominv}.

The theory of locally constant sheaves is really a special case of the more general theory of {\em constructible} sheaves on a stratified topological space $X$, which we review in \S \ref{coofer}. While locally constant sheaves on $X$ can be described as $\SSet$-valued functors on the Kan complex $\Sing(X)$, constructible sheaves on an $A$-stratified topological space can be described as $\SSet$-valued functors on an
$\infty$-category $\Sing^{A}(X) \subseteq \Sing(X)$, which we call the {\it $\infty$-category of exit paths} of
$X$. We will define this $\infty$-category in \S \ref{lubos} and establish its connection with constructible sheaves in \S \ref{cloop}. The proof relies on a generalization of the Seifert-van Kampen theorem to stratified spaces, which we prove in \S \ref{sec5sub8}, and on a general formalism for analyzing ``stratified'' $\infty$-categories, which we discuss in \S \ref{sec5sub9}. In \S \ref{sec5sub7}, we will give a detailed description of
the $\infty$-category of exit paths in the case of a simplicial complex (stratified by its simplices), in which case
the $\infty$-category $\Sing^{A}(X)$ is equivalent to the nerve of the partially ordered set $A$ of simplices of $X$ (Theorem \ref{ikos}).

The theory of factorizable (co)sheaves developed in \S \ref{sec3} relies heavily on understanding moduli spaces of embeddings between manifolds of the same dimension. For this reason, we collect together (with proofs) some basic facts about these embedding spaces in \S \ref{cowpy}.
Finally, in \S \ref{kopo}, we sketch a version of Verdier duality which can be applied to sheaves of spectra (or sheaves with values in any other stable $\infty$-category) on a locally compact topological space $X$.
This result has a simple consequence (Corollary \ref{staffer}) that plays an essential role in our discussion of nonabelian Poincare duality in \S \ref{nonpn}.

\subsection{The Seifert-van Kampen Theorem}\label{sec5sub1}

Let $X$ be a topological space covered by a pair of open sets $U$ and $V$, such that
$U$, $V$, and $U \cap V$ are path-connected. The Seifert-van Kampen theorem asserts that, for any choice of base point $x \in U \cap V$, the diagram of groups
$$ \xymatrix{ \pi_1( U \cap V,x) \ar[r] \ar[d] & \pi_1( U,x) \ar[d] \\
\pi_1(V,x) \ar[r] & \pi_1(X,x) }$$
is a pushout square. In this section, we will prove a generalization of the Seifert-van Kampen theorem, which
describes the entire weak homotopy type of $X$ in terms of any sufficiently nice covering of $X$ by open sets:

\begin{theorem}\label{vankamp}
Let $X$ be a topological space, let $\calU(X)$ denote the collection of all open subsets of $X$
(partially ordered by inclusion). Let $\calC$ be a small category and let
$\chi: \calC \rightarrow \calU(X)$ be a functor. For every $x \in X$, let
$\calC_{x}$ denote the full subcategory of $\calC$ spanned by those objects
$C \in \calC$ such that $x \in \chi(C)$. Assume that $\chi$ satisfies the following condition:
\begin{itemize}
\item[$(\ast)$] For every point $x$, the simplicial set $\Nerve( \calC_{x})$ is weakly contractible.
\end{itemize}
Then the canonical map $\varinjlim_{C \in \calC} \Sing( \chi(C) ) \rightarrow \Sing(X)$
exhibits the simplicial set $\Sing(X)$ as a homotopy colimit of the diagram
$\{ \Sing( \chi(C) ) \}_{C \in \calC}$.
\end{theorem}

The proof of Theorem \ref{vankamp} will occupy our attention throughout this section. The main step will be to establish the following somewhat weaker result:

\begin{proposition}\label{swugman}
Let $X$ be a topological space, let $\calU(X)$ be the partially ordered set of all
open subsets of $X$, and let $S \subseteq \calU(X)$ be a covering sieve on $X$.
Then the canonical map
$\varinjlim_{U \in S} \Sing(U) \rightarrow \Sing(X)$ exhibits the simplicial set $\Sing(X)$
as the homotopy colimit of the diagram of simplicial sets $\{ \Sing(U) \}_{U \in S}$.
\end{proposition}

Proposition \ref{swugman} is itself a consequence of the following result, which guarantees
that $\Sing(X)$ is weakly homotopy equivalent to the simplicial subset consisting of
``small'' simplices:

\begin{lemma}\label{twop}
Let $X$ be a topological space, and let $\{ U_{\alpha} \}$ be an open covering of $X$.
Let $\Sing'(X)$ be the simplicial subset of $\Sing(X)$ spanned by those $n$-simplices
$| \Delta^n | \rightarrow X$ which factor through some $U_{\alpha}$. Then
the inclusion $i: \Sing'(X) \subseteq \Sing(X)$ is a weak homotopy equivalence of simplicial sets.
\end{lemma}

The proof of Lemma \ref{twop} will require a few technical preliminaries.

\begin{lemma}
Let $X$ be a compact topological space and let $K$ be a simplicial set.
Then every continuous map $f: X \rightarrow |K|$ factors through
$|K_0|$, for some finite simplicial subset $K_0 \subseteq K$.
\end{lemma}

\begin{proof}
Let $K_0$ be the simplicial subset of $K$ spanned by those simplices
$\sigma$ such that the interior of $| \sigma |$ intersects $f(X)$. We claim
that $K_0$ is finite. Otherwise, we can choose an infinite sequence of points
$x_0, x_1, \ldots \in X$ such that each $f(x_i)$ belongs to the interior of
a {\em different} simplex of $|K|$. Let $U = |K| - \{ f(x_0), f(x_1), \ldots, \}$,
and for each $i \geq 0$ let $U_{i} = U \cup \{ f(x_i ) \}$. Then the
collection of open sets $\{ U_i \}$ forms an open cover of $K$, so that
$\{ f^{-1} U_i \}$ forms an open covering of $X$. This open covering does not
admit a finite subcovering, contradicting our assumption that $X$ is compact.
\end{proof}

\begin{lemma}\label{capst}
Let $i: K_0 \subseteq K$ be an inclusion of simplicial sets.
Suppose that the following condition is satisfied:
\begin{itemize}
\item[$(\ast)$] For every finite simplicial subset $L \subseteq K$, 
there exists a homotopy $h: |L| \times [0,1] \rightarrow |K|$
such that $h | (|L| \times \{0\})$ is the inclusion, 
$h| (|L| \times \{1\} ) \subseteq |K_0|$, and
$h| (|L_0| \times [0,1]) \subseteq |K_0|$, where
$L_0 = L \cap K_0$.
\end{itemize}
Then the inclusion $i$ is a weak homotopy equivalence.
\end{lemma}

\begin{proof}
We first show the following:
\begin{itemize}
\item[$(\ast')$] Let $X$ be a compact topological space, $X_0$ a closed subspace, and
$f: X \rightarrow |K|$ a continuous map such that $f( X_0) \subseteq |K_0|$. Then
there exists a homotopy $h: X \times [0,1] \rightarrow |K|$ such that
$h | (X \times \{0\} ) = f$, $h( X \times \{1\}) \subseteq |K_0|$, and
$h | (X_0 \times [0,1]) \subseteq |K_0|$.
\end{itemize}
To prove $(\ast')$, we note that since $X$ is compact, the map $f$
factors through $|L|$, where $L$ is some finite simplicial subset of $K$.
Then $f| X_0$ factors through $|L_0|$, where $L_0 = L \cap K_0$. We may therefore
replace $X$ and $X_0$ by $|L|$ and $|L_0|$, in which case $(\ast')$ is equivalent to our assumption $(\ast)$.

Applying $(\ast')$ in the case where $X$ is a point and $X_0$ is empty, we deduce that
the inclusion $i$ is surjective on connected components. It will therefore suffice to show that $i$ induces a bijection $\phi: \pi_n( |K_0|, v) \rightarrow \pi_{n}( |K|, v )$ for each $n \geq 0$ and each vertex $v$ of $K$. To prove that $\phi$ is surjective, consider a homotopy class $\eta \in \pi_{n}( |K|, v)$.
This homotopy class can be represented by a pointed map $f: (S^n, \ast) \rightarrow (|K|,v)$.
Applying $(\ast')$, we deduce that $f$ is homotopic to a another map
$g: S^n \rightarrow |K_0|$, via a homotopy which, when restricted to the base point
$\ast \in S^n$, determines a path $p$ from $v$ to another point $v' \in |K_0|$. 
Then $g$ determines an element $\eta' \in \pi_{n}( |K_0|,v')$. The image of
$\eta'$ under the transport isomorphism $p_{\ast}: \pi_{n}( |K_0|, v')
\simeq \pi_{n}( |K_0|, v)$ is a preimage of $\eta$ under $\phi$.

We now prove that $\phi$ is injective. Suppose we are given a continuous map
$f_0: S^{n} \rightarrow |K_0|$ which extends to a map $f: D^{n+1} \rightarrow |K|$; we wish to show that 
$f_0$ is nullhomotopic. Applying $(\ast')$, we deduce that $f_0$ is homotopic to a map
which extends over the disk $D^{n+1}$, and is therefore itself nullhomotopic.
\end{proof}

Before we can proceed with the proof of Lemma \ref{twop}, we need to recall some properties
of the {\it barycentric subdivision} construction in the setting of simplicial sets.

\begin{notation}\label{cuspa}
Let $[n]$ be an object of $\cDelta$. We let $P[n]$ denote the collection of all nonempty subsets of $[n]$, partially ordered by inclusion. We let $\overline{P}[n]$ denote the disjoint union
$P[n] \coprod [n]$. We regard $\overline{P}[n]$ as endowed with a partial ordering which
extends the partial orderings on $P[n]$ and $[n]$, where we let
$i \nleq \sigma$ for $i \in [n]$ and $\sigma \in P[n]$, while $\sigma \leq i$ if and only if
each element of $\sigma$ is $\leq i$.

The functors $[n] \mapsto \Nerve P[n]$ and $[n] \mapsto \Nerve \overline{P}[n]$
extend to colimit-preserving functors from the category of simplicial sets to itself.
We will denote these functors by $\sd$ and $\overline{\sd}$, respectively.

Let us identify the topological $n$-simplex $| \Delta^n |$ which the set of
all maps $t: [n] \rightarrow [0,1]$ such that $t(0) + \ldots + t(n) = 1$.
For each $n \geq 0$, there is a homeomorphism
$\eta_{n}: | \overline{P}[n] | \rightarrow | \Delta^n | \times [0,1]$ which 
is linear on each simplex, carries a vertex $i \in [n]$ to $(t_{i},0)$ where
$t_i$ is given by the formula $t_{i}(j) = \begin{cases} 1 & \text{ if } i= j \\
0 & \text{ if } i \neq j, \end{cases}$
and carries a vertex $\sigma \in P[n]$ to the pair
$(t_{\sigma}, 1)$, where $$t_{\sigma}(i) = \begin{cases} \frac{ 1 }{m} & \text{ if } i \in \sigma \\
0 & \text{ if } i \notin \sigma \end{cases}$$
where $m$ is the cardinality of $\sigma$. This construction is functorial
in $[n]$, and induces a homeomorphism
$| \overline{\sd} K| \rightarrow |K| \times [0,1]$ for every simplicial set $K$.
We observe that $\overline{\sd} K$ contains $K$ and $\sd K$ as simplicial
subsets, whose geometric realizations map homeomorphically to $|K| \times \{0\}$ and $|K| \times \{1\}$, respectively.
\end{notation}

\begin{proof}[Proof of Lemma \ref{twop}]
We will show that $i$ satisfies the criterion of Lemma \ref{capst}. Let
$L \subseteq \Sing(X)$ be a finite simplicial subset, and let $L_0 = L \cap \Sing'(X)$.
Fix $n \geq 0$, let $\overline{L}$ denote the iterated pushout
$$ \overline{\sd} \sd^{n-1} L \coprod_{ \sd^{n-1} L}
\overline{\sd} \sd^{n-2} L \coprod_{ \sd^{n-2} L} \ldots \coprod_{ \sd L} L,$$
and define $\overline{L}_0$ similarly.
Using the homeomorphisms $| \overline{\sd} K | \simeq |K| \times [0,1]$ of
Notation \ref{cuspa} repeatedly, we obtain a homeomorphism
$| \overline{L} | \simeq |L| \times [0,n]$ (which restricts to a homeomorphism
$| \overline{L}_0 | \simeq |L_0| \times [0,n]$). 

The inclusion map $L \subseteq \Sing(X)$ is adjoint to a continuous map of
topological spaces $f: |L| \rightarrow X$. Let $\overline{f}$ denote the composite map
$$ | \overline{L} | \simeq |L| \times [0,n] \rightarrow |L| \stackrel{f}{\rightarrow} X.$$
Then $\overline{f}$ determines a map of simplicial sets $\overline{L} \rightarrow \Sing(X)$;
we observe that this map carries $\overline{L}_0$ into $\Sing'(X)$. Passing to geometric realizations, 
we get a map
$h: | L | \times [0,n] \simeq | \overline{L} | \rightarrow | \Sing(X) |$, which is a homotopy
from the inclusion $|L| \subseteq | \Sing(X) |$ to the map
$g = h | ( |L| \times \{n\})$ (by construction, this homotopy carries $| L_0 | \times [0,n]$ into
$| \Sing'(X) |$). We note that $g$ is the geometric realization of the map
$\sd^{n} L \rightarrow \Sing'(X)$, which is adjoint to the composition
$| \sd^{n} L| \simeq |L| \stackrel{f}{\rightarrow} X$. To complete the proof, it suffices to observe
that for $n$ sufficiently large, each simplex of the $n$-fold barycentric subdivision
$| \sd^{n} L|$ will map into one of the open sets $U_{\alpha}$, so that
$g$ factors through $| \Sing'(X)|$ as required.
\end{proof}

Armed with Lemma \ref{twop}, it is easy to finish the proof of Proposition \ref{swugman}.

\begin{proof}[Proof of Proposition \ref{swugman}]
Choose a collection of open sets $\{ U_{\alpha} \}_{ \alpha \in A}$ which generates the sieve
$S$. Let $P(A)$ denote the collection of all nonempty subsets of $A$, partially ordered
by reverse inclusion. Let $P_0(A)$ be the subset consisting of nonempty {\em finite} subsets of $A$.
For each $A_0 \in P(A)$, let $U_{A_0} = \bigcap_{ \alpha \in A_0} U_{\alpha}$ (if $A_0$ is finite,
this is an open subset of $X$, though in general it need not be). The construction
$A_0 \mapsto U_{A_0}$ determines a map of partially ordered sets
$P_0(A) \rightarrow S$. Using Theorem \toposref{hollowtt}, we deduce that
the map $\Nerve( P_0(A) ) \rightarrow \Nerve(S)$ is cofinal, so that (by virtue of Theorem \toposref{colimcomparee}) it will suffice to show that $\Sing(X)$ is a homotopy colimit of the diagram
$\{ \Sing(U_{A_0}) \}_{A_0 \in P_0(A) }$. A similar argument shows that the inclusion
$\Nerve( P_0(A) ) \subseteq \Nerve( P(A) )$ is cofinal, so we are reduced to showing that
$\Sing(X)$ is a homotopy colimit of the diagram $\psi = \{ \Sing(U_{A_0}) \}_{A_0 \in P(A) }$. 
The actual colimit of the diagram $\psi$ is the simplicial set
$\Sing'(X)$ which is weakly equivalent to $\Sing(X)$ by Lemma \ref{twop}. It will therefore suffice
to show that the diagram $\psi$ is projectively cofibrant. To prove this, we will show more generally
that for any pair of simplicial subsets $K_0 \subseteq K \subseteq \Sing(X)$, the induced map
$$ \phi: \{ \Sing( U_{A_0}) \cap K_0 \}_{A_0 \in P(A)} \hookrightarrow
\{ \Sing( U_{A_0}) \cap K \}_{ A_0 \in P(A) }$$ is a projective cofibration
of diagrams (taking $K_0 = \emptyset$ and $K = \Sing(X)$ will then yield the desired result).
Working simplex by simplex, we may assume that $K$ is obtained from $K_0$ by adjoining
a single nondegenerate simplex $\sigma: | \Delta^n | \rightarrow X$ whose boundary already belongs to $K_0$. Let $A' = \{ \alpha \in A: \sigma( | \Delta^n| ) \subseteq U_{\alpha} \}$. If
$A'$ is empty, then $\phi$ is an isomorphism. Otherwise, $\phi$ is a pushout of the projective
cofibration $F_0 \hookrightarrow F$, where 
$$F_0(A_0) = \begin{cases} \bd \Delta^n & \text{ if } A_0 \subseteq A' \\
\emptyset & \text{ otherwise }\end{cases} \quad \quad F(A_0) = \begin{cases} \Delta^n & \text{ if } A_0 \subseteq A' \\
\emptyset & \text{ otherwise. }\end{cases}$$
\end{proof}

\begin{variant}\label{varswug}
If $X$ is a paracompact topological space, we can replace $\calU(X)$ with the collection of
all open $F_{\sigma}$ subsets of $X$ in the statement of Proposition \ref{swugman}; the proof remains the same.
\end{variant}

\begin{remark}\label{kwuppus}
Let $X$ be a topological space, and let $\calU(X)$ denote the partially ordered set of all open subsets
of $X$. The construction $U \mapsto \Sing(U)$ determines a functor between
$\infty$-categories $\Nerve( \calU(X) ) \rightarrow \SSet$. Theorem \toposref{charpresheaf} implies that this functor is equivalent to a composition
$$ \Nerve( \calU(X) ) \stackrel{j}{\rightarrow} \calP( \calU(X) ) \stackrel{F}{\rightarrow} \SSet,$$
where $j$ denotes the Yoneda embedding and the functor $F$ preserves small colimits (moreover, the functor $F$ is determined uniquely up to equivalence). Proposition \ref{swugman} implies that $F$
is equivalent to the composition
$$ \calP( \calU(X) ) \stackrel{L}{\rightarrow} \Shv(X) \stackrel{ F}{\rightarrow} \SSet,$$
where $L$ denotes a left adjoint to the inclusion $\Shv(X) \subseteq \calP(\calU(X))$
and we identify $F$ with its restriction to $\Shv(X)$. In particular, the functor
$F: \Shv(X) \rightarrow \SSet$ preserves small colimits.
\end{remark}

We now explain how to deduce Theorem \ref{vankamp} from Proposition \ref{swugman}. The main technical obstacle is that the $\infty$-topos $\Shv(X)$ need not be hypercomplete. We will address this problem by showing that the functor $F$ of Remark \ref{kwuppus} factors through the hypercompletion of $\Shv(X)$: in other words, that $F$ carries $\infty$-connected morphisms in $\Shv(X)$ to equivalences in $\SSet$ (Lemma \ref{swugman2}). We first note that $\infty$-connectedness is a condition which can be tested ``stalkwise'':

\begin{lemma}\label{copus}
Let $X$ be a topological space, and let $\alpha: \calF \rightarrow \calF'$ be a morphism
in the $\infty$-category $\Shv(X)$. For each point $x \in X$, let $x^{\ast}:
\Shv(X) \rightarrow \Shv( \{x\} ) \simeq \SSet$ denote the pullback functor.
The following conditions are equivalent:
\begin{itemize}
\item[$(1)$] The morphism $\alpha$ is $\infty$-connective.
\item[$(2)$] For each $x \in X$, the morphism $x^{\ast}( \alpha)$ is an equivalence in $\SSet$.
\end{itemize}
\end{lemma}

\begin{proof}
The implication $(1) \Rightarrow (2)$ is obvious, since the pullback functors
$x^{\ast}$ preserve $\infty$-connectivity and the $\infty$-topos $\SSet$ is hypercomplete.
Conversely, suppose that $(2)$ is satisfied. We will prove by induction on $n$
that the morphism $\alpha$ is $n$-connective. Assume that $n > 0$. By virtue of Proposition \toposref{trowler}, it will suffice to show that the diagonal map $\calF \times_{ \calF'} \calF$ is $(n-1)$-connective, which follows from the inductive hypothesis. We may therefore reduce to the case
$n = 0$: that is, we must show that $\alpha$ is an effective epimorphism.
According to Proposition \toposref{pi00detects}, this is equivalent to the requirement that
the induced map $\alpha': \tau_{\leq 0} \calF \rightarrow \tau_{\leq 0} \calF'$ is
an effective epimorphism. We may therefore replace $\alpha$ by $\alpha'$
and thereby reduce to the case where $\calF, \calF' \in \Shv_{\Set}(X)$ are sheaves
of sets on $X$, in which case the result is obvious.
\end{proof}

\begin{lemma}\label{swugman2}
Let $X$ be a topological space, and let $F: \Shv(X) \rightarrow \SSet$ be as in Remark \ref{kwuppus}. Then $F$ carries $\infty$-connective morphisms of $\Shv(X)$ to equivalences in $\SSet$.
\end{lemma}

\begin{proof}
Let $\alpha$ be an $\infty$-connectivemorphism in $\Shv(X)$. We will show that
$F(\alpha)$ is an $\infty$-connective morphism in $\SSet$, hence an equivalence (since the
$\infty$-topos $\SSet$ is hypercomplete). For this, it suffices to show that for each $n \geq 0$, 
the composite functor
$$ \Shv(X) \stackrel{F}{\rightarrow} \SSet \stackrel{ \tau^{\SSet}_{\leq n}}{\rightarrow} \tau_{\leq n} \SSet$$
carries $\alpha$ to an equivalence. Since $\tau_{\leq n} \SSet$ is an $n$-category, the functor
$\tau^{\SSet}_{\leq n} \circ F$ is equivalent to a composition
$$ \Shv(X) \stackrel{ \tau_{\leq n}^{\Shv(X)} }{\rightarrow} \tau_{\leq n} \Shv(X) \stackrel{ F_{n}}{\rightarrow} \tau_{\leq n} \SSet.$$
We now observe that $\tau_{\leq n}^{\Shv(X)}(\alpha)$ is an equivalence, since $\alpha$
is assumed to be $\infty$-connective. 
\end{proof}

We now have the tools in place to complete the proof of our main result.

\begin{proof}[Proof of Theorem \ref{vankamp}]
Passing to nerves, we obtain a diagram of $\infty$-categories
$\overline{p}: \Nerve(\calC)^{\triangleright} \rightarrow \SSet$. In view of
Theorem \toposref{colimcomparee}, it will suffice to show that $\overline{p}$ is a colimit
diagram. Note that $\overline{p}$ is equivalent to the composition
$$ \Nerve(\calC)^{\triangleright} \stackrel{ \overline{\chi}}{\rightarrow}
\Nerve( \calU(X) ) \stackrel{j}{\rightarrow} \Shv(X)^{\hyp} \stackrel{ F}{\rightarrow} \SSet,$$
where $\Shv(X)^{\hyp}$ denotes the full subcategory of $\calP( \calU(X) )$
spanned by the hypercomplete sheaves on $X$, $j$ denotes the Yoneda embedding,
and $F$ is defined as in Remark \ref{kwuppus}. Using Proposition \ref{swugman} and
Lemma \ref{swugman2}, we deduce that $F$ preserves small colimits. It therefore suffices to show
that $j \circ \overline{\chi}$ is a colimit diagram. Since $\Shv(X)^{\hyp}$ is hypercomplete, it suffices to show that the composition $f^{\ast} \circ j \circ \overline{\chi}$ is a colimit diagram, where
$f: \{ x\} \hookrightarrow X$ is the inclusion of any point into $X$. This follows immediately from
assumption $(\ast)$.
\end{proof}

\subsection{Locally Constant Sheaves}\label{clome}

Let $X$ be a topological space. A sheaf of sets $\calF$ on $X$ is said to be {\it constant} if
there exists a set $A$ and a map $\eta: A \rightarrow \calF(X)$ such that, for every point
$x \in X$, the composite map $A \rightarrow \calF(X) \rightarrow \calF_{x}$ is a bijection from
$A$ to the stalk $\calF_{x}$ of $\calF$ at $x$. More generally, we say that a sheaf of sets $\calF$ is {\it locally constant} if every point $x \in X$ has an open neighborhood $U$ such that the restriction $\calF|U$ is a constant sheaf on $U$. The category of locally constant sheaves of sets on $X$ is equivalent to the category of covering spaces of $X$. If $X$ is path connected and semi-locally simply connected, then the theory of covering spaces guarantees that this category is equivalent to the category of sets with an action of the fundamental group $\pi_1(X,x)$ (where $x$ is an arbitrarily chosen point of $X$). 

Our goal in this section is to obtain an $\infty$-categorical analogue of the above picture. More precisely, we will replace the topological space $X$ by an $\infty$-topos $\calX$. Our goal is to introduce a full subcategory of $\calX$ consisting of ``locally constant'' objects (see Definition \ref{placer}). We will further show that if
$\calX$ is sufficiently well-behaved, then this full subcategory is itself an $\infty$-topos: more precisely, it
is equivalent to an $\infty$-category of the form $\SSet_{/ K}$, for some Kan complex $K$. In \S \ref{lubla}, we will show that if $\calX$ is the $\infty$-category $\Shv(X)$ of sheaves on a well-behaved topological space $X$,
then we can take $K$ to be the Kan complex $\Sing(X)$. 

The first step is to formulate a condition on an $\infty$-topos which is a counterpart to the hypothesis of semi-local simple connectivity in the usual theory of covering spaces.

\begin{definition}
Let $\calX$ be an $\infty$-topos, let $\pi_{\ast}: \calX \rightarrow \SSet$ be a functor corepresented by the final object of $\calX$, and let $\pi^{\ast}$ be a right adjoint to $\pi_{\ast}$. We will
say that $\calX$ has {\it constant shape} if the composition $\pi_{\ast} \pi^{\ast}: \SSet \rightarrow \SSet$
is corepresentable.
\end{definition}

\begin{remark}
Recall that the {\it shape} of an $\infty$-topos $\calX$ is the functor
$\pi_{\ast} \pi^{\ast}: \SSet \rightarrow \SSet$, which can be regarded as
a pro-object of the $\infty$-category $\SSet$ (see \S \toposref{shapesec}). The $\infty$-topos $\calX$ has constant shape if this pro-object can be taken to be constant.
\end{remark}

\begin{remark}\label{saplen}
According to Proposition \toposref{representableprime}, an $\infty$-topos
$\calX$ has constant shape if and only if the functor $\pi_{\ast} \pi^{\ast}$ preserves small limits.
\end{remark}

\begin{remark}\label{siba}
Let $X$ be a paracompact topological space, and let $\pi_{\ast}: \Shv(X) \rightarrow \Shv(\ast) \simeq \SSet$ be the global sections functor. It follows from the results of \S \toposref{paracompactness} that we can identify the composition $\pi_{\ast} \pi^{\ast}$ with the functor
$K \mapsto \bHom_{ \Top}( X, |K| )$. Consequently, the $\infty$-topos $\Shv(X)$ has constant shape if and only if there exists a simplicial set $K_0$ and a continuous map $f: X \rightarrow |K_0|$ such that,
for every Kan complex $K$, composition with $f$ induces a homotopy equivalence
$\bHom_{ \sSet}( K_0, K) \simeq \bHom_{ \Top}( |K_0|, |K|) \rightarrow \bHom_{ \Top}( X, |K|)$. 
This is guaranteed, for example, if $f$ is a homotopy equivalence: in other words, if
$X$ is a paracompact topological space with the homotopy type of a CW complex, then
$X$ has constant shape.
\end{remark}

\begin{definition}\label{soko}
Let $\calX$ be an $\infty$-topos. We will say that an object $U \in \calX$ has {\it constant shape}
if the $\infty$-topos $\calX_{/U}$ has constant shape. We will say that $\calX$ is
{\it locally of constant shape} if every object $U \in \calX$ has constant shape.
\end{definition}

The following result guarantees that Definition \ref{soko} is reasonable:

\begin{proposition}\label{caspa}
Let $\calX$ be an $\infty$-topos, and let $\calX'$ be the full subcategory of $\calX$ spanned by those objects which have constant shape. Then $\calX'$ is stable under small colimits in $\calX$.
\end{proposition}

\begin{proof}
For each $U \in \calX$, let $\chi_{U}: \calX \rightarrow \SSet$ be the functor corepresented by $U$,
and let $\pi^{\ast}: \SSet \rightarrow \calX$ be a geometric morphism. Then $U$ has constant shape
if and only if the functor $\chi_{U} \circ \pi^{\ast}$ is corepresentable: in other words, if and only
if $\chi_{U} \circ \pi^{\ast}$ preserves small limits (Remark \ref{saplen}). Suppose that
$U$ is the colimit of a diagram $\{ U_{\alpha} \}$. Then $\chi_{U}$ is the limit of the induced diagram
of functors $\{ \chi_{U_{\alpha}} \}$ (Proposition \toposref{yonedaprop}), so that
$\chi_{U} \circ \pi^{\ast}$ is a limit of the diagram of functors $\{ \chi_{U_{\alpha}} \circ \pi^{\ast} \}$.
If each $U_{\alpha}$ has constant shape, then each of the functors $\chi_{U_{\alpha}} \circ \pi^{\ast}$ preserves small limits, so that $\chi_{U} \circ \pi^{\ast}$ preserves small limits
(Lemma \toposref{limitscommute}). 
\end{proof}

\begin{corollary}
Let $\calX$ be an $\infty$-topos. Suppose that there exists a collection of objects
$U_{\alpha} \in \calX$ such that the projection $U = \coprod_{\alpha} U_{\alpha} \rightarrow {\bf 1}$ is an effective epimorphism, where ${\bf 1}$ denotes the final object of $\calX$. If
each of the $\infty$-topoi $\calX_{ / U_{\alpha}}$ is locally of constant shape, then
$\calX$ is locally of constant shape.
\end{corollary}

\begin{proof}
Let $V \in \calX$; we wish to show that $V$ has constant shape.
Let $V_{0} = U \times V$, and let $V_{\bigdot}$ be the \Cech nerve of the
effective epimorphism $V_0 \rightarrow V$. Since $\calX$ is an $\infty$-topos,
$V$ is equivalent to the geometric realization of the simplicial object $V_{\bigdot}$.
In view of Proposition \ref{caspa}, it will suffice to show that each $V_{n}$ has constant shape.
We note that $V_{n}$ is a coproduct of objects of the form $U_{\alpha_0} \times \ldots \times U_{\alpha_{n}} \times V$. Then $\calX_{/V_{n}}$ admits an \etale geometric morphism to
the $\infty$-topos $\calX_{ / U_{\alpha_0}}$, which is locally of constant shape by assumption. It follows that $\calX_{/V_n}$ is of constant shape.
\end{proof}

\begin{proposition}\label{postkus}
Let $\calX$ be an $\infty$-topos and let $\pi^{\ast}: \SSet \rightarrow \calX$ be a geometric morphism.
The following conditions are equivalent:
\begin{itemize}
\item[$(1)$] The $\infty$-topos $\calX$ is locally of constant shape.
\item[$(2)$] The functor $\pi^{\ast}$ admits a left adjoint $\pi_{!}$.
\end{itemize}
\end{proposition}

\begin{proof}
According to Corollary \toposref{adjointfunctor}, condition $(2)$ is equivalent to the requirement
that $\pi^{\ast}$ preserves small limits. In view of Proposition \toposref{yonedaprop}, this is
equivalent to the assertion that for each $U \in \calX$, the composition
$\chi_{U} \circ \pi^{\ast}: \SSet \rightarrow \SSet$ preserves limits, where
$\chi_{U}: \calX \rightarrow \SSet$ is the functor corepresented by $U$.
\end{proof}

Let $\calX$ be an $\infty$-topos which is locally of constant shape, and let
$\pi_{!}$ and $\pi^{\ast}$ be the adjoint functors appearing in Proposition \ref{postkus}. 
Let $X \rightarrow Y$ be a morphism in $\SSet$ and let
$Z \rightarrow \pi^{\ast} Y$ be a morphism in $\calX$. Then we have a commutative
diagram
$$ \xymatrix{ \pi_{!}( \pi^{\ast} X \times_{ \pi^{\ast} Y} Z) \ar[r] \ar[d] & \pi_! Z \ar[d] \\
\pi_{!} \pi^{\ast} X \ar[r] \ar[d] & \pi_! \pi^{\ast} Y \ar[d] \\
X \ar[r] & Y, }$$
and the outer square determines a canonical map $\pi_{!}( \pi^{\ast} X \times_{ \pi^{\ast} Y} Z)
\rightarrow X \times_{Y} \pi_{!} Z$.

\begin{proposition}\label{cenc}
Let $\calX$ be an $\infty$-topos which is locally of constant shape, let
$\pi^{\ast}: \SSet \rightarrow \calX$ be a geometric morphism and
$\pi_{!}$ a left adjoint to $\pi^{\ast}$ (so that $\calX$ is locally of constant shape).
For every morphism $\alpha: X \rightarrow Y$ in $\SSet$ and every morphism
$\beta: Z \rightarrow \pi^{\ast} Y$ in $\calX$, the associated push-pull morphism
$$\pi_{!}( \pi^{\ast} X \times_{ \pi^{\ast} Y} Z)
\rightarrow X \times_{Y} \pi_{!} Z$$
is an equivalence.
\end{proposition}

\begin{proof}
Let us first regard the morphism $\alpha$ as fixed, and consider the full subcategory
$\calY \subseteq \calX_{/ \pi^{\ast} Y}$ spanned by those objects $Z$ for which the conclusion holds.
Since both $\pi_{!}( \pi^{\ast} X \times_{ \pi^{\ast} Y} Z)$ and $X \times_{Y} \pi_{!} Z$
are colimit-preserving functors of $Z$, the full subcategory $\calY$ is stable under colimits
in $\calX_{/ \pi^{\ast} Y}$. Regard $Y$ as a Kan complex, and let $\calC$ be the category
of simplices of $Y$, so that we can identify $Y$ with the colimit $\varinjlim_{C \in \calC}(\Delta^0)$ of
the constant diagram $\calC \rightarrow \SSet$ taking the value $\Delta^0$. For
every $Z \in \calX_{ / \pi^{\ast} Y}$, we have a canonical equivalence
$Z \simeq \varinjlim_{C \in \calC}(Z \times_{ \pi^{\ast} Y} \pi^{\ast} \Delta^0)$. We may
therefore replace $Z$ by the fiber product $Z \times_{\pi^{\ast} Y} \pi^{\ast} \Delta^0$, and thereby reduce to the case where $\beta$ factors through the map $\pi^{\ast} \Delta^0 \rightarrow \pi^{\ast} Y$
determined by a point of $Y$. Replacing $Y$ by $\Delta^0$ and $X$ by
$X \times_{Y} \Delta^0$, we can reduce to the case where $Y = \Delta^0$. In this case,
we must show that the canonical map $\pi_{!}( \pi^{\ast} X \times Z) \rightarrow X \times \pi_{!} Z$
is an equivalence. Let us now regard $Z$ as fixed and consider the full subcategory
$\calZ \subseteq \SSet$ spanned by those objects for which the conclusion holds. Since the
functors $\pi_{!}( \pi^{\ast} X \times Z)$ and $X \times \pi_{!} Z$ both preserve colimits
in $X$, the full subcategory $\calZ \subseteq \SSet$ is stable under small colimits.
It will therefore suffice to show that $\Delta^0 \in \SSet$, which is obvious.
\end{proof}

Let $\calX$ be an $\infty$-topos which is locally of constant shape. Let
$\pi_{!}$ and $\pi^{\ast}$ denote the adjoint functors appearing in Proposition \ref{postkus}. Let ${\bf 1}$ be a final object of $\calX$. We have a canonical functor
$$ \calX \simeq \calX_{ / {\bf 1}} \stackrel{ \pi_{!}}{\rightarrow} \SSet_{/ \pi_{!} {\bf 1}},$$
which we will denote by $\psi_{!}$. The functor $\psi_{!}$ admits a right adjoint $\psi^{\ast}$, which can
be described informally by the formula $\psi^{\ast} X = \pi^{\ast}X \times_{ \pi^{\ast} \pi_{!} {\bf 1}} {\bf 1}$
(Proposition \toposref{curpse}). We observe that $\psi^{\ast}$ preserves small colimits, and is
therefore a geometric morphism of $\infty$-topoi.

\begin{remark}
The object $\pi_{!} {\bf 1} \in \SSet$ can
be identified with the shape of the $\infty$-topos $\calX$.
\end{remark}

\begin{proposition}\label{standy}
Let $\calX$ be an $\infty$-topos which is locally of constant shape, and let
$\psi^{\ast}: \SSet_{ / \pi_{!} {\bf 1}} \rightarrow \calX$ be defined as above.
Then $\psi^{\ast}$ is fully faithful.
\end{proposition}

\begin{proof}
Fix an object $X \rightarrow \pi_{!} {\bf 1}$ in $\SSet_{/ \pi_{!} {\bf 1}}$; we wish to show that
the counit map $v: \psi_{!} \psi^{\ast} X \rightarrow X$ is an equivalence.
Unwinding the definitions, we see that $v$ can be identified with the push-pull transformation
$$ \pi_{!}( {\bf 1} \times_{ \pi^{\ast} \pi_{!} {\bf 1}} \pi^{\ast} X) \rightarrow
\pi_{!} {\bf 1} \times_{ \pi_{!} {\bf 1}} X \simeq X,$$
which is an equivalence by virtue of Proposition \ref{cenc}.
\end{proof}

We now describe the essential image of the fully faithful embedding $\psi^{\ast}$.

\begin{definition}\label{placer}
Let $\calX$ be an $\infty$-topos, and let $\calF$ be an object of $\calX$. We will say
that $\calF$ is {\it constant} if it lies in the essential image of a geometric morphism
$\pi^{\ast}: \SSet \rightarrow \calX$ (the geometric morphism $\pi^{\ast}$ is unique up to equivalence, by virtue of Proposition \toposref{spacefinall}). We will say that $\calF$ is {\it locally constant} if there
exists a small collection of objects $\{ U_{\alpha} \in \calX \}_{ \alpha \in S}$ such that the following
conditions are satisfied:
\begin{itemize}
\item[$(i)$] The objects $U_{\alpha}$ cover $\calX$: that is, there is an effective epimorphism
$\coprod U_{\alpha} \rightarrow {\bf 1}$, where ${\bf 1}$ denotes the final object of $\calX$.
\item[$(ii)$] For each $\alpha \in S$, the product $\calF \times U_{\alpha}$ is a
constant object of the $\infty$-topos $\calX_{ / U_{\alpha} }$.
\end{itemize}
\end{definition}

\begin{remark}
Let $f^{\ast}: \calX \rightarrow \calY$ be a geometric morphism of $\infty$-topoi. Then
$f^{\ast}$ carries constant objects of $\calX$ to constant objects of $\calY$ and locally constant objects of $\calX$ to locally constant objects of $\calY$.
\end{remark}

\begin{remark}\label{slope}
Let $\calF$ be a locally constant object of $\Shv(X)$, where $X$ is a topological space.
Then there exists an open covering $\{ U_{\alpha} \subseteq X \}$ such that each
$\calF | U_{\alpha}$ is constant. Moreover, if $X$ is paracompact, we can assume that
each $U_{\alpha}$ is an open $F_{\sigma}$ set.
\end{remark}

We now come to the main result of this section, which provides an $\infty$-categorical version of the classical theory of covering spaces.

\begin{theorem}\label{san}
Let $\calX$ be an $\infty$-topos which is locally of constant shape, and let
$\psi^{\ast}: \SSet_{ / \pi_{!} { \bf 1}} \rightarrow \calX$ be the functor of Proposition \ref{standy}. Then $\psi^{\ast}$ is a fully faithful embedding, whose essential image is the full subcategory of $\calX$ spanned by the locally constant objects.
\end{theorem}

\begin{proof}
Suppose first that $X \rightarrow \pi_{!} {\bf 1}$ is an object of $\SSet_{ / \pi_{!} {\bf 1}}$;
we will prove that $\psi^{\ast}(X)$ is locally constant. Choose an effective epimorphism
$\coprod_{ \alpha \in A} K_{\alpha} \rightarrow \pi_{!} {\bf 1}$ in $\SSet$, where
each $K_{\alpha}$ is contractible. Then we obtain
an effective epimorphism $\coprod_{\alpha \in A} \psi^{\ast} K_{\alpha} \rightarrow {\bf 1}$; it
will therefore suffice to show that each $\psi^{\ast} X \times \psi^{\ast} K_{\alpha}$ is
a constant object of $\calX_{/ \psi^{\ast} K_{\alpha}}$. The composite functor
$$ \SSet_{ / \pi_{!} {\bf 1} } \stackrel{ \psi^{\ast}}{\rightarrow} 
\calX \stackrel{ \times \psi^{\ast} K_{\alpha}}{\rightarrow} \calX_{ / \psi^{\ast} K_{\alpha}}$$
is equivalent to a composition of geometric morphisms
$$ \SSet_{ / \pi_{!} {\bf 1}} \rightarrow \SSet_{ / K_{\alpha}} \simeq \SSet \rightarrow
\calX_{ / \psi^{\ast} K_{\alpha}}$$
and so its essential image consists of constant objects.

For the converse, suppose that $\calF \in \calX$ is a locally constant object; we wish to show
that $\calF$ belongs to the essential image of $\psi^{\ast}$. 
Since $\calF$ is locally constant, there exists a diagram
$\{ U_{\alpha} \}$ in $\calX$ having colimit ${\bf 1}$, such that each
product $U_{\alpha} \times \calF$ is a constant object of $\calX_{/ U_{\alpha}}$.
We observe that $\SSet_{ / \pi_{!} {\bf 1}}$ can be identified with the limit
of the diagram of $\infty$-categories $\{ \SSet_{ / \pi_{!} U_{\alpha} } \}$, and that
$\calX$ can be identified with the limit of the diagram
of $\infty$-categories $\{ \calX_{ / U_{\alpha} } \}$ (Theorem \toposref{charleschar}).
Moreover, the fully faithful embedding $\psi^{\ast}$ is the limit of fully faithful embeddings
$\psi^{\ast}_{\alpha}: \SSet_{ / \pi_{!} U_{\alpha}} \rightarrow \calX_{ / U_{\alpha}}$. 
Consequently, $\calF$ belongs to the essential image of $\psi^{\ast}$ if and only if
each product $\calF \times U_{\alpha}$ belongs to the essential image of
$\psi^{\ast}_{\alpha}$. We may therefore replace $\calX$ by $\calX_{/U_{\alpha}}$ and
thereby reduce to the case where $\calF$ is constant. In this case, $\calF$ belongs to the essential
image of {\em any} geometric morphism $\phi^{\ast}: \calY \rightarrow \calX$, since we have a homotopy
commutative diagram of geometric morphisms
$$ \xymatrix{ & \calY \ar[dr]^{ \phi^{\ast} } &  \\
\SSet \ar[ur] \ar[rr]^{\pi^{\ast}} & & \calX. }$$
\end{proof}

\begin{corollary}\label{twux2}
Let $\calX$ be an $\infty$-topos which is locally of constant shape. Then the collection of
locally constant objects of $\calX$ is stable under small colimits.
\end{corollary}

\begin{corollary}\label{tablor2}
Let $\calX$ be an $\infty$-topos which is locally of constant shape. Then
for every locally constant object $X \in \calX$, the canonical map
$X \rightarrow \varprojlim \tau_{\leq n} X$ is an equivalence; in particular,
$X$ is hypercomplete.
\end{corollary}

\begin{proof}
Let $\pi_{!}: \calX \rightarrow \SSet$ and $\psi^{\ast}: \SSet_{ / \pi_! {\bf 1}} \rightarrow \calX$
be as in Proposition \ref{standy}. According to Theorem \ref{san}, we can write
$X = \psi^{\ast} X_0$ for some $X_0 \in \SSet_{ / \pi_{!} {\bf 1}}$. Since
$\psi^{\ast}$ commutes with truncations and preserves limits (being a right adjoint),
we can replace $\calX$ by $\SSet_{/ \pi_{!} {\bf 1}}$. Since the result is local on 
$\calX$, we can reduce further to the case where $\calX = \SSet$, in which case there
is nothing to prove.
\end{proof}

\subsection{Homotopy Invariance}\label{hominv}

Let $X$ be a topological space, and let $\calF$ be a locally constant sheaf of sets on $X$.
If $p: [0,1] \rightarrow X$ is a continuous path from $x = p(0)$ to $y = p(1)$, then $p$ induces
a bijection between the stalks $\calF_{x}$ and $\calF_{y}$ of the sheaf $\calF$, given by transport along $p$. More generally, if $h: Y \times [0,1] \rightarrow X$ is any homotopy from a continuous map
$h_0: Y \rightarrow X$ to a continuous map $h_1: Y \rightarrow X$, then $h$ induces an isomorphism
of sheaves $h_0^{\ast} \calF \simeq h_1^{\ast} \calF$. Our goal in this section is to generalize these statements to the case where $\calF$ is a sheaf of spaces.

Our first step is to study locally constant sheaves on the unit interval $[0,1]$. These are characterized by the following result:

\begin{proposition}\label{sooly}
Let $X$ be the unit interval $[0,1]$, and let $\calF \in \Shv(X)$. Let
$\pi_{\ast}: \Shv(X) \rightarrow \Shv( \ast) = \SSet$ be the global sections functor,
and let $\pi^{\ast}$ be a left adjoint to $\pi_{\ast}$. The following conditions are equivalent:
\begin{itemize}
\item[$(i)$] The sheaf $\calF$ is locally constant.
\item[$(ii)$] The sheaf $\calF$ is constant.
\item[$(iii)$] The canonical map $\theta: \pi^{\ast} \pi_{\ast} \calF \rightarrow \calF$ is an equivalence.
\end{itemize}
\end{proposition}

Before giving the proof, we need an easy lemma.

\begin{lemma}\label{sooch}
Let $X$ be a contractible paracompact topological space, let
$\pi_{\ast}: \Shv(X) \rightarrow \Shv(\ast) \simeq \SSet$ be the global sections functor, and let
$\pi^{\ast}$ be a right adjoint to $\pi_{\ast}$. Then $\pi^{\ast}$ is fully faithful.
\end{lemma}

\begin{proof}
Let $K$ be a Kan complex (regarded as an object of $\SSet$; we wish to prove that the unit map $u: K \rightarrow \pi_{\ast} \pi^{\ast} K$ is an equivalence. The results of \S \toposref{paracompactness} show that $\pi_{\ast} \pi^{\ast} K$ has the homotopy type of the Kan complex of maps
$\bHom_{ \Top}( X, |K|)$. Under this identification, the map $u$ corresponds to the diagonal
inclusion $K \rightarrow \Sing |K| \simeq \bHom_{\Top}( \ast, |K|) \rightarrow \bHom_{ \Top}( X, |K| )$.
Since $X$ is contractible, this inclusion is a homotopy equivalence.
\end{proof}

\begin{proof}[Proof of Proposition \ref{sooly}]
The implications $(iii) \Rightarrow (ii) \Rightarrow (i)$ are obvious. We prove that
$(ii) \Rightarrow (iii)$. Suppose that $\calF$ is constant; then $\calF \simeq \pi^{\ast} K$ for some
$K \in \SSet$. Then $\theta$ admits a right homotopy inverse, given by applying $\pi^{\ast}$ to the unit map $u: K \rightarrow \pi_{\ast} \pi^{\ast} K$. It follows from Lemma \ref{sooch} that
$u$ is an equivalence, so that $\theta$ is an equivalence as well.

We now prove that $(i) \Rightarrow (ii)$. Assume that $\calF$ is locally constant.
Let $S \subseteq [0,1]$ be the set of real numbers $t$ such that
$\calF$ is constant in some neighborhood of the interval $[0,t] \subseteq [0,1]$.
Let $s$ be the supremum of the set $S$ (since $\calF$ is constant in a neighborhood of $0$,
we must have $s > 0$). We will show that $s \in S$. It will follow
that $s = 1$ (otherwise, since $\calF$ is locally constant on $[0, s + \epsilon]$ for
$\epsilon$ sufficiently small, we would have $s + \frac{ \epsilon}{2} \in S$) so that
$\calF$ is locally constant on $[0,1]$, as desired.

Since $\calF$ is locally constant, it is constant when restricted to some open neighborhood
$U$ of $s \in [0,1]$. Since $s$ is a limit point of $S$, we have $S \cap U \neq \emptyset$.
Consequently, we can choose some point $t \in S \cap U$, so that $\calF$ is constant on
$U$ and on $[0,t)$. We will prove that $\calF$ is constant on the neighborhood $V = U \cup [0,t)$ of
$[0,s]$, so that $s \in S$ as desired.

Since $\calF$ is constant on $[0,t)$, we have an equivalence $\alpha: (\calF | [0,t) ) \simeq ( \pi^{\ast} K | [0,t) )$ for some object $K \in \SSet$. Similarly, we have an equivalence
$\beta: ( \calF | U ) \simeq ( \pi^{\ast} K' | U)$ for some $K' \in \SSet$. Restricting to the intersection,
we get an equivalence $\gamma: ( \pi^{\ast} K | U \cap [0,t) ) \simeq (\pi^{\ast} K' | U \cap [0,t) )$.
Since the intersection $U \cap [0,t)$ is contractible, Lemma \ref{sooch} guarantees that
$\gamma$ is induced by an equivalance $\gamma_0: K \simeq K'$ in the $\infty$-category $\SSet$.
Identifying $K$ with $K'$ via $\gamma_0$, we can reduce to the case where $K = K'$ and
$\gamma'$ is homotopic to the identity. 
For every open subset $W \subseteq [0,1]$, let $\chi_{W} \in \Shv(X)$ denote the sheaf given by the formula
$$ \chi_{W}(W') = \begin{cases} \ast & \text{ if } W' \subseteq W \\
\emptyset & \text{ otherwise.} \end{cases}$$
We then have a commutative diagram
$$ \xymatrix{ \pi^{\ast} K \times \chi_{ U \cap [0,t) } \ar[r] \ar[d] & \pi^{\ast} K \times \chi_U \ar[d] \\
\pi^{\ast} K \times \chi_{[0,t)} \ar[r] & \calF. }$$
This diagram induced a map $\pi^{\ast} K \times \chi_{V} \rightarrow \calF$, which
determines the required equivalence $\pi^{\ast} K | V \simeq \calF | V$.
\end{proof}

\begin{remark}
Proposition \ref{sooly} remains valid (with essentially the same proof) if we replace the closed
unit interval $[0,1]$ by an open interval $(0,1)$ or a half-open interval $[0,1)$.
\end{remark}

Let $h_0, h_1: X \rightarrow Y$ be a pair of continuous maps from a topological space $X$ to another
topological space $Y$. If $h_0$ is homotopic to $h_1$, then there exists a continuous map
$h: X \times \R \rightarrow Y$ such that $h_0 = h | X \times \{0\}$ and $h_1 = h | X \times \{1\}$. 
In this case, we can attempt to understand the relationship between the pullbacks
$h_0^{\ast} \calF$ and $h_1^{\ast} \calF$ of a sheaf $\calF$ on $Y$ by studying the pullback
$h^{\ast} \calF \in \Shv(Y \times \R)$. If $\calF$ is locally constant, then so is $h^{\ast} \calF$. It will
be convenient for us to consider a more general situation where $\calF$ is only required to be locally constant {\em along} the paths $h | (\{y\} \times \R)$ (and, for technical reasons, hypercomplete). The following definition axiomatizes the expected properties of the pullback $h^{\ast} \calF$: 

\begin{definition}
Let $X$ be a topological space and let $\calF \in \Shv(X \times \R)$. We will say that
$\calF$ is {\it foliated} if the following conditions are satisfied:
\begin{itemize}
\item[$(i)$] The sheaf $\calF$ is hypercomplete (see \S \toposref{hyperstacks}).
\item[$(ii)$] For every point $x \in X$, the restriction $\calF | ( \{x\} \times \R )$ is constant.
\end{itemize}
\end{definition}

The main result of this section is the following result, which should be regarded as a relative version of Proposition \ref{sooly} (where we have replaced the unit interval $[0,1]$ with the entire real line):

\begin{proposition}\label{postcanca}
Let $X$ be a topological space, let $\pi: X \times \R \rightarrow X$ denote the projection, and let $\calF \in \Shv(X \times \R)$. The following conditions are equivalent:
\begin{itemize}
\item[$(1)$] The sheaf $\calF$ is foliated.
\item[$(2)$] The pushforward $\pi_{\ast} \calF$ is hypercomplete, and the counit map
$v: \pi^{\ast} \pi_{\ast} \calF \rightarrow \calF$ is an equivalence.
\end{itemize}
\end{proposition}

The proof of Proposition \ref{postcanca} will require a few preliminaries. 

\begin{lemma}\label{sceen}
Let $f^{\ast}: \calX \rightarrow \calY$ be a geometric morphism of $\infty$-topoi.
Assume that $f^{\ast}$ admits a left adjoint $f_{!}$. Then $f^{\ast}$ carries
hypercomplete objects of $\calX$ to hypercomplete objects of $\calY$.
\end{lemma}

\begin{proof}
To show that $f^{\ast}$ preserves hypercomplete objects, it will suffice to show that
the left adjoint $f_{!}$ preserves $\infty$-connective morphisms. We will show that
$f_{!}$ preserves $n$-connective morphisms for every nonnegative integer $n$.
This is equivalent to the assertion that $f^{\ast}$ preserves $(n-1)$-truncated morphisms, which follows from Proposition \toposref{eaa}.
\end{proof}

\begin{example}\label{spazt}
Every \etale map of $\infty$-topoi satisfies the hypothesis of Lemma \ref{sceen}. 
Consequently, if $X$ is a hypercomplete object of an $\infty$-topos $\calX$, then
$X \times U$ is a hypercomplete object of $\calX_{/U}$ for each $U \in \calX$.
\end{example}

\begin{example}\label{tsooch}
Let $X$ and $Y$ be topological spaces, and let $\pi: X \times Y \rightarrow X$ be the projection. Assume that $Y$ is locally compact and locally of constant shape. Then $\pi^{\ast}$ satisfies the hypothesis of
Lemma \ref{sceen}, and therefore preserves hypercompletess. To prove this, we observe
that $\Shv(X \times Y)$ can be identified with $\Shv(X) \otimes \Shv(Y)$, where $\otimes$ denotes the tensor product operation on presentable $\infty$-categories described in \S \symmetricref{comm7}: this follows from Proposition \toposref{cartmun} and Example \monoidref{swikka}. The functor $\pi^{\ast}$ can be identified with the tensor product $\id_{ \Shv(X)} \otimes {\pi'}^{\ast}$, where 
$\pi': Y \rightarrow \ast$ is the projection. Proposition \ref{postkus} guarantees that
${\pi'}^{\ast}$ admits a left adjoint $\pi'_{!}$. It follows that $\id_{ \Shv(X)} \otimes \pi'_{!}$ is a left adjoint to $\pi^{\ast}$. Moreover, if ${\pi'}^{\ast}$ is fully faithful, then the counit map
$v: \pi'_{!} {\pi'}^{\ast} \rightarrow \id$ is an equivalence, so the counit map
$\pi_{!} \pi^{\ast} \rightarrow \id_{ \Shv(X)}$ is also an equivalence: it follows that
$\pi^{\ast}$ is fully faithful.
\end{example}

\begin{lemma}\label{ander}
Let $X$ be a topological space and let $\pi: X \times (0,1) \rightarrow X$ denote the projection. Then the pullback functor $\pi^{\ast}: \Shv(X) \rightarrow \Shv(X \times (0,1))$ is fully faithful
(so that the unit map $\calF \rightarrow \pi_{\ast} \pi^{\ast} \calF$ is an equivalence for
every $\calF \in \Shv(X)$).
\end{lemma}

\begin{proof}
Let $\psi: (0,1) \rightarrow \ast$ denote the projection map, and let
$\psi^{\ast}: \SSet \rightarrow \Shv( (0,1) )$ be the associated geometric morphism.
Then $\psi^{\ast}$ admits a left adjoint $\psi_{!}$ (Proposition \ref{postkus}) and
the counit transformation $v: \psi_{!} \psi^{\ast} \rightarrow \id$ is an equivalence
of functors from $\SSet$ to itself. As in Example \ref{tsooch}, we can identify
$\Shv( X \times (0,1) )$ with the tensor product $\Shv(X) \otimes \Shv( (0,1) )$,
so that $\psi_{!}$ and $\psi^{\ast}$ induce a pair of adjoint functors
$$ \Adjoint{ F}{ \Shv(X \times (0,1) )}{\Shv(X ).}{G}$$
The functor $G$ can be identified with $\pi^{\ast}$. Since the counit map $v$
is an equivalence, the counit $F \circ G \rightarrow \id_{ \Shv(X)}$ is likewise an equivalence,
which proves that $G \simeq \pi^{\ast}$ is fully faithful.
\end{proof}

\begin{variant}\label{tsongv}
In the statement of Lemma \ref{ander}, we can replace $(0,1)$ by a closed or half-open interval.
\end{variant}

\begin{proof}[Proof of Proposition \ref{postcanca}]
Suppose first that $(2)$ is satisfied, and let $\calG = \pi_{\ast} \calF$. Then 
$\calG$ is hypercomplete, so $\pi^{\ast} \calG$ is hypercomplete (Example \ref{tsooch});
since $v: \pi^{\ast} \calG \rightarrow \calF$ is an equivalence, it follows that $\calF$ is hypercomplete.
It is clear that $\calF \simeq \pi^{\ast} \calG$ is constant along $\{x \} \times \R$, for each $x \in X$.

Conversely, suppose that $\calF$ is foliated. To prove that $\pi_{\ast} \calF$ is hypercomplete, it suffices to show that $\pi_{\ast} \calF$ is local with respect to every $\infty$-connective morphism $\alpha$ in $\Shv(X)$. This is equivalent to the requirement that $\calF$ is local with respect to $\pi^{\ast}(\alpha)$.
This follows from our assumption that $\calF$ is hypercomplete, since $\pi^{\ast}(\alpha)$ is again $\infty$-connective. To complete the proof that $(1) \Rightarrow (2)$, it will suffice to show that
the counit map $v: \pi^{\ast} \calG \rightarrow \calF$ is an equivalence.

For each positive integer $n$, let $\calF_n = \calF | (X \times (-n,n)) \in \Shv(X \times (-n,n))$, let
$\pi_{n}: X \times (-n, n) \rightarrow X$ be the projection map, and let $\calG_{n} = (\pi_{n})_{\ast} \calF$.
We have a commutative diagram
$$ \xymatrix{ (\pi^{\ast} \calG) | (X \times (-n,n) ) \ar[r]^{v} \ar[d] & \calF| ( X \times (-n,n) ) \ar[d] \\
\pi_{n}^{\ast} \calG_{n} \ar[r]^{v_n} & \calF_{n}. }$$
To prove that $v$ is an equivalence, it will suffice to show that the left vertical and lower horizontal maps in this diagram are equivalences (for each $n$). This will follow from the following
pair of assertions:
\begin{itemize}
\item[$(a)$] For each $n > 0$, the restriction map $\calG_{n+1} \rightarrow \calG_{n}$ is an equivalence
(so that $\calG \simeq \varprojlim_{n} \calG_{n}$ is equivalent to each $\calG_{n}$).
\item[$(b)$] For each $n > 0$, the map $\pi_{n}^{\ast} \calG_{n} \rightarrow \calF_{n}$ is an equivalence.
\end{itemize}
Note that assertion $(a)$ follows from $(b)$: if we let $i: X \rightarrow X \times \R$ be the map
induced by the inclusion $\{0\} \hookrightarrow \R$, then we have a commutative diagram
$$ \xymatrix{ \calG_{n} \ar[r] \ar[d] & \calG_{n+1} \ar[d] \\
i^{\ast} \pi_{n}^{\ast} \calG_{n} \ar[r] \ar[d] & i^{\ast} \pi_{n+1}^{\ast} \calG_{n+1} \ar[d] \\
i^{\ast} \calF_{n} \ar[r]^{s} & i^{\ast} \calF_{n+1} }$$
in which the upper vertical maps are equivalences, the lower horizontal maps are equivalences
by $(b)$, and the map $s$ is an equivalence by construction.

To prove $(b)$, let $\calF_{n}^{+} \in \Shv( X \times [-n,n] )$ denote the hypercompletion of the
restriction $\calF | (X \times [-n,n])$, let $\pi^{n}: X \times [-n,n] \rightarrow X$ be the projection, and let $\calG_{n}^{+} = \pi^{n}_{\ast} \calF_{n}^{+}$. Let $v': (\pi^{n})^{\ast} \calG_{n}^{+} \rightarrow \calF_{n}^{+}$ be the counit map. We claim that $v'$ is an equivalence. Since $\calF_{n}^{+}$ is hypercomplete by assumption, $\calG_{n}^{+} \simeq \pi^{n}_{\ast} \calF_{n}^{+}$ is likewise hypercomplete and
so $(\pi^{n})^{\ast} \calG_{n}^{+}$ is hypercomplete by virtue of Example \ref{tsooch}. Consequently, to prove that $v'$ is an equivalence, it will suffice to show that $v'$ is $\infty$-connective.
To prove this, choose a point
$x \in X$ and let $j: [-n,n] \rightarrow X \times [-n,n]$ be the map induced by the inclusion
$j': \{x \} \hookrightarrow X$. We will show that $j^{\ast}(v')$ is an equivalence.
Consider the diagram of $\infty$-topoi
$$ \xymatrix{ \Shv( [-n,n] ) \ar[r]^{j_{\ast}} \ar[d]^{\psi_{\ast}} & \Shv( X \times [-n,n]) \ar[d]^{\pi^{n}_{\ast} } \ar[r] & \Shv( [0,1] ) \ar[d]^{\psi_{\ast}} \\
\Shv( \ast) \ar[r]^{j'_{\ast}} & \Shv(X) \ar[r] & \Shv( \ast) . }$$
The right square and the outer rectangle are pullback diagrams (Proposition \toposref{cartmun}),
so the left square is a pullback diagram as well. Moreover, the geometric morphism $\psi_{\ast}$ is proper (Corollary \toposref{compactprop}), so that $\pi^{n}_{\ast}$ is likewise proper and the
push-pull morphism $e: {j'}^{\ast} \pi^{n}_{\ast} \rightarrow \psi_{\ast} j^{\ast}$ is an equivalence.
We have a commutative diagram
$$ \xymatrix{ \psi^{\ast} {j'}^{\ast} \pi^{n}_{\ast} \calF^{+}_{n} \ar[r] \ar[d]^{e} & j^{\ast} (\pi^{n})^{\ast} \pi^{n}_{\ast} \calF^{+}_{n} \ar[d]^{j^{\ast}(v')} \\
\psi^{\ast} \psi_{\ast} j^{\ast} \calF^{+}_{n} \ar[r]^{v'_{x}} & j^{\ast} \calF^{+}_{n}. }$$
By virtue of the above diagram (and the fact that $e$ is an equivalence), we are reduced to proving that $v'_{x}$ is an equivalence.
To prove this, it suffices to verify that $j^{\ast} \calF^{+}_{n} \in \Shv( [-n,n])$ is constant
(Proposition \ref{sooly}). We have an $\infty$-connective morphism
$\theta: \calF | ( \{x\} \times [-n,n]) \rightarrow j^{\ast} \calF^{+}_{n}$. Since every open subset of
the topological space $[-n,n]$ has covering dimension $\leq 1$, the $\infty$-topos
$\Shv( [-n,n])$ is locally of homotopy dimension $\leq 1$ (Theorem \toposref{paradimension}) and
therefore hypercomplete. It follows that $\theta$ is an equivalence. Since
$\calF$ is foliated, the restriction $\calF | ( \{x\} \times [-n,n])$ is constant, from which it follows immediately that $j^{\ast} \calF^{+}_{n}$ is constant as well.

The $\infty$-connective
morphism $\calF | ( X \times [-n,n]) \rightarrow \calF^{+}_{n}$ induces another $\infty$-connective
morphism $\alpha: \calF_{n} \rightarrow \calF^{+}_{n} | ( X \times (-n,n) )$. Since
the domain and codomain of $\alpha$ are both hypercomplete (Example \ref{spazt}), 
we deduce that $\alpha$ is an equivalence. In particular, we have
$\calF_{n} \simeq ((\pi^{n})^{\ast} \calG_{n}^{+}) | (X \times (-n,n) = \pi_{n}^{\ast} \calG^{+}_{n}$.
Thus $\calF_{n}$ lies in the essential image of the functor $\pi_{n}^{\ast}$, which is fully faithful by virtue of Lemma \ref{ander}. It follows that that the counit map $\pi_{n}^{\ast} (\pi_{n})_{\ast} \calF_{n} \rightarrow \calF_{n}$ is an equivalence as desired.
\end{proof}

\subsection{Singular Shape}\label{lubla}

In \S \ref{clome}, we defined the notion of a {\em locally constant} object of an $\infty$-topos
$\calX$. Moreover, we proved that the $\infty$-topos $\calX$ is locally of constant shape, then
the $\infty$-category of locally constant objects of $\calX$ is equivalent to the $\infty$-topos
$\SSet_{/K}$ of spaces lying over some fixed object $K \in \SSet$ (Theorem \ref{san}). This can be regarded as an analogue of the main result in the theory of covering spaces, which asserts that the category of covering spaces of a sufficiently nice topological space $X$ can be identified with the category of sets acted on by the fundamental group of $X$. If we apply Theorem \ref{san} in the special case $\calX = \Shv(X)$, then we deduce that the fundamental groups of $X$ and $K$ are isomorphic to one another. Our objective in this section is to strengthen this observation: we will show that if $X$ is a sufficiently nice topological space, then the
$\infty$-topos $\Shv(X)$ of sheaves on $X$ is locally of constant shape, and the shape $K$ of
$\Shv(X)$ can be identified with the singular complex $\Sing(X)$.

\begin{remark}
We refer the reader to \cite{toengalois} for a closely related discussion, at least in the case where
$X$ is a CW complex.
\end{remark}

Our first step is to describe a class of topological spaces $X$ for which the theory of
locally constant sheaves on $X$ is well-behaved. By definition, if $\calF$ is a locally constant
sheaf on $X$, then every point $x \in X$ has an open neighborhood $U$ such that the restriction
$\calF| U$ is constant. Roughly speaking, we want a condition on $X$ which guarantees that we can
choose $U$ to be independent of $\calF$. 

\begin{definition}\label{copely}
Let $f^{\ast}: \calX \rightarrow \calY$ be a geometric morphism of $\infty$-topoi.
We will say that $f^{\ast}$ is a {\it shape equivalence} if it induces an equivalence
of functors $\pi_{\ast} \pi^{\ast} \rightarrow \pi_{\ast} f_{\ast} f^{\ast} \pi^{\ast}$, where
$\pi^{\ast}: \SSet \rightarrow \calX$ is a geometric morphism.
\end{definition}

\begin{remark}
Let $\calX$ be an $\infty$-topos. Then $\calX$ has constant shape if and only if there
exists a shape equivalence $f^{\ast}: \SSet_{/K} \rightarrow \calX$, for some 
Kan complex $K$. The ``if'' direction is obvious (since $\SSet_{/K}$ is of constant shape).
Conversely, if $\calX$ is of constant shape, then $\pi_{\ast} \pi^{\ast}$ is corepresentable
by some object $K \in \SSet$. In particular, there is a canonical map
$\Delta^0 \rightarrow \pi_{\ast} \pi^{\ast} K$, which we can identify with a map
$\alpha: {\bf 1} \rightarrow \pi^{\ast} K$ in the $\infty$-topos $\calX$, where ${\bf 1}$ denotes the final
object of $\calX$. According to Proposition \toposref{goodking}, $\alpha$ determines a geometric morphism of $\infty$-topoi $f^{\ast}: \SSet_{/K} \rightarrow \calX$, which is easily verified to be a shape equivalence.
\end{remark}

\begin{definition}\label{underman}
Let $f: X \rightarrow Y$ be a continuous map of topological spaces. We will say
that $f$ is a {\it shape equivalence} if the associated geometric morphism
$f^{\ast} \Shv(X) \rightarrow \Shv(Y)$ is a shape equivalence, in the sense of
Definition \ref{copely}.
\end{definition}

\begin{example}
Let $f: X \rightarrow Y$ be a continuous map between paracompact topological spaces.
Then $f$ is a shape equivalence in the sense of Definition \ref{underman} if and only if,
for every CW complex $Z$, composition with $f$ induces a homotopy equivalence of Kan complexes
$\bHom_{ \Top}(Y, Z) \rightarrow \bHom_{\Top}(X,Z)$.
\end{example}

\begin{example}\label{sadm}
If $X$ is any topological space, then the projection map $\pi: X \times \R \rightarrow X$
is a shape equivalence. This follows immediately from the observation that
$\pi^{\ast}$ is fully faithful (Example \ref{tsooch}).
\end{example}

\begin{remark}\label{sumn}
It follows from Example \ref{sadm} that every homotopy equivalence of topological spaces is also a shape equivalence.
\end{remark}

\begin{warning}
For general topological spaces, Definition \ref{underman} does {\em not} recover the classical
notion of a shape equivalence (see, for example, \cite{shapetheory}). However, if $X$ and $Y$ are both paracompact then we recover the usual notion of strong shape equivalence (Remark \toposref{struke}).
\end{warning}

\begin{definition}
Let $X$ be a topological space. We will say that $X$ has {\it singular shape} if
the counit map $| \Sing(X) | \rightarrow X$ is a shape equivalence.
\end{definition}

\begin{remark}
If $X$ is a topological space with singular shape, then the $\infty$-topos
$\Shv(X)$ has constant shape: indeed, $\Shv(X)$ is shape equivalent to
$\Shv( | \Sing(X) |)$, and $| \Sing(X) |$ is a CW complex (Remark \ref{siba}).
\end{remark}

\begin{remark}\label{adlam}
Let $f: X \rightarrow Y$ be a homotopy equivalence of topological spaces.
Then $X$ has singular shape if and only if $Y$ has singular shape. This follows
immediately from Remark \ref{sumn} by inspecting the diagram
$$ \xymatrix{ | \Sing(X)|  \ar[r] \ar[d] & | \Sing(Y) | \ar[d] \\
X \ar[r] & Y. }$$ 
\end{remark}

\begin{example}
Let $X$ be a paracompact topological space. Then $X$ has singular shape if and only if,
for every CW complex $Y$, the canonical map
$$ \bHom_{ \Top}(X,Y) \rightarrow \bHom_{ \sSet}( \Sing(X), \Sing(Y)) \simeq \bHom_{ \Top}( | \Sing(X) |, Y)$$ is a homotopy equivalence of Kan complexes. 
\end{example}

\begin{remark}
Let $X$ be a paracompact topological space. There are two different ways that we might try to assign
to $X$ a homotopy type. The first is to consider continuous maps from nice spaces (such as CW complexes) {\em into} the space $X$. Information about such maps is encoded in the Kan complex
$\Sing(X) \in \SSet$, which controls the {\it weak homotopy type} of $X$. Alternatively, we can instead consider maps from $X$ into CW complexes. These are controlled by the pro-object
$\Sh(X)$ of $\SSet$ which corepresents the functor $K \mapsto \bHom_{ \Top}(X, |K|)$. 
There is a canonical map $\Sing(X) \rightarrow \Sh(X)$, and $X$ has singular shape if and only if this map is an equivalence.
\end{remark}

\begin{lemma}\label{kuil}
Let $X$ be a topological space, and let
$\{ U_{\alpha} \in \calU(X) \}_{ \alpha \in A}$ be an open covering of $X$.
Assume that for every nonempty finite subset $A_0 \subseteq A$, the intersection
$U_{A_0} = \bigcap_{\alpha \in A_0} U_{\alpha}$ has singular shape. Then $X$ has singular shape.
\end{lemma}

\begin{proof}
Let $\pi^{\ast}: \SSet \rightarrow \Shv(X)$ be a geometric morphism.
For each open set $U \subseteq X$, let $F_{U}: \SSet \rightarrow \SSet$ be
the functor given by composing $\pi^{\ast}$ with evaluation at $U$, and let
$G_{U}: \SSet \rightarrow \SSet$ be the functor given by 
$K \mapsto \Fun( \Sing(U), K)$. There is a natural transformation of
functors $\gamma_U: F_{U} \rightarrow G_U$, and $U$ has singular shape if and only if
$\gamma_U$ is an equivalence. We observe that
$F_{X}$ can be identified with a limit of the diagram $\{ F_{U_{A_0}} \}$ where $A_0$ ranges over the finite subsets of $A$,
and that $G_X$ can be identified with a limit of the diagram $\{ G_{ U_{A_0}} \}$ (since
$\Sing(X)$ is the homotopy colimit of $\{ \Sing(U_{A_0}) \}$ by Theorem \ref{vankamp}).
Under these identifications, $\gamma_{X}$ is a limit of the functors
$\{ \gamma_{ U_{A_0}} \}$. Since each of these functors is assumed to be an equivalence, we deduce that $\gamma_{X}$ is an equivalence.
\end{proof}

\begin{definition}
We will say that topological space $X$ is {\it locally of singular shape} if
every open set $U \subseteq X$ has singular shape.
\end{definition}

\begin{remark}\label{locsh}
Let $X$ be a topological space. Suppose that
$X$ admits a covering by open sets which are locally of singular shape. Then $X$ is locally of singular shape (this follows immediately from Lemma \ref{kuil}).
\end{remark}

Let $X$ be a topological space which is locally of singular shape. Then $\Shv(X)$ is locally of constant shape,
and the shape of $\Shv(X)$ can be identified with the Kan complex $\Sing(X)$. It follows from
Theorem \ref{san} that the $\infty$-category of locally constant objects of $\Shv(X)$ is equivalent to
$\SSet_{ / \Sing(X) }$. Our goal for the remainder of this section is to give a more explicit description of this equivalence.

\begin{construction}
Let $X$ be a topological space. We let
$\bfA_{X}$ denote the category $( \sSet)_{/ \Sing(X) }$, endowed
with the usual model structure.
Let $\bfA_{X}^{\degree}$ denote the full subcategory of $\bfA_{X}$ spanned by the
fibrant-cofibrant objects (these are precisely the Kan fibrations
$Y \rightarrow \Sing(X)$).

We define a functor $\theta: \calU(X)^{op} \times \bfA_{X} \rightarrow \sSet$ by
the formula $\theta(U,Y) = \Fun_{ \Sing(X)}( \Sing(U), Y)$. Restricting to
$\bfA_{X}^{\degree}$ and passing to nerves, we get a map of $\infty$-categories
$\Nerve( \calU(X)^{op} ) \times \Nerve( \bfA_{X}^{\degree}) \rightarrow \SSet$, which
we regard as a map of $\infty$-categories $\Nerve( \bfU_{X}^{\degree}) \rightarrow
\calP( \calU(X) )$. It follows from Variant \ref{varswug} on Proposition \ref{swugman} that
this functor factors through the full subcategory $\Shv(X) \subseteq \calP( \calU(X) )$ spanned by the sheaves on $X$. We will denote the underlying functor
$\Nerve( \bfA_{X}^{\degree}) \rightarrow \Shv(X)$ by $\Psi_{X}$.
\end{construction}

\begin{example}\label{complus}
Let $X$ be a topological space.
The construction $K \mapsto K \times \Sing(X)$ determines a functor from
$\sSet \simeq \bfA_{\ast}$ to $\bfA_{X}$, which restricts to a functor
$\bfA_{\ast}^{\degree} \rightarrow \bfA_{X}^{\degree}$. Passing to nerves
and composing with $\Psi_{X}$, we get a functor
$\psi: \SSet \rightarrow \Shv(X)$, which carries a Kan complex $K$ to the sheaf
$U \mapsto \bHom_{ \sSet}( \Sing(U), K)$. Let $\pi_{\ast}: \Shv(X) \rightarrow \SSet$
be the functor given by evaluation on $X$. There is an evident natural transformation
$\id_{\SSet} \rightarrow \pi_{\ast} \circ \psi$, which induces a natural transformation
$\pi^{\ast} \rightarrow \psi$. The space $X$ is locally of singular shape if and only if
this natural transformation is an equivalence.
\end{example}

We note that the object $\psi \Sing(X) \in \Shv(X)$ has a canonical global section
given by the identity map from $\Sing(X)$ to itself. If $Y \rightarrow \Sing(X)$ is
any Kan fibration, then $\Psi_X(Y)$ can be identified with the (homotopy) fiber
of the induced map $\psi(Y) \rightarrow \psi( \Sing(X))$. It follows that the functor $\Psi_{X}$ is an explicit model for the fully faithful embedding
described in Proposition \ref{standy}. Coupling this observation with
Theorem \ref{san}, we obtain the following:

\begin{theorem}\label{squ}
Let $X$ be a topological space which is locally of singular shape. Then
the functor $\Psi_{X}: \Nerve( \bfA^{\degree}_{X}) \rightarrow \Shv(X)$ is a fully faithful
embedding, whose essential image is the full subcategory of $\Shv(X)$ spanned by the locally constant sheaves on $X$.
\end{theorem}

\subsection{Constructible Sheaves}\label{coofer}

In \S \ref{clome} and \S \ref{lubla}, we studied the theory of {\em locally constant} sheaves
on a topological space $X$. In many applications, one encounters sheaves $\calF \in \Shv(X)$ which
are not locally constant but are nevertheless {\em constructible}: that is, they are locally constant along
each stratum of a suitable stratification of $X$. We begin by making this notion more precise.

\begin{definition}\label{cunning}
Let $A$ be a partially ordered set. We will regard $A$ as a topological space,
where a subset $U \subseteq A$ is open if it is {\it closed upwards}: that is,
if $x \leq y$ and $x \in U$ implies that $y \in U$.

Let $X$ be a topological space. An {\it $A$-stratification} of $X$ is a continuous
map $f: X \rightarrow A$. Given an $A$-stratification of a space $X$ and an element $a \in A$, we let
$X_{a}$, $X_{\leq a}$, $X_{< a}$, $X_{\geq a}$, and $X_{> a}$ denote the subsets of
$X$ consisting of those points $x \in X$ such that $f(x) = a$, $f(x) \leq a$, $f(x) < a$,
$f(x) \geq a$, and $f(x) > a$, respectively.
\end{definition}

\begin{definition}
Let $A$ be a partially ordered set and let $X$ be a topological space equipped with an $A$-stratification. We will say that an object $\calF \in \Shv(X)$ is {\it $A$-constructible} if, for every
element $a \in A$, the restriction $\calF | X_{a}$ is a locally constant object of
$\Shv(X_{a})$. Here $\calF | X_{a}$ denotes the image of $\calF$ under the left adjoint
to the pushforward functor $\Shv(X_{a}) \rightarrow \Shv(X)$. 

We let $\Shv^{A}(X)$ denote the full subcategory of $\Shv(X)$ spanned by the $A$-constructible objects.
\end{definition}

To ensure that the theory of $A$-constructible sheaves is well-behaved, it is often convenient to
make the introduce a suitable regularity condition on the stratification $X \rightarrow A$. 

\begin{definition}
Let $A$ be a partially ordered set, and let $A^{\triangleleft}$ be the partially ordered set obtained
by adjoining a new smallest element $- \infty$ to $A$. Let $f: X \rightarrow A$ be an $A$-stratified space.
We define a new $A^{\triangleleft}$-stratified space $C(X)$ as follows:
\begin{itemize}
\item[$(1)$] As a set $C(X)$ is given by the union $\{ \ast \} \cup (X \times \R_{> 0})$.
\item[$(2)$] A subset $U \subseteq C(X)$ is open if and only if $U \cap (X \times \R_{ > 0 })$ is open,
and if $\ast \in U$ then $X \times (0, \epsilon) \subseteq U$ for some positive real number $\epsilon$.
\item[$(3)$] The $A^{\triangleleft}$-stratification of $C(X)$ is determined by the map
$\overline{f}: C(X) \rightarrow A^{\triangleleft}$ such that $\overline{f}(\ast) = - \infty$
and $\overline{f}(x,t) = f(x)$ for $(x,t) \in X \times \R_{ > 0 }$.
\end{itemize}
We will refer to $C(X)$ as the {\it open cone} on $X$.
\end{definition}

\begin{remark}
If the topological space $X$ is compact and Hausdorff, then the open cone $C(X)$ is homeomorphic to the pushout $(X \times \R_{\geq 0}) \coprod_{ X \times \{0\} } \{ \ast \}$.
\end{remark}

\begin{definition}\label{sally}
Let $A$ be a partially ordered set, let $X$ be an $A$-stratified space, and let
$x \in X_{a} \subseteq X$ be a point of $X$. We will say that
$X$ is {\it conically stratified at the point $x$} if there exists an
$A_{ > a}$-stratified topological space $Y$, a topological space $Z$, and
an open embedding $Z \times C(Y) \hookrightarrow X$ of $A$-stratified spaces
whose image $U_{x}$ contains $x$. Here we regard $Z \times C(Y)$ as endowed with
the $A$-stratification determined by the $A_{> a}^{\triangleleft} \simeq A_{\geq a}$-stratification of
$C(Y)$.

We will say that $X$ is {\it conically stratified} if it is conically stratified at every point $x \in X$.
\end{definition}

\begin{remark}
In Definition \ref{sally}, we do not require that the space $Y$ itself be conically stratified.
\end{remark}

\begin{definition}
We will say that a partially ordered set $A$ satisfies the {\it ascending chain condition} if
every nonempty subset of $A$ has a maximal element.
\end{definition}

\begin{remark}
Equivalently, $A$ satisfies the ascending chain condition if there does not exist any infinite ascending sequence $a_0 < a_1 < \cdots$ of elements of $A$.
\end{remark}

The main goal of this section is to prove the following somewhat technical convergence result 
concerning constructible sheaves:

\begin{proposition}\label{toru}
Let $A$ be a partially ordered set, and let $X$ be an $A$-stratified space. Assume that:
\begin{itemize}
\item[$(i)$] The space $X$ is paracompact and locally of singular shape.
\item[$(ii)$] The $A$-stratification of $X$ is conical.
\item[$(iii)$] The partially ordered set $A$ satisfies the ascending chain condition.
\end{itemize}
Let $\calF \in \Shv^{A}(X)$ be an $A$-constructible sheaf. 
Then the canonical map
$\theta: \calF \rightarrow \varprojlim \tau_{\leq n} \calF$ is an equivalence. In particular,
$\calF$ is hypercomplete.
\end{proposition}

The proof of Proposition \ref{toru} will require several preliminaries, and will be given at the end of this section. 
Our first step is to consider the case of a very simple stratification of $X$: namely, a decomposition of
$X$ into an open set and its closed complement. The following result is useful for working with constructible sheaves: it allows us to reduce global questions to questions which concern individual strata.

\begin{lemma}\label{singbird}
Let $\calX$ be an $\infty$-topos and $U$ a $(-1)$-truncated object of $\calX$.
Let $i^{\ast}: \calX \rightarrow \calX/U$ and $j^{\ast}: \calX \rightarrow \calX_{/U}$ be the
canonical geometric morphisms, $j_{\ast}$ a right adjoint to $j^{\ast}$, and let $\overline{p}: K^{\triangleleft} \rightarrow \calX$
be a small diagram in $\calX$ indexed by a weakly contractible simplicial set $K$.
Suppose that $i^{\ast} \overline{p}$, $j^{\ast} \overline{p}$, and
$i^{\ast} j_{\ast} j^{\ast} \overline{p}$ are all limit diagrams. Then $\overline{p}$ is a limit diagram.
\end{lemma}

\begin{proof}
Let $\calF$ denote the image of the cone point of $K^{\triangleleft}$ under $\overline{p}$,
let $p': K \rightarrow \calX$ be the constant diagram taking the value $\calF$, and let
$p = \overline{p} | K$. Then $\overline{p}$ determines a natural transformation
of diagrams $\alpha: p' \rightarrow p$; we wish to prove that $\alpha$ induces an
equivalence $\varprojlim(p') \rightarrow \varprojlim(p)$ in $\calX$.
For this, it suffices to show that for every object $V \in \calX$, the induced map
$$ \theta: \bHom_{ \calX}( V, \varprojlim(p') ) \rightarrow \bHom_{ \calX}( V, \varprojlim(p) )$$
is a homotopy equivalence. Replacing $\calX$ by $\calX_{/V}$, we can reduce to the
case where $V$ is the final object of $\calX$. In this case, we let $\Sect$ denote the functor
$\calX \rightarrow \SSet$ corepresented by $V$ (the functor of global sections). 

Fix a point $\eta \in \Sect( \varprojlim(p) )$; we will show that the homotopy fiber
of $\theta$ over $\eta$ is contractible. Let $j_{\ast}$ denote a right adjoint to
$j^{\ast}$, let $q = j_{\ast} \circ j^{\ast} \circ p$, and let $q' = j_{\ast} \circ j^{\ast} \circ p'$.
Then $\eta$ determines a point $\eta_0 \in \Sect( \varprojlim(q) )$. 
Since $j^{\ast} \circ \overline{p}$ is a limit diagram (and the functor $j_{\ast}$ preserves limits),
the canonical map $\varprojlim(q') \rightarrow \varprojlim(q)$ is an equivalence, so we
can lift $\eta_0$ to a point $\eta_1 \in \Sect( \varprojlim(q') )$. This point determines
a natural transformation from the constant diagram $c: K \rightarrow \calX$ taking the value
$V \simeq {\bf 1}$ to the diagram $q'$. Let $p'_0 = c \times_{q'} p'$ and
let $p_0 = c \times_{ q} p$. We have a map of homotopy fiber sequences
$$ \xymatrix{ \Sect( \varprojlim(p'_0) ) \ar[r] \ar[d]^{\theta'} & \Sect( \varprojlim(p') ) \ar[r] \ar[d]^{\theta} & \Sect( \varprojlim(q') ) \ar[d]^{\theta''} \\
\Sect( \varprojlim( p_0) ) \ar[r] & \Sect( \varprojlim( p) ) \ar[r] & \Sect( \varprojlim(q) ). }$$
Here $\theta''$ is a homotopy equivalence. Consequently, to prove that
the homotopy fiber of $\theta$ is contractible, it will suffice to show that $\theta'$ is a homotopy
equivalence.

By construction, the diagrams $p'_0$ and $p_0$ take values in the full subcategory $\calX/U \subseteq \calX$, so that the localization maps $p'_0 \rightarrow i^{\ast} p'_0$ and
$p_0 \rightarrow i^{\ast} p_0$ are equivalences. It therefore suffices to show that the map
$\Sect( \varprojlim( i^{\ast} p'_0) \rightarrow \Sect( \varprojlim( i^{\ast} p_0))$ is a homotopy equivalence. We have another map of homotopy fiber sequences
$$ \xymatrix{ \Sect( \varprojlim(i^{\ast} p'_0) ) \ar[r] \ar[d]^{\psi'} & \Sect( \varprojlim(i^{\ast} p') ) \ar[r] \ar[d]^{\psi} & \Sect( \varprojlim(i^{\ast} q') ) \ar[d]^{\psi''} \\
\Sect( \varprojlim(i^{\ast} p_0) ) \ar[r] & \Sect(i^{\ast} \varprojlim( p) ) \ar[r] & \Sect( \varprojlim(i^{\ast} q) ). }$$
The map $\psi$ is a homotopy equivalence by virtue of our assumption that
$i^{\ast} \overline{p}$ is a limit diagram, and the map $\psi''$ is a homotopy equivalence by virtue of
our assumption that $i^{\ast} j_{\ast} j^{\ast} \overline{p}$ is a limit diagram. It follows
that $\psi'$ is also a homotopy equivalence, as desired.
\end{proof}

\begin{lemma}\label{sadwell}
Let $\calX$ be an $\infty$-topos and $U$ a $(-1)$-truncated object of $\calX$.
Let $i^{\ast}: \calX \rightarrow \calX/U$ and $j^{\ast}: \calX \rightarrow \calX_{/U}$ be the
canonical geometric morphisms, and let $\alpha: \calF \rightarrow \calG$ be a morphism in
$\calX$. Suppose that $i^{\ast}( \alpha)$ and $j^{\ast}(\alpha)$ are equivalences. Then
$\alpha$ is an equivalence.
\end{lemma}

\begin{proof}
Apply Lemma \ref{singbird} in the special case where $K = \Delta^0$ (note that
$i^{\ast} j_{\ast}$ automatically preserves $j$-indexed limits).
\end{proof}

\begin{lemma}\label{cuppa}
Let $X$ be a paracompact topological space, $Y$ any topological space, $V$ an open neighborhood
of $X$ in $X \times C(Y)$. Then there exists a continuous function
$f: X \rightarrow (0, \infty)$ such that $V$ contains
$$V_{f} = \{ (x,y,t): t < f(x) \} \subseteq X \times Y \times (0, \infty) \subseteq X \times C(Y).$$
\end{lemma}

\begin{proof}
For each point $x \in X$, there exists a neighborhood $U_x$ of $x$ and a real number $t_x$
such that $\{ (x',y,t): t < t_x \wedge x' \in U_x \} \subseteq V$. Since $X$ is paracompact, we can choose a 
locally finite partition of unity $\{ \psi_{x} \}_{x \in X}$ subordinate to the cover $\{ U_x \}_{x \in X}$.
We now define $f(y) = \sum_{x \in X} \psi_x(y) t_{x}$.
\end{proof}

\begin{remark}\label{snark}
In the situation of Lemma \ref{cuppa}, the collection of open sets of the form $V_{f}$ is nonempty
(take $f$ to be a constant function) and stable under pairwise intersections ($V_{f} \cap V_{g} = V_{ \inf \{f,g\} })$. The collection of such open sets is therefore cofinal in partially ordered set of
all open subsets of $X \times C(Y)$ which contain $X$ (ordered by reverse inclusion).
\end{remark}

\begin{lemma}\label{sadly}
Let $X$ be a paracompact topological space. Let
$\pi$ denote the projection $X \times [0, \infty) \rightarrow X$,
let $j$ denote the inclusion $X \times (0, \infty) \hookrightarrow X \times [0, \infty)$, and let
$\pi_0 = \pi \circ j$. Then the obvious equivalence $\pi_0^{\ast} \simeq j^{\ast} \pi^{\ast}$
is adjoint to an equivalence of functors $\alpha: \pi^{\ast} \rightarrow j_{\ast} \pi_0^{\ast}$
from $\Shv(X)$ to $\Shv(X \times [0, \infty))$.
\end{lemma}

\begin{proof}
Let $\calF \in \Shv(X)$; we wish to prove that $\alpha$ induces an equivalence
$\pi^{\ast} \calF \rightarrow j_{\ast} \pi_0^{\ast} \calF$. It is clear that this map is an
equivalence when restricted to the open set $X \times (0, \infty)$. Let
$i: X \rightarrow X \times [0, \infty)$ be the map induced by the inclusion
$\{0\} \subseteq [0, \infty)$. By Corollary \ref{sadwell}, it will suffice to show that 
the map
$$\beta: \calF \simeq i^{\ast} \pi^{\ast} \calF \rightarrow i^{\ast} j_{\ast} \pi_0^{\ast} \calF$$
determined by $\alpha$ is an equivalence. Let $U$ be an open $F_{\sigma}$ subset of
$X$; we will show that the map $\beta_U: \calF(U) \rightarrow (i^{\ast} j_{\ast} \pi_0^{\ast} \calF)(U)$
is a homotopy equivalence. Replacing $X$ by $U$, we can assume that $U = X$.

According to Corollary \toposref{snottle}, we can identify $(i^{\ast} j_{\ast} \pi_0^{\ast} \calF)(X)$
with the colimit $\varinjlim_{V \in S}(j_{\ast} \pi_0^{\ast} \calF)(V) \simeq \varinjlim_{V \in S} (\pi_0^{\ast} \calF)(V-X)$, where $V$ ranges over the collection $S$ of all open neighborhoods of $X = X \times \{0\}$ in $X \times [0, \infty)$. Let $S' \subseteq S$ be the collection of all open neighborhoods of the
form $V_{f} = \{ (x,t): t < f(x) \}$, where $f: X \rightarrow (0, \infty)$ is a continuous function (see Lemma \ref{cuppa}). 
In view of Remark \ref{snark}, we have an equivalence
$\varinjlim_{V \in S} (\pi_0^{\ast} \calF)(V - X) \simeq \varinjlim_{ V \in S'}( \pi_0^{\ast} \calF)(V-X)$.
Since $S'$ is a filtered partially ordered set (when ordered by reverse inclusion), to prove that
$\beta_X$ is an equivalence it suffices to show that the pullback map
$\calF(X) \rightarrow (\pi_0^{\ast} \calF)( V_f - X)$ is a homotopy equivalence, for
every continuous map $f: X \rightarrow (0, \infty)$. Division by $f$ determines a homeomorphism
$V_{f} - X \rightarrow X \times (0,1)$, and the desired result follows from Lemma \ref{ander}.
\end{proof}

\begin{lemma}\label{torpus}
Let $X$ be a paracompact topological space of the form $Z \times C(Y)$, and consider the (noncommuting)
diagram
$$ \xymatrix{ Z \times Y \times (0, \infty) \ar[r]^{j} \ar[dr]^{\pi_0} & Z \times Y \times [0, \infty) \ar[d]_{\pi} \ar[r]^-{k}
& X \\
& Z \times Y \ar[r]^{\psi} & Z \ar[u]^{i}. }$$
Let $i'$ denote the inclusion $Z \times Y \rightarrow Z \times Y \times [0, \infty)$ given by 
$\{0\} \hookrightarrow [0, \infty)$. 
Assume that $X$ is paracompact. Then:
\begin{itemize}
\item[$(i)$] The canonical map $\alpha: \pi^{\ast} \rightarrow j_{\ast} \pi_0^{\ast}$ is an equivalence
of functors from $\Shv(Z \times Y)$ to $\Shv( Z \times Y \times [0, \infty) )$. 
\item[$(ii)$] Let $\beta: \pi^{\ast} \rightarrow i'_{\ast}$ be the natural transformation
adjoint to the equivalence ${i'}^{\ast} \pi^{\ast} \simeq \id_{ \Shv(Z \times Y)}$. 
Then the natural transformation
$$\gamma: i^{\ast} k_{\ast} \pi^{\ast} \stackrel{\beta}{\rightarrow}
i^{\ast} k_{\ast} i'_{\ast} \simeq i^{\ast} i_{\ast} \psi_{\ast} \rightarrow \psi_{\ast}$$
is an equivalence of functors from $\Shv( Z \times Y)$ to $\Shv(Z)$.
\item[$(iii)$] The functor $i^{\ast} j_{\ast} \pi_0^{\ast}$ is equivalent to $\psi_{\ast}$.
\end{itemize}
\end{lemma}

\begin{proof}
Note that $Z \times Y \simeq Z \times Y \times \{1\}$ can be identified with a closed 
subset of $X$, and is therefore paracompact. Consequently, assertion
$(i)$ follows from Lemma \ref{sadly}. Assertion $(iii)$ follows immediately from
$(i)$ and $(ii)$. It will therefore suffice to prove $(ii)$.

Since $Z$ can be identified with a closed subset of $X$, it is paracompact. Let
$\calF \in \Shv( Z \times Y)$, and let $U$ be an open $F_{\sigma}$ subset of
$Z$. We will show that $\gamma$ induces a homotopy equivalence
$(i^{\ast} k_{\ast} \pi^{\ast} \calF)(U) \rightarrow (\psi_{\ast} \calF)(U)$.
Shrinking $Z$ if necessary, we may suppose that $Z = U$. The right hand side
can be identified with $\calF(Z \times Y)$, while the left hand side is given
(by virtue of Corollary \toposref{snottle}) by the colimit
$\varinjlim_{V \in S}(\pi^{\ast} \calF)( k^{-1} V)$, where $V$ ranges over partially ordered set $S$ of
open subsets of $Z \times C(Y)$ which contain $Z$. By virtue of Remark \ref{snark}, we can replace
$S$ by the cofinal subset $S'$ consisting of open sets of the form $V = V_{f}$, where
$f: Z \rightarrow (0, \infty)$ is a continuous function (see Lemma \ref{cuppa}). Since $S'$ is filtered, it
will suffice to show that each of the maps $(\pi^{\ast} \calF)(k^{-1} V) \rightarrow \calF(Z \times Y)$
is an equivalence. Division by $f$ allows us to identify $(\pi^{\ast} \calF)(k^{-1} V)$ with
$(\pi^{\ast} \calF)( Z \times Y \times [0,1))$, and the desired result now follows from
Variant \ref{tsongv} on Lemma \ref{ander}.
\end{proof}

\begin{lemma}\label{olit}
Let $X$ be a paracompact space equipped with a conical $A$-stratification.
Then every point $x \in X_{a}$ admits a open $F_{\sigma}$ neighborhood $V$
which is homeomorphic (as an $A$-stratified space) to $Z \times C(Y)$, where
$Y$ is some $A_{>a}$-stratified space.
\end{lemma}

\begin{proof}
Since the stratification of $X$ is conical, there
exists an open neighborhood $U$ of $x$ which is homeomorphic
(as an $A$-stratified space) to $Z \times C(Y)$, where
$Y$ is some $A_{>a}$-stratified space. The open set $U$ need not be paracompact.
However, there exists a smaller open set $U' \subseteq U$ containing $x$ such that
$U'$ is an $F_{\sigma}$ subset of $X$, and therefore paracompact. let
$Z' = U' \cap Z$. Then $Z'$ is a closed subset of the paracompact space $U'$, and therefore
paracompact. Replacing $Z$ by $Z'$, we can assume that $Z$ is paracompact.
Applying Lemma \ref{cuppa}, we deduce that there exists a continuous function
$f: Z \rightarrow (0, \infty)$ such that $V_{f} \subseteq U$ (see Lemma \ref{cuppa} for an explanation of this notation). The set $V_{f}$ is the union
of the closures in $U'$ of the open sets $\{ V_{ \frac{n}{n+1} f} \}_{n > 0}$. It is therefore
an open $F_{\sigma}$ subset of $U'$ (and so also an $F_{\sigma}$ subset of the space $X$). We conclude by observing that $V_{f}$ is again homeomorphic to the product
$Z \times C(Y)$.
\end{proof}

\begin{remark}\label{jult}
If $A$ is a partially ordered set satisfying the ascending chain condition, then we can define
an ordinal-valued rank function $\rk$ on $A$. The function $\rk$ is uniquely determined by the following
requirement: for every element $a \in A$, the rank $\rk(a)$ is the smallest ordinal not of the form
$\rk(b)$, where $b > a$. More generally, suppose that $X$ is an $A$-stratified topological space. We define
the {\it rank} of $X$ to be the supremum of the set of ordinals $\{ \rk(a): X_{a} \neq \emptyset \}$.
\end{remark}

\begin{remark}\label{juline}
Let $X$ be a paracompact topological space of the form $Z \times C(Y)$. Then
$Z$ is paracompact (since it is homeomorphic to a closed subset of $X$). Suppose
that $X$ has singular shape. Since
the inclusion $Z \hookrightarrow X$ is a homotopy equivalence, we deduce also that
$Z$ has singular shape (Remark \ref{adlam}). The same argument shows that if $X$ is locally of singular shape, then $Z$ is locally of singular shape.
\end{remark}

\begin{proof}[Proof of Proposition \ref{toru}]
The assertion that $\theta: \calF \rightarrow \varprojlim \tau_{\leq n} \calF$ is an equivalence is local on $X$. It will therefore
suffice to prove that every point $x \in X_{a}$ admits an open $F_{\sigma}$ neighborhood $U$
such that $\theta$ is an equivalence over $U$. Since $A$ satisfies the ascending chain condition,
we may assume without loss of generality that the same result holds for every point
$x' \in X_{> a}$. Using Lemma \ref{olit}, we may assume without loss of generality that
$U$ is a paracompact open set of the form $Z \times C(Y)$, where $Y$ is some
$A_{>a}$-stratified space.

Let $i: Z \rightarrow Z \times C(Y)$ and $j: Z \times Y \times (0, \infty) \rightarrow Z \times C(Y)$
denote the inclusion maps. According to Lemma \ref{singbird}, it will suffice to verify the following:
\begin{itemize}
\item[$(a)$] The canonical map $i^{\ast} \calF \rightarrow \varprojlim i^{\ast} \tau_{\leq n} \calF
\simeq \varprojlim \tau_{\leq n} i^{\ast} \calF$ is an equivalence.
\item[$(b)$] The canonical map $j^{\ast} \calF \rightarrow \varprojlim j^{\ast} \tau_{\leq n} \calF
\simeq \varprojlim \tau_{\leq n} j^{\ast} \calF$ is an equivalence.
\item[$(c)$] The canonical map $i^{\ast} j_{\ast} j^{\ast} \calF
\rightarrow \varprojlim i^{\ast} j_{\ast} j^{\ast} \tau_{\leq n} \calF$ is an equivalence.
\end{itemize}

Assertion $(a)$ follows from Corollary \ref{tablor2} (note that $Z$ is locally of singular shape by
Remark \ref{juline}), and assertion $(b)$ follows from the inductive hypothesis. To prove $(c)$, let $\pi: Z \times Y \times (0, \infty)$ denote the projection.
Using the inductive hypothesis, we deduce that $j^{\ast} \calF$ is hypercomplete.
Since each fiber $\{ z \} \times \{y\} \times (0, \infty)$ is contained in a stratum of $X$,
we deduce that $j^{\ast} \calF$ is foliated, so that the counit map
$\pi^{\ast} \pi_{\ast} j^{\ast} \calF \rightarrow j^{\ast} \calF$ is an equivalence.
The same reasoning shows that $\pi^{\ast} \pi_{\ast} j^{\ast} \tau_{\leq n} \calF
\rightarrow j^{\ast} \tau_{\leq n} \calF$ is an equivalence for each $n \geq 0$.
Consequently, $(c)$ is equivalent to the assertion that the canonical map
$$ i^{\ast} j_{\ast} \pi^{\ast} \calG \rightarrow \varprojlim i^{\ast} j_{\ast} \pi^{\ast} \calG_{n}$$
is an equivalence, where $\calG = \pi_{\ast} j^{\ast} \calF$ and
$\calG_{n} = \pi_{\ast} j^{\ast} \tau_{\leq n} \calF$. Since the functor
$\pi_{\ast}$ preserves limits, the canonical map $\calG \rightarrow \varprojlim \calG_{n}$
is an equivalence by virtue of $(b)$. The desired result now follows from the fact
that the functor $i^{\ast} j_{\ast} \pi^{\ast}$ is equivalent to $\pi_{\ast}$, and therefore
preserves limits (Lemma \ref{torpus}).
\end{proof}

\begin{remark}
Let $X$ be a paracompact topological space equipped with a conical $A$-stratification, where
$A$ is a partially ordered set which satisfies the ascending chain condition. Suppose that
each stratum $X_{a}$ is locally of singular shape. Then $X$ is locally of singular shape. To prove this,
it suffices to show that $X$ has a covering by open sets which are locally of singular shape (Remark \ref{locsh}). Using Lemma \ref{olit}, we may reduce to the case where $X = Z \times C(Y)$, where
$Y$ is some $A_{> a}$-stratified space and $Z \times C(Y)$ is endowed with the induced
$A_{ \geq a}$-stratification. Working by induction on $a$, we may suppose that
$X - Z \simeq Z \times Y \times (0, \infty)$ is locally of singular shape. Let $U$ be an open $F_{\sigma}$ subset of $X$ and let
$U_0 = U \cap Z$. We wish to prove that $U$ is locally of singular shape. Using Lemma \ref{cuppa}, we deduce that there exists a continuous map $f: U_0 \rightarrow (0, \infty)$ such that $U$ contains the open set 
$V_{f} = U_0 \cup \{ (z,y,t) \in U_0 \times Y \times (0,\infty): t < f(z) \}$. Then
$U$ is covered by the open subsets $V_{f}$ and $U - U_0$. According to
Lemma \ref{kuil}, it suffices to show that $V_{f}$, $U - U_0$, and $V_{f} \cap (U - U_0)$
are of singular shape. The open sets $U - U_0$ and $V_{f} \cap (U-U_0)$ belong to
$X_{>a}$ and are therefore of singular shape by the inductive hypothesis. The open
set $V_{f}$ is homotopy equivalent to $U_0$, and thus has singular shape by virtue of our assumption
that $X_{a}$ is locally shapely (Remark \ref{adlam}).
\end{remark}

\subsection{$\infty$-Categories of Exit Paths}\label{lubos}

If $X$ is a sufficiently nice topological space, then Theorem \ref{squ} guarantees that
the $\infty$-category of locally constant sheaves on $X$ can be identified with
the $\infty$-category $\SSet_{/ \Sing(X) } \simeq \Fun( \Sing(X), \SSet)$. Roughly speaking,
we can interpret a sheaf $\calF$ on $X$ as a functor which assigns to each
$x \in X$ the stalk $\calF_x \in \SSet$, and to each path $p: [0,1] \rightarrow X$ joining
$x=p(0)$ to $y = p(1)$ the homotopy equivalence $\calF_{x} \simeq \calF_y$ given by transport along
$p$ (see \S \ref{hominv}).

Suppose now that $\calF$ is a sheaf on $X$ which is {\em not} locally constant. In this case, a
path $p: [0,1] \rightarrow X$ from $x = p(0)$ to $y= p(1)$ does not necessarily define a transport map $\calF_{x} \rightarrow \calF_{y}$. However, every point $\eta_0$ in the stalk $\calF_x$ can be lifted to a section of
$\calF$ over some neighborhood of $x$, which determines points $\eta_{t} \in \calF_{p(t) }$ for $t$ sufficiently small. If we assume that $p^{\ast} \calF$ is locally constant on the half-open interval $(0, 1]$, then
each $\eta_{t}$ can be transported to a point in the stalk $\calF_{y}$, and we should again expect
to obtain a well-defined map $\calF_{x} \rightarrow \calF_{y}$. For example, suppose that
$\calF$ is a sheaf which is locally constant when restricted to some closed subset $X_0 \subseteq X$,
and also when restricted to the open set $X - X_0$. In this case, the above analysis should apply
whenever $p^{-1} X_0 = \{0\}$: that is, whenever $p$ is a path which is {\it exiting} the closed subset
$X_0 \subseteq X$. Following an idea proposed by MacPherson, this suggests that we might try to identify $\calF$ with an $\SSet$-valued functor defined on some subset of the Kan complex $\Sing(X)$, which allows paths to travel from $X_0$ to
$X - X_0$ but not vice-versa. 

Our objective in this section is to introduce a simplicial subset $\Sing^{A}(X)$ associated to 
any stratification $f: X \rightarrow A$ of a topological space $X$ by a partially ordered set $A$. 
Our main result, Theorem \ref{kanna}, asserts that $\Sing^{A}(X)$ is an $\infty$-category provided
that the stratification of $X$ is conical (Definition \ref{sally}). In this case, we will refer to $\Sing^{A}(X)$ as the
{\it $\infty$-category of exit paths in $X$ with respect to the stratification $X \rightarrow A$}. In \S \ref{cloop}, we will show that (under suitable hypotheses) the $\infty$-category of $A$-constructible sheaves on $X$
is equivalent to the $\infty$-category of functors $\Fun( \Sing^{A}(X), \SSet)$.

\begin{remark}
The exit path $\infty$-category $\Sing^{A}(X)$ can be regarded as an $\infty$-categorical generalization of the $2$-category of exit paths constructed in \cite{treumann}. 
\end{remark}

\begin{definition}
Let $A$ be a partially ordered set, and let $X$ be a topological space equipped with
an $A$-stratification $f: X \rightarrow A$. We $\Sing^{A}(X) \subseteq \Sing(X)$ to be the simplicial
subset consisting of those $n$-simplices $\sigma: | \Delta^n | \rightarrow X$ which satisfy the following condition:
\begin{itemize}
\item[$(\ast)$] Let $| \Delta^n | = \{ ( t_0, \ldots, t_n) \in [0,1]^{n+1}: t_0 + \ldots + t_n = 1 \}$. Then
there exists a chain $a_0 \leq \ldots \leq a_n$ of elements of $A$ such that
for each point $( t_0, \ldots, t_i, 0, \ldots, 0) \in | \Delta^n |$ where $t_i \neq 0$, we have
$f( \sigma( t_0, \ldots, t_n) ) = a_i$.
\end{itemize}
\end{definition}

\begin{remark}
Let $A$ be a partially ordered set, regarded as a topological space as in
Definition \ref{cunning}. Then there is a natural map of simplicial sets
$\Nerve(A) \rightarrow \Sing(A)$, which carries an $n$-simplex
$(a_0 \leq \ldots \leq a_n)$ of $\Nerve(A)$ to the map
$\sigma: | \Delta^n | \rightarrow A$ characterized by the formula
$$ \sigma( t_0, \ldots, t_i, 0, \ldots, 0) = a_i $$
whenever $t_i > 0$. For any $A$-stratified topological space $X$, the simplicial 
set $\Sing^{A}(X)$ can be described as the fiber product $\Sing(X) \times_{ \Sing(A)} \Nerve(A)$.
In particular, there is a canonical map of simplicial sets $\Sing^{A}(X) \rightarrow \Nerve(A)$.
\end{remark}

We can now state our main result as follows:

\begin{theorem}\label{kanna}
Let $A$ be a partially ordered set, and let $X$ be a conically $A$-stratified topological space.
Then:
\begin{itemize}
\item[$(1)$] The projection $\Sing^{A}(X) \rightarrow \Nerve(A)$ is an inner fibration of simplicial sets.
\item[$(2)$] The simplicial set $\Sing^{A}(X)$ is an $\infty$-category.
\item[$(3)$] A morphism in $\Sing^{A}(X)$ is an equivalence if and only if its image in
$\Nerve(A)$ is degenerate (in other words, if and only if the underlying path
$[0,1] \rightarrow X$ belongs to a single stratum).
\end{itemize}
\end{theorem}

\begin{proof}
The implication $(1) \Rightarrow (2)$ is obvious. The ``only if'' direction of $(3)$ is clear
(since any equivalence in $\Sing^{A}(X)$ must project to an equivalence in $\Nerve(A)$), and
the ``if'' direction follows from the observation that each fiber $\Sing^{A}(X) \times_{ \Nerve(A)} \{a\}$
is isomorphic to the Kan complex $\Sing(X_{a})$. It will therefore suffice to prove $(1)$.
Fix $0 < i < n$; we wish to prove that every lifting problem of the form
$$ \xymatrix{ \Lambda^n_i \ar[r]^-{ \sigma_0} \ar[d] & \Sing^{A}(X) \ar[d] \\
\Delta^n \ar[r] \ar@{-->}[ur]^{\sigma} & \Nerve(A) }$$
admits a solution. 

The map $\Delta^n \rightarrow \Nerve(A)$ determines a chain of elements
$a_0 \leq a_1 \leq \ldots \leq a_n$. Without loss of generality, we may replace
$A$ by $A' = \{ a_0, \ldots, a_n \}$ and $X$ by $X \times_{A} A'$. We may therefore 
assume that $A$ is a finite nonempty linearly ordered set. We now work by induction on the
number of elements of $A$. If $A$ has only a single element, then $\Sing^{A}(X) = \Sing(X)$ is
a Kan complex and there is nothing to prove. Otherwise, there exists some integer $p < n$ such
that $a_p = a_0$ and $a_{p+1} \neq a_0$. There are two cases to consider.

\begin{itemize}
\item[$(a)$] Suppose that $p < i < n$. Let $q = n - p - 1$ and let $j = i - p - 1$, so that we have
isomorphisms of simplicial sets 
$$\Delta^n \simeq \Delta^p \star \Delta^q \quad \quad \Lambda^n_{i} \simeq
( \Delta^p \star \Lambda^{q}_{q'}) \coprod_{ \bd \Delta^p \star \Lambda^{q}_{j}} ( \bd \Delta^p \star \Delta^q).$$
We will use the first isomorphism to identify $| \Delta^n |$ with the pushout
$$ | \Delta^p | \coprod_{ | \Delta^p | \times | \Delta^q| \times \{0\} }
( | \Delta^p| \times | \Delta^q | \times [0,1]) \coprod_{ | \Delta^p | \times | \Delta^q| \times \{1\} } | \Delta^q |.$$
Let $K \subseteq | \Delta^p | \times | \Delta^q |$ be the union of the closed subsets
$| \bd \Delta^p| \times | \Delta^q |$ and $| \Delta^p | \times | \Lambda^q_j |$, so that
$| \Lambda^n_i |$ can be identified with the pushout
$$ | \Delta^p | \coprod_{ | \Delta^p | \times | \Delta^q| \times \{0\} }
( K \times [0,1]) \coprod_{ | \Delta^p | \times | \Delta^q| \times \{1\} } | \Delta^q |.$$
Let $K' \subseteq | \Delta^p | \times | \Delta^q| \times [0,1]$ be the union of
$K \times [0,1]$ with $| \Delta^p | \times |\Delta^q| \times \{0,1\}$. 
Then $\sigma_0$ determines a continuous map
$F_0: K' \rightarrow X$. To construct the map $\sigma$, we must extend $F_0$ to a map
$F: | \Delta^p | \times | \Delta^q| \times [0,1] \rightarrow X$ satisfying the following
condition: for every point $s \in (| \Delta^p | \times | \Delta^q | \times [0,1]) - K'$, we
have $F(s) \in X_{ a_n}$. 

Let $F_{-}: | \Delta^p | \rightarrow X_{a_0}$ be the map obtained by restricting
$F_0$ to $| \Delta^p | \times | \Delta^q | \times \{0\}$. For every point $x \in X_{a_0}$, choose an open neighborhood $U_{x} \subseteq X$ as in Definition \ref{sally}. Choose a triangulation
of the simplex $| \Delta^p |$ with the following property: for every simplex $\tau$ of the triangulation,
the image $F_{-}( \tau)$ is contained in some $U_{x}$. Refining our triangulation if necessary, we may assume that $| \bd \Delta^p |$ is a subcomplex of $| \Delta^p |$. For every subcomplex $L$
of $| \Delta^p |$ which contains $| \bd \Delta^p |$, we let $K_{L} \subseteq | \Delta^p | \times | \Delta^q|$
denote the union of the closed subsets $L \times | \Delta^q |$ and $| \Delta^p | \times | \Lambda^q_{j} |$
and $K'_{L} \subseteq | \Delta^p | \times | \Delta^q| \times [0,1]$ denote the union of the closed
subsets $K_L \times [0,1]$ and $| \Delta^p | \times | \Delta^q | \times \{0,1\}$. We will show that
$F_0$ can be extended to a continuous map $F_{L}: K'_{L} \rightarrow X$ (satisfying
the condition that $F_{L}(s) \in X_{ a_n}$ for $s \notin K'$), using induction on the number of simplices
of $L$. If $L = | \bd \Delta^p |$, there is nothing to prove. Otherwise, we may assume without loss
of generality that $L = L_0 \cup \tau$, where $L_0$ is another subcomplex of
$| \Delta^p |$ containing $| \bd \Delta^p|$ and $\tau$ is a simplex of $L$ such that $\tau \cap L_0 = \bd \tau$.
The inductive hypothesis guarantees the existence of a map $F_{L_0}: 
K'_{L_0} \rightarrow X$ with the desired properties. 

Let $K_{\tau} \subseteq \tau \times | \Delta^q |$ be the union of the closed subsets
$\bd \tau \times | \Delta^q |$ and $\tau \times | \Lambda^q_{j} |$, and let
$K'_{\tau} \subseteq \tau \times | \Delta^q | \times [0,1]$ be the union of the closed
subsets $K_{\tau} \times [0,1]$ and $\tau \times | \Delta^q| \times \{ 0,1\}$. 
The map $F_{L_0}$ restricts to a map $G_0: K'_{\tau} \rightarrow X$. To
construct $F_{L}$, it will suffice to extend $G_0$ to a continuous map
$G: \tau \times | \Delta^q | \times [0,1] \rightarrow X$ (satisfying the condition that $G(s) \in X_{a_n}$ for
$s \notin K'_{\tau}$). 

By assumption, the map $G_0$ carries $\tau \times | \Delta^q| \times \{0\}$ into an open subset
$U_{x}$, for some $x \in X_{a_0}$. Let $U = U_x$, and choose a homeomorphism
$U \simeq Z \times C(Y)$, where $Y$ is an $A_{> a_0}$-stratified space. Since
$\tau \times | \Delta^q |$ is compact, we deduce that $G_0( \tau \times | \Delta^q| \times [0,r])
\subseteq U$ for some real number $0 < r < 1$. Let $X' = X - X_{ a_0}$ and let
$A' = A - \{ a_0 \}$, so that $X'$ is an $A'$-stratified space. Let $m$ be
the dimension of the simplex $\tau$. The restriction
$G_0 | ( \tau \times | \Delta^q | \times \{1\} )$ determines a map of simplicial sets
$\overline{h}_1: \Delta^m \times \Delta^q \rightarrow \Sing^{A'}(X')$. Let
$J$ denote the simplicial set $( \bd \Delta^m \times \Delta^q)
\coprod_{ \bd \Delta^m \times \Lambda^q_{j}} ( \Delta^m \times \Lambda^q_{j})$.
The restriction of
$G_0$ to $K_{\tau} \times [r,1]$ determines another map of simplicial sets
$h: J \times \Delta^1 \rightarrow \Sing^{A'}(X')$, which is a natural transformation from
$h_0 = h | (J \times \{0\})$ to $h_1 = h | (J \times \{1\}) = \overline{h}_1 | J$. 
It follows from the inductive hypothesis that $\Sing^{A'}(X')$ is an $\infty$-category,
and (using $(3)$) that natural transformation $h$ is an equivalence. Consequently, we can
lift $h$ to an equivalence $\overline{h}: \overline{h}_0 \rightarrow \overline{h}_1$ in
$\Fun( J, \Sing^{A'}(X'))$. This morphism determines a continuous map
$G_{+}: \tau \times | \Delta^q| \times [r, 1] \rightarrow X$ which agrees with
$G_0$ on $( \tau \times | \Delta^q | \times [r,1]) \cap K'_{\tau}$. 

Let us identify $| \Delta^q |$ with the set of tuples of real numbers
$\vec{t} = (t_0, t_1, \ldots, t_q)$ such that $0 \leq t_k \leq 1$ and
$t_0 + \cdots + t_q = 1$. In this case, we let $d( \vec{t} ) = \inf \{ t_{k} : k \neq j \}$:
note that $d( \vec{t} ) = 0$ if and only if $\vec{t} \in | \Lambda^q_{j} |$. 
If $u$ is a real number satisfying $0 \leq u \leq d( \vec{t})$, we let
$\vec{t}_{u}$ denote the tuple
$$( t_0 - u, t_1 - u, \ldots, t_{j-1} - u, t_{j} + q u, t_{j+1-u} - u, \ldots, t_{q} - u) \in | \Delta^q |. $$
Choose a continuous function $d': \tau \rightarrow [0,1]$ which vanishes on
$\bd \tau$ and is positive on the interior of $\R$. For every positive real number
$v$, let $c_{v}: \tau \times | \Delta^q | \times [r,1] \rightarrow \tau \times | \Delta^q | \times [r,1]$
given by the formula
$$ c_{v}( s, \vec{t}, r') = (s, \vec{t}_{ d(\vec{t})} \frac{ v d'(s) (1-r')}{ 1 + v d'(s) (1-r')}, r'),$$ 
and let $G_{+}^{v}$ denote the composition $G_{+} \circ c_{v}$. 
Since $G_{+}$ agrees with $G_0$ on $K_{\tau} \times \{r\}$, it carries
$K_{\tau} \times \{r\}$ into $U$. By continuity, there exists a neighborhood
$V$ of $K_{\tau}$ in $\tau \times | \Delta^q |$ such that $G_{+}(V \times \{r\}) \subseteq U$.
If the real number $v$ is sufficiently large, then $c_{v}( \tau \times | \Delta^q | \times \{r\})
\subseteq V$, so that $G_{+}^{v}( \tau \times | \Delta^q | \times \{r \}) \subseteq U$.
Replacing $G_{+}$ by $G_{+}^{v}$, we may assume that
$G_{+}( \tau \times | \Delta^q| \times \{r \}) \subseteq U$ (here we invoke the assumption
that $j < q$ to guarantee that $G_{+}$ continues to satisfy the requirement that
$G_{+}(s, \vec{t}, r') \in X_{a_{n}}$ whenever $\vec{t} \notin | \Lambda^q_{j} |$).

Let $X'' = U - X_{a_0} \simeq Z \times Y \times \R_{> 0}$. The $A'$-stratification of
$X'$ restricts to a (conical) $A'$-stratification of $X''$.
Let $g: \tau \times | \Delta^q | \times \{r\} \rightarrow X''$ be the map
obtained by restricting $G_{+}$. Then $g$ determines a map of simplicial sets
$\overline{\phi}_0: \Delta^m \times \Delta^q \rightarrow \Sing^{A'}(X'')$. Let $I$ denote the simplicial set
$$ \Delta^{ \{0, 1\} } \coprod_{ \{ 1\} } \Delta^{ \{1,2\} } \coprod_{ \{2\} } \Delta^{ \{2,3\} } \coprod_{ \{3\} } \ldots,$$
and identify the geometric realization $|I|$ with the open interval $(0, r]$. Then $G_0$
determines a map of simplicial sets $J \times I \rightarrow \Sing^{A'}(X'')$, which we can identify
with a sequence of maps $\phi_0, \phi_1, \ldots \in \Fun( J, \Sing^{A'}(X''))$ together with
natural transformations $\phi_0 \rightarrow \phi_1 \rightarrow \ldots$. We note that
$\phi_0 = \overline{\phi}_0 | J$. The inductive hypothesis guarantees that
$\Sing^{A'}(X'')$ is an $\infty$-category, and assertion $(3)$ ensures that each of the natural transformations $\phi_{k} \rightarrow \phi_{k+1}$ is an equivalence. It follows that we can
lift these natural transformations to obtain a sequence of equivalences
$$ \overline{\phi}_0 \rightarrow \overline{\phi}_1 \rightarrow \overline{\phi_2} \rightarrow \cdots$$
in the $\infty$-category $\Fun( \Delta^m \times \Delta^q, \Sing^{A'}(X''))$. This sequence of equivalences is given by a map of simplicial sets $\Delta^m \times \Delta^q \times I \rightarrow \Sing^{A'}(X'')$, which
we can identify with a continuous map $\tau \times | \Delta^q | \times (0, r] \rightarrow 
Z \times Y \times \R_{> 0}$. Let $y: \tau \times | \Delta^q | \times (0,r] \rightarrow Y$ be the projection
of this map onto the second fiber.

We observe that $G_{+}$ and $G_0$ together determine a map 
$(K_{\tau} \times [0,r]) \coprod_{ K_{\tau} \times \{0,r\}} ( \tau \times | \Delta^q | \times \{0,r\} )
\rightarrow X'$. Let $z$ denote the composition of this map with the projection $U \rightarrow Z \times \R_{\geq 0}$. Since the domain of $z$ is a retract of $\tau \times | \Delta^q | \times [0,r]$, we can extend $z$ to a continuous map $\overline{z}: \tau \times | \Delta^q | \times [0,r] \rightarrow Z \times \R_{\geq 0}$.
Let $\overline{z}_1: \tau \times | \Delta^q | \times [0,r] \rightarrow \R_{\geq 0}$ be obtained from
$\overline{z}$ by projection onto the second factor. By adding to 
$\overline{z}_1$ a function which vanishes on 
$(K_{\tau} \times [0,r]) \coprod_{ K_{\tau} \times \{0,r\}} ( \tau \times | \Delta^q | \times \{0,r\} )$ and is
positive elsewhere, we can assume that $\overline{z}_1^{-1} \{0\} = \tau \times | \Delta^q | \times \{0\}$.
Let $G_{-}: \tau \times | \Delta^q| \times [0,r] \rightarrow U \simeq Z \times C(Y)$ be the map
which is given by $\overline{z}$ on $\tau \times | \Delta^q | \times \{0\}$ and
by the pair $(\overline{z}, y)$ on $\tau \times | \Delta^q | \times (0,r]$. Then
$G_{-}$ and $G_{+}$ together determine an extension
$G: \tau \times | \Delta^q | \times [0,1] \rightarrow X$ of $G_0$ with the desired properties.

\item[$(b)$] Suppose now that $0 < i \leq p$. The proof proceeds as in case $(a)$ with some minor changes. We let $q = n - p - 1$ as before, so that we have an identification
of $| \Delta^n |$ with the pushout
$$ | \Delta^p | \coprod_{ | \Delta^p | \times | \Delta^q| \times \{0\} }
( | \Delta^p| \times | \Delta^q | \times [0,1]) \coprod_{ | \Delta^p | \times | \Delta^q| \times \{1\} } | \Delta^q |.$$
Let $K \subseteq | \Delta^p | \times | \Delta^q |$ be the union of the closed subsets
$| \Lambda^p_i | \times | \Delta^q |$ and $| \Delta^p | \times | \bd \Delta^q |$, so that
$| \Lambda^n_i |$ can be identified with the pushout
$$ | \Delta^p | \coprod_{ | \Delta^p | \times | \Delta^q| \times \{0\} }
( K \times [0,1]) \coprod_{ | \Delta^p | \times | \Delta^q| \times \{1\} } | \Delta^q |.$$
Let $K' \subseteq | \Delta^p | \times | \Delta^q| \times [0,1]$ be the union of
$K \times [0,1]$ with $| \Delta^p | \times |\Delta^q| \times \{0,1\}$. 
Then $\sigma_0$ determines a continuous map
$F_0: K' \rightarrow X$. To construct the map $\sigma$, we must extend $F_0$ to a map
$F: | \Delta^p | \times | \Delta^q| \times [0,1] \rightarrow X$ satisfying the following
condition: for every point $s \in (| \Delta^p | \times | \Delta^q | \times [0,1]) - K'$, we
have $F(s) \in X_{ a_n}$. 

We observe that there is a homeomorphism of $| \Delta^p |$ with
$| \Delta^{p-1} | \times [0,1]$ which carries $| \Lambda^{p}_{i} |$ to
$| \Delta^{p-1} | \times \{0\}$. 
Let $F_{-}: | \Delta^{p-1}| \times [0,1] \rightarrow X_{a_0}$ be the map determined by
$\sigma_0$ together with this homeomorphism.
For every point $x \in X_{a_0}$, choose an open neighborhood $U_{x} \subseteq X$ as in Definition \ref{sally}. Choose a triangulation of the simplex $| \Delta^{p-1}|$ and a large positive integer
$N$ so that the following condition is satisfied: for every simplex $\tau$ of
$| \Delta^{p-1}|$ and every nonnegative integer $k < N$, the map
$F_{-}$ carries $\tau \times [ \frac{k}{N}, \frac{k+1}{N}]$ into some
$U_{x}$. For every subcomplex $L$ of $| \Delta^{p-1} |$, we let
$K_{L} \subseteq | \Delta^p | \times | \Delta^q|$
denote the union of the closed subsets $L \times [0,1] \times | \Delta^q |$, $| \Delta^{p-1} | \times \{0\} \times | \Delta^q |$, and $| \Delta^{p-1} | \times [0,1] \times | \bd \Delta^q |$.
Let $K'_{L} \subseteq | \Delta^p | \times | \Delta^q| \times [0,1]$ denote the union of the closed
subsets $K_L \times [0,1]$ and $| \Delta^p | \times | \Delta^q | \times \{0,1\}$. We will show that
$F_0$ can be extended to a continuous map $F_{L}: K'_{L} \rightarrow X$ (satisfying
the condition that $F_{L}(s) \in X_{ a_n}$ for $s \notin K'$), using induction on the number of simplices
of $L$. If $L$ is empty there is nothing to prove. Otherwise, we may assume without loss
of generality that $L = L_0 \cup \tau$, where $\tau$ is a simplex
of $| \Delta^{p-1} |$ such that $\tau \cap L_0 = \bd \tau$.
The inductive hypothesis guarantees the existence of a map $F_{L_0}: 
K'_{L_0} \rightarrow X$ with the desired properties. 

For $0 \leq k \leq N$, let $K_{\tau,k} \subseteq \tau \times [0,1] \times | \Delta^q |$ be the union of the closed subsets $\bd \tau \times [0, 1] \times | \Delta^q |$, 
$\tau \times [0, \frac{k}{N}] \times | \Delta^q|$, and $\tau \times [0,1] \times | \bd \Delta^q |$. 
Let $K'_{\tau,k} \subseteq \tau \times [0,1] \times | \Delta^q | \times [0,1]$ be the union of the closed
subsets $K_{\tau,k} \times [0,1]$ and $\tau \times | \Delta^q| \times \{ 0,1\}$. 
The map $F_{L_0}$ restricts to a map $F[0]: K'_{\tau,0} \rightarrow X$. To
construct $F_{L}$, it will suffice to extend $G_0$ to a continuous map
$F[N]: K_{\tau,N} \times [0,1] \rightarrow X$ (satisfying the condition that $F[n](s) \in X_{a_n}$ for
$s \notin K'_{\tau}$). We again proceed by induction, constructing maps
$F[k]: K'_{\tau, k} \rightarrow X$ for $k \leq N$ using recursion on $k$. Assume that
$k > 0$ and that $F[k-1]$ has already been constructed. 

Let $\overline{\tau}$ denote the prism
$\tau \times [ \frac{k-1}{N}, \frac{k}{N}]$, and let $\overline{\tau}_0$ denote the closed subset
of $\overline{\tau}$ which is the union of $\bd \tau \times [ \frac{k-1}{N}, \frac{k}{N}]$ with
$\tau \times \{ \frac{k-1}{N} \}$. Let $K_{\overline{\tau}} \subseteq \overline{\tau} \times | \Delta^q |$
denote the union of the closed subsets $\overline{\tau} \times | \bd \Delta^q |$ and
$\overline{\tau}_0 \times | \Delta^q |$. Let $K'_{ \overline{\tau} } \subseteq
\overline{\tau} \times | \Delta^q | \times [0,1]$ be the union of the closed subsets
$K_{\overline{\tau}} \times [0,1]$ with $\overline{\tau} \times | \Delta^q| \times \{0,1\}$. 
Then $F[k-1]$ determines a map $G_0: K'_{\overline{\tau}} \rightarrow X$. To find the desired extension $F[k]$ of $F[k-1]$, it will suffice to prove that $G_0$ admits a continuous extension
$G: \overline{\tau} \times | \Delta^q| \times [0,1]$ (again satisfying the condition that
$G(s) \in X_{ a_n}$ whenever $s \notin K'_{\overline{\tau}}$). 

By assumption, the map $G_0$ carries $\overline{\tau} \times | \Delta^q| \times \{0\}$ into an open subset
$U_{x}$, for some $x \in X_{a_0}$. Let $U = U_x$, and choose a homeomorphism
$U \simeq Z \times C(Y)$, where $Y$ is an $A_{> a_0}$-stratified space. Since
$\tau \times | \Delta^q |$ is compact, we deduce that $G_0( \overline{\tau} \times | \Delta^q| \times [0,r])
\subseteq U$ for some real number $0 < r < 1$. Let $X' = X - X_{ a_0}$ and let
$A' = A - \{ a_0 \}$, so that $X'$ is an $A'$-stratified space. Let $m$ be
the dimension of the simplex $\tau$. The restriction
$G_0 | ( \overline{\tau} \times | \Delta^q | \times \{1\} )$ determines a map of simplicial sets
$\overline{h}_1: \Delta^m \times \Delta^1 \times \Delta^q \rightarrow \Sing^{A'}(X')$. Let
$J$ denote the simplicial subset of $\Delta^m \times \Delta^1 \times \Delta^q$ spanned
by $\Delta^m \times \{0\} \times \Delta^q$, $\Delta^m \times \Delta^1 \times \bd \Delta^q$, and
$\bd \Delta^m \times \Delta^1 \times \Delta^q$. 
The restriction of $G_0$ to $K_{\tau} \times [r,1]$ determines another map of simplicial sets
$h: J \times \Delta^1 \rightarrow \Sing^{A'}(X')$, which is a natural transformation from
$h_0 = h | (J \times \{0\})$ to $h_1 = h | (J \times \{1\}) = \overline{h}_1 | J$. 
It follows from the inductive hypothesis that $\Sing^{A'}(X')$ is an $\infty$-category,
and (using $(3)$) that natural transformation $h$ is an equivalence. Consequently, we can
lift $h$ to an equivalence $\overline{h}: \overline{h}_0 \rightarrow \overline{h}_1$ in
$\Fun( J, \Sing^{A'}(X'))$. This morphism determines a continuous map
$G_{+}: \overline{\tau} \times | \Delta^q| \times [r, 1] \rightarrow X$ which agrees with
$G_0$ on $( \overline{\tau} \times | \Delta^q | \times [r,1]) \cap K'_{\tau}$. 

Let $d: | \Delta^q | \rightarrow [0,1]$ be a continuous function which vanishes
precisely on $| \bd \Delta^q |$, and choose $d': \tau \rightarrow [0,1]$ similarly.
For every nonnegative real number $v$, let $c_{v}$ be the map from
$\tau \times [ \frac{k-1}{N}, \frac{k}{N}] \times | \Delta^q | \times [r,1] \rightarrow \tau \times | \Delta^q | \times [r,1]$ to itself which is given by the formula
$$ c_{v}( x, \frac{k-1}{N} + t, y, r') = (x, \frac{k-1}{N} + \frac{t}{ 1 + v d'(x) d(y) (1- r')}, y, r')$$
and let $G_{+}^{v}$ denote the composition $G_{+} \circ c_{v}$. 
Since $G_{+}$ agrees with $G_0$ on $K_{\overline{\tau}} \times \{r\}$, it carries
$K_{\tau} \times \{r\}$ into $U$. By continuity, there exists a neighborhood
$V$ of $K_{\overline{\tau}}$ in $\overline{\tau} \times | \Delta^q |$ such that 
$G_{+}(V \times \{r\}) \subseteq U$.
If the real number $v$ is sufficiently large, then $c_{v}( \overline{\tau} \times | \Delta^q | \times \{r\})
\subseteq V$, so that $G_{+}^{v}( \overline{\tau} \times | \Delta^q | \times \{r \}) \subseteq U$.
Replacing $G_{+}$ by $G_{+}^{v}$, we may assume that
$G_{+}( \overline{\tau} \times | \Delta^q| \times \{r \}) \subseteq U$.

Let $X'' = U - X_{a_0} \simeq Z \times Y \times \R_{> 0}$. The $A'$-stratification of
$X'$ restricts to a (conical) $A'$-stratification of $X''$.
Let $g: \overline{\tau} \times | \Delta^q | \times \{r\} \rightarrow X''$ be the map
obtained by restricting $G_{+}$. Then $g$ determines a map of simplicial sets
$\overline{\phi}_0: \Delta^m \times \Delta^1 \times \Delta^q \rightarrow \Sing^{A'}(X'')$. Let $I$ denote the simplicial set
$$ \Delta^{ \{0, 1\} } \coprod_{ \{ 1\} } \Delta^{ \{1,2\} } \coprod_{ \{2\} } \Delta^{ \{2,3\} } \coprod_{ \{3\} } \ldots,$$
and identify the geometric realization $|I|$ with the open interval $(0, r]$. Then $G_0$
determines a map of simplicial sets $J \times I \rightarrow \Sing^{A'}(X'')$, which we can identify
with a sequence of maps $\phi_0, \phi_1, \ldots \in \Fun( J, \Sing^{A'}(X''))$ together with
natural transformations $\phi_0 \rightarrow \phi_1 \rightarrow \ldots$. We note that
$\phi_0 = \overline{\phi}_0 | J$. The inductive hypothesis guarantees that
$\Sing^{A'}(X'')$ is an $\infty$-category, and assertion $(3)$ ensures that each of the natural transformations $\phi_{k} \rightarrow \phi_{k+1}$ is an equivalence. It follows that we can
lift these natural transformations to obtain a sequence of equivalences
$$ \overline{\phi}_0 \rightarrow \overline{\phi}_1 \rightarrow \overline{\phi_2} \rightarrow \cdots$$
in the $\infty$-category $\Fun( \Delta^m \times \Delta^1 \times \Delta^q, \Sing^{A'}(X''))$. This sequence of equivalences is given by a map of simplicial sets $\Delta^m \times \Delta^1 \times \Delta^q \times I \rightarrow \Sing^{A'}(X'')$, which
we can identify with a continuous map $\overline{\tau} \times | \Delta^q | \times (0, r] \rightarrow 
Z \times Y \times \R_{> 0}$. Let $y: \overline{\tau} \times | \Delta^q | \times (0,r] \rightarrow Y$ be the projection of this map onto the second fiber.

We observe that $G_{+}$ and $G_0$ together determine a map 
$(K_{\overline{\tau}} \times [0,r]) \coprod_{ \overline{K_{\tau}} \times \{0,r\}} ( \overline{\tau} \times | \Delta^q | \times \{0,r\} )
\rightarrow X'$. Let $z$ denote the composition of this map with the projection $U \rightarrow Z \times \R_{\geq 0}$. Since the domain of $z$ is a retract of $\overline{\tau} \times | \Delta^q | \times [0,r]$, we can extend $z$ to a continuous map $\overline{z}: \overline{\tau} \times | \Delta^q | \times [0,r] \rightarrow Z \times \R_{\geq 0}$.
Let $\overline{z}_1: \overline{\tau} \times | \Delta^q | \times [0,r] \rightarrow \R_{\geq 0}$ be obtained from
$\overline{z}$ by projection onto the second factor. By adding to 
$\overline{z}_1$ a function which vanishes on 
$(K_{\overline{\tau}} \times [0,r]) \coprod_{ K_{\overline{\tau}} \times \{0,r\}} ( \overline{\tau} \times | \Delta^q | \times \{0,r\} )$ and is
positive elsewhere, we can assume that $\overline{z}_1^{-1} \{0\} = \overline{\tau} \times | \Delta^q | \times \{0\}$.
Let $G_{-}: \overline{\tau} \times | \Delta^q| \times [0,r] \rightarrow U \simeq Z \times C(Y)$ be the map
which is given by $\overline{z}$ on $\overline{\tau} \times | \Delta^q | \times \{0\}$ and
by the pair $(\overline{z}, y)$ on $\overline{\tau} \times | \Delta^q | \times (0,r]$. Then
$G_{-}$ and $G_{+}$ together determine an extension
$G: \overline{\tau} \times | \Delta^q | \times [0,1] \rightarrow X$ of $G_0$ with the desired properties.
\end{itemize}
\end{proof}

\subsection{Exit Paths in a Simplicial Complex}\label{sec5sub7}

Every simplicial complex $X$ admits a canonical stratification, whose strata are the interiors of the simplices of $X$. In this section, we will show that the $\infty$-category of exit paths associated to a stratified space of this type is particularly simple: it is equivalent to the partially ordered set of simplices of $X$ (Theorem \ref{ikos}). We begin by reviewing some definitions.

\begin{definition}
An {\it abstract simplicial complex} consists of the following data:
\begin{itemize}
\item[$(1)$] A set $V$ (the set of {\it vertices} of the complex).
\item[$(2)$] A collection $S$ of nonempty finite subsets of $V$ satisfying the following condition:
\begin{itemize}
\item[$(\ast)$] If $\emptyset \neq \sigma \subseteq \sigma' \subseteq V$ and $\sigma' \in S$, then $\sigma \in S$.
\end{itemize}
\end{itemize}
We will say that $(V,S)$ is {\it locally finite} if each element $\sigma \in S$ is contained in only finitely
many other elements of $S$.

Let $(V,S)$ be an abstract simplicial complex, and choose a linear ordering on $V$.
We let $\Delta^{(V,S)}$ denote the simplicial subset of $\Delta^{V}$ spanned by those
simplices $\sigma$ of $\Delta^{V}$ such that the set of vertices of $\sigma$ belongs to $S$.
Let $| \Delta^{(V,S)} |$ denote the geometric realization of $\Delta^{(V,S)}$. This topological
space is independent of the choice of linear ordering on $S$, up to canonical homeomorphism.
As a set, $| \Delta^{V,S} |$ can be identified with the collection of maps
$w: V \rightarrow [0,1]$ such that $\supp(w) = \{ v \in V: w(v) \neq 0 \} \in S$ and 
$\sum_{ v \in V} w(v) = 1$.
\end{definition}

\begin{definition}
Let $(V,S)$ be an abstract simplicial complex. We regard $S$ as a partially ordered
set with respect to inclusions. Then $| \Delta^{(V,S)} |$ is equipped with a natural
$S$-stratification, given by the map 
$$(t \in | \Delta^{(V,S)} |) \mapsto ( \supp(t) \in S).$$
\end{definition}

\begin{proposition}\label{silkman}
Let $(V,S)$ be a locally finite abstract simplicial complex. Then the
$S$-stratification of $| \Delta^{(V,S)} |$ is conical.
\end{proposition}

\begin{proof}
Consider an arbitrary
$\sigma \in S$. Let $V' = V - \sigma$, and let
$S' = \{ \sigma' - \sigma: \sigma \subset \sigma' \in S \}$. Then $(V',S')$ is another
abstract simplicial complex. Let $Z = | \Delta^{(V,S)} |_{\sigma}$ and let
$Y = | \Delta^{(V',S') }|$. Then the inclusion $Z \hookrightarrow | \Delta^{(V,S)} |$ extends
to an open embedding $h: Z \times C(Y) \rightarrow | \Delta^{(V,S)} |$, which is given on
$Z \times Y \times (0, \infty)$ by the formula
$$ h( w_{Z}, w_{Y}, t)(v) = \begin{cases} \frac{w_{Z}(v)}{t+1} & \text{ if } v \in \sigma \\
\frac{t w_{Y}(v)}{t+1} & \text{ if } v \notin \sigma \end{cases}$$
If $(V,S)$ is locally finite, then $h$ is an open embedding whose image is
$| \Delta^{(V,S)} |_{ > \sigma}$, which proves that the $S$-stratification of
$| \Delta^{(V,S)} |$ is locally conical.
\end{proof}

\begin{corollary}
Let $(V,S)$ be an abstract simplicial complex. Then the simplicial set
$\Sing^{S} | \Delta^{(V,S)} |$ is an $\infty$-category.
\end{corollary}

\begin{proof}
For every subset $V_0 \subseteq V$, let $S_0 = \{ \sigma \in S: \sigma \subseteq V_0 \}$.
Then $\Sing^{S} | \Delta^{(V,S)} |$ is equivalent to the filtered colimit
$\varinjlim_{V_0} \Sing^{S_0} | \Delta^{ (V_0, S_0) } |$, where the colimit is taken over all
finite subsets $V_0 \subseteq V$. It will therefore suffice to prove that each
$\Sing^{S_0} | \Delta^{ (V_0, S_0) } |$ is an $\infty$-category. Replacing $V$ by $V_0$, we may assume that $V$ is finite so that $(V,S)$ is locally finite. In this case, the desired result follows immediately
from Proposition \ref{silkman} and Theorem \ref{kanna}.
\end{proof}

\begin{theorem}\label{ikos}
Let $(V,S)$ be an abstract simplicial complex. Then the projection
$q: \Sing^{S} | \Delta^{(V,S)} | \rightarrow \Nerve(S)$ is an equivalence of $\infty$-categories.
\end{theorem}

\begin{proof}
Since each stratum of $| \Delta^{(V,S)} |$ is nonempty, the map $q$ is essentially surjective.
To prove that $q$ is fully faithful, fix points $x \in | \Delta^{(V,S)} |_{\sigma}$
and $y \in | \Delta^{(V,S)} |_{\sigma'}$. It is clear that
$M = \bHom_{ \Sing^{S} | \Delta^{(V,S)} |}( x,y)$ is empty unless $\sigma \subseteq \sigma'$.
We wish to prove that $M$ is contractible if $\sigma \subseteq \sigma'$. We can
identify $M$ with $\Sing P$, where $P$ is the space of paths
$p: [0,1] \rightarrow | \Delta^{(V,S)}|$ such that $p(0) = x$, $p(1) = y$, and
$p(t) \in | \Delta{(V,S)} |_{\sigma'}$ for $t > 0$. It now suffices to observe that
there is a contracting homotopy $h: P \times [0,1] \rightarrow P$, given by the formula
$$ h( p, s)(t) = (1-s) p(t) + s(1-t) x + st y.$$
\end{proof}

\begin{remark}
Let $(V,S)$ be an abstract simplicial complex. It is possible to construct an explicit
homotopy inverse to the equivalence of $\infty$-categories
$q: \Sing^{S} | \Delta^{(V,S)} | \rightarrow \Nerve(S)$ of Theorem \ref{ikos}. 
For each $\sigma \in S$ having cardinality
$n$, we let $w_{\sigma} \in | \Delta^{(V,S)} |$ be the point described by the formula
$$ w_{\sigma}(v) = \begin{cases} \frac{1}{n} & \text{ if } v \in \sigma \\
0 & \text{ if } v \notin \sigma. \end{cases}$$
For every chain of subsets $\emptyset \neq \sigma_0  \subseteq \sigma_1 \subseteq \ldots \subseteq \sigma_k \in S,$ we define a map $| \Delta^{k} | \rightarrow | \Delta^{ (V,S) } |$ by the formula
$$ (t_0, \ldots, t_k) \mapsto t_0 w_{ \sigma_0} + \cdots + t_k w_{ \sigma_{k} }.$$
This construction determines section $\phi: \Nerve(S) \rightarrow \Sing^{S} | \Delta^{(V,S)} |$ of
$q$, and is therefore an equivalence of $\infty$-categories.
The induced map of topological spaces $| \Nerve(S) | \rightarrow | \Delta^{(V,S)} |$ is a homeomorphism:
it is given by the classical process of {\em barycentric subdivision} of the simplicial complex
$| \Delta^{(V,S)} |$. 
\end{remark}

\subsection{A Seifert-van Kampen Theorem for Exit Paths}\label{sec5sub8}

Our goal in this section is to prove the following generalization of Theorem \ref{vankamp}:

\begin{theorem}\label{superkamp}
Let $A$ be a partially ordered set, let $X$ be an $A$-stratified topological space,
and let $\calC$ be a category equipped with a functor $U: \calC \rightarrow \calU(X)$, where
$\calU(X)$ denotes the partially ordered set of all open subsets of $X$.
Assume that the following conditions are satisfied:
\begin{itemize}
\item[$(i)$] The $A$-stratification of $X$ is conical.
\item[$(ii)$] For every point $x \in X$, the full subcategory $\calC_{x} \subseteq \calC$
spanned by those objects $C \in \calC$ such that $x \in U(C)$ has weakly contractible nerve.
\end{itemize}
Then $U$ exhibits the $\infty$-category
$\Sing^{A}(X)$ as the colimit (in the $\infty$-category $\Cat_{\infty}$) of the diagram
$$\{ \Sing^{A}( U(C)) \}_{C \in \calC}.$$ 
\end{theorem}

\begin{remark}
Theorem \ref{superkamp} reduces to Theorem \ref{vankamp} in the special case where
$A$ has only a single element.
\end{remark}

The proof of Theorem \ref{superkamp} will occupy our attention throughout this section.
We begin by establishing some notation.

\begin{definition}
Let $A$ be a partially ordered set and $X$ an $A$-stratified topological space.
Given a chain of elements $a_0 \leq \ldots \leq a_n$ in $A$ (which we can identify
with an $n$-simplex $\vec{a}$ in $\Nerve(A)$), we let $\Sing^{A}(X)[ \vec{a}]$ denote
the fiber product $\Fun( \Delta^n, \Sing^{A}(X)) \times_{ \Fun( \Delta^n, \Nerve(A))} \{ \vec{a} \}$. 
\end{definition}

\begin{remark}
Suppose that $X$ is a conically $A$-stratified topological space. It follows immediately
from Theorem \ref{kanna} that for every $n$-simplex $\vec{a}$ of $\Nerve(A)$, the simplicial set
$\Sing^{A}(X)[\vec{a}]$ is a Kan complex.
\end{remark}

\begin{example}
Let $a \in A$ be a $0$-simplex of $\Nerve(A)$, and let $X$ be an $A$-stratified topological space.
Then $\Sing^{A}(X)[a]$ can be identified with the Kan complex $\Sing(X_a)$.
\end{example}

In the special case where $\vec{a} = ( a_0 \leq a_1)$ is an edge of $\Nerve(A)$, the simplicial set
$\Sing^{A}(X)[\vec{a}]$ can be viewed as the space of paths $p: [0,1] \rightarrow X$ such that
$p(0) \in X_{a_0}$ and $p(t) \in X_{a_1}$ for $t \neq 0$. The essential information is encoded in the behavior of the path $p(t)$ where $t$ is close to zero. To make this more precise, we need to introduce a bit of notation.

\begin{definition}
Let $A$ be a partially ordered set, let $X$ be an $A$-stratified topological space, and let 
$a \leq b$ be elements of $A$. We define a simplicial set $\Sing^{A}_{a \leq b}(X)$ as follows:
\begin{itemize}
\item[$(\ast)$] An $n$-simplex of $\Sing^{A}_{a \leq b}(X)$ consists of an equivalence class
of pairs $(\epsilon, \sigma)$, where $\epsilon$ is a positive real number and 
$\sigma: |\Delta^n| \times [0, \epsilon] \rightarrow X$ is a continuous map such that
$\sigma( | \Delta^n | \times \{0\} ) \subseteq X_{a}$ and $\sigma( | \Delta^n | \times (0, \epsilon]) \subseteq X_{b}$. Here we regard $(\epsilon, \sigma)$ and $(\epsilon', \sigma')$ as equivalent if
there exists a positive real number $\epsilon'' < \epsilon, \epsilon'$ such that
$\sigma | (| \Delta^n | \times [0, \epsilon'']) = \sigma' | ( | \Delta^n | \times [ 0, \epsilon''])$.
\end{itemize}
\end{definition}

More informally, we can think of $\Sing^{A}_{a \leq b}(X)$ as the space of germs of paths in $X$
which begin in $X_{a}$ and then pass immediately into $X_{b}$. There is an evident map
$\Sing^{A}(X)[a \leq b] \rightarrow \Sing^{A}_{a \leq b}(X)$, which is given by passing from paths to germs of paths.

\begin{lemma}\label{tankus}
Let $A$ be a partially ordered set, $X$ an $A$-stratified topological space, and
$a \leq b$ elements of $A$. Then the map $\phi: \Sing^{A}(X)[a \leq b] \rightarrow \Sing^{A}_{a \leq b}(X)$ is a weak homotopy equivalence of simplicial sets.
\end{lemma}

\begin{proof}
For every positive real number $\epsilon$, let $S[\epsilon]$ denote the simplicial
set whose $n$-simplces are maps $\sigma: | \Delta^n | \times [0, \epsilon] \rightarrow X$
such that $\sigma( | \Delta^n | \times \{0\} ) \subseteq X_a$ and
$\sigma( | \Delta^n| \times (0, \epsilon]) \subseteq X_{b}$. There are evident restriction maps
$$ \Sing^{A}(X)[a \leq b] = S[1] \rightarrow S[ \frac{1}{2}] \rightarrow S[ \frac{1}{4}] \rightarrow \cdots$$
and the colimit of this sequence can be identified with $\Sing^{A}_{a \leq b}(X)$.
Consequently, to prove that $\phi$ is a weak homotopy equivalence, it will suffice to show that
each of the restriction maps $\psi: S[ \epsilon] \rightarrow S[ \frac{ \epsilon}{2}]$ is a 
weak homotopy equivalence. It now suffices to observe that $\psi$ is a pullback of the
trivial Kan fibration $\Fun( \Delta^1, \Sing(X_{b})) \rightarrow \Fun( \{0\}, \Sing(X_{b}))$. 
\end{proof}

The space of germs $\Sing^{A}_{a \leq b}(X)$ enjoys a formal advantage over the space of paths of fixed  length:

\begin{lemma}\label{tanker}
Let $A$ be a partially ordered set, $X$ a conically $A$-stratified topological space, and $a \leq b$ elements of $A$. Then the restriction map $\Sing^{A}_{a \leq b}(X) \rightarrow \Sing(X_a)$ is a Kan fibration.
\end{lemma}

\begin{proof}
We must show that every lifting problem of the form
$$ \xymatrix{ \Lambda^{n+1}_i \ar[r]^-{F^{0}_{+}} \ar[d] & \Sing^{A}_{a \leq b}(X) \ar[d] \\
\Delta^{n+1} \ar[r]^-{F^{0}_{-}} \ar@{-->}[ur] & \Sing(X_a) }$$
admits a solution. Let us identify $| \Delta^{n+1} |$ with a product
$| \Delta^n | \times [0,1]$ in such a way that the closed subset $| \Lambda^{n+1}_i |$ is identified
with $| \Delta^n | \times \{0\}$. We can identify $F^{0}_{+}$ with a continuous map
$| \Delta^n | \times \{0\} \times [0, \epsilon] \rightarrow X$ for some positive real number $\epsilon$,
and $F^{0}_{-}$ with a continuous map $| \Delta^n | \times [0,1] \times \{0\} \rightarrow X_{a}$.
To solve the lifting problem, we must construct a positive real number
$\epsilon' \leq \epsilon$ and a map $F: | \Delta^n | \times [0,1] \times [0, \epsilon']
\rightarrow X$ compatible with $F^{0}_{-}$ and $F^{0}_{+}$ with the following additional property:
\begin{itemize}
\item[$(\ast)$] For $0 < t$, we have $F(v,s,t) \in X_b$. 
\end{itemize}

For each point $x \in X_{a}$, choose a neighborhood $U_{x}$ of $x$ as in
Definition \ref{sally}. Choose a triangulation of $| \Delta^n |$ and a nonnegative integer
$N \gg 0$ with the property that for each simplex $\tau$ of $| \Delta^n |$ and $0 \leq k < N$, the map
$F^{0}_{-}$ carries $\tau \times [ \frac{k}{N}, \frac{k+1}{N}]$ into some $U_x$ for some point
$x \in X_a$. For each subcomplex $L$ of $| \Delta^n |$, we will prove that there exists
a map $F_{L}: L \times [0,1] \times [0,\epsilon] \rightarrow X$ (possibly after shrinking $\epsilon$)
compatible with $F^{0}_{-}$ and $F^{0}_{+}$ and satisfying condition $\ast$. Taking $L = | \Delta^n |$ we will obtain a proof of the desired result.

The proof now proceeds by induction on the number of simplices of $L$. If $L = \emptyset$ there is nothing to prove. Otherwise, we can write $L = L_0 \cup \tau$, where $\tau$ is a simplex of
$| \Delta^n |$ such that $L_0 \cap \tau = \bd \tau$. By the inductive hypothesis, we may assume
that the map $F_{L_0}$ has already been supplied; let $F_{ \bd \tau}$ be its restriction to
$\bd \tau \times [0,1] \times [0, \epsilon]$. To complete the proof, it will suffice to show that
we can extend $F_{ \bd \tau}$ to a map $F_{\tau}: \tau \times [0,1] \times [0,\epsilon] \rightarrow X$
compatible with $F^{0}_{-}$ and $F^{0}_{+}$ and satisfying $(\ast)$ (possibly after shrinking
the real number $\epsilon$). We again proceed in stages by defining a compatible sequence of maps
$F^{k}_{\tau}: \tau \times [ 0, \frac{k}{N}] \times [0, \epsilon] \rightarrow X$ using induction on $k \leq N$.
The map $F^{0}_{\tau}$ is determined by $F^{0}_{+}$. Assume that
$F^{k-1}_{\tau}$ has already been constructed. Let $K = \tau \times [\frac{k-1}{N}, \frac{k}{N}]$
and let $K_0$ be the closed subset of $K$ given by the union of $\bd \tau \times [ \frac{k-1}{N}, \frac{k}{N}]$ and $\tau \times \{ \frac{k-1}{N} \}$. Then $F^{k-1}_{\tau}$ determines a continuous map
$$g_0: (K \times \{0\} ) \coprod_{ K_0 \times \{0\} } ( K_0 \times [0, \epsilon] ) \rightarrow X.$$
To construct $F^{k}_{\tau}$, it will suffice to extend $g_0$ to a continuous map
$g: K \times [0, \epsilon] \rightarrow X$ satisfying $(\ast)$ (possibly after shrinking $\epsilon$).

By assumption, the map $g_0$ carries $K \times \{0\}$ into some open set $U = U_x$ of
the form $Z \times C(Y)$ described in Definition \ref{sally}. Shrinking $\epsilon$ if necessary,
we may assume that $g_0$ also carries $K_0 \times [0, \epsilon]$ into $U$. 
Let $g'_{0}: (K \times \{0\} ) \coprod_{ K_0 \times \{0\} } ( K_0 \times [0, \epsilon] )
\rightarrow C(Y)$ be the composition of $g_0$ with the projection to $C(Y)$, and let
$g''_0:  (K \times \{0\} ) \coprod_{ K_0 \times \{0\} } ( K_0 \times [0, \epsilon] )
\rightarrow Z$ be defined similarly. Let $r$ be a retraction of
$K$ onto $K_0$, and let
$g'$ be the composition $K \times [0,\epsilon] \rightarrow K_0 \times [0, \epsilon]
\stackrel{g'_0}{\rightarrow} C(Y)$; we observe that $g'$ is an extension of
$g'_0$ (since $g'_0$ is constant on $K \times \{0\}$). Let
$r'$ be a retraction of $K \times [0, \epsilon]$ onto
$(K \times \{0\} ) \coprod_{ K_0 \times \{0\} } ( K_0 \times [0, \epsilon] )$, and let
$g''$ be the composition $g''_0 \circ r'$. The pair $(g', g'')$ determines a map
$g: K \times [0, \epsilon] \rightarrow X$ with the desired properties.
\end{proof}

\begin{proposition}\label{kamma}
Let $A$ be a partially ordered set, let $X$ be a conically $A$-stratified space,
let $U$ be an open subset of $X$ (which inherits the structure of a conically $A$-stratified space),
and let $\vec{a} = (a_0 \leq a_1 \leq \ldots \leq a_n)$ be an $n$-simplex of
$\Nerve(A)$. Then the diagram of Kan complexes
$$ \xymatrix{ \Sing^{A}(U)[ \vec{a} ] \ar[r] \ar[d] & \Sing^{A}(X)[\vec{a}] \ar[d] \\
\Sing(U_{a_0}) \ar[r] & \Sing(X_{a_0}) }$$
is a homotopy pullback square.
\end{proposition}

\begin{proof}
The proof proceeds by induction on $n$. If $n = 0$ the result is obvious.
If $n > 1$, then let $\vec{a}'$ denote the truncated chain $(a_0 \leq a_1)$
and $\vec{a}''$ the chain $(a_1 \leq \ldots \leq a_{n-1} \leq a_{n})$. We have a commutative diagram
$$ \xymatrix{ \Sing^{A}(U)[ \vec{a} ] \ar[r] \ar[d] & \Sing^{A}(X)[\vec{a}] \ar[d] \\
\Sing^{A}(U)[ \vec{a}'] \times_{ \Sing(U_{a_1}) } \Sing^{A}(U)[ \vec{a}''] \ar[r] \ar[d] & 
\Sing^{A}(X)[ \vec{a}'] \times_{ \Sing(X_{a_1}) } \Sing^{A}(X)[ \vec{a}''] \ar[d] \\
\Sing^{A}(U)[ \vec{a}' ] \ar[r] \ar[d] & \Sing^{A}(X)[\vec{a}'] \ar[d] \\
\Sing(U_{a_0}) \ar[r] & \Sing(X_{a_0}). }$$
The upper square is a homotopy pullback because the vertical maps
are weak homotopy equivalences (since $\Sing^{A}(U)$ and $\Sing^{A}(X)$ are
$\infty$-categories, by virtue of Theorem \ref{kanna}). The lower square is a homotopy pullback
by the inductive hypothesis. The middle square is
a (homotopy) pullback of the diagram
$$ \xymatrix{ \Sing^{A}(U)[\vec{a}''] \ar[r] \ar[d] & \Sing^{A}(X)[\vec{a}''] \ar[d] \\
\Sing(U_{a_{1}}) \ar[r] & \Sing(X_{a_1}),}$$
and therefore also a homotopy pullback square by the inductive hypothesis.
It follows that the outer rectangle is a homotopy pullback as required.

It remains to treat the case $n=1$. We have a commutative diagram
$$ \xymatrix{ \Sing^{A}(U)[a_0 \leq a_1] \ar[r] \ar[d] & \Sing^{A}(X)[a_0 \leq a_1] \ar[d] \\
\Sing^{A}_{a_0 \leq a_1}(U) \ar[r] \ar[d] & \Sing^{A}_{a_0 \leq a_1}(X) \ar[d] \\
\Sing( U_{a_0}) \ar[r] & \Sing( X_{a_0} ). }$$
The lower square is a homotopy pullback since it is a pullback square in which the
vertical maps are Kan fibrations (Lemma \ref{tanker}). The upper square is a homotopy
pullback since the upper vertical maps are weak homotopy equivalences (Lemma \ref{tankus}).
It follows that the outer square is also a homotopy pullback, as desired.
\end{proof}

We are now ready to establish our main result.

\begin{proof}[Proof of Theorem \ref{superkamp}]
Let $f: \Cat_{\infty} \rightarrow \Fun( \Nerve(\cDelta)^{op}, \SSet)$ be the functor
described in Corollary \bicatref{presquare}. Then $f$ is a fully faithful embedding, whose essential
image consists of the complete Segal objects in the $\infty$-category
$\Fun( \Nerve(\cDelta)^{op}, \SSet)$ of simplicial spaces. It will therefore suffice to prove
that the composite functor
$$ \theta: \Nerve(\calC)^{\triangleright} \rightarrow \Cat_{\infty} \rightarrow \Fun( \Nerve(\cDelta)^{op}, \SSet)$$
is a colimit diagram. Since colimits in $\Fun( \Nerve(\cDelta)^{op}, \SSet)$ are computed pointwise,
it will suffice to show that $\theta$ determines a colimit diagram in $\SSet$ after evaluation at
each object $[n] \in \cDelta$. Unwinding the definitions, we see that this diagram is given by
the formula
$$ C \mapsto \coprod_{ a_0 \leq a_1 \leq \cdots \leq a_n} \Sing^{A}( U(C) )[ a_0 \leq \ldots \leq a_n].$$
Since the collection of colimit diagrams is stable under coproducts (Lemma \toposref{limitscommute}), 
it will suffice to show that for every $n$-simplex $\vec{a} = (a_0 \leq \ldots \leq a_n)$ of $\Nerve(A)$, the functor
$$ \theta_{ \vec{a}}: \Nerve(\calC)^{\triangleright} \rightarrow \SSet$$
given by the formula $C \mapsto \Sing^{A}(U(C))[ \vec{a} ]$ is a colimit diagram in $\SSet$.

We have an evident natural tranformation $\alpha: \theta_{ \vec{a}} \rightarrow \theta_{a_0}$.
The functor $\theta_{a_0}$ is a colimit diagram in $\SSet$: this follows by applying Theorem \ref{vankamp} to the stratum $X_{a}$. Proposition \ref{kamma} guarantees that $\alpha$ is a Cartesian natural transformation. Since $\SSet$ is an $\infty$-topos, Theorem \toposref{mainchar} guarantees
that $\theta_{ \vec{a}}$ is also a colimit diagram, as desired.
\end{proof}

\subsection{Digression: Complementary Colocalizations}\label{sec5sub9}

In this section, we will describe a technical device which can be used to show that a functor
$F: \calC \rightarrow \calD$ is an equivalence of $\infty$-categories, assuming that 
$\calC$ and $\calD$ can be decomposed into ``pieces'' on which $F$ is known be an equivalence.
To make this idea more precise, we need to introduce some terminology.

\begin{definition}\label{colpus}
Let $\calC$ be an $\infty$-category which admits pushouts and which contains a pair of full subcategories
$\calC_{0}, \calC_{1} \subseteq \calC$. Assume that each of the inclusions
$\calC_{i} \subseteq \calC$ admits a right adjoint $L_{i}$. We will say that
{\it $L_0$ is complementary to $L_{1}$} if the following conditions are satisfied:
\begin{itemize}
\item[$(1)$] The $\infty$-category $\calC$ admits pushouts.
\item[$(2)$] The functors $L_0$ and $L_1$ preserve pushouts.
\item[$(3)$] The functor $L_1$ carries every morphism in $\calC_0$ to an equivalence.
\item[$(4)$] If $\alpha$ is a morphism in $\calC$ such that $L_0(\alpha)$ and
$L_1(\alpha)$ are both equivalences, then $\alpha$ is an equivalence.
\end{itemize}
\end{definition}

\begin{warning}
The relation of complementarity is not symmetric in $L_0$ and $L_1$.
\end{warning}

Our main result is the following:

\begin{proposition}\label{inkor}
Let $\calC$ and $\calC'$ be $\infty$-categories which admit pushouts.
Suppose that we are given inclusions of full subcategories
$\calC_0, \calC_1 \subseteq \calC$, $\calC'_0, \calC'_{1} \subseteq \calC'$, which
admit right adjoints $L_0$, $L_1$, $L'_0$, and $L'_{1}$. Let $F: \calC \rightarrow \calC'$
be a functor satisfying the following conditions:
\begin{itemize}
\item[$(1)$] The colocalizations $L_0$ and $L_1$ are complementary in $\calC$, and
the colocalizations $L'_0$ and $L'_1$ are complementary in $\calC'$.
\item[$(2)$] The functor $F$ restricts to equivalences $\calC_0 \rightarrow \calC'_0$ and
$\calC_1 \rightarrow \calC'_{1}$.
\item[$(3)$] The functor $F$ preserves pushout squares.
\item[$(4)$] Let $C \in \calC_1$, and let $\alpha: C_0 \rightarrow C$ be a map which exhibits
$C_0$ as a $\calC_0$-colocalization of $C$. Then $F(\alpha)$ exhibits
$F(C_0) \in \calC_0$ as a $\calC'_0$-localization of $F(C) \in \calC'_{1} \subseteq \calC'$. 
\end{itemize}
Then $F$ is an equivalence of $\infty$-categories.
\end{proposition}

Before proving Proposition \ref{inkor}, we describe a mechanism by which an
$\infty$-category $\calC$ can be recovered from a pair of complementary colocalizations.

\begin{lemma}\label{kuro}
Let $\calC$ be an $\infty$-category which admits pushouts, and let
$\calC_0, \calC_1 \subseteq \calC$ be full subcategories. Assume that
the inclusions $\calC_i \subseteq \calC$ admit right adjoints $L_i$, and that
$L_0$ is complementary to $L_1$. Let $\calD$ be the full subcategory of
$\Fun( \Delta^1 \times \Delta^1, \calC)$ spanned by those diagrams $\sigma$:
$$ \xymatrix{ C_{01} \ar[r] \ar[d] & C_0 \ar[d] \\
C_{1} \ar[r] & C }$$
satisfying the following conditions:
\begin{itemize}
\item[$(i)$] The diagram $\sigma$ is a pushout square.
\item[$(ii)$] The object $C_0$ belongs to $\calC_0$, and the object
$C_1$ belongs to $\calC_{1}$.
\item[$(iii)$] The map $C_{01} \rightarrow C_{1}$ exhibits $C_{01}$ as
a $\calC_0$-colocalization of $C_1$.
\end{itemize}
Then the evaluation functor $\sigma \mapsto C$ determines a trivial Kan fibration $\calD \rightarrow \calC$.
\end{lemma}

\begin{proof}
Let $\calD'$ denote the full subcategory of $\calD$ spanned by those diagrams
which satisfy the following additional conditions:
\begin{itemize}
\item[$(iv)$] The map $C_0 \rightarrow C$ exhibits $C_0$ as a $\calC_0$-colocalization of $C$.
\item[$(v)$] The map $C_1 \rightarrow C$ exhibits $C_1$ as a $\calC_1$-colocalization of $C$.
\end{itemize}
Let $\calD''$ denote the full subcategory of $\Fun( \Delta^1 \times \Delta^1, \calC)$ spanned by those functors which satisfy conditions $(ii)$ through $(v)$. Let $\overline{\calC}$ denote the full subcategory of $\calC \times \Delta^1 \times \Delta^1$
spanned by those objects $(C, i,j)$ such that $C \in \calC_0$ if $i=0$ and
$C \in \calC_1$ if $0 = i < j = 1$. Let $p: \overline{\calC} \rightarrow \Delta^1 \times \Delta^1$ denote the projection map. We observe that $\calD''$ can be identified with the $\infty$-category of functors
$F \in \Fun_{ \Delta^1 \times \Delta^1}( \Delta^1 \times \Delta^1, \overline{\calC})$ which are $p$-right Kan extensions of $F | \{ (1,1) \}$. It follows from Proposition \toposref{lklk} that the evaluation
functor $\calD'' \rightarrow \calC$ is a trivial Kan fibration. We will complete the proof by showing that 
$\calD = \calD' = \calD''$.

To prove that $\calD' = \calD''$, consider a diagram $\sigma:$
$$ \xymatrix{ L_0 L_1 C \ar[r] \ar[d] & L_0 C \ar[d] \\
L_1 C \ar[r] & C }$$
which belongs to $\calD''$. This diagram induces a map
$\alpha: L_0 C \coprod_{ L_0 L_1 C} L_1 C \rightarrow C$; to prove that $\sigma \in \calD'$
we must show that $\alpha$ is an equivalence. For this, it suffices to show that both $L_0(\alpha)$
and $L_1(\alpha)$ are equivalences. Since $L_0$ and $L_1$ preserve pushout squares, we
are reduced to proving that the diagrams $L_0(\sigma)$ and $L_1(\sigma)$ are pushout
squares. This is clear: in the diagram $L_0(\sigma)$, the vertical maps are both equivalences;
in the diagram $L_1(\sigma)$, the horizontal maps are both equivalences.

To show that $\calD = \calD'$, consider an arbitrary diagram
$$ \xymatrix{ C_{01} \ar[r] \ar[d] & C_0 \ar[d] \\
C_{1} \ar[r] & C }$$
satisfying $(i)$ through $(iii)$. Since $L_0$ preserves pushouts, we have a pushout
diagram
$$ \xymatrix{ L_0 C_{01} \ar[r] \ar[d] & L_0 C_0 \ar[d] \\
L_0 C_{1} \ar[r] & L_0 C. }$$
The left vertical map is an equivalence by assumption $(iii)$, so the right
vertical map is also an equivalence. Since $C_0 \in \calC_0$ by $(ii)$, this proves $(iv)$.
Similarly, since the functor $L_1$ preserves pushouts we have a pushout diagram
$$ \xymatrix{ L_1 C_{01} \ar[r] \ar[d] & L_1 C_0 \ar[d] \\
L_1 C_{1} \ar[r] & L_1 C. }$$
Since $L_0$ is complementary to $L_1$, the upper horizontal map is an equivalence.
It follows that the lower horizontal map is an equivalence. Since $C_1 \in \calC_1$ by
$(ii)$, we conclude that assertion $(v)$ holds.
\end{proof}

\begin{proof}[Proof of Proposition \ref{inkor}]
Let $\calD$ and $\calD'$ be defined as in Lemma \ref{kuro}. Let
$\calE$ be the full subcategory of $\Fun( \Lambda^2_0, \calC)$ spanned by those diagrams
$$ C_0 \leftarrow C_{01} \stackrel{\alpha}{\rightarrow} C_{1}$$
where $C_0 \in \calC_0$, $C_1 \in \calC_1$, and $\alpha$ exhibits $C_{01}$ as a 
$\calC_{0}$-colocalization of $C_1$. Proposition \toposref{lklk} guarantees that the restriction
functor $\calD \rightarrow \calE$ is a trivial Kan fibration. Similarly, we have a trivial
Kan fibration $\calD' \rightarrow \calE$. These maps fit into a commutative diagram
$$ \xymatrix{ \calC \ar[r]^{F} & \calC' \\
\calD \ar[u] \ar[r] \ar[d] & \calD' \ar[u] \ar[d]  \\
\calE \ar[r]^{F_0} & \calE }$$
Using Lemma \ref{kuro}, we are reduced to the problem of showing that $F_0$ is an
equivalence of $\infty$-categories. The map $F_0$ extends to a map of (homotopy) pullback diagrams
$$ \xymatrix{ \calE \ar[r] \ar[d] & \Fun( \Delta^1, \calC_0) \ar[d] & \calE' \ar[r] \ar[d] & \Fun( \Delta^1, \calC'_0) \ar[d] \\
\calM \ar[r] & \Fun( \{0\}, \calC_0) & \calM' \ar[r] & \Fun( \{0\}, \calC'_0), }$$
where $\calM$ is denotes the full subcategory of $\calC$ spanned by those morphisms $f: C_{01} \rightarrow C_1$ such that $C_1 \in \calC_1$ and $f$ exhibits $C_{01}$ as a $\calC_0$-colocalization of $C_1$, and $\calM'$ is defined similarly. Since $F$ induces an equivalence
$\calC_0 \rightarrow \calC'_0$ by assumption, it suffices to show that the map
$\calM \rightarrow \calM'$ (which is well-defined by virtue of $(3)$) is an equivalence of
$\infty$-categories. This follows from the assumption that $F$ restricts to an equivalence
$\calC_1 \rightarrow \calC'_{1}$, since we have a commutative diagram
$$ \xymatrix{ \calM \ar[r] \ar[d] & \calM' \ar[d] \\
\calC_1 \ar[r] & \calC'_{1} }$$
in which the vertical maps are trivial Kan fibrations.
\end{proof}

In order to apply Proposition \ref{inkor}, we will need to be able to recognize the existence of complementary colocalizations. The following result provides a useful criterion:

\begin{proposition}\label{soake}
Let $p: \calM \rightarrow \Delta^1$ be a correspondence between $\infty$-categories.
Assume that there exists a retraction $r$ from $\calM$ onto the full subcategory $\calM_0$.
Let $\calA$ be an $\infty$-category which admits finite colimits, and let
$\calC = \Fun( \calM, \calA)$. We define full subcategories
$\calC_0, \calC_1 \subseteq \calC$ as follows:
\begin{itemize}
\item[$(a)$] A functor $f: \calM \rightarrow \calA$ belongs to 
$\calC_0$ if $f$ is a left Kan extension of $f| \calM_1$ (that is,
if $f(M)$ is an initial object of $\calA$, for each $M \in \calM_0$).
\item[$(b)$] A functor $f: \calM \rightarrow \calA$ belongs to
$\calC_1$ if $f(\alpha)$ is an equivalence, for every $p$-coCartesian
morphism in $\calM$.
\end{itemize}
Then:
\begin{itemize}
\item[$(1)$] The inclusion functors $\calC_0 \subseteq \calC$ and
$\calC_1 \subseteq \calC$ admit
right adjoints $L_0$ and $L_1$.
\item[$(2)$] The functor $L_0$ is complementary to $L_1$.
\end{itemize}
\end{proposition}

\begin{proof}
We prove $(1)$ by explicit construction. The functor $L_0$ is given by
composing the restriction functor $\Fun( \calM, \calA) \rightarrow \Fun( \calM_1, \calA)$ with
a section of the trivial Kan fibration $\calC_0 \rightarrow \Fun( \calM_1, \calA)$.
The functor $L_1$ is given by composing the restriction functor
$\Fun( \calM, \calA) \rightarrow \Fun( \calM_0, \calA)$ with the retraction
$r: \calM \rightarrow \calM_0$.

To prove $(2)$, we must verify that the conditions of Definition \ref{colpus} are satisfied.
It is clear that the $\infty$-category $\calC$ admits pushouts (these are computed pointwise),
and the above constructions show that $L_0$ and $L_1$ preserve pushouts.
The restriction $L_1 | \calC_0$ factors through $\Fun( \calM_0, \calA')$, where
$\calA' \subseteq \calA$ is the contractible Kan complex spanned by the final objects of $\calA$.
It follows that $L_1 | \calC_0$ is essentially constant (and, in particular, that $L_1$ carries every morphism in $\calC_0$ to an equivalence). Finally, if $\alpha$ is a morphism in $\calC$ such that
$L_0(\alpha)$ and $L_1(\alpha)$ are both equivalences, then $\alpha$ is a natural
transformation of functors from $\calM$ to $\calA$ which induces an equivalence after evaluation
at every object $\calM_1$ and every object of $\calM_0$. Since every object of $\calM$ belongs to $\calM_0$ or $\calM_1$, we conclude that $\alpha$ is an equivalence.
\end{proof}

\subsection{Exit Paths and Constructible Sheaves}\label{cloop}

Let $A$ be a partially ordered set and let $X$ be a space equipped with an $A$-stratificatin
$f: X \rightarrow A$. Our goal in this section is to prove that, if $X$ is sufficiently well-behaved, then
the $\infty$-category of $A$-constructible objects of $\Shv(X)$ can be identified with
the $\infty$-category $\Fun( \Sing^{A}(X), \SSet)$, where $\Sing^{A}(X)$ is the $\infty$-category
of exit paths defined in \S \ref{lubos}. In fact, we will give an explicit construction of this equivalence,
generalizing the analysis we carried out for locally constant sheaves in \S \ref{lubla}. First, we need to establish a bit of terminology.

\begin{notation}
Let $A$ be a partially ordered set and let $X$ be a paracompact $A$-stratified space.
We let $\bfA_{X}$ denote the category $( \sSet)_{ / \Sing^{A}(X) }$, which we regard as
endowed with the covariant model structure described in \S \toposref{contrasec}.
Let $\calB(X)$ denote the partially ordered collection of all open $F_{\sigma}$ subsets of $X$. 
We let $\Shv(X)$ denote the full subcategory of $\calP( \calB(X) )$ spanned by those objects
which are sheaves with respect to the natural Grothendieck topology on $\calB(X)$.
\end{notation}

Proposition \toposref{gumby444} and Theorem \toposref{struns} furnish a chain of
equivalences of $\infty$-categories
$$ \Fun( \Sing^{A}(X), \SSet) \leftarrow \Nerve( ( \Set_{\Delta}^{ \sCoNerve[ \Sing^{A}(X)] })^{\degree})
\rightarrow \Nerve( \bfA^{\degree}_{X} ).$$ 

\begin{construction}
We define a functor $\theta: \calB(X)^{op} \times \bfA_{X} \rightarrow \sSet$
by the formula $$\theta(U,Y) = \Fun_{ \Sing^{A}(X)}( \Sing^{A}(U), Y).$$
Note that if $Y \in \bfA_X$ is fibrant, then $Y \rightarrow \Sing^{A}(X)$ is a left fibration
so that each of the simplicial sets $\theta(U,Y)$ is a Kan complex.
Passing to the nerve, $\theta$ induces a map of $\infty$-categories
$\Nerve( \calB(X)^{op} ) \times \Nerve( \bfA^{\degree}_{X}) \rightarrow \SSet$, which
we will identify with a map of $\infty$-categories
$$ \Psi_{X}: \Nerve( \bfA^{\degree}_{X} ) \rightarrow \calP( \calB(X) ).$$
\end{construction}

We are now ready to state the main result of this section.

\begin{theorem}\label{mainstrat}
Let $X$ be a paracompact topological space which is locally of singular shape and is equipped with a conical $A$-stratification, where
$A$ is a partially ordered set satisfying the ascending chain condition.
Then the functor $\Psi_{X}$ induces an equivalence $\Nerve( \bfA^{\degree}_{X}) \rightarrow \Shv^{A}(X)$.
\end{theorem}

The proof of Theorem \ref{mainstrat} will be given at the end of this section, after we have developed a number of preliminary ideas. For later use, we record the following easy consequence of Theorem \ref{mainstrat}:

\begin{corollary}
Let $X$ be a paracompact topological space which is locally of singular shape and is equipped with a conical $A$-stratification, where $A$ is a partially ordered set satisfying the ascending chain condition. Then the inclusion
$i: \Sing^{A}(X) \hookrightarrow \Sing(X)$ is a weak homotopy equivalence of simplicial sets.
\end{corollary}

\begin{proof}
Let $X'$ denote the topological space $X$ equipped with the trivial stratification. 
The inclusion $i$ induces a pullback functor $i^{\ast}: \Nerve(\bfA_{X'}^{\degree}) \rightarrow \Nerve(\bfA^{\degree}_{X})$, and we have an evident natural transformation
$\alpha: \Psi_{X} \circ i^{\ast} \rightarrow \Psi_{X'}$ from $\Nerve( \bfA_{X'}^{\degree})$ to
$\Shv(X)$. We claim that $\alpha$ is an equivalence. Since both functors
take values in the full subcategory of hypercomplete objects of $\Shv(X)$ (Lemma \ref{curin} and Proposition \ref{toru}), it suffices to show that $\alpha(Y)$ is $\infty$-connective for
each $Y \in \Nerve( \bfA_{X'}^{\degree})$. For this, it suffices to show that $x^{\ast} \alpha(Y)$
is an equivalence for every point $x \in X$ (Lemma \ref{copus}). Using Proposition \ref{swilker}, we can reduce to the case $X = \{x\}$ where the result is obvious. Applying the functor of global sections
to $\alpha$, we deduce that for every Kan fibration $Y \rightarrow \Sing(X)$ the restriction map
$$ \Fun_{ \Sing(X)}( \Sing(X), Y) \rightarrow \Fun_{ \Sing(X)}( \Sing^{A}(X), Y)$$ is a
homotopy equivalence of Kan complexes, which is equivalent to the assertion that $i$
is a weak homotopy equivalence.
\end{proof}

We now turn to the proof of Theorem \ref{mainstrat} itself. Our first objective 
is to show that the functor $\Psi_X$ takes values in the the full subcategory
$\Shv(X) \subseteq \calP( \calB(X) )$. 

\begin{lemma}
Let $A$ be a partially ordered set, let $X$ be a paracompact topological space
equipped with a conical $A$-stratification. The functor $\Psi_{X}: \Nerve( \bfA^{\degree}_{X}) \rightarrow \calP( \calB(X) )$
factors through the full subcategory $\Shv(X) \subseteq \calP(\calB(X))$.
\end{lemma}

\begin{proof}
Let $U \in \calB(X)$, and let $S \subseteq \calB(U)$ be a covering sieve on $U$. In view of
Theorem \toposref{colimcomparee}, it will suffice to show that
for every left fibration $Y \rightarrow \Sing^{A}(X)$, the canonical map
$$ \Fun_{ \Sing^{A}(X)}(\Sing^{A}(U) , Y) \rightarrow \varprojlim_{V \in S} \Fun_{ \Sing^{A}(X)}( \Sing^{A}(V), Y)$$ 
exhibits the Kan complexes $\Fun_{ \Sing^{A}(X) }( \Sing^{A}(U), Y)$ as a homotopy limit of the diagram
of Kan complexes $\{ \Fun_{ \Sing^{A}(X)}( \Sing^{A}(V), Y) \}_{V \in S}$. For this, it suffices to show
that $\Sing^{A}(U)$ is a homotopy colimit of the simplicial sets $\{ \Sing^{A}(V) \}_{V \in S}$ in
the category $(\sSet)_{ / \Sing^{A}(X)}$, endowed with the covariant model structure. This follows
from the observation that the covariant model structure on $( \sSet)_{/ \Sing^{A}(X)}$ is a localization
of the Joyal model structure, and $\Sing^{A}(U)$ is a homotopy colimit of
$\{ \Sing^{A}(V) \}_{V \in S}$ with respect to the Joyal model structure (by Theorems \ref{superkamp}
and \toposref{colimcomparee}).
\end{proof}

\begin{remark}\label{sablo}
Let $X$ be a paracompact space equipped with an $A$-stratification. For each open $F_{\sigma}$ subset
$U$ of $X$, the composition of $\Psi_{X}: \Nerve( \bfA^{\degree}) \rightarrow
\Shv(X)$ with the evaluation functor $\calF \mapsto \calF(U)$ from $\Shv(X)$ to $\SSet$ is equivalent
to the functor $\Nerve( \bfA^{\degree}_{X}) \rightarrow \SSet$ corepresented by (a fibrant replacement for) the object $\Sing^{A}(U) \in \bfA_{X}$. It follows that $\Psi_{X}$ preserves small limits.
\end{remark}

\begin{remark}\label{sailo}
Combining Remark \ref{sablo} with Proposition \toposref{eaa}, we deduce that the functor
$$\Psi_{X}: \Fun( \Sing^{A}(X), \SSet) \simeq \Nerve( \bfA^{\degree}_{X}) \rightarrow
\Shv(X)$$ preserves $n$-truncated objects for each $n \geq -1$.
Since every object $F \in \Fun( \Sing^{A}(X), \SSet)$ equivalent to a limit
of truncated objects (since Postnikov towers in $\SSet$ are convergent),
we deduce from Remark \ref{sablo} that $\Psi_{X}(F)$ is also equivalent to a limit
of truncated objects, and therefore hypercomplete.
\end{remark}

We now discuss the functorial behavior of the map $\Psi_{X}$.
Let $f: X' \rightarrow X$ be a continuous map of paracompact spaces. Let
$A$ be a partially ordered set such that $X$ is endowed with an $A$-stratification.
Then $X'$ inherits an $A$-stratification. The map $f$ determines a morphism of simplicial
sets $\Sing^{A}(X') \rightarrow \Sing^{A}(X)$; let $r: \bfA_{X} \rightarrow \bfA_{X'}$ be the
associated pullback functor and $R: \Nerve( \bfA_{X}^{\degree}) \rightarrow
\Nerve( \bfA_{X'}^{\degree})$ the induced map of $\infty$-categories.
For each $U \in \calB(X)$, we have $f^{-1} U \in \calB(X')$. The canonical map
$\Sing^{A}( f^{-1} U) \rightarrow \Sing^{A}(U)$ induces a map
$\theta_{X}(U,Y) \rightarrow \theta_{X'}( f^{-1} U, r(Y))$. These maps together
determine a natural transformation of functors
$\Psi_{X} \rightarrow f_{\ast} \Psi_{X'} R$
from $\Nerve( \bfA_{X}^{\degree})$ to $\Shv(X)$. We let
$\phi_{X',X}: f^{\ast} \Psi_{X} \rightarrow \Psi_{X'} R$ denote the adjoint transformation (which is well-defined up to homotopy).

\begin{example}\label{ulj}
If $X'$ is an open $F_{\sigma}$ subset of $X$, then the pullback functor
$f^{\ast}: \Shv(X) \rightarrow \Shv(X')$ can be described as the restriction
along the inclusion of partially ordered sets $\calB(X') \subseteq \calB(X)$.
In this case, the natural transformation $\phi_{X',X}$ can be chosen to be an isomorphism of simplicial sets, since the maps 
$\theta_{X}(U,Y) \rightarrow \theta_{X'}(U, r(Y) )$ are isomorphisms for $U \subseteq X'$.
\end{example}

\begin{lemma}\label{olba}
Let $X$ be a paracompact topological space equipped with an $A$-stratification.
Let $a \in A$, let $X' = X_{a}$, and let $f: X' \rightarrow X$ denote the inclusion map. Assume that
$X_{a}$ is paracompact. Then the natural transformation $\phi_{X',X}$ defined above is an equivalence.
\end{lemma}

\begin{proof}
Fix a left fibration $M \rightarrow \Sing^{A}(X)$, and let $M' = M \times_{ \Sing^{A}(X)} \Sing(X_a)$.
We wish to show that $\phi_{X',X}$ induces an equivalence of sheaves $f^{\ast} \Psi_{X}(M)
\rightarrow \Psi_{X_a}(M')$. This assertion is local on $X_{a}$. We may therefore use
Lemma \ref{olit} (and Example \ref{ulj}) to reduce to the case where
$X$ has the form $Z \times C(Y)$, where $Y$ is an $A_{> a}$-stratified space.
Corollary \toposref{snottle} implies that the left hand side can be identified with the (filtered) colimit
$\varinjlim_{V} (\Psi_{X}(Y))(V)$, where $V$ ranges over the collection of all open neighborhoods of
$Z$ in $Z \times C(Y)$. In view of Lemma \ref{cuppa}, it suffices to take the same limit indexed
by those open neighborhoods of the form $V_{g}$, where $g: Z \rightarrow (0, \infty)$ is a continuous function. It will therefore suffice to show that each of the maps $\Psi_{X}(Y)(V_{g}) \rightarrow
\Psi_{X_a}(Y')(Z)$ is a homotopy equivalence. This map is given by the restriction
$$ \Fun_{ \Sing^{A}(X)}( \Sing^{A}(V_g), Y) \rightarrow \Fun_{ \Sing^{A}(X)}( \Sing(Z), Y).$$
To show that this map is a homotopy equivalence, it suffices to show that the inclusion
$i: \Sing(Z) \hookrightarrow \Sing^{A}( V_g)$ is a covariant equivalence in $\Sing^{A}(X)$.
We will show that $i$ is left anodyne. Let $h: C(Y) \times [0,1] \rightarrow C(Y)$
be the map which carries points$(y,s,t) \in Y \times (0, \infty) \times (0,1]$ to $(y,st) \in Y \times (0,\infty)$,
and every other point to the cone point of $C(Y)$. Then $h$ induces a homotopy 
$H: V_{g} \times [0,1] \rightarrow V_{g}$ from the projection $V_{g} \rightarrow Z \subseteq V_g$
to the identity map on $V_{g}$. The homotopy $H$ determines a natural transformation
from the projection $\Sing^{A}(V_g) \rightarrow \Sing(Z)$ to the identity map from
$\Sing^{A}(V_g)$ to itself, which exhibits the map $i$ as a retract of the left anodyne inclusion
$$ (\Sing(Z) \times \Delta^1) \coprod_{ \Sing(Z) \times \{0\} } ( \Sing^{A}(V_g) \times \{0\} )
\subseteq \Sing^{A}( V_g) \times \Delta^1.$$
\end{proof}

\begin{lemma}\label{curin}
Let $X$ be a paracompact topological space which is locally of singular shape and is equipped with a conical $A$-stratification. Then the functor $\Psi_{X}: \Nerve( \bfA^{\degree}_{X}) \rightarrow \Shv(X)$ factors through
the full subcategory $\Shv^{A}(X) \subseteq \Shv(X)$ spanned by the $A$-constructible sheaves on $X$.
\end{lemma}

\begin{proof}
Choose a left fibration $Y \rightarrow \Sing^{A}(X)$ and an element $a \in A$; we wish to prove
that $(\Psi_X(Y) | X_{a}) \in \Shv(X_a)$ is locally constant. The assertion is local on $X$,
so we may assume without loss of generality that $X$ has the form $Z \times C(Y)$
(Lemma \ref{olit}), so that $X_{a} \simeq Z$ is locally of singular shape (Remark \ref{juline}). 
Using Lemma \ref{olba}, we can replace $X$ by $Z$, and thereby reduce to the case where
$X$ consists of only one stratum. In this case, the desired result follows from
Theorem \ref{squ}.
\end{proof}

\begin{lemma}\label{sooit}
Let $X$ be a paracompact topological space of the form $Z \times C(Y)$, and let
$\pi: X \rightarrow Z$ denote the projection map. Then the pullback functor
$\pi^{\ast}: \Shv(Z) \rightarrow \Shv(X)$ is fully faithful.
\end{lemma}

\begin{proof}
Fix an object $\calF \in \Shv(Z)$; we will show that the unit map
$\calF \rightarrow \pi_{\ast} \pi^{\ast} \calF$ is an equivalence.
In view of Corollary \toposref{wamain}, we may suppose that there exists
a map of topological spaces $Z' \rightarrow Z$ such that $\calF$ is 
given by the formula $U \mapsto \bHom_{ \Top_{/Z}}( U, Z')$. Using
the results of \S \toposref{dooky}, we may suppose also that
$\pi^{\ast} \calF$ is given by the formula $V \mapsto \bHom_{ \Top_{/X}}( V, Z' \times_{ Z} X)$.
It will suffice to show that the induced map
$\calF(U) \rightarrow (\pi^{\ast} \calF)( \pi^{-1} U)$ is a homotopy equivalence
for each $U \in \calB(Z)$. Replacing $Z$ by $U$, we may assume that $U = Z$.
In other words, we are reduced to proving that the map
$$ \bHom_{ \Top_{/Z}}( Z, Z') \rightarrow \bHom_{ \Top_{/Z}}( X, Z')$$
is a homotopy equivalence of Kan complexes. This follows from the observation
that there is a deformation retraction from $X$ onto $Z$ (in the category $\Top_{/Z}$ of topological spaces over $Z$).
\end{proof}

\begin{lemma}\label{corpwell}
Let $X$ be a paracompact space of the form $Z \times C(Y)$, let
$\pi: X \rightarrow Z$ denote the projection map, and let $i: Z \rightarrow X$ be the inclusion.
Let $\calF \in \Shv(X)$ be a sheaf whose restriction to $Z \times Y \times (0, \infty)$ is foliated.
Then the canonical map $\pi_{\ast} \calF \rightarrow i^{\ast} \calF$ is an equivalence.
\end{lemma}

\begin{proof}
It will suffice to show that for every $U \in \calB(Z)$, the induced map
$\calF( \pi^{-1}(U) ) \rightarrow (i^{\ast} \calF)(U)$ is a homotopy equivalence. Replacing
$Z$ by $U$, we can assume $U = Z$. Using Corollary \toposref{snottle}, we can identify
$( i^{\ast} \calF)(Z)$ with the filtered colimit $\varinjlim_{V} \calF(V)$, where $V$ ranges
over all open neighborhoods of $Z$ in $X$. In view of Lemma \ref{cuppa}, it suffices to take
the colimit over the cofinal collection of open sets of the form $V_{f}$, where $f: Z \rightarrow (0, \infty)$ is a continuous map. To prove this, it suffices to show that each of the restriction maps
$\theta: \calF(X) \rightarrow \calF( V_{f})$ is an equivalence. Let $W \subseteq Z \times Y \times (0, \infty)$
be the set of triples $(z,y,t)$ such that $t > \frac{ f(z)}{2}$, so that we have a pullback diagram
$$ \xymatrix{ \calF(X) \ar[r]^{\theta} \ar[d] & \calF(V_{f}) \ar[d] \\
\calF(W) \ar[r]^-{\theta'} & \calF( W \cap V_{f} ). }$$
To prove that $\theta$ is a homotopy equivalence, it suffices to show that $\theta'$ is a homotopy
equivalence. The map $\theta'$ fits into a commutative diagram
$$ \xymatrix{ \calF(W) \ar[rr]^{\theta'} \ar[dr] & & \calF( W \cap V_{f} ) \ar[dl] \\
& (s^{\ast} \calF)( Z \times Y), & }$$
where $s: Z \times Y \rightarrow W \cap V_{f}$ is the section given by the continuous map
$\frac{3}{4} f: Z \rightarrow (0, \infty)$. Since $\calF$ is foliated, Proposition \ref{postcanca} and Lemma \ref{ander} guarantee that the vertical maps in this diagram are both equivalences, so that
$\theta'$ is an equivalence as well.
\end{proof}

\begin{lemma}\label{sou}
Let $A$ be a partially ordered set containing an element $a$. 
Let $X$ be a paracompact $A_{\geq a}$-stratified topological space of the form 
$Z \times C(Y)$, where $Y$ is an $A_{> a}$-stratified space.
Let $\calC = \Shv^{A}(X)$. Let $j: Z \times Y \times (0,\infty) \rightarrow X$ denote the inclusion and
let $\calC_0$ denote the intersection of $\calC$ with the essential image of
the left adjoint $j_{!}: \Shv( Z \times Y \times (0,\infty)) \rightarrow \Shv(X)$ to the pullback functor
$j^{\ast}$. Let $\pi: X \rightarrow Z$ be the projection map, and let $\calC_1$ denote the intersection
of $\calC$ with the essential image of $\pi^{\ast}$ (which is fully faithful by Lemma \ref{sooit}). Then:
\begin{itemize}
\item[$(1)$] The inclusion functors $\calC_0 \subseteq \calC$ and
$\calC_1 \subseteq \calC$ admit
right adjoints $L_0$ and $L_1$.
\item[$(2)$] The functor $L_0$ is complementary to $L_1$.
\end{itemize}
\end{lemma}

\begin{proof}
Let $i: Z \rightarrow X$ be the inclusion map. 
The functor $L_0$ is given by the composition $j_{!} j^{\ast}$, and the functor
$L_1$ is given by the composition $\pi^{\ast} \pi_{\ast}$ (which is equivalent
to $\pi^{\ast} i^{\ast}$ by Lemma \ref{corpwell}, and therefore preserves constructibility
and pushout diagrams). Since the composition $i^{\ast} j_{!}$ is equivalent to the constant
functor $\Shv( Z \times Y \times (0, \infty)) \rightarrow \Shv(Z)$ (taking value equal to the initial object of $\Shv(Z)$), the functor $L_1$ carries every morphism in $\calC_0$ to an equivalence.
Finally, suppose that $\alpha$ is a morphism in $\calC$ such that $L_0(\alpha)$ and
$L_1(\alpha)$ are equivalences. Since $j_{!}$ and $\pi^{\ast}$ are fully faithful,
we conclude that $j^{\ast}(\alpha)$ and $i^{\ast}(\alpha)$ are equivalences, so that
$\alpha$ is an equivalence (Corollary \ref{sadwell}).
\end{proof}

\begin{lemma}\label{sogg}
Let $X$ be a paracompact topological space which is locally of singular shape and is equipped with a conical $A$-stratification. 
Then the full subcategory $\Shv^{A}(X) \subseteq \Shv(X)$ is stable under finite colimits in $\Shv(X)$.
\end{lemma}

\begin{proof}
Let $\calF \in \Shv(X)$ be a finite colimit of $A$-constructible sheaves; we wish to show that
$\calF | X_{a}$ is constructible for each $a \in A$. The assertion is local; we may therefore
assume that $X$ has the form $Z \times C(Y)$ (Lemma \ref{olit}). Then
$X_{a} \simeq Z$ is paracompact and locally of singular shape (Remark \ref{juline}) so the desired result follows from Corollary \ref{twux2}.
\end{proof}

\begin{lemma}\label{joblek}
Let $X$ be a paracompact topological space which is locally of singular shape and equipped with a conical $A$-stratification,
where $A$ satisfies the ascending chain condition. Then the functor
$\Psi_{X}: \Nerve( \bfA^{\degree}_{X}) \rightarrow \Shv(X)$ preserves finite colimits.
\end{lemma}

\begin{proof}
Fix a diagram $p: K \rightarrow \Nerve( \bfA^{\degree}_{X})$ having a colimit $Y$, where
$K$ is finite.
We wish to prove that the induced map $\alpha: \varinjlim(\Psi_{X} \circ p) \rightarrow
\Psi_{X}(Y)$ is an equivalence. Lemma \ref{curin} implies that
$\Psi_{X}(Y) \in \Shv^{A}(X)$, and is therefore hypercomplete (Proposition \ref{toru}). Similarly,
$\varinjlim( \Psi_{X} \circ p)$ is a finite colimit in $\Shv(X)$ of $A$-constructible sheaves, hence
$A$-constructible (Lemma \ref{sogg}) and therefore hypercomplete. Consequently, to
prove that $\alpha$ is an equivalence, it will suffice to show that $\alpha$ is $\infty$-connective.
This condition can be tested pointwise (Lemma \ref{copus}); we may therefore reduce to the problem of showing that $\alpha$ is an equivalence when restricted to each stratum $X_{a}$. Shrinking
$X$ if necessary, we may suppose that $X$ has the form $Z \times C(Y)$ (Lemma \ref{olit})
so that $X_{a} \simeq Z$ is paracompact and locally of singular shape (Remark \ref{juline}). 
Using Lemma \ref{olba} we can replace $X$ by $X_{a}$ and thereby reduce to the case of a trivial stratification. In this case, the functor $\Psi_{X}$ is a fully faithful embedding (Theorem \ref{squ})
those essential image is stable under finite colimits (Corollary \ref{twux2}), and therefore
preserves finite colimits.
\end{proof}

We can use the same argument to prove a sharpened version of Lemma \ref{olba} (at least in case where $A$ satisfies the ascending chain condition):

\begin{proposition}\label{swilker}
Let $A$ be a partially ordered set which satisfies the ascending chain condition, and let
$f: X' \rightarrow X$ be a continuous map between paracompact topological spaces which are locally of singular shape.
Suppose that $X$ is endowed with a conical $A$-stratification, and that the induced
$A$-stratification of $X'$ is also conical. Then the natural transformation $\phi_{X',X}$ is an equivalence of functors from $\Nerve( \bfA^{\degree}_{X})$ to $\Shv(X')$.
\end{proposition}

\begin{lemma}\label{singsong}
Let $X$ be a topological space of singular shape. For every point $x \in X$, there
exists an open neighborhood $U$ of $x$ such that the inclusion
of Kan complexes $\Sing(U) \rightarrow \Sing(X)$ is nullhomotopic.
\end{lemma}

\begin{proof}
Let $K = \Sing(X) \in \SSet$, and let $\pi: X \rightarrow \ast$ denote the projection map.
Since $X$ is of singular shape, there exists a morphism ${\bf 1} \rightarrow \pi^{\ast} K$
in $\Shv(X)$ 
The geometric realization $| \Sing(X) |$ is a CW complex. Since $X$ is of singular shape, composition
with the counit map $v: | \Sing(X) | \rightarrow X$ induces a homotopy equivalence of Kan complexes
$\bHom_{ \Top}( X, | \Sing(X) | ) \rightarrow \bHom_{ \Top}( | \Sing(X) | , | \Sing(X) | )$.
In particular, there exists a continuous map $s: X \rightarrow | \Sing(X) |$ such that
$s \circ v$ is homotopic to the identity. Choose a contractible open subset
$V \subseteq | \Sing(X) |$ containing $s(x)$, and let $U = s^{-1}(V)$. We claim
that the inclusion $i: \Sing(U) \rightarrow \Sing(X)$ is nullhomotopic. To prove this,
it suffices to show that $|i|: | \Sing(U) | \rightarrow | \Sing(X) |$ is nullhomotopic. This map
is homotopic to the composition
$s \circ v \circ |i|$, which factors through the contractible open subset $V \subseteq | \Sing(X) |$.
\end{proof}

\begin{proof}[Proof of Proposition \ref{swilker}]
Let $Y \in \Nerve( \bfA^{\degree}_{X})$, and let $Y' = Y \times_{ \Sing^{A}(X)} \Sing^{A}(X')$. 
We wish to prove that the map $\alpha: f^{\ast} \Psi_{X}(Y) \rightarrow \Psi_{X'}(Y')$ is an equivalence
in $\Shv(X')$. Lemma \ref{curin} implies that $\Psi_{X}(Y) \in \Shv^{A}(X)$, so that
$f^{\ast} \Psi_{X}(Y) \in \Shv^{A}(X')$. Similarly, $\Psi_{X'}(Y') \in \Shv^{A}(X')$, so that
both $f^{\ast} \Psi_{X}(Y)$ and $\Psi_{X'}(Y')$ are hypercomplete (Proposition \ref{toru}).
To prove that $\alpha$ is an equivalence, it will suffice to show that $\alpha$ is $\infty$-connective.
Since this condition can be tested pointwise, it will suffice to show that $\alpha$ induces
an equivalence after restricting to each stratum $X'_{a}$ of $X'$. Using Lemma
\ref{olit} and Remark \ref{juline}, we can shrink $X$ and $X'$ so that $X_{a}$ and
$X'_{a}$ are again paracompact and locally of singular shape. Applying
Lemma \ref{olba}, we can reduce to the case where $X = X_{a}$ and $X' = X'_{a}$.
Shrinking $X$ further (using Lemma \ref{singsong}), we may assume that $Y \simeq \Sing(X) \times K$ for some
Kan complex $K \in \SSet$. In this case, Example \ref{complus} allows us to identify
$\Psi_{X}(Y)$ with the pullback $\pi^{\ast} K$ and $\Psi_{X'}(Y')$ with ${\pi'}^{\ast} K$, where
$\pi: X \rightarrow \ast$ and $\pi': X' \rightarrow \ast$ denote the projection maps. 
Under these identifications, the natural transformation $\phi_{X',X}(Y)$ is induced by the canonical equivalence
$f^{\ast} \circ \pi^{\ast} \simeq (\pi \circ f)^{\ast} = {\pi'}^{\ast}$.
\end{proof}

\begin{proof}[Proof of Theorem \ref{mainstrat}]
We will prove more generally that for every $U \in \calB(X)$, the
functor $\Psi_{U}: \Nerve( \bfA^{\degree}_{U}) \rightarrow \Shv^{A}(U)$ is an equivalence
of $\infty$-categories. The proof proceeds by induction on $\rk(U)$, where the rank
functor $\rk$ is defined in Remark \ref{jult}.

Let $S$ denote the partially ordered set of all open sets $V \in \calB(U)$
which are homeomorphic to a product $Z \times C(Y)$, where $Y$ is an $A_{> a}$-stratified space, and $Z \times C(Y)$ is endowed with the induced
$A_{ \geq a}$-stratification. For every such open set $V$, let $\chi_{V} \in \Shv(X)$ be
the sheaf determined by the formula
$$ \chi_{V}(W) = \begin{cases} \ast & \text{if } W \subseteq V \\
\emptyset & \text{ otherwise.} \end{cases}$$
Let $\alpha$ denote the canonical map $\varinjlim_{V} \chi_V \rightarrow \chi_{U}$.
For each point $x \in U$, the stalk of the colimit $\varinjlim_{V} \chi_{V}$ at $x$
is homotopy equivalent to the nerve of the partially ordered set $S_{x} = \{ V \in S: x \in V \}$.
It follows from Lemma \ref{olit} that the partially ordered set $S_{x}^{op}$ is filtered,
so that $| S_{x} | $ is contractible: consequently, the map $\alpha$ is $\infty$-connective.
Consequently, $\alpha$ induces an equivalence
$\varinjlim_{V} \chi_{V} \rightarrow \chi_{U}$ in the hypercomplete $\infty$-topos $\Shv(X)^{\hyp}$. 
Applying Theorem \toposref{charleschar} to the $\infty$-topos $\Shv(X)^{\hyp}$, we conclude
that $\Shv(U)^{\hyp} \simeq \Shv(X)^{\hyp}_{/ \chi_U}$ is equivalent to the homotopy limit of the diagram
of $\infty$-categories $\{ \Shv(V)^{\hyp} \simeq \Shv(X)^{\hyp}_{/ \chi_V} \}_{V \in S}$.
Proposition \ref{toru} guarantees that $\Shv^{A}(U) \subseteq \Shv(U)^{\hyp}$
(and similarly $\Shv^{A}(V) \subseteq \Shv(V)^{\hyp})$ for each $V \in S$).
Since the property of being constructible can be tested locally, we obtain an equivalence
$$ \Shv^{A}(U) \simeq \varprojlim \{ \Shv^{A}(V) \}_{V \in S}.$$

We next show that the restriction maps $\Nerve( \bfA^{\degree}_{U} )
\rightarrow \Nerve( \bfA^{\degree}_{V} )$ exhibit $\Nerve( \bfA^{\degree}_{U})$ as the homotopy
limit of the diagram of $\infty$-categories $\{ \Nerve( \bfA^{\degree}_{V} ) \}_{V \in S}$.
In view of the natural equivalences
$$ \Fun( \Sing^{A}(V), \SSet) \leftarrow \Nerve(( \Set_{\Delta}^{\sCoNerve[ \Sing^{A}(V)]})^{\degree})
\rightarrow \Nerve( \bfA^{\degree}_{V}),$$
it will suffice to show that the canonical map
$$ \Fun( \Sing^{A}(U), \SSet) \rightarrow \varprojlim \{ \Fun( \Sing^{A}(V), \SSet) \}_{V \in S}$$
is an equivalence. This follows immediately from Theorem \ref{superkamp}.

We have a commutative diagram
$$ \xymatrix{ \Nerve( \bfA^{\degree}_{U} ) \ar[r] \ar[d] & \varprojlim_{V \in S} \Nerve( \bfA^{\degree}_{V} ) \ar[d] \\
 \Shv^{A}(U) \ar[r] & \varprojlim_{V \in S} \Shv^{A}(V) }$$
where the vertical maps are equivalences. Consequently, to prove that $\Psi_{U}$ is an equivalence,
it will suffice to show that $\Psi_{V}$ is an equivalence for each $V \in S$. Replacing $U$ by
$V$, we can assume that $U$ has the form $Z \times C(Y)$. We will also assume that $Z$ is nonempty (otherwise there is nothing to prove).

Let $U' = Z \times Y \times (0, \infty)$, which we regard as an open subset of $U$.
Let $\calC_0 \subseteq \Nerve( \bfA^{\degree}_{U})$ be the full subcategory spanned by
the left fibrations $Y \rightarrow \Sing^{A}(U)$ which factor through $\Sing^{A}(U')$, and let
$\calC_1 \subseteq \Nerve( \bfA^{\degree}_{U})$ be the full subcategory spanned
by the Kan fibrations $Y \rightarrow \Sing^{A}(U)$. Under the equivalence
$\Nerve( \bfA^{\degree}_{U}) \simeq \Fun( \Sing^{A}(U), \SSet)$, these
correspond to the full subcategories described in Proposition \ref{soake} (where
$\calA = \SSet$ and $p: \Sing^{A}(U) \rightarrow \Delta^1$ is characterized by
the requirements that $p^{-1} \{0\} = \Sing(U_a)$ and $p^{-1} \{1\} = \Sing^{A}( U')$). 
It follows that the inclusions $\calC_0, \calC_1 \subseteq \Nerve( \bfA^{\degree}_{U})$
admit right adjoints $L_0$ and $L_1$, and that $L_0$ is complementary to $L_1$.
Let $\calC'_0, \calC'_{1} \subseteq \Shv^{A}(U)$ be defined as in Lemma \ref{sou}, so that 
we again have right adjoints $L'_{0}: \Shv^{A}(U) \rightarrow \calC'_{0}$ and
$L'_{1}: \Shv^{A}(U) \rightarrow \calC'_{1}$ which are complementary.
We will prove that the functor $\Psi_{U}$ is an equivalence of $\infty$-categories by
verifying the hypotheses of Proposition \ref{inkor}:

\begin{itemize}
\item[$(2)$] The functor $\Psi_{U}$ restricts to an equivalence $\calC_0 \rightarrow \calC'_0$.
Let $Y \rightarrow \Sing^{A}(U')$ be an object of $\calC_0$. Then
$(\Psi_{U}(Y))(W)$ is empty if $W$ is not contained in $U'$, so that $\Psi_{U}(Y) \in \calC'_0$.
Moreover, the composition of $\Psi_{U} | \calC_0$ with the equivalence
$\calC'_0 \simeq \Shv^{A}(U')$ coincides with the functor $\Psi_{U'}$. Since
the strata $U'_{b}$ are empty unless $b > a$, while $U_{a}$ is nonempty (since $Z \neq \emptyset$),
we have $\rk(U') < \rk(U)$ so that $\Psi_{U'}$ is an equivalence of $\infty$-categories by the inductive hypothesis.

\item[$(2')$] We must show that the functor $\Psi_{U}$ restricts to an equivalence $\calC_1 \rightarrow \calC'_{1}$. Let $\pi: U \rightarrow Z$ denote the projection. We have a diagram of $\infty$-categories
$$ \xymatrix{ \Nerve( \bfA_{Z}^{\degree}) \ar[r] \ar[d]^{\Psi_{Z}} & \Nerve( \bfA_{U}^{\degree}) \ar[d]^{ \Psi_{U}} \\
\Shv(Z) \ar[r]^{\pi^{\ast}} & \Shv(U) }$$
which commutes up to homotopy (Proposition \ref{swilker}). The upper horizontal arrow is
fully faithful, and its essential image is precisely the $\infty$-category $\calC_{1}$.
Consequently, it suffices to show that the composite map $\pi^{\ast} \Psi_{Z}$
is a fully faithful embedding whose essential image is precisely $\calC'_{1}$.
Theorem \ref{squ} implies that $\Psi_{Z}$ is fully faithful, and that its essential image
is the full subcategory of $\Shv(Z)$ spanned by the locally constant sheaves.
The desired result now follows from the definition of $\calC'_{1}$.

\item[$(3)$] The functor $\Psi_{U}$ preserves pushouts. This follows from Lemma \ref{joblek}.

\item[$(4)$] If $\alpha: Y_0 \rightarrow Y$ is a morphism which exhibits
$Y_0 \in \calC_0$ as a $\calC_0$-colocalization of $Y \in \calC_1$, then
$\Psi_{U}(\alpha)$ exhibits $\Psi_{U}(Y_0)$ as a $\calC'_0$-localization
of $\Psi_{U}(Y_1)$. Unwinding the definitions, $\alpha$ induces an equivalence
of left fibrations $Y_0 \rightarrow Y \times_{ \Sing^{A}(U)} \Sing^{A}(U')$, and we must
show that for each $W \in \calB( U')$ that the induced map
$\Fun_{ \Sing^{A}(U)}( \Sing^{A}(W), Y_0 ) \rightarrow \Fun_{ \Sing^{A}(U)}( \Sing^{A}(W), Y)$
is a homotopy equivalence. This is clear, since the condition that $W \subseteq U'$ guarantees
that any map $\Sing^{A}(W) \rightarrow Y$ factors uniquely through the fiber product
$Y \times_{ \Sing^{A}(U)} \Sing^{A}(U')$.
\end{itemize}
\end{proof}

\subsection{Embeddings of Topological Manifolds}\label{cowpy}

Let $M$ and $N$ be topological manifolds of the same dimension.
We let $\Emb(M, N)$ denote the set of all open embeddings $M \hookrightarrow N$.
We will regard $\Emb(M,N)$ as a topological space: it is a subspace of the collection of all
continuous maps from $M$ to $N$, which we endow with the compact-open topology.
We let $\Homeo(M,N)$ denote the set of all homeomorphisms of $M$ with $N$, regarded as a subspace of $\Emb(M,N)$. For $k \geq 0$, we let $\Top(k)$ denote the topological group $\Homeo( \R^{k}, \R^{k})$ of homeomorphisms from $\R^{k}$ to itself.

\begin{warning}
We will regard Convention \ref{koson} as in force throughout this section: the word {\it manifold} will always refer to a paracompact, Hausdorff, topological manifold of some fixed dimension $k$.
\end{warning}

We begin by recalling some classical results from point-set topology. The following is a parametrized version of Brouwer's invariance of domain theorem (for the reader's convenience, we reproduce a proof at the end of this section).

\begin{theorem}[Brouwer]\label{lowcu}
Let $M$ and $N$ be manifolds of dimension $k$, and let $S$ be an arbitrary topological space. Suppose we are given a continuous map $f: M \times X \rightarrow N \times X$
satisfying the following pair of conditions:
\begin{itemize}
\item[$(i)$] The diagram
$$ \xymatrix{ M \times X \ar[rr]^{f} \ar[dr] & & N \times X \ar[dl] \\
& X & }$$
is commutative.
\item[$(ii)$] The map $f$ is injective.
\end{itemize}
Then $f$ is an open map.
\end{theorem}

\begin{corollary}
Let $M$ and $N$ be manifolds of the same dimension, and let $f: M \times X \rightarrow N \times X$
be a continuous bijection which commutes with the projection to $X$. Then $f$ is a homeomorphism.
\end{corollary}

\begin{remark}\label{kum}
Let $M$ and $N$ be topological manifolds of the same dimension, and let $\bHom(M,N)$ denote the set of all continuous maps from $M$ to $N$, endowed with the compact-open topology. Since 
$M$ is locally compact, $\bHom(M,N)$ classifies maps of topological spaces from $M$ to $N$:
that is, for any topological space $X$, giving a continuous map $X \rightarrow \bHom(M,N)$ is equivalent to giving a continuous map $M \times X \rightarrow N$, which is in turn equivalent to giving a commutative diagram
$$ \xymatrix{ M \times X \ar[rr]^{f} \ar[dr] & & N \times X \ar[dl] \\
& X. & }$$
Under this equivalence, continuous maps from $X$ to $\Emb(M,N)$ correspond to commutative diagrams as above where $f$ is injective (hence an open embedding, by Theorem \ref{lowcu}), and
continuous maps from $X$ to $\Homeo(M,N)$ correspond to commutative diagrams as above where
$f$ is bijective (and therefore a homeomorphism, by Theorem \ref{lowcu}). It follows that
the space of embeddings $\Emb(M,M)$ has the structure of a topological monoid, and that
$\Homeo(M,M)$ has the structure of a topological group.
\end{remark}

In \S \ref{colee}, we need the to know that the topological monoid $\Emb( \R^{k}, \R^{k})$ is grouplike: that is, the set of path components $\pi_0 \Emb( \R^{k}, \R^{k})$ forms a group under composition. This is an immediate consequence of the following version of the Kister-Mazur theorem (again, we reproduce a proof at the end of this section):

\begin{theorem}[Kister-Mazur]\label{scen}
For each $k \geq 0$, the inclusion $\Top(k) \hookrightarrow \Emb( \R^{k}, \R^{k})$ is
a homotopy equivalence.
\end{theorem}

We now describe some variants on the embedding spaces $\Emb(M,N)$ and their homotopy types.

\begin{definition}
Let $M$ be a topological manifold of dimension $k$, let $S$ be a finite set, and for every positive
real number $t$ let $B(t) \subset \R^{k}$ be as in Lemma \ref{cookra}. We let
$\Germ(S, M)$ denote the simplicial set
$\varinjlim_{n} \Sing \Emb( B( \frac{1}{2^n}) \times S, M).$
We will refer to $\Germ(S,M)$ as the {\it simplicial set of $S$-germs in $M$}.
\end{definition}

\begin{lemma}\label{cookra}
Let $M$ be a topological manifold of dimension $k$ and $S$ a finite set. For every
positive real number $t$, let $B(t) \subset \R^{k}$ denote the open ball of radius $t$. 
For every pair of positive real numbers $s < t$, the restriction map
$r: \Emb( B(t) \times S, M) \rightarrow \Emb(B(s) \times S, M)$ is a homotopy
equivalence.
\end{lemma}

\begin{proof}
This follows from the observation that the embedding $B(s) \hookrightarrow B(t)$ is isotopic
to a homeomorphism.
\end{proof}

By repeated application of Lemma \ref{cookra} we deduce the following:

\begin{proposition}\label{splay}
Let $M$ be a topological manifold of dimension $k$ and let $S$ be a finite set. Then
the obvious restriction map
$\Sing \Emb( \R^{k} \times S, M) \rightarrow \Germ(S, M)$ is a homotopy equivalence of
Kan complexes.
\end{proposition}

\begin{notation}
Let $M$ be a topological manifold of dimension $k$. Evaluation at the origin
$0 \in \R^{k}$ induces a map $\theta: \Emb( \R^{k}, M) \rightarrow M$. We will denote the fiber
of this map over a point $x \in M$ by $\Emb_{x}( \R^{k}, M)$. The map $\theta$ is a Serre fibration, so we have a fiber sequence of topological spaces
$$ \Emb_{x}( \R^{k},M) \rightarrow \Emb( \R^{k}, M) \rightarrow M.$$

We let $\Germ(M)$ denote the simplicial set $\Germ( \{ \ast \}, M)$. Evaluation at $0$
induces a Kan fibration $\Germ(M) \rightarrow \Sing M$; we will denote the fiber of this map
over a point $x \in M$ by $\Germ_{x}(M)$. We have a map of fiber sequences
$$ \xymatrix{ \Sing \Emb_{x}( \R^{k}, M) \ar[r] \ar[d]^{\psi} & \Sing \Emb(\R^{k}, M) \ar[r] \ar[d]^{\psi'} & \Sing M \ar[d]^{\psi''} \\
\Germ_{x}(M) \ar[r] & \Germ(M) \ar[r] & \Sing M. }$$
Since $\psi'$ is a homotopy equivalence (Proposition \ref{splay}) and $\psi''$ is an isomorphism, we conclude that $\psi$ is a homotopy equivalence.

The simplicial set $\Germ_0( \R^{k})$ forms a simplicial group with respect to the operation of composition of germs. Since $\R^{k}$ is contractible, we have homotopy equivalences of simplicial
monoids
$$ \Germ_0(\R^{k}) \leftarrow \Sing \Emb_{x}( \R^{k},\R^{k}) \rightarrow \Sing \Emb( \R^{k}, \R^{k})
\leftarrow \Sing \Top(k)$$
(see Theorem \ref{scen}): in other words, $\Germ_0( \R^{k})$ can be regarded as a model
for the homotopy type of the topological group $\Top(k)$. 
\end{notation}

\begin{remark}\label{geck}
For any topological $k$-manifold $M$, the group $\Germ_0( \R^{k})$ acts on
$\Germ(M)$ by composition. This action is free, and we have a canonical
isomorphism of simplicial sets $\Germ(M) / \Germ_0( \R^{k}) \simeq \Sing M$. 
\end{remark}

\begin{remark}\label{gec}
Let $j: U \rightarrow M$ be an open embedding of topological $k$-manifolds and $S$ a finite set.
Then evaluation at $0$ determines a diagram of simplicial sets
$$ \xymatrix{ \Sing \Emb( \R^{k} \times S, U) \ar[r] \ar[d] & \Sing \Emb( \R^{k} \times S, M) \ar[d] \\
\Conf(S,U) \ar[r] & \Conf(S, M). }$$
We claim that this diagram is homotopy Cartesian. In view of Proposition \ref{splay}, it suffices
to show that the equivalent diagram
$$ \xymatrix{ \Germ(S,U) \ar[r] \ar[d] & \Germ(S,M) \ar[d] \\
\Conf(S,U) \ar[r] & \Conf(S,M), }$$
is homotopy Cartesian. This diagram is a pullback square and the vertical maps are Kan fibrations:
in fact, the vertical maps are principal fibrations with structure group $\Germ_0(\R^{k})^{S}$. 

Taking $U = \R^{k}$, and $S$ to consist of a single point, we have a larger diagram
$$ \xymatrix{ \Sing \Emb_0( \R^{k}, \R^{k}) \ar[r] \ar[d] & \Sing \Emb( \R^{k}, \R^{k}) \ar[r] \ar[d] & \Sing \Emb( \R^{k}, M) \ar[d] \\
\{0\} \ar[r] & \Sing \R^{k} \ar[r] & \Sing M. }$$
Since the horizontal maps on the left are homotopy equivalences of Kan complexes, we obtain
a homotopy fiber sequence of Kan complexes
$$ \Sing \Emb_0( \R^{k}, \R^{k} ) \rightarrow \Sing \Emb( \R^{k}, M) \rightarrow \Sing M.$$
\end{remark}

We conclude this section with the proofs of Theorems \ref{lowcu} and \ref{scen}. 

\begin{proof}[Proof of Theorem \ref{lowcu}]
Fix a continuous map $f: M \times X \rightarrow N \times X$ and an open set $U \subseteq M \times X$; we wish to show that $f(U)$ is open in $N \times X$.
In other words, we wish to show that for each $u = (m,x) \in U$, the set $f(U)$ contains a neighborhood of $f(u)=(n,x)$ in $N \times S$. Since $N$ is a manifold, there exists an open neighborhood
$V \subseteq N$ containing $n$ which is homeomorphic to Euclidean space $\R^{k}$.
Replacing $N$ by $V$ (and shrinking $M$ and $X$ as necessary), we may assume that
$N \simeq \R^{k}$. Similarly, we can replace $M$ and $X$ by small neighborhoods of $m$ and
$s$ to reduce to the case where $M \simeq \R^{k}$ and $U = M \times X$. 

We first treat the case where $X$ consists of a single point. Let $D \subseteq M$ be a closed neighborhood of $m$ homeomorphic to a (closed) $k$-dimensional disk, and regard $N$ as
an open subset of the $k$-sphere $S^{k}$. We have a long exact sequence of compactly supported cohomology groups
$$ 0 \simeq \HH^{k-1}_{c}( S^{k}; \Z) \rightarrow \HH^{k-1}_{c}( f( \bd D); \Z) \rightarrow \HH^{k}_c( S^{k} - f( \bd D); \Z) \rightarrow \HH^{k}_c(S^{k}; \Z) \rightarrow \HH^{k}_{c}( f(\bd D); \Z) \simeq 0.$$
Since $f$ is injective, $f( \bd D)$ is homeomorphic to a $(k-1)$-sphere.
It follows that $\HH^{k}_{c}( S^k- f(\bd D); \Z)$ is a free $\Z$-module of rank $2$, so that
(by Poincare duality) the ordinary cohomology $\HH^{0}( S^{k}- f(\bd D); \Z)$ is also free of rank
$2$: in other words, the open set $S^{k} - f(\bd D)$ has exactly two connected components.
We have another long exact sequence
$$ 0 \simeq \HH^{k-1}_{c}( f(D) ; \Z) \rightarrow \HH^{k}_c( S^{k} - f(D) ; \Z) \rightarrow \HH^{k}_c(S^{k}; \Z) \rightarrow \HH^{k}_{c}( f(D); \Z) \simeq 0.$$
This proves that $\HH^{k}_{c}( S^k - f(D); \Z)$ is free of rank $1$ so that (by Poincare duality)
$S^{k} - f(D)$ is connected. The set $S^{k} - f( \bd D)$ can be written as a union of connected sets
$f( D - \bd D)$ and $S^{k} - f(D)$, which must therefore be the connected components of
$S^{k} - f( \bd D)$. It follows that $f( D - \bd D)$ is open $S^{k}$ so that $f(M)$ contains
a neighborhood of $f(m)$ as desired.

Let us now treat the general case. Without loss of generality, we may assume that
$f(u) = (0, x)$, where $x \in X$ and $0$ denotes the origin of $\R^{k}$.
Let $f_{x}: M \rightarrow N$ be the restriction of $f$ to $M \times \{x\}$. The above argument
shows that $f_{x}$ is an open map, so that $f_x(M)$ contains a closed ball $\overline{B(\epsilon)} \subseteq \R^{k}$ for some positive radius $\epsilon$. Let $S \subseteq M - \{m\}$ be the inverse
image of the boundary $\bd \overline{ B(\epsilon)}$, so that $S$ is homeomorphic to the
$(k-1)$-sphere. In particular, $S$ is compact. Let $\pi: M \times X \rightarrow \R^{k}$ denote the composition of $f$ with the projection map $N \times X \rightarrow N \simeq \R^{k}$.
Shrinking $X$ if necessary, we may suppose that the distance $d( f(s,x), f(s,y) ) < \frac{ \epsilon}{2}$ for 
all $s \in S$ and all $y \in X$. We will complete the proof by showing that $B( \frac{\epsilon}{2}) \times X$ is contained in the image of $f$. Supposing otherwise; then there exists $v \in B( \frac{\epsilon}{2})$
and $y \in X$ such that $(v,y) \notin f( M \times X)$. Then $f_{y}$ defines a map from $M$ to
$\R^{k} - \{v\}$, so the restriction $f_{y} | S$ is nullhomotopic when regarded as a map from
$S$ to $\R^{k} - \{v\}$. However, this map is homotopic (via a straight-line homotopy) to
$f_{x} | S$, which carries $S$ homeomorphically onto $\bd \overline{ B( \epsilon)} \subseteq
\R^{k} - \{v\}$. It follows that the inclusion $\bd \overline{ B( \epsilon)} \subseteq \R^{k} - \{v\}$ is nullhomotopic, which is impossible.
\end{proof}

We now turn to the proof of Theorem \ref{scen}. The main step is the following technical result:

\begin{lemma}\label{prescen}
Let $X$ be a paracompact topological space, and suppose that there exists a continuous map
$f_0: \R^{k} \times X \rightarrow \R^{k}$ such that, for each $x \in X$, the restriction
$f_{0,x} = f_0 | \R^{k} \times \{x\}$ is injective. Then there exists an isotopy
$f: \R^{k} \times X \times [0,1] \rightarrow \R^{k}$ with the following properties:
\begin{itemize}
\item[$(i)$] The restriction $f | \R^{k} \times X \times \{0\}$ coincides with $f_0$.
\item[$(ii)$] For every pair $(x,t) \in X \times [0,1]$, the restricted map
$f_{t,x}= f| \R^{k} \times \{x\} \times \{t\}$ is injective.
\item[$(iii)$] For each $x \in X$, the map $f_{1,x}$ is bijective.
\item[$(iv)$] Suppose $x \in X$ has the property that $f_{0,x}$ is bijective.
Then $f_{t,x}$ is bijective for all $t \in [0,1]$.
\end{itemize}
\end{lemma}

\begin{proof}
Let $w: X \rightarrow \R^{k}$ be given by the formula $w(x) = f_0(0,x)$. Replacing
$f_0$ by the map $(v,x) \mapsto f_0(v,x) - w(x)$, we can reduce to the case where
$w=0$: that is, each of the maps $f_{0,x}$ carries the origin of $\R^{k}$ to itself. 

For every continuous positive real-valued function $\epsilon: X \rightarrow \R_{> 0}$, we let
$B(\epsilon)$ denote the open subset of $\R^{k} \times X$ consisting of those pairs
$(v,x)$ such that $|v| < \epsilon(x)$. If $r$ is a real number, we let $B(r) = B(\epsilon)$, where
$\epsilon: X \rightarrow \R_{ > 0}$ is the constant function taking the value $r$.

Let $g^1: \R^{k} \times X \rightarrow \R^{k} \times X$ be given by the formula
$g^{1}(v,x) = (f_1(v,x),x)$. The image $g^1( B(1))$ is an open subset
of $\R^{k} \times X$ (Theorem \ref{lowcu}) which
contains the zero section $\{0\} \times X$; it follows that $g^1(B(1))$ 
contains $B( \epsilon )$ for some positive real-valued continuous function
$\epsilon: X \rightarrow \R_{>0}$. Replacing $f_0$ by the funciton
$(v,x) \mapsto \frac{ f_0(v,x)}{\epsilon(x)}$, we can assume that $B(1) \subseteq
g^1( B(1) )$. 

We now proceed by defining a sequence of open embeddings
$\{ g^{i}: \R^{k} \times X \rightarrow \R^{k} \times X \}_{i \geq 2}$ and isotopies
$\{ h^{i}_{t} \}_{0 \leq t \leq 1}$ from $g^{i}$ to $g^{i+1}$, so that the following conditions are satisfied:
\begin{itemize}
\item[$(a)$] Each of the maps $g^{i}$ is compatible with the projection to $X$.
\item[$(b)$] Each isotopy $\{ h^{i}_{t} \}_{0 \leq t \leq 1}$ consists of open embeddings
$\R^{k} \times X \rightarrow \R^{k} \times X$ which are compatible with the projection to $X$.
Moreover, this isotopy is constant on the open set $B(i) \subseteq \R^{k} \times X$.
\item[$(c)$] For $i \geq 1$, we have $B(i) \subseteq g^{i}( B(i) )$. 
\item[$(d)$] Let $x \in X$ be such that the map
$g^{i}_{x}: \R^{k} \rightarrow \R^{k}$ is a homeomorphism. Then
$h^{i}_{t,x}: \R^{k} \rightarrow \R^{k}$ is a homeomorphism for all $t \in [0,1]$.
\end{itemize}

Assuming that these requirements are met, we can obtain the desired isotopy $f_{t}$ by
the formula
$$ f_{t}(v,x) = \begin{cases} \pi g^{i}(v,x) & \text{ if } (|v| < i) \wedge ( t > 1 - \frac{1}{2^{i-1}}) \\
\pi h^{i}_{s}(v,x) & \text{ if } t = 1 + \frac{s-2}{2^{i}}, \end{cases}$$
where $\pi$ denotes the projection from $\R^{k} \times X$ onto $\R^{k}$. We now proceed by induction on $i$. Assume that $g^{i}$ has already been constructed; we will construct an isotopy
$h^{i}$ from $g^{i}$ to another open embedding $g^{i+1}$ to satisfy the above conditions.
First, we need to establish a bit of notation.

For every pair of real numbers $r < s$, let $\{ H(r,s)_{t}: \R^{k} \rightarrow \R^{k} \}_{0 \leq t \leq 1}$ be
a continuous family of homeomorphisms satisfying the following conditions:
\begin{itemize}
\item[$(i)$] The isotopy $\{ H(r,s)_{t} \}$ is constant on $\{ v \in \R^k: |v| < \frac{r}{2} \}$ and
$\{ v \in \R^k: |v| > s+1 \}$.
\item[$(ii)$] The map $H(r,s)$ restricts to a homeomorphism of $B(r)$ with $B(s)$.
\end{itemize}
We will assume that the homeomorphisms $\{ H(r,s)_{t} \}$ are chosen to depend continuously
on $r$, $s$, and $t$. Consequently, if $\epsilon < \epsilon'$ are positive real-valued functions on
$X$, we obtain an isotopy $\{ H(\epsilon, \epsilon')_{t}: \R^{k} \times X \rightarrow
\R^{k} \times X \}$ by the formula $H(\epsilon, \epsilon')_{t}( v, x) =
( H(\epsilon(x), \epsilon'(x))_{t}(v), x)$.

Since $g^i$ is continuous and $\{0\} \times X \subseteq (g^i)^{-1} B( \frac{1}{2})$, there
exists a real-valued function $\delta: X \rightarrow (0,1)$ such that
$g^i( B( \delta) ) \subseteq B( \frac{1}{2})$. We define a homeomorphism $c: \R^{k} \times X \rightarrow \R^{k} \times X$ as follows:
$$ c(v,x) = \begin{cases} (v,x) & \text{if } (v,x) \notin g^{i} ( \R^{k} \times X ) \\
g^i( H( \delta(x), i)_{1}^{-1}(w),x) & \text{ if } (v,x) = g^{i}(w,x). \end{cases}$$
Since $g^i$ carries $B( \delta)$ into $B( \frac{1}{2} )$, we
deduce that $c( g^{i}(v,x) ) \in B( \frac{1}{2} )$ if $(v,x) \in B(i)$. 
Note that $c$ is the identity outside of the image $g^{i} B(i+1)$; we can therefore
choose a positive real valued function $\epsilon: X \rightarrow (i+1, \infty)$ such that
$c$ is the identity outside of $B(\epsilon)$.

We now define $h^{i}_{t}$ by the formula
$h^{i}_{t} = c^{-1} \circ H(1, \epsilon)_{t} \circ c \circ g^{i}$ (here we
identify the real number $1 \in \R$ with the constant function
$X \rightarrow \R$ taking the value $1$). 
It is clear that $h^{i}_{t}$ is an isotopy from $g^{i} = g^{i}_{0}$ to another
map $g^{i+1} = g^{i}_{1}$, satisfying conditions $(a)$ and $(d)$ above.
Since $H(1,\epsilon)_{t}$ is
the identity on $B( \frac{1}{2})$ and $c \circ g^{i}$ carries
$B(i)$ into $B( \frac{1}{2})$, we deduce that $h^{i}_{t}$ is constant
on $B(i)$ so that $(b)$ is satisfied. It remains only to verify
$(c)$: we must show that $g^{i+1} B(i+1)$ contains $B(i+1)$.
In fact, we claim that $g^{i+1} B(i+1)$ contains $B(\epsilon)$. Since
$c$ is supported in $B(\epsilon)$, it suffices to show that
$(c g^{i+1})B(i+1) = (H( 1, \epsilon)_{1} \circ c \circ g^{i}) B(i+1)$
contains $B(\epsilon)$. For this we need only show that
$(c \circ g^{i}) B(i+1)$ contains $B(1) \subseteq B(i) \subseteq g^{i} B(i) \subseteq g^{i} B(i+1)$.
This is clear, since $H( \delta(x), i)_{1}$ induces a homeomorphism of
$B(i+1)$ with itself.
\end{proof}

\begin{proof}[Proof of Theorem \ref{scen}]
For every compact set $K \subseteq \R^{k}$, the compact open topology on the set
of continuous maps $\bHom(K, \R^{k})$ agrees with the topology induced by the metric
$d_K(f,g) = \sup \{ | f(v) - g(v) |, v \in K \}$. Consequently, the compact open topology on
the entire mapping space $\bHom ( \R^{k}, \R^{k})$ is defined by the countable sequence of metrics
$\{ d_{ \overline{B(n)}} \}_{n \geq 0}$ (here $\overline{ B(n)}$ denotes the closed ball of
radius $n$), or equivalently by the single metric
$$ d(f,g) = \sum_{n \geq 0} \frac{1}{2^n} \inf \{ 1, d_{ \overline{B(n)}}(f,g) \}.$$
It follows that $\Emb( \R^{k}, \R^{k}) \subseteq \bHom( \R^{k}, \R^{k})$ is metrizable
and therefore paracompact. Applying Lemma \ref{prescen} to the canonical pairing
$$f_0: \R^{k} \times \Emb( \R^{k}, \R^{k}) \hookrightarrow
\R^{k} \times \bHom( \R^{k}, \R^{k}) \rightarrow \R^{k},$$
we deduce the existence of an map $f: \R{k} \times \Emb( \R^{k}, \R^{k}) \times [0,1] \rightarrow \R^{k}$ which is classified by a homotopy $\chi: \Emb( \R^{k}, \R^{k}) \times [0,1] \rightarrow \Emb( \R^{k}, \R^{k})$ from $\id_{ \Emb( \R^{k}, \R^{k})}$ to some map $s: \Emb( \R^{k}, \R^{k}) \rightarrow
\Homeo( \R^{k}, \R^{k})$. We claim that $s$ is a homotopy inverse to the inclusion
$i: \Homeo(\R^{k}, \R^{k}) \rightarrow \Emb( \R^{k}, \R^{k})$. The homotopy
$\chi$ shows that $i \circ s$ is homotopy to the identity on $\Emb( \R^{k}, \R^{k})$, and
the restriction of $\chi$ to $\Homeo( \R^{k}, \R^{k}) \times [0,1]$ shows that
$s \circ i$ is homotopic to the identity on $\Homeo( \R^{k}, \R^{k})$.
\end{proof}

\subsection{Verdier Duality}\label{kopo}

Our goal in this section is to prove the following result:

\begin{theorem}[Verdier Duality]\label{vedu}
Let $\calC$ be a stable $\infty$-category which admits small limits and colimits, and let $X$ be a locally compact topological space. There is a canonical equivalence of $\infty$-categories
$$ \mathbb{D}: \Shv( X; \calC )^{op} \simeq \Shv(X; \calC^{op}).$$
\end{theorem}

\begin{remark}
Let $k$ be a field and let $\bfA$ denote the category of chain complexes of $k$-vector spaces. Then
$\bfA$ has the structure of a simplicial category; we let $\calC = \Nerve(\bfA)$ denote the nerve of $\bfA$
(that is, the derived $\infty$-category of the abelian category of $k$-vector spaces; see Definition \stableref{smucky}). Vector space duality induces a simplicial functor $\bfA^{op} \rightarrow \bfA$, which
in turn gives rise to a functor $\calC^{op} \rightarrow \calC$. This functor preserves limits, and therefore
induces a functor $\Shv( X; \calC^{op}) \rightarrow \Shv(X; \calC)$ for any locally compact topological space $X$. Composing this map with the equivalence $\mathbb{D}$ of Theorem \ref{vedu}, we obtain a functor $\mathbb{D}': \Shv( X; \calC)^{op} \rightarrow \Shv(X,\calC)$: that is, a contravariant functor from $\Shv(X; \calC)$ to itself. It is the functor $\mathbb{D}'$ which is usually called Verdier duality. Note that $\mathbb{D}'$ is not an equivalence of $\infty$-categories: it is obtained by composing the equivalence $\mathbb{D}$ with vector space duality, which fails to be an equivalence unless suitable finiteness restrictions are imposed.
\end{remark}

The first step in the proof of Theorem \ref{vedu} is to choose a convenient model for the
$\infty$-category $\Shv(X; \calC)$ of $\calC$-valued sheaves on $X$. Let $\calK(X)$ denote the collection of all compact subsets of $X$, regarded as a partially ordered set with respect to inclusion. Recall (Definition \toposref{labsos}) that a {\it $\calK$-sheaf} on $X$ (with values in an $\infty$-category $\calC$) is a functor
$\calF: \Nerve( \calK(X) )^{op} \rightarrow \calC$ with the following properties:
\begin{itemize}
\item[$(i)$] The object $\calF(\emptyset) \in \calC$ is final.
\item[$(ii)$] For every pair of compact sets $K, K' \subseteq X$, the diagram
$$ \xymatrix{ \calF( K \cup K' ) \ar[r] \ar[d] & \calF(K) \ar[d] \\
\calF(K') \ar[r] & \calF(K \cap K') }$$
is a pullback square in $\calC$.
\item[$(iii)$] For every compact set $K \subseteq X$, the canonical map
$\varinjlim_{K'} \calF(K') \rightarrow \calF(K)$ is an equivalence, where $K'$ ranges over
all compact subsets of $X$ which contain a neighborhood of $K$.
\end{itemize}
We let $\Shv_{\calK}(X; \calC)$ denote the full subcategory of $\Fun( \Nerve( \calK(X)^{op}, \calC)$ spanned by the $\calK$-sheaves. We now have the following:

\begin{lemma}\label{boffly}
Let $X$ be a locally compact topological space and $\calC$ a stable $\infty$-category which admits small limits and colimits. Then there is a canonical equivalence of $\infty$-categories $\Shv(X; \calC) \simeq \Shv_{\calK}(X; \calC)$. 
\end{lemma}

\begin{proof}
Since $\calC$ is stable, filtered colimits in $\calC$ are left exact. The desired result is now a consequence of Theorem \toposref{kuku} (note that Theorem \toposref{kuku} is stated under the hypothesis that $\calC$ is presentable, but this hypothesis is used only to guarantee the existence of small limits and colimits in $\calC$).
\end{proof}

Using Lemma \ref{boffly}, we can reformulate Theorem \ref{vedu} as follows:

\begin{theorem}\label{vedu2}
Let $X$ be a locally compact topological space and let $\calC$ be a stable $\infty$-category which admits small limits and colimits. Then there is a canonical equivalence of $\infty$-categories
$$ \Shv_{\calK}( X; \calC)^{op} \simeq \Shv_{\calK}( X; \calC^{op} ).$$
\end{theorem}

We will prove Theorem \ref{vedu2} by introducing an $\infty$-category which
is equivalent to both $\Shv_{\calK}(X; \calC)^{op}$ and $\Shv_{\calK}(X; \calC^{op} )$.

\begin{notation}\label{buxy}
Fix a locally compact topological space $X$. We define a partially ordered set $M$ as follows:
\begin{itemize}
\item[$(1)$] The objects of $M$ are pairs $(i,S)$ where $0 \leq i \leq 2$ and $S$ is a subset of $X$ such that
$S$ is compact if $i=0$ and $X-S$ is compact if $i=2$.
\item[$(2)$] We have $(i, S) \leq (j, T)$ if either $i \leq j$ and $S \subseteq T$, or $i = 0$ and $j = 2$.
\end{itemize}
\end{notation}

\begin{remark}
The projection $(i,S) \mapsto i$ determines a map of partially ordered sets $\phi: M \rightarrow [2]$. For
$0 \leq i \leq 2$, we let $M_i$ denote the fiber $\phi^{-1} \{i\}$. We have canonical isomorphisms
$M_0 \simeq \calK(X)$ and $M_2 \simeq \calK(X)^{op}$, while $M_1$ can be identified with the partially ordered set of {\em all} subsets of $X$. 
\end{remark}

The proof of Theorem \ref{vedu2} rests on the following:

\begin{proposition}\label{cobine}
Let $X$ be a locally compact topological space, $\calC$ a stable $\infty$-category which admits small limits and colimits, and let $M$ be the partially ordered set of Notation \ref{buxy}. Let $F: \Nerve(M) \rightarrow \calC$ be a functor. The following conditions are equivalent:
\begin{itemize}
\item[$(1)$] The restriction $( F | \Nerve(M_0) )^{op}$ determines a $\calK$-sheaf
$\Nerve( \calK(X) )^{op} \rightarrow \calC^{op}$, the restriction $F| \Nerve( M_1)$ is zero, 
and $F$ is a left Kan extension of the restriction $F | \Nerve( M_0 \cup M_1)$.
\item[$(2)$] The restriction $F | \Nerve(M_2)$ determines a $\calK$-sheaf
$\Nerve( \calK(X))^{op} \rightarrow \calC$, the restriction $F| \Nerve( M_1)$ is zero, and
$F$ is a right Kan extension of $F | \Nerve( M_1 \cup M_2)$. 
\end{itemize}
\end{proposition}

Assuming Proposition \ref{cobine} for the moment, we can give the proof of Theorem \ref{vedu2}.

\begin{proof}[Proof of Theorem \ref{vedu2}]
Let $\calE(\calC)$ be the full subcategory of $\Fun( \Nerve(M), \calC)$ spanned by those functors which satisfy the equivalent conditions of Proposition \ref{cobine}. The inclusions $M_0 \hookrightarrow M \hookleftarrow M_2$ determine restriction functors
$$ \Shv_{\calK}( X; \calC^{op}) \stackrel{\theta}{\leftarrow} \calE(\calC)^{op} 
\stackrel{\theta'}{\rightarrow} \Shv_{\calK}( X; \calC)^{op}.$$
Note that a functor $F \in \Fun( \Nerve(M), \calC)$ belongs to $\calE(\calC)$ if and only if
$F | \Nerve(M_0)$ belongs to $\Shv_{\calK}( X; \calC^{op})$, $F | \Nerve( M_0 \cup M_1)$ is a right Kan
extension of $F| \Nerve(M_0)$, and $F$ is a left Kan extension of $F| \Nerve( M_0 \cup M_1)$. Applying
Proposition \toposref{lklk}, we deduce that $\theta$ is a trivial Kan fibration. The same argument shows that
$\theta'$ is a trivial Kan fibration, so that $\theta$ and $\theta'$ determine an equivalence
$\Shv_{\calK}(X; \calC^{op}) \simeq \Shv_{\calK}( X; \calC)^{op}$.
\end{proof}

\begin{remark}\label{spain}
The construction $(i, S) \mapsto (2-i, X-S)$ determines an order-reversing bijection from the partially ordered
set $M$ to itself. Composition with this involution induces an isomorphism $\calE(\calC)^{op} \simeq
\calE(\calC)^{op}$, which interchanges the restriction functors $\theta$ and $\theta'$ appearing in
the proof of Theorem \ref{vedu2}. It follows that the equivalence of Theorem \ref{vedu2} is symmetric in $\calC$ and $\calC^{op}$ (up to coherent homotopy).
\end{remark}

We will give the proof of Proposition \ref{cobine} at the end of this section. For the moment, we will concentrate on the problem of making the equivalence of Theorem \ref{vedu} more explicit.  

\begin{definition}\label{spose}
Let $X$ be a locally compact topological space and let $\calC$ be a pointed $\infty$-category which admits small limits and colimits. Let $\calF$ be a $\calC$-valued sheaf on $X$. For every compact set $K \subseteq X$, we let $\Gamma_{K}(X; \calC)$ denote the fiber product
$\calF(X) \times_{ \calF(X-K)} 0$, where $0$ denotes a zero object of $\calC$.
For every open set $U \subseteq X$, we let $\Gamma_c(U; \calF)$ denote the filtered colimit
$\varinjlim_{K \subseteq U} \Gamma_{K}(M;\calF)$, where $K$ ranges over all compact subsets of $U$.
The construction $U \mapsto \Gamma_c(U; \calF)$ determines a functor
$\Nerve( \calU(X)) \rightarrow \calC$, which we will denote by $\Gamma_c( \bigdot; \calF)$.
\end{definition}

\begin{proposition}\label{capsu}
In the situation of Definition \ref{spose}, suppose that the $\infty$-category $\calC$ is stable.
Then the equivalence $\mathbb{D}$ of Theorem \ref{vedu} is given by the formula
$\mathbb{D}(\calF)(U) = \Gamma_c(U; \calF)$.
\end{proposition}

\begin{remark}
Proposition \ref{capsu} is an abstract formulation of the following more classical fact: conjugation by Verdier duality exchanges cohomology with compactly supported cohomology.
\end{remark}

\begin{proof}
It follows from the proof of Theorem \toposref{kuku} that the equivalence
$$\theta: \Shv_{\calK}(X; \calC^{op})^{op} \simeq \Shv(X; \calC^{op})^{op}$$ of Lemma \ref{boffly} is given by
the formula $\theta( \calG)(U) = \varinjlim_{K \subseteq U} \calG(K)$. Consequently, it will suffice
to show that the composition of the equivalence
$\psi: \Shv(X; \calC) \rightarrow \Shv_{\calK}(X; \calC)$ of Lemma \ref{boffly} with the equivalence
$\psi': \Shv_{\calK}(X; \calC) \rightarrow \Shv_{\calK}( X; \calC^{op})^{op}$ is given by the formula
$(\psi' \circ \psi)(\calF)(K) = \Gamma_K(X; \calF)$. To prove this, we need to introduce a bit of notation.

Let $M'$ denote the partially ordered set of pairs $(i, S)$, where $0 \leq i \leq 2$ and
$S$ is a subset of $X$ such that $S$ is compact if $i=0$ and $X - S$ is either open or compact if $i=2$;
we let $(i, S) \leq (j, T)$ if $i \leq j$ and $S \subseteq T$ or if $i =0$ and $j=2$. We will regard 
the set $M$ of Notation \ref{buxy} as a partially ordered subset of $M'$. For $0 \leq i \leq 2$, let
$M'_i$ denote the subset $\{ (j, S) \in M': j = i \} \subseteq M'$. Let $\calD$ denote the full subcategory
of $\Fun( \Nerve(M'), \calC)$ spanned by those functors $F$ which satisfy the following conditions:
\begin{itemize}
\item[$(i)$] The restriction $F | \Nerve(M_2)$ is a $\calK$-sheaf on $X$.
\item[$(ii)$] The restriction $F | \Nerve(M'_2)$ is a right Kan extension of $F | \Nerve(M_2)$.
\item[$(iii)$] The restriction $F| \Nerve( M'_1)$ is zero.
\item[$(iv)$] The restriction $F| \Nerve( M')$ is a right Kan extension of $F| \Nerve(M'_1 \cup M'_2)$.
\end{itemize}

Note that condition $(ii)$ is equivalent to the requirement that $F | \Nerve( M'_1 \cup M'_2)$ is a right
Kan extension of $F| \Nerve(M_1 \cup M_2)$. It follows from Proposition \toposref{acekan} that condition
$(iv)$ is equivalent to the requirement that $F| \Nerve(M)$ is a right Kan extension of
$F| \Nerve( M_1 \cup M_2)$. Consequently, the inclusion $M \hookrightarrow M'$ induces a
restriction functor $\calD \rightarrow \calE$, where $\calE \subseteq \Fun( \Nerve(M), \calC)$ is
defined as in the proof of Theorem \ref{vedu2}. Using Theorem \toposref{kuku} and
Proposition \toposref{lklk}, we deduce that the restriction functor
$\calD \rightarrow \Fun( \Nerve( \calU(X))^{op}, \calC)$ is a trivial Kan fibration onto the full subcategory
$\theta: \Shv(X; \calC) \subseteq \Fun( \Nerve(\calU(X))^{op}, \calC)$; moreover, the composition
$\psi' \circ \psi$ is given by composing a homotopy inverse to $\theta$ with the restriction functor
$\calD \rightarrow \Fun( \Nerve(M_0), \calC) \simeq \Fun( \Nerve( \calK(X))^{op}, \calC^{op})^{op}$. 

We define a map of simplicial sets $\phi: \Nerve(M_0) \rightarrow \Fun( \Delta^1 \times \Delta^1, \Nerve(M') )$ so that $\phi$ carries an object $(0,K) \in M_0$ to the diagram
$$ \xymatrix{ (0,K) \ar[r] \ar[d] & (1,K) \ar[d] \\
(2, \emptyset) \ar[r] & (2, K). }$$
It follows from Theorem \toposref{hollowtt} that for each $(0, K) \in M_0$, the image $\phi( 0, K)$ can
be regarded as a cofinal map $\Lambda^2_2 \rightarrow \Nerve( M')_{ (0,K)/ } \times_{ \Nerve(M')}
\Nerve( M'_1 \cup M'_2)$. Consequently, if $F \in \calD$ then condition $(iv)$ is equivalent to the requirement that the composition of $F$ with each $\phi(0,K)$ yields a pullback diagram
$$ \xymatrix{ F(0,K) \ar[r] \ar[d] & F(1, K) \ar[d] \\
F(2, \emptyset) \ar[r] & F(2, K) }$$
in the $\infty$-category $\calC$. Since $F(1,K)$ is a zero object of $\calC$ (condition $(iii)$), 
we can identify $F(0,K)$ with the kernel of the map $F(2, \emptyset) \rightarrow F(2,K)$. 
Taking $F$ to be a preimage of $\calF \in \Shv(X; \calC)$ under the functor $\theta$, we obtain
the desired equivalence
$$ (\psi' \circ \psi)(\calF)(K) \simeq \ker( \calF(X) \rightarrow \calF(X-K)) = \Gamma_{K}(X; \calF).$$
\end{proof}

\begin{corollary}\label{sloke}
Let $X$ be a locally compact topological space, let $\calC$ be a stable
$\infty$-category which admits small limits and colimits, and let $\calF \in \Shv(X; \calC)$ be a $\calC$-valued sheaf on $X$. Then the functor $\Gamma_c( \bigdot; \calF)$ is a $\calC$-valued cosheaf on $X$.
\end{corollary}

We will need the following consequence of Corollary \ref{sloke} in the next section.

\begin{corollary}\label{staffer}
Let $M$ be a manifold and let $\calF \in \Shv(M; \Spectra)$ be a spectrum-valued sheaf on $M$.
Then:
\begin{itemize}
\item[$(1)$] The functor $\calF$ exhibits $\Gamma_c( M; \calF)$ as a colimit of the diagram
$\{ \Gamma_c(U; \calF) \}_{U \in \Disk{M} }$. 
\item[$(2)$] The functor $\calF$ exhibits $\Gamma_c(M; \calF)$ as a colimit of the diagram
$\{ \Gamma_c(U; \calF) \}_{U \in \Disj{M} }$. 
\end{itemize}
\end{corollary}

\begin{proof}
We will give the proof of $(1)$; the proof of $(2)$ is similar.
According to Corollary \ref{sloke}, the functor $U \mapsto \Gamma_c(U; \calF)$ is a cosheaf
of spectra on $M$. Since every open subset of $M$ is a paracompact topological space of finite covering dimension, the $\infty$-topos
$\Shv(M)$ is hypercomplete so that $\calF$ is automatically hypercomplete. According to
Remark \ref{splait}, it will suffice to show that for every point $x \in M$, the category
$\Disk{M}_x = \{ U \in \Disk{M}: x \in U \}$ has weakly contractible nerve. This follows
from the observation that $\Disk{M}_x^{op}$ is filtered (since every open neighborhood of $M$ contains
an open set $U \in \Disk{M}_x$). 
\end{proof}

We conclude this section by giving the proof of Proposition \ref{cobine}.

\begin{proof}[Proof of Proposition \ref{cobine}]
We will prove that condition $(2)$ implies $(1)$; the converse follows by symmetry, in view of
Remark \ref{spain}. Let $F: \Nerve(M) \rightarrow \calC$ be a functor satisfying condition $(2)$, and
let $M'$ and $\calD \subseteq \Fun( \Nerve(M'), \calC)$ be defined as in the proof of Proposition \ref{capsu}. Using Proposition \toposref{lklk}, we deduce that $F$ can be extended to a functor $F': \Nerve(M') \rightarrow \calC$ belonging to $\calD$. It follows from Theorem \toposref{kuku} that the inclusion
$\calU(X)^{op} \subseteq M'_2$ determines a restriction functor $\calD \rightarrow \Shv(X; \calC)$;
let $\calF \in \Shv(X; \calC)$ be the image of $F'$ under this restriction functor. The proof of
Proposition \ref{capsu} shows that $\calG = F | \Nerve(M_0)$ is given informally by the formula
$\calG(K) = \Gamma_K(X; \calF)$.  

We first show that $\calG^{op}$ is a $\calC^{op}$-valued $K$-sheaf on $X$. For this, we must verify the following:
\begin{itemize}
\item[$(i)$] The object $\calG( \emptyset) \simeq \Gamma_{\emptyset}( X; \calF)$ is zero. This is
clear, since the restriction map $\calF(X) \rightarrow \calF(X - \emptyset)$ is an equivalence. 
\item[$(ii)$] Let $K$ and $K'$ be compact subsets of $X$. Then the diagram $\sigma$: 
$$ \xymatrix{ \calG( K \cap K') \ar[r] \ar[d] & \calG( K ) \ar[d] \\
\calG(K') \ar[r] & \calG(K \cup K') }$$
is a pushout square in $\calC$. Since $\calC$ is stable, this is equivalent to the requirement that $\sigma$ is a pullback square. This follows from the observation that $\sigma$ is the fiber
of a map between the squares
$$ \xymatrix{ \calF(X) \ar[r] \ar[d] & \calF(X) \ar[d] & \calF( X - (K \cap K') ) \ar[r] \ar[d] & \calF( X- K) \ar[d] \\
\calF(X) \ar[r] & \calF(X) & \calF( X- K') \ar[r] & \calF( X - (K \cup K')). }$$
The left square is obviously a pullback, and the right is a pullback since $\calF$ is a sheaf.
\item[$(iii)$] For every compact subset $K \subseteq X$, the canonical map
$\theta: \calG(K) \rightarrow \varprojlim_{K'} \calG(K')$ is an equivalence in $\calC$, where
$K'$ ranges over the partially ordered set $A$ of all compact subsets of $X$ which contain a neighborhood of $K$. We have a map of fiber sequences
$$ \xymatrix{ \calG(K) \ar[r]^-{\theta} \ar[d] & \varprojlim_{K' \in A} \calG(K') \ar[d] \\
\calF(X) \ar[r]^-{\theta'} \ar[d] & \varprojlim_{K' \in A} \calF(X) \ar[d] \\
\calF( X - K) \ar[r]^-{\theta''} & \varprojlim_{ K' \in A} \calF( X - K' ). }$$
It therefore suffices to show that $\theta'$ and $\theta''$ are equivalences. The map $\theta'$ is an equivalence
because the partially ordered set $A$ has weakly contractible nerve (in fact, both $A$ and $A^{op}$ are filtered). The map $\theta''$ is an equivalence because $\calF$ is a sheaf and the collection
$\{ X- K' \}_{K' \in A}$ is a covering sieve on $X-K$. 
\end{itemize}

To complete the proof, we will show that $F$ is a left Kan extension of $F| \Nerve(M_0 \cup M_1)$. 
Let $M'' \subseteq M_0 \cup M_1$ be the subset consisting of objects of the form
$(i, S)$, where $0 \leq i \leq 1$ and $S \subseteq X$ is compact. We note that $F| \Nerve(M_0 \cup M_1)$ is a left Kan extension of $F| \Nerve(M'')$. In view of Proposition \toposref{acekan}, it will suffice to show that
$F$ is a left Kan extension of $F| \Nerve(M'')$ at every element $(2, S) \in M_2$. We will prove the stronger
assertion that $F' | \Nerve( M'' \cup M'_2)$ is a left Kan extension of $F | \Nerve(M'')$. To prove this,
we let $B$ denote the subset of $M'_2$ consisting of pairs $(2, X-U)$ where $U \subseteq X$ is an open set with compact closure. In view of Proposition \toposref{acekan}, it suffices to prove the following:
\begin{itemize}
\item[$(a)$] The functor $F' | \Nerve( M'' \cup M'_2)$ is a left Kan extension of $F' | \Nerve( M'' \cup B)$.
\item[$(b)$] The functor $F' | \Nerve( M'' \cup B)$ is a left Kan extension of $F | \Nerve(M'')$.
\end{itemize}
To prove $(a)$, we note that Theorem \toposref{kuku} guarantees that $F' | \Nerve(M'_2)$ is a left Kan
extension of $F' | \Nerve(M''')$ (note that, if $K$ is a compact subset of $X$, then the collection of open
neighborhoods of $U$ of $K$ with compact closure is cofinal in the collection of all open neighborhoods of $K$ in $X$). To complete the proof, it suffices to observe that for every
object $(2, X-K) \in M'_2 - B$, the inclusion
$\Nerve( M''')_{ / (2, X-K) } \subseteq \Nerve( M'' \cup M''')_{ / (2, X -K)}$ is cofinal.
In view of Theorem \toposref{hollowtt}, this is equivalent to the requirement that for every
object $(i, S) \in M''$, the partially ordered set $P = \{ (2, X - U) \in B: (i,S) \leq (2, X-U) \leq (2, X-K) \}$
has weakly contractible nerve. This is clear, since $P$ is nonempty and stable under finite unions
(and therefore filtered). This completes the proof of $(a)$. 

To prove $(b)$, fix an open subset $U \subseteq X$ with compact closure; we wish to prove that
$F'(2,X-U)$ is a colimit of the diagram $F' | \Nerve( M'')_{/ (2, X-U)}$. For every compact set $K \subseteq X$,
let $M''_{K}$ denote the subset of $M''$ consisting of those pairs $(i,S)$ with
$(0,K) \leq (i,S) \leq (2, X-U)$. Then $\Nerve( M'')_{/ (2, X-U) }$ is a filtered colimit of
the simplicial sets $\Nerve( M''_{K} )$, where $K$ ranges over the collection of compact subsets of
$X$ which contain $U$. It follows that $\colim( F'| \Nerve(M'')_{/ (2,X-U)} )$ can be identified with the filtered colimit of the diagram $\{ \colim( F' | \Nerve(M''_K) \}_{K}$ (see \S \toposref{quasilimit1}). Consequently,
it will suffice to prove that for every compact set $K$ containing $U$, the diagram $F'$ exhibits
$F'(2,X-U)$ as a colimit of $F' | \Nerve( M''_K)$. Theorem \toposref{hollowtt} guarantees that the diagram
$(K, 0) \leftarrow (K-U,0) \rightarrow (K-U,1)$ is cofinal in $\Nerve(M''_{K})$. Consequently, we are reduced to proving that the diagram
$$ \xymatrix{ F'(0, K-U) \ar[r] \ar[d] & F'(1, K-U) \ar[d] \\
F'(0, K) \ar[r]  & F'(2, X-U) }$$ is a pushout square in $\calC$. Form a larger commutative diagram
$$ \xymatrix{ F'(0,K-U) \ar[r] \ar[d] & F(1, K-U) \ar[d] & \\
F'(0, K) \ar[r] \ar[d] & Z \ar[r] \ar[d] & F(1, K) \ar[d] \\
F(2, \emptyset) \ar[r] & F(2, K-U) \ar[r] \ar[d] & F(2, K) \ar[d] \\
& F(2, X -U ) \ar[r] & F( 2, X), }$$
where the middle right square is a pullback. Since $F'$ is a right Kan extension of
$F' | \Nerve( M_1 \cup M'_2)$, the proof of Proposition \ref{capsu} shows that the middle horizontal rectangle is also a pullback square. It follows that the lower middle square is a pullback. Since the left vertical rectangle is a pullback diagram (Proposition \ref{capsu} again), we deduce that the upper left square is a pullback.
Since $\calC$ is stable, we deduce that the upper left square is a pushout diagram. To complete
the proof of $(b)$, it suffices to show that the composite map $Z \rightarrow F(2, K-U) \rightarrow F(2,X-U)$
is an equivalence. We note that $F(1,K-U)$ and $F(2,X) \simeq \calF(\emptyset)$ are zero objects of $\calC$,
so the composite map $F(1,K-U) \rightarrow F(2,K) \rightarrow F(X)$ is an equivalence. It will therefore
suffice to show that the right vertical rectangle is a pullback square. Since the middle right square is a pullback by construction, we are reduced to proving that the lower right square is a pullback. This is the diagram
$$ \xymatrix{ \calF( (X - K ) \cup U ) \ar[r] \ar[d] & \calF(X-K) \ar[d] \\
\calF( U) \ar[r] & \calF( \emptyset ), }$$ 
which is a pullback square because $\calF$ is a sheaf and the open sets $U, X-K \subseteq X$ are disjoint.
\end{proof}

\section{Generalities on $\infty$-Operads}\label{secB}

In this appendix, we collect some general results about $\infty$-operads which are needed for the study of little cubes $\infty$-operads undertaken in the body of this paper. In \S \ref{justac}, we will describe mutually inverse ``assembly'' and ``disintegration'' processes which allow us to decompose an arbitrary unital $\infty$-operad $\calO^{\otimes}$ into {\em reduced} pieces, provided that $\calO$ is a Kan complex (for a precise statement, see Theorem \ref{assum}). The proof makes use of the notion of an {\it ornamental map} between
$\infty$-operads, which plays an important role throughout \S \ref{sec3}. It also makes use of the process of
{\it unitalization}: that is, the process of transforming an arbitrary $\infty$-operad into a unital $\infty$-operad, which we describe in \S \ref{unitos}.

In \S \symmetricref{iljest}, we introduced the notion of an {\it operadic left Kan extension}. If $\calC$ is a symmetric monoidal $\infty$-category and $\calM^{\otimes} \rightarrow \Delta^1 \times \Nerve(\FinSeg)$ is a correspondence from an $\infty$-operad
$\calM^{\otimes}_0 = \calM^{\otimes} \times_{ \Delta^1} \{0\}$ to another $\infty$-operad
$\calM^{\otimes}_{1} = \calM^{\otimes} \times_{ \Delta^1} \{1\}$, then (in good cases) operadic left Kan extension gives rise to a functor $\Alg_{ \calM_0}(\calC) \rightarrow \Alg_{ \calM_1}(\calC)$. If
we are given instead a family of $\infty$-operads $\calN^{\otimes} \rightarrow \Delta^2 \times \Nerve(\FinSeg)$, then we obtain a diagram of operadic left Kan extension functors
$$ \xymatrix{ & \Alg_{ \calN_1}(\calC) \ar[dr] & \\
\Alg_{\calN_0}(\calC) \ar[ur] \ar[rr] & & \Alg_{ \calN_2}(\calC). }$$
In \S \ref{tensor1}, we will show that this diagram commutes up to homotopy provided that the map $\calN^{\otimes} \rightarrow \Delta^2$ is a flat categorical fibration (Corollary \ref{sabbit}). This transitivity result
will play a crucial role in our analysis of tensor products of $\OpE{k}$-algebras in \S \ref{tensor2}.

In \S \symmetricref{moduldef}, we introduced the notion of a {\em coherent} $\infty$-operad, and showed that the coherence of an $\infty$-operad $\calO^{\otimes}$ guarantees the existence of a good theory of modules over arbitrary $\calO$-algebras. However, our definition of coherence was somewhat cumbersome and difficult to verify. Our goal in \S \ref{cohcrit} is to reformulate this definition in a more conceptual way. We use this reformulation in \S \ref{slabba} to prove that the $\infty$-operad $\OpE{k}$ is coherent for each $k \geq 0$ (Theorem \ref{cubecoh}). 

The final three sections of this appendix are devoted to generalizing some basic constructions of higher category theory to the setting of $\infty$-operads:
\begin{itemize}
\item[$(a)$] If $\calC$ and $\calD$ are $\infty$-categories, then the disjoint union $\calC \coprod \calD$ is again an $\infty$-category: moreover, it is the coproduct of $\calC$ and $\calD$ in the setting of $\infty$-categories.
The $\infty$-category of $\infty$-operads also admits coproducts, but these are a bit more difficult to describe: we will give an explicit construction of these coproducts in \S \ref{coprodinf}.
\item[$(b)$] If $\calC$ is an $\infty$-category and $p: K \rightarrow \calC$ is a diagram, then we can define
an overcategory $\calC_{/p}$ and an undercategory $\calC_{p/}$. This operation also has an analogue in the $\infty$-operadic setting, which we will describe in \S \ref{cluper}.
\item[$(c)$] If $\calC$ and $\calD$ are $\infty$-categories, then the product $\calC \times \calD$ is also an $\infty$-category. This operation has more than one analogue in the $\infty$-operadic setting. In addition to
the Cartesian product $\calO^{\otimes} \times_{ \Nerve(\FinSeg)} {\calO'}^{\otimes}$ of $\infty$-operads, there
is also the {\em tensor product} of $\infty$-operads (induced by the monoidal structure $\odot$ on
the model category $\PreOp$ of $\infty$-preoperads discussed in \S \symmetricref{comm1.8}). This tensor product is difficult to describe directly, but can often be analyzed using the closely related {\em wreath product} construction described in \S \ref{wreath}. Our comparison between wreath products and tensor products
(given by Theorem \ref{kuj}) will play a key role in our proof of the additivity theorem (Theorem \ref{slide}) of \S \ref{sass1}.
\end{itemize}

\subsection{Unitalization}\label{unitos}

In \S \symmetricref{comm1.9}, we introduced the notion of a {\em unital} $\infty$-operad. The $\infty$-category
of unital $\infty$-operads is a localization of the $\infty$-category of all $\infty$-operads: that is,
the inclusion from the $\infty$-category of unital $\infty$-operads to the $\infty$-category of all
$\infty$-operads admits a left adjoint. Our goal in this section is to give an explicit construction of this left adjoint. We begin by introducing some terminology.

\begin{definition}
Let $f: {\calO'}^{\otimes} \rightarrow \calO^{\otimes}$ be a map of $\infty$-operads.
We will say that {\it $f$ exhibits ${\calO'}^{\otimes}$ as a unitalization of
$\calO^{\otimes}$} if the following conditions are satisfied:
\begin{itemize}
\item[$(1)$] The $\infty$-operad ${\calO'}^{\otimes}$ is unital.
\item[$(2)$] For every unital $\infty$-operad $\calC^{\otimes}$, composition with $f$ induces an
equivalence of $\infty$-categories $\Alg_{\calC}(\calO') \rightarrow \Alg_{\calC}(\calO)$.
\end{itemize}
\end{definition}

It is clear that unitalizations of $\infty$-operads are unique up to equivalence, provided that they exist. For existence, we have the following result:

\begin{proposition}\label{amat}
Let $\calO^{\otimes}$ be an $\infty$-operad, and let $\calO^{\otimes}_{\ast}$
be the $\infty$-category of pointed objects of $\calO^{\otimes}$. Then:
\begin{itemize}
\item[$(1)$] The forgetful map $p: \calO^{\otimes}_{\ast} \rightarrow \calO^{\otimes}$ is
a fibration of $\infty$-operads (in particular, $\calO^{\otimes}_{\ast}$ is an $\infty$-operad).
\item[$(2)$] The $\infty$-operad $\calO^{\otimes}_{\ast}$ is unital.
\item[$(3)$] For every unital $\infty$-operad $\calC^{\otimes}$, composition with
$p$ induces a trivial Kan fibration
$\theta: \Alg_{\calC}( \calO_{\ast} ) \rightarrow \Alg_{\calC}(\calO)$
(here $\Alg_{\calC}(\calO_{\ast})$ denotes the $\infty$-category of
$\calC$-algebra objects in the $\infty$-operad $\calO^{\otimes}_{\ast}$).
\item[$(4)$] The map $p$ exhibits $\calO^{\otimes}_{\ast}$ as a unitalization of
the $\infty$-operad $\calO^{\otimes}$.
\end{itemize}
\end{proposition}

\begin{lemma}\label{asm}
Let $\calC$ be a pointed $\infty$-category, and let $\calD$ be an $\infty$-category with a final
object. Let $\Fun'(\calC, \calD)$ be the full subcategory of $\Fun(\calC, \calD)$ spanned by those functors which preserve final objects, and let $\Fun'( \calC, \calD_{\ast})$ be defined similarly.
Then the forgetful functor
$$ \Fun'( \calC, \calD_{\ast}) \rightarrow \Fun'( \calC, \calD)$$
is a trivial Kan fibration.
\end{lemma}

\begin{proof}
Let $\calE \subseteq \calC \times \Delta^1$ be the full subcategory spanned by objects
$(C,i)$, where either $C$ is a zero object of $\calC$ or $i = 1$. Let $\Fun'( \calE, \calD)$
be the full subcategory of $\Fun(\calE, \calD)$ spanned by those functors $F$ such that
$F(C,i)$ is a final object of $\calD$, whenever $C \in \calC$ is a zero object. 
We observe that a functor $F \in \Fun( \calE, \calD)$ belongs to $\Fun'( \calE, \calD)$ if
and only if $F_0 = F | \calC \times \{1\}$ belongs to $\Fun'( \calC, \calD)$, and $F$ is
a right Kan extension of $F_0$. We can identify $\Fun'( \calC, \calD_{\ast})$ with the full subcategory of
$\Fun( \calC \times \Delta^1, \calD)$ spanned by those functors $G$ such that
$G_0 = G| \calE \in \Fun'(\calE, \calD)$ and $G$ is a left Kan extension of $G_0$.
It follows from Proposition \toposref{lklk} that the restriction maps
$$ \Fun'( \calC, \calD_{\ast}) \rightarrow \Fun'( \calE, \calD) \rightarrow \Fun'( \calC, \calD)$$
are trivial Kan fibrations, so that their composition is a trivial Kan fibration as desired.
\end{proof}

\begin{proof}[Proof of Proposition \ref{amat}]
We first prove $(1)$.
Fix an object $X_{\ast} \in \calO^{\otimes}_{\ast}$ lying over $X \in \calO^{\otimes}$,
and let $\alpha: X \rightarrow Y$ be an inert morphism in $\calO^{\otimes}$.
Since the map $q: \calO^{\otimes}_{\ast} \rightarrow \calO^{\otimes}$ is a left fibration,
we can lift $\alpha$ to a morphism $X_{\ast} \rightarrow Y_{\ast}$, which is automatically
$q$-coCartesian. Let $\seg{n}$ denote the image of $X$ in $\FinSeg$, and choose
inert morphisms $\alpha^{i}: X \rightarrow X^i$ covering the maps $\Colp{i}: \seg{n} \rightarrow \seg{1}$ for $1 \leq i \leq n$. We claim that the induced functors $\alpha^{i}_{!}$ induce an equivalence
$( \calO^{\otimes}_{\ast})_{X} \rightarrow \prod_{1 \leq i \leq n} (\calO^{\otimes}_{\ast})_{X_i}$. 
Fix a final object $1$ in $\calO^{\otimes}$, so that $\calO^{\otimes}_{\ast}$ is equivalent
to $\calO^{\otimes}_{1/}$. The desired assertion is not equivalent to the assertion that the maps $\alpha^{i}$ induce a homotopy equivalence
$$ \bHom_{ \calO^{\otimes}}( 1, X) \rightarrow \prod_{1 \leq i \leq n} \bHom_{ \calO^{\otimes}}( 1, X_i ),$$
which follows immediately from our assumptions that $\calO^{\otimes}$ is an $\infty$-operad and that each $\alpha^{i}$ is inert.

To complete the proof that $p$ is an $\infty$-operad fibration,
let $X_{\ast}$ be as above, let $\seg{n}$ be its image in $\FinSeg$, and suppose we have chosen morphisms $X_{\ast} \rightarrow X_{\ast}^{i}$ in $\calO^{\otimes}_{\ast}$ whose images in
$\calO^{\otimes}$ are inert and which cover the inert morphisms
$\Colp{i}: \seg{n} \rightarrow \seg{1}$ for $1 \leq i \leq n$; we wish to show that the
induced diagram $\overline{\delta}: \nostar{n}^{\triangleleft} \rightarrow \calO^{\otimes}_{\ast}$ is
a $p$-limit diagram. Let $\delta = \overline{\delta} | \nostar{n}$; we wish to prove that the map
$$ (\calO^{\otimes}_{\ast})_{ / \overline{\delta}} \rightarrow (\calO^{\otimes}_{\ast})_{ / \delta}
\times_{ \calO^{\otimes}_{/ p \delta}} \calO^{\otimes}_{ / p \overline{\delta}}$$
is a trivial Kan fibration. Since $\calO^{\otimes}_{\ast}$ is equivalent to $\calO^{\otimes}_{1/}$, this is equivalent to the requirement that every extension problem of the form
$$ \xymatrix{ \bd \Delta^m \star \nostar{n} \ar[r]^{f} \ar[d] & \calO^{\otimes} \\
\Delta^m \star \nostar{n} \ar@{-->}[ur] & }$$
admits a solution, provided that $m \geq 2$, $f$ carries the initial vertex of
$\Delta^m$ to $1 \in \calO^{\otimes}$, and $f| \{m\} \star \nostar{n} = p \circ \overline{\delta}$.
Let $\pi: \calO^{\otimes} \rightarrow \Nerve(\FinSeg)$. The map $\pi \circ f$ admits a unique
extension to $\Delta^{m} \star \nostar{n}$: this is obvious if $m > 2$, and for $m=2$ it follows
from the observation that $\pi(1) = \seg{0}$ is an initial object of $\Nerve(\FinSeg)$.
The solubility of the relevant lifting problem now follows from the observation that
$p \circ \overline{\delta}$ is a $\pi$-limit diagram.

Assertion $(2)$ is clear (since $\calO^{\otimes}_{\ast}$ has a zero object), assertion
$(3)$ follows from the observation that $\theta$ is a pullback of the morphism
$\Fun'( \calC^{\otimes}, \calO^{\otimes}_{\ast}) \rightarrow \Fun'( \calC^{\otimes}, \calO^{\otimes})$
described in Lemma \ref{asm}, and assertion $(4)$ follows immediately from $(2)$ and $(3)$.
\end{proof}

We conclude this section with two results concerning the behavior of unitalization in families.

\begin{proposition}
Let $p: \calC^{\otimes} \rightarrow \calO^{\otimes}$ be a coCartesian fibration of $\infty$-operads,
where $\calO^{\otimes}$ is unital.
The following conditions are equivalent:
\begin{itemize}
\item[$(1)$] The $\infty$-operad $\calC^{\otimes}$ is unital.
\item[$(2)$] For every object $X \in \calC$, the unit object of $\calC_{X}$ (see \S \symmetricref{unitr})
is an initial object of $\calC_{X}$.
\end{itemize}
\end{proposition}

\begin{proof}
Choose an object $1 \in \calC^{\otimes}_{\seg{0}}$. Assertion $(1)$ is equivalent to the requirement that $1$ be an initial object of $\calC^{\otimes}$. Since $p(1)$ is an initial object of $\calO^{\otimes}$, this is equivalent
to the requirement that $\emptyset$ is $p$-initial (Proposition \toposref{basrel}).
Since $p$ is a coCartesian fibration, $(1)$ is equivalent to the requirement that for every morphism
$\beta: p(1) \rightarrow X$ in $\Ass$, the object $\beta_{!}(1)$ is an initial object of
$\calC^{\otimes}_{X}$ (Proposition \toposref{relcolfibtest}). Write $X = \bigoplus X_{i}$, where
each $X_i \in \calO$. Using the equivalence $\calC^{\otimes}_{ X} \simeq \prod_{i} \calC_{X_i}$, we see that it suffices to check this criterion when $X \in \calO$, in which case we are reduced to assertion $(2)$.
\end{proof}

\begin{proposition}
Let $p: \calC^{\otimes} \rightarrow \calO^{\otimes}$ be a coCartesian fibration of $\infty$-operads, where
$\calO^{\otimes}$ is unital. Then:
\begin{itemize}
\item[$(1)$] Let $q: {\calC'}^{\otimes} \rightarrow \calC^{\otimes}$ be a categorical fibration which exhibits $\calC^{\otimes}$ as a unitalization of ${\calC'}^{\otimes}$. 
Then the map $p \circ q: {\calC'}^{\otimes} \rightarrow \calO^{\otimes}$ is a coCartesian fibration of $\infty$-operads.
\item[$(2)$] For every map of unital $\infty$-operads ${\calO'}^{\otimes} \rightarrow \calO^{\otimes}$,
the map $q$ induces an equivalence of $\infty$-categories $\theta: \Alg_{\calO'}(\calC') \rightarrow
\Alg_{\calO'}(\calC)$.
\end{itemize} 
\end{proposition}

\begin{proof}
By virtue of Proposition \ref{amat}, we may assume without loss of generality that
${\calC'}^{\otimes} = \calC^{\otimes}_{\ast}$. In this case, the map $p \circ q$ factors as a composition
$\calC^{\otimes}_{\ast} \rightarrow \calO^{\otimes}_{\ast} \rightarrow \calO^{\otimes}.$
The functor $\calC^{\otimes}_{\ast} \rightarrow \calO^{\otimes}_{\ast}$ is equivalent to
$\calC^{\otimes}_{1/} \rightarrow \calO^{\otimes}_{p(1)/}$, where $1 \in \calC^{\otimes}_{\seg{0}}$ is a final object of $\calC^{\otimes}$, and therefore a coCartesian fibration (Proposition \toposref{werylonger}), and the map $\calO^{\otimes}_{\ast} \rightarrow \calO^{\otimes}$ is a trivial Kan fibration by virtue of our assumption that $\calO^{\otimes}$ is unital. This proves $(1)$.
To prove $(2)$, it suffices to observe that $\theta$ is a pullback of the map
$\Fun^{\lax}( {\calO'}^{\otimes}, \calC^{\otimes}_{\ast}) \rightarrow \Fun^{\lax}( {\calO'}^{\otimes}, \calC^{\otimes})$, which is a trivial Kan fibration by Proposition \ref{amat}.
\end{proof}

\subsection{Disintegration of $\infty$-Operads}\label{justac}

Let $A$ be an associative ring. Recall that an {\it involution} on $A$ is a map
$\sigma: A \rightarrow A$ satisfying the conditions
$$ (a+b)^{\sigma} = a^{\sigma} + b^{\sigma} \quad \quad (ab)^{\sigma} = b^{\sigma} a^{\sigma}
\quad \quad (a^{\sigma})^{\sigma} = a.$$
Let $\Ring$ denote the category of associative rings, and let $\Ring^{\sigma}$ denote the
category of associative rings equipped with an involution (whose morphisms are ring homomorphisms that are compatible with the relevant involutions). To understand the relationship between these two categories, we observe that the construction $A \mapsto A^{op}$ defines an action of the symmetric group $\Sigma_{2}$ on the category $\Ring$. The category $\Ring^{\sigma}$ can be described as
the category of (homotopy) fixed points for the action of $\Sigma_{2}$ on $\Ring$. In particular, we can 
reconstruct the category $\Ring^{\sigma}$ by understanding the category $\Ring$ together with its action of $\Sigma_{2}$.

Note that the category $\Ring$ can be described as the category of algebras over the associative
operad $\calO$ in the (symmetric monoidal) category of abelian groups. Similarly, we can describe
$\Ring^{\sigma}$ as the category of algebras over a larger operad $\calO'$ in the category of abelian groups. The relationship between $\Ring$ and $\Ring^{\sigma}$ reflects a more basic relationship between the operads $\calO$ and $\calO'$: namely, the operad $\calO$ carries an action of the group
$\Sigma_{2}$, and the operad $\calO'$ can be recovered as a kind of semidirect product
$\calO \rtimes \Sigma_{2}$. This description of $\calO'$ is potentially useful because $\calO$
is a simpler object. For example, the operad $\calO$ is {\it reduced}: that is, it contains
only a single unary operation.

Our goal in this section is to show that the paradigm described above is quite general.
Namely, if $\calO^{\otimes}$ is a unital $\infty$-operad whose underlying $\infty$-category
$\calO$ is a Kan complex, then $\calO^{\otimes}$ can be ``assembled'' (in a precise sense to be defined below) from a family of reduced $\infty$-operads parametrized by $\calO$.
We begin by introducing some of the relevant terminology.

\begin{remark}\label{coas}
Let $\calC$ be a Kan complex, and let $q: \calO^{\otimes} \rightarrow \calC \times \Nerve(\FinSeg)$
be a $\calC$-family of $\infty$-operads. Every object $X \in \calO^{\otimes}_{\seg{0}}$ is
$q$-final, so that we have a trivial Kan fibration
$\calO^{\otimes}_{/X} \rightarrow \calC_{/C} \times \Nerve(\FinSeg)$, where $C$ denotes the image of
$X$ in $\calC$. Since $\calC$ is a Kan complex, the $\infty$-category $\calC_{/C}$ is a contractible Kan complex, so that $\calO^{\otimes}_{/X}$ is equivalent to the $\infty$-operad $\calO^{\otimes}
\times_{ \calC} \{C\}$.

Let $\calO^{\otimes}$ be an arbitrary $\infty$-operad family, and suppose that
$\calC = \calO^{\otimes}_{\seg{0}}$ is a Kan complex. Then there is an equivalence of
$\infty$-operad families $e: \calO^{\otimes} \rightarrow {\calO'}^{\otimes}$, where
${\calO'}^{\otimes}$ is a $\calC$-family of $\infty$-operads (see Remark \symmetricref{sope}).
For each object $C \in \calC$, we have equivalences of $\infty$-operads
$$ \calO^{\otimes}_{/C} \rightarrow {\calO'}^{\otimes}_{/e(C)} \leftarrow {\calO'}^{\otimes}_{C}.$$
In other words, we can identify $\calO^{\otimes}$ with a $\calC$-family of $\infty$-operads whose
fibers are given by $\calO^{\otimes}_{/C}$.
\end{remark}

\begin{definition}
Let $\calO^{\otimes}$ be an $\infty$-operad. We will say that $\calO^{\otimes}$ is
{\it reduced} if $\calO^{\otimes}$ is unital and the underlying $\infty$-category
$\calO$ is a contractible Kan complex. More generally, we will say that an $\infty$-operad family $\calO^{\otimes}$ is {\it reduced} if $\calO^{\otimes}_{\seg{0}}$ is a Kan complex and, for each object
$X \in \calO^{\otimes}_{\seg{0}}$, the $\infty$-operad $\calO^{\otimes}_{/X}$ is reduced.
\end{definition}

\begin{example}
For $0 \leq k \leq \infty$, the little cubes $\infty$-operad $\OpE{k}$ is reduced.
This follows from the observation that the space $\Rect( \Cube{k}, \Cube{k})$ of
rectilinear embeddings from $\Cube{k}$ to itself is contractible.
\end{example}

\begin{definition}\label{stubly}
Let $\calO^{\otimes}$ be an $\infty$-operad family and ${\calO'}^{\otimes}$ an $\infty$-operad.
We will say that a map $\gamma: \calO^{\otimes} \rightarrow {\calO'}^{\otimes}$
{\it assembles $\calO^{\otimes}$ to ${\calO'}^{\otimes}$} if, for every $\infty$-operad
${\calO''}^{\otimes}$, composition with
$\gamma$ induces an equivalence of $\infty$-categories
$\Alg_{\calO'}(\calO'') \rightarrow \Alg_{ \calO}(\calO'')$. In this case we will also say that
${\calO'}^{\otimes}$ is an {\it assembly} of $\calO^{\otimes}$, or that $\gamma$ {\it exhibits
${\calO'}^{\otimes}$ as an assembly of $\calO^{\otimes}$}.
\end{definition}

\begin{remark}
In the situation of Definition \ref{stubly}, suppose that $\calO^{\otimes} \rightarrow \calC \times \Nerve(\FinSeg)$ is a $\calC$-family of $\infty$-operads. We can think of an object of
$\Alg_{\calO}(\calO'')$ as a family of $\infty$-operad maps $\calO^{\otimes}_{C} \rightarrow {\calO''}^{\otimes}$ parametrized by the objects $C \in \calC$. The map $\gamma$ assembles $\calO^{\otimes}$ if this is equivalent to the data
of a single $\infty$-operad map ${\calO'}^{\otimes} \rightarrow {\calO''}^{\otimes}$. In this case,
we can view ${\calO'}^{\otimes}$ as a sort of colimit of the family of $\infty$-operads $\{ \calO^{\otimes}_{C} \}_{C \in \calC}$. This description is literally correct in the case where
$\calC$ is a Kan complex.
\end{remark}

If $\calO^{\otimes}$ is an $\infty$-operad family, then
an assembly of $\calO^{\otimes}$ is clearly well-defined up to equivalence, provided that
it exists. To verify the existence, let $M$ be the collection of inert morphisms in $\calO^{\otimes}$,
so that $(\calO^{\otimes},M)$ is an object in the category $\PreOp$ of $\infty$-preoperads
(see \S \symmetricref{comm1.8}). We can then take ${\calO'}^{\otimes}$ to be the underlying
simplicial set of any fibrant replacement for $(\calO^{\otimes},M)$ with respect to the $\infty$-operadic model structure on $\PreOp$ (see Proposition \symmetricref{cannwell}).

We now describe the process of assembly in more precise terms. Let $\FamOp_{\Delta}$ denote
the simplicial category whose objects are $\infty$-operad families, where $\bHom_{ \FamOp_{\Delta}}( \calO^{\otimes}, {\calO'}^{\otimes})$ is the largest Kan complex contained in
the $\infty$-category $\Alg_{\calO}(\calO')$ of $\infty$-operad family maps from $\calO^{\otimes}$
to ${\calO'}^{\otimes}$. Let $\FamOp = \Nerve( \FamOp_{\Delta})$ be the associated
$\infty$-category. We can regard the $\infty$-category $\Cat_{\infty}^{\lax}$ of $\infty$-operads
as a full subcategory of $\FamOp$. The process of assembly can be regarded as a left
adjoint $\Assem: \FamOp \rightarrow \Cat_{\infty}^{\lax}$ to the inclusion functor. Our main result is the following:

\begin{theorem}\label{assum}
Let $\FamOpRed$ denote the full subcategory of $\FamOp$ spanned by the reduced $\infty$-operad families. Then the assembly functor $\Assem: \FamOp \rightarrow \Cat_{\infty}^{\lax}$ induces an
equivalence from $\FamOp$ to the full subcategory of $\Cat_{\infty}^{\lax}$ spanned by those
$\infty$-operads $\calO^{\otimes}$ such that the underlying $\infty$-category $\calO$ is a Kan complex.
\end{theorem}

In other words, if $\calO^{\otimes}$ is a unital $\infty$-operad such that $\calO$ is a Kan complex, then
$\calO^{\otimes}$ can be obtained (in an essentially unique way) as the assembly of a reduced $\infty$-operad family. The proof of Theorem \ref{assum} will be given at the end of this section. First, we need to establish
a criterion for testing when a map $\gamma: \calO^{\otimes} \rightarrow {\calO'}^{\otimes}$ assembles $\calO^{\otimes}$ into ${\calO'}^{\otimes}$. 

\begin{definition}\label{castu}
Let $\gamma: \calO^{\otimes} \rightarrow {\calO'}^{\otimes}$ be a map of $\infty$-operads.
We will say that $\gamma$ is {\it ornamental} if, for every object $X \in \calO^{\otimes}$, 
and every active morphism $\seg{n} \rightarrow \seg{1}$ in $\Nerve(\FinSeg)$, $\gamma$
induces a weak homotopy equivalence of simplicial sets
$$ \calO^{\otimes}_{/X} \times_{ \Nerve(\FinSeg)_{ / \seg{1} }} \{ \seg{n} \}
\rightarrow 
{\calO'}^{\otimes}_{/\gamma(X)} \times_{ \Nerve(\FinSeg)_{/ \seg{1} }} \{ \seg{n} \}.$$
\end{definition}

\begin{remark}\label{kuz}
Suppose we are given maps of $\infty$-operads $\calO^{\otimes} \stackrel{f}{\rightarrow} {\calO'}^{\otimes} \stackrel{g}{\rightarrow} {\calO''}^{\otimes}$, where $g$ is ornamental. Then
$f$ is ornamental if and only if $g \circ f$ is ornamental.
\end{remark}

\begin{remark}\label{spazzy}
Let $\gamma: \calO^{\otimes} \rightarrow {\calO'}^{\otimes}$ be a map of $\infty$-operads,
and let $\gamma': \calO^{\otimes}_{\nunit} \rightarrow {\calO'}^{\otimes}_{\nunit}$ be the induced map.
If $\gamma$ is ornamental, then $\gamma'$ is ornamental. Conversely, if $\gamma'$ is ornamental
and both $\calO^{\otimes}$ and ${\calO'}^{\otimes}$ are unital, then $\gamma$ is ornamental.
\end{remark}

\begin{definition}\label{camba}
Let $\calO^{\otimes}$ be an $\infty$-operad family and ${\calO'}^{\otimes}$ an $\infty$-operad.
We will say that an $\infty$-operad family map $\gamma: \calO^{\otimes} \rightarrow {\calO'}^{\otimes}$ is {\it ornamental} if, for every object $X \in \calO^{\otimes}_{\seg{0}}$, the induced map
$\calO^{\otimes}_{/X} \rightarrow {\calO'}^{\otimes}$ is an ornamental map of $\infty$-operads.
\end{definition}

\begin{example}\label{clupe}
Let $\calC$ be an $\infty$-category, and let $\gamma: \calC \times \Nerve(\FinSeg) \rightarrow \calC^{\amalg}$ be the canonical map. Then $\gamma$ is ornamental. Unwinding the definitions, this is equivalent to the assertion that for every
object $C \in \calC$ and each $n \geq 0$, the $\infty$-category
$\calC_{/C}^{n}$ is weakly contractible, which is clear (since $\calC_{/C}^{n}$ has a final object
$(C,C, \ldots, C)$).
\end{example}

\begin{definition}\label{plux}
Let $\calO^{\otimes}$ be an $\infty$-operad family and $\calC^{\otimes}$ an $\infty$-operad map.
We will say that a map of $\infty$-operad families $A: \calO^{\otimes} \rightarrow \calC^{\otimes}$
is {\it locally constant} if $A$ carries every morphism in $\calO$ to an equivalence in $\calC$.
We let $\Alg_{\calO}^{\loc}(\calC)$ denote the full subcategory of $\Alg_{\calO}(\calC)$ spanned
by the locally constant maps of $\infty$-operad families.
\end{definition}

The following result describes a connection between the process of assembly and the theory of ornamental maps.

\begin{proposition}\label{loam}
Let $\gamma: {\calO}^{\otimes} \rightarrow {\calO'}^{\otimes}$ be a map between small $\infty$-operad families, where ${\calO'}^{\otimes}$ is an $\infty$-operad, and the $\infty$-categories
$\calO^{\otimes}_{\seg{0}}$ and $\calO'$ are Kan complexes.
Then:
\begin{itemize}
\item[$(1)$] Suppose that $\gamma$ is ornamental and induces a weak homotopy equivalence $\calO \rightarrow \calO'$. Then, for every $\infty$-operad $\calC^{\otimes}$, composition with 
$\gamma$ induces an equivalence of $\infty$-categories
$\psi: \Alg_{\calO'}(\calC) \rightarrow \Alg_{\calO}^{\loc}(\calC)$. In particular, if
$\calO$ is a Kan complex (so that $\Alg_{\calO}^{\loc}(\calC) = \Alg_{\calO}(\calC)$), then
$\gamma$ exhibits ${\calO'}^{\otimes}$ as an assembly of $\calO^{\otimes}$.

\item[$(2)$] Conversely, suppose that $\calO^{\otimes}$ is a unital $\infty$-operad family and
that $\calO$ is a Kan complex. If $\gamma$ exhibits ${\calO'}^{\otimes}$ as an assembly of $\calO^{\otimes}$, then $\gamma$ is ornamental and the underlying map $\calO \rightarrow \calO'$ is a homotopy equivalence of Kan complexes. Moreover, the $\infty$-operad ${\calO'}^{\otimes}$ is also unital.
\end{itemize}
\end{proposition}

The proof of Proposition \ref{loam} will require some preliminaries.

\begin{lemma}\label{snoke}
Let $\gamma: \calO^{\otimes} \rightarrow {\calO'}^{\otimes}$ be a map of $\infty$-operads.
The following conditions are equivalent:
\begin{itemize}
\item[$(1)$] The map $\gamma$ is ornamental.
\item[$(2)$] For every object $X \in \calO^{\otimes}_{\seg{n}}$ and every active map
$\alpha: \seg{n} \rightarrow \seg{m}$ in $\Nerve(\FinSeg)$, $\gamma$ induces a weak homotopy
equivalence
$$ \calO^{\otimes}_{/X} \times_{ \Nerve(\FinSeg)_{ / \seg{m} }} \{ \seg{n} \}
\rightarrow 
{\calO'}^{\otimes}_{/\gamma(X)} \times_{ \Nerve(\FinSeg)_{/ \seg{m} }} \{ \seg{n} \}.$$
\end{itemize}
If $\calO$ and $\calO'$ are Kan complexes, then $(1)$ and $(2)$ are equivalent to the following other conditions:
\begin{itemize}
\item[$(3)$] For every object $X \in \calO$ the map $\gamma$ induces an equivalence
of $\infty$-categories from 
$(\calO^{\otimes})^{\acti}_{/X}$ and $( {\calO'}^{\otimes})^{\acti}_{/\gamma(X)}$.
\item[$(4)$] For every object $X \in \calO^{\otimes}$ the map $\gamma$ induces an equivalence
of $\infty$-categories from
$(\calO^{\otimes})^{\acti}_{/X}$ to $( {\calO'}^{\otimes})^{\acti}_{/\gamma(X)}$.
\end{itemize}
\end{lemma}

\begin{proof}
The implication $(2) \Rightarrow (1)$ is obvious, and the converse implication follows from
the observation that if $\alpha: \seg{n} \rightarrow \seg{m}$ is an active map inducing a decomposition $n = n_1 + \cdots + n_m$ and $X \simeq X_1 \oplus \ldots \oplus X_m$, then we have canonical equivalences
$$  \calO^{\otimes}_{/X} \times_{ \Nerve(\FinSeg)_{ / \seg{m} }} \{ \seg{n} \}
 \simeq \prod_{1 \leq i \leq m} \calO^{\otimes}_{ /X_{i}} \times_{ \Nerve( \FinSeg)_{ / \seg{1}}} \{ \seg{n_i} \}$$
$$ {\calO'}^{\otimes}_{/\gamma(X)} \times_{ \Nerve(\FinSeg)_{ / \seg{m} } } \{ \seg{n} \}
 \simeq \prod_{1 \leq i \leq m} {\calO'}^{\otimes}_{ /\gamma(X_{i})} \times_{ \Nerve( \FinSeg)_{ / \seg{1}}} \{ \seg{n_i} \}.$$
The proof of the equivalence $(3) \Leftrightarrow (4)$ is similar. Suppose now that
$\calO$ and $\calO'$ are Kan complexes. The implication
$(3) \Rightarrow (1)$ follows from the observation that if
$X \in \calO$, then $\coprod_{n} (\calO^{\otimes}_{/X} \times_{ \Nerve(\FinSeg)_{ / \seg{1} }} \{ \seg{n} \})$ and $\coprod_{n} {\calO'}^{\otimes}_{/\gamma(X)} \times_{ \Nerve(\FinSeg)_{/ \seg{1} }} \{ \seg{n} \}$
are the largest Kan complexes contained in the $\infty$-categories
$(\calO^{\otimes})^{\acti}_{/X}$ and $( {\calO'}^{\otimes})^{\acti}_{ /\gamma(X)}$, respectively.
We will complete the proof by showing that $(2) \Rightarrow (4)$.

We wish to show that for each $X \in \calO^{\otimes}$, the induced map
$\phi: (\calO^{\otimes})^{\acti}_{/X} \rightarrow ({\calO'}^{\otimes})^{\acti}_{/\gamma(X)}$ is an
equivalence of $\infty$-categories. Since $\phi$ induces a homotopy equivalence between the underlying Kan complexes, it is essentially surjective. It therefore suffices to show that
$\phi$ is fully faithful. Fix active morphisms $Y \rightarrow X \leftarrow Z$ in $\calO^{\otimes}$.
We wish to show that $\gamma$ induces a homotopy equivalence
$$ \bHom_{ \calO^{\otimes}_{/X}}( Y, Z) \rightarrow \bHom_{ {\calO'}^{\otimes}_{/\gamma(X)}}( \gamma(Y), \gamma(Z)).$$ For every $\infty$-category $\calC$, let $\calC^{\sim}$ denote the largest Kan complex contained in $\calC$. We have a map between homotopy fiber sequences
$$ \xymatrix{ \bHom_{ \calO^{\otimes}_{/X} }(Y,Z) \ar[r] \ar[d]^{\theta} & (( \calO^{\otimes})^{\acti}_{/Z})^{\sim} \ar[r] \ar[d]^{\theta'} & (( \calO^{\otimes})^{\acti}_{/X})^{\sim} \ar[d]^{\theta''} \\
\bHom_{ {\calO'}^{\otimes}_{/\gamma(X)}}(\gamma(Y), \gamma(Z)) \ar[r] & (( {\calO'}^{\otimes})^{\acti}_{/\gamma(Z)})^{\sim} \ar[r] & (( {\calO'}^{\otimes})^{\acti}_{/ \gamma(X)})^{\sim}. }$$
Assumption $(2)$ guarantees that $\theta'$ and $\theta''$ are homotopy equivalences, so that $\theta$ is a homotopy equivalence as well.
\end{proof}

\begin{remark}\label{poas}
Let $\gamma: \calC \rightarrow \calD$ be a categorical fibration of $\infty$-categories.
For every object $C \in \calC$, the induced map $\gamma_{C}: \calC_{/C} \rightarrow \calD_{/ \gamma(C)}$ is also a categorical fibration, so that $\gamma_{C}$ is an equivalence of
$\infty$-categories if and only if it is a trivial Kan fibration. Consequently, the requirement that
$\gamma_{C}$ be an equivalence for each object $C \in \calC$ is equivalent to the requirement
that $\gamma$ have the right lifting property with respect to the inclusion $\Lambda^{n}_{n} \subset \Delta^n$ for each $n \geq 0$. Since $\gamma$ is an inner fibration by assumption, this is equivalent to the requirement that $\gamma$ be a right fibration.

Combining this observation with Lemma \ref{snoke}, we obtain the following result:
if $\gamma: \calO^{\otimes} \rightarrow {\calO'}^{\otimes}$ is a fibration of $\infty$-operads
where $\calO$ and $\calO'$ are Kan complexes, then 
$\gamma$ is ornamental if and only if the induced map $\gamma': (\calO^{\otimes})^{\acti} \rightarrow
({\calO'}^{\otimes})^{\acti}$ is a right fibration. Note that the ``if'' direction is valid without the assumption that $\calO$ and $\calO'$ are Kan complexes: if $\gamma'$ is a right fibration, then
each of the maps $\calO^{\otimes}_{/X} \times_{ \Nerve(\FinSeg)_{ / \seg{1} }} \{ \seg{n} \}
\rightarrow 
{\calO'}^{\otimes}_{/\gamma(X)} \times_{ \Nerve(\FinSeg)_{/ \seg{1} }} \{ \seg{n} \}$ is a trivial
Kan fibration.
\end{remark}

\begin{remark}\label{pullbo}
Suppose we are given a homotopy pullback diagram of $\infty$-operads
$$ \xymatrix{ \calC^{\otimes} \ar[r]^{\beta} \ar[d] & {\calC'}^{\otimes} \ar[d] \\
{\calO^{\otimes}} \ar[r]^{\gamma} & {\calO'}^{\otimes}. }$$
If $\gamma$ is ornamental and the $\infty$-categories $\calO$ and $\calO'$ are
Kan complexes, then $\beta$ is ornamental. To prove this, we may assume without loss
of generality that $\gamma$ is a categorical fibration and that $\calC^{\otimes}
= {\calC'}^{\otimes} \times_{ {\calO'}^{\otimes} } \calO^{\otimes}$. Then
$\gamma$ induces a right fibration $\gamma': (\calO^{\otimes})^{\acti} \rightarrow
({\calO'}^{\otimes})^{\acti}$ (Remark \ref{poas}); it follows that the map
$\beta': (\calC^{\otimes})^{\acti} \rightarrow ({\calC'}^{\otimes})^{\acti}$ is also
a right fibration so that $\beta$ is ornamental by Remark \ref{poas}.
\end{remark}

\begin{lemma}\label{spade}
Let $f: X \rightarrow Y$ be a map of simplicial sets. If $f$ is a weak homotopy equivalence and $Y$ is a Kan complex, then $f$ is cofinal.
\end{lemma}

\begin{proof}
The map $f$ factors as a composition $X \stackrel{f'}{\rightarrow} X' \stackrel{f''}{\rightarrow} Y$, where
$f'$ is a categorical equivalence and $f''$ is a categorical fibration. Replacing $f$ by $f''$, we can reduced to the case where $f$ is a categorical fibration so that $X$ is an $\infty$-category.
According to Theorem \toposref{hollowtt}, it suffices to show that for every vertex $y \in Y$, the
fiber product $X \times_{Y} Y_{y/}$ is weakly contractible. Consider the pullback diagram
$$ \xymatrix{ X \times_{Y} Y_{y/} \ar[r]^{f'} \ar[d] & Y_{y/} \ar[d]^{g} \\
X \ar[r]^{f} & Y. }$$
The map $g$ is a left fibration over a Kan complex, and therefore a Kan fibration (Lemma \toposref{toothie2}). Since the usual model structure on simplicial sets is right proper, our diagram
is a homotopy pullback square. Because $f$ is a weak homotopy equivalence, we deduce that
$f'$ is a weak homotopy equivalence. Since $Y_{y/}$ is weakly contractible, we deduce
that $X \times_{Y} Y_{y/}$ is weakly contractible, as desired.
\end{proof}

\begin{lemma}\label{spude}
Let $f: X \rightarrow Y$ be a weak homotopy equivalence of simplicial sets, let $\calC$ be an $\infty$-category, and let $\overline{p}: Y^{\triangleright} \rightarrow \calC$ be a colimit diagram. Suppose
that $\overline{p}$ carries every edge of $Y$ to an equivalence in $\calC$. Then the composite map
$X^{\triangleright} \rightarrow Y^{\triangleright} \rightarrow \calC$ is a colimit diagram.
\end{lemma}

\begin{proof}
Let $C \in \calC$ be the image under $\overline{p}$ of the cone point of $Y^{\triangleright}$.
Let $\calC'$ be the largest Kan complex contained in $\calC$, so that $\overline{p}$ induces a map
$p: Y \rightarrow \calC'$. Factor the map $p$ as a composition
$$ Y \stackrel{p'}{\rightarrow} Z \stackrel{p''}{\rightarrow} \calC',$$
where $p'$ is anodyne and $p''$ is a Kan fibration (so that $Z$ is a Kan complex).
Lemma \ref{spade} guarantees that the inclusion $Y \rightarrow Z$ is cofinal and therefore right anodyne (Proposition \toposref{cofbasic}). Applying this observation to the lifting problem
$$ \xymatrix{ Y \ar[d] \ar[r] & \calC_{/C} \ar[d] \\
Z \ar[r]^{p''} \ar@{-->}[ur] & \calC, }$$
we deduce that $\overline{p}$ factors as a composition
$$ Y^{\triangleright} \rightarrow Z^{\triangleright} \stackrel{\overline{q}}{\rightarrow} \calC.$$
Since $p'$ is cofinal, the map $\overline{q}$ is a colimit diagram. Lemma \ref{spade} also guarantees
that the composition $f \circ p': X \rightarrow Z$ is cofinal, so that
$$ X^{\triangleright} \rightarrow Z^{\triangleright} \stackrel{ \overline{q}}{\rightarrow} \calC$$
is also a colimit diagram.
\end{proof}

\begin{proposition}\label{sorpi}
Let $\gamma: \calO^{\otimes} \rightarrow {\calO'}^{\otimes}$ be a map between small $\infty$-operads, and let $\calC^{\otimes}$ be a symmetric monoidal $\infty$-category. Assume that $\calC$ admits small colimits, and that the tensor product on $\calC$ preserves small colimits in each variable, and let $F: \Fun( \calO, \calC) \rightarrow \Alg_{\calO}(\calC)$ and $F': \Fun( \calO', \calC) \rightarrow
\Alg_{\calO'}(\calC)$ be left adjoints to the forgetful functors (Example \symmetricref{lad}).
The commutative diagram of forgetful functors
$$ \xymatrix{ \Alg_{\calO'}(\calC) \ar[r]^{\theta} \ar[d] & \Alg_{\calO}(\calC) \ar[d] \\
\Fun(\calO, \calC) \ar[r]^{\theta'} & \Fun( \calO', \calC) }$$
induces a natural transformation $\alpha: F \circ \theta' \rightarrow \theta \circ F'$
from $\Fun(\calO', \calC)$ to $\Alg_{\calO}(\calC)$.
\begin{itemize}
\item[$(1)$] If the map $\gamma$ is ornamental and $A_0: \calO' \rightarrow \calC$ is a map which
carries every morphism in $\calO'$ to an equivalence in $\calC$, then $\alpha$ induces
an equivalence $F( \theta'( A_0 ) ) \rightarrow \theta F'( A_0 )$. In particular, if
$\calO'$ is a Kan complex, then $\alpha$ is an equivalence.

\item[$(2)$] Conversely, suppose that $\alpha$ is an equivalence in the special case where $\calC = \SSet$ (equipped with the Cartesian monoidal structure) and when evaluated on the constant functor $\calO' \rightarrow \calC$ taking the value $\Delta^0$. Then $\gamma$ is ornamental.
\end{itemize}
\end{proposition}

\begin{proof}
Fix a map $A_0 \in \Fun( \calO', \calC)$ and let $X \in \calO$. Let
$\calX$ be the subcategory of $\calO^{\otimes}_{/X}$ whose objects are active maps
$Y \rightarrow X$ in $\calO^{\otimes}$ and whose morphisms are maps which induce equivalences
in $\Nerve(\FinSeg)$, and let $\calX' \subseteq {\calO'}^{\otimes}_{/\gamma(X)}$ be defined similarly.
Then $A_0$ determines diagrams $\chi: \calX \rightarrow \calC$ and
$\chi': \calX' \rightarrow \calC$ (here $\chi$ is given by composing $\chi'$ with the map
$\calX \rightarrow \calX'$ induced by $\gamma$). Using the characterization of free algebras given
in \S \symmetricref{comm33}, we deduce that $\alpha(A_0)(X): (F \circ \theta')(A_0)(X)
\rightarrow (\theta \circ F')(A_0)(X)$ is given by the evident map
$\colim_{ \calX} \chi \rightarrow \colim_{\calX'} \chi'$. If $A_0$ carries every morphism in
$\calO'$ to an equivalence in $\calC$, then $\chi'$ carries every morphism in $\calX'$ to an equivalence
in $\calC$. If $\gamma$ is ornamental, then the evident
map $\calX \rightarrow \calX'$ is a weak homotopy equivalence, so
that $\alpha$ is an equivalence by Lemma \ref{spude}: this proves $(1)$.

Conversely, suppose that the hypotheses of $(2)$ are satisfied. Taking
$A_0$ to be the constant functor taking the value $\Delta^0 \in \SSet$, we deduce from Corollary \toposref{needka} that the map $\calX \rightarrow \calX'$ is a weak homotopy equivalence
for each $X \in \calO$, so that $\gamma$ is ornamental.
\end{proof}

\begin{proposition}\label{lebmag}
Let $S$ be a Kan complex, let $\calO^{\otimes} \rightarrow S \times \Nerve(\FinSeg)$ be
an $S$-family of $\infty$-operads, and let $\calC^{\otimes}$ be a symmetric monoidal
$\infty$-category. Suppose that, for each $s \in S$, the restriction functor
$\Alg_{ \calO_{s}}( \calC) \rightarrow \Fun( \calO_{s}, \calC)$ admits a left adjoint $F_{s}$.
Then:
\begin{itemize}
\item[$(1)$] The restriction functor $\theta: \Alg_{ \calO}(\calC) \rightarrow \Fun( \calO, \calC)$
admits a left adjoint $F$.
\item[$(2)$] Let $A \in \Alg_{\calO}(\calC)$, let $B \in \Fun(\calO, \calC)$, and let
$\alpha: B \rightarrow \theta(A)$ be a morphism in $\Fun(\calO,\calC)$. Then the adjoint
map $F(B) \rightarrow A$ is an equivalence in $\Alg_{\calO}(\calC)$ if and only if, for each
$s \in S$, the underlying map $F_{s}( B | \calO_{s}) \rightarrow A | \calO^{\otimes}_{s}$
is an equivalence in $\Alg_{\calO_{s}}(\calC)$.
\end{itemize}
\end{proposition}

\begin{proof}
Fix $B \in \Fun( \calO, \calC)$. For every map of simplicial sets $\psi: T \rightarrow S$, let
$\calO_{T} = \calO \times_{S} T$, $B_{T} = B | \calO_{T}$, and 
$X(T)$ denote the full subcategory of $\Alg_{ \calO_{T} }( \calC) \times_{
\Fun( \calO_{T}, \calC) } \Fun( \calO_{T}, \calC)_{B_T/}$ spanned by those
objects $(A_{T} \in \Alg_{ \calO_{T}}(\calC), \phi: B_{T} \rightarrow A_{T} | \calO_{T})$ such that, for each
vertex $t \in T$, the induced map $F_{\psi(t)}( B_T | \calO_{\psi(t)} ) \rightarrow A_{T} | \calO^{\otimes}_{ \psi(t)}$ is an equivalence. We claim that every inclusion of simplicial sets
$i: T' \hookrightarrow T$ in $(\sSet)_{/S}$, the restriction map $X(T) \rightarrow X(T')$ is a trivial
Kan fibration. The collection of maps $i$ for which the conclusion holds is clearly weakly saturated;
it therefore suffices to prove the claim in the case where $i$ is an inclusion of the form
$\bd \Delta^n \subset \Delta^n$. The proof proceeds by induction on $n$. The inductive hypothesis
implies that the restriction map $X( \bd \Delta^n) \rightarrow X(\emptyset) \simeq \Delta^0$ is a
trivial Kan fibration, so that $X( \bd \Delta^n)$ is a contractible Kan complex. The map
$X( \Delta^n) \rightarrow X(\bd \Delta^n)$ is evidently a categorical fibration; it therefore suffices to show
that it is a categorical equivalence. In other words, it suffices to show that $X( \Delta^n)$ is also a contractible Kan complex. Let $s \in S$ denote the image of the vertex $\{0\} \in \Delta^n$ in
$S$. Since the inclusion $\calO^{\otimes}_{s} \hookrightarrow \calO^{\otimes}_{\Delta^n}$ is
a categorical equivalence, it induces a categorical equivalence $X(\Delta^n) \rightarrow X( \{s\})$.
We are therefore reduced to proving that $X( \{s\})$ is a contractible Kan complex, which is obvious.

The above argument shows that $X(S)$ is a contractible Kan complex; in particular, $X(S)$ is nonempty.
Consequently, there exists a map $\phi: B \rightarrow \theta(A)$ satisfying the condition described in $(2)$. We will prove $(1)$ together with the ``if'' direction of $(2)$ by showing that that $\phi$ induces
a homotopy equivalence $\rho: \bHom_{ \Alg_{\calO}(\calC)}( A, C) \rightarrow \bHom_{ \Fun(\calO, \calC)}( B, \theta(C) )$ for each $C \in \Alg_{\calO}(\calC)$. The ``only if'' direction of $(2)$ will then follow by the usual uniqueness argument. We proceed as before: for every map of simplicial sets $T \rightarrow S$,
let $Y(T)$ denote the $\infty$-category $\Alg_{\calO_{T}}(\calC)_{(A | \calO^{\otimes}_{T})/} \times_{ \Fun(\calO_{T}, \calC)_{(A_{T} | \calO_{T})/}}
\Fun( \calO_{T}, \calC)_{ \phi_{T}/}$ and $Y'(T) = \Fun( \calO_{T}, \calC)_{ (B | \calO_{T})/}$. 
The map $\rho$ can be regarded as a pullback of the restriction map $Y(S) \rightarrow Y'(S)$.
To complete the proof, it will suffice to show that $Y(S) \rightarrow Y'(S)$ is a trivial Kan fibration.
We will prove the following stronger assertion: for every inclusion $T' \hookrightarrow T$ in
$(\sSet)_{/S}$, the restriction map $\pi: Y(T) \rightarrow Y(T') \times_{ Y'(T')} Y'(T)$ is a trivial Kan fibration. As before, the collection of inclusions which satisfy this condition is weakly saturated, so we may reduce to the case where $T = \Delta^n$, $T' = \bd \Delta^n$, and the result holds for inclusions of simplicial sets having dimension $< n$. Moreover, since $\pi$ is easily seen to be a categorical fibration, it suffices to show that $\pi$ is a categorical equivalence. Using the inductive hypothesis, we deduce that
$Y(T') \rightarrow Y'(T')$ is a trivial Kan fibration, so that the pullback map
$Y(T') \times_{ Y'(T')} Y'(T) \rightarrow Y'(T)$ is a categorical equivalence. By a two-out-of-three argument, we are reduced to proving that the restriction map $Y(T) \rightarrow Y'(T)$ is a categorical equivalence. If we define $s$ to be the image of $\{0\} \subseteq \Delta^n \simeq T$ in $S$, then
we have a commutative diagram
$$ \xymatrix{ Y(T) \ar[r] \ar[d] & Y'(T) \ar[d] \\
Y( \{s\}) \ar[r] & Y'(\{s\} ) }$$
in which the vertical maps are categorical equivalences. We are therefore reduced to showing
that $Y(\{s\}) \rightarrow Y'( \{s\})$ is a categorical equivalence, which is equivalent to the requirement
that the map $F_{s}( B | \calO_{s}) \rightarrow A | \calO^{\otimes}_{s}$ be an equivalence in
$\Alg_{ \calO_{s}}(\calC)$.
\end{proof}

\begin{remark}\label{pullbo2}
In the situation of Definition \ref{camba}, suppose that $\calO$, $\calO'$, and $\calO^{\otimes}_{\seg{0}}$ are Kan complexes. A map $\gamma: \calO^{\otimes} \rightarrow {\calO'}^{\otimes}$
is ornamental if and only if it induces an equivalence of $\infty$-categories
$(\calO^{\otimes})^{\acti}_{/X} \rightarrow ({\calO'}^{\otimes})^{\acti}_{/ \gamma(X)}$ for every
$X \in \calO^{\otimes}$. The alternative characterizations given in Remark \ref{poas}, Remark \ref{pullbo}, and Lemma \ref{snoke} remain valid in this context.
\end{remark}


\begin{proof}[Proof of Proposition \ref{loam}]
We first prove $(1)$. Note that since $\gamma$ induces a homotopy equivalence
$\gamma_0: \calO \rightarrow \calO'$ and $\calO'$ is a Kan complex, the map $\gamma_0$ is essentially surjective. We may assume withoout loss of generality that the $\infty$-operad $\calC^{\otimes}$ is small. Let $\calD$ be the $\infty$-category of presheaves
$\calP((\calC)^{\otimes}_{\acti})$.
The small $\infty$-category $(\calC)^{\otimes}_{\acti} \simeq \Envv(\calC)$ has the structure
of a symmetric monoidal $\infty$-category, and there is a fully faithful embedding of
$\infty$-operads ${\calC}^{\otimes} \hookrightarrow \Envv(\calC)^{\otimes}$ (see Remark \symmetricref{lazier}). Combining this with Corollary \symmetricref{banpan}, we obtain
a symmetric monoidal structure on the $\infty$-category $\calD$ (such that the tensor product
preserves colimits separately in each variable) and a fully faithful embedding of $\infty$-operads
${\calC}^{\otimes} \hookrightarrow \calD^{\otimes}$. Let ${\calD'}^{\otimes}$ denote the essential image of this embedding; it will suffice to show that the restriction map $\Alg_{\calO'}(\calD') \rightarrow
\Alg^{\loc}_{\calO}(\calD')$ is an equivalence of $\infty$-categories. We have a commutative diagram
$$ \xymatrix{ \Alg_{\calO'}(\calD') \ar[r] \ar[d] & \Alg^{\loc}_{\calO'}(\calD) \ar[d] \\
\Alg_{\calO}( \calD') \ar[r] & \Alg^{\loc}_{\calO}(\calD). }$$
Since $\gamma_0$ is essentially surjective, this diagram is a homotopy pullback square.
Consequently, it will suffice to show that the lower horizontal map is a categorical equivalence.
We may therefore replace $\calC^{\otimes}$ with $\calD^{\otimes}$, and thereby reduce to the
case where $\calC$ is a symmetric monoidal $\infty$-category which admits small colimits, where
the tensor product on $\calC$ preserves small colimits separately in each variable.

We may assume without loss of generality that $\calO^{\otimes}$ is an $S$-family of $\infty$-operads
for some Kan complex $S$. Using Corollary \symmetricref{spaltwell} (and Proposition \ref{lebmag}), we deduce that the forgetful functors
$$ \theta: \Alg_{\calO}(\calC) \rightarrow \Fun(\calO, \calC) \quad \quad \theta': \Alg_{\calO'}(\calC) \rightarrow \Fun( \calO', \calC)$$
admit left adjoints $F: \Fun( \calO, \calC) \rightarrow \Alg_{\calO}(\calC)$, $F': \Fun( \calO', \calC)
\rightarrow \Alg_{\calO'}(\calC)$. Let $\Fun^{\loc}(\calO, \calC)$ denote the full subcategory
of $\Fun( \calO, \calC)$ spanned by those functors which carry each morphism in $\calO$ to an equivalence in $\calC$. Since $\gamma_0$ is a weak homotopy equivalence and
$\calO'$ is a Kan complex, composition with $\gamma_0$ induces a categorical equivalence
$\phi: \Fun(\calO', \calC) \rightarrow \Fun^{\loc}(\calO, \calC)$; we let $\phi^{-1}$ denote any homotopy inverse to $\phi$. Propositions \ref{lebmag} and \ref{sorpi} guarantee that the canonical natural
transformation $F \circ \phi \rightarrow \psi \circ F'$ is an equivalence of functors. In particular,
$F$ carries the essential image $\Fun^{\loc}(\calO, \calC)$ into $\Alg^{\loc}_{\calO}(\calC)$.
It follows that $F$ and $\theta$ restrict to a pair of adjoint functors
$\Adjoint{ F^{\loc}}{ \Fun^{\loc}(\calO,\calC) }{ \Alg^{\loc}_{\calO}(\calC).}{\theta^{\loc}}$

We wish to prove that $\psi$ is an equivalence of $\infty$-categories. To this end, we consider
the diagram
$$ \xymatrix{ \Alg_{\calO'}(\calC) \ar[dr] \ar[rr]^{\psi} & & \Alg_{\calO}^{\loc}(\calC) \ar[dl]^{\theta^{\loc}} \\
& \Fun^{\loc}(\calO, \calC). & }$$
Using Corollaries \symmetricref{jumunj22} and \symmetricref{fillfemme}, we deduce
that the functors $\theta'$ and $\theta' \circ \psi \simeq \phi \circ \psi' \circ \phi^{-1}$ are conservative and preserve geometric realizations of simplicial objects. Consequently, to prove that
$\psi$ is an equivalence of $\infty$-categories, it will suffice to show that the map of
monads $\theta^{\loc} \circ F^{\loc} \rightarrow \phi \circ \theta' \circ F' \circ \phi^{-1}$
is an equivalence (Corollary \monoidref{littlerbeck}), which follows from Proposition \ref{sorpi}. This completes the proof of $(1)$.

To prove $(2)$, suppose that $\calO^{\otimes}$ is unital and that $\gamma$
exhibits ${\calO'}^{\otimes}$ as an assembly of $\calO^{\otimes}$. 
It follows from Proposition \ref{amat} that for each $s \in S$, the induced map
$\Alg_{ \calO_{s}}( \calO'_{\ast}) \rightarrow \Alg_{\calO_{s}}(\calO)$ is a trivial Kan fibration.
Arguing as in Proposition \ref{lebmag}, we deduce that $\Alg_{ \calO}(\calO'_{\ast}) \rightarrow
\Alg_{\calO}(\calO')$ is a trivial Kan fibration. Since $\gamma$ exhibits ${\calO'}^{\otimes}$
as an assembly of $\calO^{\otimes}$, we deduce that the map $\Alg_{\calO'}(\calO'_{\ast}) \rightarrow \Alg_{\calO'}(\calO')$ is an equivalence of $\infty$-categories, and therefore (since it is a categorical fibration) a trivial Kan fibration. In particular, the projection map ${\calO'}^{\otimes}_{\ast} \rightarrow {\calO'}^{\otimes}$ admits a section, so the final object of ${\calO'}^{\otimes}$ is initial and
${\calO'}^{\otimes}$ is also unital. Let $\calC$ be an arbitrary
$\infty$-category, which we regard as the underlying $\infty$-category of the $\infty$-operad
$\calC^{\amalg}$. We have a commutative diagram
$$ \xymatrix{ \Alg_{\calO'}(\calC) \ar[r] \ar[d] & \Alg_{\calO}(\calC) \ar[d] \\
\Fun( \calO', \calC) \ar[r] & \Fun( \calO, \calC) }$$
where the upper horizontal map is an equivalence and the vertical maps are equivalences
by virtue of Proposition \symmetricref{kime3}. It follows that the lower horizontal map is an equivalence.
Allowing $\calC$ to vary, we deduce that $\gamma$ induces an equivalence of $\infty$-categories
$\calO \rightarrow \calO'$.

It remains to show $\gamma$ is ornamental.
To prove this, let us regard $\SSet$ as endowed with the Cartesian monoidal structure, and let
$A_0: \calO' \rightarrow \SSet$ be the constant functor taking the value $\Delta^0$. Since
$\psi: \Alg_{\calO'}(\calC) \rightarrow \Alg_{\calO}(\calC)$ is an equivalence of $\infty$-categories, the canonical map $F'(A_0) \rightarrow (F \circ \phi)(A_0)$ is an equivalence (where $F$, $F'$, and $\psi$
are defined as above). Using the characterization of $F$ given in Proposition \ref{lebmag}, we deduce that for each $s \in S$ the induced map $\gamma: \calO^{\otimes}_{s} \rightarrow {\calO'}^{\otimes}$ satisfies the criterion of Proposition \ref{sorpi} and is therefore ornamental. It follows that $\gamma$
is ornamental as desired.
\end{proof}

We now turn to the proof of Theorem \ref{assum}. We need one more preliminary result:

\begin{lemma}\label{stoag}
Let $\gamma: \calO^{\otimes} \rightarrow {\calO'}^{\otimes}$ be an ornamental map of
$\infty$-operads. Suppose that $\calO^{\otimes}$ and ${\calO'}^{\otimes}$ are reduced.
Then $\gamma$ is an equivalence of $\infty$-operads.
\end{lemma}

\begin{proof}
Since $\calO^{\otimes}$ is reduced, each of the $\infty$-categories
$\calO^{\otimes}_{\seg{n}}$ has a unique object (up to equivalence) which we will
denote by $X_n$. The image $\gamma( X_n )$ can be identified with the unique
object of ${\calO'}^{\otimes}_{\seg{n}}$. It follows that $\gamma$ is essentially surjective.
We will complete the proof by showing that $\gamma$ is fully faithful. For every morphism
$\alpha: \seg{m} \rightarrow \seg{n}$ in $\FinSeg$, let $\bHom_{\calO^{\otimes}}^{\alpha}(X_m, X_n)$ be
the summand of $\bHom_{\calO^{\otimes}}(X_m, X_n)$ consisting of those connected components lying
over $\alpha \in \Hom_{\FinSeg}( \seg{m}, \seg{n})$, and define $\bHom_{ {\calO'}^{\otimes}}^{\alpha}( \gamma X_m, \gamma X_n)$ similarly. We will prove that each of the maps
$$ \bHom_{ \calO^{\otimes}}^{\alpha}( X_m, X_n) \rightarrow \bHom_{{\calO'}^{\otimes}}^{\alpha}( \gamma X_m, \gamma X_n)$$
is a homotopy equivalence. To begin, choose a factorization of $\alpha$
as a composition $\seg{m} \stackrel{\alpha'}{\rightarrow} \seg{m'} \stackrel{\alpha''}{\rightarrow} \seg{n}$, where $\alpha'$ is inert and $\alpha''$ is active. The map $\alpha'$ lifts (in an essentially unique fashion)
to an inert morphism $X_{m} \rightarrow X_{m'}$ in $\calO^{\otimes}$, and we have a homotopy
commutative diagram
$$ \xymatrix{  \bHom_{ \calO^{\otimes}}^{\alpha''}( X_{m'}, X_n) \ar[r] \ar[d]  & \bHom_{{\calO'}^{\otimes}}^{\alpha''}( \gamma X_{m'}, \gamma X_n) \\
 \bHom_{ \calO^{\otimes}}^{\alpha}( X_m, X_n) \ar[r] &  \bHom_{{\calO'}^{\otimes}}^{\alpha}( \gamma X_m, \gamma X_n) }$$
in which the vertical maps are homotopy equivalences. We may therefore replace $\alpha$ by
$\alpha''$ and thereby reduce to the case where $\alpha$ is active. Passing to the union over
all active maps $\seg{m} \rightarrow \seg{n}$, we are reduced to proving that the map
$$ (\calO^{\otimes})_{/X_n}^{\acti} \times_{ \calO^{\otimes} } \{ X_m \}
\rightarrow ( {\calO'}^{\otimes})^{\acti}_{/ \gamma X_n} \times_{ {\calO'}^{\otimes}} \{ \gamma X_m \}$$
is a homotopy equivalence. The desired conclusion now follows by examining the commutative
diagram
$$ \xymatrix{ 
(\calO^{\otimes})_{/X_n}^{\acti} \times_{ \calO^{\otimes} } \{ X_m \}
\ar[r] \ar[d] & ( {\calO'}^{\otimes})^{\acti}_{/ \gamma X_n} \times_{ {\calO'}^{\otimes}} \{ \gamma X_m \} \ar[d] \\
 ( \calO^{\otimes})_{/X_n}^{\acti} \times_{ \Nerve(\FinSeg)} \{ \seg{m} \}
\ar[r] &  ( {\calO'}^{\otimes})^{\acti}_{/ \gamma X_n} \times_{ \Nerve(\FinSeg)} \{ \seg{m} \}. }$$
The vertical maps are categorical equivalences since $\calO^{\otimes}$ and ${\calO'}^{\otimes}$ are reduced, and the lower horizontal map is a categorical equivalence because $\gamma$ is ornamental.
\end{proof}

\begin{proof}[Proof of Theorem \ref{assum}]
It follows
from Proposition \ref{loam} that the assembly functor $\Assem$ carries $\FamOpRed$
into the full subcategory $\calX \subseteq \Cat_{\infty}^{\lax}$ spanned by those
those unital $\infty$-operads $\calO^{\otimes}$ such that $\calO$ is a Kan complex. 
We next show that $\Assem: \FamOpRed \rightarrow \calX$ is essentially surjective.
Let $\calO^{\otimes}$ be such an $\infty$-operad and
choose a homotopy equivalence $f: \calO \rightarrow S$ for some Kan complex $S$ (for
example, we can take $S = \calO$ and $f$ to be the identity map).
Using Proposition \symmetricref{kime1}, we can extend $f$ to an $\infty$-operad map
$f': \calO^{\otimes} \rightarrow S^{\amalg}$. Replacing $\calO^{\otimes}$ by an
equivalent $\infty$-operad if necessary, we may suppose that $f'$ is a fibration of $\infty$-operads.
Let ${\calO'}^{\otimes}$ be the fiber product $\calO^{\otimes} \times_{ S^{\amalg} } (S \times \Nerve(\FinSeg))$. Then ${\calO'}^{\otimes}$ is an $S$-family of $\infty$-operads equipped with a map
$\gamma: {\calO'}^{\otimes} \rightarrow \calO^{\otimes}$ which induces an isomorphism
$\calO' \rightarrow \calO$. The map $\gamma$ is a homotopy pullback of the ornamental map
$S \times \Nerve(\FinSeg) \rightarrow S^{\amalg}$ of Example \ref{clupe}, so that
$\gamma$ is ornamental (Remarks \ref{pullbo} and \ref{pullbo2}). Invoking Proposition \ref{loam}, we deduce
that $\gamma$ exhibits $\calO^{\otimes}$ as an assembly of ${\calO'}^{\otimes}$, so that
we have an equivalence $\Assem( {\calO'}^{\otimes}) \simeq \calO^{\otimes}$.
To deduce the desired essential surjectivity, it suffices to show that ${\calO'}^{\otimes}$ is reduced.
In other words, we must show that for each $s \in S$, the $\infty$-operad
${\calO'}^{\otimes}_{s} \simeq \calO^{\otimes} \times_{ S^{\amalg}} \Nerve(\FinSeg)$ is
reduced. This is clear: the underlying $\infty$-category $\calO_{s}$ is given by
the fiber of a trivial Kan fibration $f: \calO \rightarrow S$, and ${\calO'}^{\otimes}_{s}$ is unital
because it is a homotopy fiber product of unital $\infty$-operads.

We now show that $\Assem: \FamOpRed \rightarrow \calX$ is fully faithful.
Let $\calC^{\otimes}$ and $\calD^{\otimes}$ be reduced $\infty$-operad families,
and choose assembly maps $\calC^{\otimes} \rightarrow {\calC'}^{\otimes}$ and
$\calD^{\otimes} \rightarrow \calO^{\otimes}$. We will show that the canonical map
$\Alg_{\calC}(\calD) \rightarrow \Alg_{\calC}(\calO) \simeq \Alg_{\calC'}(\calO)$ is
an equivalence of $\infty$-categories. As above, we choose a Kan complex $S \simeq \calO$
and a fibration of $\infty$-operads $\calO^{\otimes} \rightarrow S^{\amalg}$, and define
${\calO'}^{\otimes}$ to be the fiber product $(S \times \Nerve(\FinSeg)) \times_{ S^{\amalg} } \calO^{\otimes}$. Using the equivalences $\Alg_{\calC}( S^{\amalg}) \simeq \Fun(\calC, S)$
and $\Alg_{\calC}( S \times \Nerve(\FinSeg)) \simeq \Fun( \calC^{\otimes}_{\seg{0}}, S)$
provided by Propositions \symmetricref{kime3} and \symmetricref{slaper}, we obtain a homotopy pullback diagram of $\infty$-categories
$$ \xymatrix{ \Alg_{\calC}( \calO') \ar[r] \ar[d] & \Alg_{\calC}( \calO) \ar[d] \\
\Fun( \calC^{\otimes}_{\seg{0}},S) \ar[r] & \Fun( \calC, S). }$$
Here the lower horizontal map is obtained by composing with the functor
$\calC = \calC^{\otimes}_{\seg{1}} \rightarrow \calC^{\otimes}_{\seg{0}}$ induced by
the map $\seg{1} \rightarrow \seg{0}$ in $\FinSeg$. Since $\calC$ is reduced, this map is an
equivalence of $\infty$-categories, so the natural map $\Alg_{\calC}(\calO') \rightarrow
\Alg_{\calC}(\calO)$ is an equivalence. Similarly, we have an equivalence $\Alg_{\calD}(\calO') \rightarrow \Alg_{\calD}(\calO)$. We may therefore assume that the assembly map
$\calD^{\otimes} \rightarrow \calO^{\otimes}$ factors through a map of $\infty$-operad families
$\gamma: \calD^{\otimes} \rightarrow {\calO'}^{\otimes}$. To complete the proof, it will suffice to show
that $\gamma$ is an equivalence of $\infty$-operad families (and therefore induces an equivalence
of $\infty$-categories $\Alg_{\calC}(\calD) \rightarrow \Alg_{\calC}(\calO') \simeq \Alg_{\calC}(\calO)$).

Replacing $\calD^{\otimes}$ by an equivalent $\infty$-operad family if necessary, we can assume
that $\gamma: \calD^{\otimes} \rightarrow {\calO'}^{\otimes}$ is a categorical fibration, so that
the composite map $\calD^{\otimes} \rightarrow {\calO'}^{\otimes} \rightarrow S \times \Nerve(\FinSeg)$ exhibits $\calD$ as an $S$-family of $\infty$-operads. It will therefore suffice to show
that for each $s \in S$, the induced map of fibers $\gamma_{s}: \calD^{\otimes}_{s} \rightarrow {\calO'}^{\otimes}_{s}$ is an equivalence of $\infty$-operads. For each $D \in \calD^{\otimes}_{s}$ having an image $X \in \calO^{\otimes}$, we have a commutative
diagram
$$ \xymatrix{ (\calD^{\otimes}_{s})^{\acti}_{/D} \ar[r] \ar[d] & (\calD^{\otimes})^{\acti}_{/D} \ar[d] \ar[r] & 
(\calO^{\otimes})^{\acti}_{/X} \ar[d] \\
( {\calO'}^{\otimes}_{s})^{\acti}_{/ \gamma(D)} \ar[r] & ( {\calO'}^{\otimes})^{\acti}_{/\gamma(D)}
\ar[r] & (\calO^{\otimes})^{\acti}_{/X} }$$
in which the horizontal maps are categorical equivalences (Proposition \ref{loam}).
It follows that the vertical maps are also categorical equivalences, so that
$\gamma_{s}$ is an ornamental map between reduced $\infty$-operads. It follows
from Lemma \ref{stoag} that $\gamma_{s}$ is an equivalence of $\infty$-operads as desired.
\end{proof}

\subsection{Transitivity of Operadic Left Kan Extensions}\label{tensor1}

In this section, we will prove the following transitivity formula for operadic left Kan extensions:

\begin{theorem}\label{transit}
Let $\calM^{\otimes} \rightarrow \Delta^2 \rightarrow \Nerve(\FinSeg)$ be a $\Delta^2$-family of $\infty$-operads (Definition \symmetricref{sike}). Let $q: \calC^{\otimes} \rightarrow \calD^{\otimes}$ be a fibration of $\infty$-operads, and let $A: \calM^{\otimes} \rightarrow \calC^{\otimes}$ be an $\infty$-operad family map. Assume that $A | ( \calM^{\otimes} \times_{ \Delta^2} \Delta^{ \{0,1\} })$ and
$A | (\calM^{\otimes} \times_{ \Delta^2} \Delta^{ \{1,2\} })$ are operadic $q$-left Kan extensions,
and that the map $\calM^{\otimes} \rightarrow \Delta^2$ is a flat categorical fibration.
Then $A | (\calM^{\otimes} \times_{ \Delta^2} \Delta^{ \{0,2\}})$ is an operadic $q$-left Kan extension.
\end{theorem}

Theorem \ref{transit} has the following consequence:

\begin{corollary}\label{sabbit}
Let $\calM^{\otimes} \rightarrow \Delta^2 \rightarrow \Nerve(\FinSeg)$ be a $\Delta^2$-family of $\infty$-operads, $\calC^{\otimes}$ a symmetric monoidal $\infty$-category, and $\kappa$ an uncountable regular cardinal. Assume that:
\begin{itemize}
\item[$(i)$] The $\infty$-category $\calM^{\otimes}$ is essentially $\kappa$-small.
\item[$(ii)$] The $\infty$-category $\calC$ admits $\kappa$-small colimits, and the tensor product
on $\calC$ preserves $\kappa$-small colimits separately in each variable.
\item[$(iii)$] The projection map $\calM^{\otimes} \rightarrow \Delta^2$ is a flat categorical fibration.
\end{itemize}
For $i \in \{0,1,2\}$, let $\calM^{\otimes}_{i}$ denote the fiber $\calM^{\otimes} \times_{ \Delta^2} \{i\}$.
Let $f_{0,1}: \Alg_{\calM_0}(\calC) \rightarrow \Alg_{\calM_1}(\calC)$, $f_{1,2}: \Alg_{\calM_1}(\calC) \rightarrow
\Alg_{\calM_2}(\calC)$, and $f_{0,2}: \Alg_{\calM_0}(\calC) \rightarrow \Alg_{\calM_2}(\calC)$ be the functors
given by operadic $q$-left Kan extension (see below). Then there is a canonical equivalence of functors
$f_{0,2} \simeq f_{1,2} \circ f_{0,1}$.
\end{corollary}

\begin{proof}
For $0 \leq i \leq j \leq 2$, let $\Alg_{i,j}( \calC)$ denote the full subcategory of 
$\Fun_{ \Nerve(\FinSeg)}( \calM^{\otimes} \times_{ \Delta^2} \Delta^{ \{i,j\}}), \calC^{\otimes})$ spanned
by those $\infty$-operad family maps which are operadic $q$-left Kan extensions, where
$q: \calC^{\otimes} \rightarrow \Nerve(\FinSeg)$ denotes the projection. Using Lemma
\symmetricref{siwd}, Theorem \symmetricref{oplk}, and Proposition \symmetricref{slavewell}, we
see that conditions $(i)$ and $(ii)$ guarantee that the restriction map $r: \Alg_{i,j}(\calC) \rightarrow \Alg_{ \calM_{i}}(\calC)$ is a trivial Kan fibration. The map $f_{i,j}$ is defined to be the composition
$$ \Alg_{\calM_i}(\calC) \stackrel{s}{\rightarrow} \Alg_{i,j}(\calC) \rightarrow \Alg_{\calM_j}(\calC),$$
where $s$ is a section of $r$. Consequently, the composition $f_{1,2} \circ f_{0,1}$ can be defined
as a composition
$$ \Alg_{ \calM_0}(\calC) \stackrel{s'}{\rightarrow} \Alg_{0,1}(\calC) \times_{ \Alg_{\calM_1}(\calC)}
\Alg_{1,2}(\calC) \rightarrow \Alg_{\calM_2}(\calC),$$
where $s'$ is a section of the trivial Kan fibration
$\Alg_{0,1}(\calC) \times_{ \Alg_{\calM_1}(\calC)} \Alg_{1,2}(\calC) \rightarrow \Alg_{\calM_0}(\calC)$.

Let $\Alg_{0,1,2}(\calC)$ denote the full subcategory of
$\Fun_{ \Nerve(\FinSeg)}(\calM^{\otimes}, \calC^{\otimes})$ spanned by the $\infty$-operad
family maps whose restrictions to $\calM^{\otimes} \times_{ \Delta^2} \Delta^{ \{0,1\} }$ and
$\calM^{\otimes} \times_{ \Delta^2} \Delta^{ \{1,2\} }$ are operadic $q$-left Kan extensions.
Condition $(iii)$ guarantees that the inclusion $\calM^{\otimes} \times_{ \Delta^2} \Lambda^2_1 \subseteq \calM^{\otimes}$ is a categorical equivalence, so that the restriction maps
$$ \Fun_{ \Nerve(\FinSeg)}(\calM^{\otimes}, \calC^{\otimes}) \rightarrow
\Fun_{ \Nerve(\FinSeg)}( \calM^{\otimes} \times_{ \Delta^2} \Lambda^2_1, \calC^{\otimes})$$
$$ \Alg_{0,1,2}(\calC) \rightarrow \Alg_{0,1}(\calC) \times_{ \Alg_{ \calM_1}(\calC)} \Alg_{1,2}(\calC)$$
are trivial Kan fibrations. It follows that the restriction map $r'': \Alg_{0,1,2}( \calC) \rightarrow
\Alg_{ \calM_0}(\calC)$ is a trivial Kan fibration admitting a section $s''$, and that
$f_{1,2} \circ f_{0,1}$ can be identified with the composition
$$ \Alg_{\calM_0}(\calC) \stackrel{s''}{\rightarrow} \Alg_{0,1,2}(\calC) \rightarrow \Alg_{\calM_2}(\calC).$$

Using Theorem \ref{transit}, we deduce that the restriction map $\Alg_{0,1,2}(\calC) \rightarrow
\Alg_{\calM_2}(\calC)$ factors as a composition
$$ \Alg_{0,1,2}(\calC) \stackrel{\theta}{\rightarrow} \Alg_{0,2}(\calC) \stackrel{\theta'}{\rightarrow} \Alg_{ \calM_2}(\calC).$$ 
The composition $\theta \circ s''$ is a section of the trivial Kan fibration $\Alg_{0,2}(\calC) \rightarrow
\Alg_{\calM_0}(\calC)$, so that $f_{0,2}$ can be identified with the composition
$\theta' \circ (\theta \circ s'') \simeq f_{1,2} \circ f_{0,1}$ as desired.
\end{proof}

The proof of Theorem \ref{transit} rests on a more basic transitivity property of operadic
colimit diagrams. To state this property, we need to introduce a bit of terminology.
Let $q: \calC^{\otimes} \rightarrow \calD^{\otimes}$ be a fibration of $\infty$-operads,
and let $p: K \diamond \Delta^0 \rightarrow \calC^{\otimes}$ be a map of simplicial
sets which carries each edge of $K \diamond \Delta^0$ to an active morphism in $\calC^{\otimes}$.
Since the map $K \diamond \Delta^0 \rightarrow K^{\triangleright}$ is a categorical equivalence (Proposition \toposref{rub3}), there exists a map $p': K^{\triangleright} \rightarrow \calC^{\otimes}$ such that $p$ is homotopic to the composition $K \diamond \Delta^0 \rightarrow K^{\triangleright} \stackrel{p'}{\rightarrow} \calC^{\otimes}$. Moreover, the map $p'$ is unique up to homotopy. We will say that
$p$ is a {\it (weak) operadic $q$-colimit diagram} if $p'$ is a (weak) operadic $q$-colimit diagram, in the sense of
Definition \symmetricref{swm}.

\begin{lemma}\label{prelax}
Let $X \rightarrow S$ be a coCartesian fibration of simplicial sets, and let
$q: \calC^{\otimes} \rightarrow \calD^{\otimes}$ be a fibration of $\infty$-operads.
Let $$\theta: (X \diamond_{S} S) = (X \times \Delta^1) \coprod_{ X \times \{1\} } S
\rightarrow \calC^{\otimes}$$
be a map satisfying the following conditions:
\begin{itemize}
\item[$(i)$] The map $\theta$ carries every edge in $X \diamond_{S} S$ to an
active morphism in $\calC^{\otimes}$.
\item[$(ii)$] For every vertex $s \in S$, the induced map
$\theta_{s}: X_{s} \diamond \Delta^0 \rightarrow \calC^{\otimes}$ is a weak operadic $q$-colimit diagram.
\end{itemize}
Let $\theta_0 = \theta | X$. Let $\calC^{\acti}_{\theta/}$ denote the full subcategory of
$\calC^{\otimes}_{\theta/} \times_{ \calC^{\otimes} } \calC$ spanned by those objects
which correspond to maps $\overline{\theta}: ( X \diamond_{S} S)^{\triangleright} \rightarrow
\calC^{\otimes}$ which carry every edge of $(X \diamond_{S} S)^{\triangleright}$ to an inert
morphism of $\calC^{\otimes}$, and define $\calC^{\acti}_{\theta_0/}$,
$\calD^{\acti}_{q \theta/}$, and $\calD^{\acti}_{q \theta_0/}$ similarly. Then:
\begin{itemize}
\item[$(1)$] The map $\calC^{\acti}_{\theta/} \rightarrow
\calC^{\acti}_{\theta_0/} \times_{ \calD^{\acti}_{q \theta_0/} } \calD^{\acti}_{q \theta/}$
is a trivial Kan fibration.
\item[$(2)$] Let $\overline{\theta}: (X \diamond_{S} S)^{\triangleright} \rightarrow \calC^{\otimes}$
be an extension of $\theta$ which carries each edge of $(X \diamond_{S} S)^{\triangleright}$ to an active morphism in $\calC^{\otimes}$. Then $\overline{\theta}$ is a weak operadic $q$-colimit diagram
if and only if $\overline{\theta}_0 = \overline{\theta} | X^{\triangleright}$ is a weak operadic
$q$-colimit diagram.
\item[$(3)$] Assume that each $\theta_{s}$ is an operadic $q$-colimit diagram, and let
$\overline{\theta}$ be as in $(2)$. Then $\overline{\theta}$ is an operadic $q$-colimit diagram
if and only if $\overline{\theta}_0$ is an operadic $q$-colimit diagram.
\end{itemize}
\end{lemma}

\begin{proof}
Assertion $(2)$ follows immediately from $(1)$, and assertion $(3)$ follows from
$(2)$ after replacing $\theta$ by the composite functor
$$ X \diamond_{S} S \stackrel{\theta}{\rightarrow} \calC^{\otimes} \stackrel{ \oplus Y}{\rightarrow} \calC^{\otimes},$$
where $Y$ denotes an arbitrary object of $\calC^{\otimes}$. It will therefore suffice to prove
$(1)$. For every map of simplicial sets $K \rightarrow S$, let $\theta_{K}$ denote the induced map
$X \diamond_{S} K \rightarrow \calC^{\otimes}$. We will prove more generally that for
$K' \subseteq K$, the induced map
$$ \psi_{K', K}: 
\calC^{\acti}_{\theta_K/} \rightarrow
\calC^{\acti}_{\theta_{K'}/} \times_{ \calD^{\acti}_{q \theta_{K'}/} } \calD^{\acti}_{q \theta_K/}$$
is a trivial Kan fibration. We proceed by induction on the (possibly infinite) dimension $n$ of $K$.
If $K$ is empty, the result is obvious. Otherwise, working simplex-by-simplex, we can assume that
$K$ is obtained from $K'$ by adjoining a single nondegenerate $m$-simplex $\sigma$ whose boundary already belongs to $K'$. Replacing $K$ by $\sigma$, we may assume that $K = \Delta^m$ and 
$K' = \bd \Delta^m$. If $m = 0$, then the desired result follows from assumption $(ii)$.
Assume therefore that $m > 0$.

Because $\theta_{K',K}$ is clearly
a categorical fibration (even a left fibration), to prove that $\theta_{K',K}$ is a trivial Kan fibration
it suffices to show that $\theta_{K', K}$ is a categorical equivalence. 
Since $m \leq n$, $K'$ has dimension $< n$, so the inductive hypothesis guarantees that $\psi_{\{m\}, K'}$ is a trivial Kan fibration. The map $\psi_{\{m\}, K}$ is a composition of $\psi_{K',K}$ with a pullback of
$\psi_{\{m\},K'}$. Using a two-out-of-three argument, we are reduced to proving that
$\psi_{\{m\}, K}$ is a categorical equivalence. For this, it suffices to show that the inclusion
$f: X \diamond_{S} \{m\} \rightarrow X \diamond_{S} \Delta^m$ is cofinal. 

Let $X' = X \times_{S} \Delta^m$. The map $f$ is a pushout of the inclusion
$$ f': X' \coprod_{ X'_{m}} (X'_{m} \diamond \{m\}) \hookrightarrow X' \diamond_{ \Delta^m} \Delta^m.$$
It will therefore suffice to show that $f'$ is cofinal. We have a commutative diagram
$$ \xymatrix{  \{m\} \ar[r]^{f''} \ar[d]^{g} & \Delta^m \ar[dd]^{g''} \\
X'_{m} \diamond \{m\} \ar[d]^{g'} & \\
X' \coprod_{ X'_{m} } (X'_{m} \diamond \{m\}) \ar[r]^-{f'} & X' \diamond_{ \Delta^m} \Delta^m. }$$
The map $g''$ is a pushout of the inclusion $X' \times \{1\} \subseteq X' \times \Delta^1$,
and therefore cofinal; the same argument shows that $g$ is cofinal. The map
$f''$ is obviously cofinal. The map $g'$ is a pushout of the inclusion
$X'_{m} \subseteq X'$, which is cofinal because $\{m\}$ is a final object of $\Delta^m$ and the map
$X' \rightarrow \Delta^m$ is a coCartesian fibration. It now follows from Proposition \toposref{cofbasic} that $f'$ is cofinal, as required.
\end{proof}

\begin{proof}[Proof of Theorem \ref{transit}]
Fix an object $Z \in \calM^{\otimes}_2$, and let
$\calZ$ denote the full subcategory of $\calM^{\otimes}_{/Z}$ whose objects are active morphisms
$X \rightarrow Z$ where $X \in \calM^{\otimes}_{0}$. We wish to prove that the composite map
$$ \phi: \calZ^{\triangleright} \rightarrow (\calM^{\otimes}_{/Z})^{\triangleright}
\rightarrow \calM^{\otimes} \stackrel{A}{\rightarrow} \calC^{\otimes}$$
is an operadic $q$-colimit diagram. Let $\overline{\calZ}$ denote the subcategory of
$\Fun( \Delta^1, \calM^{\otimes}_{/Z})$ whose objects are diagrams
of active morphisms $$ \xymatrix{ & Y \ar[dr] & \\
X \ar[ur] \ar[rr] & & Z}$$
in $\calM^{\otimes}$ such that $X \in \calM^{\otimes}_0$ and $Y \in \calM^{\otimes}_1$.
Evaluation at $\{0\}$ induces a Cartesian fibration $\psi: \overline{\calZ} \rightarrow \calZ$. 
Let $Z'$ be an object of $\calZ$, corresponding to an active morphism
$X \rightarrow Z$ in $\calM^{\otimes}$. Then the fiber $\psi^{-1} \{ Z' \}$ is
a localization of the $\infty$-category $(\calM^{\otimes}_{/Z})^{X/}$, which
is equivalent to $(\calM^{\otimes})_{X/ \, /Z}$ and therefore weakly contractible
(since $\calM^{\otimes} \rightarrow \Delta^2$ is flat).
Note that the map $\psi': \overline{\calZ} \times_{ \calZ} \calZ_{Z'/}
\rightarrow \calZ_{Z'/}$ is a Cartesian fibration (Proposition \toposref{verylonger}).
Since $\calZ_{Z'/}$ has an initial object $\id_{Z'}$, the weakly contractible simplicial set
$\psi^{-1} \{ Z' \}  \simeq {\psi'}^{-1} \{ \id_{Z'} \}$ is weakly homotopy equivalent to 
$\overline{\calZ} \times_{ \calZ} \calZ_{Z'/}$. Applying Theorem \toposref{hollowtt}, we deduce
that $\psi$ is cofinal. Consequently, it will suffice to show that $\phi \circ \psi: \overline{\calZ}^{\triangleright} \rightarrow \calC^{\otimes}$ is an operadic $q$-colimit diagram.

Let $\calY$ denote the full subcategory of $\calM^{\otimes}_{/Z}$ spanned by active morphisms $Y \rightarrow Z$ where $Y \in \calM^{\otimes}_{1}$. Evaluation at $\{1\}$ induces a coCartesian
fibration $\rho: \overline{\calZ} \rightarrow \calY$. We observe that there is a canonical map
$\overline{ \calZ} \diamond_{ \calY} \calY \rightarrow \calM^{\otimes}_{/Z}$, which
determines a map
$$ \theta: ( \overline{\calZ} \diamond_{\calY} \calY)^{\triangleright} \rightarrow \calC^{\otimes}$$
extending $\phi \circ \psi$. Fix an object $Y' \in \calY$, corresponding to an active
morphism $Y \rightarrow Z$ in $\calM^{\otimes}$. Then $\theta$ induces a map
$\theta_{Y'}: \rho^{-1} \{Y' \} \diamond \Delta^0 \rightarrow \calC^{\otimes}$. We claim that
$\theta_{Y'}$ is an operadic $q$-colimit diagram. To prove this, let $\calX(Y)$
denote the full subcategory of $(\calM^{\otimes})^{/Y}$ spanned by the active morphisms $X \rightarrow Y$, and define $\calX'(Y) \subseteq (\calM^{\otimes})_{/Y}$ similarly. The map $\theta_{Y'}$ factors
through a map
$$\theta'_{Y'}: \calX(Y) \diamond \Delta^0 \rightarrow \calC^{\otimes}.$$
Since $\calM^{\otimes}_{/Z} \rightarrow \calM^{\otimes}$ is a left fibration, the map
$\rho^{-1} \{Y'\} \rightarrow \calX(Y)$ is a trivial Kan fibration; it therefore suffices to show that
$\theta'_{Y'}$ is an operadic $q$-colimit diagram. Since the evident map
$\calX'(Y) \rightarrow \calX(Y)$ is a categorical equivalence (Proposition \toposref{certs}),
it suffices to show that the induced map $\calX'(Y) \diamond \Delta^0 \rightarrow \calC^{\otimes}$
is an operadic $q$-colimit diagram, which is equivalent to the requirement that the composite map
$$ \calX'(Y)^{\triangleright} \rightarrow (\calM^{\otimes}_{/Y})^{\triangleright} \rightarrow \calM^{\otimes} \rightarrow \calC^{\otimes}$$
is an operadic $q$-colimit diagram. This follows from our assumption that
$A | (\calM^{\otimes} \times_{ \Delta^2} \Delta^{ \{0,1\} })$ is an operadic $q$-left Kan extension.

Since $A | ( \calM^{\otimes} \times_{ \Delta^2} \Delta^{ \{0,1\} })$, the restriction of
$\theta$ to $\calY^{\triangleright}$ is an operadic $q$-colimit diagram.
The inclusion $\calY \rightarrow \overline{\calZ} \diamond_{\calY} \calY$
is a pushout of the inclusion $\overline{\calZ} \times \{1\} \subseteq \overline{\calZ} \times \Delta^1$, and therefore cofinal. It follows that $\theta$ itself is an operadic $q$-colimit diagram.
Invoking Lemma \ref{prelax}, we conclude that $\phi \circ \psi$ is an operadic $q$-colimit diagram,
as desired.
\end{proof}

\subsection{A Coherence Criterion}\label{cohcrit}

In \S \symmetricref{moduldef}, we introduced the notion of a {\it coherent} $\infty$-operad (Definition \symmetricref{koopa}), and showed that
if $\calO^{\otimes}$ is coherent then there is a good notion of module over every $\calO$-algebra. 
However, the definition presented there is somewhat cumbersome. Our goal in this section is to give a characterization of coherent $\infty$-operads which is easy to verify in practice. The principal application is the verification that the little cubes operads $\OpE{k}$ are coherent: see \S \ref{slabba}. 
We begin with an informal sketch of the basic idea.

Fix an active morphism $f: X \rightarrow Y$ in $\calO^{\otimes}$. An {\it extension} of
$f$ consists of an object $X_0 \in \calO$ together with an active morphism
$f^{+}: X \oplus X_0 \rightarrow Y$ such that $f^{+}|X \simeq f$; here
the hypothesis that $\calO^{\otimes}$ is unital guarantees that there is an essentially
section to the projection $X \oplus X_0 \rightarrow X$, so that the restriction
is well-defined. The collection of extensions of $f$ can be organized into an $\infty$-category,
which we will denote by $\Ext(f)$ (see Definition \ref{cubulos} below for a precise definition).

If $g: Y \rightarrow Z$ is another active morphism, then there are canonical maps
$$ \Ext(f) \rightarrow \Ext( g \circ f) \leftarrow \Ext(g),$$
well-defined up to homotopy. In particular, we have canonical maps
$$ \Ext(f) \leftarrow \Ext( \id_{Y} ) \rightarrow \Ext(g)$$
which fit into a (homotopy coherent) diagram
$$ \xymatrix{ & \Ext( \id_Y) \ar[dl] \ar[dr] & \\
\Ext(f) \ar[dr] & & \Ext(g) \ar[dl] \\
& \Ext(g \circ f). & }$$
If we assume that $\calO$ is a Kan complex, then the $\infty$-categories appearing in this diagram are all Kan complexes. Our main result assert that in this case, the $\infty$-operad $\calO^{\otimes}$ is coherent if and only if this diagram is a pushout square, for every composable pair of
active morphisms $f: X \rightarrow Y$ and $g: Y \rightarrow Z$ in $\calO^{\otimes}$ (Theorem \ref{uggus}).
In fact, it suffices to check this condition in the special case where $Z \in \calO$.

We begin by giving a careful definition of the $\infty$-category $\Ext(f)$.

\begin{definition}\label{cubulos}
Let $q: \calO^{\otimes} \rightarrow \Nerve(\FinSeg)$ be a unital $\infty$-operad, and let
$\sigma: \Delta^n \rightarrow \calO^{\otimes}_{\acti}$ be an
$n$-simplex of $\calO^{\otimes}$ corresponding to a composable chain
$X_0 \stackrel{f_1}{\rightarrow} \ldots \stackrel{f_n}{\rightarrow} X_n$ of
active morphisms in $\calO^{\otimes}$.

If $S \subseteq [n]$ is a downward-closed subset, we let $\Ext(\sigma,S)$ denote the full subcategory of 
$\Fun( \Delta^n, \calO^{\otimes})_{\sigma/}$ spanned by those diagrams
$$ \xymatrix{ X_0 \ar[r]^{f_1} \ar[d]^{g_0} & \cdots \ar[r]^{f_n} & X_{n} \ar[d]^{g_n} \\
X'_0 \ar[r]^{f'_1} & \cdots \ar[r]^{f'_n} & X'_{n} }$$
with the following properties:
\begin{itemize}
\item[$(a)$] If $i \notin S$, the map $g_{i}$ is an equivalence.
\item[$(b)$] If $i \in S$, the map $g_{i}$ is semi-inert and
$q(g_{i})$ is an inclusion $\seg{n_{i}} \rightarrow \seg{n_i+1}$ which omits a single
value $a_{i} \in \seg{n_i+1}$.
\item[$(c)$] If $0 < i \in S$, then the map $q( f_{i})$ carries $a_{i-1} \in q(X'_{i-1})$
to $a_{i} \in q( X'_{i})$.
\item[$(d)$] Each of the maps $f'_{i}$ is active.
\end{itemize}

If $f: \Delta^1 \rightarrow \calO^{\otimes}_{\acti}$ is an active morphism in
$\calO^{\otimes}$, we will denote $\Ext(f, \{0\})$ by $\Ext(f)$.
\end{definition}

\begin{remark}
Let $\Ext(\sigma, S)$ be as in Definition \ref{cubulos}. If $\calO$ is a Kan complex, then it is easy to
see that every morphism in $\Ext(\sigma,S)$ is an equivalence, so that $\Ext(\sigma,S)$ is also a Kan complex.
\end{remark}

\begin{remark}\label{sadlus}
Let $\sigma: \Delta^n \rightarrow \calO^{\otimes}_{\acti}$ correspond to a sequence of active morphisms
$$ X_0 \stackrel{f_1}{\rightarrow} X_1 \stackrel{f_2}{\rightarrow} \cdots \stackrel{f_n}{\rightarrow} X_n$$
and let $S \subseteq [n]$. For every morphism $j: \Delta^{m} \rightarrow \Delta^n$, composition with $j$ induces a restriction map $$ \Ext( \sigma, S) \rightarrow \Ext( \sigma \circ j, j^{-1}(S)).$$

In particular, if $S$ has a largest element $i < n$, then we obtain a canonical map
$\Ext( \sigma, S) \rightarrow \Ext( f_{i+1}).$ If $\calO$ is a Kan complex, then this map is a trivial Kan fibration.
\end{remark}

\begin{remark}\label{sped}
Let $q: \calO^{\otimes} \rightarrow \Nerve(\FinSeg)$ be an $\infty$-operad, 
let $\sigma: \Delta^n \rightarrow \calO^{\otimes}_{\acti}$ correspond to a sequence of active morphisms
$$ X_0 \stackrel{f_1}{\rightarrow} X_1 \stackrel{f_2}{\rightarrow} \cdots \stackrel{f_n}{\rightarrow} X_n,$$
let $\seg{m} = q( X_n)$, and let $S$ be a nonempty proper subset of $[n]$. Then $\Ext(\sigma,S)$ decomposes naturally as a disjoint union $\coprod_{ 1 \leq i \leq m} \Ext(\sigma,S)_{i}$, where
$\Ext(\sigma,S)_{i}$ denotes the full subcategory of $\Ext(\sigma, S)$ spanned by those diagrams
$$ \xymatrix{ X_0 \ar[r]^{f_1} \ar[d]^{g_0} & \cdots \ar[r] & X_{n} \ar[d]^{g_n} \\
X'_0 \ar[r]^{f'_1} & \cdots \ar[r] & X'_{n} }$$
where $q( f'_{n} \cdots f'_{1})$ carries the unique element of $q(X'_0) - q(X_0)$
to $i \in \seg{m} \simeq q( X'_{n})$. In this case, the diagram $\sigma$ is equivalent
to an amalgamation $\bigoplus_{1 \leq i \leq m} \sigma_{i}$, and we have canonical
equivalences $\Ext(\sigma, S)_{i} \simeq \Ext(\sigma_i,S)$.
\end{remark}

\begin{remark}\label{preug}
Let $f: \calO^{\otimes} \rightarrow {\calO'}^{\otimes}$ be an ornamental map between unital $\infty$-operads, where
$\calO$ and $\calO'$ are Kan complexes. Let $\sigma$ be an $n$-simplex of $\calO^{\otimes}_{\active}$ and
let $S$ be a downward-closed subset of $[n]$. Then $f$ induces a homotopy equivalence
$\Ext( \sigma,S) \rightarrow \Ext( f(\sigma), S)$. This follows from
Remark \ref{sadlus} and Lemma \ref{snoke}.
\end{remark}

We can now state our main result precisely as follows:

\begin{theorem}\label{uggus}
Let $\calO^{\otimes}$ be a unital $\infty$-operad such that $\calO$ is a Kan complex. The following conditions are equivalent:
\begin{itemize}
\item[$(1)$] The $\infty$-operad $\calO^{\otimes}$ is coherent.
\item[$(2)$] Suppose we are given a degenerate 
$3$-simplex $\sigma:$
$$ \xymatrix{ & Y \ar[dr]^{\id_Y} \ar[rr]^{g} & & Z \\
X \ar[ur]^{f} \ar[rr]^{f} & & Y \ar[ur]^{g}}$$
in $\calO^{\otimes}$, where $f$ and $g$ are active. Then the diagram
$$ \xymatrix{ \Ext( \sigma, \{0,1\} ) \ar[r] \ar[d] & \Ext( \sigma | \Delta^{ \{0,1,3 \} }, \{0,1\} ) \ar[d] \\
\Ext( \sigma | \Delta^{ \{0,2,3\} }, \{0\} ) \ar[r] & \Ext( \sigma| \Delta^{ \{0,3\} }, \{0\} ) }$$
is a homotopy pushout square.
\item[$(3)$] Conclusion $(2)$ holds whenever $Z \in \calO$.
\end{itemize}
\end{theorem}

\begin{remark}
In view of Remark \ref{sadlus}, we can think of the diagram appearing in the statement
of Theorem \ref{uggus} as giving a homotopy coherent diagram
$$ \xymatrix{ \Ext(\id_Y) \ar[r] \ar[d] & \Ext(g) \ar[d] \\
\Ext(f) \ar[r] & \Ext( g \circ f). }$$
\end{remark}

\begin{remark}\label{saym}
The implication $(2) \Rightarrow (3)$ in Theorem \ref{uggus} is immediate, and the converse
implication follows from Remark \ref{sped}.
\end{remark}

\begin{corollary}\label{kisu}
Let $f: \calO^{\otimes} \rightarrow {\calO'}^{\otimes}$ be an ornamental map between
unital $\infty$-operads, where both $\calO$ and $\calO'$ are Kan complexes. If
${\calO'}^{\otimes}$ is coherent, then $\calO^{\otimes}$ is coherent. The converse holds if
the underlying map $\pi_0 \calO \rightarrow \pi_0 \calO'$ is surjective.
\end{corollary}

\begin{proof}
We may assume without loss of generality that $f$ is a categorical fibration.
The first assertion follows immediately from Theorem \ref{uggus} and Remark \ref{preug}.
To prove the second, it will suffice (by virtue of Theorem \ref{uggus} and Remark \ref{preug}) to show
that every $3$-simplex $\sigma: \Delta^3 \rightarrow {\calO'}^{\otimes}_{\acti}$ can be lifted
to a $3$-simplex of $\calO^{\otimes}_{\active}$. Let $Z = \sigma(3) \in \calO'$. Since $f$ is a categorical fibration, the induced map $\calO \rightarrow \calO'$ is a Kan fibration. Since this Kan fibration induces a surjection on connected components, it is surjective on vertices and we may write $Z = f( \overline{Z})$ for
some $\overline{Z} \in \calO$. Lemma \ref{snoke} guarantees that the induced map
$f': (\calO^{\otimes}_{/ \overline{Z} })^{\acti} \rightarrow ( {\calO'}^{\otimes}_{/Z})^{\acti}$ is a categorical equivalence. Since $f$ is a categorical fibration, $f'$ is also a categorical fibration and therefore a trivial Kan fibration. We can interpret the $3$-simplex $\sigma$ as a $2$-simplex $\tau: \Delta^2 \rightarrow
( {\calO'}^{\otimes}_{/Z})^{\acti}$, which can be lifted to a $2$-simplex $\overline{\tau}: \Delta^2
\rightarrow (\calO^{\otimes}_{/ \overline{Z} })^{\acti}$. This map determines a $3$-simplex of
$\calO^{\otimes}_{\acti}$ lifting $\sigma$, as desired.
\end{proof}

\begin{remark}
Let $\calO^{\otimes}$ be a unital $\infty$-operad such that $\calO$ is a Kan complex. According to Theorem \ref{assum}, the $\infty$-operad $\calO^{\otimes}$ can be obtained as the assembly of a $\calO$-family of
reduced unital $\infty$-operads ${\calO'}^{\otimes} \rightarrow \calO \times \Nerve( \FinSeg)$. Corollary
\ref{kisu} (and Proposition \ref{loam}) imply that $\calO^{\otimes}$ is coherent if and only if, for each
$X \in \calO$, the $\infty$-operad ${\calO'}^{\otimes}_{X}$ is coherent. In this case, there is a good theory of
modules associated to $\calO$-algebras and to $\calO'_{X}$-algebras, for each $X \in \calO$. 
One can show that these module theories are closely related to one another. To describe the relationship, suppose that $\calC$ is another $\infty$-operad and $A \in \Alg_{\calO}(\calC)$ is a $\calO$-algebra object of $\calC$, corresponding to a family of ${\calO'}^{\otimes}_{X}$-algebra objects $\{ A_X \in \Alg_{ \calO'_X}(\calC) \}_{X \in \calO}$. Then giving an $A$-module object $M \in \Mod^{\calO}_{A}(\calC)$ is equivalent to giving a family $\{ M_X \in \Mod^{\calO'_{X}}_{A_X}(\calC) \}_{X \in \calO'}$. We leave the precise formulation to the reader. 
\end{remark}

The rest of this section is devoted to the proof of Theorem \ref{uggus}. Our first step is to introduce an apparently weaker notion of coherence.

\begin{definition}
Let $q: \calO^{\otimes} \rightarrow \Nerve(\FinSeg)$ be a unital $\infty$-operad.
We will say that a morphism $f: X \rightarrow X'$ in $\calO^{\otimes}$ is
{\it $m$-semi-inert} if $f$ is semi-inert and the underlying map $q(f): \seg{n} \rightarrow \seg{n'}$
is such that the cardinality of the set $\seg{n'} - f( \seg{n} )$ is less than or equal to $m$.
\end{definition}

By definition, a unital $\infty$-operad $\calO^{\otimes}$ is coherent if, for every map
$\Delta^2 \rightarrow \calO^{\otimes}$ corresponding to a diagram
$X \rightarrow Y \rightarrow Z$ and every morphism $\overline{X} \rightarrow \overline{Z}$ in
$\calK_{\calO}$ lifting the underlying map $X \rightarrow Z$, the $\infty$-category
$( \calK_{\calO})_{ \overline{X} / \, / \overline{Z} } \times_{ \calO^{\otimes}_{X/ \, /Z} } \{Y\}$ is
weakly contractible. We will say that $\calO^{\otimes}$ is {\it $m$-coherent} if
this condition holds whenever $\overline{X}$ is $m$-semi-inert. Note that
$\calO^{\otimes}$ is coherent if and only if it is $m$-coherent for all $m \geq 0$.

\begin{lemma}\label{stokums}
Let $\calO^{\otimes}$ be a unital $\infty$-operad. The following conditions are equivalent:
\begin{itemize}
\item[$(1)$] The $\infty$-operad $\calO^{\otimes}$ is $m$-coherent.
\item[$(2)$] Consider a diagram $\sigma:$
$$ \xymatrix{ X \ar[d]^{f} \ar[r] & Y \\
X' & }$$ in $\calO^{\otimes}$, where $f$ is $m$-semi-inert. Let
$\calA[\sigma]$ denote the full subcategory of $\calO^{\otimes}_{\sigma/}$ spanned by those commutative diagrams
$$ \xymatrix{ X \ar[d]^{f} \ar[r] & Y \ar[d]^{g} \\
X' \ar[r] & Y' }$$
where $g$ is also semi-inert. Then the inclusion 
$\calA[\sigma]^{op} \subseteq ( \calO^{\otimes}_{\sigma/})^{op}$
is cofinal.

\item[$(3)$] Let $\sigma$ be as in $(2)$, and let $Z'$ be an object of
$\calO^{\otimes}_{\sigma/}$, and let $\calA[\sigma]_{/Z'}$ denote the full subcategory of
$\calO^{\otimes}_{\sigma/ \, /Z'}$ spanned by those diagrams
$$ \xymatrix{ X \ar[d]^{f} \ar[r] & Y \ar[d]^{g} \ar[dr] & \\
X' \ar[r] & Y' \ar[r] & Z'}$$
such that $g$ is semi-inert. Then $\calA[\sigma]_{/Z'}$ is a weakly contractible simplicial set.
\end{itemize}
\end{lemma}

\begin{proof}
The equivalence of $(2)$ and $(3)$ follows immediately from Theorem \toposref{hollowtt}.
We next prove that $(1) \Rightarrow (3)$. We can extend $\sigma$ and $Z$ to a commutative
diagram
$$ \xymatrix{ X \ar[r] \ar[d]^{f} & Y \ar[r] & Z \ar[d]^{h} \\
X' \ar[rr] & & Z' }$$
where $h$ is semi-inert (for example, we can take $Z = Z'$ and $h = \id_{Z'}$.
The upper line of this diagram determines a diagram $\Delta^2 \rightarrow \calO^{\otimes}$.
Let $\calK$ denote the fiber product $\calK_{\calO} \times_{ \calO^{\otimes}} \Delta^2$.
Since $\calO^{\otimes}$ is coherent, the projection map $\pi: \calK \rightarrow \Delta^2$ is a flat inner fibration. The maps $f$ and $h$ determine objects of $\calK \times_{ \Delta^2} \{0\}$ and
$\calK \times_{ \Delta^2} \{2\}$, which we will denote by $\overline{X}$ and $\overline{Z}$.
Since $\pi$ is flat, the $\infty$-category $\calK_{ \overline{X}/ \, / \overline{Z} } \times_{ \Delta^2 } \{1\}$ is weakly contractible. We now observe that there is a trivial Kan fibration $\psi: \calK_{ \overline{X}/ \, /\overline{Z}} \rightarrow \calA[\sigma]_{/Z'}$, so that $\calA[\sigma]_{/Z'}$ is likewise weakly contractible.

Now suppose that $(3)$ is satisfied. We wish to show that evaluation at $\{0\} \subseteq \Delta^1$
induces a flat categorical fibration $\calK_{\calO} \rightarrow \calO^{\otimes}$. Fix a map
$\Delta^2 \rightarrow \calO^{\otimes}$, and let $\calK$ be the fiber product
$\calK_{\calO} \times_{ \calO^{\otimes} } \Delta^2$. Suppose we are given objects
$\overline{X} \in \calK \times_{ \Delta^2 } \{0\}$ and $\overline{Z} \in \calK \times_{ \Delta^2 } \{2\}$; we wish to prove that if $\overline{X}$ is $m$-semi-inert, then $\calK_{ \overline{X}/ \, /\{Z\} } \times_{ \Delta^2} \{1\}$ is weakly contractible. The above data determines a commutative diagram
$$ \xymatrix{ X \ar[r] \ar[d]^{f} & Y \ar[r] & Z \ar[d]^{h} \\
X' \ar[rr] & & Z' }$$
in $\calO^{\otimes}$. If we let $\sigma$ denote the left part of this diagram, then
we can define a simplicial set $\calA[\sigma]_{/Z'} \subseteq ( \calO^{\otimes}_{\sigma/} )_{/Z'}$ as in $(3)$,
which is weakly contractible by assumption. Once again, we have a trivial Kan fibration
$\psi: \calK_{ \overline{X}/ \, /\overline{Z}} \rightarrow \calA[\sigma]_{/Z'}$, so that
$\calK_{ \overline{X}/ \, /\overline{Z}}$ is weakly contractible as desired.
\end{proof}

\begin{remark}\label{funly}
Let $\calO^{\otimes}$, $\sigma$, and $Z'$ be as in the part $(3)$ of Lemma \ref{stokums}.
Let $\calA[\sigma]'_{/Z}$ denote the full subcategory of $\calA[\sigma]_{/Z}$ spanned by those diagrams
$$ \xymatrix{ X \ar[d]^{f} \ar[r] & Y \ar[d]^{g} \ar[dr] & \\
X' \ar[r] & Y' \ar[r]^{j} & Z'}$$
for which the map $j$ is active. The inclusion $\calA[\sigma]'_{/Z} \subseteq \calA[\sigma]_{/Z}$ admits a left adjoint, and is therefore a weak homotopy equivalence. Consequently, condition $(3)$ of
Lemma \ref{stokums} is equivalent to the requirement that $\calA[\sigma]'_{/Z}$ is weakly contractible.
\end{remark}

\begin{remark}\label{saigr}
In the situation of Lemma \ref{stokums}, let $Z' \in \calO^{\otimes}$ and let
$\calB$ denote the full subcategory of $\calO^{\otimes}_{/Z'}$ spanned by the active morphisms
$W \rightarrow Z'$. The inclusion $\calB \subseteq \calO^{\otimes}_{/Z'}$ admits a left adjoint $L$.
For any diagram $\sigma:$
$$ \xymatrix{ X \ar[r] \ar[d] & Y \\
X' & }$$
in $\calO^{\otimes}_{/Z}$, $L$ induces a functor $\calA[\sigma]_{/Z'} \rightarrow
\calA[ L \sigma]_{/Z'}$, which restricts to an equivalence
$\calA[\sigma]'_{/Z'} \rightarrow \calA[ L \sigma]'_{/Z'}$,
where $\calA[\sigma]'_{/Z'}$ and $\calA[L \sigma]'_{/Z'}$ are defined as in Remark
\ref{funly}. Consequently, to verify condition $(3)$ of Lemma \ref{stokums}, we are free to replace $\sigma$ by $L \sigma$ and thereby reduce to the case where the maps
$X \rightarrow Z'$, $Y \rightarrow Z'$, and $X' \rightarrow Z'$ are active.
\end{remark}

\begin{remark}\label{saig2}
Let $Z' \in \calO^{\otimes}$ be as in Lemma \ref{stokums}, and choose an equivalence
$Z' \simeq \bigoplus Z'_{i}$, where $Z'_{i} \in \calO$. Let $\calB \subseteq \calO^{\otimes}_{/Z'}$ be defined as in Remark \ref{saigr}, and let $\calB_{i} \subseteq \calO^{\otimes}_{/Z'_{i}}$ be defined similarly. Every diagram $\sigma: \Lambda^2_0 \rightarrow \calB$ can be identified with
an amalgamation $\bigoplus_{i} \sigma_i$ of diagrams $\sigma_i: \Lambda^2_0 \rightarrow \calB_i$.
We observe that $\calA[ \sigma]'_{/Z'}$ is equivalent to the product of the $\infty$-categories
$\calA[ \sigma_i]'_{/Z'_i}$. Consequently, to verify that $\calA[\sigma_i]'_{/Z'}$ is weakly contractible,
we may replace $Z'$ by $Z'_i$ and thereby reduce to the case where $Z' \in \calO$.
\end{remark}

\begin{lemma}\label{clapsos}
Let $q: X \rightarrow S$ be an inner fibration of $\infty$-categories. Suppose that the following conditions are satisfied:
\begin{itemize}
\item[$(a)$] The inner fibration $q$ is flat.
\item[$(b)$] Each fiber $X_{s}$ of $q$ is weakly contractible.
\item[$(c)$] For every vertex $x \in X$, the induced map
$X_{x/} \rightarrow S_{q(x)/}$ has weakly contractible fibers.
\end{itemize}
Then for every map of simplicial sets $S' \rightarrow S$, the pullback map
$X \times_{S} S' \rightarrow S'$ is weak homotopy equivalence. In particular, $q$ is a weak homotopy equivalence.
\end{lemma}

\begin{proof}
We will show more generally that for every map of simplicial sets $S' \rightarrow S$, the induced map
$q_{S'}: X \times_{S} S' \rightarrow S'$ is a weak homotopy equivalence. Since the collection of
weak homotopy equivalences is stable under filtered colimits, we may suppose that $S'$ is finite.
We now work by induction on the dimension $n$ of $S'$ and the number of simplices of $S'$ of maximal dimension. If $S'$ is empty the result is obvious; otherwise we have a homotopy pushout diagram
$$ \xymatrix{ \bd \Delta^n \ar[r] \ar[d] & \Delta^n \ar[d] \\
S'_0 \ar[r] & S'. }$$
By the inductive hypothesis, the maps $q_{S'_0}$ and $q_{ \bd \Delta^n}$ are weak homotopy equivalences. Since $q_{S'}$ is a homotopy pushout of the morphisms
$q_{S'_0}$ with $q_{\Delta^n}$ over $q_{ \bd \Delta^n}$, we can reduce to proving that
$q_{ \Delta^n}$ is a weak homotopy equivalence. Note that assumption $(a)$ guarantees that $q_{ \Delta^n}$ is a flat categorical fibration. If $n > 1$, then we have a commutative diagram
$$ \xymatrix{ X \times_{S} \Lambda^n_1 \ar[d]^{q_{ \Lambda^n_1}} \ar[r] & X \times_{S} \Delta^n \ar[d]^{q_{\Delta^n}} \\
\Lambda^n_1 \ar[r] & \Delta^n }$$
where the upper horizontal map is a categorical equivalence (Corollary \bicatref{silman}) and
therefore a weak homotopy equivalence; the lower horizontal map is obviously a weak homotopy equivalence. Since $q_{ \Lambda^n_1}$ is a weak homotopy equivalence by the inductive hypothesis, we conclude that $q_{ \Delta^n}$ is a weak homotopy equivalence.

It remains to treat the cases where $n \leq 1$. If $n = 0$, the desired result follows from $(b)$.
Suppose finally that $n=1$. Let $X' = X \times_{S} \Delta^1$; we wish to prove that $X'$ is
weakly contractible. Let $X'_0$ and $X'_1$ denote the fibers of the map $q_{ \Delta^1}$, and let $Y = \Fun_{ \Delta^1}( \Delta^1, X')$. According to Proposition \bicatref{balder}, the natural map
$$ X'_0 \coprod_{ Y \times \{0\}} (Y \times \Delta^1) \coprod_{ Y \times \{1\} } X'_1 \rightarrow X'$$
is a categorical equivalence. Since $X'_0$ and $X'_1$ are weakly contractible, we are reduced
to showing that $Y$ is weakly contractible. Let $p: Y \rightarrow X'_0$ be the map given by evaluation
at $\{0\}$. Let $x' \in X'_0$, and let $x$ denote its image in $x$; the fiber $p^{-1} \{x'\}$ is isomorphic
to $(X')^{x'/} \times_{ \Delta^1} \{1\}$ and therefore categorically equivalent to
$X'_{x'/} \times_{ \Delta^1} \{1\} \simeq X_{x/} \times_{ S_{q(x)/} } \{f\}$, where
$f$ denotes the edge $\Delta^1 \rightarrow S$ under consideration. Assumption
$(c)$ guarantees that $p^{-1} \{x'\}$ is weakly contractible. Since $p$ is a Cartesian fibration,
Lemma \toposref{trull6prime} guarantees that $p$ is cofinal and therefore a weak homotopy equivalence. Since $X'_0$ is weakly contractible by $(b)$, we deduce that $Y$ is weakly contractible as desired.
\end{proof}

\begin{example}\label{stublos}
Let $q: X \rightarrow S$ be a flat inner fibration of simplicial sets, and let
$f: x \rightarrow y$ be an edge of $X$. Then the induced map
$X_{x/ \, /y} \rightarrow S_{ q(x)/ \, / q(y) }$ satisfies the hypotheses of Lemma \ref{clapsos} (see Proposition \bicatref{sliceflat}), and is therefore a categorical equivalence.
\end{example}

\begin{proposition}\label{scuzly}
Let $q: \calO^{\otimes} \rightarrow \Nerve(\FinSeg)$ be a unital $\infty$-operad. The following conditions
are equivalent:
\begin{itemize}
\item[$(1)$] The $\infty$-operad $\calO^{\otimes}$ is coherent.
\item[$(2)$] Let $Z \in \calO$, and suppose we are given a diagram
$\sigma:$
$$ \xymatrix{ X \ar[r]^{f} \ar[d] & Y \\
X' & }$$
in $\calO^{\otimes}_{/Z}$ where $f$ is semi-inert and the maps
$X \rightarrow Z$, $Y \rightarrow Z$, and $X' \rightarrow Z$ are active.
Let $\calB[\sigma,Z]$ denote the full subcategory of
$\calO^{\otimes}_{\sigma/ \, /Z}$ spanned by those diagrams
$$ \xymatrix{ X \ar[r]^{f} \ar[d] & Y \ar[d]^{g} \\
X' \ar[r] & Y' }$$
in $\calO^{\otimes}$ where the map $Y' \rightarrow Z$ is active, the map
$g$ is semi-inert, and the map $q(X') \coprod_{ q(X)}  q(Y) \rightarrow q(Y')$ is a surjective
map of pointed finite sets. Then $\calB[\sigma,Z]$ is weakly contractible.
\item[$(3)$] Condition $(2)$ holds in the special case where $f$ is required to
be $1$-semi-inert.
\end{itemize}
\end{proposition}

\begin{proof}
In the situation of $(2)$, let $\calA[\sigma]'_{/Z}$ be defined as in Remark \ref{funly}.
There is a canonical inclusion $\calB[\sigma,Z] \subseteq \calA[ \sigma]'_{/Z}$.
Our assumption that $\calO^{\otimes}$ is unital implies that this inclusion admits a right adjoint,
and is therefore a weak homotopy equivalence. The equivalence
$(1) \Leftrightarrow (2)$ now follows by combining Lemma \ref{stokums} with
Remarks \ref{funly}, \ref{saigr}, and \ref{saig2}. The same argument shows that
$(3)$ is equivalent to the condition that $\calO^{\otimes}$ is $1$-coherent.
The implication $(2) \Rightarrow (3)$ is obvious; we will complete the proof by showing that
$(3) \Rightarrow (2)$. 

Let $(\sigma,Z)$ be as in $(2)$; we wish to show that $\calB[ \sigma, Z]$ is weakly contractible.
The image $q(f)$ is a semi-inert morphism $\seg{n} \rightarrow \seg{n+m}$ in $\Nerve(\FinSeg)$, for
some $m \geq 0$. If $m = 0$, then $\calB[ \sigma, Z]$ has an initial object and there is nothing to prove.
We assume therefore that $m > 0$, so that (since $\calO^{\otimes}$ is unital) $f$ admits a factorization
$$ X \stackrel{ f'}{\rightarrow} X_0 \stackrel{ f''}{\rightarrow} X'$$
such that $q(f')$ is an inclusion $\seg{n} \hookrightarrow \seg{n+m-1}$ and
$q(f'')$ is an inclusion $\seg{n+m-1} \hookrightarrow \seg{n+m}$.
Let $\tau: \Delta^1 \coprod_{ \{0\}} \Delta^2 \rightarrow \calO^{\otimes}_{/Z}$ be the diagram
given by $Y \leftarrow X \rightarrow X_0 \rightarrow X'$, and let $\tau_0$ be the restriction of
$\tau$ to $\Delta^1 \coprod_{ \{0\} } \Delta^1$. Let $\calC$ denote the $\infty$-category
$$ \Fun( \Delta^1, \calO^{\otimes}_{ \tau_0 / \, / Z}) \times_{ \Fun( \{1\},
\calO^{\otimes}_{ \tau_0 / \, /Z}) } \Fun( \{1\}, \calO^{\otimes}_{ \tau/ \, / Z})$$
whose objects are commutative diagrams
$$ \xymatrix{ X \ar[r] \ar[d]^{f'} & Y \ar[d]^{g'} & \\
X_0 \ar[d]^{f''} \ar[r] & Y_0 \ar[d]^{g''} & \\
X' \ar[r] & Y' \ar[r] & Z }$$
in $\calO^{\otimes}$. Let $\calC_0$ denote the full subcategory of
$\calC$ spanned by those diagrams where the maps $g'$ and $g''$ are semi-inert, the maps $Y_0 \rightarrow Z$ and $Y' \rightarrow Z$ are active, and the maps $q(Y) \coprod_{ q(X)} q(X_0) \rightarrow q(Y_0)$ and $q(Y_0) \coprod_{ q(X_0)} q(X') \rightarrow q(Y')$ are surjective.
There are evident forgetful functors $\calB[\sigma],Z] \stackrel{\phi}{\leftarrow}
\calC_0 \stackrel{\psi}{\rightarrow} \calB[ \tau_0,Z]$. The map
$\phi$ admits a right adjoint and is therefore a weak homotopy equivalence.
The simplicial set $\calB[\tau_0,Z]$, and therefore
weakly contractible by the inductive hypothesis. To complete the proof, it will suffice to show that
$\psi$ is a weak homotopy equivalence. For this, we will show that $\psi$ satisfies the hypotheses
of Lemma \ref{clapsos}:

\begin{itemize}
\item[$(a)$] The map $\psi$ is a flat inner fibration. Fix a diagram
$$ \xymatrix{ & B \ar[dr] & \\
A \ar[rr]^{j} \ar[ur] & & C }$$
in $\calA[\tau_0]'_{/Z}$ and a morphism $\overline{j}: \overline{A} \rightarrow \overline{C}$ in
$\calC_0$ lifting $j$; we wish to show that the $\infty$-category
$$ \calZ = (\calC_0)_{ \overline{A} / \, / \overline{C} } \times_{ ( \calB[\tau_0, Z])_{A/ \, /C} } \{B\}$$
is weakly contractible. We have a commutative diagram
$$ \xymatrix{ X \ar[r] \ar[d] & Y \ar[d] & & & \\
X_0 \ar[r] \ar[d] & Y_A \ar[d] \ar[r] & Y_{B} \ar[r] & Y_{C} \ar[d] & \\
X' \ar[r] & Y'_{A} \ar[rr] & & Y'_{C} \ar[r] & Z }$$
in $\calO^{\otimes}$. Restricting our attention to the rectangle in the lower right,
we obtain a commutative diagram
$$ \xymatrix{ & Y_{B} \ar[dr] & \\
Y_{A} \ar[rr]^{ j_0} \ar[ur] & & Y_{C} }$$
in $\calO^{\otimes}$ and a morphism $\overline{j}_0: \overline{Y}_A \rightarrow \overline{Y}_C$
in $\calK_{\calO}$ lifting $j$. Let $\calZ_0 = ( \calK_{\calO})_{ \overline{Y}_A / \, / \overline{Y}_C }
\times_{ \calO^{\otimes}_{ Y_A/ \, /Y_{C} }} \{ Y_B \}$ be the $\infty$-category whose objects
are diagrams
$$ \xymatrix{ Y_A \ar[r] \ar[d]^{\alpha} & Y_B \ar[r] \ar[d] & Y_{C} \ar[d] \\
Y'_{A} \ar[r] & Y'_{B} \ar[r] & Y'_{C} }$$
in $\calO^{\otimes}$. Let $\calZ_{1}$ be the full subcategory of $\calZ_{0}$ spanned
by those diagrams for which the map $Y'_{B} \rightarrow Y'_{C}$ is active, and let
$\calZ_{2}$ be the full subcategory of $\calZ_1$ spanned by those maps
for which $q( Y_B) \coprod_{ q(Y_A)} q(Y'_{A}) \rightarrow q(Y'_{B})$ is surjective.
Since the map $q(Y_A) \coprod_{ q(X_0)} q(X') \rightarrow q(Y'_A)$ is surjective
and $f''$ is $1$-semi-inert, we deduce that $\alpha$ is $1$-semi-inert.
Condition $(3)$ guarantees that $\calO^{\otimes}$ is $1$-coherent, so the
$\infty$-category $\calZ_0$ is weakly contractible. The inclusion
$\calZ_1 \subseteq \calZ_0$ admits a left adjoint, and the inclusion
$\calZ_{2} \subseteq \calZ_{1}$ admits a right adjoint. It follows that both of these
inclusions are weak homotopy equivalences, so that $\calZ_{2}$ is weakly contractible.
There is an evident restriction map $\calZ \rightarrow \calZ_{2}$, which is easily shown
to be a trivial Kan fibration. It follows that $\calZ$ is weakly contractible as desired.

\item[$(b)$] The fibers of $\psi$ are weakly contractible. To prove this, we observe
an object $b \in \calB[ \tau_0, Z]$ determines a commutative diagram
$$ \xymatrix{ X \ar[r] \ar[d] & Y \ar[d] \\
X_0 \ar[r] \ar[d] & Y_0 \\
X' & }$$
in $\calO^{\otimes}_{/Z}$. If we let $\sigma'$ denote the lower part of this diagram,
then we have a trivial Kan fibration $\psi^{-1}(b) \rightarrow \calB[ \sigma',Z]$,
which is weakly contractible by virtue of $(3)$.

\item[$(c)$] For every object $c \in \calC_0$ and every morphism
$\beta: \psi(c) \rightarrow b$ in $\calB[ \tau_0, Z]$, the $\infty$-category
$$\calY = ( \calC_0)_{c/} \times_{ \calB[ \tau_0, Z]_{ \psi(c)/ } } \{ \beta \}$$ is weakly contractible.
The pair $(c, \beta)$ determines a diagram
$$ \xymatrix{ X \ar[r] \ar[d] & Y \ar[d] &  \\
X_0 \ar[r] \ar[d] & Y_0 \ar[r] \ar[d]^{g''} & Y_1 \\
X' \ar[r] & Y' & }$$
in $\calO^{\otimes}_{/Z}$. Let $\sigma''$ denote the lower right corner of this diagram.
Then we have a trivial Kan fibration $\calY \rightarrow \calB[ \sigma'', Z]$.
Since the map $q(Y_0) \coprod_{ q(X_0)} q(X') \rightarrow q(Y')$ is surjective,
we deduce that $g''$ is $1$-semi-inert, so that $\calB[ \sigma'', Z]$ is weakly contractible by virtue of $(3)$.
\end{itemize}
\end{proof}

\begin{proof}[Proof of Theorem \ref{uggus}]
In view of Remark \ref{saym}, it will suffice to show that conditions
$(1)$ and $(3)$ of Theorem \ref{uggus} are equivalent.
Fix a pair of active morphisms $f: X \rightarrow Y$ and $g: Y \rightarrow Z$ in
$\calO^{\otimes}$, where $Z \in \calO$. Let $\sigma: \Delta^3 \rightarrow \calO^{\otimes}$ be
as in the formulation of condition $(3)$, and consider the diagram
$$ \xymatrix{ \Ext( \sigma, \{0,1\} ) \ar[r] \ar[d] & \Ext( \sigma | \Delta^{ \{0,1,3 \} }, \{0,1\} ) \ar[d] \\
\Ext( \sigma | \Delta^{ \{0,2,3\} }, \{0\} ) \ar[r] & \Ext( \sigma| \Delta^{ \{0,3\} }, \{0\} ). }$$
Each of the maps in this diagram is a Kan fibration between Kan complexes.
Consequently, condition $(3)$ is satisfied if and only if, for every vertex
$v$ of $\Ext( \sigma | \Delta^{ \{0,3\} }, \{0\})$, the induced diagram of fibers
$$  \Ext( \sigma | \Delta^{ \{0,1,3 \} }, \{0,1\} )_{v} \leftarrow \Ext( \sigma, \{0,1\})_{v}
\rightarrow \Ext( \sigma| \Delta^{ \{0,2,3\}}, \{0\})_{v}$$
has a contractible homotopy pushout. 
Without loss of generality, we may assume that
$v$ determines a diagram
$$ \xymatrix{ X \ar[r] \ar[d] & Y \ar[r] & Z \ar[d]^{\id_{Z} } \\
X' \ar[rr] & & Z }$$ 
in $\calO^{\otimes}$, where the left vertical map is $1$-semi-inert. 
Let $\tau$ denote the induced diagram $X' \leftarrow X \rightarrow Y$ in
$\calO^{\otimes}_{/Z}$, and let $\calB[\tau,Z]$ be the $\infty$-category
defined in Proposition \ref{scuzly}. Let $\calB[\tau,Z]_1$ denote the full subcategory
of $\calB[\tau,Z]$ spanned by those objects which correspond to diagrams
$$ \xymatrix{ X \ar[r] \ar[d] & Y \ar[d] \\
X' \ar[r] & Y' }$$
where the right vertical map is an equivalence. Then there is a unique
map $p: \calB[ \tau, Z] \rightarrow \Delta^1$ such that $p^{-1} \{1\} \simeq \calB[ \tau,Z]_1$.
Let $\calB[ \tau, Z]_0 = p^{-1} \{0\}$, and let $\calZ = \Fun_{ \Delta^1}( \Delta^1, \calB[\tau,Z])$ be the
$\infty$-category of sections of $p$. We have a commutative diagram
$$ \xymatrix{ \Ext( \sigma| \Delta^{ \{0,1,3 \}}, \{0,1\})_{v} \ar[r] \ar[d] & \Ext( \sigma, \{0,1\})_{v} \ar[l] \ar[r] \ar[d] &
\Ext( \sigma| \Delta^{ \{0,2,3\}}, \{0\})_{v} \ar[d] \\
\calB[\tau,Z]_0 & \calZ \ar[l] \ar[r] & \calB[\tau,Z]_1. }$$
in which the vertical maps are trivial Kan fibrations. Consequently, condition $(3)$ is satisfied
if and only if each homotopy pushout
$$ \calB[ \tau,Z]_0 \coprod_{ \calZ \times \{0\} } ( \calZ \times \Delta^1) \coprod_{ \calZ \times \{1\}}
\calB[\tau,Z]_{1}$$
is weakly contractible. According to Proposition \bicatref{balder}, this homotopy pushout is
categorically equivalent to $\calB[\tau,Z]$, so the equivalence of
$(1)$ and $(3)$ follows from Proposition \ref{scuzly}.
\end{proof}

\subsection{Coproducts of $\infty$-Operads}\label{coprodinf}

Let $\Cat_{\infty}^{\lax}$ denote the $\infty$-category of $\infty$-operads.
Then $\Cat_{\infty}^{\lax}$ can be realized as the underlying $\infty$-category of a combinatorial
simplicial model category $\PreOp$ (see \S \symmetricref{comm1.8}), and therefore admits
small limits and colimits (Corollary \toposref{limitsinmodel}). The limit of a diagram $\sigma$ in
$\Cat_{\infty}^{\lax}$ can usually be described fairly explicitly: namely, choose an injectively fibrant diagram $\overline{\sigma}$ in $\PreOp$ representing $\sigma$, and then take the limit of 
$\overline{\sigma}$ in the ordinary category of $\infty$-preoperads. The case of colimits is more difficult: we can apply the same procedure to construct an $\infty$-preoperad which represents $\varinjlim(\sigma)$, but this representative will generally not be fibrant and the process of ``fibrant replacement'' is fairly inexplicit. Our goal in this section is to give a more direct construction of colimits in a special case: namely, the case of coproducts.

\begin{definition}\label{tsongd}
Given an object $\seg{n} \in \FinSeg$ and a subset $S \subseteq \seg{n}$ which contains
the base point, there is a unique integer $k$ and bijection $\seg{k} \simeq S$ whose
restriction to $\nostar{k}$ is order-preserving; we will denote the corresponding object
of $\FinSeg$ by $[S]$. We define a category $\Subd$ as follows:
\begin{itemize}
\item[$(1)$] The objects of $\Subd$ are triples $( \seg{n}, S, T)$ where
$\seg{n} \in \FinSeg$, $S$ and $T$ are subsets of $\seg{n}$ such that
$S \cup T = \seg{n}$ and $S \cap T = \{ \ast \}$.
\item[$(2)$] A morphism from $( \seg{n}, S, T)$ to $( \seg{n'}, S', T')$ in
$\Subd$ is a morphism $f: \seg{n} \rightarrow \seg{n'}$ in $\FinSeg$ 
such that $f(S) \subseteq S'$ and $f(T) \subseteq T'$.
\end{itemize}

There is an evident triple of functors $\pi, \pi_{-}, \pi_{+}: \Subd \rightarrow \FinSeg$,
given by the formulas
$$\pi_{-}( \seg{n}, S, T) = [S] \quad \quad  \pi( \seg{n}, S, T) = \seg{n} \quad \quad \pi_{+}( \seg{n}, S, T) = [T]. $$

For any pair of $\infty$-operads $\calC^{\otimes}$ and $\calD^{\otimes}$, we define
a new simplicial set $\calC^{\otimes} \boxplus \calD^{\otimes}$ so that we have a pullback diagram
$$ \xymatrix{ \calC^{\otimes} \boxplus \calD^{\otimes} \ar[r] \ar[d] & \calC^{\otimes} \times \calD^{\otimes} \ar[d] \\
\Nerve( \Subd) \ar[r]^-{ \pi_{-} \times \pi_{+}} & \Nerve(\FinSeg) \times \Nerve(\FinSeg). }$$
We regard $\calC^{\otimes} \boxplus \calD^{\otimes}$ as equipped with a map to
$\Nerve(\FinSeg)$, given by the composition
$$ \calC^{\otimes} \boxplus \calD^{\otimes} \rightarrow \Nerve(\Subd) \stackrel{\pi}{\rightarrow} \Nerve(\FinSeg).$$
\end{definition}

\begin{remark}
The product functor $( \pi_{-} \times \pi_{+}): \Subd \rightarrow \FinSeg \times \FinSeg$ is an
equivalence of categories. Consequently, $\calC^{\otimes} \boxplus \calD^{\otimes}$ is equivalent
(as an $\infty$-category) to the product $\calC^{\otimes} \times \calD^{\otimes}$. 
However, it is slightly better behaved in the following sense: the composite map
$$ \calC^{\otimes} \boxplus \calD^{\otimes} \rightarrow \Nerve(\Subd) \rightarrow \Nerve(\FinSeg)$$
is a categorical fibration, since it is the composition of a pullback of the categorical fibration
$\calC^{\otimes} \times \calD^{\otimes} \rightarrow \Nerve(\FinSeg) \times \Nerve(\FinSeg)$
with the categorical fibration $\Nerve(\Subd) \rightarrow \Nerve(\FinSeg)$.
\end{remark}

\begin{lemma}\label{scak}
Let $\calC^{\otimes}$ and $\calD^{\otimes}$ be $\infty$-operads. Then
$\calC^{\otimes} \boxplus \calD^{\otimes}$ is an $\infty$-operad.
\end{lemma}

\begin{proof}
We can identify an object of $( \calC^{\otimes} \boxplus \calD^{\otimes})_{ \seg{n}}$ with a
quintuple $X = (\seg{n}, S, T, C, D)$, where $(\seg{n}, S, T) \in \Subd$,
$C \in \calC^{\otimes}_{[S]}$, and $D \in \calD^{\otimes}_{[T]}$. Suppose we are given
an inert map $\alpha: \seg{n} \rightarrow \seg{n'}$ in $\FinSeg$. Let
$S' = \alpha(S)$ and $T' = \alpha(T)$. Then $\alpha$ induces inert morphisms
$\alpha_{-}: [S] \rightarrow [S']$ and $\alpha_{+}: [T] \rightarrow [T']$. Choose an inert
morphism $f_{-}: C \rightarrow C'$ in $\calC^{\otimes}$ lifting $\alpha_{-}$ and an inert morphism $f_{+}: D \rightarrow D'$ in $\calD^{\otimes}$ lifting $\alpha_{+}$. These maps together determine
a morphism $f: ( \seg{n}, S,T, C, D) \rightarrow ( \seg{n'}, S', T', C', D')$. 
Since $(f_{-}, f_{+})$ is coCartesian with respect to the projection
$\calC^{\otimes} \times \calD^{\otimes} \rightarrow \Nerve(\FinSeg) \times \Nerve(\FinSeg)$, the map
$f$ is $p$-coCartesian, where $p$ denotes the map $\calC^{\otimes} \boxtimes \calD^{\otimes}
\rightarrow \Nerve(\Subd)$. Let $\pi: \Nerve(\Subd) \rightarrow \Nerve(\FinSeg)$ be as
in Definition \ref{tsongd}. It is easy to see that $p(f)$ is $\pi$-coCartesian, so that
$f$ is $(\pi \circ p)$-coCartesian by virtue of Proposition \toposref{stuch}.

Choose $(\pi \circ p)$-coCartesian morphisms $X \rightarrow X_i$ covering the inert morphisms
$\Colp{i}: \seg{n} \rightarrow \seg{1}$ for $1 \leq i \leq n$. We claim that these maps exhibit
$X$ as a $(\pi \circ p)$-product of the objects $X_i$. Using our assumption that
$\calC^{\otimes}$ and $\calD^{\otimes}$ are $\infty$-operads, we deduce that
these maps exhibit $X$ as a $p$-product of the objects $\{ X_i \}$. It therefore suffices to show that
they exhibit $p(X)$ as a $\pi$-product of the objects $p(X_i)$ in the $\infty$-category
$\Nerve(\Subd)$, which follows immediately from the definitions.

It remains only to show that for each $n \geq 0$, the canonical functor
$\phi: ( \calC^{\otimes} \boxplus \calD^{\otimes})_{ \seg{n}} \rightarrow
( \calC^{\otimes} \boxplus \calD^{\otimes})_{\seg{1}}^{n}$ is essentially surjective.
We can identify the latter with $( \calC \coprod \calD)^{n}$. Given a collection of objects
$X_1, \ldots, X_n$ of $\calC \coprod \calD$, we let $S = \{ \ast \} \cup \{ i: X_i \in \calC \}$
and $T = \{ \ast \} \cup \{ i: X_i \in \calD \}$. Let $C = \bigoplus_{ X_i \in \calC} X_i
\in \calC^{\otimes}_{[S]}$ and let $D = \bigoplus_{X_i \in \calD} X_i \in \calD^{\otimes}_{[T]}$. Then 
$X = ( \seg{n}, S, T, C, D)$ is an object of $(\calC^{\otimes} \boxplus \calD^{\otimes})_{\seg{n}}$ such
that $\phi(X) \simeq (X_1, \ldots, X_n)$.
\end{proof}

\begin{remark}
The operation $\boxplus$ of Definition \ref{tsongd} is commutative and associative up to coherent isomorphism, and determines a symmetric monoidal structure on the category
$(\sSet)_{/ \Nerve(\FinSeg)}$ of simplicial sets $X$ endowed with a map $X \rightarrow \Nerve(\FinSeg)$. This restricts to a symmetric monoidal structure on the (ordinary) category of $\infty$-operads and
maps of $\infty$-operads.
\end{remark}

We next show that the $\infty$-operad $\calC^{\otimes} \boxplus \calD^{\otimes}$
can be identified with a coproduct of $\calC^{\otimes}$ with $\calD^{\otimes}$ in
$\Cat_{\infty}^{\lax}$. 

\begin{theorem}\label{cross}
Let $\calC^{\otimes}$ and $\calD^{\otimes}$ be $\infty$-operads. We
let $( \calC^{\otimes} \boxplus \calD^{\otimes})_{-}$ denote the full subcategory
of $\calC^{\otimes} \boxplus \calD^{\otimes}$ spanned by those objects whose
image in $\Subd$ has the form $( \seg{n}, \seg{n}, \{ \ast \})$, and
let $( \calC^{\otimes} \boxplus \calD^{\otimes})_{+}$ denote the full subcategory spanned by those objects whose image in $\Subd$ has the form $( \seg{n}, \{ \ast \}, \seg{n} )$. Then:
\begin{itemize}
\item[$(1)$] The projection maps $(\calC^{\otimes} \boxplus \calD^{\otimes})_{-} \rightarrow \calC^{\otimes}$ and $( \calC^{\otimes} \boxplus \calD^{\otimes})_{+} \rightarrow
\calD^{\otimes}$ are trivial Kan fibrations.
\item[$(2)$] The map $\calC^{\otimes} \boxplus \calD^{\otimes} \rightarrow \Nerve(\FinSeg)$ exhibits both $( \calC^{\otimes} \boxplus \calD^{\otimes})_{-}$ and $( \calC^{\otimes} \boxplus \calD^{\otimes})_{+}$ as $\infty$-operads.
\item[$(3)$] For any $\infty$-operad $\calE^{\otimes}$, the inclusions
$$ ( \calC^{\otimes} \boxplus \calD^{\otimes})_{-} \stackrel{i}{\hookrightarrow}
\calC^{\otimes} \boxplus \calD^{\otimes} \stackrel{j}{\hookleftarrow} ( \calC^{\otimes} \boxplus \calD^{\otimes})_{+}$$
induce an equivalence of $\infty$-categories
$$ \Fun^{\lax}( \calC^{\otimes} \boxplus \calD^{\otimes}, \calE^{\otimes})
\rightarrow \Fun^{\lax}( (\calC^{\otimes} \boxplus \calD^{\otimes})_{-}, \calE^{\otimes}) 
\times \Fun^{\lax}( ( \calC^{\otimes} \boxplus \calD^{\otimes})_{+}, \calE^{\otimes}).$$
In particular, $i$ and $j$ exhibit $\calC^{\otimes} \boxplus \calD^{\otimes}$ as a 
coproduct of $( \calC^{\otimes} \boxplus \calD^{\otimes})_{-} \simeq \calC^{\otimes}$ and
$( \calC^{\otimes} \boxplus \calD^{\otimes})_{+} \simeq \calD^{\otimes}$ in the $\infty$-category
$\Cat_{\infty}^{\lax}$. 
\end{itemize}
\end{theorem}

\begin{proof}
Assertion $(1)$ follows from the evident isomorphisms
$$ (\calC^{\otimes} \boxplus \calD^{\otimes})_{-} \simeq \calC^{\otimes} \times \calD^{\otimes}_{\seg{0}} \quad \quad (\calC^{\otimes} \boxplus \calD^{\otimes})_{+} \simeq \calC^{\otimes}_{ \seg{0}} \times \calD^{\otimes},$$
together with the observation that $\calC^{\otimes}_{\seg{0}}$ and $\calD^{\otimes}_{\seg{0}}$ are contractible Kan complexes. Assertion $(2)$ follows immediately from $(1)$. To prove
$(3)$, let $\calX = ( \calC^{\otimes} \boxplus \calD^{\otimes})_{-} \cap (\calC^{\otimes} \boxplus \calD^{\otimes})_{+} \simeq \calC^{\otimes}_{\seg{0}} \times \calD^{\otimes}_{\seg{0}}$ and
let $\calY =  ( \calC^{\otimes} \boxplus \calD^{\otimes})_{-} \cup (\calC^{\otimes} \boxplus \calD^{\otimes})_{+}$. Let $\calA$ denote the full subcategory of $\Fun_{ \Nerve(\FinSeg)}( \calY, \calE^{\otimes})$ spanned
by those functors whose restriction to both $( \calC^{\otimes} \boxplus \calD^{\otimes})_{-}$ and
$( \calC^{\otimes} \boxplus \calD^{\otimes})_{+}$ are $\infty$-operad maps. We have homotopy
pullback diagram
$$ \xymatrix{ \calA \ar[r] \ar[d] &  \Fun^{\lax}( (\calC^{\otimes} \boxplus \calD^{\otimes})_{-}, \calE^{\otimes}) \ar[d] \\
 \Fun^{\lax}( ( \calC^{\otimes} \boxplus \calD^{\otimes})_{+}, \calE^{\otimes}) \ar[r] & \Fun_{ \Nerve(\FinSeg)}( \calX, \calE^{\otimes}). }$$
Since $\calE^{\otimes}_{\seg{0}}$ is a contractible Kan complex, the
simplicial set $\Fun_{ \Nerve(\FinSeg)}( \calX, \calE^{\otimes})$ is also a contractible
Kan complex, so the map 
$$\calA \rightarrow \Fun^{\lax}( (\calC^{\otimes} \boxplus \calD^{\otimes})_{-}, \calE^{\otimes}) 
\times \Fun^{\lax}( ( \calC^{\otimes} \boxplus \calD^{\otimes})_{+}, \calE^{\otimes})$$
is a categorical equivalence. We will complete the proof by showing that the map $\Fun^{\lax}( \calC^{\otimes} \boxplus \calD^{\otimes}, \calE^{\otimes})
\rightarrow \calA$ is a trivial Kan fibration.

Let $q: \calE^{\otimes} \rightarrow \Nerve(\FinSeg)$ denote the projection map. 
In view of Proposition \toposref{lklk}, it will suffice to show the following:
\begin{itemize}
\item[$(a)$] An arbitrary map $A \in \Fun_{ \Nerve(\FinSeg)}( \calC^{\otimes} \boxplus \calD^{\otimes}, \calE^{\otimes})$ is an $\infty$-operad map if and only if it satisfies the following conditions:
\begin{itemize}
\item[$(i)$] The restriction $A_0 = A | \calY$ belongs to $\calA$.
\item[$(ii)$] The map $A$ is a $q$-right Kan extension of $A_0$.
\end{itemize}
\item[$(b)$] For every object $A_0 \in \calA$, there exists an extension
$A \in \Fun_{ \Nerve(\FinSeg)}( \calC^{\otimes} \boxplus \calD^{\otimes}, \calE^{\otimes})$
of $A_0$ which satisfies the equivalent conditions of $(a)$.
\end{itemize}

To prove $(a)$, consider an object $A \in \Fun_{ \Nerve(\FinSeg)}( \calC^{\otimes} \boxplus \calD^{\otimes}, \calE^{\otimes})$ and an object $X = ( \seg{n}, S, T, C, D) \in \calC^{\otimes} \boxplus \calD^{\otimes}$. Choose morphisms $\alpha: C \rightarrow C_0$ and $\beta: D \rightarrow D_0$, where
$C_0 \in \calC^{\otimes}_{\seg{0}}$ and $D_0 \in \calD^{\otimes}_{\seg{0}}$. Set
$$X_{-} = ( [S], [S], \{ \ast \}, C, D_0) \quad \quad
X_0 = ( \seg{0}, \{ \ast \}, \{ \ast \}, C_0, D_0 ) \quad \quad X_{+} = ( [T], \{ \ast \}, [T], C_0, D).$$
Then $\alpha$ and $\beta$ determine a commutative diagram
$$ \xymatrix{ X \ar[r] \ar[d] & X_{-} \ar[d] \\
X_{+} \ar[r] & X_0. }$$
Using Theorem \toposref{hollowtt}, we deduce that this diagram determines a map
$$ \phi: \Lambda^2_{2} \rightarrow \calY \times_{ \calC^{\otimes} \boxplus \calD^{\otimes}}
( \calC^{\otimes} \boxplus \calD^{\otimes})_{X/}$$
such that $\phi^{op}$ is cofinal. It follows that $A$ is a $q$-right Kan extension of
$A_0$ at $X$ if and only if the diagram
$$ \xymatrix{ A(X) \ar[r] \ar[d] & A(X_{-}) \ar[d] \\
A(X_{+}) \ar[r] & A(X_0) }$$ is
a $q$-pullback diagram. Since $q: \calE^{\otimes} \rightarrow \Nerve(\FinSeg)$ is an
$\infty$-operad, this is equivalent to the requirement that the maps
$A(X_{-}) \leftarrow A(X) \rightarrow A(X_{+})$ are inert. In other words, we obtain
the following version of $(a)$:

\begin{itemize}
\item[$(a')$] A map $A \in \Fun_{ \Nerve(\FinSeg)}( \calC^{\otimes} \boxplus \calD^{\otimes}, \calE^{\otimes})$ is a $q$-right Kan extension of $A_0 = A | \calY$ if and only if,
for every object $X \in \calC^{\otimes} \boxplus \calD^{\otimes}$ as above, the maps
$A(X_{-}) \leftarrow A(X) \rightarrow A(X_{+})$ are inert.
\end{itemize}

We now prove $(a)$ Suppose first that $A$ is an $\infty$-operad map; we wish to prove
that $A$ satisfies conditions $(i)$ and $(ii)$. Condition $(i)$ follows immediately from the
description of the inert morphisms in $\calC^{\otimes} \boxplus \calD^{\otimes}$ provided
by the proof of Lemma \ref{scak}, and condition $(ii)$ follows from $(a')$. Conversely,
suppose that $(i)$ and $(ii)$ are satisfied. We wish to prove that $A$ preserves inert morphisms.
In view of Remark \symmetricref{casper}, it will suffice to show that $A(X) \rightarrow A(Y)$ is inert whenever $X \rightarrow Y$ is an inert morphism such
that $Y \in ( \calC^{\otimes} \boxplus \calD^{\otimes})_{\seg{1}}$. It follows
that $Y \in \calY$; we may therefore assume without loss of generality that
$Y \in ( \calC^{\otimes} \boxplus \calD^{\otimes})_{-}$. The map $X \rightarrow Y$
then factors as a composition of inert morphisms $X \rightarrow X_{-} \rightarrow Y$,
where $X_{-}$ is defined as above. Then $A(X) \rightarrow A(X_{-})$ is inert by
virtue of $(ii)$ and $(a')$, while $A(X_{-}) \rightarrow A(Y)$ is inert by virtue of assumption
$(i)$.

To prove $(b)$, it will suffice (by virtue of Lemma \toposref{kan2}) to show that
for each $X \in \calC^{\otimes} \boxtimes \calD^{\otimes}$, the induced diagram
$$ \calY_{X/} \rightarrow \calY \stackrel{A_0}{\rightarrow} \calE^{\otimes}$$
admits a $q$-limit. Since $\phi^{op}$ is cofinal, it suffices to show that there
exists a $q$-limit of the diagram
$$ A_0( X_{-}) \rightarrow A_0(X_0) \leftarrow A_0( X_{+}).$$
The existence of this $q$-limit follows immediately from the assumption that
$\calE^{\otimes}$ is an $\infty$-operad.
\end{proof}

\begin{remark}\label{grill}
We can informally summarize Theorem \ref{cross} as follows: for every
triple of $\infty$-operads $\calO_{-}^{\otimes}$, $\calO_{+}^{\otimes}$, and $\calC^{\otimes}$, we have
a canonical equivalence of $\infty$-categories
$$ \Alg_{\calO}(\calC) \rightarrow \Alg_{ \calO_{-}}(\calC) \times \Alg_{\calO_{+}}(\calC),$$
where $\calO^{\otimes} = \calO_{-}^{\otimes} \boxplus \calO^{\otimes}_{+}.$
\end{remark}

As an application of Theorem \ref{cross}, we will explain how to view understand the formation of tensor products of algebras as an instance of the theory of operadic left Kan extensions. First, we need to introduce a bit of notation.

\begin{notation}\label{skaz}
We let $\Subd'$ denote the categorical mapping cylinder of the forgetful functor
$\pi: \Subd \rightarrow \FinSeg$ of Definition \ref{tsongd}. More precisely, the category
$\Subd$ can be described as follows:
\begin{itemize}
\item[$(i)$] An object of $\Subd'$ is either an object of $\Subd$ or an object of $\FinSeg$.
\item[$(ii)$] Morphisms in $\Subd'$ are given by the formulas
$$\Hom_{ \Subd'}( \seg{m}, ( \seg{n}, S, T) ) = \emptyset.$$
$$ \Hom_{ \Subd'}( \seg{m}, \seg{n} ) = \Hom_{ \FinSeg}( \seg{m}, \seg{n}) \quad \quad
\Hom_{ \Subd'}( (\seg{m}, S, T), \seg{n}) = \Hom_{ \FinSeg}( \seg{m}, \seg{n})$$
$$\Hom_{ \Subd'}( (\seg{m}, S,T), ( \seg{n}, S', T')) = \Hom_{ \Subd}( (\seg{m}, S, T)
( \seg{n}, S', T')).$$
\end{itemize}
The functors $\pi_{-}$, $\pi_{+}$, and $\pi$ of Definition \ref{tsongd} extend naturally
to retractions of $\Subd'$ onto the full subcategory $\FinSeg \subseteq \Subd'$; we will denote these retractions by $\pi'_{-}$, $\pi'_{+}$, and $\pi'$.

Let $\calO^{\otimes}$ be an $\infty$-operad.
We define a simplicial set
$\calT_{\calO}$ equipped with a map $\calT_{\calO} \rightarrow \Nerve(\Subd')$ so
that the following universal property is satisfied: for every simplicial set
$K$ equipped with a map $K \rightarrow \Nerve(\Subd')$, the set
$\Hom_{ \Nerve(\Subd')}( K, \calT_{\calO})$ can be identified with the set of
pairs of maps $e_{-},e_{+}: K \rightarrow \calO^{\otimes}$ satisfying the following conditions:
\begin{itemize}
\item[$(1)$] The diagram
$$ \xymatrix{ K \ar[r]^{e_{-}} \ar[d] & \calO^{\otimes} \ar[d] & K \ar[l]_{e_{+}} \ar[d] \\
\Nerve(\Subd') \ar[r]^{\pi'_{-}} & \Nerve(\FinSeg) & \Nerve(\Subd') \ar[l]_{ \pi'_{+}} }$$ 
is commutative.
\item[$(2)$] The maps $e_{-}$ and $e_{+}$ coincide on $K \times_{ \Nerve( \Subd')} \Nerve(\FinSeg)$.
\end{itemize}
We regard $\calT_{\calO}$ as an object of $( \sSet)_{/ \Delta^1 \times \Nerve(\FinSeg)}$ via
the map
$$ \calT_{\calO} \rightarrow \Nerve(\Subd') \stackrel{\psi \times \pi'}{\rightarrow} \Delta^1 \times \Nerve(\FinSeg),$$
where $\psi: \Nerve(\Subd') \rightarrow \Delta^1$ is the map given by $0$ on
the subcategory $\Subd \subseteq \Subd'$ and $1$ on $\FinSeg \subseteq \Subd'$.
\end{notation}

We observe that the fiber product $\calT_{\calO} \times_{ \Delta^1} \{0\}$ is isomorphic with
$\calO^{\otimes} \boxplus \calO^{\otimes}$, while $\calT_{\calO} \times_{ \Delta^1} \{1\}$ is isomorphic to $\calO^{\otimes}$.

\begin{lemma}\label{kwop}
Let $\calO^{\otimes}$ be an $\infty$-operad. Then the map
$q: \calT_{\calO} \rightarrow \Delta^1 \times \Nerve(\FinSeg)$ of Notation \ref{skaz} exhibits $\calT_{\calO}$ as a $\Delta^1$-family of $\infty$-operads.
\end{lemma}

\begin{proof}
It is easy to see that $\calT_{\calO}$ is an $\infty$-category, and that the fibers $\calT_{\calO} \times_{ \Delta^1} \{i\}$ are $\infty$-operads (Lemma \ref{scak}). To complete the proof, we must show the following:
\begin{itemize}
\item[$(a)$] Every inert morphism in $\calT_{\calO} \times_{ \Delta^1} \{i\}$ is $q$-coCartesian.
\item[$(b)$] Given an object $X \in \calT_{\calO} \times_{ \Delta^1 } \{i\}$ lying over
$\seg{n} \in \Nerve(\FinSeg)$, any collection of inert morphisms
$X \rightarrow X_j$ covering $\Colp{j}: \seg{n} \rightarrow \seg{1}$ for $1 \leq j \leq n$
exhibits $X$ as a $q$-product of the objects $\{ X_j \}_{1 \leq j \leq n}$.
\end{itemize}
Assertion $(a)$ is obvious when $i=1$, and follows easily from the characterization
of inert morphisms given in the proof of Lemma \ref{scak} when $i=0$. Assertion
$(b)$ is obvious when $i=0$. Let us prove that $(b)$ holds when $i=1$. 
Fix an object $Y \in \calT_{\calO} \times_{ \Delta^1} \{0\}$ lying over
$\seg{m} \in \Nerve(\FinSeg)$, and for every map $\alpha: \seg{m} \rightarrow \seg{n}$ in $\FinSeg$
let $\bHom_{ \calT_{\calO}}^{\alpha}( Y, X)$ denote the summand of
$\bHom_{ \calT_{\calO}}(Y, X)$ lying over $\alpha$. We wish to prove that the map
$$\phi: \bHom_{ \calT_{\calO}}^{\alpha}(Y,X) \rightarrow \prod_{1 \leq j \leq n} \bHom_{\calT_{\calO}}
^{ \Colp{j} \alpha}(Y, X_j)$$
is a homotopy equivalence. We can identify $Y$ with a quintuple
$( \seg{m}, S, T, Y_{-}, Y_{+}) \in \calO^{\otimes} \boxplus \calO^{\otimes}$.
The map $\alpha$ determines a pair of maps $\alpha_{-}: [S] \rightarrow \seg{n}$ and
$\alpha_{+}: [T] \rightarrow \seg{n}$. It now suffices to observe that
$\phi$ is equivalent to a product of the maps
$$ \bHom_{ \calO^{\otimes} }^{\alpha_{-}}(Y_{-}, q_{-}(X)) \rightarrow
\prod_{1 \leq j \leq n} \bHom_{\calO^{\otimes}}^{ \Colp{j} \alpha_{-}}(Y_{-}, X_j) $$
$$ \bHom_{ \calO^{\otimes} }^{\alpha_{+}}(Y_{+}, q_{+}(X)) \rightarrow
\prod_{1 \leq j \leq n} \bHom_{\calO^{\otimes}}^{ \Colp{j} \alpha_{+}}(Y_{+}, X_j),$$
which are homotopy equivalences by virtue of our assumption that $\calO^{\otimes}$ is an $\infty$-operad.
\end{proof}

\begin{construction}\label{urmouth}
Let $\calO^{\otimes}$ be an $\infty$-operad. There is a canonical map
$H: \calO^{\otimes} \times \Delta^1 \rightarrow \calT_{ \calO}$,
where $H_1 = H | \calO^{\otimes} \times \{1\}$ is the canonical isomorphism
$\calO^{\otimes} \simeq \calT_{\calO} \times_{ \Delta^1} \{1\}$
and $H_0 = H | \calO^{\otimes} \times \{0\}$ is given on objects by the formula
$(X \in \calO^{\otimes}_{\seg{n}}) \mapsto ( \seg{n}, \seg{n}, \seg{n}, X, X)$.

Suppose we are given an $\infty$-operad family map $\psi: \calT_{\calO}
\rightarrow {\calD}^{\otimes}$, where $\calD^{\otimes}$ is an $\infty$-operad. Restricting $\psi$ to $\calT_{\calO} \times_{ \Delta^1} \{1\}$, we
obtain an $\infty$-operad map $\psi_{\pm}: \calO^{\otimes} \rightarrow {\calD}^{\otimes}$.
Similarly, the restriction of $\psi$ to $\calO^{\otimes} \boxplus \calO^{\otimes}$ yields a pair
of $\infty$-operad maps $\psi_{-}, \psi_{+}: \calO^{\otimes} \rightarrow {\calD'}^{\otimes}$.
Let $p: \calC^{\otimes} \rightarrow {\calD}^{\otimes}$ be a coCartesian fibration of $\infty$-operads.
We let $\Alg^{-}_{\calO}(\calC)$, $\Alg^{+}_{\calO}(\calC)$, and $\Alg^{\pm}_{\calO}(\calC)$
be the $\infty$-categories of $\calO$-algebras in $\calC$, defined using the $\infty$-operad
morphisms $\psi_{-}$, $\psi_+$, and $\psi_{\pm}$, respectively.
For every algebra object $A \in \Fun_{{\calD}^{\otimes}(\calO^{\otimes} \boxplus \calO^{\otimes}},\calC)$,
we can choose a $p$-coCartesian natural transformation
$\alpha: A \circ H_0 \rightarrow A'$ in $\Fun( \calO^{\otimes}, \calC^{\otimes})$ lifting
the natural transformation $\psi \circ H$. This construction determines a functor
$$ \Fun^{\lax}_{ {\calD}^{\otimes}}( \calO^{\otimes} \boxplus \calO^{\otimes}, \calC)
\rightarrow \Alg_{\calO}^{\pm}(\calC).$$
Theorem \ref{cross} allows us to identify the $\infty$-category on the left hand side
with $\Alg_{\calO}^{-}(\calC) \times \Alg_{\calO}^{+}(\calC)$. Under this identification,
we obtain a functor
$$\Alg_{\calO}^{-}(\calC) \times \Alg_{\calO}^{+}(\calC) \rightarrow \Alg^{\pm}_{\calO}(\calC);$$
we will refer to this functor as the {\it tensor product} and denote it by $\otimes$.

In the special case where ${\calD}^{\otimes} = \Nerve(\FinSeg)$, we recover the usual tensor product operation on algebra objects of a symmetric monoidal $\infty$-category (see Proposition \symmetricref{slape}).
\end{construction}

\begin{proposition}\label{alug}
Suppose we are given an $\infty$-operad $\calO^{\otimes}$, an $\infty$-operad family map
$\psi: \calT_{\calO} \rightarrow {\calD}^{\otimes}$, and a coCartesian
fibration of $\infty$-operads $q: \calC^{\otimes} \rightarrow {\calD}^{\otimes}$.
Then the tensor product functor
$\otimes: \Alg_{\calO}^{-}(\calC) \times \Alg_{\calO}^{+}(\calC) 
\simeq \Fun^{\lax}_{ {\calD}^{\otimes}}( \calO^{\otimes} \boxplus
\calO^{\otimes}, \calC^{\otimes}) \rightarrow \Alg^{\pm}_{\calO}(\calC)$
described in Construction \ref{urmouth} is induced by operadic $q$-left Kan extension along
the correspondence of $\infty$-operads
$\calT_{\calO} \rightarrow \Delta^1 \times \Nerve(\FinSeg)$.
\end{proposition}

\begin{proof}
Let
$\calA$ denote the full subcategory of $\Fun_{ {\calD}^{\otimes}}( \calT_{\calO}, \calC^{\otimes})$ spanned by the $\infty$-operad family maps which are operadic $q$-left Kan extensions. The assertion of Proposition \ref{alug} can be stated more precisely as follows:
\begin{itemize}
\item[$(a)$] Every algebra object $A_0 \in 
 \Fun^{\lax}_{ {\calD}^{\otimes}}( \calO^{\otimes} \boxplus
\calO^{\otimes}, \calC^{\otimes})$ admits an operadic
$q$-left Kan extension $A \in \calA$.

\item[$(b)$] The restriction map $\calA \rightarrow  \Fun^{\lax}_{ {\calD}^{\otimes}}( \calO_{-}^{\otimes} \boxplus
\calO_{+}^{\otimes}, \calC^{\otimes})$ is a trivial Kan fibration, and therefore admits a section $s$.

\item[$(c)$] The composition
$$  \Fun^{\lax}_{ {\calD}^{\otimes}}( \calO_{-}^{\otimes} \boxplus
\calO_{+}^{\otimes}, \calC^{\otimes}) \stackrel{s}{\rightarrow} \calA \rightarrow \Alg^{\pm}_{\calO}(\calC)$$
is equivalent to the tensor product functor described in Construction \ref{urmouth}.
\end{itemize}

To prove assertion $(a)$, we must show that for every object $X \in \calO$,
the induced diagram
$$ (\calT_{\calO}^{\acti})_{/X} \times_{ \calT_{\calO}} ( \calO_{-}^{\otimes} \boxplus \calO_{+}^{\otimes})
\rightarrow \calC^{\otimes}$$ can be extended to an operadic $q$-colimit diagram lying over the obvious map
$$ (\calT^{\acti}_{\calO})_{/X} \times_{ \calT} ( \calO^{\otimes} \boxplus \calO^{\otimes}))^{\triangleright}
\rightarrow {\calD}^{\otimes}$$
(Theorem \symmetricref{oplk}). We observe that the $\infty$-category
$(\calT^{\acti}_{\calO})_{/X} \times_{ \calT_{\calO}} ( \calO_{-}^{\otimes} \boxplus \calO_{+}^{\otimes}) $
has a final object, given by the quintuple
$X': ( \seg{2}, \{ 0, \ast \}, \{ 1, \ast \}, X, X)$. It therefore suffices to show that $A_0(X')$ can be extended to 
an operadic $q$-colimit diagram $\{ X' \}^{\triangleright} \rightarrow \calC^{\otimes}$ lying
over the active map $\seg{2} \rightarrow \seg{1}$. The existence of this extension follows from
our assumption that $q$ is a coCartesian fibration of $\infty$-operads (Proposition \symmetricref{justica}). Assertion $(b)$ follows immediately from $(a)$, by virtue of
Lemma \symmetricref{siwd}.

To prove $(c)$, consider the functor $H: \calO^{\otimes} \times \Delta^1 \rightarrow \calT_{\calO}$ of Construction \ref{urmouth}. Composition with $H$ induces a map of simplicial sets $h: \calA \times \Delta^1 \rightarrow \Fun( \calO^{\otimes}, \calC^{\otimes})$. The proof of $(a)$ shows that $h$ can be regarded
as a coCartesian natural transformation from the composite functor 
$$h_0: \calA \rightarrow \Fun^{\lax}_{ {\calD}^{\otimes}}( \calO^{\otimes} \boxplus
\calO^{\otimes}, \calC^{\otimes}) \stackrel{\circ H_0}{\rightarrow} \Fun( \calO^{\otimes}, \calC^{\otimes} )$$
to a functor $h_1: \calA \rightarrow \Alg^{\pm}_{\calO}(\calC)$, so that $h_1$ can be identified with the tensor product of the forgetful functors $\calA \rightarrow \Alg^{-}_{ \calO}( \calC)$ and $\Alg^{+}_{ \calO}(\calC)$. Composing this identification with the section $s$, we obtain the desired result.
\end{proof}

\subsection{Slicing $\infty$-Operads}\label{cluper}

Let $\calC$ be a symmetric monoidal category and let $A$ be a commutative algebra object of $\calC$.
Then the overcategory $\calC_{/A}$ inherits the structure of a symmetric monoidal category: the tensor product of a map $X \rightarrow A$ with a map $Y \rightarrow A$ is given by the composition
$$ X \otimes Y \rightarrow A \otimes A \stackrel{m}{\rightarrow} A,$$
where $m$ denotes the multiplication on $A$. Our goal in this section is to establish an $\infty$-categorical analogue of this observation (and a weaker result concerning undercategories). Before we can state our result, we need to introduce a bit of notation.

\begin{definition}\label{staid}
Let $q: X \rightarrow S$ be a map of simplicial sets, and suppose we are given a commutative diagram
$$ \xymatrix{ & X \ar[d]^{q} \\
S \times K \ar[ur]^{p} \ar[r] & S. }$$
We define a simplicial set $X_{p_S/}$ equipped with a map $q': X_{p_S/} \rightarrow S$
so that the following universal property is satisfied: for every map of simplicial sets
$Y \rightarrow S$, there is a canonical bijection of $\Fun_{S}( Y, X_{p_S/})$ with
the collection of commutative diagrams
$$\xymatrix{ Y \times K \ar[r] \ar[d] & Y \times K^{\triangleright} \ar[d] \ar[r] & Y \ar[d] \\
S \times K \ar[r]^{p} & X \ar[r] & S. }$$
Similarly, we define a map of simplicial sets $X_{/ p_S} \rightarrow S$ so
that $\Fun_{S}(Y, X_{/ p_S})$ is in bijection with the set of diagrams
$$\xymatrix{ Y \times K \ar[r] \ar[d] & Y \times K^{\triangleleft} \ar[d] \ar[r] & Y \ar[d] \\
S \times K \ar[r]^{p} & X \ar[r] & S. }$$
\end{definition}

\begin{remark}
If $S$ consists of a single point, then $X_{p_S/}$ and $X_{/p_S}$ coincide with the usual
overcategory and undercategory constructions $X_{p/}$ and $X_{/p}$. In general,
the fiber of the morphism $X_{p_S/} \rightarrow S$ over a vertex $s \in S$ can be identified with
$(X_s)_{p_s/}$, where $X_{s} = X \times_{S} \{s\}$ and $p_s: K \rightarrow X_s$ is the induced map;
similarly, we can identify $X_{/p_S} \times_{S} \{s\}$ with $(X_{s})_{/ p_s}$.
\end{remark}

\begin{notation}\label{kurplex}
Let $q: \calC^{\otimes} \rightarrow \calO^{\otimes}$ be a fibration of $\infty$-operads, and let
$p: K \rightarrow \Alg_{\calO}(\calC)$ be a diagram. We let $\calC^{\otimes}_{p_{\calO}/}$
and $\calC^{\otimes}_{/ p_{\calO}}$ denote the simplicial sets $(\calC^{\otimes})_{p_{\calO^{\otimes}}/}$ 
and $(\calC^{\otimes})_{/ p_{\calO^{\otimes}}}$ described in Definition \ref{staid}. 

In the special case where $K = \Delta^0$, the diagram $p$ is simply given by a $\calO$-algebra object
$A \in \Alg_{\calO}(\calC)$; in this case, we will denote $\calC^{\otimes}_{p_{\calO}/}$ and
$\calC^{\otimes}_{/ p_{\calO}}$ by $\calC^{\otimes}_{A_{\calO}/}$ and $\calC^{\otimes}_{/A_{\calO}}$, respectively.
\end{notation}

We can now state the main result of this section.

\begin{theorem}\label{eli}
Let $q: \calC^{\otimes} \rightarrow \calO^{\otimes}$ be a fibration of $\infty$-operads,
and let $p: K \rightarrow \Alg_{\calO}(\calC)$ be a diagram. Then:
\begin{itemize}
\item[$(1)$] The maps $\calC^{\otimes}_{p_{\calO}/} \stackrel{q'}{\rightarrow} \calO^{\otimes}
\stackrel{q''}{\leftarrow} \calC^{\otimes}_{/p_{\calO}}$ are fibrations of $\infty$-operads.
\item[$(2)$] A morphism in $\calC^{\otimes}_{p_{\calO}/}$ is inert if and only if its
image in $\calC^{\otimes}$ is inert; similarly, a morphism in $\calC^{\otimes}_{/p_{\calO}}$ is inert if and only if its image in $\calC^{\otimes}$ is inert.
\item[$(3)$] If $q$ is a coCartesian fibration of $\infty$-operads, then $q''$ is a coCartesian fibration
of $\infty$-operads. If, in addition, $p(k): \calO^{\otimes} \rightarrow \calC^{\otimes}$ is a $\calO$-monoidal functor for each vertex $k \in K$, then $q'$ is also a coCartesian fibration of $\infty$-operads.
\end{itemize}
\end{theorem}

The remainder of this section is devoted to the proof of Theorem \ref{eli}. We will need a few lemmas.

\begin{lemma}\label{staffly}
Suppose we are given a diagram of simplicial sets
$$ \xymatrix{ & X \ar[d]^{q} \\
S \times K \ar[ur]^{p} \ar[r] & S }$$
where $q$ is an inner fibration, and let $q': X_{p_S/} \rightarrow S$ be the induced map. 
Then $q'$ is a inner fibration. Similarly, if $q$ is a categorical fibration, then $q'$ is a categorical fibration.
\end{lemma}

\begin{proof}
We will prove the assertion regarding inner fibrations; the case of categorical fibrations is handled similarly. We wish to show that every lifting problem of the form
$$ \xymatrix{ A \ar[r] \ar[d]^{j} & X_{p_S/} \ar[d]^{q'} \\
B \ar[r] \ar@{-->}[ur] & S}$$
admits a solution, provided that $j$ is inner anodyne. Unwinding the definitions, we arrive at an equivalent lifting problem
$$ \xymatrix{ (A \times K^{\triangleright}) \coprod_{A \times K} (B \times K) \ar[r] \ar[d]^{j'} & X \ar[d]^{q} \\
B \times K^{\triangleright} \ar[r] & S, }$$
which admits a solution by virtue of the fact that $q$ is an inner fibration and $j'$ is
inner anodyne (Corollary \toposref{prodprod2}).
\end{proof}

\begin{lemma}\label{cuffly}
Let $q: X \rightarrow S$ be an innert fibration of simplicial sets and let
$K$ and $Y$ be simplicial sets. Suppose that $\overline{h}: K \times Y^{\triangleright} \rightarrow X$
is a map such that, for each $k \in K$, the induced map $\{k\} \times Y^{\triangleright} \rightarrow X$ is a $q$-colimit diagram. Let $h = \overline{h} | K \times Y$. Then the map
$$ X_{\overline{h}/} \rightarrow X_{h/} \times_{ S_{qh/}} S_{q \overline{h}/}$$ is a trivial Kan fibration.
\end{lemma}

\begin{proof}
We will prove more generally that if $K_0 \subseteq K$ is a simplicial subset and
$\overline{h}_0 = \overline{h} | ( K \times Y) \coprod_{ K_0 \times Y} (K \times Y^{\triangleright})$, then
the induced map
$\theta: X_{\overline{h}/} \rightarrow X_{\overline{h}_0/} \times_{ S_{q \overline{h}_0/}} S_{q \overline{h}/}$
is a trivial Kan fibration. Working simplex-by-simplex, we can reduce to the case where
$K = \bd \Delta^n$ and $K' = \bd \Delta^n$. Let us identify 
$Y \star \Delta^n$ with the full simplicial subset of $Y^{\triangleright} \times \Delta^n$ spanned
by $\Delta^n$ and $Y^{\triangleright} \times \{0\}$. Let $\overline{g} = \overline{h} | Y \star \Delta^n$, and
let $g = \overline{g} | Y \star \bd \Delta^n$. Then $\theta$ is a pullback of the map
$$ \theta': X_{ \overline{g}/} \rightarrow X_{g/} \times_{ S_{qg/}} S_{q \overline{g}/}.$$
It will now suffice to show that $\theta'$ has the right lifting property with respect to every
inclusion $\bd \Delta^m \subseteq \Delta^m$. Unwinding the definition, this is equivalent to solving a lifting problem of the form
$$ \xymatrix{ Y \star \bd \Delta^{n+m+1} \ar[r] \ar[d] & X \ar[d]^{q} \\
Y \star \Delta^{n+m+1} \ar[r] \ar@{-->}[ur] & S. }$$
This lifting problem admits a solution by virtue of our assumption that $\overline{h} | \{0\} \times Y^{\triangleright}$ is a $q$-colimit diagram.
\end{proof}

\begin{lemma}\label{stuffly}
Let 
$$ \xymatrix{ & X \ar[d]^{q} \\
S \times K \ar[ur]^{p} \ar[r] & S }$$
be a diagram of simplicial sets, where $q$ is an inner fibration, 
let $q': X^{p_S/} \rightarrow S$ be the induced map, and suppose we are given a commutative diagram
$$ \xymatrix{ Y \ar[r]^{f} \ar[d] & X^{p_S/} \ar[d]^{q'} \\
Y^{\triangleright} \ar[r]^{g} \ar@{-->}[ur]^{\overline{f}} & S }$$
satisfying the following conditions:
\begin{itemize}
\item[$(i)$] For each vertex $k \in K$, the diagram
$$ Y^{\triangleright} \stackrel{g}{\rightarrow} S \simeq S \times \{k\}
\hookrightarrow S \times K \stackrel{p}{\rightarrow} X$$
is a $q$-colimit diagram.
\item[$(ii)$] The composite map $Y \stackrel{f}{\rightarrow} X^{p_S/} \rightarrow X$
can be extended to a $q$-colimit diagram $Y^{\triangleright} \rightarrow X$ lying over $g$.
\end{itemize}
Then:
\begin{itemize}
\item[$(1)$] Let $\overline{f}: Y^{\triangleright} \rightarrow X^{p_S/}$ be a map rendering
the diagram commutative. Then $\overline{f}$ is a $q'$-colimit diagram if and only if the composite map
$Y^{\triangleright} \stackrel{ \overline{f}}{\rightarrow} X^{p_S/} \rightarrow X$ is a
$q$-colimit diagram.
\item[$(2)$] There exists a map $\overline{f}$ satisfying the equivalent conditions of $(1)$.
\end{itemize}
\end{lemma}

\begin{proof}
Let $Z$ be the full simplicial subset of $K^{\triangleright} \times Y^{\triangleright}$ obtained by removing the final object, so we have a canonical isomorphism $Z^{\triangleright} \simeq K^{\triangleright} \times Y^{\triangleright}$. The maps $f$ and $g$ determine a diagram $h: Z \rightarrow X$.
We claim that $h$ can be extended to a $q$-colimit diagram $\overline{h}: Z^{\triangleright} \rightarrow X$
lying over the map $Z^{\triangleright} \rightarrow Y^{\triangleright} \stackrel{g}{\rightarrow} S$.
To prove this, let $h_0 = h | K^{\triangleright} \times Y$, $h_1 = h | K \times Y^{\triangleright}$, and
$h_2 = h | K \times Y$. Using $(i)$ we deduce that the map 
$\theta: X_{h_1/} \rightarrow X_{h_2/} \times_{ S_{q h_2/ } } S_{q h_1/}$ is a trivial Kan fibration
(Lemma \ref{cuffly}). The map $X_{h/} \rightarrow X_{h_0/} \times_{ S_{q h_0/}} S_{qh/}$ is a pullback
of $\theta$, and therefore also a trivial Kan fibration. Consequently, to show that
$h$ admits a $q$-colimit diagram compatible with $g$, it suffices to show that $h_0$ admits a $q$-colimit diagram compatible with $g$.
Since the inclusion $Y \hookrightarrow K^{\triangleright} \times Y$ is cofinal, this follows
immediately from $(ii)$. This proves the existence of $\overline{h}$: moreover, it shows that an arbitrary extension $\overline{h}$ of $h$ (compatible with $g$) is a $p$-colimit diagram if and only if
it restricts to a $p$-colimit diagram $Y^{\triangleright} \rightarrow Z$.

The map $\overline{h}$ determines an extension $\overline{f}: Y^{\triangleright} \rightarrow
X_{p_S/}$ of $f$. We will show that $\overline{f}$ is a $q'$-colimit diagram. This will prove
the ``if'' direction of $(1)$ and $(2)$; the ``only if'' direction of $(1)$ will then follow from the
uniqueness properties of $q'$-colimit diagrams. 

We wish to show that every lifting problem of the form
$$ \xymatrix{ Y \star \bd \Delta^n \ar[r]^{F} \ar[d] & X_{p_S/} \ar[d]^{q'} \\
Y \star \Delta^n \ar@{-->}[ur] \ar[r] & S}$$
admits a solution, provided that $n > 0$ and $F| Y \star \{0\}$ coincides with $\overline{f}$.
This is equivalent to a lifting problem of the form
$$ \xymatrix{ ((Y \star \bd \Delta^n) \times K^{\triangleright})
\coprod_{ (Y \star \bd \Delta^n) \times K} ((Y \star \Delta^n) \times K) \ar[r] \ar[d]^{j} & X \ar[d]^{q} \\
(Y \star \Delta^n) \times K^{\triangleright} \ar[r] \ar@{-->}[ur] & S}$$
It now suffices to observe that the map $j$ is a pushout of the inclusion
$Z \star \bd \Delta^n \hookrightarrow Z \star \Delta^n$, so the desired lifting problem can
be solved by virtue of our assumption that $\overline{h}$ is a $q$-colimit diagram.
\end{proof}

The following result is formally similar to Lemma \ref{stuffly} but requires a slightly different proof:

\begin{lemma}\label{spad}
Let 
$$ \xymatrix{ & X \ar[d]^{q} \\
S \times K \ar[ur]^{p} \ar[r] & S }$$
be a diagram of simplicial sets, where $q$ is an inner fibration, 
let $q': X^{p_S/} \rightarrow S$ be the induced map, and suppose we are given a commutative diagram
$$ \xymatrix{ Y \ar[r]^{f} \ar[d] & X^{p_S/} \ar[d]^{q'} \\
Y^{\triangleleft} \ar[r]^{g} \ar@{-->}[ur]^{\overline{f}} & S }$$
satisfying the following condition:
\begin{itemize}
\item[$(\ast)$] The composite map $Y \stackrel{f}{\rightarrow} X^{p_S/} \rightarrow X$
can be extended to a $q$-limit diagram $g': Y^{\triangleleft} \rightarrow X$ lying over $g$.
\end{itemize}
Then:
\begin{itemize}
\item[$(1)$] Let $\overline{f}: Y^{\triangleleft} \rightarrow X^{p_S/}$ be a map rendering
the diagram commutative. Then $\overline{f}$ is a $q'$-limit diagram if and only if the composite map
$Y^{\triangleleft} \stackrel{ \overline{f}}{\rightarrow} X^{p_S/} \rightarrow X$ is a
$q$-limit diagram.
\item[$(2)$] There exists a map $\overline{f}$ satisfying the equivalent conditions of $(1)$.
\end{itemize}
\end{lemma}

\begin{proof}
Let $v$ be the cone point of $K^{\triangleright}$ and $v'$ the cone point of $Y^{\triangleleft}$.
Let $Z$ be the full subcategory of $K^{\triangleright} \times Y^{\triangleleft}$ obtained
by removing the vertex $(v,v')$. The maps $f$ and $g$ determine a map
$h: Z \rightarrow X$. Choose any map $g'$ as in $(i)$, and let $g'_0 = g' | Y$. We claim that
there exists an extension $\overline{h}: K^{\triangleright} \times Y^{\triangleleft} \rightarrow X$
of $h$ which is compatible with $g$, such that $\overline{h} | \{v\} \times Y^{\triangleleft} = g'$.
Unwinding the definitions, we see that providing such a map $\overline{h}$ is equivalent to
solving a lifting problem of the form
$$ \xymatrix{ \emptyset \ar[r] \ar[d] & X_{/g'} \ar[d] \\
K \ar[r] & X_{g'_0} \times_{ S_{/ q g'_0} } S_{/g}, }$$
which is possible since the left vertical map is a trivial Kan fibration (since $g'$ is a $q$-limit diagram).

The map $\overline{h}$ determines a diagram $\overline{f}: Y^{\triangleleft} \rightarrow X_{p_S/}$.
We will prove that $\overline{f}$ is a $q'$-limit diagram. This will prove the ``if'' direction of $(1)$ and $(2)$; the ``only if'' direction of $(1)$ will then follow from the uniqueness properties of $q$-limit diagrams.

To show that $\overline{f}$ is a $q$-limit diagram, we must show that every lifting problem of the form
$$ \xymatrix{ \bd \Delta^n \star Y \ar[r]^{F} \ar[d] & X_{p_S/} \ar[d]^{q'} \\
\Delta^n \star Y \ar@{-->}[ur] \ar[r] & S}$$
admits a solution, provided that $n > 0$ and $F | \{n\} \star Y = \overline{f}$. Unwinding the definitions,
we obtain an equivalent lifting problem
$$ \xymatrix{ ( \bd \Delta^n \star Y) \times K^{\triangleright})
\coprod_{ (\bd \Delta^n \star Y) \times K} ((\Delta^n \star Y) \times K) \ar[r] \ar[d]^{j} & X \ar[d]^{q} \\
(\Delta^n \star Y) \times K^{\triangleright} \ar[r] \ar@{-->}[ur] & S.}$$
It now suffices to observe that $j$ is a pushout of the inclusion 
$K \star \bd \Delta^n \star Y
\hookrightarrow K \star \Delta^n \star Y$, so that the desired extension exists because
$\overline{h}| \{v\} \times Y^{\triangleleft} = g'$ is a $q$-limit diagram.
\end{proof}

\begin{proof}[Proof of Theorem \ref{eli}]
We will prove $(1)$, $(2)$, and $(3)$ for the simplicial set $\calC^{\otimes}_{p_{\calO}/}$; the analogous
assertions for $\calC^{\otimes}_{/ p_{\calO}}$ follow by the same reasoning.
We first observe that $q'$ is a categorical fibration (Lemma \ref{staffly}). Let
$\overline{X} \in \calC^{\otimes}_{p_{\calO}/}$, and suppose we are given an inert morphism
$\alpha: q(\overline{X}) \rightarrow Y$ in $\calO^{\otimes}$; we wish to show that there exists
a $q'$-coCartesian morphism $\overline{X} \rightarrow \overline{Y}$ in $\calC^{\otimes}_{p_{\calO}/}$ lifting $\alpha$. This follows immediately from Lemma \ref{stuffly}.

Suppose next that we are given an object $X \in \calO^{\otimes}$ lying over $\seg{n} \in \FinSeg$, and a collection of inert morphisms $\alpha^{i}: X \rightarrow X_i$ lying over $\Colp{i}: \seg{n} \rightarrow \seg{1}$ for $1 \leq i \leq n$. We wish to prove that the maps $\alpha^{i}$ induce an equivalence
$$ \theta: (\calC^{\otimes}_{p_{\calO}/})_{X} \simeq \prod_{1 \leq i \leq n} (\calC^{\otimes}_{p_{\calO}/})_{X_i}.$$
Let $p_X: K \rightarrow \calC^{\otimes}_{X}$ be the map induced by $p$, and define
maps $p_{X_i}: K \rightarrow \calC_{X_i}$ similarly. We observe that $p_{X_i}$ can be identified with the composition of $p_X$ with $\alpha^{i}_{!}: \calC^{\otimes}_{X} \rightarrow \calC_{X_i}$. Since
$q$ is a fibration of $\infty$-operads, we have an equivalence of $\infty$-categories
$$ \calC^{\otimes}_{X} \rightarrow \prod_{1 \leq i \leq n} \calC_{X_i}.$$
Passing to the $\infty$-categories of objects under $p$, we deduce that $\theta$ is also an
equivalence.

Now suppose that $X$ is as above, that $\overline{X} \in \calC^{\otimes}_{p_{\calO}/}$ is a preimage
of $X$, and that we are given $q'$-coCartesian morphisms $\overline{X} \rightarrow \overline{X}_i$ lying over the maps $\alpha^{i}$. We wish to show that the induced map $\delta: \nostar{n}^{\triangleleft} \calC^{\otimes}_{p_{\calO}/}$ is a $q'$-limit diagram. This follows from Lemma \ref{spad}, since the image of $\delta$ in $\calC^{\otimes}$ is a $q$-limit diagram. This completes the proof of $(1)$.
Moreover, our characterization of $q'$-coCartesian morphisms immediately implies $(2)$.
Assertion $(3)$ follows immediately from Lemma \ref{stuffly}.
\end{proof}

\subsection{Wreath Products of $\infty$-Operads}\label{wreath}

\begin{definition}\label{ahmad}
Suppose we are given a bifunctor $F: \calO^{\otimes} \times {\calO'}^{\otimes}
\rightarrow {\calO''}^{\otimes}$ of $\infty$-operads. Let $\calO^{\otimes, \natural}$
denote the marked simplicial set $(\calO^{\otimes}, M)$, where $M$ is the collection of all inert
morphisms in $\calO^{\otimes}$, and define ${\calO'}^{\otimes, \natural}$ and
${\calO''}^{\otimes, \natural}$ similarly. We will say that $F$ {\it exhibits
${\calO''}^{\otimes}$ as a tensor product of the $\infty$-operads $\calO^{\otimes}$
and ${\calO'}^{\otimes}$} if the underlying map
$\calO^{\otimes, \natural} \odot {\calO'}^{\otimes, \natural} \rightarrow
{\calO''}^{\otimes, \natural}$ is a weak equivalence of $\infty$-preoperads (with respect to the model structure of Proposition \symmetricref{cannwell}).
\end{definition}

\begin{remark}
In other words, a bifunctor of $\infty$-operads $F: \calO^{\otimes} \times {\calO'}^{\otimes}
\rightarrow {\calO''}^{\otimes}$ exhibits ${\calO''}^{\otimes}$ as a tensor product of
$\calO^{\otimes}$ and ${\calO'}^{\otimes}$ if and only if, for every $\infty$-operad
$\calC^{\otimes}$, composition with $F$ induces an equivalence of $\infty$-categories
$\Alg_{ \calO''}(\calC) \rightarrow \Alg_{ \calO}( \Alg_{\calO'}(\calC) )$. 
\end{remark}

For every pair of $\infty$-operads $\calO^{\otimes}$ and ${\calO'}^{\otimes}$, there exists a bifunctor
$F: \calO^{\otimes} \times {\calO'}^{\otimes}
\rightarrow {\calO''}^{\otimes}$ which exhibits ${\calO''}^{\otimes}$ as a tensor product of the
$\infty$-operads $\calO^{\otimes}$ and ${\calO'}^{\otimes}$. Moreover, the $\infty$-operad
${\calO''}^{\otimes}$ (and the bifunctor $F$) are determined up to equivalence. However, it can be quite difficult to describe ${\calO''}^{\otimes}$ directly. The product $\calO^{\otimes,\natural} \odot {\calO'}^{\otimes, \natural}$
is essentially {\em never} a fibrant $\infty$-preoperad, and the process of fibrant replacement is
fairly inexplicit. Our goal in this section is to partially address this difficulty by introducing a map of
$\infty$-preoperads $\calO^{\otimes, \natural} \times {\calO'}^{\otimes, \natural} \rightarrow
( \calO^{\otimes} \wreath {\calO'}^{\otimes}, M)$. Our main result, Theorem \ref{kuj}, asserts that this map is a weak equivalence. This is not really a complete answer, since the codomain $(\calO^{\otimes} \wreath {\calO'}^{\otimes},M)$ is still generally not a fibrant $\infty$-preoperad. However, it is in many ways more convenient than the Cartesian product $\calO^{\otimes, \natural} \odot {\calO'}^{\otimes, \natural}$, and will play an important technical role in analyzing the tensor products of little cubes $\infty$-operads in \S \ref{sass1}. 

\begin{construction}
If $\calC$ is an $\infty$-category, we let $\calC^{\amalg}$ be defined as in Construction
\symmetricref{jiu}. Note that if $\calC$ is the nerve of a category $\calJ$, then
$\calC^{\amalg}$ can be identified with the nerve of the category
$\calJ^{\amalg}$ defined as follows:
\begin{itemize}
\item[$(i)$] The objects of $\calJ^{\amalg}$ are finite sequences $(J_1, \ldots, J_n)$ of
objects in $\calJ$.
\item[$(ii)$] A morphism from $(I_1, \ldots, I_m)$ to $(J_1, \ldots, J_n)$ in
$\calJ^{\amalg}$ consists of a map $\alpha: \seg{m} \rightarrow \seg{n}$ in
$\FinSeg$ together with a collection of maps $\{ I_{i} \rightarrow J_{j} \}_{ \alpha(i) = j }$.
\end{itemize}

There is an evident functor $\FinSeg^{\amalg} \rightarrow \FinSeg$, given on objects by the formula
$( \seg{k_1}, \ldots, \seg{k_n} ) \mapsto \seg{ k_1 + \cdots + k_n }$. This functor induces a map
$$\Phi: \Nerve(\FinSeg)^{\amalg} \rightarrow \Nerve(\FinSeg).$$

Let $\calC^{\otimes}$ and $\calD^{\otimes}$ be $\infty$-operads. We let
$\calC^{\otimes} \wreath \calD^{\otimes}$ denote the simplicial set
$$ \calC^{\otimes} \times_{ \Nerve(\FinSeg)} ( \calD^{\otimes})^{\amalg}.$$
We define a map of simplicial sets $\pi: \calC^{\otimes} \wreath \calD^{\otimes} \rightarrow \Nerve(\FinSeg)$ by considering the composition
\begin{eqnarray*}
\calC^{\otimes} \wreath \calD^{\otimes} & = & \calC^{\otimes} \times_{ \Nerve(\FinSeg)} (\calD^{\otimes})^{\amalg} \\
& \rightarrow & (\calD^{\otimes})^{\amalg} \\
& \rightarrow & ( \Nerve(\FinSeg))^{\amalg} \\
& \stackrel{\Phi}{\rightarrow} & \Nerve(\FinSeg).
\end{eqnarray*}

We can identify a morphism $f$ in $\calC^{\otimes} \wreath \calD^{\otimes}$ with
a map $g: ( D_1, \ldots, D_m) \rightarrow (D'_1, \ldots, D'_n)$ in $(\calD^{\otimes})^{\amalg}$
lying over $\alpha: \seg{m} \rightarrow \seg{n}$ in $\Nerve(\FinSeg)$, together with a
map $h: C \rightarrow C'$ in $\calC^{\otimes}$ lying over $\alpha$. We will say that
$f$ is {\it inert} if $h$ is an inert morphism in $\calD^{\otimes}$ and
$g$ determines a set of inert morphisms $\{ D_{i} \rightarrow D'_{j} \}_{ \alpha(i) = j }$
in $\calD^{\otimes}$. Note that the map $\pi$ carries inert morphisms of
$\calC^{\otimes} \wreath \calD^{\otimes}$ to inert morphisms in $\Nerve(\FinSeg)$.
\end{construction}

\begin{remark}\label{hulker}
Let $\calC^{\otimes}$ and $\calD^{\otimes}$ be $\infty$-operads. 
The map $\calD^{\otimes} \times \Nerve(\FinSeg) \rightarrow (\calD^{\otimes})^{\amalg}$
of Example \symmetricref{stuly} induces a monomorphism of simplicial sets
$\calC^{\otimes} \times \calD^{\otimes} \rightarrow \calC^{\otimes} \wreath \calD^{\otimes}$.
\end{remark}

\begin{theorem}\label{kuj}
Let $\calC^{\otimes}$ and $\calD^{\otimes}$ be $\infty$-operads, and let
$M$ be the collection of inert morphisms in $\calC^{\otimes} \wreath \calD^{\otimes}$.
Then the inclusion $\calC^{\otimes} \times \calD^{\otimes} \rightarrow \calC^{\otimes} \wreath \calD^{\otimes}$ of Remark \ref{hulker} induces a weak equivalence of $\infty$-preoperads
$$ \calC^{\otimes, \natural} \odot \calD^{\otimes, \natural} \rightarrow ( \calC^{\otimes} \wreath \calD^{\otimes}, M).$$
\end{theorem}

\begin{proof}
The proof is nearly identical to that of Theorem \symmetricref{kinema}. Using Proposition
\toposref{minimod}, we may assume without loss of generality that $\calD^{\otimes}$ is minimal.
Let $\calE^{\otimes}$ be an $\infty$-operad. We let
$\calX$ denote the full subcategory of
$\Fun_{ \Nerve(\FinSeg)}( \calC^{\otimes} \wreath \calD^{\otimes}, \calE^{\otimes})$ spanned
functors $F$ which carry inert morphisms in $\calC^{\otimes} \wreath \calD^{\otimes}$ to inert morphisms in $\calE^{\otimes}$, and define $\calY \subseteq \Fun_{ \Nerve(\FinSeg)}( \calC^{\otimes} \times \calD^{\otimes}, \calE^{\otimes})$ similarly. We will show that the restriction functor
$\calX \rightarrow \calY$ is a trivial Kan fibration.

We now introduce a bit of terminology.
Recall that a morphism $\overline{\alpha}$ from $(D_1, \ldots, D_m)$ to
$(D'_1, \ldots, D'_n)$ in $(\calD^{\otimes})^{\amalg}$ consists of a map of pointed sets
$\alpha: \seg{m} \rightarrow \seg{n}$ together with a morphism $f_i: D_i \rightarrow
D'_{\alpha(i)}$ for each $i \in \alpha^{-1} \nostar{n}$. We will say that $\overline{\alpha}$ is
{\it quasidegenerate} if each of the morphisms $f_i$ is a degenerate edge of $\calD^{\otimes}$.
Let $\sigma$ be an $n$-simplex of $(\calD^{\otimes})^{\amalg}$ given by a sequence of morphisms
$$ \sigma(0) \stackrel{ \overline{\alpha}(1)}{\rightarrow} \sigma(1)
\rightarrow \cdots \stackrel{ \overline{\alpha}(n)}{\rightarrow} \sigma(n)$$
and let
$$ \seg{k_0} \stackrel{\alpha(1)}{\rightarrow} \seg{k_1} \rightarrow \cdots \stackrel{\alpha(n)}{\rightarrow} \seg{k_n}$$
be the underlying $n$-simplex of $\Nerve(\FinSeg)$. We will say that $\sigma$ is
{\it closed} if $k_n = 1$, and {\it open} otherwise. If $\sigma$ is closed, we define the
{\it tail length} of $\sigma$ to be the largest integer $m$ such that the maps 
$\alpha(k)$ are isomorphisms for $n-m < k \leq n$. We will denote the tail length of $\sigma$ by $t(\sigma)$. We define the {\it break point} of a closed simplex $\sigma$ to be smallest nonnegative integer $m$ such that the maps $\overline{\alpha}(k)$ are
active and quasidegenerate for $m < k \leq n- t(\sigma)$. We will denote the break point of $\sigma$ by
$b(\sigma)$. Let $S = \coprod_{0 \leq i \leq n} \nostar{k_i}$. We will say that
an element $j \in \nostar{k_i} \subseteq S$ is a {\it leaf} if $i=0$ or if $j$ does not lie in the image
of the map $\alpha(i)$, and we will say that $j$ is a {\it root} if $i=n$ or if $\alpha(i+1)(j) = \ast$.
We define the {\it complexity} $c(\sigma)$ of $\sigma$ to be $2l - r$, where $l$ is the number of
leaves of $\sigma$ and $r$ is the number of roots of $\sigma$. We will say that $\sigma$
is {\it flat} if it belongs to the image of the embedding $\Nerve(\FinSeg) \times \calD^{\otimes} \rightarrow (\calD^{\otimes})^{\amalg}$. Since $\calD^{\otimes}$ is minimal, Proposition \toposref{minstrict} implies that if
$\sigma$ is closed and $b(\sigma) = 0$, then $\sigma$ is flat.

We now partition the nondegenerate, nonflat simplices of $(\calD^{\otimes})^{\amalg}$ into six groups:

\begin{itemize}
\item[$(A)$] An $n$-dimensional nonflat nondegenerate simplex $\sigma$ of
$(\calD^{\otimes})^{\amalg}$ belongs to $A$ if $\sigma$ is closed and
the map $\alpha(b(\sigma))$ is not inert.

\item[$(A')$] An $n$-dimensional nonflat nondegenerate simplex $\sigma$ of
$(\calD^{\otimes})^{\amalg}$ belongs to $A'$ if $\sigma$ is closed, $b(\sigma) < n - t(\sigma)$,
and the map $\alpha(b(\sigma))$ is inert.

\item[$(B)$] An $n$-dimensional nonflat nondegenerate simplex $\sigma$ of
$(\calD^{\otimes})^{\amalg}$ belongs to $B$ if $\sigma$ is closed, $b(\sigma) = n - t(\sigma)$,
the map $\alpha( b(\sigma) )$ is inert, and $\overline{\alpha}( b(\sigma) )$ is 
not quasidegenerate.

\item[$(B')$] An $n$-dimensional nonflat nondegenerate simplex $\sigma$ of
$(\calD^{\otimes})^{\amalg}$ belongs to $B$ if $\sigma$ is closed, $b(\sigma) = n - t(\sigma) < n$,
the map $\alpha(b(\sigma))$ is inert, and $\overline{\alpha}( b(\sigma))$ is quasidegenerate.

\item[$(C)$] An $n$-dimensional nonflat nondegenerate simplex $\sigma$ of
$(\calD^{\otimes})^{\amalg}$ belongs to $C$ if it is open.

\item[$(C')$] An $n$-dimensional nonflat nondegenerate simplex $\sigma$ of
$(\calD^{\otimes})^{\amalg}$ belongs to $C'$ is it is closed, $b(\sigma) = n-t(\sigma) = n$, the
map $\alpha(b(\sigma))$ is inert, and $\overline{\alpha}(b(\sigma))$ is quasidegenerate.
\end{itemize}

If $\sigma$ belongs to $A'$, $B'$, or $C'$, then we define the {\it associate} $a(\sigma)$ of
$\sigma$ to be the face of $\sigma$ opposite the $b(\sigma)$th vertex. It follows from
Proposition \toposref{minstrict} that $a(\sigma)$ belongs to $A$ if $\sigma \in A'$,
$B$ if $\sigma \in B'$, and $C$ if $\sigma \in C'$. In this case, we will say that
$\sigma$ is an {\it associate} of $a(\sigma)$. We note that every simplex belonging to
$A$ or $B$ has a unique associate, while a simplex $\sigma$ of $C$ has precisely $k$ associates,
where $\seg{k}$ is the image of the final vertex of $\sigma$ in $\Nerve(\FinSeg)$.

For each $n \geq 0$, let $K(n) \subseteq (\calD^{\otimes})^{\amalg}$ be the simplicial subset
generated by those nondegenerate simplices which are either flat, have dimension $\leq n$, 
or have dimension $n+1$ and belong to either $A'$, $B'$, or $C'$. We observe that 
$K(0)$ is generated by $\calD^{\otimes} \times \Nerve(\FinSeg)$ together with the collection of
$1$-simplices belonging to $C'$. Let $\calX(n)$ denote the full subcategory of
$\bHom_{\Nerve(\FinSeg)}( \calC^{\otimes} \times_{\Nerve(\FinSeg)}
K(n), \calE^{\otimes} )$ spanned by those maps $F$ with the following properties:
\begin{itemize}
\item[$(i)$] The restriction of $F$ to $\calC^{\otimes} \times \calD^{\otimes}$ belongs to
$\calY$.

\item[$(ii)$] Let $f$ be an edge of $\calC^{\otimes} \times_{ \Nerve(\FinSeg)} K(0)$ whose
image in $\calC^{\otimes}$ is inert and whose image in $K(0)$ belongs to $C'$. Then
$F(f)$ is an inert morphism in $\calE^{\otimes}$.
\end{itemize}

To complete the proof, it will suffice to show that the restriction maps
$$ \calX \stackrel{\theta'}{\rightarrow} \calX(0) \stackrel{\theta''}{\rightarrow} \calY$$
are trivial Kan fibrations. For the map $\theta''$, this follows from repeated application
of Lemma \symmetricref{stealthwick}.
To prove that $\theta'$ is a trivial Kan fibration, we define 
$\calX(n)$ to be the full subcategory of 
$\bHom_{\Nerve(\FinSeg)}( \calC^{\otimes} \times_{\Nerve(\FinSeg)}
K(n), \calE^{\otimes} )$ spanned by those functors $F$ whose restriction
to $\calC^{\otimes} \times_{ \Nerve(\FinSeg)} K(0)$ belongs to $\calX(0)$.
We will prove the following:

\begin{itemize}
\item[$(a)$] A functor $F \in \Fun_{ \Nerve(\FinSeg)}( \calC^{\otimes} \wreath \calD^{\otimes}, \calE^{\otimes})$ carries inert morphisms to inert morphisms if
and only if $F$ satisfies conditions $(i)$ and $(ii)$. Consequently, the $\infty$-category
$\calX$ can be identified with the inverse limit of the tower
$$ \cdots \rightarrow \calX(2) \rightarrow \calX(1) \rightarrow \calX(0).$$

\item[$(b)$] For $n > 0$, the restriction map $\calX(n) \rightarrow \calX(n-1)$ is a trivial Kan fibration.
\end{itemize}

We first prove $(a)$. The ``only if'' direction is obvious. For the converse, suppose that
an object $F$ of $\Fun_{ \Nerve(\FinSeg)}( \calC^{\otimes} \wreath \calD^{\otimes}, \calE^{\otimes})$ satisfies conditions $(i)$ and $(ii)$ above. We wish to prove that $F$ preserves inert morphisms. Let
$f: X \rightarrow X'$ be an inert morphism in $\calC^{\otimes} \wreath \calD^{\otimes}$
covering the map $f_0: ( \seg{k_1}, \ldots, \seg{k_m}) \rightarrow
( \seg{k'_1}, \ldots, \seg{k'_{m'}})$ in $\Nerve(\FinSeg)^{\amalg}$; we wish to prove that
$F(f)$ is an inert morphism in $\calE^{\otimes}$. If $m' = k'_{1} = 1$, then $f_0$ factors as a composition of inert morphisms
$$ ( \seg{k_1}, \ldots, \seg{k_m}) \stackrel{f'_0}{\rightarrow} ( \seg{k_i} ) \stackrel{f''_0}{\rightarrow} ( \seg{1} )$$
for some $i \in \nostar{m}$, which we can lift to a factorization $f \simeq f'' \circ f'$ of $f$
where $f'$ is quasidegenerate. Condition $(ii)$ guarantees that $F( f')$ is inert, and condition
$(i)$ guarantees that $F(f'')$ is inert. In the general case, we consider for each
$j \in \nostar{k'_{i}}$ an inert morphism $g_{i,j}: X' \rightarrow X''$ lifting the composite map
$$ ( \seg{k'_1}, \ldots, \seg{k'_{m'}}) \rightarrow ( \seg{k'_{i}}) \rightarrow ( \seg{1} ).$$
The above argument shows that $F( g_{i,j})$ and $F( g_{i,j} \circ f)$ are inert morphisms
in $\calE^{\otimes}$. The argument of Remark \symmetricref{casper} shows that $F(f)$ is inert, as desired.

We now prove $(b)$. For each integer $c \geq 0$, let $K(n,c)$ denote the simplicial subset
$K(n)$ spanned by those simplices which either belong to $K(n-1)$ or have complexity $\leq c$.
Let $\calX(n,c)$ denote the full subcategory of
$\Fun_{\Nerve(\FinSeg)}( \calC^{\otimes} \times_{\Nerve(\FinSeg)}
K(n,c), \calE^{\otimes} )$
spanned by those maps $F$ whose restriction to $K(0)$ satisfies conditions $(i)$ and $(ii)$.
We have a tower of simplicial sets
$$ \cdots \rightarrow \calX(n,2) \rightarrow \calX(n,1) \rightarrow \calX(n,0) \simeq \calX(n-1)$$
with whose inverse limit can be identified with $\calX(n)$. It will therefore suffice to show that for
each $c > 0$, the restriction map $\calX(n,c) \rightarrow \calX(n,c-1)$ is a trivial Kan fibration.

We now further refine our filtration as follows. Let $K(n,c)_{A}$ denote the 
simplicial subset of $K(n,c)$ spanned by $K(n,c-1)$ together with those simplices
of $K(n,c)$ which belong to $A$ or $A'$ and let $K(n,c)_{B}$ denote the simplicial subset
of $K(n,c)$ spanned by $K(n,c-1)$ together with those simplices which belong to
$A$, $A'$, $B$, or $B'$. Let
$\calX(n,c)_{A}$ denote the full subcategory of
$\Fun_{\Nerve(\FinSeg)}( \calC^{\otimes} \times_{ \Nerve(\FinSeg)} K(n,c)_{A}, \calE^{\otimes} )$
spanned by those maps $F$ satisfying $(i)$ and $(ii)$, and define
$\calX(n,c)_{B}$ similarly. To complete the proof, it will suffice to prove the following:

\begin{itemize}
\item[$(A)$] The restriction map $\calX(n,c)_{A} \rightarrow \calX(n,c-1)$ is a trivial Kan fibration.
To prove this, it suffices to show that the inclusion
$\calC^{\otimes} \times_{ \Nerve(\FinSeg)} K(n,c-1) \rightarrow \calC^{\otimes} \times_{ \Nerve(\FinSeg)} K(n,c)_{A}$ is a categorical equivalence.
Let $A_{n,c}$ denote the collection of all $n$-simplices belonging to $A$ having
complexity $c$. Choose a well-ordering of $A_{n,c}$ with the following properties:
\begin{itemize}
\item If $\sigma, \sigma' \in A_{n,c}$ and $t(\sigma) < t(\sigma')$, then $\sigma < \sigma'$.
\item If $\sigma, \sigma' \in A_{n,c}$, $t(\sigma) = t(\sigma')$, and $b(\sigma) < b(\sigma')$, then
$\sigma < \sigma'$.
\end{itemize}
For each $\sigma \in A_{n,c}$, let $K(n,c)_{\leq \sigma}$ denote the simplicial subset
of $K(n,c)$ generated by $K(n,c-1)$, all simplices $\tau \leq \sigma$ in $A_{n,c}$, and
all of the simplices in $A'$ which are associated to simplices of the form $\tau \leq \sigma$. Define $K(n,c)_{< \sigma}$ similarly. Using transfinite induction on $A_{n,c}$, we are reduced to proving
that for each $\sigma \in A_{n,c}$, the inclusion
$$ i: \calC^{\otimes} \times_{ \Nerve(\FinSeg)} K(n,c)_{< \sigma} \rightarrow 
\calC^{\otimes} \times_{ \Nerve(\FinSeg)} K(n,c)_{\leq \sigma}$$
is a categorical equivalence. Let $\sigma': \Delta^{n+1} \rightarrow (\calD^{\otimes})^{\amalg}$ be the
unique $(n+1)$-simplex of $A'$ associated to $\sigma$. We observe that $\sigma'$ determines
a pushout diagram
$$ \xymatrix{ \Lambda^{n+1}_{b(\sigma')} \ar[r] \ar[d] & K(n,c)_{< \sigma} \ar[d] \\
\Delta^{n+1} \ar[r] & K(n,c)_{\leq \sigma}. }$$
Consequently, the map $i$ is a pushout of an inclusion
$$ i': \calC^{\otimes} \times_{ \Nerve(\FinSeg)} \Lambda^{n+1}_{b(\sigma')} 
\rightarrow \calC^{\otimes} \times_{ \Nerve(\FinSeg)} \Delta^{n+1}_{b(\sigma)}.$$
Since the Joyal model structure is left proper, it suffices to show that $i'$ is a
categorical equivalence, which follows from Lemma \symmetricref{skope}.

\item[$(B)$] The map $\calX(n,c)_{B} \rightarrow \calX(n,c)_{A}$ is a trivial Kan fibration.
To prove this, it suffices to show that the inclusion 
$\calC^{\otimes} \times_{ \Nerve(\FinSeg)} K(n,c)_{A}
\subseteq \calC^{\otimes} \times_{ \Nerve(\FinSeg)} K(n,c)_{B}$
is a categorical equivalence of simplicial sets. Let $B_{n,c}$ denote the collection
of all $n$-simplices belonging to $B$ having complexity $c$. Choose
a well-ordering of $B_{n,c}$ such that the function $\sigma \mapsto t(\sigma)$ is nonstrictly
decreasing. For each $\sigma \in B_{n,c}$, we let
$K(n,c)_{\leq \sigma}$ be the simplicial subset of $K(n,c)$ generated by
$K(n,c)_{A}$, those simplices $\tau$ of $B_{n,c}$ such that $\tau \leq \sigma$, and
those simplices of $B'$ which are associated to $\tau \leq \sigma \in B_{n,c}$.
Let $K(n,c)_{< \sigma}$ be defined similarly. Using a induction on $B_{n,c}$, we
can reduce to the problem of showing that each of the inclusions
$$ \calC^{\otimes} \times_{ \Nerve(\FinSeg)} K(n,c)_{< \sigma} \rightarrow
\calC^{\otimes} \times_{ \Nerve(\FinSeg)} K(n,c)_{\leq \sigma}$$
is a categorical equivalence. 
Let $\sigma': \Delta^{n+1} \rightarrow (\calD^{\otimes})^{\amalg}$ be the
unique $(n+1)$-simplex of $B'$ associated to $\sigma$. We observe that $\sigma'$ determines
a pushout diagram
$$ \xymatrix{ \Lambda^{n+1}_{b(\sigma')} \ar[r] \ar[d] & K(n,c)_{< \sigma} \ar[d] \\
\Delta^{n+1} \ar[r] & K(n,c)_{\leq \sigma}. }$$
Consequently, the map $i$ is a pushout of an inclusion
$$ i': \calC^{\otimes} \times_{ \Nerve(\FinSeg)} \Lambda^{n+1}_{b(\sigma')}
\rightarrow \calC^{\otimes} \times_{ \Nerve(\FinSeg)} \Delta^{n+1}_{b(\sigma)}.$$
Since the Joyal model structure is left proper, it suffices to show that $i'$ is a
categorical equivalence, which follows from Lemma \symmetricref{skope}.

\item[$(C)$] The map $\calX(n,c) \rightarrow \calX(n,c)_{B}$ is a trivial
Kan fibration. To prove this, let $C_{n,c}$ denote the subset of
$C$ consisting of $n$-dimensional simplices of complexity $c$, and choose
a well-ordering of $C_{n,c}$. For each $\sigma \in C_{n,c}$, let
$K(n,c)_{\leq \sigma}$ denote the simplicial subset of $K(n,c)$ generated
by $K(n,c)_{B}$, those simplices $\tau \in C_{n,c}$ such that $\tau \leq \sigma$,
and those simplices of $C'$ which are associated to $\tau \in C_{n,c}$ with
$\tau \leq \sigma$. Let $\calX(n,c)_{\leq \sigma}$ be the
full subcategory of $\Fun_{\Nerve(\FinSeg)}( \calC^{\otimes} \times_{ \Nerve(\FinSeg)} K(n,c)_{\leq \sigma}, \calE^{\otimes} )$ spanned by those maps $F$ satisfying $(i)$ and $(ii)$.
We define $K(n,c)_{< \sigma}$ and
$\calX(n,c)_{< \sigma}$ similarly. Using transfinite induction on $C_{n,c}$, we
are reduced to the problem of showing that for each $\sigma \in C_{n,c}$, the map
$\psi: \calX(n,c)_{\leq \sigma} \rightarrow \calX(n,c)_{< \sigma}$ is a trivial Kan fibration.

Let $\seg{k}$ denote the image of the final vertex of $\sigma$ in $\FinSeg$.
For $1 \leq i \leq k$, let $\sigma_i \in C'$ denote the unique $(n+1)$-simplex
associated to $\sigma$ such that $\sigma_{i}$ carries $\Delta^{ \{n, n+1\} }$ to
the morphism $\Colp{i}$ in $\FinSeg$. The simplices
$\{ \sigma_{i} \}_{1 \leq i \leq k}$ determine a map of simplicial sets
$\Delta^n \star \nostar{k} \rightarrow K(n,c)_{\leq \sigma}$. We have a pushout diagram of simplicial sets
$$ \xymatrix{ (\bd \Delta^n) \star \nostar{k} \ar[r] \ar[d] & K(n,c)_{< \sigma} \ar[d] \\
\Delta^n \star \nostar{k} \ar[r] & K(n,c)_{\leq \sigma}. }$$
The map $\psi$ fits into a pullback diagram
$$ \xymatrix{ \calE(n,c)_{\leq \sigma} \ar[r] \ar[d]^{\psi} & \Fun'_{ \Delta^n \star \nostar{k}}( \calC^{\otimes} \times_{ \Nerve(\FinSeg)} (\Delta^n \star \nostar{k}),
\calE^{\otimes} \times_{ \Nerve(\FinSeg)} ( \Delta^n \star \nostar{k})) \ar[d]^{\psi'} \\
\calE(n,c)_{< \sigma} \ar[r] & \Fun'_{ \bd \Delta^n \star \nostar{k}}( \calC^{\otimes} \times_{ \Nerve(\FinSeg)} ( \bd \Delta^n \star \nostar{k}), \calE^{\otimes} \times_{ \Nerve(\FinSeg)} (\bd \Delta^n \star \nostar{k}))}$$
where $\psi'$ denotes the trivial Kan fibration of Lemma \symmetricref{sleuth}.
\end{itemize}
\end{proof}

\end{document}